\documentclass[12pt,reqno]{amsart}
\usepackage{amsfonts,amsrefs,latexsym,amsmath, amssymb, mathrsfs, verbatim, slashed}
\usepackage{url,color}
\usepackage{pifont}
\usepackage{upgreek}
\usepackage{fancyhdr}
\usepackage{hyperref}
\usepackage{calligra}
\usepackage{marvosym}
\usepackage[percent]{overpic}
\usepackage{pict2e}
\usepackage[hypcap=false]{caption}
\usepackage{xcolor}
\usepackage{wasysym}
\usepackage{accents}

\topmargin - .23 in

\newtheorem{theorem}{Theorem}[section]
\newtheorem{proposition}{Proposition}[section]
\newtheorem{lemma}[proposition]{Lemma}
\newtheorem{corollary}[proposition]{Corollary}

\theoremstyle{definition}
\newtheorem{definition}{Definition}[section]
\newtheorem{remark}{Remark}[section]

\theoremstyle{plain}

\DeclareMathAlphabet{\mathcalligra}{T1}{calligra}{m}{n}
\DeclareFontShape{T1}{calligra}{m}{n}{<->s*[2.2]callig15}{}

\newcommand{\Tboot}{T_{(Boot)}}

\newcommand{\threePsi}{\vec{\Psi}}

\newcommand{\contr}{\diamond}

\newcommand{\Trandatasize}[1]{\mathring{\updelta}}

\newcommand{\TranminusdatasizeWithFactor}{\mathring{\updelta}_{\ast}}

\newcommand{\gt}{\underline{g}}

\newcommand{\gsphere}{g \mkern-8.5mu / }
\newcommand{\ginversesphere}{\gsphere^{-1}}

\newcommand{\gtancomp}{\upsilon}

\newcommand{\mytr}{{\mbox{\upshape{tr}}_{\mkern-2mu \gsphere}}}
\newcommand{\myspacetimetr}{{\mbox{\upshape{tr}}_g}}

\newcommand{\D}{\mathscr{D}}
\newcommand{\angD}{ {\nabla \mkern-14mu / \,} }

\newcommand{\angDsquared}{ {\angD_{\mkern-3mu}^2} }

\newcommand{\angdiv}{\mbox{\upshape{div} $\mkern-17mu /$\,}}
\newcommand{\angLap}{ {\Delta \mkern-12mu / \, } }

\newcommand{\Speed}{c_s}

\newcommand{\Transport}{B}

\newcommand{\Densrenormalized}{\uprho}

\newcommand{\Vortrenormalized}{\upomega}

\newcommand{\GdVar}{\upgamma}
\newcommand{\BadVar}{\underline{\upgamma}}

\newcommand{\Fullset}{\mathscr{Z}}
\newcommand{\Tanset}{\mathscr{P}}

\newcommand{\Singletan}{P}

\newcommand{\angdiff}{ {{d \mkern-9mu /} }}

\newcommand{\angdiffuparg}[1]{ {d \mkern-9mu /}^{#1} }

\newcommand{\Lie}{\mathcal{L}}
\newcommand{\SigmatLie}{\underline{\mathcal{L}}}
\newcommand{\angLie}{ { \mathcal{L} \mkern-10mu / } }

\newcommand{\Lgeo}{L_{(Geo)}}
\newcommand{\Lunit}{L}
\newcommand{\uLunit}{\underline{L}}
\newcommand{\uLgood}{\breve{\underline{L}}}
\newcommand{\Radunit}{X}
\newcommand{\Rad}{\breve{X}}

\newcommand{\NonRadialRad}{\Xi}
\newcommand{\XiCoordComp}{\upxi}

\newcommand{\CoordAng}{\Theta}

\newcommand{\Mult}{T}
\newcommand{\GeoAng}{Y}

\newcommand{\GeoAngFlatRadComponent}{y}

\newcommand{\angJ}{ {\mathscr{J} \mkern-14mu / \, } }

\newcommand{\angpi}{ { \pi \mkern-10mu / }}

\newcommand{\angk}{ { {k \mkern-10mu /} \, } }

\newcommand{\angkdoublearg}[2]{ {{k \mkern-10mu /}_{#1 #2} \, } }

\newcommand{\angG}{ {{\vec{G} \mkern-12mu /} \, }}
\newcommand{\angGarg}[1]{ {{\vec{G} \mkern-12mu /}_{\mkern 1mu #1} \, }}

\newcommand{\angGmixedarg}[2]{ {{\vec{G} \mkern-12mu /}_{#1}^{\ #2} \, }}
\newcommand{\angGnospacemixedarg}[2]{ {{\vec{G} \mkern-12mu /}_{#1}^{\mkern 2mu #2}}}

\newcommand{\NovecangGarg}[2]{ {{G \mkern-12mu /}_{\mkern 1mu #1}^{\mkern 6mu #2}}}

\newcommand{\angH}{ {{\vec{H} \mkern-12mu /} \, }}
\newcommand{\angHarg}[1]{ {{\vec{H} \mkern-12mu /}_{\mkern 1mu #1} \, }}

\newcommand{\NovecangH}[2]{ {{H \mkern-12mu /} \, }^{#1 #2}}
\newcommand{\NovecangHarg}[3]{ {{H \mkern-12mu /}_{\mkern 1mu #1}^{\mkern 6mu #2 #3}}}

\newcommand{\angxi}{ { {\xi \mkern-9mu /}  \, }}

\newcommand{\angxiarg}[1]{ {{\xi \mkern-9mu /}_{#1}  \, }}

\newcommand{\deform}[1]{{^{(#1)} \mkern-1mu \pi}}
\newcommand{\deformarg}[3]{{^{(#1)} \mkern-1mu \pi_{#2 #3}}}
\newcommand{\deformuparg}[3]{{^{(#1)} \mkern-1mu \pi^{#2 #3}}}

\newcommand{\angdeform}[1]{{^{(#1)} \mkern-2mu \angpi}}
\newcommand{\angdeformoneformarg}[2]{{^{(#1)} \mkern-2mu \angpi_{#2}}}
\newcommand{\angdeformoneformupsharparg}[2]{{^{(#1)} \mkern-2mu {\pi \mkern-10mu /}_{#2}^{\#}}}
\newcommand{\angdeformarg}[3]{{^{(#1)} \mkern-2mu \angpi_{#2 #3}}}

\newcommand{\Lineproject}{{\Pi \mkern-12mu / } \, }
\newcommand{\Sigmatproject}{\underline{\Pi}}

\newcommand{\vol}{\varpi}
\newcommand{\tvol}{\underline{\varpi}}
\newcommand{\conevol}{\overline{\varpi}}

\newcommand{\spherevol}{\uplambda_{{g \mkern-8.5mu /}}}
\newcommand{\argspherevol}[1]{\uplambda_{{g \mkern-8.5mu /}#1}}

\newcommand{\enzero}{\mathbb{E}^{(Wave)}}
\newcommand{\flzero}{\mathbb{F}^{(Wave)}}
\newcommand{\coercivespacetime}{\mathbb{K}}
\newcommand{\totmax}[1]{\mathbb{Q}_{#1}}
\newcommand{\easytotmax}[1]{\mathbb{Q}_{#1}^{(Partial)}}
\newcommand{\coercivespacetimemax}[1]{\mathbb{K}_{#1}}
\newcommand{\easycoercivespacetimemax}[1]{\mathbb{K}_{#1}^{(Partial)}}

\newcommand{\Vortenzero}{\mathbb{E}^{(Vort)}}
\newcommand{\Vortflzero}{\mathbb{F}^{(Vort)}}

\newcommand{\Vorttotmax}[1]{\mathbb{V}_{#1}}
\newcommand{\VorttotTanmax}[1]{\mathbb{V}_{#1}}

\newcommand{\upchifullmodarg}[1]{{^{(#1)} \mkern-4mu \mathscr{X}}}
\newcommand{\upchifullmodinhom}{\mathfrak{X}}

\newcommand{\upchipartialmodarg}[1]{{^{(#1)} \mkern-4mu \widetilde{\mathscr{X}}}}
\newcommand{\upchipartialmodinhom}{\widetilde{\mathfrak{X}}}
\newcommand{\upchipartialmodinhomarg}[1]{{^{(#1)} \mkern-2mu \widetilde{\mathfrak{X}}}}

\newcommand{\waveinhom}{\mathfrak{F}}
\newcommand{\vortinhom}{\mathfrak{F}}

\newcommand{\inhomleftexparg}[2]{{^{(#1)} \mkern-.5mu #2 }}

\newcommand{\upchipartialmodsourcearg}[1]{\inhomleftexparg{#1}{\mathfrak{B}}}

\newcommand{\Jcurrent}[1]{^{(#1)} \mkern-10mu \mathscr{J}}

\newcommand{\Vplus}[2]{{^{(+)} \mkern-.5mu  \mathcal{V}_{#1}^{#2}}}
\newcommand{\Vminus}[2]{{^{(-)} \mkern-.5mu  \mathcal{V}_{#1}^{#2}}}
\newcommand{\Sigmaplus}[3]{{^{(+)} \mkern-.5mu  \Sigma_{#1;#2}^{#3}}}
\newcommand{\Sigmaminus}[3]{{^{(-)} \mkern-.5mu  \Sigma_{#1;#2}^{#3}}}

\newcommand{\Contwo}{b}

\newcommand{\LateTimeLUnitMu}{\upkappa}

\newcommand{\smoothfunction}{\mathrm{f}}

\newcommand{\ThirdSmoothFunction}{\upeta}

\newcommand{\basicenergyerror}[1]{{^{(#1)} \mkern-.5mu \mathfrak{P}}}
\newcommand{\basicenergyerrorarg}[2]{{^{(#1)} \mkern-.5mu \mathfrak{P}}_{(#2)}}

\newcommand{\Cur}{\mathscr{R}}

\newcommand{\enmomtensor}{Q}

\newcommand{\BigConst}{C'}

\newcommand{\myarray}[2][]{\left(
		\begin{array}{lr}
    	 #1 \\
    	 #2 
     \end{array} \right)}

\newcommand{\threemyarray}[3][]{\left(
		\begin{array}{lr}
    	 #1 \\
    	 #2 \\
    	 #3
     \end{array} \right)}

\newcommand{\fourmyarray}[4][]{\left(
		\begin{array}{lr}
    	 #1 \\
    	 #2 \\
    	 #3 \\
    	 #4
     \end{array} \right)}

\newcommand{\f}{\frac}
\newcommand{\rd}{\partial}

\numberwithin{equation}{subsection}

\begin{document}
\title{
	Shock formation in solutions to the $2D$ compressible Euler equations 
	in the presence of non-zero vorticity
}
\author[JL,JS]{Jonathan Luk$^{*}$ and Jared Speck$^{** \dagger}$}

\thanks{$^{\dagger}$JS gratefully acknowledges support from NSF grant \# DMS-1162211,
from NSF CAREER grant \# DMS-1454419,
from a Sloan Research Fellowship provided by the Alfred P. Sloan foundation,
and from a Solomon Buchsbaum grant administered by the Massachusetts Institute of Technology.
}

\thanks{$^{*}$Stanford University, Palo Alto, CA, USA.
\texttt{jluk@stanford.edu}}

\thanks{$^{**}$Massachusetts Institute of Technology, Cambridge, MA, USA.
\texttt{jspeck@math.mit.edu}}

\begin{abstract}
We study the Cauchy problem for the
compressible Euler equations 
in two spatial dimensions under any physical barotropic equation of state except 
that of a Chaplygin gas. We prove that the 
well-known phenomenon of shock formation in simple plane wave solutions,
starting from smooth initial conditions,
is stable under perturbations that break the plane symmetry.
Moreover, we provide a sharp asymptotic description of the singularity formation.
The new feature of our work is that the perturbed solutions
are allowed to have small but non-zero vorticity, even at the location of the 
shock. Thus, our results provide the first constructive description
of the vorticity near a singularity formed from compression:
relative to a system of geometric coordinates
adapted to the acoustic characteristics, the vorticity 
remains many times differentiable, all the way up to the shock.
In addition, relative to the Cartesian coordinates, the vorticity remains bounded, 
and the specific vorticity remains uniformly Lipschitz, up to the shock.

To control the vorticity,
we rely on a coalition of new geometric and analytic insights
that complement the ones used by Christodoulou 
in his groundbreaking, sharp proof of shock formation in vorticity-free regions.
In particular, we rely on a new formulation of the compressible Euler equations 
(derived in a companion article) exhibiting remarkable structures.
To derive estimates, we construct an eikonal function adapted to the acoustic
characteristics (which correspond to sound wave propagation) 
and a related set of geometric coordinates and differential operators. 
Thanks to the remarkable structure of the equations, 
the same set of coordinates and differential operators 
can be used to analyze the vorticity, 
whose characteristics are transversal to the acoustic characteristics. 
In particular, our work provides the first constructive description of shock formation
without symmetry assumptions in a system with multiple speeds.

\bigskip

\bigskip

\end{abstract}

\maketitle

\centerline{\today}

\tableofcontents
\setcounter{tocdepth}{2}

\section{Introduction}
\label{S:INTRO}
In two spatial dimensions,
the isentropic compressible Euler equations 
are evolution equations for
the velocity $v:\mathbb{R} \times \Sigma \rightarrow \mathbb{R}^2$
and the density $\rho:\mathbb{R} \times \Sigma \rightarrow [0,\infty)$,
where $\Sigma$ is the (two-dimensional) space manifold, 
which we assume throughout to be $\Sigma = \mathbb{R} \times \mathbb{T}$.
Here and throughout, $\mathbb{T} := [0,1)$ 
(with the endpoints identified)
denotes the standard one-dimensional torus.
We now fix a constant $\bar{\rho} > 0$ corresponding to a constant background density.\footnote{In this paper, we will study solutions with density close to $\bar{\rho}$.}
Under a barotropic\footnote{A barotropic equation of state is one in which the pressure $p$
can be expressed as a function of the density $\rho$ alone.} 
equation of state
and in terms of the logarithmic density 
$\displaystyle \Densrenormalized := \ln \left(\frac{\rho}{\bar{\rho}} \right)$,
the equations
take\footnote{Throughout, if $V$ is a vectorfield
and $f$ is a function, then $V f := V^{\alpha} \partial_{\alpha} f$
denotes the derivative of $f$ in the direction $V$. Lower case Latin indices 
correspond to the Cartesian spatial coordinates and lower case Greek indices 
correspond to the Cartesian spacetime coordinates. 
We also use Einstein's summation convention.
\label{F:VECTORFIELDSACTONFUNCTIONS}}
the following form relative to the usual Cartesian coordinates,\footnote{Throughout, $\lbrace x^{\alpha} \rbrace$
are the usual Cartesian coordinates with corresponding partial derivative vectorfields
$
\displaystyle
\partial_{\alpha}
:=
\frac{\partial}{\partial x^{\alpha}}
$.
We also set $t := x^0$ and $\partial_t := \partial_0$.
}
$(i=1,2)$:
\begin{subequations}
\begin{align} \label{E:TRANSPORTDENSRENORMALIZEDRELATIVETORECTANGULAR}
	\Transport \Densrenormalized
	& = - \partial_a v^a,
		\\
	\Transport v^i 
	& = - \Speed^2 \delta^{ia} \partial_a \Densrenormalized,
	\label{E:TRANSPORTVELOCITYRELATIVETORECTANGULAR}
\end{align}
\end{subequations}
where $\Transport = \partial_t + v^a \partial_a$
is the material derivative vectorfield (see \eqref{E:MATERIALVECTORVIELDRELATIVETORECTANGULAR})
and $\Speed = \Speed(\Densrenormalized)$ 
is the speed of sound, which depends on the equation of state
(see \eqref{E:SOUNDSPEED}). Throughout this paper, we assume the normalization condition\footnote{As we explain in Sect.~\ref{S:GEOMETRICSETUP}, this can always be achieved by a change of variables.}
\begin{equation}\label{speed.normalization.intro}
\Speed(0) = 1,
\end{equation}
which simplifies some aspects of the analysis and presentation.

As has been known since the foundational work of Riemann \cite{bR1860} in one spatial
dimension, initially smooth solutions to the compressible Euler equations can form 
shock singularities in finite time, even though the solutions enjoy a conserved energy.\footnote{More precisely, 
solutions to
	the compressible Euler equations
	\eqref{E:TRANSPORTDENSRENORMALIZEDRELATIVETORECTANGULAR}-\eqref{E:TRANSPORTVELOCITYRELATIVETORECTANGULAR}
	enjoy the conserved energy
	\begin{align} \label{E:CONSERVEDENERGY}
		\int_{\Sigma_t}
			\rho e 
			+
			\frac{1}{2}
			\rho
			\delta_{ab} v^a v^b
		\, d^2 x,
	\end{align}
	where $e$, the specific internal energy,
	is given by $\displaystyle e = h - \frac{p}{\rho}$,
	where $p$ is the pressure and
	the enthalpy $h$ is as in Subsect.~\ref{sec.ideas.Euler}.
	However, the energy \eqref{E:CONSERVEDENERGY}
	is far too weak to prevent the formation of singularities.
	For this reason, it plays no role in our analysis.}
We recall that a shock singularity is such that the velocity and density
remain bounded while some first partial derivative of these quantities with respect
to the Cartesian coordinates blows up in finite time. This phenomenon is also
known in the literature as wave breaking.
Our main result is a proof of finite-time shock formation 
for solutions generated by an open set of regular Sobolev-class initial data 
in two spatial dimensions
verifying suitable relative smallness assumptions. 
The solutions that we study here are perturbations of simple plane wave solutions 
that are close, as measured by suitable Sobolev norms, to constant state solutions;
see Subsect.~\ref{sec.data} for further discussion.
The main new feature of our work is that 
the vorticity of the perturbed solutions,\footnote{Plane symmetric solutions have
vanishing vorticity.} 
defined to be 
$\omega := \partial_1 v^2 - \partial_2 v^1$,
is allowed to be non-zero
in the region where the shock forms.  
Actually, in our analysis, it is more convenient
to work with the \emph{specific vorticity} $\Vortrenormalized$, 
defined by
\[
\Vortrenormalized 
:= \omega/\exp(\Densrenormalized),
\]
since it satisfies a simpler evolution equation and better estimates.
As we describe in more detail later in this section,
\emph{our proof applies in particular 
to data such that the solution's vorticity 
is provably non-zero at the location
of the first shock singularity.} Therefore, to close the proof, 
we in particular have to control the vorticity 
(and, as it turns out, many of its derivatives)
in a past neighborhood of the first singularity.
To this end, we rely on a new formulation of the 
compressible Euler equations (see Prop.~\ref{P:GEOMETRICWAVETRANSPORTSYSTEM}), 
which we describe below in detail. 

We now provide a rough summary of our main results.
We plan to extend our results to the case of
three spatial dimensions in forthcoming work \cite{jLjS2017}.
As we briefly describe below, the case of three spatial dimensions requires
substantial new technical innovations compared to the case of two spatial dimensions.
The new innovations are tied to the need to derive, in three spatial dimensions,
elliptic estimates to control the vorticity at the top order. 
In contrast, elliptic estimates are not necessary
in two spatial dimensions; see \cite{jLjS2016a} for a more substantial overview of this issue.

\begin{theorem}[\textbf{Rough version}]
\label{T:ROUGH}
For any physical equation of state except that 
of the Chaplygin gas,\footnote{ The equation of state of a Chaplygin gas is $p = p(\rho)=C_0-\f{C_1}{\rho}$, where $C_0\in \mathbb R$ and $C_1>0$; see \eqref{E:EOSCHAPLYGINGAS}.}
there exists an \underline{open} set of regular initial data,
with elements close to the data of a subset of simple plane wave solutions, 
that leads to \underline{stable} finite-time shock formation. 
The specific vorticity, which is \underline{provably non-vanishing at the shock} for some of our solutions, 
remains uniformly Lipschitz relative to the Cartesian coordinates, all the way up to the shock.
Moreover, the dynamics are ``well-described'' by the irrotational Euler equations. 
\end{theorem}
\begin{remark}[\textbf{Assumption on the spatial manifold}]
Our assumption that $\Sigma = \mathbb{R} \times \mathbb{T}$
is mainly for technical convenience and is not of fundamental 
importance. For instance, the case $\Sigma = \mathbb{R}^2$ could be treated 
with a similar approach, though the set of initial data to which our methods apply
might be quite different (see also the discussion at the end of Subsect.~\ref{sec.data}).
\end{remark}

\begin{remark}[\textbf{Maximal classical development}]
	\label{R:MAXIMALDEVELOPMENT}
	Our main results provide information about the solution only
	up to the constant-time hypersurface of first blowup.
	However, thanks to the sharp estimates of Theorem~\ref{T:MAINTHEOREM},
	our results could in principle be extended to give a detailed description of a portion of 
	the maximal classical development\footnote{Roughly, the maximal classical development 
	is the largest possible classical solution that is uniquely determined by the data; see, for example, 
	\cites{jSb2016,wW2013}
	for further discussion.} of the data corresponding to times up to approximately twice 
	the time\footnote{Roughly, we could propagate our bootstrap assumptions for this amount of time.} 
	of first blowup, including the shape of the boundary and the behavior of the solution along it.
	More precisely, the estimates that we prove are similar to the 
	ones used by Christodoulou in his work \cite{dC2007}*{Ch.~15}, 
	in which he revealed the structure of a large irrotational portion
	of the maximal classical development 
	of solutions to the relativistic Euler equations;
	see also \cite{dCsM2014} for a similar picture of an irrotational portion of the maximal classical development 
	for solutions to the non-relativistic compressible Euler equations.
	However, for the sake of brevity, we have chosen not to carry out those arguments.
	We clarify that in obtaining their sharp picture
	of the boundary of the maximal classical development,
	the authors of \cites{dC2007,dCsM2014} relied on 
	technical non-degeneracy assumptions on the behavior 
	of the solution at the boundary;
	we would have to make similar non-degeneracy assumptions
	if we were to study the boundary of the maximal classical development in regions with vorticity.
	\end{remark}

We will give a precise version of Theorem~\ref{T:ROUGH} in Theorem~\ref{T:MAINTHEOREM}.
The estimates in  Theorem~\ref{T:MAINTHEOREM} give a precise sense in which the dynamics are ``well-described'' by the irrotational Euler equations. In particular, our proof shows that\footnote{We note here that in the irrotational case, 
Christodoulou--Miao have already proved \cite{dCsM2014}
that the quantities $\lbrace \partial_j v^i \rbrace_{i,j=1,2,3}$ blow up,
while $\lbrace v^i \rbrace_{i=1,2}$ and $\ln \left(\frac{\rho}{\bar{\rho}} \right)$
remain uniformly bounded, all the way up to the shock.} some of
the quantities $\lbrace \partial_j v^i \rbrace_{i,j=1,2,3}$ blow up,
while $\lbrace v^i \rbrace_{i=1,2}$, 
$
\displaystyle
\ln \left(\frac{\rho}{\bar{\rho}} \right)
$,
and $\omega$ remain uniformly bounded, all the way up to the shock.
Another important part of our proof is that the 
specific vorticity and vorticity
are more regular than the velocity estimates would suggest, 
both in terms of Cartesian coordinates and geometric\footnote{In order 
to capture the geometry of shock formation, we define geometric coordinates 
similar to the ones used by Christodoulou in \cite{dC2007}; 
see Subsect.~\ref{Chr.intro}.} coordinates.
In particular, to close our estimates, 
we must show that 
relative to the geometric coordinates,
the specific vorticity and vorticity are exactly as differentiable
as the velocity and density, which represents a gain of one derivative compared to viewing
vorticity as a first derivative of velocity.

\begin{remark}[\textbf{Assumption on spatial dimensions and the regularity of the vorticity}] 
In our proof, we rely on the assumption of two spatial dimensions to control the 
specific
vorticity and the top-order derivatives of the eikonal function 
(see Subsubsect.~\ref{SSS:IDENTIFYINGBLOWUPMECHANISM} for its definition). 
In particular, 
in two spatial dimensions, 
the specific vorticity equation is homogeneous (see \eqref{E:RENORMALIZEDVORTICTITYTRANSPORTEQUATION}), 
which allows for 
a relatively straightforward proof that the specific vorticity gains regularity.
In three spatial dimensions,
the specific vorticity equation contains an 
additional ``vorticity stretching'' term, which
introduces significant technical complications into
the analysis.
As we overviewed in our companion article \cite{jLjS2016a}, 
the vorticity stretching term can be controlled
using additional elliptic estimates,\footnote{We also note that one encounters 
other new technical difficulties in three spatial dimensions compared to the case of two spatial dimensions.
In particular, in three spatial dimensions, one must also derive elliptic estimates, 
distinct from those mentioned above for the specific vorticity, 
to control the top-order derivatives of the eikonal function.
However, an approach to implementing the elliptic estimates for the eikonal function in three spatial dimensions
has been well understood since the Christodoulou--Klainerman proof \cite{dCsK1993}
of the stability of Minkowski spacetime, and, in the context of shock formation,
since Christodoulou's work \cite{dC2007}.}
and a similar gain in regularity for the specific vorticity with respect to the geometric coordinates can be achieved;
we will treat this in detail in a future work. 
Nevertheless, since the case of two spatial dimensions already requires substantial new ideas, 
of interest in themselves, we have chosen to treat it separately here.

We also note that the gain in regularity for the vorticity
is familiar to the community of researchers who have proved well-posedness results
for the compressible Euler equations in the presence of a physical vacuum boundary; see, for example,
\cites{dChLsS2010,dCsS2011,dCsS2012,jJnM2011,jJnM2009}. 
However, in those works, 
the proof of the gain in regularity relied on the special properties of Lagrangian
coordinate partial derivative vectorfields. In the study of shock formation,
Lagrangian coordinates are entirely inadequate for measuring regularity
since they are not adapted to the acoustic characteristics (which we describe in detail later on),
whose intersection corresponds to the singularity. 
For this reason, in three spatial dimensions, 
we need to rely on a different approach, tied to our new formulation of the equations
(see below and Prop.~\ref{P:GEOMETRICWAVETRANSPORTSYSTEM}), 
which allows us to realize the gain in regularity 
relative to vectorfields adapted to the acoustic characteristics.
In fact, if regularity were the only consideration,
then our approach for gaining a derivative in the vorticity 
could be implemented with \emph{any} sufficiently smooth 
spanning set of vectorfields,
not just the geometric ones (described in Subsect.~\ref{sec.data})
that we use to study shock formation.
\end{remark}

The aforementioned results \cites{dC2007,dCsM2014} on shock formation for compressible fluids,
though foundational, crucially relied on the assumption that the fluid is irrotational, 
at least in a neighborhood near the shock.\footnote{A vorticity-free region near the shock can be achieved, for instance, in the small data dispersive setting by exploiting the fact that the characteristic speed for the vorticity is slower than the sound speed. See also the discussions in Subsect.~\ref{sec.history}.}
In the irrotational case, the dynamics are completely determined by a fluid potential
$\Phi$ and the Euler equations can be written as a single quasilinear scalar wave equation; this is
a big simplification compared to the structure of the compressible Euler equations with vorticity.
Moreover, we note that the assumption of irrotationality is very restrictive from a physical point of view. 
In particular, irrotational data constitute only a very small (infinite co-dimension with empty interior!) subset of all initial data. 
It is therefore of interest to prove, at the very least,
that previous shock formation results still hold 
under perturbations with small vorticity.
As we explain below, substantial new ideas are needed to  
accommodate the presence of even small amounts of vorticity near the singularity.
In this context, let us note that accommodating vorticity
is particularly relevant when one is interested in extending the solution (in a weak sense) after the shock is formed.
The reason is that
even if the fluid is initially irrotational, 
vorticity may be generated after a shock has formed \cite{dC2007}.

Moreover, in the larger context of the study of singularity formation for evolution partial differential equations, 
our theorem appears to be the first shock formation result in more than one spatial dimension
that involves a system of quasilinear wave equations 
coupled to another quasilinear evolution equation \emph{with a different characteristic speed}. 
More precisely, in the presence of vorticity,
the Euler equations exhibit the following two kinds of characteristics: acoustic characteristics
(corresponding to the propagation of sound waves)
and the integral curves of the material derivative vectorfield
(corresponding to the transporting of vorticity).
We hope that the techniques introduced here will be relevant to other 
problems featuring multiple characteristic speeds.

As we have alluded to above, the starting
point of our proof is a new formulation of the compressible Euler equations as
a coupled system of 
covariant
wave and transport equations.
The new formulation, which we derived in the companion article \cite{jLjS2016a}, 
is a consequence of \eqref{E:TRANSPORTDENSRENORMALIZEDRELATIVETORECTANGULAR}-\eqref{E:TRANSPORTVELOCITYRELATIVETORECTANGULAR},
obtained by differentiating the equations with suitable operators and observing remarkable cancellations.
One key
advantage of the new formulation,
stated below as Prop.~\ref{P:GEOMETRICWAVETRANSPORTSYSTEM},
is that the inhomogeneous terms exhibit surprisingly good null structures that are
preserved under commutations with well-constructed geometric vectorfields, 
adapted to the acoustic characteristics.
The good null structure, which we referred to as the ``strong null condition''
in \cite{jLjS2016a},
signifies the \underline{complete nonlinear absence}
of certain quadratic and higher-order interactions   
that we would not be able to control near the shock.
The new formulation also allows for the aforementioned
\emph{gain of regularity in the specific vorticity},
both in terms of Cartesian and geometric coordinates, which is central to closing the proof.
We will further discuss this in Subsect.~\ref{SS:MODELPROBLEM}.

\subsubsection{Organization of Sect.~\ref{S:INTRO}}
\label{SSS:INTROORGANIZATION}
We have organized the remainder of Sect.~\ref{S:INTRO} as follows: In Subsect.~\ref{sec.setup}, 
we describe the setup of the problem, noting in particular that we will only study the solution in the causal future of a portion of the initial data. In Subsect.~\ref{sec.data}, we describe the solution regime that we study and 
the size parameters that we use in the analysis. In Subsect.~\ref{SS:MODELPROBLEM}, we describe some new ideas in the proof of our main theorem (although we will postpone a more detailed discussion of the main ideas to Sect.~\ref{sec.ideas}). 
Finally, in Subsect.~\ref{sec.history}, we close the introduction with a discussion on relevant previous works.

\subsection{Setup of the problem}
\label{sec.setup}
Instead of studying the solution in the entire spacetime $\mathbb{R} \times \Sigma$,
we study only the future portion of the solution 
that is completely determined by the portion of the data
lying in the subset
$\Sigma_0^{U_0} \subset \Sigma_0$
of thickness $U_0$
(as measured by the eikonal function $u$, described below)
and on the curved null hyperplane portion
$\mathcal{P}_0^t$; see Figure~\ref{F:REGION}.

\begin{center}
\begin{overpic}[scale=.2]{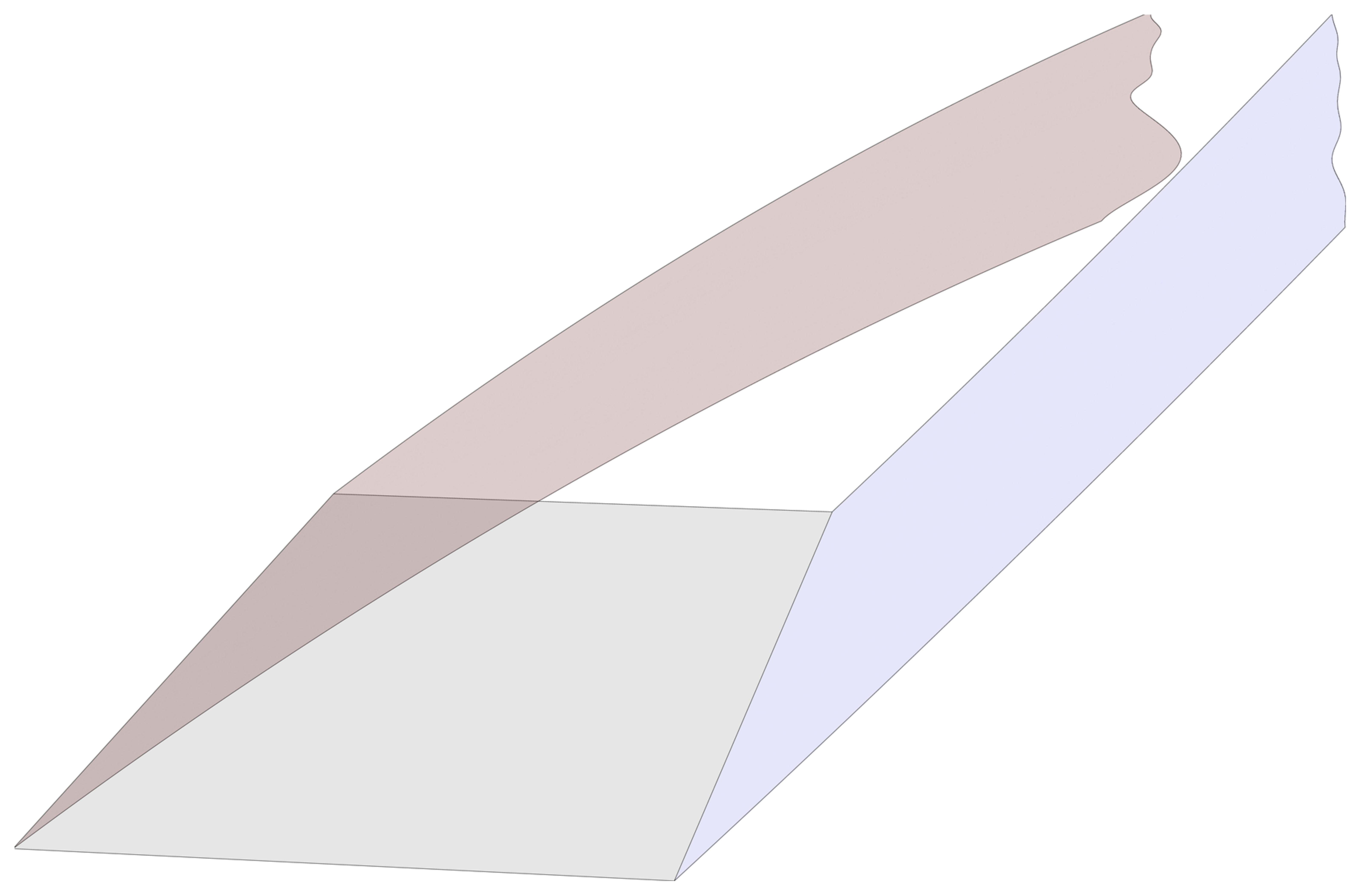} 
\put (37,33) {\large$\displaystyle \mathcal{P}_{U_0}^t$}
\put (74,33) {\large$\displaystyle \mathcal{P}_0^t$}
\put (18,5) {\large$\displaystyle \mbox{``interesting'' data}$}
\put (77,38) {\large \rotatebox{45}{$\displaystyle \mbox{very small data}$}}
\put (35,18) {\large$\Sigma_0^{U_0}$}
\put (31,10) {\large$\displaystyle U_0$}
%
\put (-.8,16) {\large$\displaystyle x^2 \in \mathbb{T}$}
\put (24,-3.2) {\large$\displaystyle x^1 \in \mathbb{R}$}
\thicklines
\put (-1.1,3){\vector(.9,1){22}}
\put (.5,1.8){\vector(100,-4.5){48}}
\put (10.5,13.9){\line(.9,1){2}}
\put (54.1,11.9){\line(.43,1){1}}
\put (11.5,15){\line(100,-4.5){43.2}}
\end{overpic}
\captionof{figure}{The spacetime region under study.}
\label{F:REGION}
\end{center}

The set $\Sigma_t$ is the level set of constant Cartesian time.
Moreover, here and throughout, 
$0 \leq t < 2 \TranminusdatasizeWithFactor^{-1}$, where
$\TranminusdatasizeWithFactor > 0$ is a data-dependent parameter
(see Def.~\ref{D:CRITICALBLOWUPTIMEFACTOR})
connected to the expected time of first shock formation,
and 
\begin{align} \label{E:FIXEDPARAMETER}
0 < U_0 \leq 1 
\end{align}
is a parameter,
fixed until Theorem~\ref{T:MAINTHEOREM}
(our main theorem).
The data that we treat are such that 
$\TranminusdatasizeWithFactor$ is large \emph{relative} to other 
parameters that control various seminorms of the data. 
That is, $\TranminusdatasizeWithFactor > 0$ can be arbitrary, 
but another parameter must be appropriately small; see Subsect.~\ref{sec.data} for further discussion.
We assume that the ``interesting, relatively large'' 
(in a sense that we explain below)
portion of the data 
lies in $\Sigma_0^{U_0}$.
Moreover,
we assume that the data 
are ``very small'' on $\mathcal{P}_0^{2 \TranminusdatasizeWithFactor^{-1}}$,
though the vorticity is allowed to be non-zero everywhere there.

More precisely, in
our main theorem, we consider perturbations of simple outgoing\footnote{Here and throughout, outgoing
simply means right-moving as is indicated in Figure~\ref{F:REGION}.} 
plane wave solutions. 
We focus our attention on perturbations of
the subset of these plane wave solutions 
that have data supported in $\Sigma_0^1$ 
and that satisfy the relative size condition mentioned above.
Domain of dependence considerations imply that
such plane wave solutions completely vanish along
$\mathcal{P}_0^{2 \TranminusdatasizeWithFactor^{-1}}$.
Thus, ``very small'' (not necessarily symmetric)
perturbations of their data on $\Sigma_0$, \emph{which we now allow
to have large spatial support}, \emph{induce}\footnote{Notice that while in principle one can attempt to directly prescribe data on the null hypersurface $\mathcal{P}_0^{2 \TranminusdatasizeWithFactor^{-1}}$, in practice, this involves solving a rather complicated system of \emph{constraint} equations, which can be difficult to implement.} 
``very small'' data along $\mathcal{P}_0^{2 \TranminusdatasizeWithFactor^{-1}}$ such that the data on 
$\Sigma_0^{U_0}\cup\mathcal{P}_0^{2 \TranminusdatasizeWithFactor^{-1}}$ is a perturbation of that of the simple plane wave solution. Moreover, it can be easily arranged that the vorticity is non-vanishing on 
$\Sigma_0^{U_0}\cup\mathcal{P}_0^{2 \TranminusdatasizeWithFactor^{-1}}$; since 
the specific vorticity is transported by the material derivative vectorfield (see the one-dimensional curves in Figure~\ref{F:FRAME} for a depiction of the integral curves of the material derivative vectorfield), if the vorticity is everywhere non-zero along 
$\Sigma_0^{U_0}\cup\mathcal{P}_0^{2 \TranminusdatasizeWithFactor^{-1}}$,
then it will be non-zero at the location of the first shock.

To summarize, the dynamic region of interest\footnote{Notice that by domain of dependence considerations, the solution in this region indeed depends only on the initial data on $\Sigma_0^{U_0}\cup\mathcal{P}_0^{2 \TranminusdatasizeWithFactor^{-1}}$.} lies in between
$\Sigma_0^{U_0}$ 
and the curved null hyperplane portions
$\mathcal{P}_{U_0}^t$
and
$\mathcal{P}_0^t$,
where $0 \leq t < 2 \TranminusdatasizeWithFactor^{-1}$.
We rigorously define these sets in Def.~\ref{D:HYPERSURFACESANDCONICALREGIONS},
but let us say a few words about them here and, at the same time, about
some other important spacetime subsets depicted in Figure~\ref{F:SOLIDREGION}.
The definitions of these sets refer to an eikonal function $u$, whose level sets 
are acoustic characteristics; we postpone our extensive discussion of $u$ until later.
For now, we simply note that the level sets of $u$ are denoted by
$\mathcal{P}_u$ or, when they are truncated at time $t$, by $\mathcal{P}_u^t$.
We refer to the $\mathcal{P}_u$ and $\mathcal{P}_u^t$ as
``null hypersurfaces,''
``null hyperplanes,''
``characteristics,''
or ``acoustic characteristics.''
We use the notation $\mathcal{M}_{t,u}$ to denote
the open-at-the-top region trapped in between 
$\Sigma_0$,
$\Sigma_t$,
$\mathcal{P}_0^t$,
and
$\mathcal{P}_u^t$.
We refer to the portion of 
$\Sigma_t$ trapped in between
$\mathcal{P}_0^t$
and
$\mathcal{P}_u^t$
as $\Sigma_t^u$.
The trace of the level sets of $u$ along $\Sigma_0$ are chosen (see condition \eqref{E:INTROEIKONALINITIALVALUE}) to be 
straight lines,
which we denote by $\ell_{0,u}$. For $t > 0$,
the trace of the level sets of $u$ along $\Sigma_t$
are (typically) curves\footnote{More precisely, the $\ell_{t,u}$ 
are diffeomorphic to the torus $\mathbb{T}$\label{FN:THEELLTUARETORI}.}  
$\ell_{t,u}$.
We restrict our attention to spacetime regions with $0 \leq u \leq U_0$.

\begin{center}
\begin{overpic}[scale=.2]{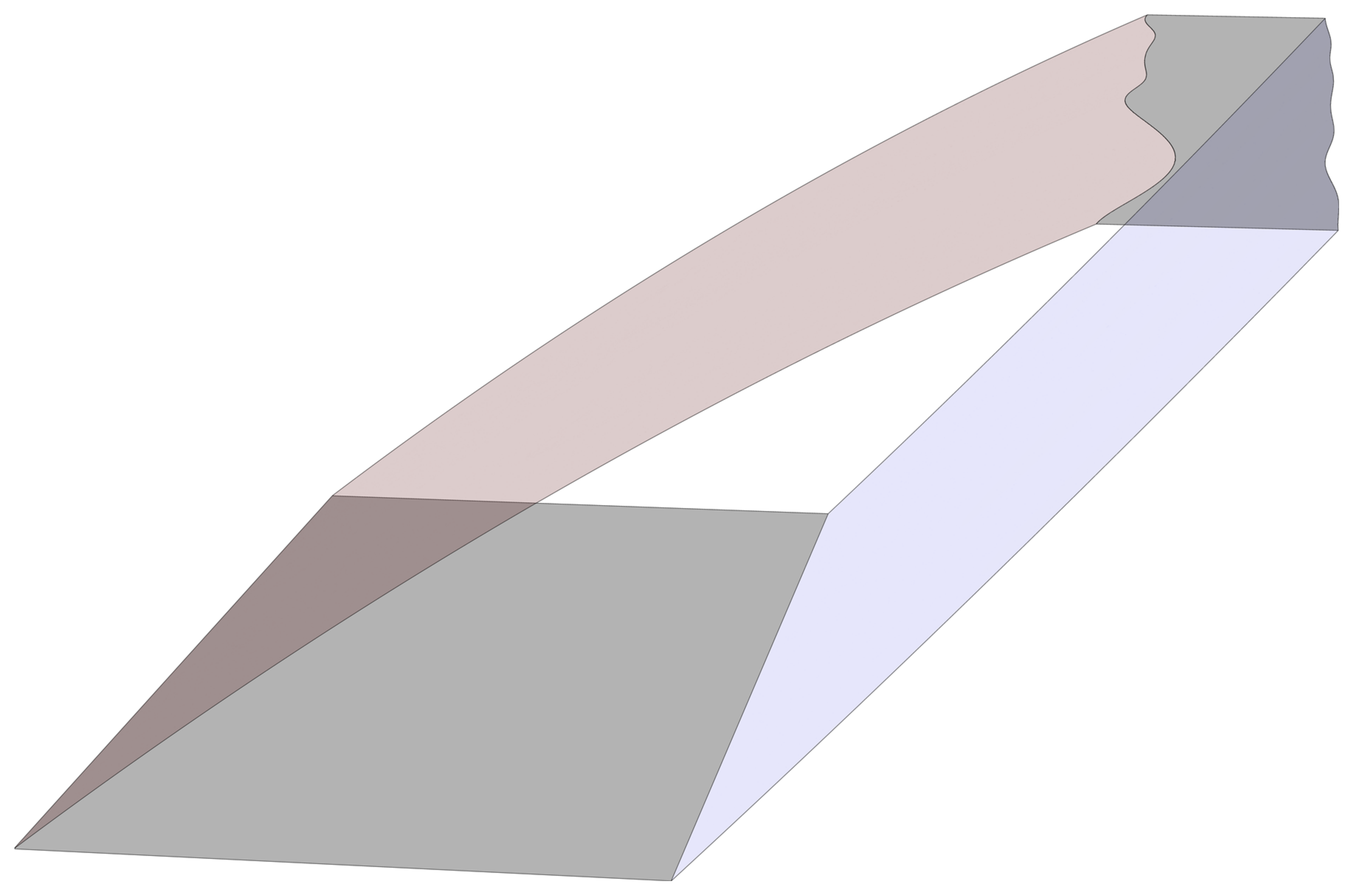} 
\put (54,34) {\large$\displaystyle \mathcal{M}_{t,u}$}
\put (38,33) {\large$\displaystyle \mathcal{P}_u^t$}
\put (74,37) {\large$\displaystyle \mathcal{P}_0^t$}
\put (59,15) {\large \rotatebox{45}{$\displaystyle \Densrenormalized,v^1,v^2,\Vortrenormalized 
\mbox{ very small}$}}
\put (32,17) {\large$\displaystyle \Sigma_0^u$}
\put (48.5,13) {\large$\displaystyle \ell_{0,0}$}
\put (12,13) {\large$\displaystyle \ell_{0,u}$}
\put (89.5,61) {\large$\displaystyle \Sigma_t^u$}
\put (93.5,56) {\large$\displaystyle \ell_{t,0}$}
\put (86,56) {\large$\displaystyle \ell_{t,u}$}
\put (-.6,16) {\large$\displaystyle x^2 \in \mathbb{T}$}
\put (22,-3) {\large$\displaystyle x^1 \in \mathbb{R}$}
\thicklines
\put (-.9,3){\vector(.9,1){22}}
\put (.7,1.8){\vector(100,-4.5){48}}
\end{overpic}
\captionof{figure}{The spacetime region and various subsets.}
\label{F:SOLIDREGION}
\end{center}

\subsection{Further description of the data and the solution regime}
\label{sec.data}
We now provide more details about
the data and solutions that we study. 
They are close to plane symmetric solutions
with data supported in the $x^1$ interval $[0,1]$
such that one Riemann invariant 
($\mathcal{R}_-:=
v^1
-
\int_{\widetilde{\Densrenormalized}=0}^{\Densrenormalized} \Speed(\widetilde{\Densrenormalized}) \, d \widetilde{\Densrenormalized}$)
completely vanishes,
while the other one 
($\mathcal{R}_+:=
v^1
+
\int_{\widetilde{\Densrenormalized}=0}^{\Densrenormalized} \Speed(\widetilde{\Densrenormalized}) \, d \widetilde{\Densrenormalized}$)
is initially small but with \emph{relatively large}\footnote{$\partial_1 \mathcal{R}_+$  
is allowed to be initially small in an absolute sense, 
as long as $\mathcal R_+$ is restricted to be even smaller.}
spatial derivatives.
Using the approach taken by Riemann in his famous work \cite{bR1860}
(in which he invented Riemann invariants), one may show that
$\partial_1 \mathcal{R}_+$ experiences 
a Riccati-type blowup
along a characteristic curve 
while $\mathcal{R}_+$ remains bounded; we stress that these phenomena occur while $\mathcal R_-$ remains identically zero.
Such a solution is known as a \emph{simple plane wave} and
arguably represents the simplest kind of symmetric shock-forming solution that
one can study in a perturbative sense.
It is for this reason that our main results apply to neighborhoods of
a subclass of simple plane waves.

To further describe the class of data and solutions under study, 
we need to introduce the geometric vectorfields $\Lunit$, $\GeoAng$ and $\Rad$, 
which are depicted in Figure~\ref{F:FRAME} below,
and also the parameters $\mathring{\updelta}_*$, 
$\mathring{\updelta}$, 
and $\mathring{\upepsilon}$ describing the sizes of different quantities.
At this point, let us just note that the vectorfields $\Lunit$ and $\GeoAng$ are chosen to be tangential to the 
acoustic characteristics, while $\Rad$ is transversal to them. Moreover, we use
$\Lunit$, $\Rad$, and $\GeoAng$ to commute the equations and obtain estimates for the solution's derivatives.
We refer the readers to Subsect.~\ref{SS:FRAMEANDRELATEDVECTORFIELDS} for 
rigorous definitions of these vectorfields and proofs of their basic properties.
The parameters $\TranminusdatasizeWithFactor$ and $\mathring{\updelta}$
are not necessarily small, 
but we require $\mathring{\upepsilon}$ to be \emph{relatively small}
in a sense explained in Subsect.~\ref{SS:SMALLNESSASSUMPTIONS}.
With the help of these parameters, we can now further describe the data
and solutions under study 
(where it may be seen that some of the conditions below are redundant):

\begin{itemize}
\item (Nearly plane symmetric perturbations of constant states) 
		The initial data and solution are $\mathring{\upepsilon}$ close 
		in $L^\infty$ to the constant state 
		$(\Densrenormalized, v^1, v^2) \equiv (0,0,0)$.
		By a ``plane symmetric'' solution, we mean 
		one such that, relative to the standard Cartesian coordinates, 
		we have $\Densrenormalized = \Densrenormalized(t,x^1)$, $v^1 = v^1(t,x^1)$, and $v^2 \equiv 0$.
		The factor $\mathbb{T}$ in the Cauchy hypersurface
		$\Sigma_0 \simeq \mathbb{R} \times \mathbb{T}$ corresponds to the direction of symmetry
		for the simple plane symmetric waves that we are perturbing.
		Thus, ``nearly plane symmetric'' means, roughly, small dependence in the
		$\mathbb{T}$ direction.
\item (Finite-time shock formation occurs) $\TranminusdatasizeWithFactor$ is defined so that $\TranminusdatasizeWithFactor^{-1}$ is 
the expected blowup time,
up to $\mathcal{O}(\mathring{\upepsilon})$ error.
$\TranminusdatasizeWithFactor$ is tied to the size of the (signed part of) $\Rad v^1|_{t=0}$;  
see \eqref{E:CRITICALBLOWUPTIMEFACTOR} for the precise definition.
\item (Boundedness of the transversal derivatives) 
$0 < \mathring{\updelta} < \infty$ bounds the initial size of the $\Rad$ (and $\Rad\Rad$ and $\Rad\Rad\Rad$, etc.) derivatives of 
$(\Densrenormalized, v^1, v^2)$.
In addition, we make smallness assumptions 
on the derivatives of $\Densrenormalized - v^1$ and $v^2$,
consistent with the behavior of the simple plane symmetric waves
that we are perturbing; see the next item.
\item (Nearly simple outgoing) Initially,\footnote{Notice that only the initial data for
$(\Densrenormalized, v^1, v^2)$ can be prescribed and that
we cannot prescribe their time derivatives. Nevertheless, 
for perturbations of appropriate simple plane wave solutions, the desired smallness can be achieved.}
 the
$\Lunit$ and $\GeoAng$ derivatives of 
$(\Densrenormalized, v^1, v^2)$ are $\mathring{\upepsilon}$-small at all
derivative levels in appropriate norms.
The same holds true for higher-order derivatives in terms of $\Lunit$, $\GeoAng$ and $\Rad$, where at least one of the derivatives is $\Lunit$ or $\GeoAng$.
Moreover, we assume that the first-order $\Rad$ derivatives of
$\Densrenormalized - v^1$ and all directional derivatives of $v^2$ are of size $\mathring{\upepsilon}$.
Roughly, these conditions correspond to a solution whose dynamics are well-described
by an outgoing (that is, moving in the direction of increasing $x^1$) and nearly simple.
\item (Smallness of the Riemann invariant $\mathcal R_-$) Initially, 
$
\displaystyle
\mathcal{R}_-
:=
v^1
-
\int_{\widetilde{\Densrenormalized}=0}^{\Densrenormalized} \Speed(\widetilde{\Densrenormalized}) \, d \widetilde{\Densrenormalized}
$ 
and \emph{all of its directional derivatives} up to top-order 
are initially $\mathring{\upepsilon}$ small in appropriate norms.
This represents a perturbation of the complete vanishing of
$\mathcal R_-$, which, in view of the above discussion, therefore 
corresponds to a perturbation of a simple outgoing plane symmetric wave.
\item (Small vorticity) This means that $\partial_1 v^2 - \partial_2 v^1$ and all of its derivatives up to top-order, 
in all directions, are initially $\mathring{\upepsilon}$ small in appropriate norms.
\item (Near-acoustic regime) In this regime, the compressible
		Euler equations are well-approximated by transport equations along 
		appropriately scaled null\footnote{By null, we mean relative to the acoustical metric of Def.~\ref{D:ACOUSTICALMETRIC}.} 
		generators $\Lunit$ of the 
		acoustic characteristics $\mathcal{P}_u$. 
		This corresponds to a flow that is dominated by sound wave 
		propagation.
\end{itemize}
Of course, one of our main tasks 
in our proof is showing that the smallness conditions stated above, 
at first only assumed for the initial data, are propagated by the 
flow. Actually, as we explain below in great detail, one of the main difficulties
in the proof is that near the top-order, the smallness can only be understood 
in terms of some singular norms. By this, we mean that the best
estimates we are able to prove allow for the possibility that the 
high-order energies might blow up.

The solution regime described above is not the only one for which we could prove
a shock formation result. However, it is perhaps the simplest one allowing
for non-zero vorticity. In particular, because the solutions are nearly plane symmetric,
there is no wave dispersion and hence no decay. Therefore, powers of $t$ and $r$
do not play a role in our analysis. We expect that a similar shock formation result
could be proved for small, nearly radially symmetric quickly decaying\footnote{By decaying, we mean in 
the Euclidean radial coordinate $r$, 
towards the data of a non-vacuum constant state.} 
data on $\mathbb{R}^2$. In this case, the analysis would involve factors of $t$
and/or $r$, which would capture the dispersive decay 
(that one expects to occur until close to the shock).
Moreover, one would have to assume that the vorticity is initially very small,
so that it and its derivatives are not able 
to become large by the time that the shock forms.

\subsection{New ideas for the proof}
\label{SS:MODELPROBLEM}
To prove our main theorem, we rely on the full strength of the technology developed in the works of Christodoulou \cite{dC2007} and Speck--Holzegel--Luk--Wong \cite{jSgHjLwW2016}. In particular,
the same size parameters $\TranminusdatasizeWithFactor$, $\mathring{\updelta}$, and $\mathring{\upepsilon}$, which are featured in the proof (and discussed in Subsect.~\ref{sec.setup}), 
are also present in \cite{jSgHjLwW2016}.
We therefore postpone a detailed discussion of our proof until Sect.~\ref{sec.ideas},
where we review the works \cites{dC2007,jSgHjLwW2016}.
Here, we will simply highlight a few key new high-level ideas
(see Subsect.~\ref{sec.ideas.vorticity} for a discussion of more technical new ideas):
\begin{enumerate}
\item In Prop.~\ref{P:GEOMETRICWAVETRANSPORTSYSTEM}, 
we reformulate the equations as a system of coupled wave and transport equations
with remarkable geometric features, including the good null structures mentioned above.
\item We prove that the transport part of the system ``interacts well'' with the wave part of the system. More precisely, one can commute geometric vectorfields 
adapted to the acoustic characteristics
$\mathcal{P}_u$
through an appropriately weighted version of the material derivative vectorfield, which is 
the principal part of the transport equation for the specific vorticity.
\item We show that the specific vorticity is uniformly
Lipschitz with respect to Cartesian coordinates up to the formation of the first shock, 
which is a much stronger estimate than what follows from
simply viewing the specific vorticity as first derivatives of $v^i$ divided by $\rho$.
\item We prove that the specific vorticity is ``one derivative better'' with respect to geometrically defined vectorfields than one naively expects, thus avoiding an apparent loss of derivatives in the new formulation of the equations.
\end{enumerate}

Our reformulation of the compressible Euler equations was derived in \cite{jLjS2016a} 
in the case of three spatial dimensions but can be easily modified so as to apply in
two spatial dimensions.
We present the two-space-dimensional version in Prop.~\ref{P:GEOMETRICWAVETRANSPORTSYSTEM} below. 
In two spatial dimensions, 
the new formulation can be modeled by
the following wave-transport system in the scalar
unknowns $\Psi$ 
(which models $v^i$ and $\Densrenormalized$) and $w$ 
(which models the specific vorticity,
defined above as
$\Vortrenormalized = \omega/\exp(\Densrenormalized)$):
\begin{align} \label{E:MODELWAVE}
	\square_{g(\Psi)} \Psi 
	& = \partial w,
		\\
	\partial_t w
	& = 0. \label{E:MODELTRANSPORT}
\end{align}
In \eqref{E:MODELWAVE}, $g = g(\Psi)$ is a Lorentzian metric
whose Cartesian components $g_{\alpha \beta}$ 
are assumed to be
explicit smooth functions of $\Psi$,
$\square_{g(\Psi)}$ is the covariant wave operator\footnote{Relative to arbitrary coordinates,
$\square_g f
= 
\frac{1}{\mbox{$\sqrt{\mbox{\upshape det} g}$}}
\partial_{\alpha}\left( \sqrt{\mbox{\upshape det} g} (g^{-1})^{\alpha \beta} \partial_{\beta} f \right)$.
\label{FN:COVWAVEOPARBITRARYCOORDS}} 
of $g(\Psi)$,
and $\partial w$ schematically denotes first-order Cartesian coordinate partial derivatives 
of $w$.
In our study of the compressible Euler equations,
$g$ is the acoustical metric (see Def.~\ref{D:ACOUSTICALMETRIC})
corresponding to the propagation of sound waves.
In Cartesian coordinates, 
the expression $\square_{g(\Psi)} \Psi$ contains (quasilinear) principal terms
of the schematic form $\smoothfunction(\Psi) \partial^2 \Psi$ and
semilinear terms of the form $\smoothfunction(\Psi) (\partial \Psi)^2$.
The precise structures of both the quasilinear and the semilinear terms are
important for our analysis.
Equation \eqref{E:MODELTRANSPORT} models the transporting of specific vorticity.
In writing down \eqref{E:MODELWAVE}-\eqref{E:MODELTRANSPORT},
we have omitted the quadratic inhomogeneous terms from Prop.~\ref{P:GEOMETRICWAVETRANSPORTSYSTEM}, 
all of which have a good null structure and
remain negligible, all the way up the shock.
The presence of this null structure,
which is available thanks 
to the special form of the equations stated in Prop.~\ref{P:GEOMETRICWAVETRANSPORTSYSTEM},
is fundamental for our proof;
see Remark~\ref{R:NULLFORMSAREIMPORTANT} for further discussion.

We also note the following aesthetically appealing
feature of the formulation: the principal parts of the system
are a wave operator and a transport operator. Thus, the two kinds
of propagation phenomena present in the compressible Euler equations,
namely the propagation of sound waves and the transporting
of vorticity, become manifest.
This stands in contrast to the usual
first-order formulation 
\eqref{E:TRANSPORTDENSRENORMALIZEDRELATIVETORECTANGULAR}-\eqref{E:TRANSPORTVELOCITYRELATIVETORECTANGULAR},
where the presence of the two kinds of propagation phenomena
are not easily visible at the level of the equations. 

Previous shock formation results, which we review in Sects.~\ref{sec.history} and \ref{sec.ideas},
apply to quasilinear wave equations. 
In contrast, in the model problem \eqref{E:MODELWAVE} and \eqref{E:MODELTRANSPORT}, 
we need to handle an extra transport equation and also additional inhomogeneous terms in the wave equation.
In previous works on shock formation in quasilinear wave equations, 
starting from 
\cites{sA1995,sA1999a,sA1999b,dC2007}, 
a crucial insight was to use geometric vectorfields 
that are adapted to the characteristics $\mathcal{P}_u$ 
and that, in directions transversal to the characteristics, 
are appropriately degenerate (with respect to the Cartesian coordinate vectorfields) near the shock. 
Morally, this is equivalent to deriving estimates relative to a system of geometric coordinates adapted to the characteristics.
To accommodate the term $\partial w$ on RHS~\eqref{E:MODELWAVE},
it is therefore important when dealing with the coupled system to ensure that 
the derivatives of the specific vorticity with respect to the \emph{same geometric vectorfields}
can be controlled. To achieve this, we rely on the fact that the transport operator is a \emph{first-order} differential operator and therefore, upon multiplying by a degeneration factor $\upmu$ (explained below in great detail), that
\emph{commuting the transport equation with the geometric vectorfields generates only controllable error terms}.

Next, we note that RHS~\eqref{E:MODELWAVE} involves a Cartesian coordinate partial derivative of $w$, 
which is therefore singular with respect to the geometric vectorfields.\footnote{This singularity is actually unimportant just from the point of view of the lower-order energy estimates. This is, however, of crucial importance at the top-order; 
see discussions in Sect.~\ref{sec.ideas}.} However, the following crucial geometric fact is available
in our formulation of the compressible Euler equations:
the transport equation has a \emph{strictly smaller speed compared to} 
the characteristic wave speed corresponding to the operator
$\square_g$. For this reason, in the actual problem under study, we can use 
the transport equation to express the transversal (to the acoustic characteristics of $g$) derivatives of $w$ 
in terms of the non-degenerate tangential derivatives of $w$. 
This can be used to show, among other things, that $w$ is in fact uniformly Lipschitz up to the shock.
The difference in the characteristic speeds for the transport operator and the wave operator 
is also important in that it leads to the availability of non-degenerate energies for
$w$ along the acoustic characteristics corresponding to $g$; see the last term on
RHS~\eqref{E:INTROVORTEN}.

Finally, we discuss the basic regularity of the solution variables,
highlighting the role of the source term $\partial w$ on the right-hand side of the wave equation \eqref{E:MODELWAVE}.
In the case of the compressible Euler equations, vorticity can be viewed as the first derivatives of the velocity 
and hence, in the context of the regularity of solutions to the model problem, 
one might be tempted to think of $\partial w$
as corresponding to the \emph{second} derivatives of $\Psi$.
However, this perspective is insufficient from the point of view of regularity since
energy estimates for the wave equation (without commutation) 
yield control of only 
one derivative of $\Psi$. Hence, this perspective leads to an apparent loss of a derivative. 
However, since \eqref{E:MODELTRANSPORT} is a homogeneous transport equation, one expects to gain a derivative 
-- this is indeed obvious\footnote{From this point of view, the model system is oversimplified 
in that one can control an arbitrarily large number of Cartesian coordinate partial derivatives of $w$. 
In the actual system, since the transport operator depends also on the Cartesian components $(v^1,v^2)$, 
only one derivative can be gained.} 
if one takes Cartesian coordinate partial derivatives of equation \eqref{E:MODELTRANSPORT}.
What is less obvious is that in fact, the loss of derivatives can also be avoided
if one differentiates the transport equation with the geometric vectorfields
which, as it turns out, depend on $\Psi$.
We note that while it is indeed possible to carry out commutations 
the transport equation with geometric derivatives,
one encounters some singular terms tied to the degenerate top-order behavior of
$\Psi$ and the acoustic geometry, which we will discuss in detail in Sect.~\ref{sec.ideas}.

\subsection{History of the problem}
\label{sec.history}
The study of the formation of shocks for the compressible Euler equations has a long history which traces back to the aforementioned foundational 
work of Riemann \cite{bR1860}. He introduced the \emph{Riemann invariants} for the  Euler equations 
in one spatial dimension
and showed that shocks often form in finite time. In the one-dimensional case, the theory, at least in the small BV regime, 
is fairly complete. In particular, it is known that there exist unique global weak solutions in the BV class. This theory in particular incorporates formation and interaction of shocks. We refer the readers to \cites{gqCdW2002,cD2010} for surveys on the one-dimensional case.

Let us mention that the theory of finite-time blowup for solutions to hyperbolic systems in one spatial dimension
has been developed way beyond the theory of the compressible Euler equations. For example, for general $2\times 2$ genuinely nonlinear hyperbolic systems, finite-time blowup has been proven by Lax \cite{pL1964}. 
For $n \times n$ genuinely nonlinear hyperbolic systems, 
even though Riemann invariants are not available, 
John \cite{fJ1974} has obtained a shock formation result in which the waves are simple by the time a shock forms.

In two or three spatial dimensions (without symmetry assumptions), the problem becomes considerably harder.
The first general breakdown result for the compressible Euler equations 
in three spatial dimensions was achieved by Sideris \cite{tS1985} for 
a polytropic gas\footnote{That is, the equation of state is given by 
$p(\rho)=k \rho^\gamma$ for constants $k$ and $\gamma$ with $k>0$.
Actually, Sideris allowed for the presence of non-constant entropy; his
breakdown result holds for the equations of state $p(\rho)=k \rho^\gamma \exp(s)$,
where $s$, the specific entropy, verifies the evolution equation
$\Transport s = 0$.}
with adiabatic index $\gamma>1$. In particular, he 
exhibited an open set of small and regular initial data for which the corresponding solutions cease to be $C^1$ in finite time. 
However, his methods
did not
provide any information on the nature of the breakdown.
 
In a different direction, Alinhac studied the two-dimensional compressible isentropic Euler equations in \emph{radial symmetry} \cite{sA1993}. He showed that a large class of small radially symmetric data (with potentially non-vanishing vorticity\footnote{However, since the initial vorticity is required to be compactly supported and the speed of the vorticity is much slower than the sound speed, the vorticity in Alinhac's solutions vanishes in a neighborhood to the past of the first singularity.}) lead to a finite time blow up. While this result only applies to radial initial data, it gives a precise estimate on the blow up time 
(at least as the size of the data tends to $0$).

Alinhac later achieved \cites{sA1995,sA1999a,sA1999b}
important breakthroughs regarding shock formation.
His works, which addressed solutions to a large class of quasilinear wave equations,
were the first instances of proofs of shock formation for solutions
to quasilinear equations in more than one spatial dimension that did not rely on 
any symmetry assumptions. In particular, his work yielded a precise description of the 
singularity and tied its formation to the intersection of the characteristics.
While he did not explicitly study the compressible Euler equations,
his works provided all of the main insights needed to extend the result to the \emph{irrotational} Euler equations.
More precisely, for 
all quasilinear wave equations $(g^{-1})^{\alpha \beta}(\partial \Phi) \partial_{\alpha} \partial_{\beta} \Phi = 0$
that fail to satisfy the null condition, 
Alinhac exhibited a set of initial data leading to finite-time shock formation. As we will discuss in Subsect.~\ref{sec.ideas.Euler}, under the irrotationality assumption, the compressible Euler equations can be written as a quasilinear wave equation in the above form. Moreover, the null condition is violated whenever the equation of state is not that of a Chaplygin gas.  The data in Alinhac's works
were small and satisfied a \emph{non-degeneracy condition}.
For this class of data, he gave precise estimates of the solution up to the first singular time. 
In his proof, he recognized the importance of 
deriving estimates relative to a geometric coordinate system tied to
an eikonal function, which captures the geometry of shock formation. 
However, Alinhac's approach to deriving energy estimates 
was based on a Nash--Moser iteration scheme featuring a free boundary,
and the iteration scheme relied in a fundamental way on his non-degeneracy condition on the initial data.

In a monumental work in 2007, Christodoulou \cite{dC2007} studied shock formation for all\footnote{Actually, there is one exceptional equation of state such that the null condition is satisfied, in which case the corresponding wave equation
admits small data global solutions \cite{hL2004}; 
see \cite{dC2007} for further discussion. The exceptional equation of state
corresponds to the equation of state of a Chaplygin gas in the non-relativistic case, for which
our main shock formation results do not apply.} 
wave equations of irrotational relativistic fluid mechanics.\footnote{Roughly speaking, these equations form a subclass of the equations studied by Alinhac and
enjoy special additional properties, such as an Euler-Lagrange structure
and invariance under the Poincar\'{e} group. However, as is shown in \cite{jS2014b}, the insights introduced in \cite{dC2007} can be applied to a much larger class of quasilinear wave equations.} Christodoulou
proved that a large class of small initial data give rise to shock formation and 
he gave a precise description of a portion of the boundary of the maximal classical development of the data.
Compared to the work of Alinhac, Christodoulou introduced a fully geometric framework such that the breakdown of the solution is completely described in terms of the vanishing of the inverse foliation density $\upmu$ (see definition \eqref{E:FIRSTMU})
of the acoustic characteristics.
As a consequence, 
his work applied to an open neighborhood of
solutions whose data are small and compactly supported perturbations of the non-vacuum constant states. 
In particular, for data that are small as measured by a high-order Sobolev norm, he
showed that (at least outside the causal future of a compact set) shocks are the only possible singularities. 
Moreover, he exhibited an open condition on the data that guarantees that a shock will form in finite time.\footnote{Given Christodoulou's result that shocks are the only possible singularities, there remains a possibility of some non-trivial global solutions arising from small data.}

The geometric framework introduced in \cite{dC2007} has 
proven to be useful for studying shock formation 
in other settings.
Most relevant to our current work is the aforementioned work of Christodoulou--Miao \cite{dCsM2014}, 
which used the geometric insights of \cite{dC2007} to study shock formation for small and compactly supported 
perturbations of non-vacuum constant state solutions 
to the non-relativistic compressible Euler equations. In particular,\footnote{Note that Sideris' proof of blowup by contradiction \cite{tS1985} applies only to adiabatic equations of state with index bigger than one, 
while the work \cite{dCsM2014} allows for an arbitrary barotropic equation of state
(except that of the Chaplygin gas).} 
the results of \cite{dCsM2014} provided a precise picture of the singularity formation exhibited by Sideris in \cite{tS1985}.

While the shock formation result of \cite{dC2007} was proved
in \emph{irrotational} regions of spacetime,
the result also applies to initial data with non-vanishing vorticity \emph{which satisfies appropriate conditions on its (compact) support}. This is because for such initial data, using that the vorticity and sound travel with different speeds, 
one can show that in the complement of the causal future of an appropriate compact set, the vorticity vanishes. In particular, 
in the small-data regime, the vorticity travels with 
small speed and hence completely vanishes in the acoustic wave zone,
where the shock forms. A similar result could be proved for the non-relativistic compressible Euler equations
using the techniques of \cite{dCsM2014}, even though such a result was not stated there.
However, we stress that the approach of \cites{dC2007,dCsM2014}
is not sufficient, in itself, for controlling solutions with non-vanishing at the first shock singularity;
for this, one seems to need all of the new structural features afforded by Prop.~\ref{P:GEOMETRICWAVETRANSPORTSYSTEM}.

The seminal work of Christodoulou also inspired some recent developments on shock formation for quasilinear wave equations in 
more than one spatial dimension.
See \cites{jS2014b,gHsKjSwW2016,sMpY2014} for a sample of such results. 
The work \cite{jS2014b} in particular generalized the results in \cite{dC2007} to a much larger class of quasilinear wave equations. We refer the readers also to \cites{dCaL2016,dCdRP2015} 
for some recent developments in \emph{symmetry-reduced} problems motivated by \cite{dC2007}.

\subsection{Outline of the paper} 
\label{SS:PAPEROUTLINE}
In Sect.~\ref{sec.ideas}, we describe the main ideas behind our proof.
Although we need many new ideas to treat the vorticity,
we extensively rely on the framework developed by Christodoulou \cite{dC2007} 
in his proof of shock formation in irrotational regions 
and on the methods of Speck--Holzegel--Luk--Wong,
who proved \cite{jSgHjLwW2016} shock formation for perturbations of
simple outgoing plane waves for general classes of quasilinear wave equations.
Hence, we review the relevant aspects of those works in detail.
Readers who are familiar with those works
might prefer to skip to Subsect.~\ref{sec.ideas.vorticity},
where we overview the main new ideas needed to handle the presence
of vorticity at the shock.

Starting in Sect.~\ref{S:GEOMETRICSETUP}, we give detailed proofs.
Specifically, in Sects.~\ref{S:GEOMETRICSETUP}-\ref{S:MODIFIED},
we construct all of the geometric quantities that we need to study the solution.
In Sect.~\ref{S:DATASIZEANDBOOTSTRAP} we describe our 
assumptions on the data and formulate suitable $L^{\infty}$-type
bootstrap assumptions. 
In Sects.~\ref{S:PRELIMINARYPOINTWISE}-\ref{S:SHARPESTIMATESFORUPMU}
and \ref{S:POINTWISEESTIMATESFORWAVEEQUATIONERRORTERMS},
we use the bootstrap assumptions to derive
$L^{\infty}$ and pointwise estimates for the solution.
In Sect.~\ref{S:FUNDAMENTALL2CONTROLLINGQUANTITIES},
we construct the $L^2$-type quantities that 
we later bound with energy estimates.
In Sect.~\ref{SS:SOBOLEVEMBEDDING},
we provide a geometric Sobolev embedding theorem,
which we will use to recover the $L^{\infty}$ bootstrap assumptions
from energy estimates.
Sect.~\ref{S:ENERGYESTIMATES} is the most important 
part of the paper. There we use the previous estimates
to derive a priori energy estimates for the $L^2$-type
quantities mentioned above. In Sect.~\ref{S:STABLESHOCKFORMATION}, 
we prove our main shock formation theorem,
including recovering the bootstrap assumptions and showing
that the shock forms.
The theorem is relatively easy to prove given the 
estimates from the prior sections.

\section{Ideas of the proof}
\label{sec.ideas}

In this section, we describe the ideas of the proof of our main theorem. While the main novelty in this paper is that we allow for non-vanishing vorticity all the way up to shock formation, in order to describe our proof, we nonetheless have to recall some of the main points in the work of Christodoulou \cite{dC2007} and the work of Speck--Holzegel--Luk--Wong \cite{jSgHjLwW2016}. In particular, in the present work, we will work with a solution regime similar to that in \cite{jSgHjLwW2016}.

We have organized Sect.~\ref{sec.ideas} as follows: 
In Subsect.~\ref{Chr.intro}, we review the work \cite{dC2007}, emphasizing the geometric insights that are relevant to our present work. In Subsect.~\ref{sec.ideas.review.plane}, we review the work \cite{jSgHjLwW2016}. In Subsect.~\ref{sec.ideas.Euler}, 
we discuss how the work \cite{jSgHjLwW2016} 
can be applied to the compressible Euler equations and how \cite{jSgHjLwW2016} is related to our present work. 
Finally, in Subsect.~\ref{sec.ideas.vorticity}, 
we discuss the main new ingredients that we use to prove our main theorem, 
which requires controlling the interaction between sound waves and vorticity
up to the first singularity caused by compression.

\subsection{Review of Christodoulou's work}\label{Chr.intro}
We begin with a review of the main ideas in \cite{dC2007}. However,
in this subsection,
we will not restrict ourselves to discussing the small-data
regime in $1+3$ dimensions as in \cite{dC2007}. Instead, we focus on
\emph{general} principles regarding the geometric structure of shock formation for quasilinear wave equations, 
which can be applied to settings beyond the original work \cite{dC2007}, for example in different solution regimes and for 
more general equations; see \cites{jS2014b,dCsM2014,sMpY2014,gHsKjSwW2016} and also a more thorough discussion in the survey article \cite{gHsKjSwW2016} in the case of small compactly supported data on $\mathbb{R}^3$.
In particular, in this subsection, we will suppress discussion of 
the precise estimates that are \emph{specific} to each problem.\footnote{In particular, 
when studying shock formation in a particular solution regime,
it is important to track the ``smallness'' in the problem. 
This plays a crucial role in \cite{dC2007}, 
which specifically considers the small-data regime in $1+3$ dimensions 
for the irrotational relativistic Euler equations. We will completely suppress this discussion in this subsection, but in later subsections, we will emphasize the importance of the role of certain kinds of smallness present
in the solution regime that we consider in the present paper.}

In this subsection, we will consider $(1+2)$ dimensional\footnote{The work \cite{dC2007} was carried out in $(1+3)$ dimensions. However, since the rest of our present paper is in $(1+2)$ dimensions, we will discuss the ideas of \cite{dC2007} as adapted to that case instead. Notice that the $(1+2)$ dimensional case already requires almost all of the new ideas introduced in \cite{dC2007}, with the exception of top-order elliptic estimates for the eikonal function.} 
covariant (see Footnote~\ref{FN:COVWAVEOPARBITRARYCOORDS}) quasilinear wave equations of the form\footnote{In \cite{dC2007},
the equations were 
in fact a subclass of equations, derivable from a Lagrangian, which take the form
$(g^{-1})^{\alpha \beta}(\partial \Phi) \partial_{\alpha} \partial_{\beta} \Phi = 0$.
In particular, the metric depends on the first derivative of the unknown.
However, by differentiation, the equation can be transformed into a system of scalar equations of type
\eqref{waveeqn.form1} that can be studied using essentially the same techniques
needed for proving shock formation in solutions to \eqref{waveeqn.form1};
see the survey article \cite{gHsKjSwW2016} for further discussion.
For this reason, we focus here on equation \eqref{waveeqn.form1}.}
\begin{equation}\label{waveeqn.form1}
\square_{g(\Psi)}\Psi=0
\end{equation}
for real-valued scalar functions $\Psi$ with regular initial data. 
By assumption, the metric $g(\Psi)$ is a Lorentzian metric 
that we will refer to as the acoustic metric since, in the context of the Euler equations,
the wave equations correspond to the propagation of sound waves.
As in Subsect.~\ref{SS:MODELPROBLEM},
we assume that the Cartesian components $g_{\alpha \beta}$ 
are explicit smooth functions of $\Psi$.
We require the nonlinearity in \eqref{waveeqn.form1} to obey certain conditions
so that a shock can form for appropriate initial conditions; see Footnote~\ref{footnote.GLL}.

\subsubsection{Identification of the blowup-mechanism}
\label{SSS:IDENTIFYINGBLOWUPMECHANISM}
In the work \cite{dC2007},  
Christodoulou studied the formation of shocks
by introducing a geometric framework 
tied to an eikonal function $u$, which is a solution to the
eikonal equation (a hyperbolic PDE)
\begin{align} \label{E:FIRSTEIK}
	(g^{-1})^{\alpha \beta}(\Psi) \partial_{\alpha} u \partial_{\beta} u 
	& = 0,
	& \partial_t u & > 0,
\end{align}
supplemented with appropriate initial conditions.
The level sets of $u$ are null hypersurfaces 
(also known as characteristics)
relative to $g(\Psi)$, which
we denoted above by $\mathcal{P}_u$.
The characteristics provide a foliation of the spacetime that is essential for understanding the shock. 
The most important quantity in the study of shock formation is the \emph{inverse foliation density}
$\upmu > 0$, defined as
\begin{align} \label{E:FIRSTMU}
	\upmu 
	& := - \frac{1}{(g^{-1})^{\alpha \beta}(\Psi) \partial_{\alpha} t \partial_{\beta} u},
\end{align}
where $t$ is the Cartesian time coordinate.
$1/\upmu$ is a measure, relative to the constant time slices $\Sigma_t$, of how ``densely packed'' the 
characteristics are. That is, $\upmu = 0$ corresponds to infinite density, the intersection of characteristics,
and the formation of a shock.
One of the key features of Christodoulou's work \cite{dC2007} 
is his proof that for small data, the only possible blowup in a certain solution regime is the formation of shocks. 
That is, he showed that the regularity of the solution is completely 
determined by $\upmu$ and that, for a class of data, one has very good control on how 
$\upmu \to 0$.
In particular, he proved the following facts for solutions generated by an open 
set\footnote{By open, we mean relative to a high-order Sobolev topology.} of data.
\begin{itemize}
\item The solution remains regular at time $t$ if $\upmu_{\star}(t) > 0$,
where $\upmu_{\star}(t) := \inf_{\Sigma_t} \upmu$.
\item Given appropriate initial conditions (consistent with the assumptions for the solution regime), 
$\upmu_{\star} \to 0$ in finite time.
\item It can be justified\footnote{This justification of course relies on a full bootstrap argument, for which the bounds for $\upmu_{\star}$ have to be obtained simultaneously with all the other estimates.} that $\upmu_{\star}$ approaches $0$ \emph{linearly},
	a fact which turns out to be crucial for deriving energy estimates.
\end{itemize}
The blowup-mechanism described above already suggests that the estimates proven for the solutions must take into account appropriate weights of $\upmu$ and $\upmu_\star$. The precise estimates, however, requires further geometric inputs, which we will explain in the next subsubsection.

\subsubsection{The geometric coordinates and the geometric vectorfields}\label{sec.ideas.geo}

A second key feature of Christodoulou's proof, which was also found in Alinhac's works
\cites{sA1995,sA1999a,sA1999b},
is that relative to a \emph{geometric coordinate system} $(t,u,\vartheta)$,
the solution and its low-order partial derivatives remain bounded.
That is, \emph{relative to the geometric coordinates, one does not see the shock ``singularity.''} 
This suggests the main paradigm for approaching the problem:
to the extent possible, prove ``long-time-existence-type'' estimates for the solution relative to the geometric
coordinates and then recover the formation of the shock singularity as a degeneration between the geometric 
coordinates and the Cartesian ones.
Above, $t$ is the Cartesian time coordinate, $u$ is the eikonal function, and $\vartheta$ is a geometrically
defined coordinate that satisfies the transport equation 
$(g^{-1})^{\alpha \beta}(\Psi) \partial_{\alpha} u \partial_{\beta} \vartheta = 0$;
we will downplay the role of $\vartheta$ here since it is better to avoid the use
of coordinates in most of the analysis.

It turns out that deriving the regularity of the solution relative to the geometric coordinates
is equivalent to proving that appropriately $\upmu$-rescaled  
derivatives of various quantities remain bounded. 
That is, one may insert factors of $\upmu$ into various estimates in such a way that 
the vanishing of $\upmu$ exactly compensates for the singularity. 
One might say that many quantities featured in the problem
``blow up like $1/\upmu$.''
More specifically, the tangential (to the characteristics) derivatives of $\Psi$ remain bounded 
\emph{without any factor of $\upmu$} while for the transversal derivative $\Radunit$ (see the 
next paragraph for further discussions), 
the $\Rad :=\upmu \Radunit$ derivatives of $\Psi$ remain bounded. 
Furthermore, it was shown that 
$| \Rad \Psi|$ is bounded from below, strictly away from $0$, when $\upmu$ becomes $0$.
Hence, at those points, $\Radunit \Psi$ blows up and the solution cannot be extended classically.

To prove that the above picture regarding shock formation holds, 
Christodoulou introduced an extensive geometric setup, 
tied to the eikonal function,
which we now adapt to the context of the present article: 
the case of two space dimensions for solutions with approximate plane symmetry. 
In addition to the geometric coordinates described above, he
also introduced geometric vectorfields $\Lunit$, $\GeoAng$ and $\Rad$ 
adapted to the characteristics.
$\Lunit$ is defined to be tangential to the null generators of the $\mathcal{P}_u$,
normalized such that $\Lunit t=1$. 
Specifically, we have
$\Lunit^{\alpha} = - \upmu (g^{-1})^{\alpha \beta}(\Psi) \partial_{\beta} u$
and moreover, 
$
\displaystyle
\Lunit = \frac{\partial}{\partial t}
$
relative to the geometric coordinates.
Let $\ell_{t,u}$ be the intersections\footnote{Note that $\vartheta$ as defined above is a local coordinate on $\ell_{t,u}$.}
$\Sigma_t \cap \mathcal{P}_u$.
Then $\GeoAng$ is the 
$g$-orthogonal
projection of\footnote{Recall that $\partial_2$ is a Cartesian coordinate partial derivative vectorfield.}
$\rd_2$ to $\ell_{t,u}$. The vectorfield $\GeoAng$ is a replacement\footnote{It turns out that $\GeoAng$
has better regularity properties 
than
$
\displaystyle
\frac{\partial}{\partial \vartheta},
$
which are essential for closing the energy estimates.} for the geometric
coordinate partial derivative vectorfield
$
\displaystyle
\frac{\partial}{\partial \vartheta}.
$
Finally, define $\Rad$ to be tangential to $\Sigma_t$
and $g$-orthogonal to $\ell_{t,u}$, normalized such that $\Rad u=1$.
That is, 
$
\displaystyle
\Rad 
= \frac{\partial}{\partial u}
$
plus a small error vectorfield that is tangent to  $\ell_{t,u}$.
Importantly, $\Rad$ becomes degenerate (with respect to the Cartesian coordinate vectorfields) as $\upmu\to 0$. 
That is, the Cartesian components $\Rad^{\alpha}$ vanish precisely at the points where $\upmu = 0$.
On the other hand, the vectorfield $\Radunit=\upmu^{-1}\Rad$ remains non-degenerate,
all the way up to the shock. See Figure \ref{F:FRAME} for a depiction of these vectorfields.

\begin{center}
\begin{overpic}[scale=.35]{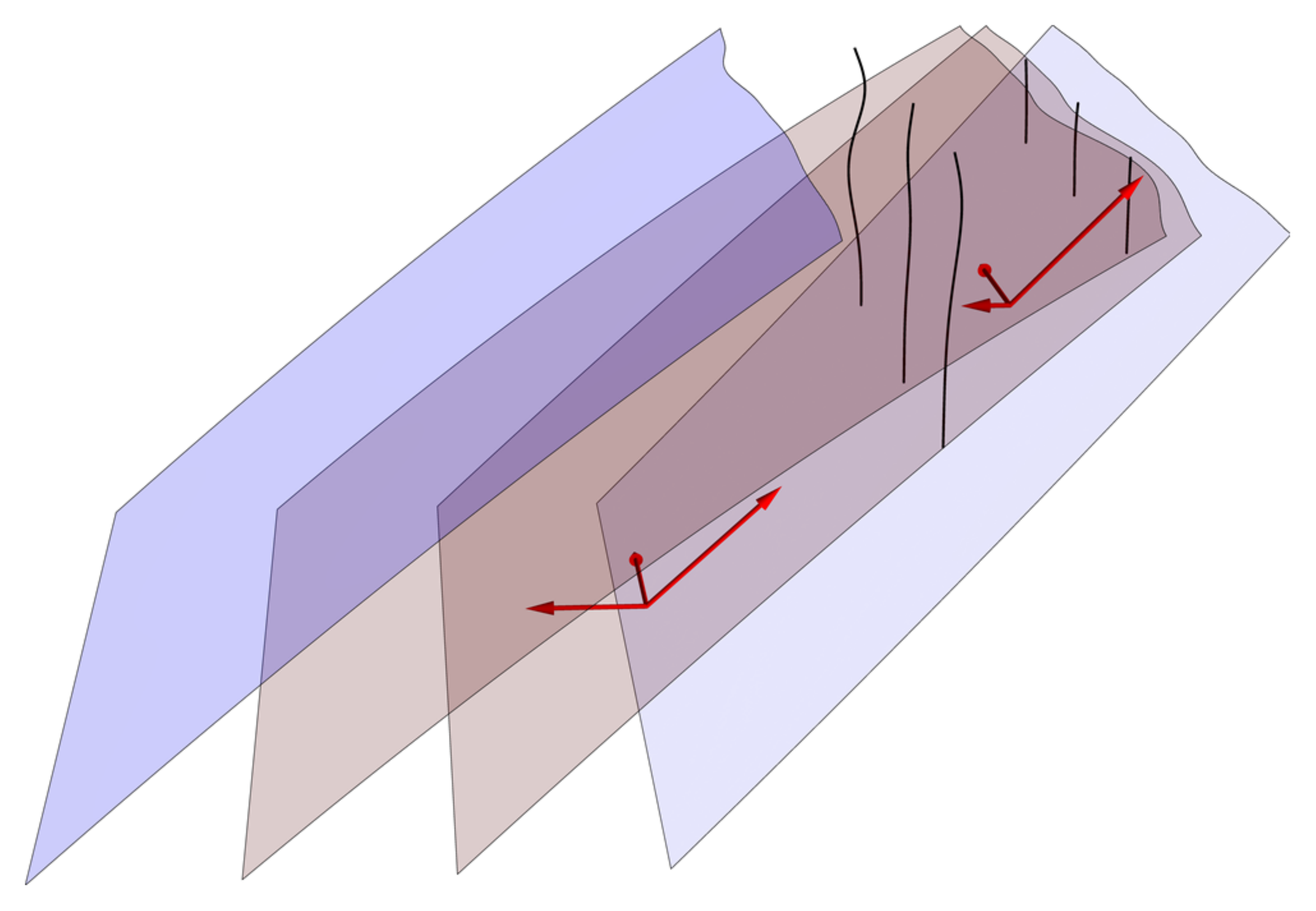} 
\put (82,48) {\large$\displaystyle \Lunit$}
\put (72.5,41) {\large$\displaystyle \Rad$}
\put (73.5,50) {\large$\displaystyle Y$}
\put (54,24) {\large$\displaystyle \Lunit$}
\put (39.5,18) {\large$\displaystyle \Rad$}
\put (47,28) {\large$\displaystyle Y$}
\put (51,13) {\large$\displaystyle \mathcal{P}_0^t$}
\put (37,13) {\large$\displaystyle \mathcal{P}_u^t$}
\put (7,13) {\large$\displaystyle \mathcal{P}_1^t$}
\put (22,26) {\large$\displaystyle \upmu \approx 1$}
\put (70,66) {\large$\displaystyle \upmu \ \mbox{\upshape small}$}
%
%
\end{overpic}
\captionof{figure}{The dynamic vectorfield frame at two distinct points in 
$\mathcal{P}_u^t$ and the integral curves of $\Transport$}
\label{F:FRAME}
\end{center}

Once these geometric vectorfields have been defined, the problem can be reduced to the following steps:
\begin{enumerate}
\item Prove that $\Psi$ and its lower-order derivatives with respect to the geometric vectorfields $\Lunit$, $\GeoAng$ and $\Rad$ are appropriately bounded, with estimates that are independent of how small $\upmu$ is.
\item Prove that the higher-order derivatives of $\Psi$ with respect to the geometric vectorfields 
are not too singular\footnote{As it turns out, the scheme in \cite{dC2007} does not show that the high-order 
derivatives of $\Psi$ with respect to the geometric vectorfields
are bounded. Of course, as we will explain in great detail below, the possible blowup of the solution's high-order derivatives is the source of many difficulties in the problem.} in terms of $\upmu_\star$. Moreover, show that these not-too-singular estimates can be used to derive the non-singular estimates for the lower-order derivatives of $\Psi$ (as described in point (1) above).
\item Justify the transport equation\footnote{In order to guarantee that shock forms, we need $G_{LL}\neq 0$. This can be viewed as a condition on the Cartesian components
$g_{\alpha \beta}$, viewed as a function of $\Psi$. \label{footnote.GLL}} 
(where $G_{LL}$ depends on $g$ and $\Psi$
(cf.\ Def.~\ref{D:GFRAMEANDHFRAMEARRAYS})
and the extra terms are small error terms by Step (1))
$$\Lunit \upmu = \frac{1}{2} G_{LL}\Rad\Psi +\dots $$
and prove that $G_{\Lunit \Lunit} \Rad \Psi$ 
can be precisely controlled in terms of its initial value. Hence, 
under appropriate negativity
assumptions on $G_{\Lunit \Lunit }\Rad\Psi|_{t=0}$,
one can guarantee that $\upmu$ approaches $0$ within the time for which the solution is controlled.
\end{enumerate}

\subsubsection{Degenerate energy estimates and the coercive spacetime bulk term}

In carrying out Step (1) of Subsubsect.~\ref{sec.ideas.geo},
the most crucial estimates are of course $L^2$-based energy estimates. 
It turns out that in order to handle the shock, one needs to incorporate degenerate weights in the energies. To obtain 
suitable 
degenerate energy estimates, 
we apply the vectorfield multiplier method 
with the help of the energy-momentum tensor (see \eqref{E:ENERGYMOMENTUMTENSOR}) and the vectorfield multiplier 
$\Mult = (1+2\upmu)L+2\Rad$.
$\Mult$ has the property that it becomes null and tangential to the 
characteristics $\mathcal{P}_u$
as $\upmu$ vanishes. Moreover, the degeneration is chosen precisely so that the energy controls the following quantities 
$\Sigma_t$
hypersurfaces
(truncated at eikonal function value $u$):
\begin{equation}\label{E.intro}
\int_{\Sigma_t^u} \left((\Rad\Psi)^2+\upmu\left((\Lunit \Psi)^2+(\GeoAng \Psi)^2\right)\right) \,d\vartheta\, du.
\end{equation}
In particular, only the estimate for $\Rad\Psi$ is non-degenerate,
by which we mean the energy becomes very weak in $\Lunit \Psi$ and $\GeoAng \Psi$ along
$\Sigma_t$ when $\upmu$ is small. On the other hand, the energy identities also yield control over
the following quantities on the characteristics $\mathcal P_u^t$ (truncated at time $t$):
\begin{equation}\label{F.intro}
\int_{\mathcal P_u^t} \left((\Lunit \Psi)^2+\upmu(\GeoAng \Psi)^2\right) \,d\vartheta\, dt,
\end{equation}
i.e., one obtains a non-degenerate control for $\Lunit \Psi$ if one considers the energy flux on constant-$u$ hypersurfaces. Notice that in both \eqref{E.intro} and \eqref{F.intro}, the control for $ \GeoAng\Psi$ is degenerate.

Naively, one might expect that in deriving energy estimates, one encounters
terms that are not controllable by the energy itself.
This is because proving degenerate energy estimates corresponds to putting $\upmu$ weights in the ``standard'' energy estimates, and the weights are differentiated during integration by parts. 
If the weights were small with large derivatives of an unfavorable sign, then 
this would lead to potentially insurmountable obstacles to closing the estimates.
However, it turns out that by obtaining detailed information about the way
that $\upmu$ behaves along the integral curves of $\Lunit$, 
one can suitably control the geometric derivatives of the $\upmu$ weights
in the energy. Moreover, 
as was first observed by Christodoulou in his work \cite{dC2007},
one of the spacetime terms in the energy identities
in fact has a \emph{good sign} and is bounded below by
\begin{equation}\label{K.intro}
\int_{\mathcal M_{t,u}} [L\upmu]_- (\GeoAng \Psi)^2 \,d\vartheta\, du\, dt.
\end{equation}
It can be proven that when $\upmu$ is sufficiently small, 
then the negative part $[L\upmu]_-$ is \emph{bounded below} and the integrated term above is non-degenerate and coercive. That is, one can quantify the following heuristic statement
for the solution regime under consideration: 
the only way that $\upmu$ can become small is for $\Lunit \upmu$ to be sufficiently negative.
It is this crucial observation that allows $\GeoAng\Psi$ to be controlled without degeneration.

One must also obtain similar energy estimates for the
higher-order derivatives of $\Psi$.
In order to prove estimates consistent with the expected shock formation picture, 
only the geometric vectorfields can be used as commutators to derive higher-order estimates;
the Cartesian coordinate partial derivative vectorfields would generate uncontrollable error terms
if they were used to commute the wave equation
since they are generally transversal to the $\mathcal{P}_u$ and are not $\upmu$-weighted.
The main technical difficulty that one encounters is that
some of the commutator error terms are exceptionally difficulty to control.
The reason is that the commutator terms depend on the derivatives of the vectorfields
and thus, in view of their connection to the characteristics, on the derivatives
of the eikonal function.  
As we describe in the next subsubsection, it
turns out that one must work hard to avoid losing
derivatives in the most difficult of these terms and, crucially, that
avoiding the derivative loss comes with a price: it introduces a dangerous factor
of $1/\upmu$ into the energy identities, which leads to energy estimates
that are allowed to blow up in terms of powers of $\upmu_{\star}^{-1}$ as the shock forms.

\subsubsection{Top-order estimates for the eikonal function}\label{eikonal.wave.intro}
While the use of tensorfields adapted the characteristics, especially geometric commutation vectorfields,
is necessary to prove shock formation, 
a naive implementation of this framework leads
to a loss of derivatives that threatens to obstruct the closure of the energy estimates.
The difficulty is that 
the commutator of\footnote{It turns out that in our proof, we must commute the weighted operator 
$\upmu \square_g$ in order to avoid generating uncontrollable error terms.
\label{FN:MUWEIGHTEDBOX}} 
$\upmu \square_g$ and the geometric vectorfields generates error terms that depend on the third derivatives of the eikonal function $u$, 
which, as is suggested by the eikonal equation \eqref{E:FIRSTEIK},
can be controlled only by obtaining control over
three derivatives of $\Psi$. On the other hand, after one commutation 
of the wave equation $\upmu \square_g \Psi = 0$
with the geometric vectorfields, only two derivatives of $\Psi$ can be estimated. Nevertheless, as is known since the works\footnote{The work \cite{dCsK1993} exploited this gain of a
derivative in the specific case of the Einstein vacuum equations, 
for which this structure is more easily seen.
Nevertheless, the ideas in \cite{dCsK1993}
already serve as a blueprint for 
gaining the derivative in the context of more
general quasilinear wave equations.} \cites{dCsK1993,sKiR2003}, one can exploit the fact that the 
Cartesian components of the metric $g(\Psi)$ 
also satisfy a wave equation\footnote{This claim is a simple consequence of the chain rule applied to the component
functions $g_{\alpha \beta}(\Psi)$.}
to gain a derivative for certain special combinations of third derivatives of $u$ and second
derivatives of $\Psi$. The gain is tensorial in nature, and it is a happy fact that
the vectorfields $\Lunit$, $\Rad$, and $\GeoAng$ generate only commutation
error terms featuring those special combinations.

One example (in fact, the most important example in the problem) 
of a controllable error term depending on the eikonal function
is the null mean curvature $\mytr \upchi$ of the characteristics $\mathcal{P}_u$.
It turns out that after commuting the wave equation one time with geometric
vectorfields, one encounters the first derivatives of $\mytr \upchi$, 
which, 
as we alluded to above, 
we must carefully treat to avoid losing a derivative.
We now explain how to avoid this derivative loss. For convenience, instead of addressing this difficulty at the level of 
one commutation of the wave equation, we consider the analogous difficulty 
at the level of zero commutations. That is, we explain how to 
control the undifferentiated quantity $\mytr \upchi$ in terms of one derivative of $\Psi$
(which is the allowed regularity for $\Psi$ without commuting).
The source of the difficulty is that $\mytr \upchi$
depends on two derivatives of the eikonal function $u$.
This seems to be incompatible with the available regularity of $\Psi$
since a
general second derivative of $u$
has to be estimated by two derivatives of the metric (and hence two derivatives of $\Psi$). However, $\mytr \upchi$ is a special combination of two derivatives of the eikonal function $u$ and one can ``gain'' in derivatives with the following procedure. First, one derives the following transport equation for $\mytr \upchi$ (well-known in general relativity as the Raychaudhuri equation):
\begin{equation}\label{Ray.eqn.intro}
\upmu \Lunit \mytr \upchi=(\Lunit \upmu)\mytr \upchi-\upmu(\mytr \upchi)^2-\upmu \mbox{\upshape Ric}_{\Lunit \Lunit}, 
\end{equation}
where $\mbox{\upshape Ric}_{\Lunit \Lunit}$ is the $\Lunit \Lunit$-component of the spacetime Ricci curvature of $g$. The second main observation is that $\mbox{\upshape Ric}_{\Lunit \Lunit}$ 
is equal
to a sum of terms controllable by only one derivative of $\Psi$ and a term which can be written as $\Lunit \upchifullmodinhom$, where $\upchifullmodinhom$ can be expressed in terms of at most one derivative of $\Psi$. That this can be achieved crucially depends on the wave equation for $\Psi$. More precisely, using the wave equation $\upmu \square_g \Psi = 0$, one can replace 
the second derivative term\footnote{Here, $\slashed\Delta$ denotes the Laplacian with respect to 
the Riemannian metric induced by $g(\Psi)$ on $\ell_{t,u} = \Sigma_t \cap \mathcal{P}_u$.}
$\upmu\slashed\Delta\Psi$, which appears in the expression of $\upmu \mbox{\upshape Ric}_{LL}$, with $\Lunit(\upmu \Lunit \Psi+2\Rad \Psi)$ (and lower-order terms). As a consequence, instead of directly studying $\mytr \upchi$, we can instead study the
\emph{modified quantity} 
$\upmu\mytr \upchi+\upchifullmodinhom$.
The key point is that the
right hand side of the transport equation $\Lunit(\upmu\mytr \upchi+\upchifullmodinhom)$ now depends on at most one derivative of 
$\Psi$. In total, this procedure allows us to control 
$\mytr \upchi$ using estimates for one derivative of $\Psi$ only, 
which is better than what one would naively expect.

On the other hand, as the shock is approached, this procedure of gaining derivatives is coupled with the difficulty of\footnote{Recall that $\upmu_{\star}$ is the minimum of $\upmu$ on a constant-$t$ hypersurface.} $\upmu_{\star}\to 0$. This is because the modified quantity is $\upmu\mytr \upchi+\upchifullmodinhom$, where $\upchifullmodinhom$ is merely bounded. Hence, to recover estimates for
(the higher-order derivatives of) $\mytr \upchi$ from estimates for (the higher-order derivatives of) the modified quantity, one faces the critical difficulty of a discrepancy factor of $1/\upmu$; this discrepancy is central to most of the difficulties that one faces in closing the problem.

We now illustrate how this difficulty enters into the top-order energy estimates
for $\Psi$ by keeping one of the most significant terms.\footnote{Notice that there are other terms which are of the same strength (from the point of view of the singularity) as the term that is shown.} This leads to an estimate of the 
following form for the top-order energy $\mathbb E_{top}$ on $\Sigma_t^u$:
\begin{equation}\label{sing.top.intro}
\mathbb E_{top}(t,u)\leq \mbox{Data}+C_{fix}\int_{t'=0}^{t'=t} \left(\sup_{\Sigma^u_{t'}} \left|\f{L\upmu}{\upmu}\right|\right)\mathbb E_{top}(t',u)\, dt +\dots,
\end{equation}
where\footnote{In the proof of our main theorem, $C_{fix}$ will in fact be an explicit numerical constant.} $C_{fix}>0$ is a \emph{fixed} constant independent of how many derivatives we choose to be the ``top level''. To proceed, 
one needs very precise estimates for
$\upmu$ and $\Lunit \upmu$ 
and to show that $\upmu_{\star}$ tends to $0$ \emph{linearly}.
This implies,\footnote{We also note that in closing the top-order energy estimates, one must perform
some crucially important integrations-by-parts in time that lead to singular boundary terms that must 
be controlled. We will suppress this technical difficulty here in order to keep the discussion short,
instead referring readers to Subsect.~\ref{SS:ERROINTEGRALSINVOLVINGIBPL} for details regarding this estimate.} 
via a difficult analog of Gronwall's inequality 
that relies on the sharp information for $\upmu_{\star}$,
the estimate\footnote{To obtain some heuristic understanding of
why \eqref{E.sing.top.intro} follows from \eqref{sing.top.intro}, one may 
replace $\upmu$ and $\upmu_{\star}$ with $1-t$ and $\Lunit \upmu$ with $-1$.
Then the standard Gronwall inequality yields \eqref{E.sing.top.intro}
with $c_{fix} = C_{fix}$.
\label{FN:HUERISTICTOPORDERBLOWUP}}
\begin{equation}\label{E.sing.top.intro}
\mathbb E_{top}(t,u)\leq \mbox{Data}\times \upmu_{\star}^{-c_{fix}}(t,u),
\end{equation}
for some \emph{universal constant $c_{fix}>0$ that is independent 
of the structure of the nonlinearities.}

\subsubsection{Energy hierarchy and the descent scheme}\label{sec.descent.intro}

In the previous subsubsection, we saw that in order not to lose derivatives, 
the top level energy estimates must degenerate in terms of $\upmu_{\star}$. To finish the argument, Christodoulou introduced a \emph{descent scheme} 
in which he showed that for every order below the top-order, 
the degeneration can be improved by a fixed amount. 
In particular, at some sufficiently low-order of derivatives, the energy can be shown to be bounded. 
This then also yields, by a geometric Sobolev embedding estimate,
the necessary low-order $L^\infty$ estimates  
that allow the argument to be closed.

Let us describe\footnote{We do not directly describe the numerology of \cite{dC2007} here as it is slightly different from that of 
the present paper and doing so might create some confusion. On the other hand, the main ideas can be traced back to \cite{dC2007}.}
the relevant numerology in the adaption of \cite{dC2007} to the present paper. One proves estimates of the type
\begin{subequations}
	\begin{align}
		\sqrt{\mathbb E_{15+K}}(t,u)
				& \leq C \mathring{\upepsilon} \upmu_{\star}^{-(K+.9)}(t,u),
			&& (0 \leq K \leq 5), \label{E.descent.intro}\\
		\sqrt{\mathbb E_{N}}(t,u)
				& \leq C \mathring{\upepsilon},
		&&	(1\leq N\leq 14),
	\end{align}
	\end{subequations}
where $\mathbb E_N$(t,u) denotes the energy on $\Sigma_t^u$
after $N$ commutations. Here, $\mathbb E_{20}$ is the top-order energy,
and it can be controlled using the approach described in the previous subsubsection. 
On the other hand, when controlling $\mathbb E_{19}$, one can control the highest order 
(that is, $19^{th}$)
derivative of $\mytr \upchi$ appearing in the energy estimates by
$\mathbb E_{20}$ instead of $\mathbb E_{19}$. Put differently,
below-top-order, one can simply allow the loss of a derivative and avoid using the modified version of $\mytr \upchi$.
In this way, one avoids introducing an explicit singular factor of $1/\upmu$ into the 
below-top-order energy estimates, at the expense of introducing a coupling to the energy at one higher level.
A key point, 
which was exploited by Christodoulou in controlling the energies at all derivative levels,
is that
$\upmu_{\star}\to 0$ at worst linearly (with precise estimates). This in particular allows one to show that
every integration in $t$ reduces the strength of the singularity by a power of $\upmu_\star$.
These ``descent estimates,'' though exceptionally technical to implement, are nothing other than a ``quasilinear version''
of the estimate $\int_{s=0}^t s^{-b} \, ds \lesssim t^{1 - b}$ (for $b>1$), 
where $s=0$ represents the ``vanishing'' of $\upmu_\star$.

To be more concrete, let us consider the most difficult
inhomogeneous term on the right-hand side of the equation $\upmu\square_{g(\Psi)}\GeoAng^{19}\Psi = \cdots$,
which modulo bounded factors is the term $(\Lunit \GeoAng^{19}\Psi) (\GeoAng^{19}\mytr \upchi)$.
Since the control for $\Lunit \GeoAng^{19} \Psi$
on a constant $t$-hypersurface in $\mathbb E_{19}^{\frac{1}{2}}$ is $\upmu^{\f12}$ degenerate (recall \eqref{E.intro}),
these considerations lead to the estimate
$$\mathbb E_{19}(t,u)\leq \mbox{Data}+ \int_0^t \upmu_{\star}^{-\frac{1}{2}} \mathbb E_{19}^{\frac{1}{2}}(t',u) \| \GeoAng^{19}\mytr \upchi \|_{L^2(\Sigma_{t'}^u)} dt'+\dots$$
To estimate $\GeoAng^{19} \mytr \upchi$, 
one can again use \eqref{Ray.eqn.intro}, but this time 
directly controlling $\mbox{\upshape Ric}_{LL}$. 
Since $\mathbb E_{20}^{\f12}$ again has a $\upmu^{\f12}$ degeneration for the $\GeoAng$ derivative, integrating \eqref{Ray.eqn.intro} yields 
\begin{equation}\label{trch.no.gain.intro}
\| \GeoAng^{19} \mytr \upchi \|_{L^2(\Sigma_{t'}^u)}\leq \mbox{Data}+\int_{t''=0}^{t''=t'} \upmu_{\star}^{-\frac{1}{2}}\mathbb E_{20}^{\f12}(t'',u)\, dt''+\dots
\end{equation}
Substituting this back into the estimate for $\mathbb E_{19}$ gives
\begin{equation}\label{E19.intro}
\mathbb E_{19}(t,u)\leq \mbox{Data}+ \int_0^t \upmu_{\star}^{-\frac{1}{2}} \mathbb E_{19}^{\frac{1}{2}}(t',u) \int_{t''=0}^{t''=t'} \upmu_{\star}^{-\frac{1}{2}} \mathbb E_{20}^{\f12}(t'',u)\, dt'' dt'+\dots
\end{equation}
Now since $\upmu_{\star}$ tends to $0$ at worst linearly, we have the following estimate, which we alluded to above:
for every $b>1$, we have
\begin{equation}\label{int.lemma.intro}
\int_{t'=0}^{t'=t} \upmu_{\star}^{-b}(t')\, dt\leq C\upmu_{\star}^{-b+1}(t).
\end{equation}
Therefore, \eqref{E19.intro} is indeed consistent\footnote{To actually close the estimates,
one needs to derive a Gronwall estimate for a coupled system featuring
$\mathbb E_{19}$ and $\mathbb E_{20}$. We refer the readers to
Subsect.~\ref{SS:PROOFOFPROPMAINAPRIORIENERGY} for the relevant details
in the context of the present article.} 
with the reduced blowup-rate for $\mathbb E_{19}$ compared to $\mathbb E_{20}$,
as stated in \eqref{E.descent.intro}.
One can continue the descent and show
that the blowup-rates for the energies continues to improve
as the number of derivatives is reduced, 
until one actually obtains boundedness of the lower-order energies.

As we can see from the above scheme, the number of derivatives that is needed for this descent scheme depends on 
the strength of the top-order singularity, represented by the constant $c_{fix}$ in \eqref{E.sing.top.intro}.
In turn, $c_{fix}$ depends on the constant in $C_{fix}$ in \eqref{sing.top.intro}.

\subsubsection{Formation of shocks}
Once one closes all the estimates, the formation of shocks follows easily. Indeed, with the estimates at hand, it is easy to conclude the solution remains regular relative to both\footnote{Indeed, the change of variables map from
geometric to Cartesian coordinates is a diffeomorphism when $\upmu > 0$.} 
geometric and Cartesian coordinates
as long as 
$\upmu>0$ and that $\upmu\to 0$ corresponds indeed to a shock. Moreover, 
the non-degenerate low-level energy estimates imply, via Sobolev embedding, non-degenerate low-level $L^{\infty}$ estimates
that lead to
$$\Lunit \upmu=\frac{1}{2} G_{\Lunit \Lunit} \Rad \Psi+\dots,$$
where $\dots$ denotes small error terms
and also that $G_{LL} \Rad\Psi$ is essentially transported along the integral curves of $\Lunit$. 
Therefore with appropriate negativity assumptions on $G_{\Lunit \Lunit} \Rad \Psi|_{t=0}$,
it is easy to prove that $\upmu$ goes to $0$ in finite time.

\subsection{Review of the stability of shock formation for nearly simple outgoing plane symmetric solutions to quasilinear wave equations}
\label{sec.ideas.review.plane}

Together with Holzegel and Wong, 
we proved \cite{jSgHjLwW2016} stable shock formation for ``nearly simple outgoing plane symmetric'' solutions to\footnote{In fact, similar methods could be used to
show that the solutions are stable under non-symmetric perturbations in \emph{three} spatial dimensions;
see the discussion in \cite{jSgHjLwW2016}.} 
a class of quasilinear wave equations in two spatial dimensions.
In this paper, 
we study a similar regime of nearly simple outgoing plane symmetric solutions.
More precisely, we extend the results of \cite{jSgHjLwW2016} 
to the compressible Euler's equations \emph{without the irrotationality assumption}.\footnote{See Subsect.~\ref{sec.ideas.Euler} 
for discussion on the relation between the class of equations discussed here and the compressible Euler equations.}

We now describe the case of \emph{exact} simple outgoing plane symmetric solutions. We
first recall that $(1+1)$-dimensional quasilinear wave equations can be greatly simplified using the conformal invariance of $\square_g$ and the fact that $(1+1)$-dimensional Lorentzian manifolds are (locally) conformally flat. Indeed, defining appropriate null functions $u$ and $w$ satisfying the eikonal equation
$$(g^{-1})^{\alpha \beta}\rd_{\alpha} u \rd_{\beta} u = (g^{-1})^{\alpha \beta}\rd_{\alpha} w \rd_{\beta} w = 0,$$
it is easy to show the quasilinear wave equation $\square_{g(\Psi)}\Psi=0$ is equivalent to 
$$\rd_w \rd_u\Psi=0.$$
Note that the equation is still quasilinear as $u$ and $w$ depend on $\Psi$. 
We say that a solution $\Psi$ is a \emph{simple outgoing}\footnote{We use the term ``outgoing'' to mean that the solution travels towards the ``right.'' That is, we have used the convention that initially $\f{\rd w}{\rd x}>0$. Of course, the restriction to outgoing waves is merely for notational convenience, as the analysis remains identical if we instead consider ``incoming'' solutions.} if $\rd_w\Psi=0$. When expressed in terms of $\Psi$ alone, $\rd_w\Psi=0$ can be viewed as a Burger's-type equation.

In \cite{jSgHjLwW2016}, the authors studied 
$(1+2)$-dimensional
quasilinear wave equations. 
The equations admit plane symmetric solutions which do not depend on the Cartesian spatial coordinate $x^2$.
Analogous to the $(1+1)$-dimensional case, they can be written as
$$\rd_w \rd_u\Psi=\mathcal N,$$
where 
$\mathcal N = \mathcal{N}(\Psi, \partial \Psi) \rd_u \Psi \cdot \rd_w\Psi$
is a null form relative to $g$, and 
a solution is said to be
\emph{simple outgoing} if $\rd_w\Psi=0$. It is not difficult to see that under a condition of genuine nonlinearity,
there exist 
simple outgoing solutions for which shocks form in finite time. 
The main result of \cite{jSgHjLwW2016} is that a subclass of such 
shock-forming solutions to equation \eqref{waveeqn.form1}
is stable under \emph{non-symmetric} perturbations. Proving this result requires, in addition to the ideas of \cite{dC2007} described in the previous subsection, a method to propagate smallness parameters relevant to this solution regime.\footnote{In particular, this is in contrast to \cite{dC2007}, where dispersion was crucially used to propagate smallness. The lack of dispersion in \cite{jSgHjLwW2016} requires the introduction of the $\TranminusdatasizeWithFactor$-$\mathring{\upepsilon}$ size hierarchy of the initial data (to be described below), but it turns out that to propagate that smallness is slightly less involved than that in \cite{dC2007}.} 

In order to achieve this, 
the authors introduced
the parameters\footnote{These parameters are closely related to those introduced in Subsect.~\ref{sec.data}. See Subsect.~\ref{sec.ideas.Euler} for further discussion.} $\TranminusdatasizeWithFactor$, $\mathring{\updelta}$ and $\mathring{\upepsilon}$ to describe the relative sizes of the derivatives of $\Psi$. Here, $\TranminusdatasizeWithFactor$ and $\mathring{\updelta}$ are not necessarily small: $\mathring{\updelta}$ describes the size of 
the \emph{transversal} (to the characteristics) derivatives
of the data of $\Psi$ and $\TranminusdatasizeWithFactor^{-1}$
is the ``expected the blow up time'', which depends on the first transversal derivative of $\Psi$ at time $0$
(compare with Definition \ref{D:CRITICALBLOWUPTIMEFACTOR}). On the other hand, $\mathring{\upepsilon}$, which, roughly speaking, 
describes the size of the initial $L^\infty$ norm of $\Psi$ as well as its initial outgoing derivatives 
and its derivatives in the direction $\rd_2$, is required to be small compared to $\TranminusdatasizeWithFactor$ and $\mathring{\updelta}^{-1}$.
In order to more precisely
describe the smallness of the initial data, 
we will again use the geometric vectorfields $\Lunit$, $\GeoAng$ and $\Rad$ described in the previous subsection.\footnote{Here, one could think of $\Lunit$ as an analogue of $\rd_w$ in the exact plane symmetric case and $\GeoAng$ as in the direction $\rd_2$.} 
In \cite{jSgHjLwW2016}, 
the geometric derivatives of $\Psi$ at time $0$
are required to be $\mathring{\upepsilon}$-small whenever 
\emph{at least one} of the differentiations is in the direction of $\Lunit$ or $\GeoAng$. 
It is straightforward to show that this initial smallness follows whenever the
initial data are $\mathring{\upepsilon}$-perturbations of simple outgoing plane symmetric solutions.

In order to control the solution, 
we do not need to explicitly subtract the simple outgoing plane symmetric solution from the full nonlinear solution. Instead, we show that the solution remains nearly plane symmetric and nearly simple outgoing up to the 
time of first shock formation
in the sense that 
\begin{quote}
\begin{center}
For the derivatives of $\Psi$ with respect to the
geometric vectorfields, if \emph{at least one} of the vectorfields is $\Lunit$ or $\GeoAng$, 
then the quantity in an appropriate norm is $\mathcal{O}(\mathring{\upepsilon})$ small, 
all the way up to the shock.\footnote{We recall from the previous subsection that some higher-order norms are allowed to blow up. Therefore, the $\mathcal{O}(\mathring{\upepsilon})$ ``smallness'' of the higher-order energies
must be carefully interpreted as smallness relative to singular norms.}
\end{center}
\end{quote}
In other words, the smallness of the $\mathcal P_u^t$-tangential derivatives of $\Psi$, which is originally assumed for the initial data, is propagated by the flow. Notice that in this process, not only do we use the geometric vectorfields $\Lunit$, 
$\GeoAng$ and $\Rad$ to capture the formation of shocks, we also use them to track the smallness in the problem. 
The following geometric and analytic properties are crucial in order to achieve this:
\begin{enumerate}
\item (Commutation properties of the geometric vectorfields) We 
stress that we have crucially used the property that even if $\Psi$
is $\Rad$-differentiated, as long  
it is also hit with one or more $\Lunit$ or
$\GeoAng$ derivatives (for instance for the quantities $\Lunit \Rad \Psi$, $\Rad L\Lunit \Psi$, etc), 
then the quantity is still $\mathcal{O}(\mathring{\upepsilon})$ small. That this holds of course relies on good commutation properties of the geometric vectorfields, in particular that the commutator of any two of 
$\lbrace \Lunit, \Rad, \GeoAng \rbrace$
is \emph{tangential to the characteristics $\mathcal{P}_u$} 
(in fact, the commutators may be seen to be tangent to 
$\ell_{t,u} = \Sigma_t \cap \mathcal{P}_u$)!
\item (Null structure of the nonlinear terms when decomposed with respect to geometric vectorfields) The wave equation  \eqref{waveeqn.form1}
is equivalent to (see Proposition \ref{P:GEOMETRICWAVEOPERATORFRAMEDECOMPOSED})
\begin{equation}\label{transport.intro.sch}
-L(\upmu L \Psi + 2\Rad \Psi)+\upmu\angLap \Psi=\mathcal N,
\end{equation}
where $\mathcal N$ denotes nonlinear terms with 
\emph{at most one factor transversal}
to $\mathcal{P}_u$, that is, with at most one factor equal to $\Rad\Psi$.
Consequently, under appropriate
bootstrap assumptions, 
the term $\mathcal N$ can be shown to be $\mathcal{O}(\mathring{\upepsilon})$-small.\footnote{The 
implicit constants are allowed to depend on
$\mathring{\updelta}$.} Notice that this smallness partly comes from the geometry associated to the problem. For instance, one of the terms in $\mathcal N$ is $\mytr \upchi \Rad\Psi$
(where as before, $\mytr \upchi$ is the null mean curvature of the $\mathcal{P}_u$).
It was shown
that $\mytr \upchi$ is $\mathcal{O}(\mathring{\upepsilon})$-small, 
which is a consequence of approximate plane symmetry
the solution.\footnote{Notice also that this smallness is tied to our foliation of spacetime by 
the nearly flat characteristics $\mathcal{P}_u$.
One might say that we made an ``educated'' guess about
how to construct a foliation that allows us to propagate the smallness.}
\item (Commutation properties between the geometric vectorfields and $\upmu \square_g$)
)In order to prove size $\mathcal{O}(\mathring{\upepsilon})$ estimates for the (higher-order) $\Lunit$ and $\GeoAng$ derivatives of $\Psi$, we rely on the fact that the commutators terms\footnote{The factor of $\upmu$ generates important cancellations.} 
$[\upmu \square_g, \Lunit]\Psi$ and $[\upmu \square_g,\GeoAng]\Psi$
are $\mathcal{O}(\mathring{\upepsilon})$-small. Moreover, we also use the fact that
$[\upmu \square_g, \Lunit]\Psi$ and $[\upmu \square_g, \GeoAng]\Psi$ do not generate 
$\Rad\Rad \Psi$ terms.\footnote{Dimensional considerations imply that these terms, if present, would
be multiplied by an uncontrollable factor of $1/\upmu$. Furthermore, the absence of these terms in the commutators is also useful for the higher-order energy estimates; see point (1) below.}
These can be viewed as a consequence of the commutation properties described in the first point above.
\end{enumerate}
Using the above properties, we carry out our estimates as follows:
\begin{enumerate}
\item (Higher-order energy estimates) For the energy estimates, we only use $\Lunit$ and $\GeoAng$ as commutators. We also only carry out the energy estimates after at least one commutation. Notice that this is sufficient from the point of view of regularity since 
$[\upmu \square_g, \Lunit]$ and $[\upmu \square_g, \GeoAng]$
do not generate $\Rad\Rad$ terms! Moreover, 
the energy corresponding to commuting the wave equation with one or more
factors of $\Lunit$ or $\GeoAng$ is initially $O(\mathring{\upepsilon}^2)$-small,
which is convenient\footnote{The energy of the non-commuted equation is lower-bounded by the square of the 
$L^2$ norm of $\Rad\Psi$, which can be of a relatively large size $\mathring{\updelta}^2$.}
for deriving estimates.
\item (Estimates for the eikonal function) The estimates for the eikonal function $u$ up to the highest order are intimately tied to the energy estimates. For example, we show that $\mytr \upchi$
(and its higher-order $\Lunit$ and $\GeoAng$ derivatives) inherits
the $\mathcal{O}(\mathring{\upepsilon})$ smallness from the energy estimates, 
as is expected since the solution is nearly outgoing simple plane symmetric.\footnote{Note that for exact outgoing simple plane symmetric solutions, we have $\mytr \upchi=0$.}
Note that in order to implement steps (1) and (2), it is important
that one can close the energy estimates 
and the estimates for the eikonal function (which, as we described above, are highly coupled!) 
by commuting only with $\Lunit$ and $\GeoAng$.
\item (Lower-order estimates for $\Rad \Psi$, $\Rad \Lunit \Psi$, $\Rad \Rad \Psi$, $\Rad \Rad \Rad \Psi$, etc.) 
Since we only derive
energy estimates after commuting with at least one factor of $\Lunit$ or $\GeoAng$, our energies cannot be
directly combined with Sobolev embedding to yield pointwise control of $\Rad \Psi$.
To obtain pointwise control of $\Rad \Psi$,
we use the wave equation in the form \eqref{transport.intro.sch} 
as a transport equation in the unknown $\Rad\Psi$.
The pointwise estimates for $\Rad \Lunit \Psi$, $\Rad \Rad \Psi$, 
$\Rad \Rad \Rad\ Psi$, etc.\ are obtained in a similar manner,
after commuting the wave equation.
Notice that some of these terms, for instance $\Rad \Psi$, are of relatively large size $\mathring{\updelta}$, but this size can be propagated since the error terms in \eqref{transport.intro.sch} are all of smaller size $\mathcal{O}(\mathring{\upepsilon})$. In other words, in the nonlinear error terms, we never encounter, say, quadratic terms of size $\mathring{\updelta}^2$.
\item (Sharp control of $\upmu$) Using the smallness above, we show that $\upmu$ satisfies the transport equation
$$\Lunit \upmu=\frac{1}{2} G_{LL}\Rad\Psi+\mathcal{O}(\mathring{\upepsilon}),$$
where $G_{LL}$ is a function of $\Psi$ depending on the Cartesian components $g_{\alpha \beta}$.
This equation,
together with a precise estimate for $G_{LL}\Rad\Psi$,
are crucial for obtaining
sharp control\footnote{Let us recall from the previous subsection that we crucially need to show that $\upmu\to 0$ at worst linearly and also to prove that $[\Lunit \upmu]_-$ is bounded from below whenever $\upmu$ is sufficiently small.} of $\upmu$ and $\Lunit \upmu$ and 
for showing that a shock indeed forms in finite time.
\end{enumerate}

\subsection{Nearly simple outgoing plane symmetric solutions to irrotational Euler equations}\label{sec.ideas.Euler}
As was already discussed in the paper \cite{jSgHjLwW2016}
on quasilinear wave equations of the type $\square_{g(\Psi)}\Psi=0$, 
with very few modifications,
the same methods can be used to prove shock formation in solutions to
quasilinear equations of the form 
\begin{equation}\label{waveeqn.form2}
(g^{-1})^{\alpha \beta}(\rd\Phi) \rd_{\alpha} \rd_{\beta} \Phi=0.
\end{equation}
This can be seen by considering the vector 
$\threePsi:= (\Psi_{\nu})_{\nu=0,1,2} := (\rd_\nu \Phi)_{\nu=0,1,2}$, 
differentiating \eqref{waveeqn.form2} 
and deriving the system of equations
\begin{equation}\label{waveeqn.form3}
\square_{g(\threePsi)} \Psi_\nu=\mathcal Q(\rd \threePsi, \rd \threePsi_\nu),
\end{equation}
where $\square_{g(\threePsi)}$ is to be understood as acting on scalar functions and 
the inhomogeneous terms $\mathcal Q$ are quadratic \emph{null forms} (relative to $g$). 
Both the vectorial nature of the unknown 
and the additional nonlinear terms pose almost no additional challenge and \eqref{waveeqn.form3} can be treated with essentially the same methods as the scalar equation \eqref{waveeqn.form1}. 
This in particular crucially relies on the structure of the null forms,
which have only a negligible influence on the solution, all the way up to the shock.

To handle the equation \eqref{waveeqn.form2}, Speck--Holzegel--Luk--Wong \cite{jSgHjLwW2016} 
considered initial data such that each of the components $\Psi_\nu$ (which are viewed as scalar functions)
obey the 
$\mathring{\upepsilon}$-$\mathring{\updelta}$
size estimates as described in Subsect.~\ref{sec.ideas.review.plane}. Notice that,
in view of the identity $\rd_{\alpha} \Psi_{\beta}  =\rd_{\beta} \Psi_{\alpha}$, 
these assumptions imply smallness estimates for the $\Rad$ derivative of certain combinations
of the $\Psi_\nu$.
We note that the result in \cite{jSgHjLwW2016} can be applied to the 
irrotational compressible Euler equations.\footnote{There is an explicit justification of this fact in \cite{jSgHjLwW2016} for the \emph{relativistic} Euler equations. It can easily be seen that this also applies to the non-relativistic case.}
More precisely,
by introducing a potential function $\Phi$ for the flow,
we obtain an equation of the form \eqref{waveeqn.form2}.  

Let us clarify the connection between the Riemann invariants 
(see Subsect.~\ref{sec.data}),
our assumption that we are studying perturbations of simple plane waves,
and the size assumptions
on the potential $\Phi$ described in the previous two paragraphs.
We first note that in the non-relativistic case, we have, relative to the Cartesian
spatial coordinates, $$\partial_i \Phi = - v^i.$$ 
Next, we note that in
one spatial dimension, we have (see \cite{dCsM2014})
\begin{equation}\label{Phi.1D}
\rd_t\Phi- \frac{1}{2}(\partial_1 \Phi)^2=h,
\end{equation}
where $h$ is the enthalpy, defined such that 
$
\displaystyle
dh 
= 
\frac{\Speed^2}{\rho} \, d \rho
=
\Speed^2 \,d \Densrenormalized$. 
The assumption (stated in the previous paragraph)
that $\Lunit \Psi_1$ should initially be of small size $\mathcal{O}(\mathring{\upepsilon})$
(where $\Psi_1 :=\rd_1\Phi$)
can be stated as
$\Lunit_{(plane)} \Psi_1
=\mathcal{O}(\mathring{\upepsilon})$
where, in one spatial dimension,
$\Lunit_{(plane)} = \Lunit = \rd_t+(v+c_s)\rd_1$ is the outgoing null vectorfield with 
$\Lunit_{(plane)} t=1$.
This smallness assumption implies, via \eqref{Phi.1D} 
and the formula
$
\displaystyle
\mathcal{R}_-
=
v^1-\int_{\widetilde{\Densrenormalized}=0}^{\Densrenormalized} \Speed(\widetilde{\Densrenormalized}) \, d \widetilde{\Densrenormalized}
$,
that, at time $0$, we have
$
\displaystyle
-\rd_1 \mathcal R_-
= \rd_1\left(-v_1 + 
	\int_{\widetilde{\Densrenormalized}=0}^{\Densrenormalized} 
		\Speed(\widetilde{\Densrenormalized})
	\, d \widetilde{\Densrenormalized}\right)
= \rd_1^2 \Phi + \Speed \rd_1 \Densrenormalized
= \rd_1^2 \Phi + \Speed^{-1} \partial_1 h
= \rd_1^2 \Phi + \Speed^{-1} \rd_t \partial_1 \Phi - \Speed^{-1} (\rd_1 \Phi) \partial_1^2 \Phi
=  \Speed^{-1} L_{(plane)} \Psi_1
= \mathcal{O}(\mathring{\upepsilon})
$.
Therefore, the smallness assumption implies that
the derivatives of the Riemann invariant $\mathcal R_-$  
are initially small, which is a perturbation of the simple plane wave case
$\mathcal{R}_- \equiv 0$ described in Subsect.~\ref{sec.data}.

Finally, recalling the discussions of the initial data in Subsect.~\ref{sec.data}, we see that the data considered in this paper 
are a generalization of that 
considered in \cite{jSgHjLwW2016} such that in this paper, the \emph{vorticity is not required to vanish identically}.

\subsection{New ideas in the case of non-vanishing vorticity}\label{sec.ideas.vorticity}

We are now ready to
discuss the main new ideas in the present paper, 
which are needed to handle the interaction between the vorticity and the sound waves. 
We recall that we described our assumptions on the initial data in Subsect.~\ref{sec.data}.
In this paper, we also need all of the ideas as described in Subsects.~\ref{Chr.intro} and \ref{sec.ideas.review.plane}, although for the sake of brevity we will often not repeat them. In particular, we mostly
suppress in this subsection the issue of propagating $L^{\infty}$-type
smallness at the low derivative levels; most of the ideas in that regard are similar 
to those discussed in Subsect.~\ref{sec.ideas.review.plane}.
Instead, we focus on the crucially important issue of closing 
the energy estimates.

As we mentioned earlier, the starting point of our proof is 
the following reformulation of the compressible Euler equations,
valid in two spatial dimensions:
\begin{subequations}
	\begin{align}
		\upmu \square_g v^i
		& = - [ia] (\exp \Densrenormalized) \Speed^2 (\upmu \partial_a \Vortrenormalized)
			+ 2 [ia] (\exp \Densrenormalized) \Vortrenormalized (\upmu \Transport v^a)
			+ \upmu \mathscr{Q}^i,
				\label{intro.wave.1} \\
	\upmu \square_g \Densrenormalized
	& = \upmu \mathscr{Q},
		\label{intro.wave.2} \\
	\upmu \Transport \Vortrenormalized
	& = 0.
	\label{intro.transport}
	\end{align}
	\end{subequations}
Here, $\Transport =\rd_t+v^a\rd_a$ (as before), $g$ is the acoustical metric depending on $v^i$
and $\Densrenormalized$ (see Definition \ref{D:ACOUSTICALMETRIC} for the precise definition), 
$\Vortrenormalized = (\partial_1 v^2 - \partial_2 v^1)/\exp(\Densrenormalized)
$
is the specific vorticity, 
and $\mathscr{Q}$ and $\mathscr{Q}^i$ are \emph{null forms relative to $g$}
(which we sometimes refer to as $g$-null forms);
see Proposition \ref{P:GEOMETRICWAVETRANSPORTSYSTEM} for precise definitions. 
It turns out that the $g$-null form structure is crucially important.
In contrast to a $g$-null form, a typical quadratic term could 
severely distort the dynamics near the shock and could in principle prevent
it from forming; see Remark~\ref{R:NULLFORMSAREIMPORTANT}
for further comments.
Thus, under this new formulation, 
we need to consider a coupled system of three quasilinear wave equations and one transport equation.\footnote{Notice that this system is in principle over-determined, but its local-in-time well-posedness follows from that of the 
Euler equations \eqref{E:TRANSPORTDENSRENORMALIZEDRELATIVETORECTANGULAR}-\eqref{E:TRANSPORTVELOCITYRELATIVETORECTANGULAR}.}

\begin{remark}[\textbf{Avoiding vacuum regions}]
	\label{R:AVOIDINGVACUUMREGIONS}
	In this article, we show that the solutions under study have densities that are from bounded from below,
	strictly away from $0$. We therefore 
	avoid the difficult problem of studying the dynamics of a fluid containing vacuum
	regions and hence there is no difficulty in dividing
	by the density to form the specific vorticity.
\end{remark}

Our approach is to 
treat the wave part of the system 
using the ideas from \cite{jSgHjLwW2016}
and to handle the additional terms involving the specific vorticity within the same geometric framework.
In particular, we prove estimates for the geometric vectorfield derivatives of the specific vorticity.
The following are the main tasks that we must accomplish:
\begin{itemize}
	\item Make sure that we can control all terms at a consistent level of derivatives
		(that is, without derivative loss). As we have mentioned, 
		as part of this scheme, we must show that the specific vorticity
		has the same differentiability as the velocity and density, representing
		a gain of one derivative.
	\item Understand the expected blowup-rate of the $L^2$ norms of
		all quantities at all orders
		in terms of powers of $\upmu_{\star}^{-1}$;
		see Subsubsect.~\ref{sec.together} for a summary of the blowup-rates.
		As before, we must distinguish between the top-order and the 
		below-top-order energy
		estimates. 
		An notable feature of the present work is that 
		the top-order derivatives of the specific
		vorticity are allowed to blow up at a worse rate
		than any of the terms that arise in the irrotational case;
		see \eqref{E:INTROTOPVORT}. However, in the coupling to the wave equation,
		the top-derivatives of the specific vorticity appear as a source term
		multiplied by a \emph{critically important factor of $\upmu$}
		(see the term $\upmu \partial_a \Vortrenormalized$ on RHS~\eqref{intro.wave.1}).
		This factor of $\upmu$ turns out to be enough
		to compensate for the especially singular behavior
		of the top-order derivatives of $\Vortrenormalized$.
		We also note that we must ensure the viability
		of the energy descent scheme (see Subsubsect.~\ref{sec.descent.intro}) 
		so that, in particular, we can obtain
		non-degenerate energy and $L^{\infty}$ estimates at the low derivative levels.

		\begin{remark}[\textbf{An alternate approach to controlling the top-order derivatives of $\Vortrenormalized$}]
			\label{R:ALTERNATEAPPROACHTOSPECIFICVORTICITY}
			Although we do not use it in the present article,
			there is alternate approach to controlling the top-order derivatives of $\Vortrenormalized$.
			Specifically, one could 
			differentiate the transport equation $\Transport \Vortrenormalized = 0$
			with the spatial Cartesian coordinate partial derivative vectorfields 
			$\partial_i$ to obtain the evolution equation
			$\Transport \partial_i \Vortrenormalized = - (\partial_i v^a) \partial_a \Vortrenormalized$.
			One could then think of the quantities $\partial_i \Vortrenormalized$ as new variables
			that need to be controlled, in addition to $\Vortrenormalized$.
			Although this approach would involve some 
			additional analysis compared to analysis carried out here,
			the advantage would be that we could close the transport equation
			energy estimates by commuting the transport equations only up to $20$ times
			with geometric vectorfields, as opposed to the approach of 
			the present article, which relies on
			commuting the equation $\Transport \Vortrenormalized = 0$ up to $21$ times.
			In carrying out this alternate strategy, one would avoid generating error
			terms in the transport equations that depend on the top-order derivatives of the 
			eikonal function. In particular, this would allow us to avoid the most singular terms
			and to therefore derive less degenerate estimates for $\Vortrenormalized$ at the top-order
			compared to the estimates that we obtain in this article.
			In our forthcoming work \cite{jLjS2017} on shock formation with vorticity in three spatial dimensions, 
			it turns out that we are forced to employ a closely related strategy
			and, as we mentioned earlier, to complement it with elliptic estimates.
			The reason is that in three spatial dimensions, the evolution equation verified by $\Vortrenormalized$ is
			no longer homogeneous, but rather $\Transport \Vortrenormalized^i = \Vortrenormalized^a \partial_a v^i$
			(recall that in three spatial dimensions, $\Vortrenormalized$ is a $\Sigma_t$-tangent vectorfield).
			Thus, the simplified approach of the present article, which is based on commuting the \emph{homogeneous}
			transport equation for $\Vortrenormalized$ up to top-order with geometric derivative vectorfields,
			would result in the loss of a derivative in three spatial dimensions
			(coming from the term generated when all derivatives fall on the factor $\partial_a v^i$ 
			in the product $\Vortrenormalized^a \partial_a v^i$).
	\end{remark}
	\item Make sure that the expected blowup-rates are consistent in the sense that the coupling
		does not spoil the expectation.
		By coupling, we mean coupling between the 
		``wave variables''
		$v^i$ and $\Densrenormalized$, the specific vorticity
		$\Vortrenormalized$, and the acoustic geometry (that is, the eikonal function), 
		which enters into the analysis in particular through the term $\mytr \upchi$.
	\item Go beyond ``consistency'' by actually closing 
		the energy estimates. For this, it is important to exploit 
		various kinds of
		smallness in the problem (in addition to those that are already present in the irrotational case).
		For instance, in the energy estimates for the specific vorticity, 
		the wave variables and the acoustic geometry enter with an extra smallness constant $\mathring{\upepsilon}^2$ (see \eqref{vorticity.lower.intro} and \eqref{vorticity.top.intro}) so that the latter variables couple only weakly\footnote{This is a big difference from the case 
		of three spatial dimensions, where the coupling is much stronger. 
		We will discuss this issue in our forthcoming work \cite{jLjS2017} in the three spatial dimensional case.} 
		to the specific vorticity. This allows the energy estimates for the 
		specific vorticity to be closed semi-independently with the help of appropriate
		bootstrap assumptions for the behavior of the wave variables and the acoustic geometry.
		Another useful but more subtle source of smallness is 
		tied to the fact that we are treating perturbations of simple 
		outgoing (that is, right-moving)
		plane waves. For example, $\Rad(v^1 - \Densrenormalized)$
		is $\mathcal{O}(\mathring{\upepsilon})$ small even though
		$\Rad v^1$ and $\Rad \Densrenormalized$ are not.
		This smallness allows us to exploit effective decoupling
		between different solution variables, which turns out to be important
		for minimizing the size of certain key coefficients and therefore minimizing\footnote{The size of the coefficients
		is tied to the blowup-rate of the top-order energies which is in turn tied to the number of derivatives needed to close;
		see, for example, the ``$6$'' on RHS~\eqref{E:KEYCOEFFICIENT}.}
		the number of derivatives needed to close the problem; see the discussion in Subsubsect.~\ref{sec.gd.comp}.
\end{itemize}

In order to close the estimates, we will commute the wave equations with up to $20$ 
geometric vectorfields and the transport equation with $21$
geometric vectorfields;\footnote{With additional effort, we could slightly reduce the number of derivatives
that we need to close.} 
see, however, Remark~\ref{R:ALTERNATEAPPROACHTOSPECIFICVORTICITY}.
As we described above, since $\Vortrenormalized$ is at the level of one
derivative of $v^i$, this represents a gain of one derivative for $\Vortrenormalized$. 
Define\footnote{In the proof, we will denote the boundary energy norms by $\mathbb Q_N$ and the bulk spacetime norm by $\mathbb K_N$. In this subsection, in order to simplify the exposition, we will not make this distinction.} 
\[
\mathbb{W}_N
\]
to be the energy norm for $v^i$ and $\Densrenormalized$
corresponding to $N$ commutations of the wave equations with geometric vectorfields, 
where we require\footnote{Note that $N=1$ corresponds to controlling \emph{two} derivatives of $v^i$ and $\Densrenormalized$.}
$N\geq 1$ and allow at most one of them to be $\Rad$. Moreover, 
the case of a single pure $\Rad$ commutation is excluded.
Notice that for technical reasons,\footnote{In commuting the specific vorticity equation, 
we encounter a new term that forces us to commute the wave equations with one copy of $\Rad$,
namely the terms $\Rad \GeoAng^{N-1} \mytr \upchi$
on RHSs \eqref{E:VORTICITYLISTHEFIRSTCOMMUTATORIMPORTANTTERMS}-\eqref{E:VORTICITYGEOANGANGISTHEFIRSTCOMMUTATORIMPORTANTTERMS}.
We will downplay this issue here.}
we have slightly modified the approach to commuting the wave equations taken in \cite{jSgHjLwW2016}.
In particular, unlike in \cite{jSgHjLwW2016}, 
we now commute with up to one $\Rad$ (recall 
that in Subsect.~\ref{sec.ideas.review.plane}, only $\Lunit$ and $\GeoAng$ were used as commutators).
We also note the energy includes a term on 
$\Sigma_t^u$ (cf. \eqref{E.intro}), a term on 
the characteristics $\mathcal{P}_u^t$ (cf. \eqref{F.intro}) and a spacetime bulk term (cf. \eqref{K.intro}). Our choice of the structure of the strings of commutation
vectorfields ensures that the energy is $\mathcal{O}(\mathring{\upepsilon})$-small.\footnote{Recall that the derivatives of $(v^1,v^2,\Densrenormalized)$ 
are small if at least one of the geometric vectorfields is $\Lunit$ or $\GeoAng$.} 

To control the specific vorticity $\Vortrenormalized$, we will use an energy norm $\mathbb{V}$, 
which we define below (see \eqref{E:INTROVORTEN}). Since $\Vortrenormalized$ satisfies a \emph{homogeneous} transport equation, the main challenge is to control the commutators of the weighted transport operator $\upmu \Transport$ 
and the geometric commutation fields, which are adapted to the acoustic characteristics.
In order to minimize the number of $\Rad$ commutators needed 
to control the ``wave variables'' $(v^1,v^2,\Densrenormalized)$ (which satisfy the wave equations) 
and $\mytr \upchi$, 
we commute the transport equation only with the $\mathcal{P}_u$-tangent vectorfields
$\Lunit$ and $\GeoAng$. Because the material derivative vectorfield is transversal\footnote{The transversality follows from a simple geometric fact: 
in all solution regimes, $\Transport$ is a $g$-timelike vectorfield (that is, $g(\Transport,\Transport) < 0$);
thus, $\Transport$ 
cannot be tangent to any $g$-null hypersurface. In fact, we have $\Transport = \Lunit + \Radunit$ 
and hence the $\Radunit$ component of $\Transport$ is bounded below all the way up to the shock. This then allows us to use equation \eqref{intro.transport} to algebraically express 
$\Radunit \Vortrenormalized = - \Lunit \Vortrenormalized$. Similarly, higher $\Rad$ derivatives of $\Radunit \Vortrenormalized$ can be expressed in terms of derivatives tangential to the $g$-null hypersurfaces. \label{FN:TRANSVERSALEQUALSTANGENTIAL}} 
to the acoustic characteristics 
$\mathcal{P}_u$, this is sufficient for obtaining
estimates for all directional derivatives of $\Vortrenormalized$ and closing the argument.

In the next few subsubsections, we will discuss
the various energy estimates needed 
to control
the specific vorticity, 
the eikonal function and the wave variables. 
Let us already note at this point that the main difficulty comes in the \emph{top-order} derivative, 
where the singular behaviors of
the wave variables, the specific vorticity,
and the geometry of the null hypersurfaces are all coupled.

\subsubsection{Lower-order energy estimates for the specific vorticity}
The energy norm that we use to control the specific vorticity at the lowest order is
\begin{align} \label{E:INTROVORTEN}
	\mathbb{V}(t,u)
	& := \int_{\Sigma_t^u} \upmu \Vortrenormalized^2 \,d\vartheta\, du+ \int_{\mathcal P_u^t} \Vortrenormalized^2\, d\vartheta\, dt.
\end{align}
In other words, the energy for $\Vortrenormalized$ on a constant-$t$ hypersurface $\Sigma_t^u$
is ``degenerate'' in $\upmu$, while that for $\Vortrenormalized$ on a constant-$u$ hypersurface $\mathcal{P}_u^t$
is ``non-degenerate''.\footnote{The energies for $V(t,u)$
are in fact very natural.
If one changes variables and expresses the forms $\upmu \,d\vartheta\, du$ and 
$d\vartheta\, dt$ relative to the Cartesian coordinates, then one sees that,
up to $\mathcal{O}(1)$ multiplicative factors, these forms agree with the usual forms
induced on the corresponding hypersurfaces
by the Euclidean metric on $\mathbb{R}^{1+2}$.}

In deriving energy estimates, we also control the derivatives of $\Vortrenormalized$ with respect to $\Lunit$ and $\GeoAng$ and use the notation 
$\mathbb{V}_{N}(t,u)$ to denote the corresponding energy norm after $N$ commutations. In particular, 
for $N$ sufficiently large,
the non-degenerate control
of $\mathbb{V}_0$, $\cdots$, $\mathbb{V}_{N}$ on the acoustic characteristics $\mathcal{P}_u^t$,
when combined with Sobolev embedding, gives rise to\footnote{Let us note that pointwise control can alternatively be derived directly using the transport equation itself.} pointwise 
control of $\Vortrenormalized$ and its lower-order $\Lunit$ and $\GeoAng$ derivatives.

While $\Vortrenormalized$ itself satisfies a homogeneous transport equation (see \eqref{intro.transport}), 
to control its derivatives, 
we need to bound the commutator terms and 
derive estimates for solutions to
inhomogeneous transport equations. 
For the general inhomogeneous equation $\upmu \Transport \Vortrenormalized = \vortinhom$,
we have the following estimate (see Proposition \ref{P:ENERGYIDENTITYRENORMALIZEDVORTICITY}):
\begin{equation}\label{transport.inho.intro}
\int_{\Sigma_t^u} \upmu \Vortrenormalized^2 \,d\vartheta\, du+ \int_{\mathcal P_u^t} \Vortrenormalized^2\, d\vartheta\, dt\lesssim \mbox{Data}+\int_{\mathcal M_{t,u}} \vortinhom^2 \,d\vartheta\, du\, dt.
\end{equation}
Both of the commutators $[\upmu \Transport, \Lunit]$ and $[\upmu \Transport, \GeoAng]$ 
generate controllable error terms that are regular with respect to $\upmu$;
this of course is the main reason to commute the geometric vectorfields with $\upmu \Transport$ (instead of, say, $\Transport$).
The following equation exhibits a typical difficult inhomogeneous term that we have to control
after $N$ commutations:
\begin{equation}\label{Vr.commuted.intro}
\upmu \Transport \GeoAng^N \Vortrenormalized=(\GeoAng \Vortrenormalized)\Rad \GeoAng^{N-2} \mytr \upchi+\dots, 
\end{equation}
where $\dots$ denotes terms that are easier to handle. Therefore, except for the top-order case $N=21$, 
one can control the term 
$(\GeoAng \Vortrenormalized)\Rad \GeoAng^{N-2} \mytr \upchi$ by first showing\footnote{Actually, 
the smallness of $\GeoAng \Vortrenormalized$ is one of our bootstrap assumptions.}
that $\GeoAng \Vortrenormalized$ is $\mathcal{O}(\mathring{\upepsilon})$-small in $L^\infty$ and that $\Rad \GeoAng^{N-2} \mytr \upchi$ can be controlled by an analogue of \eqref{trch.no.gain.intro}. Using \eqref{transport.inho.intro}, this roughly yields
the following inequality for $N<21$:
\begin{equation}\label{vorticity.lower.intro}
\mathbb{V}_N(t,u) \lesssim  \mathring{\upepsilon}^2
	+ \mathring{\upepsilon}^2 \int_{t'=0}^t\left(\int_{s=0}^{t'} \upmu^{-\frac{1}{2}}_{\star}(s,u) \mathbb{W}_N^{\frac{1}{2}}(s,u) \, ds\right)^2\, dt' +\dots.
\end{equation}
As we will see, the coupling with $\mathbb{W}_N$ 
in equation \eqref{vorticity.lower.intro} is quite weak. More precisely,
due to the small factor $\mathring{\upepsilon}^2$ and the
large number of time integrations on RHS~\eqref{vorticity.lower.intro},
the influence of RHS~\eqref{vorticity.lower.intro} on $\mathbb{V}_N$ is easy to control.

On the other hand, for $N=21$, one does not have 
the luxury 
of using the $\mathbb{W}_{21}$ norm on the right hand side, since $\mathbb{W}_{20}$ is top-order. 
Hence,
as we will later see, 
at the top-order, we have to take a different approach
to controlling certain terms in the top-order inhomogeneous transport equation, 
an approach which avoids relying on $\mathbb{W}_{21}$.
The different approach forces us to confront 
the most singular terms in 
the $21$-times-commuted transport equation:
$\Rad \GeoAng^{19} \mytr \upchi$
and
$\GeoAng^{20} \mytr \upchi$.
Specifically, we have to account for the singular behavior of 
the $L^2$ norms of
$\Rad \GeoAng^{19} \mytr \upchi$
and
$\GeoAng^{20} \mytr \upchi$
in terms of powers of $\upmu_{\star}^{-1}$.
Note that, as we described in Subsubsect.~\ref{eikonal.wave.intro},
the singular behavior of the top-order derivatives of $\mytr \upchi$,
which is tied to the necessity of using modified quantities to avoid derivative loss,
is already present in the irrotational case as the primary source of degeneracy.

\subsubsection{Top-order estimates for the eikonal function}\label{sec.top.trch}
Before we discuss the top-order estimates for $\Vortrenormalized$,
it makes sense to first consider the top-order derivatives of the eikonal function
(in particular the top-order derivatives of the mean curvature $\mytr \upchi$ of $\mathcal{P}_u$), 
as they are the main source terms in the vorticity estimates
(see equation \eqref{Vr.commuted.intro}). In our setting, 
we again need to 
use modified quantities as described in Subsubsect.~\ref{eikonal.wave.intro} 
in order to obtain sufficient top-order estimates for $\mytr \upchi$.
However, since the ``gain of a derivative'' 
that one achieves with modified quantities
relies on the wave equations 
satisfied by the Cartesian metric components,\footnote{By \eqref{E:ACOUSTICALMETRIC},
the Cartesian metric components depend on $v^1$, $v^2$ and $\Densrenormalized$.} 
which feature source terms
depending on the specific vorticity, 
this procedure is now \emph{coupled} with the estimates for the specific vorticity. 
As a consequence, at the top-order, \emph{the specific vorticity is directly coupled to the evolution
of $\mytr \upchi$}.
Indeed, we recall from the discussion\footnote{Let us note that while in Subsubsect.~\ref{eikonal.wave.intro} we were dealing with a scalar equation, the system case can be dealt with similarly. We discuss here only the estimates involving $v^i$ (as it is slightly harder) and suppress those involving the wave equation for $\Densrenormalized$.} in Subsubsect.~\ref{eikonal.wave.intro} that in order to use \eqref{Ray.eqn.intro} to gain a derivative for $\mytr \upchi$,
we need to use the wave equation to exchange $\upmu\slashed\Delta v^i$ with an exact $\Lunit$-derivative. 
Since the wave equation features the inhomogeneous terms $\upmu \rd\Vortrenormalized$ and $\upmu\Vortrenormalized\rd \Densrenormalized$,
these terms will couple into the estimates for $\mytr \upchi$. 
In the next paragraph, we describe the effect of this coupling.\footnote{Note that at the same time, equation \eqref{Vr.commuted.intro}
shows that the top-order derivatives of $\mytr \upchi$
couple into the top-order transport equation for $\Vortrenormalized$.
However, we will postpone the discussion of the effect of the top-order derivatives
of $\mytr \upchi$ on the top-order derivatives of $\Vortrenormalized$
until Subsubsect.~\ref{SSS:TOPORDERMODVORT}.}

Let us focus on the term $\upmu \rd \Vortrenormalized$ in equation \eqref{intro.wave.1}
since the second term $\upmu\Vortrenormalized\rd \Densrenormalized$, 
though it gives rise to some singular estimates, 
is easier to handle. A crucial observation, which we already made in Footnote~\ref{FN:TRANSVERSALEQUALSTANGENTIAL}, 
is that by algebraically using the transport equation for $\Vortrenormalized$, 
one can express $\upmu \rd \Vortrenormalized$ as linear combinations of $\upmu L\Vortrenormalized$ and $\upmu Y\Vortrenormalized$.
As a consequence, the $\mathbb{V}_N$ norms suffice\footnote{Recall that by definition,
the norms $\mathbb{V}_N$ control only the $\Lunit$ and $\GeoAng$ derivatives of $\Vortrenormalized$.}
to control these terms and we obtain the following bound for the top-order\footnote{The top-order derivatives of $\mytr \upchi$ involving at least
one $\Lunit$ differentiation are much easier to control since one 
can directly bound it by estimating the RHS of \eqref{Ray.eqn.intro} and 
hence does not need to use the modified quantities of Subsect.~\ref{eikonal.wave.intro}
to handle them.} 
derivatives of $\mytr \upchi$:
\begin{equation}\label{trch.top.intro}
\|\upmu \Rad \GeoAng^{19} \mytr \upchi\|_{L^2(\Sigma_t^u)}, \,\|\upmu \GeoAng^{20} \mytr \upchi\|_{L^2(\Sigma_t^u)}\lesssim \mathbb{W}^{\frac{1}{2}}_{20}+\int_{t'=0}^t \mathbb{V}_{21}^{\frac{1}{2}}(t',u)\, dt'+\dots,
\end{equation}
where $\dots$ are similar or less singular
terms.
We stress that the time integral term on RHS~\eqref{trch.top.intro} 
is exactly the term that accounts for the influence
of the top-order derivatives of the specific vorticity on the 
acoustic geometry.

\subsubsection{Top-order energy estimates for the specific vorticity}
\label{SSS:TOPORDERMODVORT}

We now return to the discussions for the estimates for $\Vortrenormalized$, but this time at the top-order derivative. According to \eqref{transport.inho.intro} and \eqref{Vr.commuted.intro}, at the top-order, we need to bound\footnote{There is in fact a similar term which features $\GeoAng^{20} \mytr \upchi$ that we have suppressed in \eqref{Vr.commuted.intro}. It is as difficult as the term featuring $\Rad \GeoAng^{19} \mytr \upchi$, although in view of \eqref{trch.top.intro}, it can be estimated in a similar manner.} $\Rad \GeoAng^{19} \mytr \upchi$ and $\GeoAng^{20} \mytr \upchi$
in a suitable spacetime $L^2$ norm. 
With the help of the estimate \eqref{trch.top.intro} for 
$\Rad \GeoAng^{19} \mytr \upchi$ and $\GeoAng^{20} \mytr \upchi$,
we can obtain the following top-order estimate 
(see Prop.~\ref{P:VORTICITYENERGYINTEGRALINEQUALITIES} for the details):
\begin{equation}\label{vorticity.top.intro}
\begin{split}
\mathbb{V}_{21}(t,u)\lesssim  \mathring{\upepsilon}^2 \left(\underbrace{\int_{t'=0}^t \upmu_{\star}^{-2}\mathbb{W}_{20}(t',u)\, dt'}_{VT_1}+\underbrace{\int_{t'=0}^t\upmu_{\star}^{-2}(t',u)\left( \int_{s=0}^{t'} \mathbb{V}_{21}^{\frac{1}{2}}(s,u)\, ds\right)^2\, dt'}_{VT_2}\right)+\dots
\end{split}
\end{equation}
Notice that $VT_2$  
features the singular factor
$\upmu_{\star}^{-2}$ and a total of three time integrations. Since $\upmu_{\star} \to 0$ at worst linearly
(as we described in Subsubsect.~\ref{SSS:IDENTIFYINGBLOWUPMECHANISM}),
as long as we are willing to settle for proving
a sufficiently singular bound,\footnote{That is, as long as we are proving that $\mathbb{V}_{21}(t,u)$ is 
bounded from above by some negative powers of $\upmu_{\star}$.} the term $VT_2$
can be treated with a Gronwall-type argument. On the other hand, as we will later see, the
term $VT_1$ determines the blowup-rate of $\mathbb{V}_{21}$ in terms of negative powers of $\upmu_{\star}$.
We now recall a crucial feature of the estimate \eqref{vorticity.top.intro}
mentioned earlier, namely that
the coupling to $\Rad \GeoAng^{19} \mytr \upchi$ and $\GeoAng^{20} \mytr \upchi$
is weak in that there is a small factor $\mathring{\upepsilon}^2$
on RHS~\eqref{vorticity.top.intro}.
For this reason, one can actually derive sufficient estimates
for $\mathbb{V}_{21}(t,u)$
using only bootstrap assumptions for the energy norms,
based on having a good guess for the blowup-rates of all quantities. 
One is aided in this endeavor by the fact, justified later on, 
that $\mathbb{W}_{20}$ blows up at the same rate as in the irrotational case.
This basic fact allows one to control the specific vorticity in a relatively
straightforward fashion.

\subsubsection{Lower-order energy estimates for the wave variables}
To close the argument, we need to 
derive estimates for the 
``wave variables'' $v^1$, $v^2$, and $\Densrenormalized$, that is, for
solutions to the wave equations \eqref{intro.wave.1}-\eqref{intro.wave.2},
and in particular to estimate the vorticity terms arising on the right hand side of equation 
\eqref{intro.wave.1}. 
When we are bounding the below-top-order derivatives of the wave variables, 
we do not need to rely on modified quantities to control the eikonal function.
For this reason, the below-top-order estimates are relatively easy to derive, as we now describe.
For this discussion, we suppress most of the terms that do not involve $\Vortrenormalized$ except 
for one that is analogous to the term on RHS~\eqref{E19.intro}; we denote this analogous term by  
$WL_1$ below in \eqref{wave.lower.intro}. The inhomogeneous terms not involving $\Vortrenormalized$ 
can be bounded by using the same arguments as in the irrotational case, so we do not discuss them in detail.
By equation \eqref{intro.wave.1},
the terms in the equation 
$\upmu \square_{g} v^i = \cdots$
involving $\Vortrenormalized$ can be expressed in the form\footnote{Here, we have used the observation discussed in Footnote~\ref{FN:TRANSVERSALEQUALSTANGENTIAL},
namely that by using the transport equation for $\Vortrenormalized$, the term $\Rad \Vortrenormalized$ can 
expressed as $- \upmu \Lunit \Vortrenormalized$.} 
$\upmu L\Vortrenormalized$, $\upmu Y\Vortrenormalized$ or $\Vortrenormalized(\upmu L v^i+\Rad v^i)$. Since the terms $\upmu L\Vortrenormalized$ and $\upmu Y\Vortrenormalized$ contain factors of $\upmu$, when estimating $\mathbb{W}_N$, these terms can be controlled by the degenerate energy 
on constant-$t$ hypersurfaces (that is, the analog of the first term on RHS~\eqref{E:INTROVORTEN}).
On the other hand, since there are 
no extra factors of $\upmu$ in the product $\Vortrenormalized \Rad v^i$,
the $\Vortrenormalized$ factor cannot be bounded by the degenerate energy. Instead, we control
it using the non-degenerate flux 
on $\mathcal{P}_u^t$ (that is, the analog of the second term on RHS~\eqref{E:INTROVORTEN}).
Since the factor $\Vortrenormalized$ is not top-order,
when estimating $\mathbb{W}_N$, 
one needs only to use $\mathbb{V}_N$ to control its up-to-order $N$ derivatives.  
In total, we roughly obtain the following estimate for $N<20$:
\begin{equation}\label{wave.lower.intro}
\begin{split}
\mathbb{W}_N(t,u)\lesssim &\mbox{Data}+\underbrace{\int_0^t \upmu_{\star}^{-\frac{1}{2}} \mathbb{W}_{N}^{\frac{1}{2}}(t',u) \int_{t''=0}^{t''=t'} \upmu_{\star}^{-\frac{1}{2}} \mathbb{W}_{N+1}^{\f12}(t'',u)\, dt'' dt'}_{WL_1}\\
&+\underbrace{\int_{t'=0}^t \mathbb{V}_{N+1}(t',u)\, dt'}_{WL_2}+\underbrace{\int_{u'=0}^u \mathbb{V}_{N}(t,u')\, du'}_{WL_3}
+
\dots.
\end{split}
\end{equation}
The key point here is that the $WL_2$ term has a time integration and thus one can gain\footnote{Let us recall again that $\upmu_{\star}$ goes to $0$ at worst linearly and that \eqref{int.lemma.intro} holds.} a power of $\upmu_{\star}$. 
Such gain cannot be achieved in $WL_3$, but on the other hand, 
the term only features the lower-order $\mathbb{V}_N$ norm
and no singular factor of $\upmu_{\star}^{-1}$.

\subsubsection{Top-order estimates for the wave variables}
In deriving estimates for
the top derivative norm $\mathbb{W}_{20}$, 
we again encounter terms that 
are analogous to the terms
$WL_2$ and $WL_3$ from \eqref{wave.lower.intro}, which are 
respectively denoted by $WT_1$ and $WT_2$ below in \eqref{wave.top.intro}. 
There are also additional terms involving the top derivatives of $\mytr \upchi$, which cannot be
treated like the term $WL_1$ from \eqref{wave.lower.intro}.
These additional terms are in fact precisely the ones described in Subsubsect.~\ref{eikonal.wave.intro}, which 
need to be bounded with the help of modified quantities. We stated an estimate for them in \eqref{trch.top.intro}. 
To proceed, we use the estimate \eqref{trch.top.intro}, 
but
this time carefully tracking the precise
numerical coefficient of the $\mathbb{W}_{20}^{\frac{1}{2}}$ term. 
These give rise to $WT_{main}$ and $WT_3$ in \eqref{wave.top.intro} below. 
In total, we obtain the following estimate
(see Prop.~\ref{P:WAVEENERGYINTEGRALINEQUALITIES} for the details):
\begin{equation}\label{wave.top.intro}
\begin{split}
\mathbb{W}_{20}(t,u)\leq &\mbox{Data}+\underbrace{C_{fix}\int_{t'=0}^t \left(\sup_{\Sigma^u_{t'}} \left|\f{L\upmu}{\upmu}\right|\right)\mathbb{W}_{20}(t',u)\, dt'}_{WT_{main}}\\
&+\underbrace{C\int_{t'=0}^t \mathbb{V}_{21}(t',u)\, dt'}_{WT_1}+\underbrace{C\int_{u'=0}^u \mathbb{V}_{20}(t,u')\, du'}_{WT_2}\\
&+\underbrace{C\int_{t'=0}^t  \f{\mathbb{W}^{\frac{1}{2}}_{20}(t',u)}{\upmu_{\star}(t',u)}\int_{s=0}^{t'} \mathbb{V}_{21}^{\frac{1}{2}}(s,u)\,ds \, dt'}_{WT_3}+\dots.
\end{split}
\end{equation}
Notice that the term $WT_{main}$ is analogous to the term in \eqref{sing.top.intro} and is the 
main term driving the blowup-rate of $\mathbb{W}_{20}(t,u)$. 
As we described in Sects.~\ref{eikonal.wave.intro} and \ref{sec.descent.intro},
the constant $C_{fix}$ in $WT_{main}$ is intimately tied to the number of derivatives needed 
to close the proof.

\subsubsection{Independent bounds for the ``good'' components}\label{sec.gd.comp}
The next ingredient of the proof is to 
derive independent estimates for $v^2$ and $v^1-\Densrenormalized$.
This is crucial 
for obtaining a good estimate for the constant
$C_{fix}$ in \eqref{wave.top.intro}. 
More precisely, we show that $v^2$ and $v^1-\Densrenormalized$ 
obey better bounds than either $v^1$ or $\Densrenormalized$ and 
that all geometric derivatives of $v^2$ and $v^1-\Densrenormalized$,
\emph{including their $\Rad$ derivatives}, are small.
This is of course tied to the assumption that the solution is nearly simple outgoing plane symmetric. 
Indeed, for plane symmetric solutions, 
we have $v^2=0$. Moreover, for the simple outgoing plane symmetric solutions 
described in Subsect.~\ref{sec.ideas.Euler},
we have $0 = \mathcal R_-=v^1-\int_{\Densrenormalized'=0}^\Densrenormalized c_s(\Densrenormalized')\, d\Densrenormalized'$.
Hence, it follows that under our assumed normalization condition $\Speed(\Densrenormalized = 0)=1$ 
from \eqref{speed.normalization.intro} 
and the $L^\infty$-smallness conditions for 
$\Densrenormalized$ and $v^1$, 
all geometric derivatives of $v^1-\Densrenormalized$ are small
for the perturbations of simple outgoing plane symmetric solutions under study.

In our proof, we take advantage of this smallness as follows. 
First, we explicitly prove that
$\Rad (\Densrenormalized - v^1)$ 
and 
$\Rad v^2$
are $\mathcal{O}(\mathring{\upepsilon})$-small in the $L^\infty$ sense; in the next
paragraph, it will become clear why this is important.
Next, we derive independent estimates for the top-order energy norms of $v^2$ and $v^1-\Densrenormalized$.
Let us momentarily\footnote{In the proof of the main theorem, 
we define $\mathbb Q_{20}^{(Partial)}$ and $\mathbb K_{20}^{(Partial)}$ in an analogous manner to 
respectively
denote the boundary terms and the bulk terms in the energy norms.} 
denote the top-order energy norms of $v^2$ and $v^1-\Densrenormalized$ by
$\mathbb{W}_{20}^{(Partial)}$
in order to distinguish it from the energy norm $\mathbb{W}_{20}$ 
(which controls all three of $v^1$, $v^2$ and $\Densrenormalized$).

Roughly speaking,  
if we just track the $WT_{main}$ term in \eqref{wave.top.intro}, 
then we obtain the following system of energy inequalities:
\begin{equation} \label{E:KEYCOEFFICIENT}
\begin{split}
\mathbb{W}_{20}(t,u)\leq & (6+C\mathring{\upepsilon})\int_{t'=0}^t \left(\sup_{\Sigma^u_{t'}} \left|\f{L\upmu}{\upmu}\right|\right)\mathbb{W}_{20}(t',u)\, dt'\\
& + C_* \int_{t'=0}^t \left(\sup_{\Sigma^u_{t'}} \left|\f{L\upmu}{\upmu}\right|\right)\sqrt{\mathbb{W}_{20}}(t',u)
\sqrt{\mathbb{W}^{(Partial)}_{20}}(t',u)\, dt'+\dots,
\\
\mathbb{W}_{20}^{(Partial)}(t,u)\leq & C\mathring{\upepsilon}\int_{t'=0}^t \left(\sup_{\Sigma^u_{t'}} \left|\f{L\upmu}{\upmu}\right|\right)\mathbb{W}_{20}(t',u)\, dt'+\dots,
\end{split}
\end{equation}
where $C_*$ is a possibly large constant that, unlike $C_{fix}$, depends on the equation of state.
A crucial feature of the above system is that
the main term in the inequality for 
$\mathbb{W}_{20}^{(Partial)}(t,u)$ is multiplied by a small factor $\mathring{\upepsilon}$,
which is available thanks to the $L^\infty$ estimates mentioned in the previous paragraph.
This small factor limits the contribution of the main term to the blowup-rate 
for $\mathbb{W}_{20}^{(Partial)}(t,u)$,
which in turn allows us to obtain
semi-independent control of $\mathbb{W}^{(Partial)}_{20}$ and 
thus show that the product 
$
C_*
\displaystyle
\int_{t'=0}^t \left(\sup_{\Sigma^u_{t'}} \left|\f{L\upmu}{\upmu}\right|\right)\sqrt{\mathbb{W}_{20}}(t',u)
\sqrt{\mathbb{W}^{(Partial)}_{20}}(t',u)\, dt'
$ 
does not significantly influence the blowup-rate of $\mathbb{W}_{20}(t,u)$.
In total, this structure allows us to show that the 
constant $C_{fix}$ in \eqref{wave.top.intro} is essentially $6$.
This fact, together with similar estimates for
a few other related terms that we have suppressed, 
determines the total number of derivatives that we need in the argument.
We clarify that if we did not split the energies into $\mathbb{W}_{20}$ and $\mathbb{W}_{20}^{(Partial)}$,
then the constant $C_*$ could in principle increase the blowup-rate of $\mathbb{W}_{20}$,
which would in turn increase number of derivatives we need to close the problem.
Thanks to the splitting, we are able to close the estimates 
by differentiating $v$ and $\Densrenormalized$ up to\footnote{Note that we take up to
$21$ derivatives of $\Vortrenormalized$, which corresponds to up to $22$ derivatives of $v$ and $21$ derivatives of $\Densrenormalized$.} 
$22$ times.

\subsubsection{Putting everything together}\label{sec.together}

We now combine the estimates discussed in the previous subsubsections and show, at least heuristically, that they can close. 
The detailed proof is based on a lengthy Gronwall argument that is
located in 
Subsects.~\ref{SS:MAINVORTICITYAPRIORIENERGYESTIMATES} and \ref{SS:PROOFOFPROPMAINAPRIORIENERGY}.
As is already clear from the discussions above, the estimates for the lower-order derivatives and for the top derivatives are rather different and the most difficult terms 
are found in the estimates for the top-order energies.

We first consider the lower-order estimates, 
where the blowup-rates are determined by \eqref{vorticity.lower.intro} and \eqref{wave.lower.intro}. Recall the discussions of the descent scheme in Subsubsect.~\ref{sec.descent.intro}: 
\emph{For every order of descent, one gains two powers of $\upmu_{\star}$ until one shows that the energy is bounded.} Moreover, we recall that the descent scheme is based mainly on the fact that each time integration reduces the power of the singularity by one:
\begin{equation}\label{descent.sec.2.4}
\int_{s=0}^t
	\upmu_{\star}^{-b}(s,u)
\, ds
\lesssim
\upmu_{\star}^{1-b}(t,u), \quad\mbox{for }b>1.
\end{equation}
For this reason, inequality \eqref{vorticity.lower.intro} suggests that when $\mathbb{W}_N$ is sufficiently singular, one can prove that $\mathbb{V}_N$ is 
less singular than $\mathbb{W}_N$ 
by a factor of $\upmu_{\star}^2$.
This suggests proving the following estimates:\footnote{Here, 
we emphasize the relative singularity between different norms. 
The precise absolute strength of the singularity depends on the estimates at the top level, 
which we will discuss immediately below.}
	{\large \begin{center} \underline{Below-top-order energy hierarchy} \end{center}}
	\begin{subequations}
	\begin{align}
		\sqrt{\mathbb{W}_{15+K}}(t,u)
				& \leq C \mathring{\upepsilon} \upmu_{\star}^{-(K+.9)}(t,u),
			&& (0 \leq K \leq 4),
				 \notag \\
		\sqrt{\mathbb{W}_{N}}(t,u)
				& \leq C \mathring{\upepsilon},
		&&	(0 \leq N\leq 14)				\notag		\\
		\sqrt{\VorttotTanmax{16+K}}(t,u)
		& \leq C \mathring{\upepsilon} \upmu_{\star}^{-(K+.9)}(t,u),
			&& (0 \leq K \leq 4),
				\notag \\
		\sqrt{\VorttotTanmax{N}}(t,u)
		& \leq C \mathring{\upepsilon}
		&& (0\leq N\leq 15). \notag
	\end{align}
	\end{subequations}
Notice that the above hierarchy is consistent in the following sense: 
when one substitutes the 
hierarchy estimates for $\mathbb{V}_{N+1}$ and $\mathbb{V}_N$ into the terms $WL_2$ and $WL_3$ on RHS~\eqref{wave.lower.intro}
and uses \eqref{descent.sec.2.4},
one finds that these terms contribute to the
blowup-rate for the term $\mathbb{W}_N$ on the LHS 
in a manner that is compatible with the estimates for $\mathbb{W}_N$
stated in the hierarchy.
In fact, there is even extra room in these estimates.

Finally, we consider the top-order estimates, 
which are determined by
\eqref{vorticity.top.intro}
and
\eqref{wave.top.intro},
modulo the discussion surrounding equation \eqref{E:KEYCOEFFICIENT}.
As in Christodoulou's work \cite{dC2007}, 
the top-order blowup-rate is determined mainly by the term
$WT_{main}$ on RHS~\eqref{wave.top.intro}.
As we mentioned in the previous subsubsection, the 
blowup-rate of the top-order energies
depends\footnote{Let us emphasize again that this is a slight simplification, as the strength of the singularity in fact depends on the constant in front of \emph{all} of the singular terms, only one of which is written in \eqref{wave.top.intro}.} on the constant $C_{fix}$, which can be precisely estimated, independent of the equation of state;\footnote{In particular, for \emph{all} 
equations of state other than that of the Chaplygin gas, the estimates can be closed with a total of $22$ derivatives of $v^i$ and $21$ derivatives of $\Densrenormalized$. On the other hand, the relative smallness that is required for $\mathring{\upepsilon}$ \emph{does} depend on the equation of state.} see also Footnote~\ref{FN:HUERISTICTOPORDERBLOWUP}.

To see how the top-order estimates for $\mathbb{W}_{20}$ and $\mathbb{V}_{21}$ couple, we combine \eqref{vorticity.top.intro} and 
\eqref{wave.top.intro} to see that $\mathbb{W}_{20}$ is 
better than $\mathbb{V}_{21}$
by a single factor of $\upmu_{\star}$.
Notice that such an estimate is \emph{borderline} in the sense that if either 
the term $VT_1$ in \eqref{vorticity.top.intro} or the term $WT_1$ in \eqref{wave.top.intro} 
involved a slightly worse power of $\upmu_{\star}^{-1}$, 
then the estimates could not close.
We are therefore led to prove the following estimate
(see Subsect.~\ref{SS:PROOFOFPROPMAINAPRIORIENERGY} for the precise details
concerning the blowup-rates in \eqref{E:INTROTOPWAVE}-\eqref{E:INTROTOPVORT}):
{\large \begin{center} \underline{Top-order energy estimates} \end{center}}
\begin{subequations}
	\begin{align}
		\sqrt{\mathbb{W}_{20}}(t,u)
				& \leq C \mathring{\upepsilon} \upmu_{\star}^{-5.9}(t,u),
							&& 
				 \label{E:INTROTOPWAVE} \\
				\sqrt{\VorttotTanmax{21}}(t,u)
		& \leq C \mathring{\upepsilon} \upmu_{\star}^{-6.4}(t,u).
		 &&  \label{E:INTROTOPVORT}
	\end{align}
	\end{subequations}
We clarify that \eqref{E:INTROTOPWAVE} should be viewed as the main estimate determining the blowup-rates of 
not only for $v^1$, $v^2$, and $\Densrenormalized$, but also 
$\mytr \upchi$ (in view of the discussion of Subsubsect.~\ref{eikonal.wave.intro}, 
where we described how to control $\mytr \upchi$ at the top-order by using modified quantities).
More precisely, 
one can trace through the above logic to discover that the blowup-rate 
$\upmu_{\star}^{-5.9}(t,u)$
on RHS~\eqref{E:INTROTOPWAVE}, 
if taken as given, controls the blowup-rates of all other energy quantities.
We also note that \eqref{E:INTROTOPWAVE}-\eqref{E:INTROTOPVORT}
are consistent for the terms $VT_2$ and $WT_3$, 
neither of which are borderline. 
This concludes our discussion of the main ideas of the proof.

\section{Geometric setup}
\label{S:GEOMETRICSETUP}
In this section, we construct most of the geometric objects
that we use to the shock formation and exhibit their basic
properties. We postpone our construction of energies
and the corresponding integration measures
until Sect.~\ref{S:FORMSANDENERGY}.
We postpone our construction of modified quantities,
which are needed for top-order energy estimates,
until Sect.~\ref{S:MODIFIED}.

\subsection{Notational conventions and shorthand notation}
\label{SS:NOTATION}
We start by summarizing some of our notational conventions;
the precise definitions of some of the concepts referred to
here are provided later in the article.

\begin{itemize}
	\item Lowercase Greek spacetime indices 
	$\alpha$, $\beta$, etc.\
	correspond to the Cartesian spacetime coordinates 
	defined in Sect.~\ref{S:FORMULATIONOFEQUATIONS}
	and vary over $0,1,2$.
	Lowercase Latin spatial indices
	$a$,$b$, etc.\ 
	correspond to the Cartesian spatial coordinates and vary over $1,2$.
	All lowercase Greek indices are lowered and raised with the spacetime metric
	$g$ and its inverse $g^{-1}$, and \emph{not with the Minkowski metric}.
\item We sometimes use $\cdot$ to denote the natural contraction between two tensors
		(and thus raising or lowering indices with a metric is not needed). 
		For example, if $\xi$ is a spacetime one-form and $V$ is a 
		spacetime vectorfield,
		then $\xi \cdot V := \xi_{\alpha} V^{\alpha}$.
\item If $\xi$ is an $\ell_{t,u}$-tangent one-form
	(as defined in Sect.~\ref{SS:PROJECTIONTENSORFIELDANDPROJECTEDLIEDERIVATIVES}),
	then $\xi^{\#}$ denotes its $\gsphere$-dual vectorfield,
	where $\gsphere$ is the Riemannian metric induced on $\ell_{t,u}$ by $g$.
	Similarly, if $\xi$ is a symmetric type $\binom{0}{2}$ $\ell_{t,u}$-tangent tensor, 
	then $\xi^{\#}$ denotes the type $\binom{1}{1}$ $\ell_{t,u}$-tangent tensor formed by raising one index with $\ginversesphere$
	and $\xi^{\# \#}$ denotes the type $\binom{2}{0}$ $\ell_{t,u}$-tangent tensor formed by raising both indices with $\ginversesphere$.
\item Unless otherwise indicated, 
	all quantities in our estimates that are not explicitly under
	an integral are viewed as functions of 
	the geometric coordinates $(t,u,\vartheta)$
	of Def.~\ref{D:GEOMETRICCOORDINATES}.
	Unless otherwise indicated, quantities
	under integrals have the functional dependence 
	established below in
	Def.~\ref{D:NONDEGENERATEVOLUMEFORMS}.
\item If $Q_1$ and $Q_2$ are two operators, then
	$[Q_1,Q_2] = Q_1 Q_2 - Q_2 Q_1$ denotes their commutator.
\item $A \lesssim B$ means that there exists $C > 0$ such that $A \leq C B$.
\item $A \approx B$ means that $A \lesssim B$ and $B \lesssim A$.
\item $A = \mathcal{O}(B)$ means that $|A| \lesssim |B|$.
\item Constants such as $C$ and $c$ are free to vary from line to line.
	\textbf{Explicit and implicit constants are allowed to depend in an increasing, 
	continuous fashion on the data-size parameters 
	$\mathring{\updelta}$
	and $\TranminusdatasizeWithFactor^{-1}$
	from
	Sect.~\ref{SS:FLUIDVARIABLEDATAASSUMPTIONS}.
	However, the constants can be chosen to be 
	independent of the parameters $\mathring{\upepsilon}$ 
	and $\varepsilon$ whenever $\mathring{\upepsilon}$ and $\varepsilon$
	are sufficiently small relative to 
	$\mathring{\updelta}^{-1}$
	and $\TranminusdatasizeWithFactor$.
	}
\item $\lfloor \cdot \rfloor$
	and $\lceil \cdot \rceil$
	respectively denote the floor and ceiling functions. 
\end{itemize}

\subsection{A caveat on citations}
\label{SS:REMARKSONCITATIONSOFEARLIERWORK}
	We often cite \cite{jSgHjLwW2016}
	for equations and identities
	We now point out some minor discrepancies between the work
	\cite{jSgHjLwW2016} and the present work; we will not explicitly comment
	on them again, even though they occur throughout our work.
	Some of the concepts referred to here are defined later in the article.

\begin{itemize}
	\item In citing \cite{jSgHjLwW2016}, we sometimes adjust the formulas
	involving Cartesian metric components to take into account the explicit
	form of $g_{\alpha \beta}$ and $(g^{-1})^{\alpha \beta}$ stated in Def.~\ref{D:ACOUSTICALMETRIC}.
	\item 
	In \cite{jSgHjLwW2016},
	the metric components $g_{\alpha \beta}$ were functions of a scalar-valued function $\Psi$,
	as opposed to the array $\threePsi$ (defined in Def.~\ref{D:VECPSI}).
	For this reason, we must make minor adjustments 
	to many of the formulas
	from \cite{jSgHjLwW2016} 
	to account for the fact that in the present article, $\threePsi$ is an array.
	In all cases, our minor adjustments can easily be verified by examining
	the corresponding proof in \cite{jSgHjLwW2016}.
\item In \cite{jSgHjLwW2016}, the array $\BadVar$ 
	(see definition \eqref{E:ABBREIVATEDVARIABLES})
	was defined to contain the entry
	$\upmu - 1$ rather than $\upmu$.
	However, that difference is not important
	and in particular does not affect the validity of any of the 
	schematic formulas that
	we cite from \cite{jSgHjLwW2016}.
\end{itemize}

\subsection{Formulation of the equations}
\label{S:FORMULATIONOFEQUATIONS}
We now formulate the evolution equations, the main result being Prop.~\ref{P:GEOMETRICWAVETRANSPORTSYSTEM}.
As we mentioned at the beginning, 
we assume that the space manifold, 
on which the equations are posed, is  
\begin{align} \label{E:SPACEMANIFOLD}
	\Sigma 
	& := \mathbb{R} \times \mathbb{T},
\end{align}
where $\mathbb{R}$ corresponds to time
and $\Sigma$ to space. 
We fix a standard Cartesian coordinate system 
$\lbrace x^{\alpha} \rbrace_{\alpha = 0,1,2}$
on $\mathbb{R} \times \Sigma$, where $x^0 \in \mathbb{R}$
is the time coordinate and $(x^1,x^2) \in \mathbb{R} \times \mathbb{T}$
are the spatial coordinates. The coordinate $x^2$ corresponds to perturbations 
away from plane symmetry. We denote the corresponding
Cartesian coordinate partial derivative vectorfields by
$
\displaystyle
\partial_{\alpha} 
:= \frac{\partial}{\partial x^{\alpha}} 
$.
The coordinate $x^2$ is only locally defined even though
$\partial_2$ can be extended to a globally defined vectorfield
on $\mathbb{T}$. We often use the alternate notation
$t = x^0$ and $\partial_t = \partial_0$.
 
The compressible Euler equations are evolution equations for
the velocity $v:\mathbb{R} \times \Sigma \rightarrow \mathbb{R}^2$
and the density $\rho:\mathbb{R} \times \Sigma \rightarrow [0,\infty)$.
To close the system, we assume a barotropic equation of state
\begin{align} \label{E:BARATROPICEOS}
	p = p(\rho),
\end{align}
where $p$ is the pressure. 
To the equation of state, we associate 
the quantity $\Speed$, known as the \emph{speed of sound} 
\begin{align} \label{E:SOUNDSPEED}
	\Speed 
	& := \sqrt{\frac{dp}{d \rho}}.
\end{align}
Physical equations of state are such that
\begin{itemize}
	\item $\Speed \geq 0$.
	\item $\Speed > 0$ when $\rho > 0$.
\end{itemize}
We study solutions with $\rho > 0$,
which, under the above assumptions,
ensures the hyperbolicity of the system.
In particular, we avoid the study of fluid-vacuum boundaries,
which is accompanied by technical difficulties tied to the degeneracy
of the hyperbolicity along the boundary.

Our shock formation results apply to all equations of state except for those
corresponding to a Chaplygin gas, which are of the form
\begin{align} \label{E:EOSCHAPLYGINGAS}
	p = C_0 - \frac{C_1}{\rho}
\end{align}
for constants $C_0 \in \mathbb{R}$ and $C_1 > 0$.

\subsubsection{Vorticity, modified variables, the speed of sound and its derivative with respect to 
$\Densrenormalized$}
\label{SS:VORTICITYANDMODIFIEDVARIABLES}
In two spatial dimensions, the vorticity $\omega$ is the scalar-valued function
\begin{align} \label{E:VORTICITYDEFINITION}
	\omega 
	& := \partial_1 v^2 - \partial_2 v^1.
\end{align}
Although $\omega$ is an auxiliary variable, it plays a fundamental role
in our analysis.

Rather than directly studying the
density and the vorticity,
we find it convenient to instead study the
logarithmic density
and the specific vorticity.

\begin{definition}[\textbf{Modified variables}]
\label{D:MODIFIEDVARIABLES}
We define the 
\emph{logarithmic density} $\Densrenormalized$
and the \emph{specific vorticity} $\Vortrenormalized$
as follows:
\begin{align} \label{E:MODIFIEDVARIABLES}
	\Densrenormalized
	& := \ln \rho,
		\qquad
	\Vortrenormalized
	:= \frac{\omega}{\rho}
	= \frac{\omega}{\exp \Densrenormalized}.
\end{align}
\end{definition}

From now on, we view $\Speed$, defined in \eqref{E:SOUNDSPEED}, as a function of $\Densrenormalized$:
\begin{align} \label{E:SPEEDOFSOUNDISAFUNCTIONOFRENORMALIZEDDENSITY}
	\Speed 
	& = \Speed(\Densrenormalized).
\end{align}
Moreover, we set
\begin{align} \label{E:DEFINITIONOFSPEEDPRIME}
	\Speed'
	= \Speed'(\Densrenormalized)
	& := \frac{d}{d \Densrenormalized} \Speed(\Densrenormalized).
\end{align}

\subsubsection{Geometric tensorfields associated to the flow}
\label{SS:GEOMETRICTENSORFIELDSASSOCITEDTOTHEFLOW}
To derive our main results, we rely on a geometric formulation
of the Euler equations, derived in the companion article \cite{jLjS2016a},
which exhibit remarkable structures. Before stating the equations,
we define some tensorfields that lie at the heart of our analysis.

We start by defining the material derivative vectorfield,
which transports the specific vorticity.

\begin{definition}[\textbf{Material derivative vectorfield}]
\label{D:MATERIALVECTORVIELDRELATIVETORECTANGULAR}
 The \emph{material derivative vectorfield}
$\Transport$ is defined as follows
	relative to the Cartesian coordinates:
\begin{align} \label{E:MATERIALVECTORVIELDRELATIVETORECTANGULAR}
	\Transport 
	& := \partial_t + v^a \partial_a.
\end{align}
\end{definition}

Next, we define the acoustical metric $g$. It is the Lorentzian spacetime metric
corresponding to the propagation of sounds waves.

\begin{definition}[\textbf{The acoustical metric and its inverse}] 
\label{D:ACOUSTICALMETRIC}
We define the \emph{acoustical metric} $g$ and the 
\emph{inverse acoustical metric} $g^{-1}$ relative
to the Cartesian coordinates as follows:
\begin{subequations}
	\begin{align}
		g 
		& := 
		-  dt \otimes dt
			+ 
			\Speed^{-2} \sum_{a=1}^2(dx^a - v^a dt) \otimes (dx^a - v^a dt),
				\label{E:ACOUSTICALMETRIC} \\
		g^{-1} 
		& := 
			- \Transport \otimes \Transport
			+ \Speed^2 \sum_{a=1}^2 \partial_a \otimes \partial_a.
			\label{E:INVERSEACOUSTICALMETRIC}
	\end{align}
\end{subequations}
\end{definition}

\begin{remark}
	\label{R:GINVERSEISTHEINVERSE}
	It is straightforward to verify that
	$g^{-1}$ is the matrix inverse of $g$, that
	is, we have $(g^{-1})^{\mu \alpha} g_{\alpha \nu} = \delta_{\nu}^{\mu}$,
	where $\delta_{\nu}^{\mu}$ is the standard Kronecker delta.
\end{remark}

\subsubsection{Statement of the geometric form of the equations}
\label{SS:STATEMENTOFEQUATIONS}
We now state the form of the equations that we use to analyze solutions.
The equations were essentially derived in \cite{jLjS2016a}, 
up to the following three remarks: \textbf{i)} We have multiplied the equations by a weight
$\upmu > 0$ that we explain in great detail below. The reason is that the 
$\upmu$-weighted equations have better commutation properties with various
differential operators compared to the unweighted equations.
\textbf{ii)} In \cite{jLjS2016a}, the equations were derived in three space dimensions,
in which case the specific vorticity is a vectorfield. In that context, the 
analog of equation \eqref{E:RENORMALIZEDVORTICTITYTRANSPORTEQUATION}
is a vector equation for 
the Cartesian components $\Vortrenormalized^i$, $(i=1,2,3)$.
The equation features the non-zero ``vorticity stretching'' source term 
$\sum_{a=1}^3 \Vortrenormalized^a \partial_a v^i$.
In $2D$, this source term vanishes, as we now explain.
We may view the $2D$ Euler equations as a special case of the $3D$ Euler equations
in which $v^3 \equiv 0$, $\partial_3 v^i \equiv 0$, and 
the vectorfield $\Vortrenormalized$ is proportional to 
$(\partial_2 v^3 - \partial_3 v^2) \partial_3$.
It follows that $\sum_{a=1}^3 \Vortrenormalized^a \partial_a v^i \equiv 0$ 
in $2D$, that is, the vorticity stretching term vanishes. 
Hence, in the remainder of the article, we view $\Vortrenormalized$
to be the scalar-valued function in \eqref{E:MODIFIEDVARIABLES}.
\textbf{iii)} In \cite{jLjS2016a}, an additional term
$
-
\Speed^{-1} \Speed' 
(g^{-1})^{\alpha \beta} \partial_{\alpha} \Densrenormalized \partial_{\beta} v^i
$
appeared on the analog of RHS~\eqref{E:VELOCITYNULLFORM}
and the coefficient of the first product of the analog of RHS~\eqref{E:DENSITYNULLFORM}
was $-3$ instead of $-2$. The reason for the discrepancy 
is that relative to Cartesian coordinates,
$\mbox{\upshape det} g = \Speed^{-6}$ in three space dimensions
while $\mbox{\upshape det} g = \Speed^{-4}$ in the present case of 
two space dimensions; in view of this fact and Footnote~\ref{FN:COVWAVEOPARBITRARYCOORDS}, 
we see that the form of $\square_g$ relative to the Cartesian coordinates
depends on the number of spatial dimensions. In turn, this affects 
the coefficients of the semilinear terms present 
on the RHS of the wave equations. However, this is a minor point
that has no substantial bearing on the analysis; the 
products under discussion are null forms and thus have only a negligible effect on the dynamics;
see Remark~\ref{R:NULLFORMSAREIMPORTANT}.

\begin{proposition}[\textbf{The geometric wave-transport formulation of the compressible Euler equations}]
	\label{P:GEOMETRICWAVETRANSPORTSYSTEM}
	Let $\square_g$ denote the covariant wave operator
	of the acoustic metric $g$ defined by \eqref{E:ACOUSTICALMETRIC}.
	In two spatial dimensions,
	classical solutions to
	the compressible Euler equations 
	\eqref{E:TRANSPORTDENSRENORMALIZEDRELATIVETORECTANGULAR}-\eqref{E:TRANSPORTVELOCITYRELATIVETORECTANGULAR}
	verify the following equations,
	where the Cartesian components 
	$v^i$, $(i=1,2)$, are viewed as scalar-valued functions 
	under covariant differentiation:\footnote{Here, we use the square bracket $[\cdot]$ to denote the anti-symmetrization of the indices.}
	\begin{subequations}
	\begin{align}
		\upmu \square_g v^i
		& = - [ia] (\exp \Densrenormalized) \Speed^2 (\upmu \partial_a \Vortrenormalized)
			+ 2 [ia] (\exp \Densrenormalized) \Vortrenormalized (\upmu \Transport v^a)
			+ \upmu \mathscr{Q}^i,
			\label{E:VELOCITYWAVEEQUATION}	\\
	\upmu \square_g \Densrenormalized
	& = \upmu \mathscr{Q},
		\label{E:RENORMALIZEDDENSITYWAVEEQUATION} \\
	\upmu \Transport \Vortrenormalized
	& = 0.
	\label{E:RENORMALIZEDVORTICTITYTRANSPORTEQUATION}
	\end{align}
	\end{subequations}
	In \eqref{E:VELOCITYWAVEEQUATION}-\eqref{E:RENORMALIZEDVORTICTITYTRANSPORTEQUATION},
	$\mathscr{Q}^i$ and $\mathscr{Q}$ are the \textbf{null forms relative to g}, defined by
	\begin{subequations}
		\begin{align}
		\mathscr{Q}^i
		& := 
			-(g^{-1})^{\alpha \beta} \partial_{\alpha} \Densrenormalized \partial_{\beta} v^i,
			\label{E:VELOCITYNULLFORM} \\
		\mathscr{Q}
		& := 
		- 
		2 \Speed^{-1} \Speed' 
		(g^{-1})^{\alpha \beta} \partial_{\alpha} \Densrenormalized \partial_{\beta} \Densrenormalized
		+ 
		2 \left\lbrace
				\partial_1 v^1 \partial_2 v^2
				-
				\partial_2 v^1 \partial_1 v^2 
			\right\rbrace.
			\label{E:DENSITYNULLFORM}
		\end{align}
	\end{subequations}
\end{proposition}

\begin{remark}[\textbf{The importance of the null forms relative to $g$}]
	\label{R:NULLFORMSAREIMPORTANT}
	For the proof of our main theorem, 
	it is critically important that $\mathscr{Q}^i$ and $\mathscr{Q}$ are null forms relative to $g$.
	The reason is that, due to their special structure,
	$\upmu \mathscr{Q}^i$ and $\upmu \mathscr{Q}$
	remain uniformly small, all the way up to the shock. Thus, they do not
	interfere with the singularity formation mechanisms.
	In contrast, a general quadratic term
	$\upmu (\partial v,\partial \Densrenormalized) \cdot (\partial v,\partial \Densrenormalized)$
	could become large near the expected singularity 
	and dominate the dynamics;
	had such a term been present in the equations,
	it would have completely obstructed our approach.
\end{remark}

\subsection{Constant state background solutions and the array of solution variables}
We will study perturbations of the following constant state background solution
to the system \eqref{E:VELOCITYWAVEEQUATION}-\eqref{E:RENORMALIZEDVORTICTITYTRANSPORTEQUATION}:
\begin{align} \label{E:BACKGROUNDSOLUTION}
	(\Densrenormalized,v^1,v^2,\Vortrenormalized)
	\equiv
	(0,0,0,0).
\end{align}
The solution \eqref{E:BACKGROUNDSOLUTION}
corresponds to a motionless fluid of constant
density $\bar{\rho}$,
where $\bar{\rho} > 0$ is a constant.
Note that a more general constant state
$(\Densrenormalized,v^1,v^2,\Vortrenormalized) \equiv (0,\bar{v}^1,\bar{v}^2,0)$,
the $\bar{v}^i$ are constants,
may be brought into the form \eqref{E:BACKGROUNDSOLUTION}
via a Galilean transformation.\footnote{By this, we mean the change of coordinates
$\widetilde{t} := t$ and
$\widetilde{x}^i := x^i - \bar{v}^i t$,
which implies that 
$
\displaystyle
\frac{\partial}{\partial \widetilde{t}} = \partial_t + \bar{v}^a \partial_a$
and
$
\displaystyle
\frac{\partial}{\partial \widetilde{x}^i} = \partial_i
$.
Note that the expression
$x^2 - \bar{v}^2 t$
should be interpreted as
the translation of the point $x^2 \in \mathbb{T}$
by the flow of of the vectorfield 
$-\bar{v}^2 \partial_2$ 
for $t$ units of time.
\label{FN:GALILEAN}} 
Let 
\begin{align} \label{E:SPEEDOFSOUNDOFBACKGROUND}
	\bar{\Speed}
	& := \Speed(\Densrenormalized = 0)
\end{align}
denote the speed of sound \eqref{E:SOUNDSPEED}
evaluated at the background solution \eqref{E:BACKGROUNDSOLUTION}.
Without loss of generality, we assume\footnote{We can always ensure the condition 
\eqref{E:SPEEDOFSOUNDISUNITY}
by making the following changes of variables:
\[
	\widetilde{v}^i 
	= \frac{v^i}{\bar{\Speed}},
		\qquad
	\widetilde{t}
	= \bar{\Speed} t,
		\qquad
	\widetilde{g} = \bar{\Speed}^2 g,
		\qquad
	\widetilde{\Speed} = \frac{\Speed}{\bar{\Speed}}.
\]
These changes of variables
leave the expressions
\eqref{E:ACOUSTICALMETRIC}-\eqref{E:INVERSEACOUSTICALMETRIC}
and the Euler equations 
\eqref{E:TRANSPORTDENSRENORMALIZEDRELATIVETORECTANGULAR}-\eqref{E:TRANSPORTVELOCITYRELATIVETORECTANGULAR}
invariant and are such that the desired normalization
$\widetilde{\Speed}(\Densrenormalized = 0) = 1$
holds.} 
that
\begin{align} \label{E:SPEEDOFSOUNDISUNITY}
	\bar{\Speed}
	& = 1.
\end{align}
The advantage of the assumption \eqref{E:SPEEDOFSOUNDISUNITY}
is that it simplifies many of our formulas.

Many of our estimates will apply uniformly to 
the ``wave variables'' $\Densrenormalized$,
$v^1$, and $v^2$. For this reason, we collect them
into an array.
\begin{definition}[\textbf{The array} $\threePsi$ \textbf{of wave variables}]
	\label{D:VECPSI}
		\begin{align} \label{E:THREEVECPSI}
			\threePsi
			= 
			(\Psi_0,\Psi_1,\Psi_2)
			:=
			(\Densrenormalized,
				v^1,
				v^2).
			\end{align}
\end{definition}
Since we are studying perturbations of the solution \eqref{E:BACKGROUNDSOLUTION},
we may think of $\threePsi$ as small.
However, for the solutions under study, some of the derivatives 
of $\threePsi$ are relatively large.

\begin{remark}[$\threePsi$ \textbf{is not a tensor}]
	\label{R:PSINOTTENSOR}
	Throughout, we view $\threePsi$ to be an array of
	scalar-valued functions; we will not attribute any 
	tensorial structure to the labeling index of $\Psi_{\imath}$
	besides simple contractions, denoted by $\contr$,
	corresponding to the chain rule; see Def.~\ref{D:BIGGANDBIGH}.
\end{remark}

\subsection{The metric components and their derivatives with respect to the solution}
\label{SS:METRICCOMPONENTSANDTHEIRDERIVATIVE}
Throughout the paper, we often view 
the Cartesian metric component functions $g_{\alpha \beta}$
(see \eqref{E:ACOUSTICALMETRIC})
to be (explicitly known) functions of $\threePsi$: 
$g_{\alpha \beta} = g_{\alpha \beta}(\threePsi)$.
From the expression \eqref{E:ACOUSTICALMETRIC}
and the assumption \eqref{E:SPEEDOFSOUNDISUNITY},
it follows that we can decompose
\begin{align} \label{E:LITTLEGDECOMPOSED}
	g_{\alpha \beta}(\threePsi)
	& = m_{\alpha \beta} 
		+ g_{\alpha \beta}^{(Small)}(\threePsi),
	&& (\alpha, \beta = 0,1,2),
\end{align}
where 
\begin{align} \label{E:STMINKOWKMET}
	m_{\alpha \beta} = \mbox{\upshape diag}(-1,1,1)
\end{align}
is the standard Minkowski metric and
$g_{\alpha \beta}^{(Small)}(\threePsi)$ is a smooth function of $\threePsi$ with
\begin{align} \label{E:METRICPERTURBATIONFUNCTION}
	g_{\alpha \beta}^{(Small)}(\vec{0})
	& = 0.
\end{align}
Specifically, we have the formula
\begin{align} \label{E:GSMALLTENSOR}
g^{(Small)} 
		& = 
				\sum_{a=1}^2
				\Speed^{-2} (v^a)^2
				dt \otimes dt
			- 
			\Speed^{-2}
			\sum_{a=1}^2 v^a dt \otimes dx^a
			- 
			\Speed^{-2}
			\sum_{i=1}^2 v^a dx^a \otimes dt
				\\
		& \ \
			+
			\left\lbrace
				\Speed^{-2}
				-
				1
			\right\rbrace
			\sum_{a=1}^2 dx^a \otimes dx^a.
			\notag
	\end{align}

The following quantities arise in many of the equations that we study.

\begin{definition}[\textbf{Derivatives of $g_{\alpha \beta}$ with respect to $\threePsi$}]
	\label{D:BIGGANDBIGH}
	For $\alpha,\beta = 0,1,2$ and $\imath,\jmath = 0,1,2$, we define
	\begin{subequations}
	\begin{align}
		G_{\alpha \beta}^{\imath}(\threePsi)
		& := \frac{\partial}{\partial \Psi_{\imath}} 
					g_{\alpha \beta}(\threePsi),
			\label{E:BIGGDEF} \\
		\vec{G}_{\alpha \beta}
		=
		\vec{G}_{\alpha \beta}(\threePsi)
		& := 
		\left(
			G_{\alpha \beta}^0(\threePsi),
			G_{\alpha \beta}^1(\threePsi),
			G_{\alpha \beta}^2(\threePsi)
		\right),
			\label{E:BIGGARRAY} \\
		H_{\alpha \beta}^{\imath \jmath}(\threePsi)
		& := \frac{\partial}{\partial \Psi_{\imath}} 
				 \frac{\partial}{\partial \Psi_{\jmath}}
					g_{\alpha \beta}(\threePsi),
			\label{E:BIGHDEF} \\
		\vec{H}_{\alpha \beta}
		=
		\vec{H}_{\alpha \beta}(\threePsi)
		&:= 
		\left(
			H_{\alpha \beta}^{00}(\threePsi),
			H_{\alpha \beta}^{01}(\threePsi),
			\cdots,
			H_{\alpha \beta}^{22}(\threePsi)
		\right).
		\label{E:BIGHARRAY}
	\end{align}
	\end{subequations}
\end{definition}
For $\imath = 0,1,2$, we think of
the $\lbrace G_{\alpha \beta}^{\imath} \rbrace_{\alpha,\beta = 0,1,2}$,
as the Cartesian components of a spacetime tensorfield. 
Similarly, we think of
$\lbrace \vec{G}_{\alpha \beta} \rbrace_{\alpha,\beta = 0,1,2}$
as the Cartesian components of an array-valued spacetime tensorfield.
Similar remarks apply to $H_{\alpha \beta}^{\imath \jmath}$
and $\vec{H}_{\alpha \beta}$.

The following operators naturally arise in our analysis of solutions.

\begin{definition}[\textbf{Operators involving the array $\threePsi$}]
	\label{D:ARRAYNOTATION}
	Let $U,V,U_1,U_2,V_1,V_2$ be vectorfields. We define
	\begin{subequations}
	\begin{align} \label{E:VAPPLIEDTOVECPSI}
		V \threePsi
		& := 
		(V \Psi_0, V \Psi_1, V \Psi_2),
				\\
	\vec{G}_{U_1 U_2}\contr V \threePsi
	& := \sum_{\imath=0}^2 G_{\alpha \beta}^{\imath}U_1^{\alpha} U_2^{\beta} V \Psi_{\imath},
		\label{E:VECGVECTORFIELDCONTRACTEDANDCONTRACEDWITHVVECPSI} \\
		\vec{H}_{U_1 U_2} \contr \contr(V_1 \threePsi) V_2 \threePsi
	& := \sum_{\imath,\jmath=0}^2 H_{\alpha \beta}^{\imath \jmath} 
				U_1^{\alpha} U_2^{\beta} (V_1 \Psi_{\imath}) V_2 \Psi_{\jmath}.
				\label{E:VECHVECTORFIELDCONTRACTEDANDCONTRACEDWITHVVECPSI}
	\end{align}
	\end{subequations}

	We use similar notation with other differential operators in place of vectorfield differentiation.
	For example, 
	$\vec{G}_{U_1 U_2} \contr\angLap \threePsi 
	:= \sum_{\imath=0}^2 G_{\alpha \beta}^{\imath} U_1^{\alpha} U_2^{\beta} \angLap \Psi_{\imath}
	$.

\end{definition}

\subsection{The eikonal function and related constructions}
\label{SS:EIKONALFUNCTIONANDRELATED}
To track the solution all the way to the shock,
we construct a new set of geometric coordinates, one of
which is the eikonal function.
\begin{definition}[\textbf{Eikonal function}]
\label{D:INTROEIKONAL}
The eikonal function $u$ 
solves the eikonal equation initial value problem
\begin{align} \label{E:INTROEIKONAL}
	(g^{-1})^{\alpha \beta}
	\partial_{\alpha} u \partial_{\beta} u
	& = 0, 
	\qquad \partial_t u > 0,
		\\
	u|_{\Sigma_0}
	& = 1 - x^1,
	\label{E:INTROEIKONALINITIALVALUE}
\end{align}
where $\Sigma_0$ is the hypersurface of constant Cartesian time $0$.
\end{definition}

Using $u$, we can construct many
geometric quantities that can be used to derive sharp information about
the solution. We start by defining the most important quantity in the study of shock formation:
the \emph{inverse foliation density}.

\begin{definition}[\textbf{Inverse foliation density}]
\label{D:FIRSTUPMU}
We define the inverse foliation density $\upmu$ as follows:
\begin{align} \label{E:FIRSTUPMU}
	\upmu 
	& := \frac{-1}{(g^{-1})^{\alpha \beta}(\threePsi) \partial_{\alpha} t \partial_{\beta} u} > 0,
\end{align}
where $t$ is the Cartesian time coordinate.
We note that the identity \eqref{E:MATERIALRECTANGULAR} below implies that
$
\displaystyle
\upmu 
= \frac{1}{\Transport u}
$,
where $\Transport$ is the material derivative vectorfield defined in 
\eqref{E:MATERIALVECTORVIELDRELATIVETORECTANGULAR}.
\end{definition}

The quantity $1/\upmu$ measures the density of the level sets of $u$
relative to the constant-time hypersurfaces $\Sigma_t$. When $\upmu$
becomes $0$, the density becomes infinite and the level sets of $u$
intersect. For the initial data under consideration, 
$\upmu$ starts out near unity.
It turns out that the formation of the shock, the blowup
of the eikonal function's first Cartesian coordinate partial derivatives,
and the blowup of the first derivatives of 
$v$ and $\Densrenormalized$ with respect to the Cartesian coordinate partial derivatives
are all simultaneously tied to the vanishing of $\upmu$.
We also note that the vanishing of $\upmu$
is equivalent to the blowup of $\Transport u$.

We now define the spacetime subsets on which we analyze solutions.
They are depicted in Fig.~\ref{F:SOLIDREGION}
on pg.~\pageref{F:SOLIDREGION}.

\begin{definition} [\textbf{Subsets of spacetime}]
\label{D:HYPERSURFACESANDCONICALREGIONS}
For $1 \leq t'$ and $0 \leq u \leq U_0$, we define the following subsets of spacetime:
\begin{subequations}
\begin{align}
	\Sigma_{t'} & := \lbrace (t,x^1,x^2) \in \mathbb{R} \times \mathbb{R} \times \mathbb{T}  
		\ | \ t = t' \rbrace, 
		\label{E:SIGMAT} \\
	\Sigma_{t'}^{u'} & := \lbrace (t,x^1,x^2) \in \mathbb{R} \times \mathbb{R} \times \mathbb{T} 
		 \ | \ t = t', \ 0 \leq u(t,x^1,x^2) \leq u' \rbrace, 
		\label{E:SIGMATU} 
		\\
	\mathcal{P}_{u'}^{t'} 
	& := 
		\lbrace (t,x^1,x^2) \in \mathbb{R} \times \mathbb{R} \times \mathbb{T} 
			\ | \ 1 \leq t \leq t', \ u(t,x^1,x^2) = u' 
		\rbrace, 
		\label{E:PUT} \\
	\ell_{t',u'} 
		&:= \mathcal{P}_{u'}^{t'} \cap \Sigma_{t'}^{u'}
		= \lbrace (t,x^1,x^2) \in \mathbb{R} \times \mathbb{R} \times \mathbb{T} 
			\ | \ t = t', \ u(t,x^1,x^2) = u' \rbrace, 
			\label{E:LTU} \\
	\mathcal{M}_{t',u'} & := \cup_{u \in [0,u']} \mathcal{P}_u^{t'} \cap \lbrace (t,x^1,x^2) 
		\in \mathbb{R} \times \mathbb{R} \times \mathbb{T}  \ | \ 1 \leq t < t' \rbrace.
		\label{E:MTUDEF}
\end{align}
\end{subequations}
\end{definition}
We refer to the $\Sigma_t$ and $\Sigma_t^u$ as ``constant time slices,'' 
the $\mathcal{P}_u^t$ as ``null hyperplanes,''
and the $\ell_{t,u}$ as ``curves'' or ``tori.'' 
We sometimes use the notation $\mathcal{P}_u$ in place of $\mathcal{P}_u^t$ 
when we are not concerned with the truncation time $t$.
Note that $\mathcal{M}_{t,u}$ is ``open at the top'' by construction.

We now construct a local coordinate function on the tori $\ell_{t,u}$.

\begin{definition}[\textbf{Geometric torus coordinate}]
\label{D:GEOMETRICTORUSCOORDINATE}
We define the geometric torus coordinate $\vartheta$
to be the solution to the following transport equation:
\begin{align} \label{E:APPENDIXGEOMETRICTORUSCOORD}
	(g^{-1})^{\alpha \beta}
	\partial_{\alpha} u \partial_{\beta} \vartheta
	& = 0,
		\\
	\vartheta|_{\Sigma_0}
	& = x^2. \label{E:APPENDIXGEOMETRICTORUSCOORDINITIALCOND}
\end{align}
\end{definition}

\begin{definition}[\textbf{Geometric coordinates and partial derivatives}]
	\label{D:GEOMETRICCOORDINATES}
	We refer to $(t,u,\vartheta)$ as the geometric coordinates,
	where $t$ is the Cartesian time coordinate.
	We denote the corresponding geometric coordinate partial derivative vectorfields by
	\begin{align} \label{E:GEOCOORDPARTDERIVVECTORFIELDS}
		\left\lbrace 
			\frac{\partial}{\partial t}, 
			\frac{\partial}{\partial u},
			\CoordAng := \frac{\partial}{\partial \vartheta}
		\right\rbrace.
	\end{align}
\end{definition}

\begin{remark}
	$\CoordAng$ is globally defined even though
	$\vartheta$ is only locally defined along $\ell_{t,u}$.
\end{remark}

\begin{definition}\label{D:CHOVMAP}
	We define $\Upsilon: [0,T) \times [0,U_0] \times \mathbb{T} \rightarrow \mathcal{M}_{T,U_0}$,
	$\Upsilon(t,u,\vartheta) := (t,x^1,x^2)$,
	to be the change of variables map from geometric to Cartesian coordinates.
\end{definition}

\begin{remark}[\textbf{$C^1$-equivalent differential structures until shock formation}]
We often identify
spacetime regions of the form
$\mathcal{M}_{t,U_0}$ 
(see \eqref{E:MTUDEF})
with the region
$[0,t) \times [0,U_0] \times \mathbb{T}$
corresponding to the geometric coordinates.
This identification is 
justified by the fact that 
during the classical lifespan of the solutions under consideration,
the differential structure on 
$\mathcal{M}_{t,U_0}$ corresponding to the 
geometric coordinates is $C^1$-equivalent
to the differential structure on 
$\mathcal{M}_{t,U_0}$
corresponding to the Cartesian coordinates.
The reason is that $\Upsilon$
is $C^1$ with a $C^1$ inverse until a shock forms;
this fact was proved in \cite{jSgHjLwW2016}*{Theorem 15.1}
and is revisited in limited form 
in the proof of Theorem~\ref{T:MAINTHEOREM}.
In contrast, at points where $\upmu$ vanishes,
the partial derivatives of 
$\Densrenormalized$
and
$v^1$ with respect to the Cartesian coordinates blow up,
the inverse map $\Upsilon^{-1}$ becomes singular,
and the equivalence of the differential structures breaks down.
\end{remark}

\subsection{Important vectorfields, the rescaled frame, and the unit frame}
\label{SS:FRAMEANDRELATEDVECTORFIELDS}
In this section, we define some 
vectorfields that we use in our analysis
and exhibit their basic properties.

We start by defining the (negative) gradient vectorfield
associated to the eikonal function:
\begin{align} \label{E:LGEOEQUATION}
	\Lgeo^{\nu} 
	& := - (g^{-1})^{\nu \alpha} \partial_{\alpha} u.
\end{align}
It is easy to see that $\Lgeo$ is future-directed\footnote{Here and throughout, 
a vectorfield $V$ is ``future-directed'' if its Cartesian component $V^0$ is positive.} 
with
\begin{align}  \label{E:LGEOISNULL}
g(\Lgeo,\Lgeo) 
:= g_{\alpha \beta} \Lgeo^{\alpha} \Lgeo^{\beta}
= 0,
\end{align}
that is, $\Lgeo$ is $g$-null.
Moreover, we can differentiate the eikonal equation with
$\D^{\nu} := (g^{-1})^{\nu \alpha} \D_{\alpha}$
and use the torsion-free property of the connection $\D$ to deduce that
$
0
= 
(g^{-1})^{\alpha \beta} \D_{\alpha} u \D_{\beta} \D^{\nu} u
= - \D^{\alpha} u \D_{\alpha} \Lgeo^{\nu}
= \Lgeo^{\alpha} \D_{\alpha} \Lgeo^{\nu}
$. That is, $\Lgeo$ is geodesic:
\begin{align} \label{E:LGEOISGEODESIC}
	\D_{\Lgeo} \Lgeo & = 0.
\end{align}
In addition, since $\Lgeo$ is proportional to the metric dual of the one-form $d u$,
which is co-normal to the level sets $\mathcal{P}_u$ of the eikonal function,
it follows that $\Lgeo$ is $g$-orthogonal to $\mathcal{P}_u$.
Hence, the $\mathcal{P}_u$ have null normals. 
Such hypersurfaces are known as \emph{null hypersurfaces}
or \emph{characteristics}.
Our analysis will show that the Cartesian components of $\Lgeo$
blow up when the shock forms. 

In our analysis, we work with a rescaled version of
$\Lgeo$ that we denote by $\Lunit$.  
Our proof reveals that the Cartesian components of $\Lunit$
remain near those of 
$\Lunit_{(Flat)}
	:= 
	\partial_t + \partial_1
$ 
all the way up to the shock.

\begin{definition}[\textbf{Rescaled null vectorfield}]
	\label{D:LUNITDEF}
	We define the rescaled null (see \eqref{E:LGEOISNULL}) vectorfield $\Lunit$ as follows:
	\begin{align} \label{E:LUNITDEF}
		\Lunit
		& := \upmu \Lgeo.
	\end{align}
\end{definition}
Note that $\Lunit$ is $g$-null since $\Lgeo$ is.
We also note that by \eqref{E:APPENDIXGEOMETRICTORUSCOORD}, we have
\begin{align} \label{E:LUNITTHETA}
	\Lunit \vartheta 
	& = 0.
\end{align}

We now define the vectorfields $\Rad$ and $\Radunit$, 
which are transversal to the characteristics
$\mathcal{P}_u$. It is critically important for our
work that $\Rad$ is rescaled by a factor of $\upmu$.

\begin{definition}[$\Radunit$ \textbf{and} $\Rad$]
	\label{D:RADANDXIDEFS}
	We define $\Radunit$ to be the unique vectorfield
	that is $\Sigma_t$-tangent, $g$-orthogonal 
	to the $\ell_{t,u}$, and normalized by
	\begin{align} \label{E:GLUNITRADUNITISMINUSONE}
		g(\Lunit,\Radunit) = -1.
	\end{align}
	We define
	\begin{align} \label{E:RADDEF}
		\Rad := \upmu \Radunit.
	\end{align}
\end{definition}

We use the following two vectorfield frames in our analysis.

\begin{definition}[\textbf{Two frames}]
	\label{D:RESCALEDFRAME}
	We define, respectively, the rescaled frame and the non-rescaled frame as follows:
	\begin{subequations}
	\begin{align} \label{E:RESCALEDFRAME}
		& \lbrace \Lunit, \Rad, \CoordAng \rbrace,
		&& \mbox{Rescaled frame},	\\
		&\lbrace \Lunit, \Radunit, \CoordAng \rbrace,
		&& \mbox{Non-rescaled frame}.
			\label{E:UNITFRAME}
	\end{align}
	\end{subequations}
\end{definition}

In the next lemma, we exhibit the basic properties
of some of the vectorfields that we have defined.

\begin{lemma}[\textbf{Basic properties of} $\Radunit$, $\Rad$, $\Lunit$, \textbf{and} $\Transport$]
\label{L:BASICPROPERTIESOFFRAME}
The following identities hold:
\begin{subequations}
\begin{align} \label{E:LUNITOFUANDT}
	\Lunit u & = 0, \qquad \Lunit t = \Lunit^0 = 1,
		\\
	\Rad u & = 1,
	\qquad \Rad t = \Rad^0 = 0, \label{E:RADOFUANDT}
\end{align}
\end{subequations}

\begin{subequations}
\begin{align} \label{E:RADIALVECTORFIELDSLENGTHS}
	g(\Radunit,\Radunit)
	& = 1,
	\qquad
	g(\Rad,\Rad)
	= \upmu^2,
		\\
	g(\Lunit,\Radunit)
	& = -1,
	\qquad
	g(\Lunit,\Rad) = -\upmu.
	\label{E:LRADIALVECTORFIELDSNORMALIZATIONS}
\end{align}
\end{subequations}

Moreover, relative to the geometric coordinates, we have
\begin{align} \label{E:LISDDT}
	\Lunit = \frac{\partial}{\partial t}.
\end{align}

In addition, there exists an $\ell_{t,u}$-tangent vectorfield
$\NonRadialRad = \XiCoordComp \CoordAng$ 
(where $\XiCoordComp$ is a scalar-valued function)
such that
\begin{align} \label{E:RADSPLITINTOPARTTILAUANDXI}
	\Rad 
	& = \frac{\partial}{\partial u} - \NonRadialRad
	= \frac{\partial}{\partial u} - \XiCoordComp \CoordAng.
\end{align}

The material derivative vectorfield 
$\Transport$ defined in \eqref{E:MATERIALVECTORVIELDRELATIVETORECTANGULAR}
is future-directed, $g$-orthogonal
to $\Sigma_t$ and is normalized by
\begin{align} \label{E:TRANSPROTUNITLENGTH}
	g(\Transport,\Transport) 
	& = - 1.
\end{align}
In addition, relative to Cartesian coordinates, we have
(for $\nu = 0,1,2$):
\begin{align} \label{E:MATERIALRECTANGULAR}
	\Transport^{\nu} = - (g^{-1})^{0 \nu}.
\end{align}
Moreover, we have
\begin{align} \label{E:TRANSPORTVECTORFIELDINTERMSOFLUNITANDRADUNIT}
	\Transport
	& = \Lunit + \Radunit.
\end{align}

	Finally, the following identities hold relative to the Cartesian coordinates 
	(for $\nu = 0,1,2$):
	\begin{align}  \label{E:DOWNSTAIRSUPSTAIRSSRADUNITPLUSLUNITISAFUNCTIONOFPSI} 
		\Radunit_{\nu} 
		& = - \Lunit_{\nu} - \delta_{\nu}^0,
		\qquad
		\Radunit^{\nu}
		= - \Lunit^{\nu}
			- (g^{-1})^{0\nu},
	\end{align}
	where $\delta_{\nu}^0$ is the standard Kronecker delta.
\end{lemma}
\begin{proof}
The identity \eqref{E:MATERIALRECTANGULAR}
follows trivially from \eqref{E:INVERSEACOUSTICALMETRIC}.
The remaining statements in the lemmas were proved in
\cite{jSgHjLwW2016}*{Lemma 2.1},
where the vectorfield $\Transport$ was denoted by $N$.
\end{proof}

\subsection{Projection tensorfields, \texorpdfstring{$\vec{G}_{(Frame)}$}{frame components}, and projected Lie derivatives}
\label{SS:PROJECTIONTENSORFIELDANDPROJECTEDLIEDERIVATIVES}
Many of our constructions involve projections onto $\Sigma_t$ and $\ell_{t,u}$.

\begin{definition}[\textbf{Projection tensorfields}]
We define the $\Sigma_t$ projection tensorfield $\Sigmatproject$
and the $\ell_{t,u}$ projection tensorfield
$\Lineproject$ relative to Cartesian coordinates as follows:
\begin{subequations}
\begin{align} 
	\Sigmatproject_{\nu}^{\ \mu} 
	&:=	\delta_{\nu}^{\ \mu}
			- \Transport_{\nu} \Transport^{\mu} 
		= \delta_{\nu}^{\ \mu}
			+ \delta_{\nu}^{\ 0} \Lunit^{\mu}
			+ \delta_{\nu}^{\ 0} \Radunit^{\mu},
			\label{E:SIGMATPROJECTION} \\
	\Lineproject_{\nu}^{\ \mu} 
	&:=	\delta_{\nu}^{\ \mu}
			+ \Radunit_{\nu} \Lunit^{\mu} 
			+ \Lunit_{\nu} (\Lunit^{\mu} + \Radunit^{\mu}) 
		= \delta_{\nu}^{\ \mu}
			- \delta_{\nu}^{\ 0} \Lunit^{\mu} 
			+  \Lunit_{\nu} \Radunit^{\mu}.
			\label{E:LINEPROJECTION}
	\end{align}
\end{subequations}
In \eqref{E:SIGMATPROJECTION}-\eqref{E:LINEPROJECTION},
$\delta_{\nu}^{\ \mu}$ is the standard Kronecker delta.
\end{definition}

\begin{definition}[\textbf{Projections of tensorfields}]
Given any spacetime tensorfield $\xi$,
we define its $\Sigma_t$ projection $\Sigmatproject \xi$
and its $\ell_{t,u}$ projection $\Lineproject \xi$
as follows:
\begin{subequations}
\begin{align} 
(\Sigmatproject \xi)_{\nu_1 \cdots \nu_n}^{\mu_1 \cdots \mu_m}
& :=
	\Sigmatproject_{\widetilde{\mu}_1}^{\ \mu_1} \cdots \Sigmatproject_{\widetilde{\mu}_m}^{\ \mu_m}
	\Sigmatproject_{\nu_1}^{\ \widetilde{\nu}_1} \cdots \Sigmatproject_{\nu_n}^{\ \widetilde{\nu}_n} 
	\xi_{\widetilde{\nu}_1 \cdots \widetilde{\nu}_n}^{\widetilde{\mu}_1 \cdots \widetilde{\mu}_m},
		\\
(\Lineproject \xi)_{\nu_1 \cdots \nu_n}^{\mu_1 \cdots \mu_m}
& := 
	\Lineproject_{\widetilde{\mu}_1}^{\ \mu_1} \cdots \Lineproject_{\widetilde{\mu}_m}^{\ \mu_m}
	\Lineproject_{\nu_1}^{\ \widetilde{\nu}_1} \cdots \Lineproject_{\nu_n}^{\ \widetilde{\nu}_n} 
	\xi_{\widetilde{\nu}_1 \cdots \widetilde{\nu}_n}^{\widetilde{\mu}_1 \cdots \widetilde{\mu}_m}.
	\label{E:STUPROJECTIONOFATENSOR}
\end{align}
\end{subequations}
\end{definition}
We say that a spacetime tensorfield $\xi$ is $\Sigma_t$-tangent 
(respectively $\ell_{t,u}$-tangent)
if $\Sigmatproject \xi = \xi$
(respectively if $\Lineproject \xi = \xi$).
Alternatively, we say that $\xi$ is a
$\Sigma_t$ tensor (respectively $\ell_{t,u}$ tensor).


\begin{definition}[\textbf{$\ell_{t,u}$ projection notation}]
	\label{D:STUSLASHPROJECTIONNOTATION}
	If $\xi$ is a spacetime tensor, then we define
	\begin{align} \label{E:TENSORSTUPROJECTED}
		\angxi := \Lineproject \xi.
	\end{align}

	If $\xi$ is a symmetric type $\binom{0}{2}$ spacetime tensor and $V$ is a spacetime
	vectorfield, then we define
	\begin{align} \label{E:TENSORVECTORANDSTUPROJECTED}
		\angxiarg{V} 
		& := \Lineproject (\xi_V),
	\end{align}
	where $\xi_V$ is the spacetime one-form with 
	Cartesian components $\xi_{\alpha \nu} V^{\alpha}$, $(\nu = 0,1,2)$.
\end{definition}

Throughout, $\Lie_V \xi$ denotes the Lie derivative of the tensorfield
$\xi$ with respect to the vectorfield $V$.
We often use the Lie bracket notation
$[V,W] := \Lie_V W$ when $V$ and $W$ are vectorfields.

In our analysis, we will apply the Leibniz rule for Lie derivatives to 
contractions of tensor products of $\ell_{t,u}$-tensorfields.
Due in part to the special properties 
(such as \eqref{E:LIELANDLIERADOFELLTUTANGENTISELLTUTANGENT})
of the vectorfields that we use to differentiate,
the non-$\ell_{t,u}$ components of the differentiated factor in the products typically cancel.
This motivates the following definition.

\begin{definition}[$\ell_{t,u}$ \textbf{and} $\Sigma_t$-\textbf{projected Lie derivatives}]
\label{D:PROJECTEDLIE}
Given a tensorfield $\xi$
and a vectorfield $V$,
we define the $\Sigma_t$-projected Lie derivative
$\SigmatLie_V \xi$ of $\xi$
and the $\ell_{t,u}$-projected Lie derivative
$\angLie_V \xi$ of $\xi$ as follows:
\begin{align} 
	\SigmatLie_V \xi
	& := \Sigmatproject \Lie_V \xi,
		\qquad
	\angLie_V \xi
	:= \Lineproject \Lie_V \xi.
	\label{E:PROJECTIONS}
\end{align}
\end{definition}

\begin{definition}[\textbf{Components of $\vec{G}$ and $\vec{H}$ relative to the non-rescaled frame}]
	\label{D:GFRAMEANDHFRAMEARRAYS}
	We define
	\begin{align}
		\vec{G}_{(Frame)} 
		& := 
			\left\lbrace 
				\vec{G}_{\Lunit \Lunit}, \vec{G}_{\Lunit \Radunit}, \angGarg{\Lunit}, \angGarg{\Radunit}, \angG 
			\right\rbrace,
		\qquad
		\vec{H}_{(Frame)} 
		:= \left\lbrace 
					\vec{H}_{\Lunit \Lunit}, \vec{H}_{\Lunit \Radunit}, \angHarg{\Lunit}, \angHarg{\Radunit}, \angH 
				\right\rbrace
				\label{E:HFRAME}
	\end{align}
	to be the arrays of components of the tensorfield arrays 
	\eqref{E:BIGGARRAY} and \eqref{E:BIGHARRAY}
	relative the non-rescaled frame \eqref{E:UNITFRAME}.
\end{definition}

We adopt the convention that when we differentiate $\vec{G}_{(Frame)}$ 
or $\vec{H}_{(Frame)}$,
we by definition form a new array consisting of the differentiated components.
For example,
\begin{align}\label{E:GFRAMEDIFFERENTIATEDEXAMPLE}
		\angLie_{\Lunit} \vec{G}_{(Frame)} 
		& := 
			\left\lbrace 
				\Lunit(\vec{G}_{\Lunit \Lunit}), 
				\Lunit (\vec{G}_{\Lunit \Radunit}), 
				\angLie_{\Lunit} (\angGarg{\Lunit}), 
				\angLie_{\Lunit} (\angGarg{\Radunit}), 
				\angLie_{\Lunit} \angG
			\right\rbrace,
\end{align}
where
$
\displaystyle
\Lunit(\vec{G}_{\Lunit \Lunit})
:= \left\lbrace
					\Lunit(G_{\Lunit \Lunit}^0),
					\Lunit(G_{\Lunit \Lunit}^1),
					\Lunit(G_{\Lunit \Lunit}^2)
			 \right\rbrace
$,
$
\displaystyle
	\angLie_{\Lunit} (\angGarg{\Radunit})
	:= \left\lbrace
					\angLie_{\Lunit} (\NovecangGarg{\Radunit}{0}),
					\angLie_{\Lunit} (\NovecangGarg{\Radunit}{1}),
					\angLie_{\Lunit} (\NovecangGarg{\Radunit}{2})
			 \right\rbrace,
$
etc.

\subsection{First and second fundamental forms and covariant differential operators}
\label{SS:FUNDFORMSANDCOVDERIVOPS}

\begin{definition}[\textbf{First fundamental forms}] \label{D:FIRSTFUND}
	We define the first fundamental form $\gt$ of $\Sigma_t$ and the 
	first fundamental form $\gsphere$ of $\ell_{t,u}$ as follows:
	\begin{align}
		\gt
		:= \Sigmatproject g,
		\qquad
		\gsphere
		:= \Lineproject g.
		\label{E:GTANDGSPHERESPHEREDEF}
	\end{align}

	We define the inverse first fundamental forms by raising the indices with $g^{-1}$:
	\begin{align}
		(\gt^{-1})^{\mu \nu}
		:= (g^{-1})^{\mu \alpha} (g^{-1})^{\nu \beta} \gt_{\alpha \beta},
		\qquad 
		(\gsphere^{-1})^{\mu \nu}
		:= (g^{-1})^{\mu \alpha} (g^{-1})^{\nu \beta} \gsphere_{\alpha \beta}.
		\label{E:GGTINVERSEANDGSPHEREINVERSEDEF}
	\end{align}
\end{definition}
Note that $\gt$ is the Riemannian metric on $\Sigma_t$ induced by $g$
and that $\gsphere$ is the Riemannian metric on $\ell_{t,u}$ induced by $g$.
Moreover, simple calculations yield
$(\gt^{-1})^{\mu \alpha} \gt_{\alpha \nu} = \Sigmatproject_{\nu}^{\ \mu}$
and $(\gsphere^{-1})^{\mu \alpha} \gsphere_{\alpha \nu} = \Lineproject_{\nu}^{\ \mu}$.

\begin{remark}
	\label{E:ELLTUTENSORSAREPURETRACE}
	Because the $\ell_{t,u}$ are one-dimensional manifolds, it follows that
	symmetric type $\binom{0}{2}$ $\ell_{t,u}$-tangent tensorfields $\xi$
	satisfy $\xi = (\mytr \xi) \gsphere$,
	where $\mytr \xi := \ginversesphere \cdot \xi$.
	This simple fact simplifies some of our formulas
	compared to the case of higher space dimensions.
	In the remainder of the article, we often use this fact without
	explicitly mentioning it. 
\end{remark}



\begin{definition}[\textbf{Differential operators associated to the metrics}] 
\label{D:CONNECTIONS}
	We use the following notation for various differential operators associated 
	to the spacetime metric $g$
	and the Riemannian metric $\gsphere$ induced on $\ell_{t,u}$.
	\begin{itemize}
		\item $\D$ denotes the Levi-Civita connection of the acoustical metric $g$.
		\item $\angD$ denotes the Levi-Civita connection of $\gsphere$.
		\item If $\xi$ is an $\ell_{t,u}$-tangent one-form,
			then $\angdiv \xi$ is the scalar-valued function
			$\angdiv \xi := \ginversesphere \cdot \angD \xi$.
		\item Similarly, if $V$ is an $\ell_{t,u}$-tangent vectorfield,
			then $\angdiv V := \ginversesphere \cdot \angD V_{\flat}$,
			where $V_{\flat}$ is the one-form $\gsphere$-dual to $V$.
		\item If $\xi$ is a symmetric type $\binom{0}{2}$ 
		 $\ell_{t,u}$-tangent tensorfield, then
		 $\angdiv \xi$ is the $\ell_{t,u}$-tangent 
		 one-form $\angdiv \xi := \ginversesphere \cdot \angD \xi$,
		 where the two contraction indices in $\angD \xi$
		 correspond to the operator $\angD$ and the first index of $\xi$.
		\item  $\angLap := \ginversesphere \cdot \angD^2$ denotes the covariant Laplacian 
			corresponding to $\gsphere$.
		\end{itemize}
\end{definition}

\begin{definition}[\textbf{Geometric torus differential}]
	\label{D:ANGULARDIFFERENTIAL}
	If $f$ is a scalar-valued function on $\ell_{t,u}$, 
	then $\angdiff f := \angD f = \Lineproject \D f$,
	where $\D f$ is the gradient one-form associated to $f$.
\end{definition}

Def.~\ref{D:ANGULARDIFFERENTIAL} allows us to avoid potentially confusing notation
such as $\angD \Lunit^i$ by instead writing
$\angdiff \Lunit^i$; the latter notation 
signifies to view $\Lunit^i$ as a scalar function under differentiation.

\begin{definition}[\textbf{Second fundamental forms}]
\label{D:SECONDFUND}
We define the second fundamental form $k$ of $\Sigma_t$, 
by
\begin{align} \label{E:SECONDFUNDSIGMATDEF}
	k 
	&:= \frac{1}{2} \SigmatLie_{\Transport} \gt.
\end{align}

	We define the null second fundamental form $\upchi$ of $\ell_{t,u}$ 
	by
\begin{align} \label{E:CHIDEF}
	\upchi
	& := \frac{1}{2} \angLie_{\Lunit} \gsphere.
\end{align}

\end{definition}

As was shown in \cite{jSgHjLwW2016}*{Subsection 2.6},
we have the following alternate expressions:
\begin{align} \label{E:ALTERNATESECONDFUND}
	k 
	&= \frac{1}{2} \SigmatLie_{\Transport} g,
		\qquad
	\upchi
	= \frac{1}{2} \angLie_{\Lunit} g.
\end{align}

\begin{lemma}\cite{jSgHjLwW2016}*{Lemma 2.3; \textbf{Alternate expressions for the second fundamental forms}}
We have the following identities:
\begin{align} \label{E:SECONDFUNDUSEFULID}
	\upchi_{\CoordAng \CoordAng}
	& = g(\D_{\CoordAng} \Lunit, \CoordAng),
		\qquad
	\angkdoublearg{\Radunit}{\CoordAng}
	= g(\D_{\CoordAng} \Lunit, \Radunit).
\end{align}
\end{lemma}

\begin{lemma}\cite{jSgHjLwW2016}*{Lemma 2.13; \textbf{Decompositions of some $\ell_{t,u}$ tensorfields into}
$\upmu^{-1}$-\textbf{singular and} $\upmu^{-1}$-\textbf{regular pieces}}
	\label{L:CONNECTIONLRADFRAME}
	Let $\upzeta$ be the $\ell_{t,u}$-tangent one-form defined by (see \eqref{E:SECONDFUNDUSEFULID})
	\begin{align} \label{E:ZETADEF}
		\upzeta_{\CoordAng} 
		& := \angkdoublearg{\Radunit}{\CoordAng}
		= g(\D_{\CoordAng} \Lunit, \Radunit) 
		= \upmu^{-1} g(\D_{\CoordAng} \Lunit, \Rad).
	\end{align}

	Then we can decompose the frame components of the 
	$\ell_{t,u}$-tangent
	tensorfields $\angk$ and $\upzeta$ into 
	$\upmu^{-1}$-singular and $\upmu^{-1}$-regular
	pieces as follows: 
	\begin{subequations}
	\begin{align}
	\upzeta & 
		= \upmu^{-1} \upzeta^{(Trans-\threePsi)}
		+ \upzeta^{(Tan-\threePsi)},
		 \label{E:ZETADECOMPOSED} \\
	\angk 
	& = \upmu^{-1} \angk^{(Trans-\threePsi)}
		+ \angk^{(Tan-\threePsi)},
		\label{E:ANGKDECOMPOSED}
	\end{align}
\end{subequations}
where
\begin{subequations}
	\begin{align}
	\upzeta^{(Trans-\threePsi)} 
		& :=
			- \frac{1}{2} \angGarg{\Lunit}\contr\Rad \threePsi,
			\label{E:ZETATRANSVERSAL} \\
		\angk^{(Trans-\threePsi)} 
		& := \frac{1}{2} \angG\contr\Rad \threePsi,
			\label{E:KABTRANSVERSAL}
				\\
		\upzeta^{(Tan-\threePsi)}
	& := \frac{1}{2} \angGarg{\Radunit}\contr\Lunit \threePsi
			- \frac{1}{2} \vec{G}_{\Lunit \Radunit}\contr\angdiff \threePsi
			- \frac{1}{2} \vec{G}_{\Radunit \Radunit}\contr\angdiff \threePsi,
		\label{E:ZETAGOOD} \\
	\angk^{(Tan-\threePsi)} 
	& := \frac{1}{2} \angG\contr\Lunit \threePsi
			- \frac{1}{2} \angGarg{\Lunit} \overset{\contr}{\otimes} \angdiff \threePsi
			- \frac{1}{2} \angdiff \threePsi \overset{\contr}{\otimes} \angGarg{\Lunit} 
			- \frac{1}{2} \angGarg{\Radunit} \overset{\contr}{\otimes} \angdiff \threePsi
			- \frac{1}{2} \angdiff \threePsi \overset{\contr}{\otimes} \angGarg{\Radunit}.
			\label{E:KABGOOD}
	\end{align}
\end{subequations}
\end{lemma}

\subsection{Pointwise norms}
\label{SS:POINTWISENORMS}
We always measure the magnitude of $\ell_{t,u}$ tensors using
the Riemannian metric $\gsphere$, as is captured by the following definition.

\begin{definition}[\textbf{Pointwise norms}]
	\label{D:POINTWISENORM}
	If $\xi_{\nu_1 \cdots \nu_n}^{\mu_1 \cdots \mu_m}$ 
	is a type $\binom{m}{n}$ $\ell_{t,u}$ tensor,
	then we define the norm $|\xi| \geq 0$ by
	\begin{align} \label{E:POINTWISENORM}
		|\xi|^2
		:= 
		\gsphere_{\mu_1 \widetilde{\mu}_1} \cdots \gsphere_{\mu_m \widetilde{\mu}_m}
		(\ginversesphere)^{\nu_1 \widetilde{\nu}_1} \cdots (\ginversesphere)^{\nu_n \widetilde{\nu}_n}
		\xi_{\nu_1 \cdots \nu_n}^{\mu_1 \cdots \mu_m}
		\xi_{\widetilde{\nu}_1 \cdots \widetilde{\nu}_n}^{\widetilde{\mu}_1 \cdots \widetilde{\mu}_m}.
	\end{align}
	In \eqref{E:POINTWISENORM}, 
	$\gsphere$ is the Riemannian metric on $\ell_{t,u}$ induced by 
	$g$, as given by Def.~\ref{D:FIRSTFUND}.
\end{definition}

\subsection{Expressions for the metrics}
\label{SS:METRICEXPRESSIONS}

\begin{lemma}\cite{jSgHjLwW2016}*{Lemma 2.4; \textbf{Expressions for $g$ and $g^{-1}$ in terms of the non-rescaled frame}}
\label{L:METRICDECOMPOSEDRELATIVETOTHEUNITFRAME}
We have the following identities:
\begin{subequations}
\begin{align}
	g_{\mu \nu} 
	& = - \Lunit_{\mu} \Lunit_{\nu}
			- (
					\Lunit_{\mu} \Radunit_{\nu} 
					+ 
					\Radunit_{\mu} \Lunit_{\nu}
				)
			+ \gsphere_{\mu \nu} 
			\label{E:METRICFRAMEDECOMPLUNITRADUNITFRAME},
			\\
	(g^{-1})^{\mu \nu} 
	& = 
			- \Lunit^{\mu} \Lunit^{\nu}
			- (
					\Lunit^{\mu} \Radunit^{\nu} 
					+ \Radunit^{\mu} \Lunit^{\nu}
				)
			+ (\ginversesphere)^{\mu \nu}.
			\label{E:GINVERSEFRAMEWITHRECTCOORDINATESFORGSPHEREINVERSE}
\end{align}
\end{subequations}

\end{lemma}

The following scalar-valued function captures the $\ell_{t,u}$ part of $g$.

\begin{definition}[\textbf{The metric component} $\gtancomp$]
\label{D:METRICANGULARCOMPONENT}
We define the function $\gtancomp > 0$ by
\begin{align} \label{E:METRICANGULARCOMPONENT}
	\gtancomp^2
	& :=  g(\CoordAng,\CoordAng) = \gsphere(\CoordAng,\CoordAng).
\end{align}
\end{definition}

It follows that relative to the geometric coordinates, we have
$
\ginversesphere
	= \gtancomp^{-2} \CoordAng \otimes \CoordAng
$.



\begin{lemma}\cite{jSgHjLwW2016}*{Corollary 2.6; 
\textbf{The geometric volume form factors of} $g$ \textbf{and} $\gt$}
\label{L:SPACETIMEVOLUMEFORMWITHUPMU}
The following identity is verified by the acoustcial metric $g$:
\begin{align} \label{E:SPACETIMEVOLUMEFORMWITHUPMU}
	|\mbox{\upshape{det}} g| 
	& = \upmu^2 \gtancomp^2,
\end{align}
where the determinant on the LHS is taken
relative to the geometric coordinates
$(t,u,\vartheta)$.

Furthermore, the following identity is verified by 
the first fundamental form $\gt$ of $\Sigma_t^{U_0}$: 
\begin{align} \label{E:SIGMATVOLUMEFORMWITHUPMU}
	\mbox{\upshape{det}} \gt|_{\Sigma_t^{U_0}}
	& = \upmu^2 \gtancomp^2,
\end{align}
where the determinant on the LHS is taken
relative to the geometric coordinates
$(u,\vartheta)$ induced on $\Sigma_t^{U_0}$.
\end{lemma}

\subsection{Commutation vectorfields}
\label{SS:COMMUTATIONVECTORFIELDS}
To derive estimates for the solution's higher-order derivatives,
we commute the equations with the elements of
$\lbrace \Lunit, \Rad, \GeoAng \rbrace$, where
$\GeoAng$ is the $\ell_{t,u}$-tangent vectorfield 
given in the next definition. We use $\GeoAng$ rather than
$\CoordAng$ because commuting $\CoordAng$ through
$\square_g$ seems to produce error terms 
that are uncontrollable because they lose a derivative.

\begin{definition}[\textbf{The vectorfields} $\GeoAng_{(Flat)}$ \textbf{and} $\GeoAng$]
\label{D:ANGULARVECTORFIELDS}
We define the Cartesian components of
the $\Sigma_t$-tangent vectorfields 
$\GeoAng_{(Flat)}$ and $\GeoAng$ 
as follows ($i=1,2$):
\begin{align}
	\GeoAng_{(Flat)}^i
	&: = \delta_2^i,
		\label{E:GEOANGEUCLIDEAN} \\
	\GeoAng^i
	& :=
	\Lineproject_a^{\ i} \GeoAng_{(Flat)}^a
		= \Lineproject_2^{\ i},
		\label{E:GEOANGDEF}
\end{align}
where $\Lineproject$ is the
$\ell_{t,u}$ projection tensorfield 
defined in \eqref{E:LINEPROJECTION}.
\end{definition}

When commuting the equations, we
use elements of the commutation sets 
$\Fullset$ and $\Tanset$.

\begin{definition}[\textbf{Commutation vectorfields}]
	\label{D:COMMUTATIONVECTORFIELDS}
	We define the commutation set $\Fullset$ as follows:
	\begin{align} \label{E:COMMUTATIONVECTORFIELDS}
		\Fullset
		:= \lbrace \Lunit, \Rad, \GeoAng \rbrace,
	\end{align}
	where $\Lunit$, $\Rad$, and $\GeoAng$ are respectively defined by
	\eqref{E:LUNITDEF}, \eqref{E:RADDEF}, and \eqref{E:GEOANGDEF}.

	We define the $\mathcal{P}_u$-tangent commutation set $\Tanset$ as follows:
	\begin{align} \label{E:TANGENTIALCOMMUTATIONVECTORFIELDS}
		\Tanset
		:= \lbrace \Lunit, \GeoAng \rbrace.
	\end{align}
\end{definition}

The Cartesian spatial components of
$\Lunit$,
$\Radunit$,
and $\GeoAng$
deviate from their
flat values by a small amount
that we denote by
$\Lunit_{(Small)}^i$,
$\Radunit_{(Small)}^i$,
and $\GeoAng_{(Small)}^i$.

\begin{definition}[\textbf{Perturbed part of various vectorfields}]
\label{D:PERTURBEDPART}
For $i=1,2$, we define the following scalar-valued functions:
\begin{align} \label{E:PERTURBEDPART}
	\Lunit_{(Small)}^i
	& := \Lunit^i
		- \delta_1^i,
	\qquad
	\Radunit_{(Small)}^i
	:= \Radunit^i
		+ \delta_1^i,
	\qquad
	\GeoAng_{(Small)}^i
	:= \GeoAng^i - \delta_2^i.
\end{align}
\end{definition}

\begin{lemma}[\textbf{Identity connecting} $\Lunit_{(Small)}^i$, $\Radunit_{(Small)}^i$, \textbf{and} $v^i$]
	\label{L:RADUNITSMALLLUNITSMALLRELATION}
	The following identity holds:
	\begin{align} \label{E:RADUNITSMALLLUNITSMALLRELATION}
		\Radunit_{(Small)}^i 
		& = -\Lunit_{(Small)}^i + v^i.
	\end{align}
\end{lemma}
\begin{proof}
	The identity \eqref{E:RADUNITSMALLLUNITSMALLRELATION}
	follows from
	\eqref{E:LITTLEGDECOMPOSED},
	\eqref{E:METRICPERTURBATIONFUNCTION},
	\eqref{E:DOWNSTAIRSUPSTAIRSSRADUNITPLUSLUNITISAFUNCTIONOFPSI}, 
	\eqref{E:PERTURBEDPART}, 
	and the fact that $(g^{-1})^{0i} = - v^i$. 
\end{proof}

In the next lemma, we characterize the
discrepancy between $\GeoAng_{(Flat)}$ and $\GeoAng$.

\begin{lemma}\cite{jSgHjLwW2016}*{Lemma 2.8; \textbf{Decomposition of} $\GeoAng_{(Flat)}$}
\label{L:GEOANGDECOMPOSITION}
We can decompose $\GeoAng_{(Flat)}$ into an $\ell_{t,u}$-tangent vectorfield
and a vectorfield parallel to $\Radunit$ as follows:
since $\GeoAng$ is $\ell_{t,u}$-tangent, 
there exists a scalar-valued function $\GeoAngFlatRadComponent$ such that
\begin{subequations}
\begin{align} \label{E:GEOANGINTERMSOFEUCLIDEANANGANDRADUNIT}
	\GeoAng_{(Flat)}^i
	& = \GeoAng^i 
		+ \GeoAngFlatRadComponent \Radunit^i,
			\\
	\GeoAng_{(Small)}^i
	& = - \GeoAngFlatRadComponent \Radunit^i.
	\label{E:GEOANGSMALLINTERMSOFRADUNIT}
\end{align}
\end{subequations}
Moreover, we have
\begin{align} \label{E:FLATYDERIVATIVERADIALCOMPONENT}
	\GeoAngFlatRadComponent 
	= g(\GeoAng_{(Flat)},\Radunit)
	= g_{ab} \GeoAng_{(Flat)}^a \Radunit^b
	= g_{2a} \Radunit^a
	= - \Speed^{-2} \Radunit_{(Small)}^2.
\end{align}
\end{lemma}

\subsection{Deformation tensors and basic vectorfield commutator properties}
\label{SS:BASICVECTORFIELDCOMMUTATOR}
In this section, we recall the standard definition of the deformation
of a vectorfield $V$. We then provide some simple commutator lemmas.

\begin{definition}[\textbf{Deformation tensor of a vectorfield} $V$]
\label{D:DEFORMATIONTENSOR}
If $V$ is a spacetime vectorfield,
then its deformation tensor $\deform{V}$
(relative to $g$)
is the symmetric type $\binom{0}{2}$ tensorfield
\begin{align} \label{E:DEFORMATIONTENSOR}
	\deformarg{V}{\alpha}{\beta}
	:= \Lie_V g_{\alpha \beta}
	= \D_{\alpha} V_{\beta} 
		+
		\D_{\beta} V_{\alpha},
\end{align}
where the second equality follows from
the torsion-free property of $\D$.

\end{definition}

\begin{lemma}[\textbf{Basic vectorfield commutator properties}]
\label{L:CONNECTIONBETWEENCOMMUTATORSANDDEFORMATIONTENSORS}
The vectorfields 
$[\Lunit, \Rad]$,
$[\Lunit, \GeoAng]$,
and
$[\Rad, \GeoAng] $
are $\ell_{t,u}$-tangent,
and the following identities hold:
\begin{align}
	[\Lunit, \Rad] 
	& = \angdeformoneformupsharparg{\Rad}{\Lunit},
		\qquad
	[\Lunit, \GeoAng] 
	= \angdeformoneformupsharparg{\GeoAng}{\Lunit},
		\qquad
	[\Rad, \GeoAng] 
	= \angdeformoneformupsharparg{\GeoAng}{\Rad}.
		\label{E:CONNECTIONBETWEENCOMMUTATORSANDDEFORMATIONTENSORS}
\end{align}

In addition, we have
\begin{subequations}
\begin{align}
	[\upmu \Transport, \Lunit]
	& = - (\Lunit \upmu) \Lunit
		+ \angdeformoneformupsharparg{\Lunit}{\Rad},
			\label{E:COMMUTATORFORMULAFORTRANSPORTRENORMANDLUNIT} \\
	[\upmu \Transport, \GeoAng]
	& = - (\GeoAng \upmu) \Lunit
		+ \upmu \angdeformoneformupsharparg{\GeoAng}{\Lunit}
		+ \angdeformoneformupsharparg{\GeoAng}{\Rad}.
		\label{E:COMMUTATORFORMULAFORTRANSPORTRENORMANDGEOANG}
\end{align}
\end{subequations}

Furthermore, if $Z \in \Fullset$, then
\begin{align}  \label{E:CONNECTIONBETWEENANGLIEOFGSPHEREANDDEFORMATIONTENSORS}
		\angLie_Z \gsphere 
		& = \angdeform{Z},
			\qquad
		\angLie_Z \ginversesphere
		= - \angdeform{Z}^{\# \#}.
\end{align}

Finally, if $V$ is an $\ell_{t,u}$-tangent vectorfield, then
\begin{align} \label{E:LIELANDLIERADOFELLTUTANGENTISELLTUTANGENT}
	[\Lunit, V] \mbox{ and } [\Rad, V] \mbox{ are } \ell_{t,u}-\mbox{tangent}.
\end{align}
\end{lemma}

\begin{proof}
		All aspects of the proposition except for
		\eqref{E:COMMUTATORFORMULAFORTRANSPORTRENORMANDLUNIT}-\eqref{E:COMMUTATORFORMULAFORTRANSPORTRENORMANDGEOANG}
		were derived in \cite{jSgHjLwW2016}*{Lemma 2.9}.
		The identities
		\eqref{E:COMMUTATORFORMULAFORTRANSPORTRENORMANDLUNIT}-\eqref{E:COMMUTATORFORMULAFORTRANSPORTRENORMANDGEOANG}
		are straightforward consequences of 
		the decomposition \eqref{E:TRANSPORTVECTORFIELDINTERMSOFLUNITANDRADUNIT}
		and the identities in \eqref{E:CONNECTIONBETWEENCOMMUTATORSANDDEFORMATIONTENSORS}.
\end{proof}

\begin{lemma}\cite{jSgHjLwW2016}*{Lemma 2.10; 
$\Lunit$, $\Rad$, $\GeoAng$ \textbf{commute with} $\angdiff$}
\label{L:LANDRADCOMMUTEWITHANGDIFF}
If $V \in \lbrace \Lunit, \Rad, \GeoAng \rbrace$
and $f$ is a scalar-valued function,
then
\begin{align} \label{E:ANGLIECOMMUTESWITHANGDIFF}
	\angLie_V \angdiff f
	& = \angdiff V f.
\end{align}
\end{lemma}

\subsection{Transport equations for the eikonal function quantities}
\label{SS:TRANSPORT}
We now provide transport equations verified by 
the scalar-valued functions
$\upmu$
and
$\Lunit_{(Small)}^i$.
These are the main equations we use to estimate
the eikonal function quantities below-top-order.
For top-order estimates, 
we use the modified quantities of Sect.~\ref{S:MODIFIED}.

\begin{lemma}\cite{jSgHjLwW2016}*{Lemma 2.12; \textbf{The transport equations verified by} $\upmu$ \textbf{and} $\Lunit^i$}
 \label{L:UPMUANDLUNITIFIRSTTRANSPORT}
The inverse foliation density $\upmu$ defined in \eqref{E:FIRSTUPMU} 
verifies the following transport equation:
\begin{align} \label{E:UPMUFIRSTTRANSPORT}
	\Lunit \upmu 
	& =
		\frac{1}{2} \vec{G}_{\Lunit \Lunit}\contr \Rad \threePsi
		- \frac{1}{2} \upmu \vec{G}_{\Lunit \Lunit}\contr\Lunit \threePsi
		- \upmu \vec{G}_{\Lunit \Radunit}\contr\Lunit \threePsi.
\end{align}
The scalar-valued Cartesian component functions $\Lunit_{(Small)}^i$,
($i=1,2$), defined in \eqref{E:PERTURBEDPART},
verify the following transport equation:
\begin{align} \label{E:LLUNITI} 
	\Lunit \Lunit_{(Small)}^i
	&  = - \frac{1}{2} \vec{G}_{\Lunit \Lunit}\contr(\Lunit \threePsi) \Lunit^i
			+ \frac{1}{2} \vec{G}_{\Lunit \Lunit}\contr(\Lunit \threePsi) v^i
				\\
	& \ \
			- \angGmixedarg{\Lunit}{\#}\contr(\Lunit \threePsi) \cdot (\angdiff x^i)
			+ \frac{1}{2} \vec{G}_{\Lunit \Lunit}\contr(\angdiffuparg{\#} \threePsi) \cdot \angdiff x^i.
			\notag
\end{align}
\end{lemma}

\subsection{Calculations connected to the failure of the null condition}
\label{SS:CALCSFORFAILUREOFNULLCONDITION}
Many of our most important estimates are tied to
the coefficients $\vec{G}_{\Lunit \Lunit}$.
In the next lemma, we derive expressions for them.
Then, in the subsequent lemma, we derive an expression for
the product
$
\frac{1}{2} \vec{G}_{\Lunit \Lunit} \contr \Rad \threePsi
$.
This presence of this product is tied to the failure of
Klainerman's null condition \cite{sK1984} and thus one expects
that the product must be non-zero for shocks to form;
this is explained in
the survey article \cite{gHsKjSwW2016} in 
a slightly different context.

\begin{lemma}[\textbf{Formula for} $G_{\Lunit \Lunit}^{\imath}$]
\label{L:FORMULAFORGLL}
Let $G_{\alpha \beta}^{\iota}$ be as in Def.~\ref{D:BIGGANDBIGH}.
Then for $i=1,2$, we have
\begin{subequations}
\begin{align}
	G_{\Lunit \Lunit}^0
	& = - 2 \Speed^{-1} \Speed',
		\label{E:GLUNITLUNIT0} \\
	G_{\Lunit \Lunit}^i
	& = 2 \Speed^{-2} (v^i - \Lunit^i)
		= 2 \Speed^{-2} \Radunit^i.
		\label{E:GLUNITLUNITi}
\end{align}
\end{subequations}
\end{lemma}

\begin{proof}
	We first prove \eqref{E:GLUNITLUNITi}.
	From the formula \eqref{E:ACOUSTICALMETRIC},
	Defs.~\ref{D:BIGGANDBIGH} and \ref{D:ARRAYNOTATION},
	the fact that $\Lunit^0 = 1$, 
	the identity \eqref{E:DOWNSTAIRSUPSTAIRSSRADUNITPLUSLUNITISAFUNCTIONOFPSI},
	and the fact that $\Lunit^i + \Radunit^i = v^i$
	(see \eqref{E:MATERIALVECTORVIELDRELATIVETORECTANGULAR} 
	and \eqref{E:TRANSPORTVECTORFIELDINTERMSOFLUNITANDRADUNIT}),
	we compute the desired identity as follows:
	\begin{align}
		G_{\Lunit \Lunit}^i
		= \left(
				\frac{\partial}{\partial v^i} g_{\alpha \beta} 
			\right)
			\Lunit^{\alpha} \Lunit^{\beta} \left(
				\frac{\partial}{\partial v^i} g_{\alpha \beta} 
			\right)
			\Lunit^{\alpha} \Lunit^{\beta}
		& = 2 \Speed^{-2} v^i (\Lunit^0)^2
			-2 \Speed^{-2} \Lunit^0 \Lunit^i
			 \\
		& = 2 \Speed^{-1} (v^i - \Lunit^i)
			= 2 \Speed^{-2} \Radunit^i.
			\notag 
	\end{align}

	We now prove \eqref{E:GLUNITLUNIT0}. Since 
	$g_{\alpha \beta}\Lunit^{\alpha}\Lunit^{\beta} = 0$, it suffices to 
	prove
	$
	 \left(
			\frac{\partial}{\partial \Densrenormalized} (\Speed^2 g_{\alpha \beta})
		\right)
			\Lunit^{\alpha} \Lunit^{\beta}
	=  - 2 \Speed \Speed'
	$.
	Since, among the components
	$\lbrace \Speed^2 g_{\alpha \beta} \rbrace_{\alpha,\beta = 0,1,2}$,
	only $\Speed^2 g_{00}$ depends on $\Densrenormalized$
	(see \eqref{E:ACOUSTICALMETRIC}),
	the desired identity is a simple consequence of the fact that
	$\Lunit^0 = 1$.

\end{proof}

\begin{lemma}[\textbf{Formula for} $\frac{1}{2} \vec{G}_{\Lunit \Lunit}\contr\Rad \threePsi$]
	\label{L:FORMULAFORGLLRADPSI}
	For solutions to the compressible Euler equations
	\eqref{E:TRANSPORTDENSRENORMALIZEDRELATIVETORECTANGULAR}-\eqref{E:TRANSPORTVELOCITYRELATIVETORECTANGULAR}
the following identity holds for the first product on RHS~\eqref{E:UPMUFIRSTTRANSPORT}:
\begin{align} \label{E:KEYLARGETERMEXPANDED}
	\frac{1}{2} \vec{G}_{\Lunit \Lunit}\contr\Rad \threePsi
	&
		=
		\Speed^{-2}
		\left\lbrace
			\Speed^{-1} \Speed' 
			+
			1
		\right\rbrace
		\delta_{ab} \Radunit^a \Rad v^b
		+ 
		\upmu \Speed^{-3} \Speed' \delta_{ab} \Lunit v^a \Radunit^b.
\end{align}
\end{lemma}

\begin{proof}
		Using Lemma~\ref{L:FORMULAFORGLL},
		we deduce
		$
		\frac{1}{2} \vec{G}_{\Lunit \Lunit}\contr\Rad \threePsi
		= - \Speed^{-1} \Speed' \Rad \Densrenormalized
		+ \Speed^{-2} \delta_{ab} \Radunit^a \Rad v^b
		$.
		Contracting \eqref{E:TRANSPORTVELOCITYRELATIVETORECTANGULAR}
		against $\delta_{ij} \Rad^j$ and inserting the resulting
		identity into the first product in the previous expression,
		we can rewrite it as
		$
		\Speed^{-3} \Speed' \delta_{ab} \Transport v^a \Rad^b
		+ \Speed^{-2} \delta_{ab} \Radunit^a \Rad v^b
		$. Using \eqref{E:TRANSPORTVECTORFIELDINTERMSOFLUNITANDRADUNIT}
		to substitute $\Lunit + \Radunit$ for $\Transport$ in the
		previous expression and recalling that $\Rad = \upmu \Radunit$, we 
		conclude \eqref{E:KEYLARGETERMEXPANDED}.
\end{proof}

Note that for the equation of state 
$p = C_0 - C_1 \exp(- \Densrenormalized)$
of a Chaplygin gas, 
we have
$
\Speed^{-1} \Speed' + 1 = 0
$.
For such a gas, the product 
$
\frac{1}{2} \vec{G}_{\Lunit \Lunit} \contr \Rad \threePsi
$
vanishes and our main shock formation results do not apply.
In fact, even in the plane symmetric case, it is not known whether shocks form in Chaplygin gas. In that case, only a very different type of singularity (where in particular the density itself blows up) is known to form \cite{dxKcWqZ14}. Moreover, in the case of the Chaplygin gas without vorticity, 
the wave equations 
\eqref{E:VELOCITYWAVEEQUATION}-\eqref{E:RENORMALIZEDDENSITYWAVEEQUATION}
verify Klainerman's null condition. While it is not directly related to the regime we study, we point out that in that case small-data global existence is known\footnote{Note that the equation for the irrotational Chaplygin gas is equivalent to that of a Minkowskian minimal surface equation, which is treated in \cite{hL2004}.} \cite{hL2004} when the data are given on the Cauchy hypersurface $\mathbb{R}^2$.

\subsection{Deformation tensor calculations}
\label{SS:DEFORMATIONTENSORCALCULATIONS}
In the next lemma, we provide explicit expressions
for the frame components of the deformation tensors of the commutation
vectorfields. 

\begin{lemma}\cite{jSgHjLwW2016}*{Lemma 2.18; \textbf{The frame components of} $\deform{Z}$}
\label{L:DEFORMATIONTENSORFRAMECOMPONENTS}
The following identities are verified by the deformation
tensors (see Def.~\ref{D:DEFORMATIONTENSOR}) of the elements of
$\Fullset$ (see \eqref{E:COMMUTATIONVECTORFIELDS}):
\begin{subequations}
\begin{align}
	\deformarg{\Rad}{\Lunit}{\Lunit} & = 0, 
		\qquad
	\deformarg{\Rad}{\Rad}{\Radunit} 
		= 2 \Rad \upmu,
		\qquad
	\deformarg{\Rad}{\Lunit}{\Rad} 
		= - \Rad \upmu,
			\label{E:RADDEFORMSCALARS}	\\
	\angdeformoneformarg{\Rad}{\Lunit}
	& = - \angdiff \upmu 
		- 2 \upzeta^{(Trans-\threePsi)}
		- 2 \upmu \upzeta^{(Tan-\threePsi)},
	\qquad
	\angdeformoneformarg{\Rad}{\Rad}
		= 0,
		\label{E:RADDEFORMSPHERERAD} \\
	\angdeform{\Rad}
	& = - 2 \upmu \mytr \upchi \gsphere
		+ 2 \angk^{(Trans-\threePsi)}
		+ 2 \upmu \angk^{(Tan-\threePsi)},
			\label{E:RADDEFORMSPHERE}
\end{align}
\end{subequations}

\begin{subequations}
\begin{align}
	\deformarg{\Lunit}{\Lunit}{\Lunit}
	& = 0, 
		\qquad
	\deformarg{\Lunit}{\Rad}{\Radunit} 
		= 2 \Lunit \upmu,
		\qquad
	\deformarg{\Lunit}{\Lunit}{\Rad} 
		= - \Lunit \upmu,
		\label{E:LUNITDEFORMSCALARS} \\
	\angdeformoneformarg{\Lunit}{\Lunit}
	& = 0, 
		\qquad
	\angdeformoneformarg{\Lunit}{\Rad}
		= \angdiff \upmu
	 		+ 2 \upzeta^{(Trans-\threePsi)}
			+ 2 \upmu \upzeta^{(Tan-\threePsi)}, 
	 		\label{E:LUNITDEFORMSPHERELUNITANDLUNITDEFORMSPHERERAD} \\
	\angdeform{\Lunit}
	& = 2 \mytr \upchi \gsphere,
	\label{E:LUNITDEFORMSPHERE}
\end{align}
\end{subequations}

\begin{subequations}
\begin{align}
	\deformarg{\GeoAng}{\Lunit}{\Lunit}
	& = 0, 
		\qquad
	\deformarg{\GeoAng}{\Rad}{\Radunit} 
		= 2 \GeoAng \upmu,
		\qquad
	\deformarg{\GeoAng}{\Lunit}{\Rad} 
		= - \GeoAng \upmu,
		\label{E:GEOANGDEFORMSCALARS} \\
	\angdeformoneformarg{\GeoAng}{\Lunit}
	& = - \mytr \upchi \GeoAng_{\flat}
		+ \frac{1}{2} (\angG \cdot \GeoAng)\contr\Lunit \threePsi
		+ \GeoAngFlatRadComponent \angGarg{\Radunit}\contr\Lunit \threePsi
			\label{E:GEOANGDEFORMSPHEREL} \\
	& \ \
		+ \frac{1}{2} (\angGarg{\Lunit} \cdot \GeoAng)\contr\angdiff \threePsi
		- \GeoAngFlatRadComponent \vec{G}_{\Lunit \Radunit}\contr\angdiff \threePsi
		- \frac{1}{2} \GeoAngFlatRadComponent \vec{G}_{\Radunit \Radunit}\contr\angdiff \threePsi, 
		\notag \\
	\angdeformoneformarg{\GeoAng}{\Rad}
	& = \upmu \mytr \upchi \GeoAng_{\flat}
			+ \GeoAngFlatRadComponent \angdiff \upmu
			+ \GeoAngFlatRadComponent \angGarg{\Radunit}\contr\Rad \threePsi
			- \frac{1}{2} \upmu \GeoAngFlatRadComponent \vec{G}_{\Radunit \Radunit}\contr\angdiff \threePsi
			\label{E:GEOANGDEFORMSPHERERAD} \\
	& \ \
			- \frac{1}{2} \upmu (\angG \cdot \GeoAng)\contr\Lunit \threePsi
			+ \upmu (\angGarg{\Lunit} \cdot \GeoAng)\contr\angdiff \threePsi
			+ \upmu (\angGarg{\Radunit} \cdot \GeoAng)\contr\angdiff \threePsi, 
	 		\notag \\
	\angdeform{\GeoAng}
	& = 2 \GeoAngFlatRadComponent \mytr \upchi \gsphere
		+ \frac{1}{2} (\angG \cdot \GeoAng) \overset{\contr}{\otimes} \angdiff \threePsi
		+ \frac{1}{2} \angdiff \threePsi \overset{\contr}{\otimes} (\angG \cdot \GeoAng) 
		- \GeoAngFlatRadComponent \angG\contr\Lunit \threePsi
			\label{E:GEOANGDEFORMSPHERE} \\
	& \ \
		+ \GeoAngFlatRadComponent \angGarg{\Lunit} \overset{\contr}{\otimes} \angdiff \threePsi
		+ \GeoAngFlatRadComponent \angdiff \threePsi \overset{\contr}{\otimes} \angGarg{\Lunit} 
		+ \GeoAngFlatRadComponent \angGarg{\Radunit} \overset{\contr}{\otimes} \angdiff \threePsi
		+ \GeoAngFlatRadComponent \angdiff \threePsi \overset{\contr}{\otimes} \angGarg{\Radunit}.
		\notag 
\end{align}
\end{subequations}
	The scalar-valued function $\GeoAngFlatRadComponent$ from above
	is as in Lemma~\ref{L:GEOANGDECOMPOSITION}, while
	the $\ell_{t,u}$-tangent tensorfields
	$\upchi$,
	$\upzeta^{(Trans-\threePsi)}$,
	$\angk^{(Trans-\threePsi)}$,
	$\upzeta^{(Tan-\threePsi)}$,
	and
	$\angk^{(Tan-\threePsi)}$
from above 
are as in
\eqref{E:CHIDEF},
\eqref{E:ZETATRANSVERSAL},
\eqref{E:KABTRANSVERSAL},
\eqref{E:ZETAGOOD},
and \eqref{E:ZETAGOOD}.
In \eqref{E:GEOANGDEFORMSPHERE},
$
\displaystyle
\angGarg{\Lunit} \overset{\contr}{\otimes} \angdiff \threePsi
:= 
\sum_{\imath=0}^2 
\NovecangGarg{\Lunit}{\imath} \otimes \angdiff \threePsi_{\imath} 
$,
and similarly for the other terms involving $\overset{\contr}{\otimes}$.
\end{lemma}

\subsection{Useful expressions for the null second fundamental form}
The next lemma provides explicit formulas for 
$\upchi$,
$\mytr \upchi$,
and $\Lunit \ln \gtancomp$.

\begin{lemma}\cite{jSgHjLwW2016}*{Lemma 2.15; \textbf{Identities involving} $\upchi$}
	\label{L:IDENTITIESINVOLVING}
	Let $\upchi$ be the $\ell_{t,u}$ tensorfield
	defined in \eqref{E:CHIDEF} and let
	$\gtancomp$ be the metric component
	from Def.~\ref{D:METRICANGULARCOMPONENT}.
	We have the following identities:
	\begin{subequations}
	\begin{align} \label{E:CHIINTERMSOFOTHERVARIABLES}
	\upchi 
	& = g_{ab} (\angdiff \Lunit^a) \otimes \angdiff x^b
		+ \frac{1}{2} \angG\contr\Lunit\threePsi,
			\\
	\mytr \upchi 
	& = g_{ab} \ginversesphere \cdot \left\lbrace (\angdiff \Lunit^a) \otimes \angdiff x^b \right\rbrace
		+ \frac{1}{2} \ginversesphere \cdot \angG\contr\Lunit\threePsi,
			 \label{E:TRCHIINTERMSOFOTHERVARIABLES}
			\\
	\Lunit \ln \gtancomp
	& = \mytr \upchi.
	\label{E:LDERIVATIVEOFVOLUMEFORMFACTOR}
\end{align}
\end{subequations}
\end{lemma}

\subsection{Decomposition of differential operators}
\label{SS:DECOMPOFDIFFERENTIALOPERATORS}
We start by decomposing $\upmu \square_{g(\threePsi)}$ relative to the rescaled frame.
The factor of $\upmu$ is important for our decompositions.

\begin{proposition} [\textbf{Frame decomposition of $\upmu \square_{g(\threePsi)} f$}]
	\label{P:GEOMETRICWAVEOPERATORFRAMEDECOMPOSED}
	Let $f$ be a scalar-valued function.
	Then $\upmu \square_{g(\threePsi)} f$ can be expressed in either of the following two forms:
	\begin{subequations}
	\begin{align} \label{E:LONOUTSIDEGEOMETRICWAVEOPERATORFRAMEDECOMPOSED}
		\upmu \square_{g(\threePsi)} f 
			& = - \Lunit(\upmu \Lunit f + 2 \Rad f)
				+ \upmu \angLap f
				- \mytr \upchi \Rad f
				- \upmu \mytr \angk \Lunit f
				- 2 \upmu \upzeta^{\#} \cdot \angdiff f,
					\\
			& = - (\upmu \Lunit + 2 \Rad) (\Lunit f) 
				+ \upmu \angLap f
				- \mytr \upchi \Rad f
				- (\Lunit \upmu) \Lunit f
				+ 2 \upmu \upzeta^{\#} \cdot \angdiff f
				+ 2 (\angdiffuparg{\#} \upmu) \cdot \angdiff f,
				\label{E:LONINSIDEGEOMETRICWAVEOPERATORFRAMEDECOMPOSED}
	\end{align}
	\end{subequations}
	where the $\ell_{t,u}$-tangent tensorfields
	$\upchi$,
	$\upzeta$,
	and
	$\angk$
	can be expressed via
	\eqref{E:CHIINTERMSOFOTHERVARIABLES},
	\eqref{E:ZETADECOMPOSED},
	and
	\eqref{E:ANGKDECOMPOSED}.
\end{proposition}

\begin{lemma}[\textbf{Expression for} $\partial_{\nu}$ \textbf{in terms of geometric vectorfields}]
	\label{L:CARTESIANVECTORFIELDSINTERMSOFGEOMETRICONES}
	We can express the Cartesian coordinate partial derivative 
	vectorfields in terms of 
	$\Lunit$, $\Radunit$, and $\GeoAng$ as follows,
	$(i=1,2)$:
	\begin{subequations}
	\begin{align} 
		\partial_t
		& = 
			\Lunit
			-  
			(g_{\alpha 0} \Lunit^{\alpha}) \Radunit
			+
			\left(
				\frac{g_{a0} \GeoAng^a}{g_{cd} \GeoAng^c \GeoAng^d}
			\right)
			\GeoAng,
			\label{E:PARTIALTINTERMSOFLUNITRADUNITANDGEOANG} \\
		\partial_i
		& = (g_{ai} \Radunit^a) \Radunit 
			+ 
			\left(
				\frac{g_{ai} \GeoAng^a}{g_{cd} \GeoAng^c \GeoAng^d}
			\right)
			\GeoAng.
			\label{E:PARTIALIINTERMSOFRADUNITANDGEOANG}
	\end{align}
	\end{subequations}
\end{lemma}

\begin{proof}
	We expand 
	$\partial_i = \upalpha_i \Radunit + \upbeta_i \GeoAng$
	for scalar-valued functions $\upalpha_i$ and $\upbeta_i$.
	Taking the $g$-inner product of each side with respect to
	$\Radunit$, we obtain
	$\upalpha_i = g(\Radunit,\partial_i) = g_{ab} \Radunit^a \delta_i^b = g_{ai} \Radunit^a$.
	Similarly, we take the inner product with respect to $\GeoAng$ 
	to deduce $\upbeta_i g_{cd} \GeoAng^c \GeoAng^d = g_{ai} \GeoAng^a$.
	Using these identities to substitute for
	$\upalpha_i$ and $\upbeta_i$,
	we conclude \eqref{E:PARTIALIINTERMSOFRADUNITANDGEOANG}.
	A similar argument yields \eqref{E:PARTIALTINTERMSOFLUNITRADUNITANDGEOANG}, 
	though in this
	case we must use an expansion of the form
	$\partial_t = \upalpha \Lunit + \upbeta \Radunit + \upgamma \GeoAng$;
	we omit the details.
	\end{proof}

	With the help of Lemma~\ref{L:CARTESIANVECTORFIELDSINTERMSOFGEOMETRICONES},
	we can now express the products
	on RHS~\eqref{E:VELOCITYWAVEEQUATION}
	involving $\partial_a \Vortrenormalized$ in terms of 
	\textbf{$\mathcal{P}_u$-tangent}
	geometric derivatives of $\Vortrenormalized$.

	\begin{corollary}[\textbf{Decomposition of the vorticity derivatives in equation 
	\eqref{E:VELOCITYWAVEEQUATION}}]
		\label{C:VELOCITYWAVEEQUATIONDERIVATIVEOFVORTICITYINHOMOGENEOUSTERMEXPRESSION}
		We have the following identity for the vorticity derivative-involving product
		on RHS~\eqref{E:VELOCITYWAVEEQUATION}:
		\begin{align} \label{E:VELOCITYWAVEEQUATIONDERIVATIVEOFVORTICITYINHOMOGENEOUSTERMEXPRESSION}
			- [ia] (\exp \Densrenormalized) \Speed^2 (\upmu \partial_a \Vortrenormalized)
			& =  [ia] \upmu (\exp \Densrenormalized) \Speed^2 (g_{ab} \Radunit^b) \Lunit \Vortrenormalized
					 \\
			& \ \
			- [ia] \upmu (\exp \Densrenormalized) \Speed^2
			\left(
				\frac{g_{ab} \GeoAng^b}{g_{cd} \GeoAng^c \GeoAng^d}
			\right)
			\GeoAng \Vortrenormalized.
			\notag
		\end{align}
	\end{corollary}

	\begin{proof}
			We first use the formula 
			\eqref{E:PARTIALIINTERMSOFRADUNITANDGEOANG} to express
			the factor $\partial_a \Vortrenormalized$ on 
			LHS~\eqref{E:VELOCITYWAVEEQUATIONDERIVATIVEOFVORTICITYINHOMOGENEOUSTERMEXPRESSION}
			in terms of $\Radunit \Vortrenormalized$ and $\GeoAng \Vortrenormalized$.
			We then use \eqref{E:RENORMALIZEDVORTICTITYTRANSPORTEQUATION}
			and \eqref{E:TRANSPORTVECTORFIELDINTERMSOFLUNITANDRADUNIT}
			to replace $\Radunit \Vortrenormalized$ with $-\Lunit \Vortrenormalized$.
	\end{proof}

\subsection{Arrays of fundamental unknowns and schematic notation}
\label{SS:ARRAYS}
In Lemma~\ref{L:SCHEMATICDEPENDENCEOFMANYTENSORFIELDS}, we show that
many scalar-valued functions and tensorfields
that we have introduced depend on just a handful of more fundamental functions and tensorfields.
This simplifies various aspects of our analysis.
We start by introducing some convenient shorthand notation.

\begin{definition}[\textbf{Shorthand notation for the unknowns}]
\label{D:ABBREIVATEDVARIABLES}
We define the following arrays $\GdVar$ and $\BadVar$ of scalar-valued functions:
\begin{align}
	\GdVar 
	& := \left(\threePsi, \Lunit_{(Small)}^1, \Lunit_{(Small)}^2 \right),
			\qquad
	\BadVar 
		:= \left(\threePsi, \upmu, \Lunit_{(Small)}^1, \Lunit_{(Small)}^2 \right).
	\label{E:ABBREIVATEDVARIABLES}
\end{align}
\end{definition}

\begin{remark}[\textbf{Schematic functional dependence}]
\label{R:SCHEMATICTENSORFIELDPRODUCTS}
Throughout,
$\smoothfunction(\xi_{(1)},\xi_{(2)},\cdots,\xi_{(m)})$ 
schematically
denotes an expression (often tensorial and involving contractions)
that depends smoothly on the $\ell_{t,u}$-tangent tensorfields $\xi_{(1)}, \xi_{(2)}, \cdots, \xi_{(m)}$.
In general, we have $\smoothfunction(0) \neq 0$.
We sometimes use the notation $\vec{x} := (x^1,x^2)$
and $\angdiff \vec{x} := (\angdiff x^1,\angdiff x^2)$
in our schematic depictions.
\end{remark}

\begin{lemma}[\textbf{Schematic structure of various tensorfields}]
	\label{L:SCHEMATICDEPENDENCEOFMANYTENSORFIELDS}
	We have the following schematic relations for scalar-valued functions:
	\begin{subequations}
	\begin{align} \label{E:SCALARSDEPENDINGONGOODVARIABLES}
		g_{\alpha \beta},
		(g^{-1})^{\alpha \beta},
		(\gt^{-1})^{ab},
		\gsphere_{\alpha \beta},
		(\ginversesphere)^{\alpha \beta},
		G_{\alpha \beta}^{\imath},
		H_{\alpha \beta}^{\imath \jmath},
		\Lineproject_{\beta}^{\ \alpha},
		\Lunit^{\alpha}, 
		\Radunit^{\alpha},
		\GeoAng^{\alpha},
		\Speed
		& = \smoothfunction(\GdVar),
			\\
		G_{\Lunit \Lunit}^{\imath},
		G_{\Lunit \Radunit}^{\imath},
		G_{\Radunit \Radunit}^{\imath},
		H_{\Lunit \Lunit}^{\imath \jmath},
		H_{\Lunit \Radunit}^{\imath \jmath},
		H_{\Radunit \Radunit}^{\imath \jmath}
		& = \smoothfunction(\GdVar),
			\label{E:GFRAMESCALARSDEPENDINGONGOODVARIABLES} \\
		g_{\alpha \beta}^{(Small)},
		\Lineproject_{\beta}^{\ \alpha}
		- \delta_{\beta 2} \delta^{\alpha 2},
		\GeoAng_{(Small)}^{\alpha},
		\Radunit_{(Small)}^{\alpha},
		\GeoAngFlatRadComponent
		& = \smoothfunction(\GdVar) \GdVar,
			\label{E:LINEARLYSMALLSCALARSDEPENDINGONGOODVARIABLES}
			\\
		\Rad^{\alpha} 
		& = \smoothfunction(\BadVar).
			\label{E:SCALARSDEPENDINGONBADVARIABLES}
	\end{align}
	\end{subequations}

	Moreover, we have the following schematic relations for $\ell_{t,u}$-tangent tensorfields:
	\begin{subequations}
	\begin{align}
		\gsphere,
		\angGarg{\Lunit},
		\angGarg{\Radunit},
		\angG,
		\NovecangHarg{\Lunit}{\imath}{\jmath},
		\NovecangHarg{\Radunit}{\imath}{\jmath},
		\NovecangH{\imath}{\jmath}
		& = \smoothfunction(\GdVar,\angdiff \vec{x}),
			\label{E:TENSORSDEPENDINGONGOODVARIABLES} \\
	\GeoAng
	& = \smoothfunction(\GdVar,\ginversesphere,\angdiff \vec{x}),
			\label{E:TENSORSDEPENDINGONGOODVARIABLESANDGINVERSESPHERE}
			\\
	\upzeta^{(Tan-\threePsi)},
	\angk^{(Tan-\threePsi)} 
	& = \smoothfunction(\GdVar,\angdiff \vec{x}) \Singletan \threePsi,
		\label{E:TENSORSDEPENDINGONGOODVARIABLESGOODPSIDERIVATIVES}
			\\
	\angk^{(Trans-\threePsi)}
	& = \smoothfunction(\GdVar,\angdiff \vec{x}) \Rad \threePsi,
		\label{E:TENSORSDEPENDINGONGOODVARIABLESBADDERIVATIVES} \\
	\upzeta^{(Trans-\threePsi)}
	& = \smoothfunction(\Rad \threePsi,\angdiff \vec{x}) \GdVar,
		\label{E:ZETABADDERIVATIVESBETTERTHANEXPECTED} \\
	\upchi 
	& = \smoothfunction(\GdVar,\angdiff \vec{x}) \Singletan \GdVar,
		\label{E:TENSORSDEPENDINGONGOODVARIABLESGOODDERIVATIVES}
			\\
	\mytr \upchi 
	& = \smoothfunction(\GdVar,\ginversesphere,\angdiff \vec{x}) \Singletan \GdVar.
	\label{E:TENSORSDEPENDINGONGOODVARIABLESGOODDERIVATIVESANDGINVERSESPHERE}
  \end{align}
\end{subequations}

Finally, the null forms $\mathscr{Q}^i$ and $\mathscr{Q}$
defined by 
\eqref{E:VELOCITYNULLFORM} and \eqref{E:DENSITYNULLFORM},
upon being multiplied by $\upmu$,
have the following schematic structure:
\begin{align} \label{E:UPMUTIMESNULLFORMSSCHEMATIC}
	\upmu \mathscr{Q}^i,
		\,
	\upmu \mathscr{Q}
	& =
	\smoothfunction(\BadVar,\Rad \threePsi,\Singletan \threePsi) 
	\Singletan \threePsi.
\end{align}

\end{lemma}

\begin{proof}
	Except for \eqref{E:UPMUTIMESNULLFORMSSCHEMATIC}, the desired relations
	were proved as \cite{jSgHjLwW2016}*{Lemma 2.19}.
	We now prove \eqref{E:UPMUTIMESNULLFORMSSCHEMATIC}.
	The desired result for the term
	$
	\Speed^{-1} \Speed' 
	(g^{-1})^{\alpha \beta} \partial_{\alpha} \Densrenormalized \partial_{\beta} \Densrenormalized
	$
	on RHS~\eqref{E:DENSITYNULLFORM}
	is a simple consequence of the identity
	$
	\displaystyle
	g^{-1} 
	= - \Lunit \otimes \Lunit
		- \Lunit \otimes \Radunit
		- \Radunit \otimes \Lunit
		+ \frac{1}{g_{ab} \GeoAng^a \GeoAng^b} \GeoAng \otimes \GeoAng
	$,
	which is easy to verify by contracting each side against the 
	$g$-duals of the elements of 
	$\lbrace \Lunit, \Radunit, \GeoAng \rbrace$.
	We now consider the quadratic term 
	$
	\partial_1 v^1 \partial_2 v^2
	-
	\partial_2 v^1 \partial_1 v^2 
	$
	on RHS~\eqref{E:DENSITYNULLFORM}.
	Similar remarks apply to the quadratic term $-(g^{-1})^{\alpha \beta} \partial_{\alpha} \Densrenormalized \partial_{\beta} v^i$
	on RHS~~\eqref{E:VELOCITYNULLFORM}.
	We use \eqref{E:PARTIALIINTERMSOFRADUNITANDGEOANG}
	to write the Cartesian coordinate partial derivatives in the previous expression
	in terms of $\Radunit$ and $\GeoAng$ derivatives.
	In view of the antisymmetry of the expression in $v^1$ and $v^2$,
	we see that the terms proportional to 
	$(\Radunit v^1) \Radunit v^2$ cancel 
	(although it is not important
	for the main results of this paper, we note that
	the terms proportional to $(\GeoAng v^1) \GeoAng v^2$ also cancel).
	Multiplying by $\upmu$, we conclude that the quadratic term
	under consideration is of the form RHS~\eqref{E:UPMUTIMESNULLFORMSSCHEMATIC}.
	\end{proof}

\subsection{Geometric decompositions involving \texorpdfstring{$\GeoAng$}{the geometric torus commutation vectorfield}}
\label{SS:GEOMETRICDEOMPOSITIONSINVOLVINGGEOANG}
In this section, we express various $\ell_{t,u}$ tensorfields
and operators in terms of $\GeoAng$. This allows for a simplified approach
to deriving various formulas and estimates.

\begin{lemma}[\textbf{Formula for} $\ginversesphere$ \textbf{in terms of } $\GeoAng$]
	\label{L:INVERSEANGULARMETRICINTERMSOFGEOANG}
		Let $\ginversesphere$ be the inverse 
		of the first fundamental form $\gsphere$
		from Def.~\ref{D:FIRSTFUND}
		and let $\GeoAng$ be the $\ell_{t,u}$-tangent vectorfield from
		Def.~\ref{D:ANGULARVECTORFIELDS}. 
		We have the following identity:
	\begin{align} \label{E:INVERSEANGULARMETRICINTERMSOFGEOANG}
		\ginversesphere
		& = \frac{1}{g(\GeoAng,\GeoAng)} \GeoAng \otimes \GeoAng.
	\end{align}
\end{lemma}
\begin{proof}
	Since the $\ell_{t,u}$ are one-dimensional,
	$\ginversesphere$ must be a multiple of
	$\GeoAng \otimes \GeoAng$.
	Contracting \eqref{E:INVERSEANGULARMETRICINTERMSOFGEOANG}
	against $\GeoAng_{\flat} \otimes \GeoAng_{\flat}$,
	we easily obtain that the correct proportionality factor is 
	$
	\displaystyle
	\frac{1}{g(\GeoAng,\GeoAng)}
	$.
\end{proof}

\begin{lemma}[\textbf{$\xi$ in terms of $\mytr \xi$}]
	\label{L:XIINTERMSOFTRACEXI}
	Let $\GeoAng$ be the $\ell_{t,u}$-tangent vectorfield from
	Def.~\ref{D:ANGULARVECTORFIELDS}.
	We have the following identity,
	valid for symmetric type $\binom{0}{2}$ $\ell_{t,u}$ tensorfields
	$\xi$:
	\begin{align} \label{E:XIINTERMSOFTRACEXI}
		\xi 
		& = \frac{1}{g(\GeoAng,\GeoAng)} \mytr \xi \GeoAng_{\flat} \otimes \GeoAng_{\flat}.
	\end{align}
\end{lemma}

\begin{proof}
	Since $\ell_{t,u}$ is one-dimensional, we have
	$
	\displaystyle
	\xi = A \GeoAng_{\flat} \otimes \GeoAng_{\flat}
	$
	for some scalar-valued function $A$. Taking the 
	$\gsphere$-trace of this equation, we find that
	$
	\displaystyle
	A = \frac{1}{g(\GeoAng,\GeoAng)} \mytr \xi
	$
	as desired.
	\end{proof}

\section{Area and Volume Forms and Energy-Null Flux Identities}
\label{S:FORMSANDENERGY}
In this section, we first define geometric area and volume forms
and corresponding integrals. 
Using these, we construct the energies and null fluxes
that we use to control the solution and its derivatives in $L^2$.
We then exhibit the basic coercive properties of the energies 
and null fluxes and provide the fundamental energy identities
that we use to derive a priori estimates.
There are two identities: one for wave equations, which we use
to control $\threePsi$ (see Prop.~\ref{P:DIVTHMWITHCANCELLATIONS}),
and one for transport equations, which we use to
control $\Vortrenormalized$
(see Prop.~\ref{P:ENERGYIDENTITYRENORMALIZEDVORTICITY}).

\subsection{Area and volume forms and geometric integrals}
\label{SS:FORMS}
We define our geometric integrals in terms of length, area, and volume
forms that remain non-degenerate throughout the evolution,
all the way up to the shock.

\begin{definition}[\textbf{Non-degenerate forms and related integrals}]
	\label{D:NONDEGENERATEVOLUMEFORMS}
	We define the length form
	$d \spherevol$ on $\ell_{t,u}$,
	the area form $d \tvol$ on $\Sigma_t^u$,
	the area form $d \conevol$ on $\mathcal{P}_u^t$,
	and the volume form $d \vol$ on $\mathcal{M}_{t,u}$
	as follows (relative to the geometric coordinates):
	\begin{align} \label{E:RESCALEDVOLUMEFORMS}
			d \spherevol
			& = d \spherevol{(t,u,\vartheta)}
			:= \gtancomp(t,u,\vartheta) \, d \vartheta,
				&&
			d \tvol
			=
			d \tvol(t,u',\vartheta)
			:= d \spherevol(t,u',\vartheta) du',
				\\
			d \conevol 
			& = d \conevol(t',u,\vartheta)
			:= d \spherevol(t',u,\vartheta) dt',
				&&
			d \vol 
			= d \vol(t',u',\vartheta')
			:= d \spherevol(t',u',\vartheta') du' dt',
				\notag
	\end{align}
	where $\gtancomp$ is the scalar-valued function from Def.~\ref{D:METRICANGULARCOMPONENT}.

	If $f$ is a scalar-valued function, then we define
	\begin{subequations}
	\begin{align}
	\int_{\ell_{t,u}}
			f
		\, d \spherevol
		& 
		:=
		\int_{\vartheta \in \mathbb{T}}
			f(t,u,\vartheta)
		\, \gtancomp(t,u,\vartheta) d \vartheta,
			\label{E:LINEINTEGRALDEF} \\
		\int_{\Sigma_t^u}
			f
		\, d \tvol
		& 
		:=
		\int_{u'=0}^u
		\int_{\vartheta \in \mathbb{T}}
			f(t,u',\vartheta)
		\, \gtancomp(t,u',\vartheta) d \vartheta du',
			\label{E:SIGMATUINTEGRALDEF} \\
		\int_{\mathcal{P}_u^t}
			f
		\, d \conevol
		& 
		:=
		\int_{t'=0}^t
		\int_{\vartheta \in \mathbb{T}}
			f(t',u,\vartheta)
		\, \gtancomp(t',u,\vartheta) d \vartheta dt',
			\label{E:PUTINTEGRALDEF} \\
		\int_{\mathcal{M}_{t,u}}
			f
		\, d \vol
		& 
		:=
		\int_{t'=0}^t
		\int_{u'=0}^u
		\int_{\vartheta \in \mathbb{T}}
			f(t',u',\vartheta)
		\, \gtancomp(t',u',\vartheta) d \vartheta du' dt'.
		\label{E:MTUTUINTEGRALDEF}
	\end{align}
	\end{subequations}
\end{definition}

\begin{remark}
	The canonical forms associated to
	$\gt$ and $g$ are respectively
	$\upmu d \tvol$ and $\upmu d \vol$.
\end{remark}

\subsection{Basic ingredients and the definitions of the energies and null fluxes}
We construct our fundamental energies and null fluxes for scalar-valued functions $\Psi$
with the help of the \emph{energy-momentum tensor} 
\begin{align} \label{E:ENERGYMOMENTUMTENSOR}
	\enmomtensor_{\mu \nu}
	=
	\enmomtensor_{\mu \nu}[\Psi]
	& := \D_{\mu} \Psi \D_{\nu} \Psi
	- \frac{1}{2} g_{\mu \nu} (g^{-1})^{\alpha \beta} \D_{\alpha} \Psi \D_{\beta} \Psi.
\end{align}

We construct our energies and null fluxes for $\threePsi$
by contracting the following multiplier
vectorfield $\Mult$ against $\enmomtensor$.

\begin{definition}[\textbf{The timelike multiplier vectorfield} $\Mult$]
\label{D:DEFINITIONMULT} 
We define 
\begin{align} \label{E:DEFINITIONMULT} 
		\Mult 
		& := (1 + 2 \upmu) \Lunit + 2 \Rad.
\end{align}
\end{definition}
Note that $g(\Mult,\Mult) = - 4 \upmu (1 + \upmu) < 0$. 
This property leads to coercive energy identities.

\begin{definition}[\textbf{Energies and null fluxes}]
\label{D:ENERGYFLUX}
In terms of the non-degenerate forms of Def.~\ref{E:RESCALEDVOLUMEFORMS},
we define the energy functional $\enzero[\cdot]$ and null flux functional
$\flzero[\cdot]$ as follows:
\begin{align} \label{E:ENERGYFLUX}
	\enzero[\Psi](t,u)
		& := 
		\int_{\Sigma_t^u} 
			\upmu \enmomtensor_{\Transport \Mult}[\Psi]
		\, d \tvol,
	\qquad
	\flzero[\Psi](t,u)
		:=	 
		\int_{\mathcal{P}_u^t} 
			\enmomtensor_{\Lunit \Mult}[\Psi]
		\, d \conevol,
\end{align}
where $\Transport$ and $\Mult$
are the vectorfields defined in
\eqref{E:MATERIALVECTORVIELDRELATIVETORECTANGULAR}
and
\eqref{E:DEFINITIONMULT}.

We define the energy functional $\Vortenzero[\cdot]$ and null flux functional
$\Vortflzero[\cdot]$ as follows:
\begin{align} \label{E:RENORMALIZEDVORTENERGYFLUX}
	\Vortenzero[\Vortrenormalized](t,u)
		& := 
		\int_{\Sigma_t^u} 
			\upmu \Vortrenormalized^2
		\, d \tvol,
	\qquad
	\Vortflzero[\Vortrenormalized](t,u)
	 := 
		\int_{\mathcal{P}_u^t} 
			\Vortrenormalized^2
		\, d \conevol.
\end{align}

\end{definition}

\begin{lemma}\cite{jSgHjLwW2016}*{Lemma 3.4; \textbf{Coerciveness of the energy and null flux}}
\label{L:ORDERZEROCOERCIVENESS}
The energy 
$\enzero[\Psi]$
and null flux 
$\flzero[\Psi]$
from Def.~\ref{D:ENERGYFLUX} enjoy the following 
coerciveness properties:
\begin{subequations}
\begin{align} \label{E:ENERGYORDERZEROCOERCIVENESS}
	\enzero[\Psi](t,u)
		& = 
		\int_{\Sigma_t^u} 
			\frac{1}{2} (1 + 2 \upmu) \upmu (\Lunit \Psi)^2
			+ 2 \upmu (\Lunit \Psi) \Rad \Psi
			+ 2 (\Rad \Psi)^2
			+ \frac{1}{2} (1 + 2 \upmu)\upmu |\angdiff \Psi|^2
		\, d \tvol,
		\\
	\flzero[\Psi](t,u)
		& = 
		\int_{\mathcal{P}_u^t} 
			(1 + \upmu)(\Lunit \Psi)^2
			+ \upmu |\angdiff \Psi|^2
		\, d \conevol.
		\label{E:NULLFLUXENERGYORDERZEROCOERCIVENESS}
\end{align}
\end{subequations}

\end{lemma}


\subsection{The main energy-null flux identities for wave and transport equations}
We now provide the fundamental energy-null flux identity for solutions to 
$\upmu \square_{g(\threePsi)} \Psi = \waveinhom$.

\begin{remark}[\textbf{Picture of the regions of integration for the energy identities}]
\label{R:REGIONSOFINTEGRATIONFORENERGYIDENTITIES}
See Fig.~\ref{F:SOLIDREGION} on pg.~\pageref{F:SOLIDREGION} 
for a picture of the regions of integration.
Note that the (unlabeled) front and back boundaries in that figure 
should be identified.
\end{remark}

\begin{proposition}\cite{jSgHjLwW2016}*{Proposition 3.5; \textbf{Fundamental energy-null flux identity for
	the wave equation}}
	\label{P:DIVTHMWITHCANCELLATIONS}
		For scalar-valued functions $\Psi$ that solve the covariant wave equation
		\begin{align*}
			\upmu \square_{g(\threePsi)} \Psi & = \waveinhom,
		\end{align*}
		the following identity involving
		the energy and null flux from Def.~\ref{D:ENERGYFLUX}
		holds for $t \geq 1$ and $u \in [0,U_0]$:
	\begin{align} \label{E:E0DIVID}
				\enzero[\Psi](t,u)
				+
				\flzero[\Psi](t,u)
				& 
				= 
				\enzero[\Psi](0,u)
				+
				\flzero[\Psi](t,0)
					\\
				& \ \
				- 
				\int_{\mathcal{M}_{t,u}}
					\left\lbrace
						(1 + 2 \upmu) (\Lunit \Psi)
						+ 
						2 \Rad \Psi 
					\right\rbrace
				\waveinhom 
				\, d \vol
					\notag \\
		& \ \
				- 
				\frac{1}{2} 
				\int_{\mathcal{M}_{t,u}}
					\upmu \enmomtensor^{\alpha \beta}[\Psi] \deformarg{\Mult}{\alpha}{\beta}
				\, d \vol.
				\notag
			\end{align}
Furthermore, with $f_+: = \max \lbrace f,0 \rbrace$ and $f_- := \max \lbrace -f, 0 \rbrace$, we have
\begin{align}
\basicenergyerror{\Mult}[\Psi] 
	& := - \frac{1}{2} \upmu \enmomtensor^{\alpha \beta}[\Psi] \deformarg{\Mult}{\alpha}{\beta}
	 = - \frac{1}{2} \upmu |\angdiff \Psi|^2 \frac{[\Lunit \upmu]_-}{\upmu} 
	 	+ \sum_{i=1}^5 \basicenergyerrorarg{\Mult}{i}[\Psi],
		\label{E:MULTERRORINT} 
\end{align}
where
	\begin{subequations}
		\begin{align}
			\basicenergyerrorarg{\Mult}{1}[\Psi] 
				& := (\Lunit \Psi)^2 
						\left\lbrace
							- \frac{1}{2} \Lunit \upmu
							+ \Rad \upmu
							- \frac{1}{2} \upmu \mytr \upchi
							- \mytr \angk^{(Trans-\threePsi)}
							- \upmu \mytr \angk^{(Tan-\threePsi)}
						\right\rbrace,
						\label{E:MULTERRORINTEG1} \\
			\basicenergyerrorarg{\Mult}{2}[\Psi] 
			& := - (\Lunit \Psi) (\Rad \Psi)
					\left\lbrace
						 \mytr \upchi
						+ 2 \mytr \angk^{(Trans-\threePsi)}
						+ 2 \upmu \mytr \angk^{(Tan-\threePsi)}
					\right\rbrace,
						\label{E:MULTERRORINTEG2} \\
		  \basicenergyerrorarg{\Mult}{3}[\Psi] 
			& := 
				\upmu |\angdiff \Psi|^2
				\left\lbrace
					\frac{1}{2} \frac{[\Lunit \upmu]_+}{\upmu}
					+ \frac{\Rad \upmu}{\upmu}
					+ 2 \Lunit \upmu
					- \frac{1}{2} \mytr \upchi
					- \mytr \angk^{(Trans-\threePsi)}
					- \upmu \mytr \angk^{(Tan-\threePsi)}
				\right\rbrace,
				\label{E:MULTERRORINTEG3} \\
			\basicenergyerrorarg{\Mult}{4}[\Psi] 
			& := 	(\Lunit \Psi)(\angdiffuparg{\#} \Psi) 
					\cdot
					\left\lbrace
						(1 - 2 \upmu) \angdiff \upmu 
						+ 
						2 \upzeta^{(Trans-\threePsi)}
						+
						2 \upmu \upzeta^{(Tan-\threePsi)}
					\right\rbrace,
					\label{E:MULTERRORINTEG4} \\
			\basicenergyerrorarg{\Mult}{5}[\Psi] 
			& := - 2 (\Rad \Psi)(\angdiffuparg{\#} \Psi)
					\cdot
					 \left\lbrace
							\angdiff \upmu 
						 	+ 
						 	2 \upzeta^{(Trans-\threePsi)}
							+
							2 \upmu \upzeta^{(Tan-\threePsi)}
					 \right\rbrace.
					\label{E:MULTERRORINTEG5} 
		\end{align}
\end{subequations}
	The tensorfields
	$\upchi$, 
	$\upzeta^{(Trans-\threePsi)}$,
	$\angk^{(Trans-\threePsi)}$,
	$\upzeta^{(Tan-\threePsi)}$,
	and
	$\angk^{(Tan-\threePsi)}$
	from above 
	are as in
	\eqref{E:CHIDEF},
	\eqref{E:ZETATRANSVERSAL},
	\eqref{E:KABTRANSVERSAL},
	\eqref{E:ZETAGOOD},
	and \eqref{E:KABGOOD}.
\end{proposition}


In the next proposition, 
we provide the fundamental energy-null flux identity for
solutions to the transport equation
$\upmu \Transport \Vortrenormalized = \vortinhom$.
The proof relies on the following divergence identity.
\begin{lemma}\cite{jSgHjLwW2016}*{Lemma 4.3; \textbf{Spacetime divergence in terms of derivatives of frame components}}
\label{L:DIVERGENCEFRAME}
	Let $\mathscr{J}$ be a spacetime vectorfield.
	Let 
	$\upmu \mathscr{J} 
	= 
	- \upmu \mathscr{J}_{\Lunit} \Lunit 
	- \mathscr{J}_{\Rad} \Lunit 
	- \mathscr{J}_{\Lunit} \Rad 
	+ \upmu \angJ$
	be its decomposition relative to the rescaled frame, where
	$\mathscr{J}_{\Lunit} = \mathscr{J}^{\alpha} \Lunit_{\alpha}$,
	$\mathscr{J}_{\Rad} = \mathscr{J}^{\alpha} \Rad_{\alpha}$,
	and $\angJ = \Lineproject \mathscr{J}$. Then
	\begin{align} \label{E:DIVERGENCEFRAME}
	\upmu \D_{\alpha} \mathscr{J}^{\alpha} 
	& = - \Lunit (\upmu \mathscr{J}_{\Lunit})
		- \Lunit (\mathscr{J}_{\Rad})
		- \Rad (\mathscr{J}_{\Lunit})
		+ \angdiv (\upmu \angJ)
		- \upmu \mytr \angk \mathscr{J}_{\Lunit}
		- \mytr \upchi \mathscr{J}_{\Rad},
\end{align}
where where the $\ell_{t,u}$-tangent tensorfields
$\angk$
and
$\upchi$ 
can be expressed via 
\eqref{E:ANGKDECOMPOSED}
and
\eqref{E:CHIINTERMSOFOTHERVARIABLES}.
\end{lemma}

\begin{proposition}[\textbf{Energy-null flux identity for the specific vorticity}]
	\label{P:ENERGYIDENTITYRENORMALIZEDVORTICITY}
	For scalar-valued functions $\Vortrenormalized$ that solve the transport equation
	\begin{align} \label{E:RENORMALIZEDTRANSPORTEQUATIONWITHFSOURCE}
		\upmu \Transport \Vortrenormalized
		& = \vortinhom,
	\end{align}
		the following identity involving
		the energy and null flux from Def.~\ref{D:ENERGYFLUX}
		holds for $t \geq 1$ and $u \in [0,U_0]$:
	\begin{align} \label{E:ENERGYIDENTITYRENORMALIZEDVORTICITY}
		\Vortenzero[\Vortrenormalized](t,u)
		+ 
		\Vortflzero[\Vortrenormalized](t,u)
		& 
		=
		\Vortenzero[\Vortrenormalized](0,u)
		+ 
		\Vortflzero[\Vortrenormalized](t,0)
			\\
		& \ \
			+
			2
			\int_{\mathcal{M}_{t,u}}
				\Vortrenormalized \vortinhom
			\, d \vol
			\notag
			\\
		& \ \
				+
			\int_{\mathcal{M}_{t,u}}
				\left\lbrace
					\Lunit \upmu
					+
					\upmu \mytr \angk
				\right\rbrace
				\Vortrenormalized^2 
			\, d \vol.
				\notag
\end{align}
\end{proposition}

\begin{proof}
	We define the vectorfield
	$J := 
	\Vortrenormalized^2 \Transport
	=\Vortrenormalized^2 \Lunit 
	+ \Vortrenormalized^2 \Radunit$
	and note that
	$
	J_{\Lunit}
	= - \Vortrenormalized^2
	$,
	$
		J_{\Radunit}
	= 
	J_{\CoordAng}
	= 0
	$.
	Thus, using Lemma \ref{L:DIVERGENCEFRAME}
	and equation \eqref{E:RENORMALIZEDTRANSPORTEQUATIONWITHFSOURCE}, 
	we compute that
	\begin{align} \label{E:DIVTRANSPORTCURRENTEXPRESSION}
		\upmu \D_{\alpha} J^{\alpha}
		& = (\Lunit \upmu) \Vortrenormalized^2
				+
				\upmu \mytr \angk \Vortrenormalized^2
			+ 2 \Vortrenormalized \vortinhom.
	\end{align}
	Next, using the identities 
	$
	\displaystyle
	\Lunit = \frac{\partial}{\partial t}
	$
	and 
	$
	\displaystyle
	\Rad = \frac{\partial}{\partial u} - \NonRadialRad
	$
	(see \eqref{E:RADSPLITINTOPARTTILAUANDXI})
	and the relations
	$
	J_{\Lunit}
	= - \Vortrenormalized^2
	$,
	$
		J_{\Radunit}
	= 
	J_{\CoordAng}
	= 0
	$
	obtained above,
	we obtain the following decomposition from straightforward computations:
	$
	\displaystyle
		J = 
		J^t \frac{\partial}{\partial t} 
		+ J^u \frac{\partial}{\partial u} 
		+ J^{\CoordAng} \CoordAng
	$,
	where 
	$J^t = \Vortrenormalized^2$,
	$J^u = \upmu^{-1} \Vortrenormalized^2$.
	Next, we note the following formula,
	which is the standard identity for the divergence of a vectorfield
	expressed relative to a coordinate frame 
	(here the geometric coordinates)
	and the formula \eqref{E:SPACETIMEVOLUMEFORMWITHUPMU}, which
	implies that $|\mbox{\upshape{det}} g|^{1/2} = \upmu \gtancomp$
	(where the determinant is taken relative to the geometric coordinates):
	$
	\displaystyle
		\upmu \gtancomp \D_{\alpha} J^{\alpha}
		= 
		\frac{\partial}{\partial t}
		 	\left(
		 		\upmu \gtancomp J^t
		 	\right)
		+
		 	\frac{\partial}{\partial u}
		 	\left(
		 		\upmu \gtancomp J^u
		 	\right)
	+
		 	\frac{\partial}{\partial \vartheta}
		 	\left(
		 		\upmu \gtancomp J^{\CoordAng}
		 	\right)
	$.
	Integrating this identity over $\mathcal{M}_{t,u}$
	with respect to 
	$dt' \, du' \, d \vartheta$
	and referring to Def.~\ref{D:NONDEGENERATEVOLUMEFORMS}, 
	we obtain
	\begin{align} \label{E:FIRSTDIVERGENCEID}
		\int_{\mathcal{M}_{t,u}}
			\upmu \D_{\alpha} J^{\alpha}
		 \, d \vol
		 & = 
		 \int_{t'=0}^t
		 \int_{u'=0}^u
		 \int_{\vartheta \in \mathbb{T}}
		 	\frac{\partial}{\partial t}
		 	\left(
		 		\upmu \gtancomp \Vortrenormalized^2
		 	\right)
		 	+
		 	\frac{\partial}{\partial u}
		 	\left(
		 		\gtancomp \Vortrenormalized^2
		 	\right)
		 	+
		 	\frac{\partial}{\partial \vartheta}
		 	\left(
		 		\upmu \gtancomp J^{\CoordAng}
		 	\right)
		 \, dt' 
		 \, du'
		 \, d \vartheta.
	\end{align}
	The desired identity \eqref{E:ENERGYIDENTITYRENORMALIZEDVORTICITY}
	now follows from 
	\eqref{E:DIVTRANSPORTCURRENTEXPRESSION},
	\eqref{E:FIRSTDIVERGENCEID},
	definition \eqref{E:RENORMALIZEDVORTENERGYFLUX},
	Fubini's theorem,
	and the fact that
	the integral of the last term
	$
	\displaystyle
	\frac{\partial}{\partial \vartheta}
		 	\left(
		 		\upmu \gtancomp J^{\CoordAng}
		 	\right)
	$
	over $\mathbb{T}$ vanishes.
\end{proof}

\subsection{Additional integration by parts identities}
\label{SS:ADDITONALIBPIDENTITIES}
In this section, we provide, for future use, some integration by parts identities.
We highlight here the identity \eqref{E:LUNITANDANGULARIBPIDENTITY}, 
which plays a critical role in our 
top-order energy estimates; see 
equation \eqref{E:TOPORDERLUNITIBPDIFFICULTTERMS} and just below it. 

\begin{lemma}\cite{jSgHjLwW2016}*{Lemma 3.6; \textbf{Identities connected to integration by parts}}
\label{L:LDERIVATIVEOFLINEINTEGRAL}
	The following identities hold for scalar-valued functions $f$:
	\begin{subequations}
	\begin{align} \label{E:LDERIVATIVEOFLINEINTEGRAL}
		\frac{\partial}{\partial t}
		\int_{\ell_{t,u}}
			f
		\, d \spherevol
		& 
		= 
		\int_{\ell_{t,u}}
			\left\lbrace
				\Lunit f
				+ 
				\mytr \upchi f
			\right\rbrace
		\, d \spherevol,
			\\
		\frac{\partial}{\partial u}
		\int_{\ell_{t,u}}
			f
		\, d \spherevol
		& 
		= 
		\int_{\ell_{t,u}}
			\left\lbrace
				\Rad f
				+ 
				\frac{1}{2} \mytr \angdeform{\Rad} 
				f
				\right\rbrace
		\, d \spherevol.
			\label{E:UDERIVATIVEOFLINEINTEGRAL}
	\end{align}
	\end{subequations}

	In addition, the following integration by parts identity holds for 
	scalar-valued functions $\Psi$ and $\ThirdSmoothFunction$
	(see Sect.~\ref{SS:STRINGSOFCOMMUTATIONVECTORFIELDS} regarding the vectorfield operator notation):
	\begin{align}	\label{E:LUNITANDANGULARIBPIDENTITY}
		& 
		\int_{\mathcal{M}_{t,u}}
			(1 + 2 \upmu) (\Rad \Psi) (\Lunit \Fullset_{\ast}^{N;\leq 1} \Psi) \GeoAng \ThirdSmoothFunction
		\, d \vol
			\\
	& = 
		\int_{\mathcal{M}_{t,u}}
			(1 + 2 \upmu) (\Rad \Psi) (\GeoAng \Fullset_{\ast}^{N;\leq 1} \Psi) \Lunit \ThirdSmoothFunction
		\, d \vol
		\notag	\\
	& \ \ 
		- \int_{\Sigma_t^u}
				(1 + 2 \upmu) (\Rad \Psi) (\GeoAng \Fullset_{\ast}^{N;\leq 1} \Psi) \ThirdSmoothFunction
			\, d \vol
		+ \int_{\Sigma_0^u}
				(1 + 2 \upmu) (\Rad \Psi) (\GeoAng \Fullset_{\ast}^{N;\leq 1} \Psi) \ThirdSmoothFunction
			\, d \vol
			\notag \\
	& \ \
		+
		\int_{\mathcal{M}_{t,u}}
			\mbox{\upshape Error}_1[\Fullset_{\ast}^{N;\leq 1} \Psi;\ThirdSmoothFunction]
		\, d \vol
		+ \int_{\Sigma_t^u}
				\mbox{\upshape Error}_2[\Fullset_{\ast}^{N;\leq 1} \Psi;\ThirdSmoothFunction]
			\, d \vol
		- \int_{\Sigma_0^u}
				\mbox{\upshape Error}_2[\Fullset_{\ast}^{N;\leq 1} \Psi;\ThirdSmoothFunction]
			\, d \vol,
			\notag
	\end{align}
	where
	\begin{subequations}
	\begin{align} \label{E:LUNITIBPSPACETIMEERROR}
		\mbox{\upshape Error}_1[\Fullset_{\ast}^{N;\leq 1} \Psi;\ThirdSmoothFunction] 
		& := 2 (\Lunit \upmu) (\Rad \Psi) (\GeoAng \Fullset_{\ast}^{N;\leq 1} \Psi) \ThirdSmoothFunction
			+ (1 + 2 \upmu) (\Lunit \Rad \Psi) (\GeoAng \Fullset_{\ast}^{N;\leq 1} \Psi) \ThirdSmoothFunction
				\\
		& \ \
			+ (1 + 2 \upmu) (\Rad \Psi) (\angdeformoneformupsharparg{\GeoAng}{\Lunit} \cdot \angdiff \Fullset_{\ast}^{N;\leq 1} \Psi) \ThirdSmoothFunction
			+  (1 + 2 \upmu) (\Rad \Psi) \mytr \upchi (\GeoAng \Fullset_{\ast}^{N;\leq 1} \Psi) \ThirdSmoothFunction
			\notag \\
		& \ \
			+ 2 (\GeoAng \upmu) (\Rad \Psi) (\Fullset_{\ast}^{N;\leq 1} \Psi) \Lunit \ThirdSmoothFunction
			+ (1 + 2 \upmu) (\GeoAng \Rad \Psi) (\Fullset_{\ast}^{N;\leq 1} \Psi) \Lunit \ThirdSmoothFunction
				\notag \\
		& \ \ + \frac{1}{2} (1 + 2 \upmu) (\Rad \Psi) \mytr \angdeform{Y} (\Fullset_{\ast}^{N;\leq 1} \Psi) \Lunit \ThirdSmoothFunction
				\notag \\
		& \ \ 
			+  2 (\Lunit \GeoAng \upmu) (\Rad \Psi) (\Fullset_{\ast}^{N;\leq 1} \Psi) \ThirdSmoothFunction
			+  2 (\GeoAng \upmu) (\Lunit \Rad \Psi) (\Fullset_{\ast}^{N;\leq 1} \Psi) \ThirdSmoothFunction
				\notag \\
		& \ \
			+  2 (\GeoAng \upmu) (\Rad \Psi) \mytr \upchi (\Fullset_{\ast}^{N;\leq 1} \Psi) \ThirdSmoothFunction
			+  (\Lunit \upmu) (\Rad \Psi) \mytr \angdeform{Y} (\Fullset_{\ast}^{N;\leq 1} \Psi) \ThirdSmoothFunction
				\notag \\
		& \ \
			+  (1 + 2 \upmu) (\Lunit \GeoAng \Rad \Psi) (\Fullset_{\ast}^{N;\leq 1} \Psi) \ThirdSmoothFunction
			+  (1 + 2 \upmu) (\Rad \Psi) (\GeoAng \mytr \upchi) (\Fullset_{\ast}^{N;\leq 1} \Psi) \ThirdSmoothFunction
			\notag \\
		& \ \
		   +  (1 + 2 \upmu) (\GeoAng \Rad \Psi) \mytr \upchi (\Fullset_{\ast}^{N;\leq 1} \Psi) \ThirdSmoothFunction
			+ \frac{1}{2} (1 + 2 \upmu) (\Lunit \Rad \Psi) \mytr \angdeform{Y} (\Fullset_{\ast}^{N;\leq 1} \Psi) \ThirdSmoothFunction 
			 \notag \\
		& \ \
			+ (1 + 2 \upmu) (\Rad \Psi) (\angdiv \angdeformoneformupsharparg{\GeoAng}{\Lunit}) (\Fullset_{\ast}^{N;\leq 1} \Psi) \ThirdSmoothFunction
			+ \frac{1}{2} (1 + 2 \upmu) (\Rad \Psi) \mytr \upchi \mytr \angdeform{Y} (\Fullset_{\ast}^{N;\leq 1} \Psi) \ThirdSmoothFunction,
				\notag
				\\
		\mbox{\upshape Error}_2[\Fullset_{\ast}^{N;\leq 1} \Psi;\ThirdSmoothFunction] 
		& := 
			- 2 (\GeoAng \upmu) (\Rad \Psi) (\Fullset_{\ast}^{N;\leq 1} \Psi) \ThirdSmoothFunction
			- (1 + 2 \upmu) (\GeoAng \Rad \Psi) (\Fullset_{\ast}^{N;\leq 1} \Psi) \ThirdSmoothFunction
				\label{E:LUNITIBPHYPERSURFACEERROR} \\
		& \ \
			- \frac{1}{2} (1 + 2 \upmu) (\Rad \Psi) \mytr \angdeform{Y} (\Fullset_{\ast}^{N;\leq 1} \Psi) \ThirdSmoothFunction.
			\notag	
	\end{align}
	\end{subequations}

\end{lemma}

\section{The commutator of the covariant wave operator and a vectorfield}
\label{S:WAVEEQUATIONONECOMMUTATION}
\setcounter{equation}{0}
In this section, we provide expressions for the
commutators $[\upmu \square_{g(\threePsi)},Z]$
for vectorfields $Z$ belonging to the commutation set
$\Fullset$ defined in \eqref{E:COMMUTATIONVECTORFIELDS}.
The following lemma provides a first decomposition.
In Prop.~\ref{P:COMMUTATIONCURRENTDIVERGENCEFRAMEDECOMP}, we
further decompose the main term from the lemma.

\begin{lemma}\cite{jSgHjLwW2016}*{Lemma 4.2; \textbf{Vectorfield-covariant wave operator commutator identity}}
\label{L:BOXZCOM}
For $Z \in \Fullset$ (see Def.~\ref{E:COMMUTATIONVECTORFIELDS}),
we have the following commutation identity
for scalar functions $\Psi$,
where $\myspacetimetr \deform{Z} := (g^{-1})^{\alpha \beta} \deformarg{Z}{\alpha}{\beta}$:
\begin{align} \label{E:BOXZCOM}
		\upmu \square_{g(\Psi)} (Z \Psi)
		& =   \upmu 
					\D_{\alpha} 
					\left\lbrace
					\deformuparg{Z}{\alpha}{\beta} \D_{\beta} \Psi
						- \frac{1}{2} \myspacetimetr \deform{Z} \D^{\alpha} \Psi 
					\right\rbrace
				+ Z (\upmu \square_{g(\Psi)} \Psi)
				+ \frac{1}{2} \mytr \angdeform{Z} (\upmu \square_{g(\Psi)} \Psi).
\end{align}

\end{lemma}

We now decompose the first term on RHS~\eqref{E:BOXZCOM}
relative to the rescaled frame.

\begin{proposition}\cite{jSgHjLwW2016}*{Proposition 4.4; \textbf{Frame decomposition of the divergence of 
the main inhomogeneous term in the commuted wave equation}} 
\label{P:COMMUTATIONCURRENTDIVERGENCEFRAMEDECOMP}
	For vectorfields $Z \in \Fullset$,
	we have the following identity for the
	first term on RHS~\eqref{E:BOXZCOM}:
	\begin{align} \label{E:DIVCOMMUTATIONCURRENTDECOMPOSITION}
				\upmu 
					\D_{\alpha} 
					\left\lbrace
						\deformuparg{Z}{\alpha}{\beta} \D_{\beta} \Psi
						- \frac{1}{2} \myspacetimetr \deform{Z} \D^{\alpha} \Psi 
					\right\rbrace
		& = \mathscr{K}_{(\pi-Danger)}^{(Z)}[\Psi]
			\\
		& \ \
			+ \mathscr{K}_{(\pi-Cancel-1)}^{(Z)}[\Psi]
			+ \mathscr{K}_{(\pi-Cancel-2)}^{(Z)}[\Psi]
			\notag	\\
		& \ \ 
			+ \mathscr{K}_{(\pi-Less \ Dangerous)}^{(Z)}[\Psi]
			+ \mathscr{K}_{(\pi-Good)}^{(Z)}[\Psi]
				\notag \\
		& \	\
			+ \mathscr{K}_{(\Psi)}^{(Z)}[\Psi]
			+ \mathscr{K}_{(Low)}^{(Z)}[\Psi],
			\notag
	\end{align}
	where
	\begin{subequations}
		\begin{align}
			\mathscr{K}_{(\pi-Danger)}^{(Z)}[\Psi]
			& := - (\angdiv \angdeformoneformupsharparg{Z}{\Lunit}) \Rad \Psi,
				\label{E:DIVCURRENTTRANSVERSAL}
				\\
			\mathscr{K}_{(\pi-Cancel-1)}^{(Z)}[\Psi]
			& := \left\lbrace 
						\frac{1}{2} \Rad \mytr  \angdeform{Z}
						- \angdiv \angdeformoneformupsharparg{Z}{\Rad}
						- \upmu \angdiv \angdeformoneformupsharparg{Z}{\Lunit}
					\right\rbrace 
						\Lunit \Psi,
					\label{E:DIVCURRENTCANEL1} \\
			\mathscr{K}_{(\pi-Cancel-2)}^{(Z)}[\Psi]
			& :=
				\left\lbrace
					- \angLie_{\Rad} \angdeformoneformupsharparg{Z}{\Lunit}
					+ \angdiffuparg{\#} \deformarg{Z}{\Lunit}{\Rad}
				\right\rbrace 
				\cdot
				\angdiff \Psi,
				\label{E:DIVCURRENTCANEL2} \\
			\mathscr{K}_{(\pi-Less \ Dangerous)}^{(Z)}[\Psi]
			& := \frac{1}{2} \upmu (\angdiffuparg{\#} \mytr \angdeform{Z}) \cdot \angdiff \Psi, 
				\label{E:DIVCURRENTELLIPTIC} \\
			\mathscr{K}_{(\pi-Good)}^{(Z)}[\Psi] 
			& := \frac{1}{2} \upmu (\Lunit \mytr \angdeform{Z}) \Lunit \Psi
				+ (\Lunit \deformarg{Z}{\Lunit}{\Rad}) \Lunit \Psi
				+ (\Lunit \deformarg{Z}{\Rad}{\Radunit}) \Lunit \Psi
				\label{E:DIVCURRENTGOOD} \\
			& \ \ + \frac{1}{2} (\Lunit \mytr \angdeform{Z}) \Rad \Psi
				- \upmu (\angLie_{\Lunit} \angdeformoneformupsharparg{Z}{\Lunit}) \cdot \angdiff \Psi
				- (\angLie_{\Lunit} \angdeformoneformupsharparg{Z}{\Rad}) \cdot \angdiff \Psi,
				\notag 
	\end{align}
	\end{subequations}
	\begin{align}
		\mathscr{K}_{(\Psi)}^{(Z)}[\Psi]
			& := \left\lbrace
							\frac{1}{2}
							\upmu \mytr \angdeform{Z}
							+ \deformarg{Z}{\Lunit}{\Rad}
							+ \deformarg{Z}{\Rad}{\Radunit}
						\right\rbrace
						\Lunit \Lunit \Psi  
					\label{E:DIVCURRENTPSI} \\
		& \ \
				+ \mytr \angdeform{Z}
					 \Lunit \Rad \Psi
					\notag \\
			& \ \ - 2 \upmu \angdeformoneformupsharparg{Z}{\Lunit} \cdot \angdiff \Lunit \Psi
					- 2 \angdeformoneformupsharparg{Z}{\Rad} \cdot \angdiff \Lunit \Psi
					- 2 \angdeformoneformupsharparg{Z}{\Lunit} \cdot \angdiff \Rad \Psi
				\notag \\
			& \ \ + \deformarg{Z}{\Lunit}{\Rad} \angLap \Psi
				+ \frac{1}{2} \upmu \mytr \angdeform{Z} \angLap \Psi,
				\notag 
	\end{align}
	and
	\begin{align}
		\mathscr{K}_{(Low)}^{(Z)}[\Psi] 
		& := \left\lbrace
				 	\frac{1}{2} (\Lunit \upmu) \mytr \angdeform{Z}
				 	+ \frac{1}{2} \upmu \mytr \angk \mytr \angdeform{Z}
				 	+ \mytr \upchi \deformarg{Z}{\Lunit}{\Rad}
				 	+ \mytr \upchi \deformarg{Z}{\Rad}{\Radunit}
				 	- \angdeformoneformupsharparg{Z}{\Lunit} \cdot \angdiff \upmu
				 \right\rbrace
				 \Lunit \Psi
				\label{E:DIVCURRENTLOW}  \\
		& \ \ 
				+
				\frac{1}{2} \mytr \upchi \mytr \angdeform{Z} \Rad \Psi
				\notag \\
		& \ \  
				+ \left\lbrace
						- (\Lunit \upmu) \angdeformoneformupsharparg{Z}{\Lunit}
						- \upmu \mytr \angk \angdeformoneformupsharparg{Z}{\Lunit}
						- \mytr \upchi \angdeformoneformupsharparg{Z}{\Rad}
						+ \mytr \angdeform{Z} \angdiffuparg{\#} \upmu
						+  \mytr \upchi \upmu \upzeta^{\#}
					\right\rbrace
					\cdot
					\angdiff \Psi.
				\notag
	\end{align}
	In the above expressions, the $\ell_{t,u}$-tangent tensorfields
	$\upchi$, 
	$\upzeta$,
	and
	$\angk$, 
	are as in
	\eqref{E:CHIINTERMSOFOTHERVARIABLES},
	\eqref{E:ZETADECOMPOSED},
	and 
	\eqref{E:ANGKDECOMPOSED}.
\end{proposition}

\section{Norms and Strings of Commutation Vectorfields}
\label{S:NORMS}
In this section, we define various norms and seminorms. We also
introduce schematic notation that succinctly captures
the most important properties of strings of commutation
vectorfields.

\subsection{Norms}
\label{SS:NORMS}
We now define some norms that we use in our analysis.
We recall that we defined the pointwise norm of $\ell_{t,u}$-tensors 
(relative to $\gsphere$) in Subsect.~\ref{SS:POINTWISENORMS}.

\subsubsection{Lebesgue norms}
\label{SSS:LEBESGUENORMS}

\begin{definition}[$L^2$ \textbf{and} $L^{\infty}$ \textbf{norms}]
In terms of the \textbf{non-degenerate} forms of Def.~\ref{D:NONDEGENERATEVOLUMEFORMS},
we define the following norms for 
$\ell_{t,u}$-tangent tensorfields:
\label{D:SOBOLEVNORMS}
	\begin{subequations}
	\begin{align}  \label{E:L2NORMS}
			\left\|
				\xi
			\right\|_{L^2(\ell_{t,u})}^2
			& :=
			\int_{\ell_{t,u}}
				|\xi|^2
			\, d \spherevol,
				\qquad
			\left\|
				\xi
			\right\|_{L^2(\Sigma_t^u)}^2
			:=
			\int_{\Sigma_t^u}
				|\xi|^2
			\, d \tvol,
				\\
			\left\|
				\xi
			\right\|_{L^2(\mathcal{P}_u^t)}^2
			& :=
			\int_{\mathcal{P}_u^t}
				|\xi|^2
			\, d \conevol,
			\notag
	\end{align}

	\begin{align} 
			\left\|
				\xi
			\right\|_{L^{\infty}(\ell_{t,u})}
			& :=
				\mbox{ess sup}_{\vartheta \in \mathbb{T}}
				|\xi|(t,u,\vartheta),
			\qquad
			\left\|
				\xi
			\right\|_{L^{\infty}(\Sigma_t^u)}
			:=
			\mbox{ess sup}_{(u',\vartheta) \in [0,u] \times \mathbb{T}}
				|\xi|(t,u',\vartheta),
			\label{E:LINFTYNORMS}
				\\
			\left\|
				\xi
			\right\|_{L^{\infty}(\mathcal{P}_u^t)}
			& :=
			\mbox{ess sup}_{(t',\vartheta) \in [0,t] \times \mathbb{T}}
				|\xi|(t',u,\vartheta).
			\notag
	\end{align}
	\end{subequations}
\end{definition}

\begin{remark}[\textbf{Subset norms}]
	\label{R:SUBSETNORMS}
	In our analysis below, we occasionally use norms 
	$\| \cdot \|_{L^2(\Omega)}$
	and
	$\| \cdot \|_{L^{\infty}(\Omega)}$,
	where $\Omega$ is a subset of $\Sigma_t^u$.
	These norms are defined by replacing 
	$\Sigma_t^u$ with $\Omega$ in 
	\eqref{E:L2NORMS} and \eqref{E:LINFTYNORMS}.
\end{remark}

\subsubsection{Norms of arrays}
\label{SSS:ARRAYNORMS}
We define the norms of the arrays $\vec{G}_{(Frame)}$ and $\vec{H}_{(Frame)}$
from Def.~\ref{D:BIGGANDBIGH} to be the sums of the norms of their $\imath,\jmath$-indexed 
entries. For example,
\begin{align}\label{E:GFRAMEPOINTWISENORMEXAPLE}
		\left|
			\vec{G}_{(Frame)} 
		\right|
		& := 
				\left|
					\vec{G}_{\Lunit \Lunit} 
				\right| 
				+
				\left|
					\vec{G}_{\Lunit \Radunit}
				\right|
				+
				\left|
					\angGarg{\Lunit}
				\right|
				+
				\left|
					\angGarg{\Radunit}
				\right|
				+
				\left|
					\angG
				\right|,
\end{align}
where
$\left|
		\vec{G}_{\Lunit \Lunit}
	\right|
	:= \sum_{\imath=0}^2
				\left|
					G_{\Lunit \Lunit}^{\imath}
				\right|
$,
$
	\left|
		\angGarg{\Radunit}
	\right| 
	:= \sum_{\iota=0}^2 \left|\NovecangGarg{\Radunit}{\imath}\right|
$, etc.
We similarly define
$\left\|
		\vec{G}_{(Frame)} 
\right\|_{L^{\infty}(\Sigma_t^u)}
$
and
$\left\|
		\vec{H}_{(Frame)} 
\right\|_{L^{\infty}(\Sigma_t^u)}
$,
and similarly 
for other norms.

\subsection{Strings of commutation vectorfields and vectorfield seminorms}
\label{SS:STRINGSOFCOMMUTATIONVECTORFIELDS}
The following shorthand notation captures the important structural features
of various differential operators corresponding to repeated differentiation
with respect to the commutation vectorfields. The notation allows us to schematically depict 
identities and estimates.

\begin{definition}[\textbf{Strings of commutation vectorfields and vectorfield seminorms}] \label{D:VECTORFIELDOPERATORS}
	\ \\
	\begin{itemize}
		\item $\Fullset^{N;M} f$ 
			denotes an arbitrary string of $N$ commutation
			vectorfields in $\Fullset$ (see \eqref{E:COMMUTATIONVECTORFIELDS})
			applied to $f$, where the string contains \emph{precisely} $M$ factors of the 
			$\mathcal{P}_t^u$-transversal
			vectorfield $\Rad$. 
			We also set $\Fullset^{0;0} f := f$.
			Similarly, we write
			$\Fullset^{N;\leq M} f$ when the string is allowed to contain
			$\leq M$ factors of $\Rad$.
		\item $\Tanset^N f$
			denotes an arbitrary string of $N$ commutation
			vectorfields in $\Tanset$ (see \eqref{E:TANGENTIALCOMMUTATIONVECTORFIELDS})
			applied to $f$.
		\item 
			For $N \geq 1$,
			$\Fullset_{\ast}^{N;M} f$
			denotes an arbitrary string of $N$ commutation
			vectorfields in $\Fullset$ 
			applied to $f$, where the string contains \emph{at least} one $\mathcal{P}_u$-tangent factor 
			and \emph{precisely} $M$ factors of $\Rad$.
			We also set $\Fullset_{\ast}^{0;0} f := f$.
			Similarly, we write
			$\Fullset_{\ast}^{N;\leq M} f$ when the string is allowed to contain
			$\leq M$ factors of $\Rad$.
	\item 
			For $N \geq 1$,
			$\Fullset_{\ast \ast}^{N;M} f$
			denotes an arbitrary string of $N$ commutation
			vectorfields in $\Fullset$ 
			applied to $f$, 
			where the string contains \emph{at least} two factors of $\Lunit$
			or \emph{at least} one factor of $\GeoAng$
			and \emph{precisely} $M$ factors of $\Rad$.
			Similarly, we write
			$\Fullset_{\ast \ast}^{N;\leq M} f$ when the string is allowed to contain
			$\leq M$ factors of $\Rad$.
	\item For $\ell_{t,u}$-tangent tensorfields $\xi$, 
					we similarly define strings of $\ell_{t,u}$-projected Lie derivatives 
					such as $\angLie_{\Fullset}^{N;M} \xi$.
	\end{itemize}

	We also define pointwise seminorms constructed out of sums of strings of vectorfields:
	\begin{itemize}
		\item $\left|\Fullset^{N;M} f \right|$ 
		simply denotes the magnitude of one of the $\Fullset^{N;M} f$ as defined above
		(there is no summation).
		Similarly, $\left|\Fullset^{N;\leq M} f \right|$
		denotes the magnitude of one of the $\Fullset^{N;\leq M} f$ as defined above.
		\item $\left|\Fullset^{\leq N;M} f \right|$ is the \emph{sum} over all terms of the form 
		$\left|\Fullset^{N';M} f \right|$ with $N' \leq N$.
		\item $\left|\Fullset^{\leq N;\leq M} f \right|$ is the sum over all terms of the form 
		$\left|\Fullset^{N';M'} f \right|$
			with $N' \leq N$ and $M' \leq M$.
		\item $\left|\Fullset^{[1,N];M} f \right|$ is the sum over all terms of the form 
		$\left|\Fullset^{N';M} f \right|$
			with $1 \leq N' \leq N$.
		\item $\left|\Fullset^{[1,N];\leq M} f \right|$ is the sum over all terms of the form 
		$\left|\Fullset^{N';M'} f \right|$ with $1 \leq N' \leq N$ and $M' \leq M$.
		\item Quantities such as
		  $\left|\Tanset^N f \right|$,
			$\left|\Fullset_{\ast}^{N;M} f \right|$,
			$\left|\Fullset_{\ast \ast}^{N;M} f \right|$,
			and
			$\left|\Fullset_{\ast \ast}^{N;\leq M} f \right|$
			are defined analogously
			(without summation).
		\item Sums such as 
			$\left|\Tanset^{\leq N} f \right|$,
			$\left|\Tanset^{[1,N]} f \right|$,
			$\left|\Fullset_{\ast}^{[1,N];M} f \right|$,
			$\left|\Fullset_{\ast}^{[1,N];\leq M} f \right|$,
			$\left|\Fullset_{\ast \ast}^{[1,N];M} f \right|$,
			$\left|\Fullset_{\ast \ast}^{[1,N];\leq M} f \right|$,
			$\left|\GeoAng^{\leq N} f \right|$,
			and $\left|\Rad^{[1,N]} f \right|$
			are defined analogously.
			For example,
			$\left|\Rad^{[1,N]} f \right| 
			= |\Rad f| 
			+ |\Rad \Rad f| 
			+ \cdots 
			+ \left|
				\overbrace{\Rad \Rad \cdots \Rad}^{N \mbox{ \upshape copies}} f
			\right|
			$.
	\end{itemize}

\end{definition}

\begin{remark}[\textbf{Operators decorated with} $\ast$ \textbf{or} $\ast \ast$]
	\label{R:DECORATEDOPERATROS}
	The purpose of the symbols 
	$\ast$ and $\ast \ast$
	in Def.~\ref{D:VECTORFIELDOPERATORS}
	is to highlight the presence of special structures in vectorfield operators,
	which helps us track smallness in the estimates.
	That is, in our analysis, we typically display operators 
	decorated with a $\ast$ and $\ast \ast$
	when they lead to
	quantities that are initially\footnote{At the high derivative levels, 
	the ``initially small'' quantities are allowed to blow up
	like $\mathring{\upepsilon} (\min_{\Sigma_t} \upmu)^{-P}$ for some power $P$
	as the shock forms.} 
	of small size $\mathcal{O}(\mathring{\upepsilon})$,
	where $\mathring{\upepsilon}$ is the data-size parameter
	defined in Sect.~\ref{S:DATASIZEANDBOOTSTRAP}.
	We note here that the quantities
	$\Fullset_{\ast}^{N;M} \GdVar$ 
	and $\Fullset_{\ast \ast}^{N;M} \BadVar$
	are always initially small,
	while $\Fullset_{\ast}^{N;M} \BadVar$ may not be.
	The reason 
	that $\Fullset_{\ast}^{N;M} \BadVar$ may not be small is: 
	for the solutions under consideration,
	$\Lunit \upmu$ and its $\Rad$ derivatives
	are large quantities.
	We also note that the notation
	$\ast$
	and
	$\ast \ast$
	is not important\footnote{We use it nonetheless for consistency.}
	for treating the specific vorticity variable $\Vortrenormalized$
	because our initial conditions are such that all directional
	derivatives of the specific vorticity are initially small.
\end{remark}

\section{Modified quantities}
\label{S:MODIFIED}
In this section, we define the modified quantities that allow us to 
avoid losing a derivative at the top-order. We also define the partially 
modified quantities that allow us to avoid some top-order
error integrals with magnitudes that are too large for us to control. 
We then provide transport-type evolution equations for these quantities.

\subsection{Curvature tensors and the key Ricci component identity}
\label{SS:CURVATURETENSORSRICCIID}
We use use curvature tensors of $g$
to help us organize the calculations 
in this section.

\begin{definition}[\textbf{Curvature tensors} \textbf{of} $g$]
\label{D:SPACETIMECURVATURE}
The Riemann curvature tensor $\Cur_{\alpha \beta \kappa \lambda}$ of the spacetime metric $g$ is
the type $\binom{0}{4}$ spacetime tensorfield defined by
\begin{align} \label{E:SPACETIMERIEMANN}
	g(\D_{UV}^2 W - \D_{VU}^2 W,Z)
	& = - \Cur(U,V,W,Z),
\end{align}
where $U$, $V$, $W$, and $Z$ are arbitrary spacetime vectors.
In \eqref{E:SPACETIMERIEMANN}, 
$\D_{UV}^2 W := U^{\alpha} V^{\beta} \D_{\alpha} \D_{\beta} W$.

The Ricci curvature tensor $\mbox{\upshape Ric}_{\alpha \beta}$ of $g$ is the following type 
$\binom{0}{2}$ tensorfield:
\begin{align} \label{E:RICCIDEF}
	\mbox{\upshape Ric}_{\alpha \beta}
	& := (g^{-1})^{\kappa \lambda} \Cur_{\alpha \kappa \beta \lambda}.
\end{align}
\end{definition}

The next lemma lies at the heart of the construction of the modified quantities.

\begin{lemma}[\textbf{The key identity verified by} $\upmu \mbox{\upshape Ric}_{\Lunit \Lunit}$]
\label{L:RICLLRENORMALIZED}
Assume that the entries of $\threePsi = (\Densrenormalized,v^1,v^2)$
verify the geometric wave equation system
\eqref{E:VELOCITYWAVEEQUATION}-\eqref{E:RENORMALIZEDDENSITYWAVEEQUATION}.
Then the following identity holds for the Ricci curvature component 
$\mbox{\upshape Ric}_{\Lunit \Lunit} := \mbox{\upshape Ric}_{\alpha \beta} \Lunit^{\alpha} \Lunit^{\beta}$:
\begin{align} \label{E:RICLLRENORMALIZED}
	\upmu \mbox{\upshape Ric}_{\Lunit \Lunit}
	& = \Lunit
			\left\lbrace
				- \vec{G}_{\Lunit \Lunit}\contr\Rad \threePsi
				- \frac{1}{2} \upmu \mytr \angG\contr\Lunit \threePsi
				- \frac{1}{2} \upmu \vec{G}_{\Lunit \Lunit}\contr\Lunit \threePsi
				+ \upmu \angGnospacemixedarg{\Lunit}{\#}\contr\cdot\angdiff \threePsi
			\right\rbrace
		+ \mathfrak{A},
\end{align}
where $\mathfrak{A}$ has the following schematic structure:
\begin{align} \label{E:RENORMALIZEDRICLLINHOMOGENEOUSTERM}
	\mathfrak{A}
	& = \smoothfunction(\BadVar,\ginversesphere,\angdiff \vec{x},\Rad \threePsi, \Singletan \threePsi)
			\Singletan \threePsi
			+ \upmu \smoothfunction(\GdVar) \Singletan \Vortrenormalized
			+ \Vortrenormalized \smoothfunction(\GdVar) \Rad \threePsi
			+ \upmu \Vortrenormalized \smoothfunction(\GdVar) \Singletan \threePsi.
\end{align}

Furthermore, without assuming 
that equations
\eqref{E:VELOCITYWAVEEQUATION}-\eqref{E:RENORMALIZEDDENSITYWAVEEQUATION}
hold, we have
\begin{align} \label{E:RICLLPARTIALRENORMALIZED}
	\mbox{\upshape Ric}_{\Lunit \Lunit}
	& = \frac{(\Lunit \upmu)}{\upmu} \mytr \upchi
		+ \Lunit
			\left\lbrace
				- \frac{1}{2} \mytr \angG\contr\Lunit \threePsi
				- \frac{1}{2} \vec{G}_{\Lunit \Lunit}\contr\Lunit \threePsi 
				+ \angGnospacemixedarg{\Lunit}{\#}\contr\cdot \angdiff \threePsi
			\right\rbrace
		- \frac{1}{2} \vec{G}_{\Lunit \Lunit}\contr\angLap \threePsi
		+ \mathfrak{B},
\end{align}
where $\mathfrak{B}$ has the following schematic structure:
\begin{align} \label{E:PARTIALRENORMALIZEDRICLLNHOMOGENEOUSTERM}
	\mathfrak{B}
	& = 	\smoothfunction(\GdVar,\ginversesphere,\angdiff \vec{x})
				(\Singletan \threePsi) 
				\Singletan \GdVar.
\end{align}

\end{lemma}
\begin{proof}[Sketch of proof]
	The identities \eqref{E:RICLLRENORMALIZED} and \eqref{E:RICLLPARTIALRENORMALIZED}
	were essentially proved in \cite{jSgHjLwW2016}*{Lemma 6.1}
	using calculations along the lines of those in \cite{dC2007}*{Chapter 8}.
	The only new feature in the present work
	is that RHS~\eqref{E:RICLLRENORMALIZED}
	depends on the inhomogeneous terms on the right-hand sides of
	the wave equations \eqref{E:VELOCITYWAVEEQUATION}-\eqref{E:RENORMALIZEDDENSITYWAVEEQUATION},
	which were absent in the previous works.
	The inhomogeneous terms appear because
	at the key point in the proof,
	one uses \eqref{E:LONOUTSIDEGEOMETRICWAVEOPERATORFRAMEDECOMPOSED}
and the wave equations
\eqref{E:VELOCITYWAVEEQUATION}-\eqref{E:RENORMALIZEDDENSITYWAVEEQUATION}
to express
	\begin{align} \label{E:SKETCHRICLLWAVEEQUATIONTERMNORMALIZED}
		- \frac{1}{2} 
		\upmu \vec{G}_{\Lunit \Lunit}\contr\angLap \threePsi 
		& = 
				-
				\frac{1}{2} 
				\Lunit
				\left\lbrace
					\vec{G}_{\Lunit \Lunit}
					\contr
					(\upmu \Lunit \threePsi + 2 \Rad \threePsi)
				\right\rbrace
			- 
			\frac{1}{2} \mytr \upchi \vec{G}_{\Lunit \Lunit}\contr\Rad \threePsi
				\\
		& \ \
			+ \mbox{\upshape Inhom}
				\notag \\
		& \ \
				+ 
				\myarray
				[\vec{G}_{(Frame)}^2 \ginversesphere]
				{\vec{H}_{(Frame)}}
				\threemyarray[\upmu \Lunit \threePsi]
					{\Rad \threePsi}
					{\upmu \angdiff \threePsi}
				\myarray
					[\Lunit \threePsi]
					{\angdiff \threePsi},
					\notag
	\end{align}
	where the last line of RHS~\eqref{E:SKETCHRICLLWAVEEQUATIONTERMNORMALIZED} is schematically depicted and
	term $\mbox{\upshape Inhom}$ on RHS~\eqref{E:SKETCHRICLLWAVEEQUATIONTERMNORMALIZED}
	denotes the inhomogeneous terms on RHSs
	\eqref{E:VELOCITYWAVEEQUATION}-\eqref{E:RENORMALIZEDDENSITYWAVEEQUATION}.
	The first term on RHS~\eqref{E:SKETCHRICLLWAVEEQUATIONTERMNORMALIZED}
	is incorporated into the perfect $\Lunit$ derivative term on the first line
	of RHS~\eqref{E:RICLLRENORMALIZED}.
	It is straightforward to see that the term $\mbox{\upshape Inhom}$
	is of the form of RHS~\eqref{E:RENORMALIZEDRICLLINHOMOGENEOUSTERM}:
	we use \eqref{E:UPMUTIMESNULLFORMSSCHEMATIC}
	to decompose the null forms on 
	RHSs \eqref{E:VELOCITYWAVEEQUATION}-\eqref{E:RENORMALIZEDDENSITYWAVEEQUATION},
	Cor.~\ref{C:VELOCITYWAVEEQUATIONDERIVATIVEOFVORTICITYINHOMOGENEOUSTERMEXPRESSION}
	to decompose the product on RHS~\eqref{E:VELOCITYWAVEEQUATION}
	depending on the first Cartesian coordinate partial derivatives of $\Vortrenormalized$,
	\eqref{E:TRANSPORTVECTORFIELDINTERMSOFLUNITANDRADUNIT}
	to decompose the material derivative vectorfield on RHS~\eqref{E:VELOCITYWAVEEQUATION},
	and Lemma~\ref{L:SCHEMATICDEPENDENCEOFMANYTENSORFIELDS}.
	In a detailed proof (see \cite{jSgHjLwW2016}*{Lemma 6.1}), one would find that
	the term $- \frac{1}{2} \mytr \upchi \vec{G}_{\Lunit \Lunit}\contr\Rad \threePsi$
	on RHS~\eqref{E:SKETCHRICLLWAVEEQUATIONTERMNORMALIZED} is canceled
	by another term and hence does not appear on RHS~\eqref{E:RICLLRENORMALIZED}.
	This completes our proof sketch of the lemma.
\end{proof}

\subsection{The definitions of the modified quantities and their transport equations}
\label{SS:MODIFIEDQUANTITIES}

\begin{definition}[\textbf{Modified versions of the derivatives of} $\mytr \upchi$]
\label{D:TRANSPORTRENORMALIZEDTRCHIJUNK}
Let $\Fullset_{\ast}^{N;\leq 1}$ be an $N^{th}$ order commutation vectorfield operator
(see Sect.~\ref{SS:STRINGSOFCOMMUTATIONVECTORFIELDS} regarding the notation).
We define the fully modified function $\upchifullmodarg{\Fullset_{\ast}^{N;\leq 1}}$ as follows:
\begin{subequations}
\begin{align}
	\upchifullmodarg{\Fullset_{\ast}^{N;\leq 1}}
	& := \upmu \Fullset_{\ast}^{N;\leq 1} \mytr \upchi 
			 + \Fullset_{\ast}^{N;\leq 1} \upchifullmodinhom,
		\label{E:TRANSPORTRENORMALIZEDTRCHIJUNK} 
			\\
	\upchifullmodinhom
	& := - \vec{G}_{\Lunit \Lunit}\contr\Rad \threePsi
				- \frac{1}{2} \upmu \mytr \angG\contr\Lunit \threePsi
				- \frac{1}{2} \upmu \vec{G}_{\Lunit \Lunit}\contr\Lunit \threePsi 
				+ \upmu \angGnospacemixedarg{\Lunit}{\#}\contr\cdot \angdiff \threePsi.
			\label{E:LOWESTORDERTRANSPORTRENORMALIZEDTRCHIJUNKDISCREPANCY}
\end{align}
\end{subequations}

We define the partially modified function $\upchipartialmodarg{\Fullset_{\ast}^{N;\leq 1}}$ as follows:
\begin{subequations}
\begin{align}
	\upchipartialmodarg{\Fullset_{\ast}^{N;\leq 1}}
	& := \Fullset_{\ast}^{N;\leq 1} \mytr \upchi 
		+ \upchipartialmodinhomarg{\Fullset_{\ast}^{N;\leq 1}},
		\label{E:TRANSPORTPARTIALRENORMALIZEDTRCHIJUNK} \\
	\upchipartialmodinhomarg{\Fullset_{\ast}^{N;\leq 1}}
	& := - \frac{1}{2} \mytr \angG\contr\Lunit \Fullset_{\ast}^{N;\leq 1} \threePsi
			- \frac{1}{2} \vec{G}_{\Lunit \Lunit}\contr\Lunit \Fullset_{\ast}^{N;\leq 1} \threePsi
			+ \angGmixedarg{\Lunit}{\#}\contr\cdot \angdiff \Fullset_{\ast}^{N;\leq 1} \threePsi.
			\label{E:TRANSPORTPARTIALRENORMALIZEDTRCHIJUNKDISCREPANCY}
\end{align}
\end{subequations}
We also define the following ``$0^{th}$-order'' version of \eqref{E:TRANSPORTPARTIALRENORMALIZEDTRCHIJUNKDISCREPANCY}:
\begin{align} \label{E:LOWESTORDERTRANSPORTPARTIALRENORMALIZEDTRCHIJUNKDISCREPANCY}
	\upchipartialmodinhom
	& := - \frac{1}{2} \mytr \angG\contr\Lunit \threePsi
			- \frac{1}{2} \vec{G}_{\Lunit \Lunit}\contr\Lunit \threePsi
			+ \angGmixedarg{\Lunit}{\#}\contr\cdot \angdiff \threePsi.
\end{align}

\end{definition}

\begin{proposition}\cite{jSgHjLwW2016}*{Proposition 6.2; \textbf{The transport equation for the fully modified version of} $\Fullset_{\ast}^{N;\leq 1} \mytr \upchi$}
\label{P:TOPORDERTRCHIJUNKRENORMALIZEDTRANSPORT}
Assume that the entries of $\threePsi = (\Densrenormalized,v^1,v^2)$
verify the geometric wave equation system
\eqref{E:VELOCITYWAVEEQUATION}-\eqref{E:RENORMALIZEDDENSITYWAVEEQUATION}.
Let $\Fullset_{\ast}^{N;\leq 1}$ be an $N^{th}$ order commutation vectorfield operator
(see Sect.~\ref{SS:STRINGSOFCOMMUTATIONVECTORFIELDS} regarding the notation)
and let $\upchifullmodarg{\Fullset_{\ast}^{N;\leq 1}}$ and 
$\upchifullmodinhom$ be the corresponding quantities 
defined in \eqref{E:TRANSPORTRENORMALIZEDTRCHIJUNK} and \eqref{E:LOWESTORDERTRANSPORTRENORMALIZEDTRCHIJUNKDISCREPANCY}.
Then the fully modified quantity $\upchifullmodarg{\Fullset_{\ast}^{N;\leq 1}}$ 
verifies the following transport equation:
\begin{align} \label{E:TOPORDERTRCHIJUNKRENORMALIZEDTRANSPORT}
\Lunit \upchifullmodarg{\Fullset_{\ast}^{N;\leq 1}}
	- \left(
			2 \frac{\Lunit \upmu}{\upmu}
		\right)
		\upchifullmodarg{\Fullset_{\ast}^{N;\leq 1}}
	& = 
	\upmu [\Lunit, \Fullset_{\ast}^{N;\leq 1}] \mytr \upchi
	- 2 \upmu \mytr \upchi \Fullset_{\ast}^{N;\leq 1} \mytr \upchi
		- 
				\left(
					2 \frac{\Lunit \upmu}{\upmu} 
				\right)
				\Fullset_{\ast}^{N;\leq 1} \upchifullmodinhom
				\\
	& \ \ 
			+ [\Lunit, \Fullset_{\ast}^{N;\leq 1}] \upchifullmodinhom
			+ [\upmu, \Fullset_{\ast}^{N;\leq 1}] \Lunit \mytr \upchi
			+ [\Fullset_{\ast}^{N;\leq 1},\Lunit \upmu] \mytr \upchi
			\notag \\
	& \ \ - \left\lbrace
					\Fullset_{\ast}^{N;\leq 1} \left(\upmu (\mytr \upchi)^2 \right)
					- 2 \upmu \mytr \upchi \Fullset_{\ast}^{N;\leq 1} \mytr \upchi
				\right\rbrace 
			- \Fullset_{\ast}^{N;\leq 1} \mathfrak{A},
			\notag
\end{align}
where the term $\mathfrak{A}$
on the last line of RHS~\eqref{E:TOPORDERTRCHIJUNKRENORMALIZEDTRANSPORT}
is the one appearing in \eqref{E:RICLLRENORMALIZED}-\eqref{E:RENORMALIZEDRICLLINHOMOGENEOUSTERM}.
\end{proposition}

\begin{proposition}\cite{jSgHjLwW2016}*{Proposition 6.3; \textbf{The transport equation for the partially modified 
version of} $\Fullset_{\ast}^{N-1;\leq 1} \mytr \upchi$}
\label{P:COMMUTEDTRCHIJUNKFIRSTPARTIALRENORMALIZEDTRANSPORTEQUATION}
Let $\Fullset_{\ast}^{N-1;\leq 1}$ be an $(N-1)^{st}$ order commutation vectorfield operator
(see Sect.~\ref{SS:STRINGSOFCOMMUTATIONVECTORFIELDS} regarding the notation)
and let $\upchipartialmodarg{\Fullset_{\ast}^{N-1;\leq 1}}$ 
be the corresponding partially modified quantity defined in \eqref{E:TRANSPORTPARTIALRENORMALIZEDTRCHIJUNK}.
Then $\upchipartialmodarg{\Fullset_{\ast}^{N-1;\leq 1}}$
verifies the following transport equation:
\begin{align} \label{E:COMMUTEDTRCHIJUNKFIRSTPARTIALRENORMALIZEDTRANSPORTEQUATION}
	\Lunit \upchipartialmodarg{\Fullset_{\ast}^{N-1;\leq 1}}
	& = \frac{1}{2} \vec{G}_{\Lunit \Lunit}\contr\angLap \Fullset_{\ast}^{N-1;\leq 1} \threePsi
			+ 
			\upchipartialmodsourcearg{\Fullset_{\ast}^{N-1;\leq 1}},
\end{align}
where the inhomogeneous term $\upchipartialmodsourcearg{\Fullset_{\ast}^{N-1;\leq 1}}$ is given by
\begin{align}  \label{E:TRCHIJUNKCOMMUTEDTRANSPORTEQNPARTIALRENORMALIZATIONINHOMOGENEOUSTERM}
\upchipartialmodsourcearg{\Fullset_{\ast}^{N-1;\leq 1}}
	& = - \Fullset_{\ast}^{N-1;\leq 1} \mathfrak{B}
			- \Fullset_{\ast}^{N-1;\leq 1} (\mytr \upchi)^2
				\\
	& \ \ 
			+ \frac{1}{2} [\Fullset_{\ast}^{N-1;\leq 1}, \vec{G}_{\Lunit \Lunit}]\contr\angLap \threePsi
			+ \frac{1}{2} \vec{G}_{\Lunit \Lunit} [\Fullset_{\ast}^{N-1;\leq 1}, \angLap]\contr\threePsi
			+ [\Lunit, \Fullset_{\ast}^{N-1;\leq 1}] \mytr \upchi
			  \notag \\
	& \ \ 
				+ 
				[\Lunit, \Fullset_{\ast}^{N-1;\leq 1}] \upchipartialmodinhom
				+ 
				\Lunit
				\left\lbrace
					\upchipartialmodinhomarg{\Fullset_{\ast}^{N-1;\leq 1}}
					- \Fullset_{\ast}^{N-1;\leq 1} \upchipartialmodinhom
				\right\rbrace,
				\notag
\end{align}
$\mathfrak{B}$ is defined in \eqref{E:PARTIALRENORMALIZEDRICLLNHOMOGENEOUSTERM},
$\upchipartialmodinhomarg{\Fullset_{\ast}^{N-1;\leq 1}}$ 
is defined in \eqref{E:TRANSPORTPARTIALRENORMALIZEDTRCHIJUNKDISCREPANCY},
and $\upchipartialmodinhom$
is defined in \eqref{E:LOWESTORDERTRANSPORTPARTIALRENORMALIZEDTRCHIJUNKDISCREPANCY}.
\end{proposition}

\subsection{Some identities connected to curvature}
\label{SS:IDENTITIESCONNECTEDTOCURVATURE}
We now show that
$\Rad \mytr \upchi$ and
$\angLap \upmu$
are equal up to simple error terms.
This fact allows for a simplified approach to various estimates
appearing later in the paper.

\begin{lemma}\cite{jSgHjLwW2016}*{Lemma 11.4; \textbf{Connection between} $\Rad \mytr \upchi$ \textbf{and} 
$\angLap \upmu$}
\label{L:RADDIRECTIONCHITRANSPORT}
$\Rad \mytr \upchi$ can be expressed as follows,
where the term $\angLap \upmu$ on
RHS~\eqref{E:RADDIRECTIONCHITRANSPORT} and $\smoothfunction(\cdots)$ is schematic:
\begin{align}  \label{E:RADDIRECTIONCHITRANSPORT}
	\Rad \mytr \upchi
	& = \angLap \upmu
			+ \smoothfunction(\GdVar,\ginversesphere,\angdiff \vec{x}) \Singletan \Rad \threePsi
			+ \smoothfunction(\BadVar,\ginversesphere, \angdiff \vec{x}) \Singletan \Singletan \threePsi
		\\
	& \ \ + 
		\smoothfunction(\BadVar,\Rad \threePsi, \Singletan \BadVar,\ginversesphere,\angdiff x) 
		\Singletan \GdVar.
		\notag
\end{align}

\end{lemma}

\begin{proof}[Discussion of proof]
	Lemma~\ref{L:RADDIRECTIONCHITRANSPORT} was essentially proved as
	\cite{jSgHjLwW2016}*{Lemma 11.4} and is based on an analysis
	of the Riemann curvature component $\mytr \Cur_{\Rad \cdot \Lunit \cdot}$.
	We remark that in the identity provided by \cite{jSgHjLwW2016}*{Lemma 11.4}, 
	one finds a term proportional to $\angD^2 \vec{x}$.
However, using \eqref{E:ANGDSQUAREDINTERMSOFGEOANGDERIVATIVES} with $f = x^i$
and Lemma~\ref{L:SCHEMATICDEPENDENCEOFMANYTENSORFIELDS},
we can write 
$
\angD^2 \vec{x} 
= \smoothfunction(\GdVar,\angdiff \vec{x}) \Singletan \GdVar
$,
and thus the corresponding error terms are part of the last
term on RHS~\eqref{E:RADDIRECTIONCHITRANSPORT}.
This completes our discussion of the lemma.
We remark that similar calculations are presented in \cite{dC2007}*{Chapter 4}.
\end{proof}

\section{Assumptions on the initial state of the solution and bootstrap assumptions}
\label{S:DATASIZEANDBOOTSTRAP}
In this section, we introduce our Sobolev norm assumptions on the data
for $\threePsi$,
$\Vortrenormalized$,
and the eikonal function quantities.
We also state the bootstrap assumptions that we use 
in analyzing solutions.
By data, we mean the state of the solution along
$\Sigma_0^1$ and a large potion of the outgoing null
hypersurface $\mathcal{P}_0$.
Our assumptions involve several size parameters,
and in Sect.~\ref{SS:SMALLNESSASSUMPTIONS}, we describe 
our assumptions on their relative sizes.
In Subsubsect.~\ref{SS:EXISTENCEOFDATA},
we show that there exists an open set of
nearly plane symmetric data verifying the 
size assumptions.

\subsection{Assumptions on the initial state of the fluid variables}
\label{SS:FLUIDVARIABLEDATAASSUMPTIONS}

\subsubsection{The quantity that controls the blowup-time}
\label{SSS:CRUCIALDELTADEF}
We start by introducing the data-dependent number 
$\TranminusdatasizeWithFactor$, 
which is of crucial importance.
Our main theorem shows that 
if $\mathring{\upepsilon}$ (defined just below)
is sufficiently small,
then the time of first shock formation is 
$(1 + \mathcal{O}(\mathring{\upepsilon}))\TranminusdatasizeWithFactor^{-1}$

\begin{definition}[\textbf{The quantity that controls the blowup-time}]
	\label{D:CRITICALBLOWUPTIMEFACTOR}
	We define
	\begin{align} \label{E:CRITICALBLOWUPTIMEFACTOR}
		\TranminusdatasizeWithFactor
		& := \frac{1}{2} 
		\sup_{\Sigma_0^1} 
		\left[
			\sum_{\imath=0}^1 G_{\Lunit \Lunit}^{\imath} \Rad v^1 \right]_-.
	\end{align}
\end{definition}

\begin{remark}[\textbf{Significance of} $\TranminusdatasizeWithFactor$]
	Equation \eqref{E:UPMUFIRSTTRANSPORT}
	and the estimates of Props.~\ref{P:IMPROVEMENTOFAUX}
	and \ref{P:IMPROVEMENTOFHIGHERTRANSVERSALBOOTSTRAP} can be used to show that
	there exist $u_* \in [0,U_0]$ and $\vartheta_* \in \mathbb{T}$
	such that for $t \geq 1$, we have
	$\Lunit \upmu(t,u_*,\vartheta_*) = - \TranminusdatasizeWithFactor + \mbox{\upshape Error}$,
	where (under suitable assumptions on the data) 
	$\mbox{\upshape Error}$ is small compared to $\TranminusdatasizeWithFactor$.
	That is, the maximal shrinking rate of $\upmu$ along the integral curves of $\Lunit$
	is determined by $\TranminusdatasizeWithFactor$.
	It is for this reason that $\TranminusdatasizeWithFactor^{-1}$
	is connected to the time of shock formation.
\end{remark}

\subsubsection{Size assumptions for the fluid variables}
\label{SSS:SIZEASSUMPTIONSFORFLUIDVARIABLES}
We make the following size assumptions along $\Sigma_0^1$ and 
$\mathcal{P}_0^{2 \TranminusdatasizeWithFactor^{-1}}$
(see Sect.~\ref{SS:STRINGSOFCOMMUTATIONVECTORFIELDS} regarding the vectorfield operator notation).

\medskip

\noindent \underline{$L^2$ \textbf{assumptions along} $\Sigma_0^1$}.
\begin{align} \label{E:L2SMALLDATAASSUMPTIONSALONGSIGMA0}
	\left\|
		\Fullset_{\ast}^{\leq 21;\leq 2} \threePsi
	\right\|_{L^2(\Sigma_0^1)},
			\,
	\myarray
	[\left\|
			\Rad (\Densrenormalized - v^1)
		\right\|_{L^2(\Sigma_0^1)}]
	{
		\left\|
			\Rad^{[0,2]} v^2
		\right\|_{L^2(\Sigma_0^1)}
		},
			\,
		\left\|
			\Tanset^{\leq 21} \Vortrenormalized
		\right\|_{L^2(\Sigma_0^1)}
		& \leq \mathring{\upepsilon}.
\end{align}

\noindent \underline{$L^{\infty}$ \textbf{assumptions along} $\Sigma_0^1$}.
\begin{subequations}
\begin{align} \label{E:LINFTYSMALLDATAASSUMPTIONSALONGSIGMA0}
	\left\|
		\Fullset_{\ast}^{\leq 13;\leq 1} \threePsi
	\right\|_{L^{\infty}(\Sigma_0^1)},
			\,
	\left\|
		\Fullset_{\ast}^{\leq 12;\leq 2} \threePsi
	\right\|_{L^{\infty}(\Sigma_0^1)},
		\\
	\myarray[
		\left\|
			\Rad (\Densrenormalized  - v^1)
		\right\|_{L^{\infty}(\Sigma_0^1)}
		]
		{
		\left\|
			\Rad^{[1,3]} v^2
		\right\|_{L^{\infty}(\Sigma_0^1)}
		},
		\,
		\left\|
		\Lunit \Rad \Rad \Rad \threePsi
	\right\|_{L^{\infty}(\Sigma_0^1)},
		\notag \\
		\left\|
			\Tanset^{\leq 13}\Vortrenormalized
		\right\|_{L^{\infty}(\Sigma_0^1)}
		& \leq \mathring{\upepsilon},
		\notag \\
	\myarray[
		\left\|
			\Rad^{[1,3]} \Densrenormalized
		\right\|_{L^{\infty}(\Sigma_0^1)}
		]
	{
	\left\|
		\Rad^{[1,3]} v^1
	\right\|_{L^{\infty}(\Sigma_0^1)}
	}
	& \leq \mathring{\updelta}.
	\label{E:V1LARGEDATAASSUMPTIONSALONGSIGMA0}
\end{align}
\end{subequations}

\noindent \underline{$L^2$ \textbf{assumptions along} $\mathcal{P}_0^{2 \TranminusdatasizeWithFactor^{-1}}$}.
\begin{align} \label{E:PSIL2SMALLDATAASSUMPTIONSALONGP0}
		\left\|
			\Fullset^{\leq 21;\leq 1} \threePsi
		\right\|_{L^2\left(\mathcal{P}_0^{2 \TranminusdatasizeWithFactor^{-1}}\right)},
			\,
		\left\|
			\Tanset^{\leq 21} \Vortrenormalized
		\right\|_{L^2\left(\mathcal{P}_0^{2 \TranminusdatasizeWithFactor^{-1}}\right)}
		& \leq \mathring{\upepsilon}.
\end{align}

\noindent \underline{$L^{\infty}$ \textbf{assumptions along} 
$\mathcal{P}_0^{2 \TranminusdatasizeWithFactor^{-1}}$}.
\begin{align} \label{E:PSILINFTYSMALLDATAASSUMPTIONSALONGP0}
		\left\|
			\Fullset^{\leq 19;\leq 1} \threePsi
		\right\|_{L^{\infty}\left(\mathcal{P}_0^{2 \TranminusdatasizeWithFactor^{-1}}\right)},
			\,
		\left\|
			\Tanset^{\leq 19} \Vortrenormalized
		\right\|_{L^{\infty}\left(\mathcal{P}_0^{2 \TranminusdatasizeWithFactor^{-1}}\right)}
		& \leq \mathring{\upepsilon}.
\end{align}

\noindent \underline{$L^2$ \textbf{assumptions along} $\ell_{1,u}$}.
\begin{align} \label{E:DATAASSUMPTIONSALONGELL1U}
	\left\|
			\Fullset_{\ast}^{\leq 20;\leq 1} \threePsi
		\right\|_{L^2(\ell_{1,u})},
			\,
	\myarray[
		\left\|
			\Rad (\Densrenormalized  - v^1)
		\right\|_{L^2(\ell_{1,u})}
		]
		{\left\|
			\Rad v^2 
		\right\|_{L^2(\ell_{1,u})}},
			\,
		\left\|
			\Tanset^{\leq 20} \Vortrenormalized
		\right\|_{L^2(\ell_{1,u})}
		& \leq \mathring{\upepsilon}.
\end{align}

\noindent \underline{$L^2$ \textbf{assumptions along} $\ell_{t,0}$}.
\begin{align} \label{E:SMALLDATAASSUMPTIONSALONGELLT0}
		\left\|
			\Fullset^{\leq 20;\leq 1} \threePsi
		\right\|_{L^2(\ell_{t,0})},
			\,
		\left\|
			\Tanset^{\leq 20} \Vortrenormalized
		\right\|_{L^2(\ell_{t,0})}
		& \leq \mathring{\upepsilon}.
\end{align}

\begin{remark}[\textbf{A concise summary of the effect of the size assumptions}]
	\label{R:SUMMARYOFDATASIZEASSUMPTIONS}
	The assumptions 
	\eqref{E:L2SMALLDATAASSUMPTIONSALONGSIGMA0}-\eqref{E:SMALLDATAASSUMPTIONSALONGELLT0}
	will allow us to prove that 
	among $\threePsi$, $\Vortrenormalized$ and their relevant derivatives,
	the only relatively large (in all relevant norms) quantities 
	in our analysis are $\Rad^{[1,3]} v^1$ 
	and $\Rad^{[1,3]} \Densrenormalized$
	along $\Sigma_t^u$.
	Moreover, even 
	$\Rad(\Densrenormalized  - v^1)$ is small along
	$\Sigma_t^u$, 
	and
	$\Rad v^1$ 
	and $\Rad \Densrenormalized$
	are small along 
	$\mathcal{P}_0^{2 \TranminusdatasizeWithFactor^{-1}}$.
	This division into small and large quantities 
	is fundamental for our analysis.
\end{remark}

To prove our main theorem,
we make assumptions on the
relative sizes of the above
parameters; see Sect.~\ref{SS:SMALLNESSASSUMPTIONS}.

\subsection{Assumptions on the initial conditions of the eikonal function quantities}
\label{SS:DATAFOREIKONALFUNCTIONQUANTITIES}
We now state our size assumptions for the initial conditions
of the eikonal function quantities
(see Sect.~\ref{SS:STRINGSOFCOMMUTATIONVECTORFIELDS} regarding the vectorfield operator notation).

\medskip

\noindent \underline{$L^2$ \textbf{assumptions along} $\Sigma_0^1$}.
We assume that there exist (implicit) constants, depending on 
$\mathring{\updelta}$, such that
\begin{subequations}
\begin{align} \label{E:LUNITISMALLDATA}
	\left\|
		\Fullset_{\ast}^{\leq 21;\leq 2} \Lunit_{(Small)}^i
	\right\|_{L^2(\Sigma_0^1)}
	& \lesssim \mathring{\upepsilon},
		\\
	\left\|
		\Rad^{[1,3]} \Lunit_{(Small)}^i
	\right\|_{L^2(\Sigma_0^1)}
	& \lesssim 1,
		\label{E:LUNITILARGEDATA}
\end{align}
\end{subequations}

\begin{subequations}
\begin{align}  \label{E:UPMUDATATANGENTIALL2CONSEQUENCES}
	\left\|
		\Fullset_{\ast \ast}^{[1,21];\leq 1} \upmu
	\right\|_{L^2(\Sigma_0^1)}
	& \lesssim \mathring{\upepsilon},
		\\
	\left\|
		\Lunit \Rad^{[0,2]} \upmu 
	\right\|_{L^2(\Sigma_0^1)},
		\,
	\left\|
		\Rad^{[0,2]} \Lunit \upmu 
	\right\|_{L^2(\Sigma_0^1)},
		\,
	\left\|
		\Rad \Lunit \Rad \upmu 
	\right\|_{L^2(\Sigma_0^1)},
		\,
	\left\|
		\Rad^{[0,2]} \upmu 
	\right\|_{L^2(\Sigma_0^1)}
	& \lesssim 1.
	\label{E:UPMUDATARADIALL2CONSEQUENCES}
\end{align}
\end{subequations}

\noindent \underline{$L^{\infty}$ \textbf{assumptions along} $\Sigma_0^1$}.
\begin{align} 
	\label{E:LUNITISMALLDATALINFTY}
	\left\|
		\Fullset_{\ast}^{\leq 11;2} \Lunit_{(Small)}^i
	\right\|_{L^{\infty}(\Sigma_0^1)}
	& \lesssim \mathring{\upepsilon},
		\\
	\left\|
		\Rad^{[0,2]} \Lunit_{(Small)}^i
	\right\|_{L^{\infty}(\Sigma_0^1)}
	& \lesssim 1,
		\label{E:LUNITILARGEDATALINFTY}
\end{align}

\begin{subequations}
\begin{align}  \label{E:UPMUDATATANGENTIALLINFINITYCONSEQUENCES}
	\left\|
		\upmu - 1
	\right\|_{L^{\infty}(\Sigma_0^1)},
		\,
	\left\|
		\Fullset_{\ast \ast}^{[1,11];\leq 1} \upmu
	\right\|_{L^{\infty}(\Sigma_0^1)}
	& \lesssim \mathring{\upepsilon},
		\\
	\left\|
		\Lunit \Rad^{[0,2]} \upmu 
	\right\|_{L^{\infty}(\Sigma_0^1)},
		\,
	\left\|
		\Rad^{[0,2]} \Lunit \upmu 
	\right\|_{L^{\infty}(\Sigma_0^1)},
		\,
	\left\|
		\Rad \Lunit \Rad \upmu 
	\right\|_{L^{\infty}(\Sigma_0^1)},
		\,
	\left\|
		\Rad^{[0,2]} \upmu 
	\right\|_{L^{\infty}(\Sigma_0^1)}
	& \lesssim 1,
	\label{E:UPMUDATARADIALLINFINITYCONSEQUENCES}
\end{align}
\end{subequations}

\begin{align} \label{E:LINFTYCOORADANGILARGEDATAALONGSIGMA1}
	\left\|
		\Fullset^{\leq 18;\leq 3} (\CoordAng^i - \delta_2^i)
	\right\|_{L^{\infty}(\Sigma_0^1)}
	& \lesssim \mathring{\upepsilon},
		\\
	\left\|
		\Fullset^{\leq 18;\leq 2} \NonRadialRad^i
	\right\|_{L^{\infty}(\Sigma_0^1)}
	& \lesssim \mathring{\upepsilon}.
	\label{E:LINFTYXILARGEDATAALONGSIGMA1}
\end{align}

\noindent \underline{$L^{\infty}$ \textbf{assumptions along} 
$\mathcal{P}_0^{2 \TranminusdatasizeWithFactor^{-1}}$}.

\begin{align} \label{E:SIMPLEUPMUALONGPOESTIMATE}
	\left\|
		\upmu - 1
	\right\|_{L^{\infty}(\mathcal{P}_0^{2 \TranminusdatasizeWithFactor^{-1}})}
	& \lesssim \mathring{\upepsilon}.
\end{align}

\subsection{\texorpdfstring{$\Tboot$}{The bootstrap time}, the positivity of 
\texorpdfstring{$\upmu$}{the inverse foliation density}, and the diffeomorphism property of 
\texorpdfstring{$\Upsilon$}{the change of variables map}}
\label{SS:SIZEOFTBOOT}
We now state some basic bootstrap assumptions.
We start by fixing a real number $\Tboot$ with
\begin{align} \label{E:TBOOTBOUNDS}
	0 <  \Tboot \leq 2 \TranminusdatasizeWithFactor^{-1}.
\end{align}

We assume that on the spacetime domain $\mathcal{M}_{\Tboot,U_0}$
(see \eqref{E:MTUDEF}), we have
\begin{align} \label{E:BOOTSTRAPMUPOSITIVITY} \tag{$\mathbf{BA} \upmu > 0$}
	\upmu > 0.
\end{align}
Inequality \eqref{E:BOOTSTRAPMUPOSITIVITY} implies that no shocks are present in
$\mathcal{M}_{\Tboot,U_0}$.

We also assume that
\begin{align} \label{E:BOOTSTRAPCHOVISDIFFEO}
	& \mbox{The change of variables map $\Upsilon$ from Def.~\ref{D:CHOVMAP}
	is a $C^1$ diffeomorphism from} \\
	& [0,\Tboot) \times [0,U_0] \times \mathbb{T}
	\mbox{ onto its image}.
	\notag
\end{align}

\subsection{Fundamental \texorpdfstring{$L^{\infty}$}{essential sup-norm} bootstrap assumptions}
\label{SS:PSIBOOTSTRAP}
 Our fundamental bootstrap assumptions for 
$\threePsi$ and $\Vortrenormalized$ are that the following inequalities hold on $\mathcal{M}_{\Tboot,U_0}$
 (see Sect.~\ref{SS:STRINGSOFCOMMUTATIONVECTORFIELDS} regarding the vectorfield operator notation):
\begin{align} \label{E:PSIFUNDAMENTALC0BOUNDBOOTSTRAP} \tag{$\mathbf{BA}\threePsi$}
	\left\| 
		\Fullset_{\ast}^{\leq 13;\leq 1} \threePsi
	\right\|_{L^{\infty}(\Sigma_t^u)},
		\,
	\left\| 
		\Fullset^{\leq 13;\leq 1} (\Densrenormalized  - v^1)
	\right\|_{L^{\infty}(\Sigma_t^u)}
	& \leq \varepsilon,
		\\
	\left\| 
		\Tanset^{\leq 13} \Vortrenormalized
	\right\|_{L^{\infty}(\Sigma_t^u)}
	& \leq \varepsilon,
	\label{E:VORTFUNDAMENTALC0BOUNDBOOTSTRAP} \tag{$\mathbf{BA}\Vortrenormalized$}
\end{align}
where $\varepsilon$ is a small positive bootstrap parameter whose smallness 
we describe in Sect.~\ref{SS:SMALLNESSASSUMPTIONS}.

\subsection{Auxiliary  \texorpdfstring{$L^{\infty}$}{essential sup-norm} bootstrap assumptions}
\label{SS:AUXILIARYBOOTSTRAP}
In deriving pointwise estimates, we find it convenient to make the
following auxiliary bootstrap assumptions.
In Prop.~\ref{P:IMPROVEMENTOFAUX}, we will derive strict improvements
of these assumptions.

\medskip

\noindent \underline{\textbf{Auxiliary bootstrap assumptions for small quantities}}.
We assume that the following inequalities hold on $\mathcal{M}_{\Tboot,U_0}$:
\begin{align} \label{E:PSIAUXLINFINITYBOOTSTRAP} \tag{$\mathbf{AUX} \threePsi \mathbf{SMALL}$}
	\left\| 
		\Fullset_{\ast}^{\leq 12;\leq 2} \threePsi 
	\right\|_{L^{\infty}(\Sigma_t^u)}
	& \leq \varepsilon^{1/2},
\end{align}

\begin{align}
	\left\| 
		\Fullset_{\ast}^{\leq 11;\leq 2} \Lunit_{(Small)}^i 
	\right\|_{L^{\infty}(\Sigma_t^u)}
	& \leq \varepsilon^{1/2},
		\label{E:FRAMECOMPONENTSBOOT} \tag{$\mathbf{AUX} \Lunit_{(Small)} \mathbf{SMALL}$} 
\end{align}

\begin{align}
	\left\| 
		\Fullset_{\ast \ast}^{[1,11];1} \upmu
	\right\|_{L^{\infty}(\Sigma_t^u)}
	& \leq \varepsilon^{1/2},
		\label{E:UPMUBOOT}  \tag{$\mathbf{AUX} \upmu \mathbf{SMALL}$} 
\end{align}

\begin{align}
	\left\| 
		\angLie_{\Fullset}^{\leq 11;\leq 1} \upchi
	\right\|_{L^{\infty}(\Sigma_t^u)},
		\,
	\left\| 
		\angLie_{\Fullset}^{\leq 10;\leq 2} \upchi
	\right\|_{L^{\infty}(\Sigma_t^u)}
	& \leq \varepsilon^{1/2}.
	\label{E:CHIBOOT} \tag{$\mathbf{AUX} \upchi$}
\end{align}

\noindent \underline{\textbf{Auxiliary bootstrap assumptions for quantities that are allowed to be large}}.
We assume that the following inequalities hold on $\mathcal{M}_{\Tboot,U_0}$ for $M=1,2$:
\begin{align}
	\left\| 
		\Rad^M \Densrenormalized 
	\right\|_{L^{\infty}(\Sigma_t^u)}
	& \leq 
	\left\| 
		\Rad^M \Densrenormalized 
	\right\|_{L^{\infty}(\Sigma_0^u)}
	+ \varepsilon^{1/2},
	&& 
		\label{E:DENSITYTRANSVERSALLINFINITYBOUNDBOOTSTRAP} \tag{$\mathbf{AUX} \Densrenormalized \mathbf{LARGE}$}
			\\
	\left\| 
		\Rad^M v^1 
	\right\|_{L^{\infty}(\Sigma_t^u)}
	& \leq 
	\left\| 
		\Rad^M v^1
	\right\|_{L^{\infty}(\Sigma_0^u)}
	+ \varepsilon^{1/2}.
	&& 
		\label{E:V1TRANSVERSALLINFINITYBOUNDBOOTSTRAP} \tag{$\mathbf{AUX} v^1 \mathbf{LARGE}$}
\end{align}

We assume that the following inequalities hold on $\mathcal{M}_{\Tboot,U_0}$ for $M=0,1$:
\begin{align}
	\left\| 
		\Lunit \Rad^M \upmu
	\right\|_{L^{\infty}(\Sigma_t^u)}
	& \leq 
		\frac{1}{2}
		\left\| 
			\Rad^M
			\left\lbrace
				\vec{G}_{\Lunit \Lunit} \contr \Rad \threePsi
			\right\rbrace
		\right\|_{L^{\infty}(\Sigma_0^u)}
		+ \varepsilon^{1/2},
		 \label{E:LUNITUPMUBOOT}  \tag{$\mathbf{AUX} \Lunit \upmu$}  \\
		\left\| 
			\Rad^M \upmu
		\right\|_{L^{\infty}(\Sigma_t^u)}
		& \leq
	 	\left\| 
			\Rad^M \upmu
		\right\|_{L^{\infty}(\Sigma_0^u)}
		+ 
		2 \TranminusdatasizeWithFactor^{-1} 
		\left\| 
			\Rad^M
			\left\lbrace
				\vec{G}_{\Lunit \Lunit} \contr \Rad \threePsi
			\right\rbrace
		\right\|_{L^{\infty}(\Sigma_0^u)}
		+ \varepsilon^{1/2}.
			\label{E:UPMUTRANSVERSALBOOT}  \tag{$\mathbf{AUX} \upmu$} 
	\end{align}

	We assume that the following inequalities hold on $\mathcal{M}_{\Tboot,U_0}$ for $M=1,2$:
	\begin{align}
		\left\| 
			\Rad^M \Lunit_{(Small)}^i
		\right\|_{L^{\infty}(\Sigma_t^u)}
		& \leq
	 	\left\| 
			\Rad^M \Lunit_{(Small)}^i
		\right\|_{L^{\infty}(\Sigma_0^u)}
		+ \varepsilon^{1/2}.
		\label{E:LUNITITRANSVERSALBOOT}  \tag{$\mathbf{AUX} \Lunit_{(Small)} \mathbf{LARGE}$} 
\end{align}

\subsection{Smallness assumptions}
\label{SS:SMALLNESSASSUMPTIONS}
For the remainder of the article, 
when we say that ``$A$ is small relative to $B$,''
we mean that there exists a continuous increasing function 
$f :[0,\infty) \rightarrow (0,\infty)$ 
such that 
$
\displaystyle
A \leq f(B)
$.
In principle, the functions $f$ could always be chosen to be 
polynomials with positive coefficients or exponential functions.\footnote{The exponential functions appear, for example, in our energy estimates, 
during our Gronwall argument; see the proof of Prop.~\ref{P:MAINAPRIORIENERGY}
given in Sect.~\ref{SS:PROOFOFPROPMAINAPRIORIENERGY}.} 
However, to avoid lengthening the paper, we typically do not 
specify the form of $f$.

Throughout the rest of the paper, we make the following
relative smallness assumptions. We
continually adjust the required smallness
in order to close our estimates.
\begin{itemize}
	\item $\varepsilon$ is small relative to $\mathring{\updelta}^{-1}$,
		where $\mathring{\updelta}$ is the data-size parameter 
		from \eqref{E:V1LARGEDATAASSUMPTIONSALONGSIGMA0}.
	\item $\varepsilon$ is small relative to 
		the data-size parameter $\TranminusdatasizeWithFactor$ 
		from \eqref{E:CRITICALBLOWUPTIMEFACTOR}.
\end{itemize}
The first assumption will allow us to control error terms that,
roughly speaking, are of size $\varepsilon \mathring{\updelta}^k$ 
for some integer $k \geq 0$. The second assumption 
is relevant because the expected blowup-time is 
approximately $\TranminusdatasizeWithFactor^{-1}$,
and the assumption will allow us to show that various
error products featuring a small factor $\varepsilon$
remain small for $t < 2 \TranminusdatasizeWithFactor^{-1}$, 
which is plenty of time for us to show that a shock forms.

Finally, we assume that
\begin{align} \label{E:DATAEPSILONVSBOOTSTRAPEPSILON}
	\varepsilon^{3/2}
	& \leq
	\mathring{\upepsilon} 
	\leq \varepsilon,
\end{align}
where $\mathring{\upepsilon}$ is the data smallness parameter from 
Sects.~\ref{SSS:SIZEASSUMPTIONSFORFLUIDVARIABLES} and \ref{SS:DATAFOREIKONALFUNCTIONQUANTITIES}.

\begin{remark}[\textbf{Relationship between} $\varepsilon$ \textbf{and} $\mathring{\upepsilon}$
\textbf{in the proof of our main theorem}]
\label{R:CHOICEOFTWOEPSILONS}
	In the proof of our main theorem,
	we will set $\varepsilon = \BigConst \mathring{\upepsilon}$,
	where $\BigConst > 1$ is chosen to be sufficiently large
	and $\mathring{\upepsilon}$ is assumed to be sufficiently small.
	This is compatible with \eqref{E:DATAEPSILONVSBOOTSTRAPEPSILON}.
\end{remark}

\subsection{The existence of initial data verifying the size assumptions}
\label{SS:EXISTENCEOFDATA}
In this section, we show that there exists an open set of data
verifying the size the assumptions of
Subsects.~\ref{SS:FLUIDVARIABLEDATAASSUMPTIONS},\ref{SS:DATAFOREIKONALFUNCTIONQUANTITIES},
and \ref{SS:SMALLNESSASSUMPTIONS}.
By ``open,'' we mean open relative to the Sobolev topologies
corresponding to the size assumptions stated in those subsections.
By Cauchy stability,\footnote{Here we mean continuous dependence of the solution on the data.}
it is enough to exhibit smooth plane symmetric data that are compactly supported in 
$\Sigma_0^1$ (which can be identified here with the unit $x^1$ interval $[0,1]$)
and that verify the size assumptions. By a plane symmetric solution, 
we mean that $\Densrenormalized = \Densrenormalized(t,x^1)$, $v^1 = v^1(t,x^1)$, and $v^2 \equiv 0$.
The data that we exhibit launch simple plane symmetric solutions. By ``simple,''
we mean that one Riemann invariant completely vanishes.

\begin{remark}[\textbf{Strictly non-zero vorticity along 
$\Sigma_0^1$ and $\mathcal{P}_0^{2 \TranminusdatasizeWithFactor^{-1}}$}]
\label{R:DATAWITHNONZEROVORTICITY}
Once we have exhibited the plane symmetric data described above,
it is easy to perturb it so that
the vorticity is everywhere non-zero along  
$\Sigma_0^1$ and $\mathcal{P}_0^{2 \TranminusdatasizeWithFactor^{-1}}$.
One can simply leave the data for $\Densrenormalized$ and $v^1$
along $\Sigma_0$ unchanged and 
set $v^2|_{\Sigma_0} = f(x^1)$,
where $\uplambda > 0$ is small
and $f$ is smooth with $f' > 0$ in an interval $I$ 
of length $|I| \gg 2 \TranminusdatasizeWithFactor^{-1}$
containing the origin.
Then $\Vortrenormalized$ will be small but
non-zero on $I \times \mathbb{T} \subset \Sigma_0$.
Hence, using the transport equation \eqref{E:RENORMALIZEDVORTICTITYTRANSPORTEQUATION},
it is easy to show\footnote{In the solution regime under consideration,
equation \eqref{E:RENORMALIZEDVORTICTITYTRANSPORTEQUATION}
reads $\partial_t \Vortrenormalized =$ Quadratically small error terms.} 
(under suitable smallness assumptions)
that the corresponding solution
``induces'' data for $\Vortrenormalized$
along $\mathcal{P}_0^{2 \TranminusdatasizeWithFactor^{-1}}$
such that $\Vortrenormalized$ is everywhere non-zero along
$\mathcal{P}_0^{2 \TranminusdatasizeWithFactor^{-1}}$.
\end{remark}

We remind the reader 
(see \eqref{E:INTROEIKONALINITIALVALUE})
that the initial condition for the eikonal function is
$u|_{\Sigma_0} := 1 - x^1$.
Our discussion relies on the following simple lemma.

\begin{lemma}	\cite{jSgHjLwW2016}*{Lemma 7.2; \textbf{Algebraic identities along} $\Sigma_0$}
	\label{L:ALGEBRAICIDALONGSIGMA0}
Consider a solution launched by data given along 
$\Sigma_0$. Assume that the datum for the eikonal function is
$u|_{\Sigma_0} := 1 - x^1$.
Then the following identities hold along $\Sigma_0$ (for $i=1,2$):
\begin{align}
	\upmu
	& = 
	\frac{1}{\Speed},
	\qquad
	\Lunit_{(Small)}^i
	 = 	(\Speed - 1)\delta^{i1}
			+ v^i,
	\qquad
	\NonRadialRad^i
	 =  0,
	 \qquad
	\CoordAng^i
	= \delta_2^i,
		\label{E:INITIALRELATIONS}
\end{align}
where
$\NonRadialRad$ is the $\ell_{t,u}$-tangent vectorfield from \eqref{E:RADSPLITINTOPARTTILAUANDXI}.
\end{lemma}

We now turn to the construction of plane symmetric 
initial data that lead to the desired size assumptions.
Our approach is based on Riemann's method of
Riemann invariants \cite{bR1860}. 
The results that we present here are standard.
Hence, for brevity, we do not provide detailed proofs.
In plane symmetry, in terms of the Riemann invariants
\begin{align} \label{E:RIEMANNINVARIANTS}
	\mathcal{R}_{\pm}
	& = v^1 \pm F(\Densrenormalized),
\end{align}
the compressible Euler equations
\eqref{E:TRANSPORTDENSRENORMALIZEDRELATIVETORECTANGULAR}-\eqref{E:TRANSPORTVELOCITYRELATIVETORECTANGULAR}
are equivalent to the system
\begin{align} \label{E:EVOLUTIONRIEMANNINVARIANT}
	\uLgood \mathcal{R}_-
	& = 0,
	&&
	\Lunit \mathcal{R}_+
	= 0,
\end{align}
where
\begin{align} \label{E:ULGOODAPPENDIXDEF}
	\uLgood 
		& := \upmu \uLunit,
		&&		\\
		\Lunit 
		& = \partial_t + (v^1 + \Speed) \partial_1,
		&&
		\uLunit = \partial_t + (v^1 - \Speed) \partial_1,
		\label{E:CHARVECPLANESYMMETRY}
\end{align}
and $\Lunit$ coincides with the vectorfield defined in Def.~\ref{D:LUNITDEF}.
The function $F$ in \eqref{E:EVOLUTIONRIEMANNINVARIANT}
solves the following initial value problem in $\Densrenormalized$:
\begin{align} \label{E:DENSITYFUNCTIONUSEDFORRIEMANNINVARIANTS}
	\frac{d}{d\Densrenormalized} F(\Densrenormalized) 
	&= \Speed(\Densrenormalized),
	&&
	F(\Densrenormalized=0) = 0,
\end{align}
where $F(\Densrenormalized=0) = 0$ is just a convenient normalization condition.
It is straightforward to show that
\begin{align} \label{E:ULUNITLUNITRADRELATION}
		\uLunit
		& = \Lunit + 2 \Radunit,
		\\
		\Radunit 
		& = - \Speed \partial_1,
		\label{E:PLANESYMMETRY}
	\end{align}
where $\Radunit$ coincides with the vectorfield defined in Def.~\ref{D:RADANDXIDEFS}.
Then by \eqref{E:INITIALRELATIONS}, we have
\begin{align} \label{E:PLANESYMMETRYRADATTIME0}
	\Rad|_{\Sigma_0} = - \partial_1.
\end{align}
Hence, by \eqref{E:RADDEF}, we have
\begin{align}  \label{E:SIMPLEULGOODIDENTITIES}
	\uLgood 
		& = \upmu \Lunit + 2 \upmu \Radunit
			= \upmu \Lunit + 2 \Rad.
\end{align}

The desired initial data can be constructed
by simply taking smooth data 
$(\mathcal{R}_-|_{\Sigma_0^1},\mathcal{R}_+|_{\Sigma_0^1})$
for the system \eqref{E:EVOLUTIONRIEMANNINVARIANT} that are supported in 
$\Sigma_0^1$
such that
$\mathcal{R}_-|_{\Sigma_0^1} \equiv 0$,
$\sum_{M=1}^3 \| \Rad^M \mathcal{R}_+ \|_{L^{\infty}(\Sigma_0^1)} \leq \mathring{\updelta}'$,
and
$\| \mathcal{R}_+ \|_{L^{\infty}(\Sigma_0^1)} \leq \mathring{\upepsilon}'$,
where $\mathring{\upepsilon}'$ and $\mathring{\updelta}'$ 
verify the same
relative size assumptions 
as the parameters $\mathring{\upepsilon}$ and $\mathring{\updelta}$
described in Subsect.~\ref{SS:SMALLNESSASSUMPTIONS}. 
As we now outline, this leads to the desired size assumptions
stated in Subsects.~\ref{SS:FLUIDVARIABLEDATAASSUMPTIONS},\ref{SS:DATAFOREIKONALFUNCTIONQUANTITIES},
and \ref{SS:SMALLNESSASSUMPTIONS},
where the smallness of $\mathring{\upepsilon}$ is induced by the smallness of
$\mathring{\upepsilon}'$
and the relative largeness of $\mathring{\updelta}$ is tied to the relative largeness of
$\mathring{\updelta}'$.

We first note that the support assumption on the data
implies that the solution completely vanishes along $\mathcal{P}_0$,
consistent with the data assumptions made in
Subsects.~\ref{SS:FLUIDVARIABLEDATAASSUMPTIONS} and \ref{SS:DATAFOREIKONALFUNCTIONQUANTITIES}.
We next note that the first evolution equation in \eqref{E:EVOLUTIONRIEMANNINVARIANT}
implies that $\mathcal{R}_- \equiv 0$.
One can derive estimates for the mixed derivatives of
$\mathcal{R}_+|_{\Sigma_0^1}$ with respect to 
$\Lunit$ and $\Rad$ by commuting the second evolution equation
in \eqref{E:EVOLUTIONRIEMANNINVARIANT}.
In view of the simple commutation relation
$[\Lunit, \Rad] = 0$, valid in plane symmetry, 
we obtain that 
$\Lunit^{M_1} \Rad^{M_2} \mathcal{R}_+ = 0$ if $M_1 \geq 1$,
from which it easily follows (see equation \eqref{E:UPMUFIRSTTRANSPORT})
that if $M_1 \geq 1$, then
$\Lunit \Lunit^{M_1} \Rad^{M_2} \upmu = 0$.

From these facts, one can
show that all of the data assumptions stated in
Subsects.~\ref{SS:FLUIDVARIABLEDATAASSUMPTIONS} and \ref{SS:DATAFOREIKONALFUNCTIONQUANTITIES}
are verified if $\mathring{\upepsilon}'$ is sufficiently small.
We do not give a the full proof here because it is straightforward but tedious;
instead, we prove four representative estimates.
First, using \eqref{E:SPEEDOFSOUNDISUNITY}, \eqref{E:RIEMANNINVARIANTS},
Taylor expansions, and the fact that $\mathcal{R}_- \equiv 0$,
we obtain 
$\Rad (v^1 - \Densrenormalized)
= \Rad \mathcal{R}_- + (1 - \Speed) \Rad \Densrenormalized
= \mathcal{O}(\mathcal{R}_+)  \Rad \mathcal{R}_+
$.
Hence, using the above estimates, we obtain
$\| \Rad (v^1 - \Densrenormalized) \|_{L^{\infty}(\Sigma_0^1)}
\lesssim \mathring{\upepsilon}'
$,
which is consistent with the smallness assumption \eqref{E:L2SMALLDATAASSUMPTIONSALONGSIGMA0}
for the first entry of the second term on the LHS.
As a second example,
we note that with the help of \eqref{E:INITIALRELATIONS},
we have 
$\upmu|_{\Sigma_0^1}
= 1 + \mathcal{O}(\mathcal{R}_+)
$.
Hence, using the above estimates, we obtain
$\| \upmu - 1 \|_{L^{\infty}(\Sigma_0^1)}
\lesssim \mathring{\upepsilon}'
$
(consistent with the smallness assumption \eqref{E:UPMUDATATANGENTIALLINFINITYCONSEQUENCES}
for the first term on the LHS)
and
$\| \Rad^M \upmu \|_{L^{\infty}(\Sigma_0^1)}
\lesssim \mathring{\updelta}'
\lesssim 1
$
for $M=1,2$,
consistent with the assumptions stated in
\eqref{E:UPMUDATARADIALLINFINITYCONSEQUENCES}
for the last term on the LHS.
As a third example, we note that with the help of \eqref{E:INITIALRELATIONS},
it is easy to show that\footnote{Here we are viewing plane symmetric solutions to be solutions on 
$\mathbb{R} \times \mathbb{R} \times \mathbb{T}$.} 
$\NonRadialRad^i \equiv 0$
and
$\CoordAng^i \equiv \delta_2^i$
in the maximal development of the data,
which in particular is consistent with 
\eqref{E:LINFTYCOORADANGILARGEDATAALONGSIGMA1}-\eqref{E:LINFTYXILARGEDATAALONGSIGMA1}.
As a last example, we note that $\Vortrenormalized \equiv 0$ in plane symmetry,
which is consistent with the smallness assumptions 
\eqref{E:L2SMALLDATAASSUMPTIONSALONGSIGMA0},
\eqref{E:LINFTYSMALLDATAASSUMPTIONSALONGSIGMA0},
and \eqref{E:SMALLDATAASSUMPTIONSALONGELLT0} for $\Vortrenormalized$.

\section{Preliminary pointwise estimates}
\label{S:PRELIMINARYPOINTWISE}
In this section, we derive preliminary pointwise estimates for the simplest error
terms that appear in the commuted equations.
Our arguments rely on the data-size assumptions 
and bootstrap assumptions stated in Sect.~\ref{S:DATASIZEANDBOOTSTRAP}
and are tedious to carry out but not too difficult.

In the remainder of the article, 
we schematically express many equations and inequalities
by stating them in terms of the arrays $\GdVar$ and $\BadVar$ 
from Def.~\ref{D:ABBREIVATEDVARIABLES}.
We also remind the reader that we often use the abbreviations
introduced Sect.~\ref{SS:STRINGSOFCOMMUTATIONVECTORFIELDS} 
to schematically indicate the structure of various differential 
operators.

\subsection{Differential operator comparison estimates}
\label{SS:DIFFOPCOMPARISONESTIMATES}
In this section, we provide quantitative comparison estimates
relating various differential operators on $\ell_{t,u}$.

We start by providing a simple lemma in which we
express $\angLap f$
and $\angD^2 f$
in terms of derivatives with respect to the vectorfield $\GeoAng$.

\begin{lemma}[\textbf{$\angLap$ and $\angDsquared$ in terms of $\GeoAng$ derivatives}] 
\label{L:ANGLAPINTERMSOFGEOANGDERIVATIVES}
	Let $\GeoAng$ be the $\ell_{t,u}$-tangent vectorfield defined 
	in Def.~\ref{D:ANGULARVECTORFIELDS}.
	We have the following differential operator identities, 
	valid for scalar-valued functions $f$ defined on $\ell_{t,u}$:
	\begin{subequations}
	\begin{align} \label{E:ANGLAPINTERMSOFGEOANGDERIVATIVES}
		\angLap f
		& = \frac{1}{g(\GeoAng,\GeoAng)} \GeoAng \GeoAng f
			- 
			\frac{1}{g(\GeoAng,\GeoAng)}
			\left\lbrace
				\GeoAng 
				\ln g(\GeoAng,\GeoAng) 
			\right\rbrace
			\GeoAng f,
				\\
		\angDsquared f
		& = \frac{1}{\left\lbrace g(\GeoAng,\GeoAng) \right\rbrace^2} 
				(\GeoAng \GeoAng f) \GeoAng_{\flat} \otimes \GeoAng_{\flat}
			- 
			\frac{1}{\left\lbrace g(\GeoAng,\GeoAng) \right\rbrace^2}
			\left\lbrace
				\GeoAng 
				\ln g(\GeoAng,\GeoAng) 
			\right\rbrace
			\GeoAng f.
			\label{E:ANGDSQUAREDINTERMSOFGEOANGDERIVATIVES}
	\end{align}
	\end{subequations}
\end{lemma}

\begin{proof}
	Using \eqref{E:INVERSEANGULARMETRICINTERMSOFGEOANG},
	we obtain
	\begin{align} \label{E:FIRSTIDANGLAPINTERMSOFGEOANGDERIVATIVES}
		\angLap f
		& = \frac{1}{g(\GeoAng,\GeoAng)} \GeoAng \GeoAng f
			- 
			\frac{1}{g(\GeoAng,\GeoAng)} (\angD_{\GeoAng} \GeoAng) \cdot \angdiff f.
	\end{align}
	Since $\angD_{\GeoAng} \GeoAng$ is $\ell_{t,u}$-tangent, it must be a scalar-valued function multiple,
	denoted by $M$, of $\GeoAng$:
	$
	M \GeoAng
	=
	\angD_{\GeoAng} \GeoAng 
	$.
	Taking the inner product of this identity with $\GeoAng$,
	we obtain
	$M g(\GeoAng,\GeoAng) 
	=
	\gsphere(\angD_{\GeoAng} \GeoAng,\GeoAng)
	=
	\frac{1}{2} 
	\angD_{\GeoAng} 
	\left\lbrace
		\gsphere (\GeoAng,\GeoAng)
	\right\rbrace
	=
	\frac{1}{2} 
	\GeoAng
	\left\lbrace
		g(\GeoAng,\GeoAng)
	\right\rbrace
	$.
	Solving for $M$ and substituting into
	\eqref{E:FIRSTIDANGLAPINTERMSOFGEOANGDERIVATIVES}, we conclude
	\eqref{E:ANGLAPINTERMSOFGEOANGDERIVATIVES}.

	\eqref{E:ANGDSQUAREDINTERMSOFGEOANGDERIVATIVES} then follows 
	from \eqref{E:ANGLAPINTERMSOFGEOANGDERIVATIVES}
	and the identity \eqref{E:XIINTERMSOFTRACEXI} with
	$\xi = \angDsquared f$.

\end{proof}

The next lemma shows that the pointwise norms of 
$\ell_{t,u}$ tensors are controlled by contractions
against $\GeoAng$.

\begin{lemma}[\textbf{The norm of $\ell_{t,u}$-tangent tensors can be measured via $\GeoAng$ contractions}]
	\label{L:TENSORSIZECONTROLLEDBYYCONTRACTIONS}
	Let $\xi_{\alpha_1 \cdots \alpha_n}$ be a type $\binom{0}{n}$ $\ell_{t,u}$-tangent tensor with $n \geq 1$
	and let $\GeoAng$ be the $\ell_{t,u}$-tangent vectorfield defined 
	in Def.~\ref{D:ANGULARVECTORFIELDS}.
	Under the data-size and bootstrap assumptions 
	of Sects.~\ref{SS:FLUIDVARIABLEDATAASSUMPTIONS}-\ref{SS:AUXILIARYBOOTSTRAP}
	and the smallness assumptions of Sect.~\ref{SS:SMALLNESSASSUMPTIONS}, 
	we have
	\begin{align} \label{E:TENSORSIZECONTROLLEDBYYCONTRACTIONS}
		|\xi| = 
			\left\lbrace 
				1 + \mathcal{O}(\varepsilon^{1/2})
			\right\rbrace
			|\xi_{\GeoAng \GeoAng \cdots \GeoAng}|.
	\end{align}
	The same result holds if 
	$|\xi_{\GeoAng \GeoAng \cdots \GeoAng}|$
	is replaced with 
	$|\xi_{\GeoAng \cdot}|$, 
	$|\xi_{\GeoAng \GeoAng \cdot}|$,
	etc., where $\xi_{\GeoAng \cdot}$
	is the type $\binom{0}{n-1}$ tensor with components
  $\GeoAng^{\alpha_1} \xi_{\alpha_1 \alpha_2 \cdots \alpha_n}$,
  and similarly for $\xi_{\GeoAng \GeoAng \cdot}$, etc.
\end{lemma}
\begin{proof}
	\eqref{E:TENSORSIZECONTROLLEDBYYCONTRACTIONS} is easy to derive 
	relative to Cartesian coordinates by using 
	the decomposition $(\ginversesphere)^{ij} = \frac{1}{|\GeoAng|^2} \GeoAng^i \GeoAng^j$
	(see \eqref{E:INVERSEANGULARMETRICINTERMSOFGEOANG})
	and the estimate $|\GeoAng| = 1 + \mathcal{O}(\varepsilon^{1/2})$,
	which follows from
	the identity 
	$|\GeoAng|^2 
	= g_{ab} \GeoAng^a \GeoAng^b
	= (\delta_{ab} + g_{ab}^{(Small)}) (\delta_2^a + \GeoAng_{(Small)}^a)(\delta_2^b + \GeoAng_{(Small)}^b)$,
	the schematic relations $g_{ab}^{(Small)}, \GeoAng_{(Small)}^a = \smoothfunction(\GdVar)\GdVar$
	(see Lemma~\ref{L:SCHEMATICDEPENDENCEOFMANYTENSORFIELDS}),
	and the bootstrap assumptions.
\end{proof}

We now establish some comparison estimates for various differential operators
on $\ell_{t,u}$.

\begin{lemma}[\textbf{Controlling $\angD$ derivatives in terms of $\GeoAng$ derivatives}]
\label{L:ANGDERIVATIVESINTERMSOFTANGENTIALCOMMUTATOR}
	Let $f$ be a scalar-valued function on $\ell_{t,u}$.
	Under the data-size and bootstrap assumptions 
	of Sects.~\ref{SS:FLUIDVARIABLEDATAASSUMPTIONS}-\ref{SS:AUXILIARYBOOTSTRAP}
	and the smallness assumptions of Sect.~\ref{SS:SMALLNESSASSUMPTIONS}, 
	the following comparison estimates hold
	on $\mathcal{M}_{\Tboot,U_0}$:
\begin{align} \label{E:ANGDERIVATIVESINTERMSOFTANGENTIALCOMMUTATOR}
		|\angdiff f|
		& \leq (1 + C \varepsilon^{1/2})\left| \GeoAng f \right|,
		\qquad
		|\angD^2 f|
		\leq (1 + C \varepsilon^{1/2})\left| \angdiff(\GeoAng f) \right|
			+ C \varepsilon |\angdiff f|.
\end{align}

\end{lemma}

\begin{proof}
	The first inequality in \eqref{E:ANGDERIVATIVESINTERMSOFTANGENTIALCOMMUTATOR}
	follows directly from Lemma~\ref{L:TENSORSIZECONTROLLEDBYYCONTRACTIONS}.
	To prove the second, we first use
	Lemma~\ref{L:TENSORSIZECONTROLLEDBYYCONTRACTIONS},
	the identity 
	$\angD_{\GeoAng \GeoAng}^2 f = \GeoAng \cdot \angdiff (\GeoAng f) - \angD_{\GeoAng} \GeoAng \cdot \angdiff f$,
	and the estimate $|\GeoAng| = 1 + \mathcal{O}(\varepsilon^{1/2})$ noted in the proof of
	Lemma~\ref{L:TENSORSIZECONTROLLEDBYYCONTRACTIONS}
	to deduce 
	\begin{align} \label{E:ANGDSQUAREFUNCTIONFIRSTBOUNDINTERMSOFGEOANG}
		|\angD^2 f| 
		& \leq (1 + C \varepsilon^{1/2})|\angD_{\GeoAng \GeoAng}^2 f|
			\leq (1 + C \varepsilon^{1/2})|\angdiff(\GeoAng f)|
			+ |\angD_{\GeoAng} \GeoAng||\angdiff f|.
	\end{align}
	Next, we use Lemma~\ref{L:TENSORSIZECONTROLLEDBYYCONTRACTIONS} 
	and the identity 
	$	\angdeformarg{\GeoAng}{\GeoAng}{\GeoAng} 
		= \angD_{\GeoAng} (\gsphere(\GeoAng, \GeoAng))
		= \GeoAng (g_{ab} \GeoAng^a \GeoAng^b)
	$
	to deduce that
	\begin{align} \label{E:ANGDGEOANGOFGEOANGINTERMSOFGEOANGDEFORMSPHERE}
		\left|
			\angD_{\GeoAng} \GeoAng
		\right|
		& \lesssim
			\left|
				g(\angD_{\GeoAng} \GeoAng,\GeoAng)
			\right|
			\lesssim
			\left|
				\angdeformarg{\GeoAng}{\GeoAng}{\GeoAng}
			\right|
			\lesssim
			\left|
				\GeoAng(g_{ab} \GeoAng^a \GeoAng^b)
			\right|.
	\end{align}
	Since Lemma~\ref{L:SCHEMATICDEPENDENCEOFMANYTENSORFIELDS} implies that
	$g_{ab} \GeoAng^a \GeoAng^b = \smoothfunction(\GdVar)$ with $\smoothfunction$ smooth,
	the bootstrap assumptions yield that 
	$\mbox{{\upshape RHS}~\eqref{E:ANGDGEOANGOFGEOANGINTERMSOFGEOANGDEFORMSPHERE}}
	\lesssim |\GeoAng \GdVar| \lesssim \varepsilon^{1/2}$.
	The desired estimate for $|\angD^2 f|$
	now follows from this estimate,
	\eqref{E:ANGDSQUAREFUNCTIONFIRSTBOUNDINTERMSOFGEOANG},
	and
	\eqref{E:ANGDGEOANGOFGEOANGINTERMSOFGEOANGDEFORMSPHERE}.
\end{proof}

\begin{lemma} [\textbf{Controlling $\angLie_V$ and $\angD$ derivatives in terms of $\angLie_{\GeoAng}$ derivatives}]
	\label{L:ANGLIEPXIINTERMSOFANGLIEGEOANGXI}
	Let $\xi_{\alpha_1 \cdots \alpha_n}$ be a type $\binom{0}{n}$ $\ell_{t,u}$-tangent tensor with $n \geq 1$
	and let $V$ be an $\ell_{t,u}$-tangent vectorfield.
	Under the data-size and bootstrap assumptions 
	of Sects.~\ref{SS:FLUIDVARIABLEDATAASSUMPTIONS}-\ref{SS:AUXILIARYBOOTSTRAP}
	and the smallness assumptions of Sect.~\ref{SS:SMALLNESSASSUMPTIONS}, 
	the following comparison estimates hold
	on $\mathcal{M}_{\Tboot,U_0}$:
	\begin{align} \label{E:ANGLIEPXIINTERMSOFANGLIEGEOANGXI}
		\left| 
			\angLie_V \xi
		\right|
		& \lesssim 
			|V| 
			 \left|
			 	\angLie_{\GeoAng} \xi
			 \right|
			+ |\xi|
				\left|
					\angLie_{\GeoAng} V
				\right|
			+ \left|
					\GeoAng \GdVar
				\right|
				|\xi| 
				|V|
					\\
		& \lesssim 
			|V| 
			 \left|
			 	\angLie_{\GeoAng} \xi
			 \right|
			+ 
			|\xi|
				\left|
					\angLie_{\GeoAng} V
				\right|
			+ 
			\varepsilon^{1/2} |\xi| |V|,
				\notag \\
		|\angD \xi|
		& \lesssim
			\left|\angLie_{\GeoAng} \xi \right|
			+ 
			\left|\GeoAng \GdVar \right| |\xi|
				\label{E:ANGDPIINTERMSOFLIEDERIVATIVES} \\
		& \lesssim
			\left|\angLie_{\GeoAng} \xi \right|
			+ 
			\varepsilon^{1/2} |\xi|.
		\notag
	\end{align}

\end{lemma}

\begin{proof}
	To prove \eqref{E:ANGLIEPXIINTERMSOFANGLIEGEOANGXI},
	we use the schematic Lie derivative identity
	$\angLie_V \xi 
	= 
	\angD_V \xi + \sum \xi \cdot \angD V
	$
	and Lemma~\ref{L:TENSORSIZECONTROLLEDBYYCONTRACTIONS}
	to deduce
	\begin{align} \label{E:FIRSTBOUNDANGLIEPXIINTERMSOFANGD}
		|\angLie_V \xi| \lesssim |V||\angD_{\GeoAng} \xi| + |\xi| |\angD_{\GeoAng} V|.
	\end{align}
	Next, we note that the torsion-free property of $\angD$ implies that
	$\angD_{\GeoAng} V 
		= \angLie_{\GeoAng} V 
		+ \angD_V \GeoAng
	$.
	Hence, using Lemma~\ref{L:TENSORSIZECONTROLLEDBYYCONTRACTIONS},
	\eqref{E:ANGDGEOANGOFGEOANGINTERMSOFGEOANGDEFORMSPHERE},
	and the estimate 
	$|\angD_{\GeoAng} \GeoAng| 
	\lesssim |\GeoAng \GdVar|
	\lesssim \varepsilon^{1/2}
	$
	shown in the proof of Lemma~\ref{L:ANGDERIVATIVESINTERMSOFTANGENTIALCOMMUTATOR},
	we find that
	\begin{align}  \label{E:ANGDGEOANGXINTERMSOFLIEDERIVATIVESBOUND}
		|\angD_{\GeoAng} V|
		& \lesssim 
			|\angLie_{\GeoAng} V|
			+ |V||\angD \GeoAng|
			\lesssim 
			|\angLie_{\GeoAng} V|
			+ |V||\angD_{\GeoAng} \GeoAng|
			\lesssim 
			|\angLie_{\GeoAng} V|
			+ 
			|\GeoAng \GdVar| |V|	\\
	& \ \
			\lesssim 
			|\angLie_{\GeoAng} V|
			+ 
			\varepsilon^{1/2} |V|.
			\notag
	\end{align}
	Similarly, we have
	\begin{align} \label{E:ANGDGEOANGXIINTERMSOFLIEDERIVATIVESBOUND}
		|\angD_{\GeoAng} \xi|
		& \lesssim |\angLie_{\GeoAng} \xi|
			+ \varepsilon^{1/2} |\xi|.
	\end{align}
	The desired estimate \eqref{E:ANGLIEPXIINTERMSOFANGLIEGEOANGXI}
	now follows from
	\eqref{E:FIRSTBOUNDANGLIEPXIINTERMSOFANGD},
	\eqref{E:ANGDGEOANGXINTERMSOFLIEDERIVATIVESBOUND},
	and
	\eqref{E:ANGDGEOANGXIINTERMSOFLIEDERIVATIVESBOUND}.

	The estimate \eqref{E:ANGDPIINTERMSOFLIEDERIVATIVES}
	follows from applying Lemma~\ref{L:TENSORSIZECONTROLLEDBYYCONTRACTIONS} to $\angD \xi$ and
	using \eqref{E:ANGDGEOANGXIINTERMSOFLIEDERIVATIVESBOUND}.
\end{proof}

\subsection{Basic facts and estimates that we use silently}
\label{SS:OFTENUSEDESTIMATES}
For the reader's convenience, we present here some basic facts and estimates
that we silently use throughout the rest of the paper when deriving
estimates.

\begin{enumerate}
	\item All quantities that we estimate can be controlled in terms of 
		$\BadVar = \lbrace \threePsi, \upmu, \Lunit_{(Small)}^1, \Lunit_{(Small)}^2 \rbrace$
		and the specific vorticity $\Vortrenormalized$.
	\item We typically use the Leibniz rule for 
		the operators $\angLie_Z$ and $\angD$ when deriving
		pointwise estimates for the $\angLie_Z$ and $\angD$ 
		derivatives of tensor products of the schematic form 
		$\prod_{i=1}^m v_i$, where the $v_i$ are scalar functions or
		$\ell_{t,u}$-tangent tensors. Our derivative counts are such that
		all $v_i$ except at most one are uniformly bounded in $L^{\infty}$
		on $\mathcal{M}_{\Tboot,U_0}$.
		Thus, our pointwise estimates often explicitly feature 
		(on the right-hand sides)
		only one factor with many derivatives on it,
		multiplied by a constant that uniformly bounds the other factors.
		In some estimates, the right-hand sides also gain a smallness factor,
		such as $\varepsilon^{1/2}$,
		generated by the remaining $v_i's$.
	\item The operators $\angLie_{\Fullset}^N$ commute through
		$\angdiff$, as shown by Lemma~\ref{L:LANDRADCOMMUTEWITHANGDIFF}.
	\item As differential operators acting on scalar functions, we have
		$\GeoAng 
		= \left(1 + \mathcal{O}(\GdVar)\right) \angdiff
		= (1 + \mathcal{O}(\varepsilon^{1/2})) \angdiff
		$,
		a fact which follows from 
		Lemma~\ref{E:ANGDERIVATIVESINTERMSOFTANGENTIALCOMMUTATOR},
		\eqref{E:GEOANGPOINTWISE}, and the bootstrap assumptions.
		Hence, for scalar functions $f$, we sometimes schematically depict $\angdiff f$ 
		as 
		$\left(1 + \mathcal{O}(\GdVar)\right) \Singletan f$
		or
		$\left(1 + \mathcal{O}(\GdVar)\right) \Fullset_{**}^{1;0} f$,
		or alternatively as
		$\Singletan f$ or
		$\Fullset_{**}^{1;0} f$
		when the factor 
		$1 + \mathcal{O}(\GdVar)$ is not important.
		Similarly, by
		Lemma~\ref{L:ANGDERIVATIVESINTERMSOFTANGENTIALCOMMUTATOR} 
		we can depict $\angLap f$ by 
		$
		\smoothfunction(\Tanset^{\leq 1} \GdVar,\ginversesphere)
		\Fullset_{**}^{[0,2];0} f
		$
		(or $\Fullset_{**}^{[0,2];0} f$
		when the factor $\smoothfunction(\Tanset^{\leq 1} \GdVar,\ginversesphere)$
		is not important).
		Similarly, by Lemma~\ref{L:ANGLIEPXIINTERMSOFANGLIEGEOANGXI}, 
		for type $\binom{0}{n}$ $\ell_{t,u}$-tangent tensorfields $\xi$,
		we can depict $\angD \xi$ by 
		$
		\smoothfunction(\Tanset^{\leq 1} \GdVar,\ginversesphere)
		\angLie_{\Tanset}^{\leq 1} \xi
		$
		(or $\angLie_{\Tanset}^{\leq 1} \xi$
		when the factor $\smoothfunction(\Tanset^{\leq 1} \GdVar,\ginversesphere)$
		is not important).
	\item We remind the reader that all constants are allowed to depend on 
		the data-size parameters
		$\mathring{\updelta}$
		and 
		$\TranminusdatasizeWithFactor^{-1}$.
\end{enumerate}

\subsection{Pointwise estimates for the Cartesian coordinates and the Cartesian components of some vectorfields}
\label{SS:RECTANGULARCOMPONENTSPOINTWISEESTIMATES}

\begin{lemma}[\textbf{Pointwise estimates for} $x^i$ \textbf{and the Cartesian components of several vectorfields}]
	\label{L:POINTWISEFORRECTANGULARCOMPONENTSOFVECTORFIELDS}
	Assume that\footnote{Throughout, we use the convention that 
	terms in our formulas and estimates involving operators that do not make sense
	are absent. $\Fullset_{\ast}^{1;1}$ is an example of such an operator.
	\label{FN:TERMSTHATDONOTMAKESENSE}} 
	$1 \leq N \leq 20$, 
	$0 \leq M \leq \min \lbrace N,2\rbrace$,
	and $V \in \lbrace \Lunit, \Radunit, \GeoAng \rbrace$.
	Let $x^i = x^i(t,u,\vartheta)$ denote the Cartesian spatial coordinate function
	and let $\mathring{x}^i = \mathring{x}^i(u,\vartheta) := x^i(0,u,\vartheta)$.
	Under the data-size and bootstrap assumptions 
	of Sects.~\ref{SS:FLUIDVARIABLEDATAASSUMPTIONS}-\ref{SS:AUXILIARYBOOTSTRAP}
	and the smallness assumptions of Sect.~\ref{SS:SMALLNESSASSUMPTIONS}, 
	the following pointwise estimates hold
	on $\mathcal{M}_{\Tboot,U_0}$,
	for $i = 1,2$
	(see Sect.~\ref{SS:STRINGSOFCOMMUTATIONVECTORFIELDS} regarding the vectorfield operator notation):
	\begin{subequations}
	\begin{align}
		\left| 
			V^i
		\right|
		& \lesssim 
			1
			+
			\left|
				\GdVar
			\right|,
		\label{E:ORDERZEROVECTORFIELDSRECTCOMPPOINTWISE} \\		 
		\left| 
			\Fullset^{[1,N];M} V^i
		\right|
		& \lesssim 
			\left|
				\Fullset^{[1,N];\leq M} \GdVar
			\right|,
		\label{E:NOSPECIALSTRUCTURECOMMMUTEDVECTORFIELDSRECTCOMPPOINTWISE} \\
		\left| 
			\Fullset_{\ast}^{[1,N];M} V^i
		\right|
		& \lesssim 
			\left|
				\Fullset_{\ast}^{[1,N];\leq M} \GdVar
			\right|.
		\label{E:STARTCOMMMUTEDVECTORFIELDSRECTCOMPPOINTWISE} 
\end{align}
\end{subequations}

Similarly, if 
$1 \leq N \leq 20$
and $0 \leq M \leq \min \lbrace N,1\rbrace$, 
then
\begin{subequations}
\begin{align}
		\left| 
			\Rad^i
		\right|
		& \lesssim
			1 + |\BadVar|,
			\label{E:RADRECTCOMPPOINTWISE} \\
		\left| 
			\Fullset^{[1,N];M} \Rad^i
		\right|
		& \lesssim
			\left|
				\Fullset^{[1,N];\leq M}
				\BadVar
			\right|,
			\label{E:NOSPECIALSTRUCTURECOMMUTEDRADIPOINTWISE} 
				\\
		\left| 
			\Fullset_{\ast}^{[1,N];M} \Rad^i
		\right|
		& \lesssim
			\left|
				\Fullset_{\ast}^{[1,N];\leq M}
				\BadVar
			\right|,
			\label{E:STARCOMMUTEDRADIPOINTWISE} 
				\\
		\left| 
			\Fullset_{\ast \ast}^{[1,N];M} \Rad^i
		\right|
		& \lesssim
			\left|
				\myarray
					[\Fullset_{\ast \ast}^{[1,N];\leq M} \BadVar]
					{\Fullset_{\ast}^{[1,N];\leq M} \GdVar}
				\right|.
			\label{E:DOUBLESTARCOMMUTEDRADIPOINTWISE} 
	\end{align}
	\end{subequations}

	Moreover, if 
	$1 \leq N \leq 20$
	and
	$0 \leq M \leq \min \lbrace N,2\rbrace$, 
	then
	\begin{subequations}
		\begin{align} \label{E:XIPOINTWISE}
			\left|
				x^i - \mathring{x}^i
			\right|
			& \lesssim 1,
				\\
			\left|
				\angdiff x^i
			\right|
			& \lesssim 
				1 
				+
				\left| 
					\GdVar 
				\right|,
				 \label{E:ANGDIFFXI} \\
			\left|
				\angdiff \Fullset^{[1,N];M} x^i
			\right|
			& \lesssim 
				\left|
				\myarray
					[\Fullset_{\ast \ast}^{[1,N];\leq (M-1)_+} \BadVar]
					{\Fullset_{\ast}^{[1,N];\leq M} \GdVar}
				\right|.
					\label{E:ANGDIFFXNOSPECIALSTRUCTUREDIFFERENTIATEDPOINTWISE} 
		\end{align} 
	\end{subequations}

	Finally, if $0 \leq N \leq 20$
	and
	$0 \leq M \leq \min \lbrace N,2\rbrace$,
	then
	\begin{subequations}
	\begin{align} 
	\left|
		\Fullset^{N;M} \GeoAng_{(Small)}^i
	\right|
	& \lesssim 
		\left|
			\Fullset^{\leq N;\leq M} \GdVar
		\right|,
		\label{E:NOSPECIALSTRUCTURECOMMMUTEDGEOANGSMALLPOINTWISE}
			\\
	\left|
		\Fullset_{\ast}^{N;M} \GeoAng_{(Small)}^i
	\right|
	& \lesssim 
		\left|
			\Fullset_{\ast}^{\leq N;\leq M} \GdVar
		\right|.
		\label{E:STARCOMMMUTEDGEOANGSMALLPOINTWISE} 
	\end{align}
	\end{subequations}
	In the case $i=2$ at fixed $u,\vartheta$,
	LHS \eqref{E:XIPOINTWISE} is to be interpreted as
	the Euclidean distance traveled
	by the point $x^2$
	in the flat universal covering space 
	$\mathbb{R}$ of $\mathbb{T}$
	along the corresponding integral curve of $\Lunit$
	over the time interval $[0,t]$.
\end{lemma}

\begin{proof}
See Sect.~\ref{SS:OFTENUSEDESTIMATES} for some comments on the analysis.
Lemma~\ref{L:SCHEMATICDEPENDENCEOFMANYTENSORFIELDS} implies
that for $V \in \lbrace \Lunit, \Radunit, \GeoAng \rbrace$,
the component $V^i = V x^i$ verifies $V^i = \smoothfunction(\GdVar)$
with $\smoothfunction$ smooth.
Similarly, $\GeoAng_{(Small)}^i$ verifies
$\GeoAng_{(Small)}^i = \smoothfunction(\GdVar)\GdVar$
with $\smoothfunction$ smooth
and
$\Rad x^i = \Rad^i$ verifies 
$\Rad^i = \upmu \smoothfunction(\GdVar)$
with $\smoothfunction$ smooth.
The estimates of the lemma therefore follow easily
from the bootstrap assumptions,
except for the estimates 
\eqref{E:XIPOINTWISE}-\eqref{E:ANGDIFFXNOSPECIALSTRUCTUREDIFFERENTIATEDPOINTWISE}.
To obtain \eqref{E:XIPOINTWISE}, 
we first argue as above to deduce $|\Lunit x^i| = |\Lunit^i| = |\smoothfunction(\GdVar)| \lesssim 1$.
Since $\Lunit = \frac{\partial}{\partial t}$, 
we may integrate along the integral curves of $\Lunit$ starting from time $1$
to deduce, via the fundamental theorem of calculus, that
\begin{align} \label{E:SIMPLEFTCID}
	x^i(t,u,\vartheta)
	& = x^i(0,u,\vartheta)
		+ 
		\int_{s=0}^t
			\Lunit x^i(s,u,\vartheta)
		\, ds.
\end{align}
Taking the absolute value of \eqref{E:SIMPLEFTCIDENTITY}
and using the estimate $|\Lunit x^i| \lesssim 1$
to bound the time integral by $\lesssim t \lesssim 1$,
we conclude \eqref{E:XIPOINTWISE}.
To derive \eqref{E:ANGDIFFXI},
we use \eqref{E:ANGDERIVATIVESINTERMSOFTANGENTIALCOMMUTATOR}
with $f=x^i$ to deduce 
$|\angdiff x^i| 
\lesssim |\GeoAng x^i| 
= |\GeoAng^i| 
= |\smoothfunction(\GdVar)| 
\lesssim 
1 
+
|\GdVar|
$
as desired.
The proof of 
\eqref{E:ANGDIFFXNOSPECIALSTRUCTUREDIFFERENTIATEDPOINTWISE}
is similar, but we also use
Lemma~\ref{L:LANDRADCOMMUTEWITHANGDIFF}
to commute vectorfields under $\angdiff$.
\end{proof}

\subsection{Pointwise estimates for various \texorpdfstring{$\ell_{t,u}$-}{}tensorfields}

\begin{lemma}[\textbf{Crude pointwise estimates for the Lie derivatives of $\gsphere$ and 
$\ginversesphere$}]
\label{L:POINTWISEESTIMATESFORGSPHEREANDITSDERIVATIVES}
	Assume that $1 \leq N \leq 20$ and $0 \leq M \leq \min \lbrace N,2\rbrace$.
	Under the data-size and bootstrap assumptions 
	of Sects.~\ref{SS:FLUIDVARIABLEDATAASSUMPTIONS}-\ref{SS:AUXILIARYBOOTSTRAP}
	and the smallness assumptions of Sect.~\ref{SS:SMALLNESSASSUMPTIONS}, 
	the following pointwise estimates hold
	on $\mathcal{M}_{\Tboot,U_0}$
	(see Sect.~\ref{SS:STRINGSOFCOMMUTATIONVECTORFIELDS} regarding the vectorfield operator notation):
	\begin{subequations}
	\begin{align} \label{E:POINTWISEESTIMATESFORGSPHEREANDITSNOSPECIALSTRUCTUREDERIVATIVES}
		\left|
			\angLie_{\Fullset}^{N;M} \gsphere
		\right|,
			\,
		\left|
			\angLie_{\Fullset}^{N;M} \ginversesphere
		\right|
		& \lesssim 
			\myarray[
				\left| 
					\Fullset_{\ast \ast}^{[1,N];\leq (M-1)_+} \BadVar
				\right|
			]
			{
			\left| 
				\Fullset^{[1,N];\leq M} \GdVar
			\right|
			},
				\\
	\left|
			\angLie_{\Fullset_{\ast}}^{N;M} \gsphere
		\right|,
			\,
		\left|
			\angLie_{\Fullset_{\ast}}^{N;M} \ginversesphere
		\right|
		& \lesssim 
			\myarray[
				\left| 
					\Fullset_{\ast \ast}^{[1,N];\leq (M-1)_+} \BadVar
				\right|
			]
			{
			\left| 
				\Fullset_{\ast}^{[1,N];\leq M} \GdVar
			\right|
			}.
			\label{E:POINTWISEESTIMATESFORGSPHEREANDITSSTARDERIVATIVES}
	\end{align}

	Moreover, if $0 \leq N \leq 20$ and $0 \leq M \leq \min \lbrace N,2\rbrace$,
	then
	\begin{align}
		\left|
			\angLie_{\Fullset}^{N;M} \upchi
		\right|,
			\,
		\left|
			\Fullset^{N;M} \mytr \upchi
		\right|
		& \lesssim 
			\myarray[
				\left| 
					\Fullset_{\ast \ast}^{[1,N+1];\leq (M-1)_+} \BadVar
				\right|
			]
			{
			\left| 
				\Fullset_{\ast}^{[1,N+1];\leq M} \GdVar
			\right|
			}.
			\label{E:POINTWISEESTIMATESFORCHIANDITSDERIVATIVES}
\end{align}
\end{subequations}
\end{lemma}

\begin{proof}
		See Sect.~\ref{SS:OFTENUSEDESTIMATES} for some comments on the analysis.
		By Lemma~\ref{L:SCHEMATICDEPENDENCEOFMANYTENSORFIELDS}, we have
		$\gsphere = \smoothfunction(\GdVar,\angdiff \vec{x})$.
		The desired estimates
		for $\angLie_{\Fullset}^{N;M} \gsphere$
		and
		$\angLie_{\Fullset_{\ast}}^{N;M} \gsphere$
		thus follow from
		Lemma~\ref{L:POINTWISEFORRECTANGULARCOMPONENTSOFVECTORFIELDS}
		and the bootstrap assumptions.
		The desired estimates for 
		$\angLie_{\Fullset}^{N;M} \ginversesphere$
		and 
		$\angLie_{\Fullset_{\ast}}^{N;M} \ginversesphere$
		then follow from repeated use of the second
		identity in \eqref{E:CONNECTIONBETWEENANGLIEOFGSPHEREANDDEFORMATIONTENSORS}
		and the estimates for $\angLie_{\Fullset}^{N;M} \gsphere$
		and
		$\angLie_{\Fullset_{\ast}}^{N;M} \gsphere$.
		The estimates for $\angLie_{\Fullset}^{N;M} \upchi$ 
		and $\Fullset^{N;M} \mytr \upchi$ 
		follow from 
		the estimates for $\angLie_{\Fullset}^{N+1;M} \gsphere$
		and $\angLie_{\Fullset}^{N+1;M} \ginversesphere$
		since $\upchi \sim \angLie_{\Singletan} \gsphere$
		(see \eqref{E:CHIDEF})
		and $\mytr \upchi \sim \ginversesphere \cdot \angLie_{\Singletan} \gsphere$.
\end{proof}

\begin{lemma}[\textbf{Pointwise estimates for the Lie derivatives of $\GeoAng$ and some deformation tensor components}]
\label{L:MOREPRECISEDEFTENSESTIMATESPOINTWISE}
Assume that $1 \leq N \leq 20$ and $0 \leq M \leq \min \lbrace N,2\rbrace$.
Under the data-size and bootstrap assumptions 
of Sects.~\ref{SS:FLUIDVARIABLEDATAASSUMPTIONS}-\ref{SS:AUXILIARYBOOTSTRAP}
and the smallness assumptions of Sect.~\ref{SS:SMALLNESSASSUMPTIONS}, 
the following pointwise estimates hold
on $\mathcal{M}_{\Tboot,U_0}$
(see Sect.~\ref{SS:STRINGSOFCOMMUTATIONVECTORFIELDS} regarding the vectorfield operator notation):
\begin{subequations}
		\begin{align} \label{E:GEOANGPOINTWISE}
		\left|
			\left|
				\GeoAng
			\right|
			- 
			1
		\right|
		& \leq
			C
			\left|
				\GdVar
			\right|,
				\\
		\left|
			\angLie_{\Fullset}^{N;M} \GeoAng
		\right|
		& \leq
			C
			\left|
				\myarray
				[\Fullset_{\ast \ast}^{[1,N];\leq (M-1)_+} \BadVar]
				{\Fullset^{[1,N+1];\leq M} \GdVar}
			\right|,
				\label{E:NOSPECIALSTRUCTUREDIFFERENTIATEDGEOANGPOINTWISE}
				\\
		\left|
			\angLie_{\Fullset_{\ast}}^{N;M} \GeoAng
		\right|
		& \leq
			C
			\left|
				\myarray
				[\Fullset_{\ast \ast}^{[1,N-1];\leq (M-1)_+} \BadVar]
				{\Fullset_{\ast}^{[1,N];\leq M} \GdVar}
			\right|.
				\label{E:STARDIFFERENTIATEDGEOANGPOINTWISE}
	\end{align}
	\end{subequations}

Similarly, if $0 \leq N \leq 20$ and $0 \leq M \leq \min \lbrace N,2\rbrace$,
then we have
\begin{subequations}
\begin{align}
	\left|
		\angLie_{\Fullset}^{N;M} \angdeformoneformupsharparg{\GeoAng}{\Lunit}
	\right|
	& \lesssim 
		\left|
			\myarray
			[\Fullset_{\ast \ast}^{[1,N];\leq (M-1)_+} \BadVar]
			{\Fullset_{\ast}^{[1,N+1];\leq M} \GdVar}
		\right|.
				\label{E:NOSPECIALSTRUCTUREDIFFERENTIATEDGEOANGDEFORMSPHERELSHARPPOINTWISE} 
	\end{align}
\end{subequations}

In addition, if $0 \leq N \leq 20$ and $0 \leq M \leq \min \lbrace N,1 \rbrace$,
then we have
\begin{subequations}
\begin{align}
	\left|
		\angLie_{\Fullset}^{N;M} \angdeformoneformupsharparg{\Rad}{\Lunit}
	\right|,
		\,
	\left|
		\angLie_{\Fullset}^{N;M} \angdeformoneformupsharparg{\GeoAng}{\Rad}
	\right|
	& \lesssim 
		\left|
			\Fullset_{\ast}^{[1,N+1];\leq M+1} \threePsi
		\right|
		+
		\left|
			\myarray
			[\Fullset_{\ast \ast}^{[1,N+1];\leq M} \BadVar]
			{\Fullset^{\leq N+1;\leq M} \GdVar}
		\right|,
				\label{E:NOSPECIALSTRUCTUREDIFFERENTIATEDGEOANGDEFORMSPHERERADPOINTWISE}
\end{align}
and if $1 \leq N \leq 20$ and $0 \leq M \leq \min \lbrace N-1,1 \rbrace$,
then 
\begin{align}
	\left|
		\angLie_{\Fullset_{\ast}}^{N;M} \angdeformoneformupsharparg{\Rad}{\Lunit}
	\right|,
		\,
	\left|
		\angLie_{\Fullset_{\ast}}^{N;M} \angdeformoneformupsharparg{\GeoAng}{\Rad}
	\right|
	& \lesssim 
		\left|
			\Fullset_{\ast}^{[1,N+1];\leq M+1} \threePsi
		\right|
		+
		\left|
			\myarray
			[\Fullset_{\ast \ast}^{[1,N+1];\leq M} \BadVar]
			{\Fullset_{\ast}^{[1,N+1];\leq M} \GdVar}
		\right|.
				\label{E:STARDIFFERENTIATEDGEOANGDEFORMSPHERERADPOINTWISE}
\end{align}
\end{subequations}

Moreover, if $0 \leq N \leq 20$
and
$0 \leq M \leq \min \lbrace N,2 \rbrace$,
then we have
\begin{align} \label{E:NOSPECIALSTRUCTUREDIFFERNTIATEDANGDEFORMTANGENTPOINTWISE}
	\left|
		\angLie_{\Fullset}^{N;M} \angdeform{\Lunit}
	\right|,
		\,
	\left|
		\angLie_{\Fullset}^{N;M} \angdeform{\GeoAng}
	\right|
	& \lesssim 
		\left|
			\myarray
			[\Fullset_{\ast \ast}^{[1,N];\leq (M-1)_+} \BadVar]
			{\Fullset_{\ast}^{[1,N+1];\leq M} \GdVar}
		\right|.
\end{align}

In addition, if $0 \leq N \leq 20$ and $0 \leq M \leq \min \lbrace N,1 \rbrace$,
then we have
\begin{subequations}
\begin{align} \label{E:NOSPECIALSTRUCTUREDIFFERNTIATEDANGDEFORMRADPOINTWISE}
	\left|
		\angLie_{\Fullset}^{N;M} \angdeform{\Rad}
	\right|
	& \lesssim 
		\left|
			\Fullset^{[1,N+1];\leq M+1} \threePsi
		\right|
		+
		\left|
			\myarray
			[\Fullset_{\ast \ast}^{[1,N];\leq (M-1)_+} \BadVar]
			{\Fullset_{\ast}^{[1,N+1];\leq M} \GdVar}
		\right|,
\end{align}
and if $1 \leq N \leq 20$ and $0 \leq M \leq \min \lbrace N-1,1 \rbrace$,
then we have
\begin{align}
	\left|
		\angLie_{\Fullset_{\ast}}^{N;M} \angdeform{\Rad}
	\right|
	& \lesssim 
		\left|
			\Fullset_{\ast}^{[1,N+1];\leq M+1} \threePsi
		\right|
		+
		\left|
			\myarray
			[\Fullset_{\ast \ast}^{[1,N];\leq (M-1)_+} \BadVar]
			{\Fullset_{\ast}^{[1,N+1];\leq M} \GdVar}
		\right|.
		\label{E:STARDIFFERNTIATEDANGDEFORMRADPOINTWISE}
\end{align}
\end{subequations}

\end{lemma}

\begin{proof}
	See Sect.~\ref{SS:OFTENUSEDESTIMATES} for some comments on the analysis.
	To prove \eqref{E:NOSPECIALSTRUCTUREDIFFERENTIATEDGEOANGDEFORMSPHERELSHARPPOINTWISE},
	we first note that by 
	Lemma~\ref{L:SCHEMATICDEPENDENCEOFMANYTENSORFIELDS}
	and
	\eqref{E:GEOANGDEFORMSPHEREL},
	we have
	$\angdeformoneformupsharparg{\GeoAng}{\Lunit}
	=
	\smoothfunction(\GdVar,\ginversesphere,\angdiff \vec{x}) \Singletan \GdVar
	$.
	We now apply $\angLie_{\Fullset}^{N;M}$ to the previous relation.
	We bound the derivatives of $\ginversesphere$ and $\angdiff x$
	with Lemmas~\ref{L:POINTWISEFORRECTANGULARCOMPONENTSOFVECTORFIELDS}
	and \ref{L:POINTWISEESTIMATESFORGSPHEREANDITSDERIVATIVES}.
	Also using the bootstrap assumptions, we conclude the desired result.

	Since Lemma~\ref{L:SCHEMATICDEPENDENCEOFMANYTENSORFIELDS} implies that
	$\GeoAng = \smoothfunction(\GdVar,\ginversesphere,\angdiff \vec{x})$,
	similar reasoning yields 
	\eqref{E:NOSPECIALSTRUCTUREDIFFERENTIATEDGEOANGPOINTWISE}-\eqref{E:STARDIFFERENTIATEDGEOANGPOINTWISE}.

	Inequality \eqref{E:GEOANGPOINTWISE}
	follows from the slightly more
	precise arguments already given in the proof of Lemma~\ref{L:TENSORSIZECONTROLLEDBYYCONTRACTIONS}.

	The proofs of 
	\eqref{E:NOSPECIALSTRUCTUREDIFFERENTIATEDGEOANGDEFORMSPHERERADPOINTWISE}-\eqref{E:STARDIFFERENTIATEDGEOANGDEFORMSPHERERADPOINTWISE}
	for 
	$\angdeformoneformupsharparg{\GeoAng}{\Rad}$
	are similar and are based on the observation that by
	Lemma~\ref{L:SCHEMATICDEPENDENCEOFMANYTENSORFIELDS}
	and
	\eqref{E:GEOANGDEFORMSPHERERAD},
	we have
	\[
	\angdeformoneformupsharparg{\GeoAng}{\Rad}
	=
	\smoothfunction(\BadVar,\ginversesphere,\angdiff \vec{x}) \Singletan \GdVar
	+ \smoothfunction(\GdVar,\ginversesphere,\angdiff \vec{x},\Rad \threePsi) \GdVar
	+ \smoothfunction(\GdVar,\ginversesphere) \angdiff \upmu.
	\]

	The proofs of 
	\eqref{E:NOSPECIALSTRUCTUREDIFFERENTIATEDGEOANGDEFORMSPHERERADPOINTWISE}-\eqref{E:STARDIFFERENTIATEDGEOANGDEFORMSPHERERADPOINTWISE}
	for 
	$\angdeformoneformupsharparg{\Rad}{\Lunit}$
	are similar and are based on the observation that by
	Lemma~\ref{L:SCHEMATICDEPENDENCEOFMANYTENSORFIELDS},
	\eqref{E:ZETADECOMPOSED},
	and 
	\eqref{E:RADDEFORMSPHERERAD},
	we have
	\[
	\angdeformoneformupsharparg{\Rad}{\Lunit}
	=
	\smoothfunction(\BadVar,\ginversesphere,\angdiff \vec{x}) \Singletan \threePsi
	+
	\smoothfunction(\GdVar,\ginversesphere,\angdiff \vec{x}, \Rad \threePsi) \GdVar
	+ 
	\ginversesphere \angdiff \upmu.
	\]

	The proof of 
	\eqref{E:NOSPECIALSTRUCTUREDIFFERNTIATEDANGDEFORMTANGENTPOINTWISE}
	is similar and is based on the fact that by
	Lemma~\ref{L:SCHEMATICDEPENDENCEOFMANYTENSORFIELDS},
	\eqref{E:LUNITDEFORMSPHERE},
	and
	\eqref{E:GEOANGDEFORMSPHERE},
	we have
	$\angdeform{\Lunit},
	\angdeform{\GeoAng}
	= \smoothfunction(\GdVar,\ginversesphere,\angdiff \vec{x}) \Singletan \GdVar
	$.

	The proofs of
	\eqref{E:NOSPECIALSTRUCTUREDIFFERNTIATEDANGDEFORMRADPOINTWISE}-\eqref{E:STARDIFFERNTIATEDANGDEFORMRADPOINTWISE}
	are similar 
	and are based on the fact that by
	Lemma~\ref{L:SCHEMATICDEPENDENCEOFMANYTENSORFIELDS}
	and \eqref{E:RADDEFORMSPHERE},
	we have
	$
		\angdeform{\Rad}
		= 
		\smoothfunction(\BadVar,\ginversesphere,\angdiff \vec{x}) \Singletan \GdVar
		+
		\smoothfunction(\GdVar,\angdiff \vec{x}) \Rad \threePsi
	$.
\end{proof}

\subsection{Multi-indices and commutator estimates}
\label{SS:COMMUTATORESTIMATES}
In this section, we establish some commutator estimates.

We start by defining some sets of multi-indices 
corresponding to repeated differentiation with respect
to the commutation vectorfields.

\begin{definition}[\textbf{Sets of multi-indices}]
\label{D:COMMUTATORMULTIINDICES}
We define
\begin{align} \label{E:COMMUTATORMULTIINDICES}
	\mathcal{I}^{N;M}
\end{align}
to be the set of $\Fullset$ multi-indices $\vec{I}$ with the following properties:
\begin{itemize}
	\item $|\vec{I}| = N$.
	\item $\Fullset^{\vec{I}}$ contains at least one factor belonging to 
	$\Tanset = \lbrace \Lunit, \GeoAng \rbrace$.
	\item $\Fullset^{\vec{I}}$ contains precisely $M$ factors of $\Rad$.
\end{itemize}
We define
\begin{align} \label{E:LEQCOMMUTATORMULTIINDICES}
	\mathcal{I}^{N;\leq M},
\end{align}
in the same way,
except the last condition above is replaced with the following one:
\begin{itemize}
	\item $\Fullset^{\vec{I}}$ contains no more than $M$ factors of $\Rad$.
\end{itemize} 
\end{definition}

We now provide two preliminary lemmas from
\cite{jSgHjLwW2016}.

\begin{lemma} \cite{jSgHjLwW2016}*{Lemma 5.1; 
\textbf{Preliminary identities for commuting $Z \in \Fullset$ with $\angD$}}
\label{L:COMMUTINGVEDCTORFIELDSWITHANGD}
For each $\Fullset$-multi-index $\vec{I}$ and integer $n \geq 1$,
the following commutator identity, correct up to constant factors, 
holds for all type $\binom{0}{n}$ $\ell_{t,u}$-tangent tensorfields $\xi$:
\begin{align} \label{E:ANGDANGLIEZELLTUTENSORFIELDCOMMUTATOR}
	[\angD, \angLie_{\Fullset}^{\vec{I}}] \xi
	& = 	\sum_{M=1}^{|\vec{I}|}
				\mathop{\sum_{\vec{I}_1 + \cdots + \vec{I}_{M+1} = \vec{I}}}_{|\vec{I}_a| \geq 1 \mbox{\upshape \ for } 1 \leq a \leq M} 
				(\ginversesphere)^M
				\underbrace{
				(\angLie_{\Fullset}^{\vec{I}_1} \gsphere)
				\cdots
				(\angLie_{\Fullset}^{\vec{I}_{M-1}} \gsphere)
				}_{\mbox{absent when $M=1$}}
				(\angD \angLie_{\Fullset}^{\vec{I}_M} \gsphere)
				(\angLie_{\Fullset}^{\vec{I}_{M+1}} \xi).
\end{align}
Moreover, with $\angdiv$ denoting the torus divergence operator from Def.~\ref{D:CONNECTIONS},
for each $\Fullset$-multi-index $\vec{I}$,
the following commutator identity,
correct up to constant factors,
holds for all symmetric type $\binom{0}{2}$ $\ell_{t,u}$-tangent tensorfields $\xi$:
\begin{align}
	[\angdiv, \angLie_{\Fullset}^{\vec{I}}] \xi
	& = 	\sum_{i_1 + i_2 = 1}
				\sum_{M=1}^{|\vec{I}|}
				\mathop{\sum_{\vec{I}_1 + \cdots + \vec{I}_{M+1} = \vec{I}}}_{|\vec{I}_a| \geq 1 \mbox{\upshape \ for } 1 \leq a \leq M} 
				(\ginversesphere)^{M+1}
				\underbrace{
				(\angLie_{\Fullset}^{\vec{I}_1} \gsphere)
				\cdots
				(\angLie_{\Fullset}^{\vec{I}_{M-1}} \gsphere)
				}_{\mbox{absent when $i_1=M=1$}}
				(\angD^{i_1} \angLie_{\Fullset}^{\vec{I}_M} \gsphere)
				(\angD^{i_2} \angLie_{\Fullset}^{\vec{I}_{M+1}} \xi).
				\label{E:ANGDIVANGLIEZELLTUTENSORFIELDCOMMUTATOR}
\end{align}

Finally, for each $\Fullset$-multi-index $\vec{I}$
and each commutation vectorfield $Z \in \Fullset$,
the following commutator identity,
correct up to constant factors,
holds for all scalar-valued functions $f$:
\begin{subequations}
\begin{align} \label{E:COMMUTINGANGDSQUAREDANDLIEZ}
	[\angD^2, \angLie_{\Fullset}^{\vec{I}}] f
	& = 	\sum_{M=1}^{|\vec{I}|}
				\mathop{\sum_{\vec{I}_1 + \cdots + \vec{I}_{M+1} = \vec{I}}}_{|\vec{I}_a| \geq 1 \mbox{\upshape \ for } 1 \leq a \leq M} 
				(\ginversesphere)^M
				\underbrace{
				(\angLie_{\Fullset}^{\vec{I}_1} \gsphere)
				\cdots
				(\angLie_{\Fullset}^{\vec{I}_{M-1}} \gsphere)
				}_{\mbox{absent when $M=1$}}
				(\angD \angLie_{\Fullset}^{\vec{I}_M} \gsphere)
				(\angdiff \Fullset^{\vec{I}_{M+1}} f),
				\\
	[\angLap, \Fullset^{\vec{I}}] f
	& = 	\sum_{i_1 + i_2 = 1}
				\sum_{M=1}^{|\vec{I}|}
				\mathop{\sum_{\vec{I}_1 + \cdots + \vec{I}_{M+1} = \vec{I}}}_{|\vec{I}_a| \geq 1 \mbox{\upshape \ for } 1 \leq a \leq M} 
					\label{E:COMMUTINGANGANGLAPANDLIEZ}
					\\
	& \ \ \ \ \ \
				(\ginversesphere)^{M+1}
				\underbrace{
				(\angLie_{\Fullset}^{\vec{I}_1} \gsphere)
				\cdots
				(\angLie_{\Fullset}^{\vec{I}_{M-1}} \gsphere)
				}_{\mbox{absent when $i_1=M=1$}}
				(\angD^{i_1} \angLie_{\Fullset}^{\vec{I}_M} \gsphere)
				(\angD^{i_2+1} \Fullset^{\vec{I}_{M+1}} f).
				\notag
\end{align}
\end{subequations}
In \eqref{E:ANGDANGLIEZELLTUTENSORFIELDCOMMUTATOR}-\eqref{E:COMMUTINGANGANGLAPANDLIEZ}, 
we have omitted all tensorial contractions to condense the presentation.
\end{lemma}

\begin{lemma}\cite{jSgHjLwW2016}*{Lemma 5.2; \textbf{Preliminary Lie derivative commutation identities}}
Let $\vec{I} = (\iota_1,\iota_2,\cdots, \iota_N)$ be an $N^{th}$-order $\Fullset$ multi-index,
	let $f$ be a scalar-valued function, and let
	$\xi$ be a type $\binom{m}{n}$ $\ell_{t,u}$-tangent tensorfield with $m + n \geq 1$.
	Let $i_1,i_2,\cdots,i_N$ be any permutation of
	$1,2,\cdots,N$ and let
	$\vec{I}' = (\iota_{i_1},\iota_{i_2},\cdots, \iota_{i_N})$.
	Then, up to omitted constant factors, we have
	\begin{subequations}
	\begin{align} \label{E:ALGEBRAICSUBTRACTINGTWOPERMUTEDORDERNVECTORFIELDS}
		\left\lbrace
			\Fullset^{\vec{I}}
			- \Fullset^{\vec{I}'}
		\right\rbrace
		f
		& = 	\mathop{\sum_{\vec{I}_1 + \vec{I}_2 + \iota_{k_1} + \iota_{k_2} = \vec{I}}}_{
					Z_{(\iota_{k_1})} \in \lbrace \Lunit, \Rad \rbrace, 
					\ Z_{(\iota_{k_2})} \in \lbrace \Rad, \GeoAng \rbrace, 
					\ Z_{(\iota_{k_1})} \neq Z_{(\iota_{k_2})}}
					\angLie_{\Fullset}^{\vec{I}_1} \angdeformoneformupsharparg{Z_{(\iota_{k_2})}}{Z_{(\iota_{k_1})}} 
					\cdot 
					\angdiff \Fullset^{\vec{I}_2} f,
					\\
		\left\lbrace
			\angLie_{\Fullset}^{\vec{I}}
			- \angLie_{\Fullset}^{\vec{I}'}
		\right\rbrace
		\xi
		& = 
			 	\mathop{\sum_{\vec{I}_1 + \vec{I}_2 + \iota_{k_1} + \iota_{k_2} = \vec{I}}}_{
			 		Z_{(\iota_{k_1})} \in \lbrace \Lunit, \Rad \rbrace, 
					\ Z_{(\iota_{k_2})} \in \lbrace \Rad, \GeoAng \rbrace, 
					\ Z_{(\iota_{k_1})} \neq Z_{(\iota_{k_2})}
			 	}
				\angLie_{\angLie_{\Fullset}^{\vec{I}_1} \angdeformoneformupsharparg{Z_{(\iota_{k_2})}}{Z_{(\iota_{k_1})}}} 
				\angLie_{\Fullset}^{\vec{I}_2} \xi.
				\label{E:TENSORFIELDACTINGALGEBRAICSUBTRACTINGTWOPERMUTEDORDERNVECTORFIELDS}
	\end{align}
	\end{subequations}
	In \eqref{E:ALGEBRAICSUBTRACTINGTWOPERMUTEDORDERNVECTORFIELDS}-\eqref{E:TENSORFIELDACTINGALGEBRAICSUBTRACTINGTWOPERMUTEDORDERNVECTORFIELDS},
	$
	\vec{I}_1 + \vec{I}_2 + \iota_{k_1} + \iota_{k_2} = \vec{I}
	$
	means that 
	$\vec{I}_1 = (\iota_{k_3}, \iota_{k_4}, \cdots, \iota_{k_m})$,
	and
	$\vec{I}_2 = (\iota_{k_{m+1}}, \iota_{k_{m+2}}, \cdots, \iota_{k_N})$,
	where
	$k_1, k_2, \cdots, k_N$ is a permutation of 
	$1,2,\cdots,N$. In particular, 
	$|\vec{I}_1| + |\vec{I}_2| = N-2$.
\end{lemma}

We now provide the main estimates of this section.
\begin{lemma}[\textbf{Commutator estimates}]
	\label{L:COMMUTATORESTIMATES}
	Assume that 
	$1 \leq N \leq 20$
	and
	$0 \leq M \leq \min \lbrace 2, N \rbrace$.
	Let $\vec{I}$
	be a multi-index belonging to the set
	$\mathcal{I}^{N+1;M}$
	from Def.~\ref{D:COMMUTATORMULTIINDICES}
	and let $\vec{I}'$ be any permutation of $\vec{I}$.
	Let $f$ be a scalar-valued function.
	Under the data-size and bootstrap assumptions 
	of Sects.~\ref{SS:FLUIDVARIABLEDATAASSUMPTIONS}-\ref{SS:AUXILIARYBOOTSTRAP}
	and the smallness assumptions of Sect.~\ref{SS:SMALLNESSASSUMPTIONS}, 
	the following commutator estimates hold
	on $\mathcal{M}_{\Tboot,U_0}$
	(see Sect.~\ref{SS:STRINGSOFCOMMUTATIONVECTORFIELDS} regarding the vectorfield operator notation):
	\begin{subequations}
	\begin{align}
		\left|
			\Fullset^{\vec{I}} f
			-
			\Fullset^{\vec{I}'} f
		\right|
		& \lesssim
			\left|
				\myarray
					[\Fullset_{\ast \ast}^{[1,\lceil N/2 \rceil];\leq (M-1)_+} \BadVar]
					{\Fullset_{\ast}^{[1,\lceil N/2 \rceil];\leq (M-1)_+} \GdVar}
			\right|
			\left|
				\Fullset_{\ast \ast}^{[1,N];\leq M} f
			\right|
			+
			\underbrace{
			\left|
				\Fullset_{\ast \ast}^{[1,N];\leq (M-1)_+} f
			\right|
			}_{\mbox{Absent if $M=0$}}
				\label{E:DETAILEDPERMUTEDVECTORFIELDSACTINGONFUNCTIONSCOMMUTATORESTIMATE}  \\
		& \ \
			+
			\left|
				\Fullset_{\ast \ast}^{[1,\lfloor N/2 \rfloor];\leq M} f
			\right|
			\left|
				\myarray
					[\Fullset_{\ast \ast}^{[1,N];\leq (M-1)_+} \BadVar]
					{\Fullset_{\ast}^{[1,N];\leq M} \GdVar}
			\right|,
			\notag
				\\
		\left|
			\Fullset^{\vec{I}} f
			-
			\Fullset^{\vec{I}'} f
		\right|
		& \lesssim
			\left|
				\Fullset_{\ast \ast}^{[1,N];\leq M} f
			\right|
		+
			\left|
				\Fullset_{\ast \ast}^{[1,\lfloor N/2 \rfloor];\leq M} f
			\right|
			\left|
				\myarray
					[\Fullset_{\ast \ast}^{[1,N];\leq (M-1)_+} \BadVar]
					{\Fullset_{\ast}^{[1,N];\leq M} \GdVar}
			\right|,
			\label{E:LESSPRECISEPERMUTEDVECTORFIELDSACTINGONFUNCTIONSCOMMUTATORESTIMATE}
		\end{align}
		 \end{subequations}
		where $(M-1)_+ := \max\lbrace 0, M-1 \rbrace$.

		Moreover, if $1 \leq N \leq 19$
		and $0 \leq M \leq \min \lbrace 2,N \rbrace$,
		then the following commutator estimates hold:
		\begin{subequations}
		\begin{align}
		\left|
			[\angD^2, \angLie_{\Fullset}^{N;M}] f
		\right|
		& \lesssim
			\left|
				\Fullset_{\ast \ast}^{[1,N];\leq M} f
			\right|
			+ 
			\left|
				\Fullset_{\ast \ast}^{[1,\lceil N/2 \rceil];\leq M} f
			\right|
			\left|
				\myarray
					[\Fullset_{\ast \ast}^{[1,N+1];\leq (M-1)_+} \BadVar]
					{\Fullset_{\ast}^{[1,N+1];\leq M} \GdVar}
			\right|,
				\label{E:ANGDSQUAREDNOSPECIALSTRUCTUREFUNCTIONCOMMUTATOR} \\
		\left|
			[\angLap, \Fullset^{N;M}] f
		\right|
		& \lesssim
			\left|
				\Fullset_{\ast \ast}^{[1,N+1];\leq M} f
			\right|
			+ 
			\left|
				\Fullset_{\ast \ast}^{[1,\lceil N/2 \rceil];\leq M} f
			\right|
			\left|
				\myarray
					[\Fullset_{\ast \ast}^{[1,N+1];\leq (M-1)_+} \BadVar]
					{\Fullset_{\ast}^{[1,N+1];\leq M} \GdVar}
			\right|.
			\label{E:ANGLAPNOSPECIALSTRUCTUREFUNCTIONCOMMUTATOR}
		\end{align}
		\end{subequations}

		Finally, if $\xi$ is an
		$\ell_{t,u}$-tangent one-form or a type $\binom{0}{2}$ $\ell_{t,u}$-tangent tensorfield,
		$1 \leq N \leq 19$, 
		$0 \leq M \leq \min \lbrace 2,N \rbrace$,
		and $\vec{I} \in \mathcal{I}^{N+1;M}$,
		then the following commutator estimates hold:
		\begin{subequations}
		\begin{align}
		\left|
			\angLie_{\Fullset}^{\vec{I}} \xi
			-
			\angLie_{\Fullset}^{\vec{I}'} \xi
		\right|
		& \lesssim
			\left|
				\angLie_{\Fullset_*}^{[1,N];\leq M} \xi
			\right|
			+
			\left|
				\angLie_{\Fullset}^{\leq \lceil N/2 \rceil;\leq M} \xi
			\right|
			\left|
				\myarray
					[\Fullset_{\ast \ast}^{[1,N+1];\leq (M-1)_+} \BadVar]
					{\Fullset_{\ast}^{[1,N+1];\leq M} \GdVar}
			\right|,
			\label{E:NOSPECIALSTRUCTURETENSORFIELDCOMMUTATORESTIMATE}
		\\
		\left|
			[\angD, \angLie_{\Fullset}^{N;M}] \xi
		\right|
		& \lesssim
			\left|
				\angLie_{\Fullset_*}^{[1,N-1];\leq M} \xi
			\right|
			+
			\left|
				\angLie_{\Fullset}^{\leq \lceil N/2 \rceil;\leq M} \xi
			\right|
			\left|
				\myarray
					[\Fullset_{\ast \ast}^{[1,N+1];\leq (M-1)_+} \BadVar]
					{\Fullset_{\ast}^{[1,N+1];\leq M} \GdVar}
			\right|,
				\label{E:ANGDNOSPECIALSTRUCTURETENSORFIELDCOMMUTATORESTIMATE} \\
		\left|
			[\angdiv, \angLie_{\Fullset}^{N;\leq M}] \xi
		\right|
		& \lesssim
			\left|
				\angLie_{\Fullset_*}^{[1,N];\leq M} \xi
			\right|
			+
			\left|
				\angLie_{\Fullset}^{\leq \lceil N/2 \rceil;\leq M} \xi
			\right|
			\left|
				\myarray
					[\Fullset_{\ast \ast}^{[1,N+1];\leq (M-1)_+} \BadVar]
					{\Fullset_{\ast}^{[1,N+1];\leq M} \GdVar}
			\right|.
			\label{E:ANGDIVANGLIETANGENTIALTENSORFIELDCOMMUTATORESTIMATE}
		\end{align}
		\end{subequations}


	\end{lemma}

\begin{proof}
	See Sect.~\ref{SS:OFTENUSEDESTIMATES} for some comments on the analysis.

	\noindent \textbf{Proof of \eqref{E:DETAILEDPERMUTEDVECTORFIELDSACTINGONFUNCTIONSCOMMUTATORESTIMATE}}:
	We consider the commutation formula \eqref{E:ALGEBRAICSUBTRACTINGTWOPERMUTEDORDERNVECTORFIELDS}.
	We will bound the products
	$
	\angLie_{\Fullset}^{\vec{I}_1} \angdeformoneformupsharparg{Z_{(\iota_{k_2})}}{Z_{(\iota_{k_1})}} 
	\cdot 
	\angdiff \Fullset^{\vec{I}_2} f
	$
	on RHS~\eqref{E:ALGEBRAICSUBTRACTINGTWOPERMUTEDORDERNVECTORFIELDS}
	on a case by case basis.
	Let $M'$ be the number of factors of $\Rad$ in $\Fullset^{\vec{I}_2}$.
	Note that $M' \leq M$ in view of
	the summation constraint $\vec{I}_1 + \vec{I}_2 + \iota_{k_1} + \iota_{k_2} = \vec{I}$
	on RHS~\eqref{E:ALGEBRAICSUBTRACTINGTWOPERMUTEDORDERNVECTORFIELDS}.

	\textit{Case i): $M'=M$ and $|\vec{I}_2| \in [\lfloor N/2 \rfloor,N-1]$.}
	Clearly we have
	$
	\left|
		\angdiff \Fullset^{\vec{I}_2} f
	\right|
	\lesssim
	\left|
		\Fullset_{\ast \ast}^{[1,N];\leq M} f
	\right|
	$.
	To bound the remaining factor
	$\angLie_{\Fullset}^{\vec{I}_1} \angdeformoneformupsharparg{Z_{(\iota_{k_2})}}{Z_{(\iota_{k_1})}}$,
	where $|\vec{I}_1| \in [0,\lfloor (N-1)/2 \rfloor]$,
	we note that since $M'=M$,
	it must be that $\angLie_{\Fullset}^{\vec{I}_1}$ comprises only 
	$\mathcal{P}_u^t$-tangent vectorfield factors
	and that $(Z_{(\iota_{k_1})}, Z_{(\iota_{k_2})}) = (\Lunit,\GeoAng)$.
	Hence, with the help of \eqref{E:NOSPECIALSTRUCTUREDIFFERENTIATEDGEOANGDEFORMSPHERELSHARPPOINTWISE},
	we see that the remaining factor under consideration
	is bounded in magnitude by 
	$\lesssim 
	\left|
		\angLie_{\Tanset}^{\leq \lfloor (N-1)/2 \rfloor} \angdeformoneformupsharparg{\GeoAng}{\Lunit}
	\right|
	\lesssim
	\left|
			\myarray
			[\Fullset_{\ast \ast}^{[1,\lceil N/2 \rceil];0} \BadVar]
			{\Fullset_{\ast}^{[1,\lceil N/2 \rceil];0} \GdVar}
		\right|
	$.
	In particular, the product 
	$
	\angLie_{\Fullset}^{\vec{I}_1} \angdeformoneformupsharparg{Z_{(\iota_{k_2})}}{Z_{(\iota_{k_1})}} 
	\cdot 
	\angdiff \Fullset^{\vec{I}_2} f
	$
	under consideration is bounded in magnitude by
	$
	\lesssim 
	$
	the first product on
	RHS~\eqref{E:DETAILEDPERMUTEDVECTORFIELDSACTINGONFUNCTIONSCOMMUTATORESTIMATE}.

	\textit{Case ii): $M'=M$ and $|\vec{I}_2| \in [0,\lfloor N/2 \rfloor - 1]$.}
	Clearly we have
	$
	\left|
		\angdiff \Fullset^{\vec{I}_2} f
	\right|
	\lesssim
	\left|
		\Fullset_{\ast \ast}^{[1,\lfloor N/2 \rfloor];\leq M} f
	\right|
	$.
	To bound the remaining factor
	$\angLie_{\Fullset}^{\vec{I}_1} \angdeformoneformupsharparg{Z_{(\iota_{k_2})}}{Z_{(\iota_{k_1})}}$,
	where $|\vec{I}_1| \in [\lceil N/2 \rceil,N-1]$,
	we note that since $M'=M$,
	it must be that $\angLie_{\Fullset}^{\vec{I}_1}$ comprises only 
	$\mathcal{P}_u^t$-tangent vectorfield factors
	and $(Z_{(\iota_{k_1})}, Z_{(\iota_{k_2})}) = (\Lunit,\GeoAng)$.
	Hence, with the help of \eqref{E:NOSPECIALSTRUCTUREDIFFERENTIATEDGEOANGDEFORMSPHERELSHARPPOINTWISE},
	we see that the remaining factor under consideration
	is bounded in magnitude by 
	$\lesssim 
	\left|
		\angLie_{\Tanset}^{\leq N-1} \angdeformoneformupsharparg{\GeoAng}{\Lunit}
	\right|
	\lesssim
	\left|
			\myarray
			[\Fullset_{\ast \ast}^{[1,N];0} \BadVar]
			{\Fullset_{\ast}^{[1,N+1];0} \GdVar}
		\right|
	$.
	In particular, the product 
	$
	\angLie_{\Fullset}^{\vec{I}_1} \angdeformoneformupsharparg{Z_{(\iota_{k_2})}}{Z_{(\iota_{k_1})}} 
	\cdot 
	\angdiff \Fullset^{\vec{I}_2} f
	$
	under consideration is bounded in magnitude by
	$
	\lesssim 
	$
	the last product on
	RHS~\eqref{E:DETAILEDPERMUTEDVECTORFIELDSACTINGONFUNCTIONSCOMMUTATORESTIMATE}.

	We note that we have now proved inequality 
	\eqref{E:DETAILEDPERMUTEDVECTORFIELDSACTINGONFUNCTIONSCOMMUTATORESTIMATE} 
	in the case $M=0$, and it holds without the second term 
	on the RHS (as is indicated in 
	\eqref{E:DETAILEDPERMUTEDVECTORFIELDSACTINGONFUNCTIONSCOMMUTATORESTIMATE}).

	\textit{Case iii): $1 \leq M \leq 2$, $M'\leq M-1$ and $|\vec{I}_2| \in [\lfloor N/2 \rfloor,N-1]$.}
	Clearly we have
	$
	\left|
		\angdiff \Fullset^{\vec{I}_2} f
	\right|
	\lesssim
	\left|
		\Fullset_{\ast \ast}^{[1,N];\leq M-1} f
	\right|
	$.
	Since $|\vec{I}_1| \in [0,\lfloor (N-1)/2 \rfloor]$,
	we may bound the remaining factor
	$\angLie_{\Fullset}^{\vec{I}_1} \angdeformoneformupsharparg{Z_{(\iota_{k_2})}}{Z_{(\iota_{k_1})}}$
	in the norm $\| \cdot \|_{L^{\infty}(\Sigma_t^u)}$ by $\lesssim 1$
	with the help of the pointwise estimates of Lemma~\ref{L:MOREPRECISEDEFTENSESTIMATESPOINTWISE}
	and the bootstrap assumptions.
	We note that in view of 
	the summation constraint $\vec{I}_1 + \vec{I}_2 + \iota_{k_1} + \iota_{k_2} = \vec{I}$
	and the assumption $M \leq 2$,
	we do not encounter terms of the form
	$
	\angLie_{\Fullset}^{|\vec{I}_1|;2}
	\angdeformoneformupsharparg{\Lunit}{\Rad}
	$
	or
	$
	\angLie_{\Fullset}^{|\vec{I}_1|;2}
	\angdeformoneformupsharparg{\GeoAng}{\Rad}
	$,
	which would contain factors involving derivatives of $\Rad \Rad \upmu$ 
	or $\Rad \Rad \Rad \vec{\Psi}$,
	which we are not able to control in $L^{\infty}$ based on the current
	bootstrap assumptions.\footnote{See, however, Sect.~\ref{S:LINFINITYESTIMATESFORHIGHERTRANSVERSAL}.}
	In total, we find that the product 
	$
	\angLie_{\Fullset}^{\vec{I}_1} \angdeformoneformupsharparg{Z_{(\iota_{k_2})}}{Z_{(\iota_{k_1})}} 
	\cdot 
	\angdiff \Fullset^{\vec{I}_2} f
	$
	under consideration is bounded in magnitude by
	$
	\lesssim 
	$
	the second term on
	RHS~\eqref{E:DETAILEDPERMUTEDVECTORFIELDSACTINGONFUNCTIONSCOMMUTATORESTIMATE}.

	\textit{Case iv): $1 \leq M \leq 2$, $M'\leq M-1$, 
	$|\vec{I}_2| \in [0,\lfloor N/2 \rfloor - 1]$, and $1 \leq N \leq 4$.}
	Clearly we have
	$
	\left|
		\angdiff \Fullset^{\vec{I}_2} f
	\right|
	\lesssim
	\left|
		\Fullset_{\ast \ast}^{[1,\lfloor N/2 \rfloor];\leq M-1} f
	\right|
	$.
	Since $1 \leq N \leq 4$,
	we may bound the remaining factor
	$\angLie_{\Fullset}^{\vec{I}_1} \angdeformoneformupsharparg{Z_{(\iota_{k_2})}}{Z_{(\iota_{k_1})}}$
	in the norm $\| \cdot \|_{L^{\infty}(\Sigma_t^u)}$ by $\lesssim 1$,
	as in Case iii).
	It follows that the product
	$
	\angLie_{\Fullset}^{\vec{I}_1} \angdeformoneformupsharparg{Z_{(\iota_{k_2})}}{Z_{(\iota_{k_1})}} 
	\cdot 
	\angdiff \Fullset^{\vec{I}_2} f
	$
	under consideration is bounded in magnitude by
	$
	\lesssim 
	$
	the second term on
	RHS~\eqref{E:DETAILEDPERMUTEDVECTORFIELDSACTINGONFUNCTIONSCOMMUTATORESTIMATE}.

	\textit{Case v): $1 \leq M \leq 2$,
	$M'\leq M-1$, $|\vec{I}_2| \in [0,\lfloor N/2 \rfloor - 1]$, and $5 \leq N \leq 20$.}
	Clearly we have
	$
	\left|
		\angdiff \Fullset^{\vec{I}_2} f
	\right|
	\lesssim
	\left|
		\Fullset_{\ast \ast}^{[1,\lfloor N/2 \rfloor];\leq M-1} f
	\right|
	$.
	We now bound the remaining factor
	$\angLie_{\Fullset}^{\vec{I}_1} \angdeformoneformupsharparg{Z_{(\iota_{k_2})}}{Z_{(\iota_{k_1})}}$,
	starting with the sub-case in which
	either $Z_{(\iota_{k_1})} = \Rad$
	or $Z_{(\iota_{k_2})} = \Rad$.
	In view of the summation constraint 
	$\vec{I}_1 + \vec{I}_2 + \iota_{k_1} + \iota_{k_2} = \vec{I}$,
	we see that it suffices to bound
	$
		\left|
			\angLie_{\Fullset}^{|\vec{I}_1|;\leq M-1}
			\angdeformoneformupsharparg{\GeoAng}{\Rad}
		\right|,
			\,
	\left|
		\angLie_{\Fullset}^{|\vec{I}_1|;\leq M-1}
		\angdeformoneformupsharparg{\Rad}{\Lunit}
	\right|
	$.
	Since $N \geq 5$, and $|\vec{I}_1| \in [\lceil N/2 \rceil,N-1]$, 
	we have $|\vec{I}_1| \geq 3$. Thus,
	since $M \leq 2$,
	at least $2$ vectorfield factors in the operator
	$\angLie_{\Fullset}^{|\vec{I}_1|;\leq M-1}$ must be 
	$\mathcal{P}_u^t$-tangent. That is, 
	$
	\angLie_{\Fullset}^{|\vec{I}_1|;\leq M-1} 
	= 
	\angLie_{\Fullset_{\ast}}^{|\vec{I}_1|;\leq M-1}
	$.
	We may therefore use
	\eqref{E:STARDIFFERENTIATEDGEOANGDEFORMSPHERERADPOINTWISE}
	with $M-1$ in the role of $M$
	to deduce that
	$
		\left|
			\angLie_{\Fullset_{\ast}}^{|\vec{I}_1|;\leq M-1}
			\angdeformoneformupsharparg{\GeoAng}{\Rad}
		\right|,
			\,
	\left|
		\angLie_{\Fullset_{\ast}}^{|\vec{I}_1|;\leq M-1}
		\angdeformoneformupsharparg{\Rad}{\Lunit}
	\right|
	\lesssim
	\left|
			\myarray
			[\Fullset_{\ast \ast}^{[1,N];\leq M-1} \BadVar]
			{\Fullset_{\ast}^{[1,N];\leq M} \GdVar}
	\right|
	$.
	It follows that the product 
	$
	\angLie_{\Fullset}^{\vec{I}_1} \angdeformoneformupsharparg{Z_{(\iota_{k_2})}}{Z_{(\iota_{k_1})}} 
	\cdot 
	\angdiff \Fullset^{\vec{I}_2} f
	$
	under consideration is bounded in magnitude by
	$
	\lesssim 
	$
	the last product on
	RHS~\eqref{E:DETAILEDPERMUTEDVECTORFIELDSACTINGONFUNCTIONSCOMMUTATORESTIMATE} as desired.
	Finally, we address the remaining sub-case in which
	$(Z_{(\iota_{k_1})},Z_{(\iota_{k_2})})
	=
	(\Lunit, \GeoAng)
	$.
	Thus, using \eqref{E:NOSPECIALSTRUCTUREDIFFERENTIATEDGEOANGDEFORMSPHERELSHARPPOINTWISE},
	we see that the factor
	$\angLie_{\Fullset}^{\vec{I}_1} \angdeformoneformupsharparg{Z_{(\iota_{k_2})}}{Z_{(\iota_{k_1})}}
	$
	is bounded in magnitude by
	$
	\lesssim
	\left|
		\angLie_{\Fullset_{\ast}}^{\leq N-1;\leq M}
		\angdeformoneformupsharparg{\GeoAng}{\Lunit}
	\right|
	\lesssim
	\left|
			\myarray
			[\Fullset_{\ast \ast}^{[1,N-1];\leq M-1} \BadVar]
			{\Fullset_{\ast}^{[1,N];\leq M} \GdVar}
	\right|
	$.
	It follows that the product 
	$
	\angLie_{\Fullset}^{\vec{I}_1} \angdeformoneformupsharparg{Z_{(\iota_{k_2})}}{Z_{(\iota_{k_1})}} 
	\cdot 
	\angdiff \Fullset^{\vec{I}_2} f
	$
	under consideration is bounded in magnitude by
	$
	\lesssim 
	$
	the last product on
	RHS~\eqref{E:DETAILEDPERMUTEDVECTORFIELDSACTINGONFUNCTIONSCOMMUTATORESTIMATE} as desired.
	We have thus proved \eqref{E:DETAILEDPERMUTEDVECTORFIELDSACTINGONFUNCTIONSCOMMUTATORESTIMATE}.

	\medskip
	\noindent \textbf{Proof of \eqref{E:LESSPRECISEPERMUTEDVECTORFIELDSACTINGONFUNCTIONSCOMMUTATORESTIMATE}}:
	The estimate \eqref{E:LESSPRECISEPERMUTEDVECTORFIELDSACTINGONFUNCTIONSCOMMUTATORESTIMATE}
	is a simplified version of \eqref{E:DETAILEDPERMUTEDVECTORFIELDSACTINGONFUNCTIONSCOMMUTATORESTIMATE}
	that follows as a simple consequence of
	\eqref{E:DETAILEDPERMUTEDVECTORFIELDSACTINGONFUNCTIONSCOMMUTATORESTIMATE}
	and the bootstrap assumptions.

	\medskip
	\noindent \textbf{Proof of \eqref{E:ANGDSQUAREDNOSPECIALSTRUCTUREFUNCTIONCOMMUTATOR}
	and
	\eqref{E:ANGLAPNOSPECIALSTRUCTUREFUNCTIONCOMMUTATOR}}:
	The proofs of these estimates
	are similar to the proof of \eqref{E:DETAILEDPERMUTEDVECTORFIELDSACTINGONFUNCTIONSCOMMUTATORESTIMATE}
	and are based on the commutation identities 
	\eqref{E:COMMUTINGANGDSQUAREDANDLIEZ}-\eqref{E:COMMUTINGANGANGLAPANDLIEZ},
	the estimates 
	\eqref{E:POINTWISEESTIMATESFORGSPHEREANDITSNOSPECIALSTRUCTUREDERIVATIVES}
	and
	\eqref{E:POINTWISEESTIMATESFORGSPHEREANDITSSTARDERIVATIVES},
	and Lemma~\ref{L:ANGDERIVATIVESINTERMSOFTANGENTIALCOMMUTATOR};
	we omit the details.

	\medskip
	\noindent \textbf{Proof of \eqref{E:NOSPECIALSTRUCTURETENSORFIELDCOMMUTATORESTIMATE},
	\eqref{E:ANGDNOSPECIALSTRUCTURETENSORFIELDCOMMUTATORESTIMATE}
	and \eqref{E:ANGDIVANGLIETANGENTIALTENSORFIELDCOMMUTATORESTIMATE}}:
	The proofs of these estimates 
	are similar to the proof of \eqref{E:DETAILEDPERMUTEDVECTORFIELDSACTINGONFUNCTIONSCOMMUTATORESTIMATE}
	and are based on the commutation identities 
	\eqref{E:ANGDANGLIEZELLTUTENSORFIELDCOMMUTATOR}-\eqref{E:ANGDIVANGLIEZELLTUTENSORFIELDCOMMUTATOR}
	and
	\eqref{E:TENSORFIELDACTINGALGEBRAICSUBTRACTINGTWOPERMUTEDORDERNVECTORFIELDS},
	the estimates 
	\eqref{E:POINTWISEESTIMATESFORGSPHEREANDITSNOSPECIALSTRUCTUREDERIVATIVES}
	and
	\eqref{E:POINTWISEESTIMATESFORGSPHEREANDITSSTARDERIVATIVES},
	and Lemma~\ref{L:ANGLIEPXIINTERMSOFANGLIEGEOANGXI}.
	We omit the details, noting only that
	the right-hand side of
	\eqref{E:NOSPECIALSTRUCTURETENSORFIELDCOMMUTATORESTIMATE}
	involves one more derivative of $\BadVar$ and $\GdVar$
	compared to the estimates \eqref{E:DETAILEDPERMUTEDVECTORFIELDSACTINGONFUNCTIONSCOMMUTATORESTIMATE}-\eqref{E:LESSPRECISEPERMUTEDVECTORFIELDSACTINGONFUNCTIONSCOMMUTATORESTIMATE}; 
	the reason is that we use the estimate
	\eqref{E:ANGLIEPXIINTERMSOFANGLIEGEOANGXI}
	when bounding the terms 
	on RHS~\eqref{E:TENSORFIELDACTINGALGEBRAICSUBTRACTINGTWOPERMUTEDORDERNVECTORFIELDS},
	which leads to the presence of one additional derivative on 
	$\angdeformoneformupsharparg{Z_{(\iota_{k_2})}}{Z_{(\iota_{k_1})}}$.
\end{proof}

\begin{corollary}
	\label{C:BOUNDSFORDERIVATIVESOFLAPLACIAN}
	Assume that 
	$1 \leq N \leq 20$
	and
	$0 \leq M \leq \min \lbrace 2, N \rbrace$.
	Let $\Psi \in \lbrace \Densrenormalized ,v^1,v^2 \rbrace$.
	Under the data-size and bootstrap assumptions 
	of Sects.~\ref{SS:FLUIDVARIABLEDATAASSUMPTIONS}-\ref{SS:AUXILIARYBOOTSTRAP}
	and the smallness assumptions of Sect.~\ref{SS:SMALLNESSASSUMPTIONS}, 
	the following estimates hold
	on $\mathcal{M}_{\Tboot,U_0}$
	(see Sect.~\ref{SS:STRINGSOFCOMMUTATIONVECTORFIELDS} regarding the vectorfield operator notation):
	\begin{align} \label{E:BOUNDSFORDERIVATIVESOFLAPLACIAN}
		\left|
			\Fullset^{N-1;M} \angLap \Psi
		\right|
		& \lesssim 
			\left|
					\Fullset_{\ast \ast}^{[1,N+1];\leq M} \Psi
			\right|
			+
			\left|
				\myarray
					[\Fullset_{\ast \ast}^{[1,N];\leq (M-1)_+} \BadVar]
					{\Fullset_{\ast}^{[1,N];\leq M} \GdVar}
			\right|.
	\end{align}
\end{corollary}
\begin{proof}
	See Sect.~\ref{SS:OFTENUSEDESTIMATES} for some comments on the analysis.
	We decompose
	$
	 \angLap \Psi
	= 
	\angLap \Fullset^{N-1;M} \Psi
	+
	[\Fullset^{N-1;M},\angLap] \Psi
	$.
	The first term in the decomposition is bounded in magnitude by $\lesssim$
	the first term on RHS~\eqref{E:BOUNDSFORDERIVATIVESOFLAPLACIAN}.
	To obtain 
	$\left|
		[\Fullset^{N-1;M},\angLap] \Psi
		\lesssim
		\mbox{RHS~\eqref{E:BOUNDSFORDERIVATIVESOFLAPLACIAN}}
	\right|
	$,
	we use the commutator estimate
	\eqref{E:ANGLAPNOSPECIALSTRUCTUREFUNCTIONCOMMUTATOR} with 
	$f=\Psi$ and
	$N-1$ in the role of $N$
	and the bootstrap assumptions.

\end{proof}

\subsection{Transport inequalities and improvements of the auxiliary bootstrap assumptions}
\label{SS:IMPROVEMENTOFAUX}

In the next proposition, 
we use the previous estimates to derive transport inequalities for the eikonal function quantities and improvements
of the auxiliary bootstrap assumptions. 
The transport inequalities form the starting point for our
derivation of $L^2$ estimates for the below-top-order derivatives of the eikonal function quantities
(see Sect.~\ref{SS:L2FOREIKONALNOMODIFIED}).
In proving the proposition, we must propagate the smallness of the
$\mathring{\upepsilon}$-sized quantities even though
some terms in their evolution equations 
involve the relatively large $\mathring{\updelta}$-sized
quantities.

\begin{proposition}[\textbf{Transport inequalities and improvements of the auxiliary bootstrap assumptions}] 
\label{P:IMPROVEMENTOFAUX}
Under the data-size and bootstrap assumptions 
of Sects.~\ref{SS:FLUIDVARIABLEDATAASSUMPTIONS}-\ref{SS:AUXILIARYBOOTSTRAP}
and the smallness assumptions of Sect.~\ref{SS:SMALLNESSASSUMPTIONS}, 
the following estimates hold
on $\mathcal{M}_{\Tboot,U_0}$
(see Sect.~\ref{SS:STRINGSOFCOMMUTATIONVECTORFIELDS} regarding the vectorfield operator notation).

\medskip
\noindent \underline{\textbf{Transport inequalities for the eikonal function quantities}.}

\medskip

\noindent \textbf{$\bullet$Transport inequalities for} $\upmu$.
	The following pointwise estimate holds:
	\begin{subequations}
	\begin{align} 
		\left|
			\Lunit \upmu
		\right|
		& \lesssim 
			\left|
				\Fullset^{\leq 1} \threePsi
			\right|.
			\label{E:LUNITUPMUPOINTWISE} 
		\end{align}

	Moreover,
	for $1 \leq N \leq 20$ and $0 \leq M \leq \min \lbrace 1, N-1 \rbrace$
	the following pointwise estimates hold:
	\begin{align} \label{E:LLUNITUPMUFULLSETSTARCOMMUTEDPOINTWISEESTIMATE}
		\left|
			\Lunit \Fullset_{\ast}^{N;M} \upmu
		\right|,
			\,
		\left|
			\Fullset_{\ast}^{N;M} \Lunit \upmu
		\right|
		& \lesssim 
			\left|
				\Fullset_{\ast}^{[1,N+1];\leq M+1} \threePsi
			\right|
			+
			\left|
				\myarray
					[\Fullset_{\ast \ast}^{[1,N];\leq M}\BadVar]
					{\Fullset_{\ast}^{[1,N];\leq M} \GdVar}
			\right|.
	\end{align}
	\end{subequations}

	\noindent \textbf{$\bullet$Transport inequalities for} $\Lunit_{(Small)}^i$ and $\mytr \upchi$.
	For $0 \leq N \leq 20$ and $0 \leq M \leq \min \lbrace 2,N \rbrace$, 
	the following pointwise estimates hold:
	\begin{align}
		\left|
			\myarray
				[\Lunit \Fullset^{N;M} \Lunit_{(Small)}^i]
				{\Lunit \Fullset^{N-1;M} \mytr \upchi}
		\right|,
			\,
		\left|
			\myarray
				[\Fullset^{N;M} \Lunit \Lunit_{(Small)}^i]
				{\Fullset^{N-1;M} \Lunit \mytr \upchi}
		\right|
		& \lesssim 
			\left|
				\Fullset_{\ast}^{[1,N+1];\leq M} \threePsi
			\right|
			+
		\underbrace{
		\left|
			\myarray
				[\Fullset_{\ast \ast}^{[1,N];\leq (M-1)_+}\BadVar]
				{\Fullset_{\ast}^{[1,N];\leq M} \GdVar}
		\right|
		}_{\mbox{Absent when $N=0$}}.
				\label{E:LUNITCOMMUTEDLUNITSMALLIPOINTWISE} 
\end{align}

\medskip
\noindent \underline{$L^{\infty}$ \textbf{estimates for} $\threePsi$ \textbf{and the eikonal function quantities}.}

\medskip

\noindent \textbf{$\bullet$$L^{\infty}$ estimates for $\threePsi$}. 
The following estimates hold for $M=1,2$:
\begin{subequations}
\begin{align} \label{E:PSITRANSVERSALLINFINITYBOUNDBOOTSTRAPIMPROVEDSMALL}
		\left\| 
			\Rad^M v^2
		\right\|_{L^{\infty}(\Sigma_t^u)}
	& \leq C \varepsilon,
	\\
	\left\| 
		\Rad^M \Densrenormalized
	\right\|_{L^{\infty}(\Sigma_t^u)}
	& \leq
		\left\| 
			\Rad^M \Densrenormalized
		\right\|_{L^{\infty}(\Sigma_0^u)}
		+ 
		C \varepsilon,
		 \label{E:DENSITYPSITRANSVERSALLINFINITYBOUNDBOOTSTRAPIMPROVEDLARGE} \\
	\left\| 
		\Rad^M v^1
	\right\|_{L^{\infty}(\Sigma_t^u)}
	& \leq
		\left\| 
			\Rad^M v^1
		\right\|_{L^{\infty}(\Sigma_0^u)}
		+ 
		C \varepsilon,
		 \label{E:PSITRANSVERSALLINFINITYBOUNDBOOTSTRAPIMPROVEDLARGE} \\
		\left\| 
			\Fullset_{\ast}^{\leq 12;\leq 2} \threePsi
		\right\|_{L^{\infty}(\Sigma_t^u)}
		& \leq C \varepsilon.
		\label{E:PSIMIXEDUPTOORDERTWOTRANSVERSALIMPROVED}
\end{align}
\end{subequations}
Moreover,
\begin{align}
	\label{E:CRUCIALPSITRANSVERSALLINFINITYBOUNDBOOTSTRAPIMPROVEDSMALL}
	\left\| 
		\Rad (\Densrenormalized - v^1)
	\right\|_{L^{\infty}(\Sigma_t^u)}
	& \leq C \varepsilon.
\end{align}

\noindent \textbf{$\bullet$$L^{\infty}$ estimates for $\upmu$}. 
	The following estimates hold for $M=0,1$:
	\begin{subequations} 
	\begin{align}
	\left\| 
		\Rad^M \upmu
	\right\|_{L^{\infty}(\Sigma_t^u)}
		& \leq
	 	\left\| 
			\Rad^M \upmu
		\right\|_{L^{\infty}(\Sigma_0^u)}
		+
		\TranminusdatasizeWithFactor^{-1} 
		\left\| 
			\Rad^M 
			\left\lbrace
				\vec{G}_{\Lunit \Lunit} \contr \Rad \threePsi
			\right\rbrace
		\right\|_{L^{\infty}(\Sigma_0^u)}
		+ C \varepsilon,
			\label{E:UPTOONETRANSVERSALDERIVATIVEUPMULINFTY} \\
		\left\| 
			\Lunit \Rad^M \upmu
		\right\|_{L^{\infty}(\Sigma_t^u)}
		& 
		= 
		\frac{1}{2}
		\left\| 
			\Rad^M
			\left\lbrace
				\vec{G}_{\Lunit \Lunit} \contr \Rad \threePsi
			\right\rbrace
		\right\|_{L^{\infty}(\Sigma_0^u)}
		+ \mathcal{O}(\varepsilon),
			\label{E:LUNITUPTOONETRANSVERSALUPMULINFINITY} 
			\\
		\left\| 
			\Fullset_{\ast \ast}^{[1,11];\leq 1} \upmu
		\right\|_{L^{\infty}(\Sigma_t^u)}
		& \leq
			C \varepsilon.
			\label{E:ZSTARSTARUPMULINFTY}
	\end{align}
	\end{subequations}

	Moreover, we have
	\begin{align} \label{E:UPMULINFINITYALONGP0TBOOT}
		\left\| 
		\upmu-1
	\right\|_{L^{\infty}\left(\mathcal{P}_0^{\Tboot}\right)}
	& \leq C \varepsilon.
	\end{align}

	\noindent \textbf{$\bullet$$L^{\infty}$ estimates for $\Lunit_{(Small)}^i$ and $\upchi$}.
The following estimates hold for $M=1,2$:
\begin{subequations}
\begin{align}  
	\left\|
		\Fullset_{\ast}^{\leq 11;\leq 2} \Lunit_{(Small)}^i
	\right\|_{L^{\infty}(\Sigma_t^u)}
	& \leq C \varepsilon,
		\label{E:ZSTARLISMALLLISMALLINFTYESTIMATE} \\
	\left\|
		\Rad^M \Lunit_{(Small)}^i
	\right\|_{L^{\infty}(\Sigma_t^u)}
	& \leq
	\left\| 
		\Rad^M \Lunit_{(Small)}^i
	\right\|_{L^{\infty}(\Sigma_0^u)}
	+  C \varepsilon,
	\label{E:LISMALLLUPTOTWORADIALINFINITYESTIMATE}
\end{align}
\end{subequations}

\begin{align}  
		\left\|
			\angLie_{\Fullset}^{\leq 10;\leq 2} \upchi
		\right\|_{L^{\infty}(\Sigma_t^u)},
			\,
		\left\|
			\angLie_{\Fullset}^{\leq 10;\leq 2} \upchi^{\#}
		\right\|_{L^{\infty}(\Sigma_t^u)},
			\,
		\left\|
			\Fullset^{\leq 10;\leq 2} \mytr \upchi
		\right\|_{L^{\infty}(\Sigma_t^u)}
		& \leq C \varepsilon.
		\label{E:TWORADIALCHICOMMUTEDLINFINITY}
\end{align}

\medskip
\noindent \underline{$L^{\infty}$ \textbf{estimates for} $\Vortrenormalized$.}
The following estimates hold:
\begin{align}
\left\| 
	\Fullset^{\leq 12;\leq 2} \Vortrenormalized
\right\|_{L^{\infty}(\Sigma_t^u)}
		& \leq
			C \varepsilon.
			\label{E:VORTICITYUPTOTWOTRANSVERSALLINFTY}
\end{align}
\end{proposition}

\begin{proof}
	See Sect.~\ref{SS:OFTENUSEDESTIMATES} for some comments on the analysis.
	Throughout the proof, we use the phrase ``conditions on the data''
	to refer to the assumptions stated 
	in Sects.~\ref{SS:FLUIDVARIABLEDATAASSUMPTIONS} and \ref{SS:DATAFOREIKONALFUNCTIONQUANTITIES}.
	Also, we often silently use inequality \eqref{E:DATAEPSILONVSBOOTSTRAPEPSILON}.

	\medskip

	\noindent \textbf{Proof of \eqref{E:VORTICITYUPTOTWOTRANSVERSALLINFTY}}:
	The estimate \eqref{E:VORTICITYUPTOTWOTRANSVERSALLINFTY} for
	$
	\left\| 
		\Tanset^{\leq 12} \Vortrenormalized
	\right\|_{L^{\infty}(\Sigma_t^u)}
	$
	is a direct consequence of the bootstrap assumptions.
	To prove \eqref{E:VORTICITYUPTOTWOTRANSVERSALLINFTY} for
	$
	\left\| 
		\Fullset^{\leq 12;1} \Vortrenormalized
	\right\|_{L^{\infty}(\Sigma_t^u)}
	$,
	we first note that by \eqref{E:TRANSPORTVECTORFIELDINTERMSOFLUNITANDRADUNIT},
	equation \eqref{E:RENORMALIZEDVORTICTITYTRANSPORTEQUATION}
	is equivalent to $\Rad \Vortrenormalized = - \upmu \Lunit \Vortrenormalized$.
	Applying $\Tanset^{\leq 11}$ to this equation
	and using the bootstrap assumptions, we deduce that
	$
	\left\| 
		\Tanset^{\leq 11} \Rad \Vortrenormalized
	\right\|_{L^{\infty}(\Sigma_t^u)}
	\lesssim \varepsilon
	$.
	Next, for $2 \leq K \leq 12$,
	we repeatedly use the commutator estimate 
	\eqref{E:LESSPRECISEPERMUTEDVECTORFIELDSACTINGONFUNCTIONSCOMMUTATORESTIMATE}
	and the bootstrap assumptions 
	to deduce 
	\begin{align} \label{E:VORTICITYONETRANSVERSALINDUCTIVECOMMUTATORBOUND}
		\left|
			\Fullset^{\leq K;1} \Vortrenormalized
		\right|
		& \lesssim
		\left|
			\Tanset^{[1,K-1]} \Rad \Vortrenormalized
		\right|
		+
		\left|
			\Tanset^{\leq K-1} \Vortrenormalized
		\right|.
	\end{align}
	We have already shown that the first term on RHS~\eqref{E:VORTICITYONETRANSVERSALINDUCTIVECOMMUTATORBOUND}
	is $\lesssim \varepsilon$, while the bootstrap assumptions
	imply that the second term is $\lesssim \varepsilon$.
	In total, we have shown that we can permute the vectorfield factors
	in $\Tanset^{\leq 11} \Rad \Vortrenormalized$ 
	up to $\mathcal{O}(\varepsilon)$ errors, 
	which yields the desired bound
	$
	\left\| 
		\Fullset^{\leq 12;1} \Vortrenormalized
	\right\|_{L^{\infty}(\Sigma_t^u)}
	\lesssim \varepsilon
	$.
	To prove \eqref{E:VORTICITYUPTOTWOTRANSVERSALLINFTY} for
	$
	\left\| 
		\Fullset^{\leq 12;2} \Vortrenormalized
	\right\|_{L^{\infty}(\Sigma_t^u)}
	$,
	we first apply $\Tanset^{\leq 10} \Rad$ to 
	the equation
	$\Rad \Vortrenormalized = - \upmu \Lunit \Vortrenormalized$
	and use the 
	already proven bound
	$
	\left\| 
		\Fullset^{\leq 12;\leq 1} \Vortrenormalized
	\right\|_{L^{\infty}(\Sigma_t^u)}
	\lesssim \varepsilon
	$,
	and the bootstrap assumptions to deduce
	$
	\left\| 
		\Tanset^{\leq 10} \Rad \Rad \Vortrenormalized
	\right\|_{L^{\infty}(\Sigma_t^u)}
	\lesssim \varepsilon
	$.
	Using this bound, the estimate
	$
	\left\| 
		\Fullset^{\leq 12;\leq 1} \Vortrenormalized
	\right\|_{L^{\infty}(\Sigma_t^u)}
	\lesssim \varepsilon
	$,
	the commutator estimate \eqref{E:LESSPRECISEPERMUTEDVECTORFIELDSACTINGONFUNCTIONSCOMMUTATORESTIMATE}
	with $M=2$,
	and the bootstrap assumptions, we may use an argument similar
	to the one given just below \eqref{E:VORTICITYONETRANSVERSALINDUCTIVECOMMUTATORBOUND}
	in order to permute the vectorfield factors
	in $\Tanset^{\leq 10} \Rad \Rad \Vortrenormalized$ 
	up to $\mathcal{O}(\varepsilon)$ errors. In total, we have shown that
	$
	\left\| 
		\Fullset^{\leq 12;2} \Vortrenormalized
	\right\|_{L^{\infty}(\Sigma_t^u)}
	\lesssim \varepsilon
	$,
	which completes the proof of \eqref{E:VORTICITYUPTOTWOTRANSVERSALLINFTY}.

	\medskip

	\noindent \textbf{Proof of \eqref{E:LUNITUPMUPOINTWISE}  and 
	\eqref{E:LLUNITUPMUFULLSETSTARCOMMUTEDPOINTWISEESTIMATE}}:
	The estimate \eqref{E:LUNITUPMUPOINTWISE} 
	is a simple consequence of the evolution equation
	$
	\Lunit \upmu 
	= \smoothfunction(\BadVar) \Singletan \threePsi
		+ \smoothfunction(\GdVar) \Rad \threePsi
	$
	(see equation \eqref{E:UPMUFIRSTTRANSPORT} and Lemma~\ref{L:SCHEMATICDEPENDENCEOFMANYTENSORFIELDS})
	and the bootstrap assumptions.

	We now prove \eqref{E:LLUNITUPMUFULLSETSTARCOMMUTEDPOINTWISEESTIMATE}.
	We show only how to obtain the estimates for
	$
	\left|
			\Lunit \Fullset_{\ast}^{N;M} \upmu
	\right|
	$
	since the estimates for
	$
		\left|
			\Fullset_{\ast}^{N;M} \Lunit \upmu
		\right|
	$
	are simpler because they do not involve commutations.
	To proceed, for $1 \leq N \leq 20$ and $\min \lbrace 1, N-1 \rbrace$,
	we commute the evolution equation from the previous paragraph 
	with $\Fullset_{\ast}^{N;M}$ to deduce the schematic identity
	\begin{align} \label{E:LUPMUSCHEMATICCOMMUTED}
		\Lunit \Fullset_{\ast}^{N;M} \upmu
		& = [\Lunit, \Fullset_{\ast}^{N;M}] \upmu
			+
			\Fullset_{\ast}^{N;M}
			\left\lbrace
				\smoothfunction(\GdVar) \Rad \Psi
				+ \smoothfunction(\BadVar) \Singletan \threePsi
			\right\rbrace.
	\end{align}
	To bound the magnitude of the second term on RHS \eqref{E:LUPMUSCHEMATICCOMMUTED} by
	$\lesssim$ RHS~\eqref{E:LLUNITUPMUFULLSETSTARCOMMUTEDPOINTWISEESTIMATE}, 
	we use the bootstrap assumptions. To derive
	$
	\left|
		[\Lunit, \Fullset_{\ast}^{N;M}] \upmu
	\right|
	\lesssim \mbox{RHS~\eqref{E:LLUNITUPMUFULLSETSTARCOMMUTEDPOINTWISEESTIMATE}}
	$,
	we use the commutator estimate \eqref{E:LESSPRECISEPERMUTEDVECTORFIELDSACTINGONFUNCTIONSCOMMUTATORESTIMATE}
	with $f = \upmu$ and the bootstrap assumptions.
	We have thus proved \eqref{E:LLUNITUPMUFULLSETSTARCOMMUTEDPOINTWISEESTIMATE}.

	\medskip

	\noindent \textbf{Proof of \eqref{E:LUNITCOMMUTEDLUNITSMALLIPOINTWISE} for} 
	$\Lunit \Fullset^{N;M} \Lunit_{(Small)}^i$ 
	\textbf{and} 
	$\Fullset^{N;M} \Lunit \Lunit_{(Small)}^i$:
	We first write the evolution equation \eqref{E:LLUNITI} 
	in the schematic form
	$\Lunit \Lunit_{(Small)}^i = \smoothfunction(\GdVar,\ginversesphere,\angdiff \vec{x}) \Singletan \threePsi$.
	For $0 \leq N \leq 20$ and $0 \leq M \leq \min \lbrace 2,N \rbrace$, 
	we commute this evolution
	equation with $\Fullset^{N;M}$ to obtain
	\begin{align} \label{E:LLUNITISCHEMATICCOMMUTED}
		\Lunit \Fullset^{N;M} \Lunit_{(Small)}^i
		& = [\Lunit, \Fullset^{N;M}] \Lunit_{(Small)}^i
			+
			\Fullset^{N;M}
			\left\lbrace
				\smoothfunction(\GdVar,\ginversesphere,\angdiff \vec{x}) \Singletan \threePsi
			\right\rbrace.
	\end{align}
	To bound the magnitude of the second term on
	RHS~\eqref{E:LLUNITISCHEMATICCOMMUTED} by
	$\lesssim \mbox{RHS~\eqref{E:LUNITCOMMUTEDLUNITSMALLIPOINTWISE}}$,
	we use the estimates
	\eqref{E:ANGDIFFXI}-\eqref{E:ANGDIFFXNOSPECIALSTRUCTUREDIFFERENTIATEDPOINTWISE}
	and
	\eqref{E:POINTWISEESTIMATESFORGSPHEREANDITSNOSPECIALSTRUCTUREDERIVATIVES}-\eqref{E:POINTWISEESTIMATESFORGSPHEREANDITSSTARDERIVATIVES}
	and the bootstrap assumptions.
	To deduce
	$
	\left|
		[\Lunit, \Fullset^{N;M}] \Lunit_{(Small)}^i
	\right|
	\lesssim \mbox{RHS~\eqref{E:LUNITCOMMUTEDLUNITSMALLIPOINTWISE}}
	$
	we use the commutator estimate \eqref{E:LESSPRECISEPERMUTEDVECTORFIELDSACTINGONFUNCTIONSCOMMUTATORESTIMATE}
	with $f = \Lunit_{(Small)}^i$ and the bootstrap assumptions.

	\medskip
	\noindent \textbf{Proof of \eqref{E:LUNITCOMMUTEDLUNITSMALLIPOINTWISE} for} 
	$\Lunit \Fullset^{N-1;M} \mytr \upchi$ \textbf{and} 
	$\Fullset^{N-1;M} \Lunit \mytr \upchi$:
	We first apply $\Lunit$ to equation \eqref{E:TRCHIINTERMSOFOTHERVARIABLES}
	and use the schematic identity 
	$\angLie_{\Lunit} \ginversesphere 
	= (\ginversesphere)^{-2} \upchi
	= \smoothfunction(\GdVar,\ginversesphere,\angdiff \vec{x}) 
	\Singletan \GdVar
	$
	(see \eqref{E:CONNECTIONBETWEENANGLIEOFGSPHEREANDDEFORMATIONTENSORS},
	\eqref{E:LUNITDEFORMSPHERE},
	and Lemma~\ref{L:SCHEMATICDEPENDENCEOFMANYTENSORFIELDS})
	to deduce that
	$\Lunit \mytr \upchi 
	= 
	\smoothfunction(\GdVar,\ginversesphere,\angdiff \vec{x}) \Singletan \Lunit \GdVar
	+ l.o.t.$,
	where
	$l.o.t.
	:=
	\left\lbrace
			\smoothfunction(\Tanset^{\leq 1} \GdVar, 
			\angLie_{\Tanset}^{\leq 1} \ginversesphere,\angdiff \Tanset^{\leq 1} \vec{x}) 
	\right\rbrace
	\Singletan \GdVar
	$.
	Applying
	$\Fullset^{N-1;M}$ to this identity 
	and using Lemmas~\ref{L:POINTWISEFORRECTANGULARCOMPONENTSOFVECTORFIELDS}
	and 
	\ref{L:POINTWISEESTIMATESFORGSPHEREANDITSDERIVATIVES}
	and the bootstrap assumptions,
	we find that
	$
	\left|
		\Fullset^{N-1;M} \Lunit \mytr \upchi
	\right|
	\lesssim
	\sum_{i=1}^2
	\left|
		\Fullset_*^{N+1;M} \Lunit_{(Small)}^i
	\right|
	+
	\mbox{{\upshape RHS}~\eqref{E:LUNITCOMMUTEDLUNITSMALLIPOINTWISE}}
	$,
	where the operator $\Fullset_*^{N+1;M}$ acting on $\Lunit_{(Small)}^i$
	contains a factor of $\Lunit$.
	By arguing as in our proof of the bound for the commutator term
	on RHS~\eqref{E:LLUNITISCHEMATICCOMMUTED}, we may commute the factor of 
	$\Lunit$ to the front (so that $\Lunit$ acts last), 
	thereby obtaining that
	$
	\sum_{i=1}^2
	\left|
		\Fullset_*^{N+1;M} \Lunit_{(Small)}^i
	\right|
	\lesssim
	\left|
		\Lunit \Fullset^{N;M} \Lunit_{(Small)}^i
	\right|
		+
	\mbox{{\upshape RHS}~\eqref{E:LUNITCOMMUTEDLUNITSMALLIPOINTWISE}}
	$.
	Moreover, we already showed in the previous paragraph that
	$
	\left|
		\Lunit \Fullset^{N;M} \Lunit_{(Small)}^i
	\right|
	\lesssim 
	\mbox{\upshape RHS}~\eqref{E:LUNITCOMMUTEDLUNITSMALLIPOINTWISE}
	$,
	which completes our proof of the estimate \eqref{E:LUNITCOMMUTEDLUNITSMALLIPOINTWISE}
	for
	$
	\left|
		\Fullset^{N-1;M} \Lunit \mytr \upchi
	\right|
	$.
	To obtain the same estimate for
	$
	\left|
		\Lunit \Fullset^{N-1;M} \mytr \upchi
	\right|
	$,
	we use the commutator estimate 
	\eqref{E:LESSPRECISEPERMUTEDVECTORFIELDSACTINGONFUNCTIONSCOMMUTATORESTIMATE} with $f = \mytr \upchi$,
	\eqref{E:POINTWISEESTIMATESFORCHIANDITSDERIVATIVES},
	and the bootstrap assumptions
	to deduce that
	$\left|
		\Lunit \Fullset^{N-1;M} \mytr \upchi
	\right|
	\lesssim 
	\left|
		 \Fullset^{N-1;M} \Lunit \mytr \upchi
	\right|
	+
	\left|
		\myarray
		[\Fullset_{\ast \ast}^{[1,N];\leq (M-1)_+} \BadVar]
		{\Fullset_{\ast}^{[1,N];\leq M} \GdVar}
	\right|
	$.
	The desired bound \eqref{E:LUNITCOMMUTEDLUNITSMALLIPOINTWISE} for
	$
	\left|
		\Lunit \Fullset^{N-1;M} \mytr \upchi
	\right|
	$
	now follows from this estimate and the
	one we established just above for
	$
	\left|
		\Fullset^{N-1;M} \Lunit \mytr \upchi
	\right|
	$.

	\medskip

	\noindent \textbf{Proof of an intermediate estimate:}
As an intermediate step, we now show that
\begin{align} \label{E:EIKONALFUNCTIONQUANTITIESGRONWALLED}
	\left|
		\myarray
			[\Fullset_{\ast \ast}^{[1,11];0} \upmu]
			{\Fullset_{\ast}^{\leq 11;\leq 1} (\Lunit_{(Small)}^1,\Lunit_{(Small)}^2)}
	\right|(t,u,\vartheta)
	& \lesssim \varepsilon.
\end{align} 
	We first recall that 
	$
	\displaystyle
	\Lunit = \frac{\partial}{\partial t}
	$.
	Hence, we can use 
	\eqref{E:LLUNITUPMUFULLSETSTARCOMMUTEDPOINTWISEESTIMATE}-\eqref{E:LUNITCOMMUTEDLUNITSMALLIPOINTWISE}
	and the bootstrap assumptions
	and integrate along the integral curves of $\Lunit$ 
	as in \eqref{E:SIMPLEFTCID}
	to deduce
	\begin{align} \label{E:EIKONALFUNCTIONQUANTITIESGRONWALLREADY}
	&
	\left|
		\myarray
			[\Fullset_{\ast \ast}^{[1,11];\leq 0} \upmu]
			{\Fullset_{\ast}^{\leq 11;\leq 1} (\Lunit_{(Small)}^1,\Lunit_{(Small)}^2)}
	\right|(t,u,\vartheta)
		\\
	& \leq
	C
	\left|
		\myarray
			[\Fullset_{\ast \ast}^{[1,11];0} \upmu]
			{\Fullset_{\ast}^{\leq 11;\leq 1} (\Lunit_{(Small)}^1,\Lunit_{(Small)}^2)}
	\right|(0,u,\vartheta)
		\notag \\
& \ \
	+
	C
	\int_{s=0}^t
		\left|
		\myarray
			[\Fullset_{\ast \ast}^{[1,11];0} \upmu]
			{\Fullset_{\ast}^{\leq 11;\leq 1} (\Lunit_{(Small)}^1,\Lunit_{(Small)}^2)}
		\right|(s,u,\vartheta)
	\, ds
	+ 
	C \varepsilon,
	\notag
\end{align}
where the $C \varepsilon$ term on RHS~\eqref{E:EIKONALFUNCTIONQUANTITIESGRONWALLREADY}
comes from the terms
$\left| \Fullset_{\ast}^{[1,12];\leq 1} \threePsi \right|$
(which are $\lesssim \varepsilon$ in view of the bootstrap assumptions)
on RHS~\eqref{E:LLUNITUPMUFULLSETSTARCOMMUTEDPOINTWISEESTIMATE}
and RHS~\eqref{E:LUNITCOMMUTEDLUNITSMALLIPOINTWISE}.
The conditions on the data imply that the first term on 
RHS~\eqref{E:EIKONALFUNCTIONQUANTITIESGRONWALLREADY}
is $\leq C \varepsilon$. Hence, from Gronwall's inequality, we conclude
$
\left|
		\myarray
			[\Fullset_{\ast \ast}^{[1,11];0} \upmu]
			{\Fullset_{\ast}^{\leq 11;\leq 1} (\Lunit_{(Small)}^1,\Lunit_{(Small)}^2)}
	\right|(t,u,\vartheta)
	\lesssim \varepsilon \exp(C \TranminusdatasizeWithFactor^{-1})
	\lesssim \varepsilon
$,
which yields \eqref{E:EIKONALFUNCTIONQUANTITIESGRONWALLED}.

\medskip

	\noindent \textbf{Proof of \eqref{E:PSIMIXEDUPTOORDERTWOTRANSVERSALIMPROVED}}:
		Using \eqref{E:LONOUTSIDEGEOMETRICWAVEOPERATORFRAMEDECOMPOSED},
		\eqref{E:TRANSPORTVECTORFIELDINTERMSOFLUNITANDRADUNIT},
		Cor.~\ref{C:VELOCITYWAVEEQUATIONDERIVATIVEOFVORTICITYINHOMOGENEOUSTERMEXPRESSION},
		Lemma~\ref{L:SCHEMATICDEPENDENCEOFMANYTENSORFIELDS},
		the identity \eqref{E:ANGLAPINTERMSOFGEOANGDERIVATIVES},
		and the schematic relation
	$
	\Lunit \upmu 
	= \smoothfunction(\BadVar) \Singletan \Psi
		+ \smoothfunction(\GdVar) \Rad \Psi
	$
	(which follows from Lemmas~\ref{E:UPMUFIRSTTRANSPORT} and  
	\ref{L:SCHEMATICDEPENDENCEOFMANYTENSORFIELDS}),
		we write the wave equations 
		\eqref{E:VELOCITYWAVEEQUATION}-\eqref{E:RENORMALIZEDDENSITYWAVEEQUATION}
		verified by $\Psi \in \lbrace \Densrenormalized ,v^1,v^2 \rbrace$
		in the following schematic form:
	\begin{align} \label{E:TRANSRENORMALIZEDAPPLIEDWAVEEQUATIONTRANSPORTINTERPRETATION}
		\Lunit \Rad \Psi 
		& 
		=   \smoothfunction(\BadVar,\ginversesphere,\angdiff \vec{x},\Fullset^{\leq 1} \threePsi) 
				\Tanset^{\leq 2} \Psi
			+ \smoothfunction(\BadVar,\ginversesphere,\angdiff \vec{x},\Fullset^{\leq 1} \threePsi)
				\Singletan \GdVar
					\\
		& \ \
			+ \smoothfunction(\BadVar,\Fullset^{\leq 1} \threePsi) \Tanset^{\leq 1} \Vortrenormalized.
			\notag
	\end{align}
	We now show that
	\begin{align} \label{E:LUNITAPPLIEDTOONCERADIALCOMMUTEDWAVEPOINTWISEGRONWALLREADY}
		\left|
			\Lunit \Fullset_{\ast}^{[1,11];1} \Rad \Psi 
		\right|
		\lesssim 
		\left|
			\Fullset_{\ast}^{[1,11];1} \Rad \Psi  
		\right|
		+
		\varepsilon.
	\end{align}
	To derive \eqref{E:LUNITAPPLIEDTOONCERADIALCOMMUTEDWAVEPOINTWISEGRONWALLREADY},
	we first apply $\Fullset_{\ast}^{[1,11];1}$
	to \eqref{E:TRANSRENORMALIZEDAPPLIEDWAVEEQUATIONTRANSPORTINTERPRETATION}.
	Using the bootstrap assumptions and the already proven estimates
	\eqref{E:VORTICITYUPTOTWOTRANSVERSALLINFTY} and
	\eqref{E:EIKONALFUNCTIONQUANTITIESGRONWALLED} 
	(to bound the derivatives of the terms $\Lunit_{(Small)}^1$ and $\Lunit_{(Small)}^2$
	found in the factor $\Singletan \GdVar$ on 
	RHS~\eqref{E:TRANSRENORMALIZEDAPPLIEDWAVEEQUATIONTRANSPORTINTERPRETATION}),
	we deduce the estimate
	$
		\left\| 
			\Fullset_{\ast}^{[1,11];1} \mbox{RHS \eqref{E:TRANSRENORMALIZEDAPPLIEDWAVEEQUATIONTRANSPORTINTERPRETATION}} 
		\right\|_{L^{\infty}(\Sigma_t^u)}
		\lesssim \varepsilon
	$.
	To finish the proof of \eqref{E:LUNITAPPLIEDTOONCERADIALCOMMUTEDWAVEPOINTWISEGRONWALLREADY},
	we consider the term $\Fullset_{\ast}^{[1,11];1} \Lunit \Rad \Psi$
	obtained by applying $\Fullset_{\ast}^{[1,11];1}$
	to LHS~\eqref{E:TRANSRENORMALIZEDAPPLIEDWAVEEQUATIONTRANSPORTINTERPRETATION}.
	We use the commutator estimate
	\eqref{E:LESSPRECISEPERMUTEDVECTORFIELDSACTINGONFUNCTIONSCOMMUTATORESTIMATE} 
	with $f = \Rad \Psi$, $1 \leq N \leq 11$, and $M=1$
	and the bootstrap assumptions
	to arbitrarily permute the vectorfield factors in 
	$\Fullset_{\ast}^{[1,11];1} \Lunit$
	up to error terms that are
	$\lesssim
	\left|
		\Fullset_{\ast \ast}^{[1,11];\leq 1} \Rad \Psi
	\right|
	\lesssim
	\left|
		\Fullset_{\ast \ast}^{[1,11];1} \Rad \Psi
	\right|
	+
	\varepsilon
	$,
	where we bounded the 
	last factor on RHS~\eqref{E:LESSPRECISEPERMUTEDVECTORFIELDSACTINGONFUNCTIONSCOMMUTATORESTIMATE}  
	as follows:
	$
	\left|
		\myarray
		[\Fullset_{\ast \ast}^{[1,11];0} \BadVar]
		{\Fullset_{\ast}^{[1,11];\leq 1} \GdVar}
	\right|
			\lesssim 1
	$.
	We have therefore proved \eqref{E:LUNITAPPLIEDTOONCERADIALCOMMUTEDWAVEPOINTWISEGRONWALLREADY}.
	We now integrate inequality \eqref{E:LUNITAPPLIEDTOONCERADIALCOMMUTEDWAVEPOINTWISEGRONWALLREADY}
	along the integral curves of $\Lunit$
	as in \eqref{E:SIMPLEFTCID},
	use the conditions on the data, 
	and apply Gronwall's inequality to 
	deduce
	\begin{align} \label{E:GRONWALLEDPOINTWISEESTIMATEPSIUPTOTWOTRANSVERALWITHONEALLTHEWAYIN}
		\left|
			\Fullset_{\ast}^{[1,11];1} \Rad \Psi
		\right|
		& \lesssim \varepsilon.
	\end{align}
	Using \eqref{E:GRONWALLEDPOINTWISEESTIMATEPSIUPTOTWOTRANSVERALWITHONEALLTHEWAYIN}, 
	the commutator estimate \eqref{E:LESSPRECISEPERMUTEDVECTORFIELDSACTINGONFUNCTIONSCOMMUTATORESTIMATE}
	with $M=2$,
	and the bootstrap assumptions
	(including 
	$
	\left\| 
		\Fullset_{\ast}^{[1,12];\leq 1} \Psi
	\right\|_{L^{\infty}(\Sigma_t^u)}
	\lesssim \varepsilon
	$),
	we use a commutator argument similar
	to the one surrounding equation \eqref{E:VORTICITYONETRANSVERSALINDUCTIVECOMMUTATORBOUND},
	which allows us to arbitrarily permute the vectorfield factors
	in the operator $\Fullset_{\ast}^{[1,11];1} \Rad$ 
	on LHS \eqref{E:GRONWALLEDPOINTWISEESTIMATEPSIUPTOTWOTRANSVERALWITHONEALLTHEWAYIN}
	up to $\mathcal{O}(\varepsilon)$ errors. 
	In total, we have obtained
	the desired bound
	$
	\left\| 
		\Fullset_{\ast}^{[1,12];\leq 2} \Psi
	\right\|_{L^{\infty}(\Sigma_t^u)}
	\lesssim \varepsilon
	$,
	which completes the proof of \eqref{E:PSIMIXEDUPTOORDERTWOTRANSVERSALIMPROVED}.

	\medskip

	\noindent \textbf{Proof of \eqref{E:ZSTARSTARUPMULINFTY}, 
	\eqref{E:ZSTARLISMALLLISMALLINFTYESTIMATE},
	and \eqref{E:TWORADIALCHICOMMUTEDLINFINITY}}:
	We first prove \eqref{E:ZSTARSTARUPMULINFTY}
	and \eqref{E:ZSTARLISMALLLISMALLINFTYESTIMATE}.
	Much like in our proof of \eqref{E:EIKONALFUNCTIONQUANTITIESGRONWALLED},
	we may use 
	\eqref{E:LLUNITUPMUFULLSETSTARCOMMUTEDPOINTWISEESTIMATE}-\eqref{E:LUNITCOMMUTEDLUNITSMALLIPOINTWISE}
	and the bootstrap assumptions
	and integrate along the integral curves of $\Lunit$ to deduce
	\begin{align} \label{E:UPTOTWOTRANSVERSALEIKONALFUNCTIONQUANTITIESGRONWALLREADY}
	&
	\left|
		\myarray
			[\Fullset_{\ast \ast}^{[1,11];\leq 1} \upmu]
			{\Fullset_{\ast}^{\leq 11;\leq 2} (\Lunit_{(Small)}^1,\Lunit_{(Small)}^2)}
	\right|(t,u,\vartheta)
		\\
	& \leq
	C
	\left|
		\myarray
			[\Fullset_{\ast \ast}^{[1,11];\leq 1} \upmu]
			{\Fullset_{\ast}^{\leq 11;\leq 2} (\Lunit_{(Small)}^1,\Lunit_{(Small)}^2)}
	\right|(0,u,\vartheta)
		\notag \\
& \ \
	+
	C
	\int_{s=0}^t
		\left|
		\myarray
			[\Fullset_{\ast \ast}^{[1,11];\leq 1} \upmu]
			{\Fullset_{\ast}^{\leq 11;\leq 2} (\Lunit_{(Small)}^1,\Lunit_{(Small)}^2)}
		\right|(s,u,\vartheta)
	\, ds
	+ 
	C \varepsilon,
	\notag
\end{align}
where the $C \varepsilon$ term on RHS~\eqref{E:UPTOTWOTRANSVERSALEIKONALFUNCTIONQUANTITIESGRONWALLREADY}
comes from the terms
$\left| \Fullset_{\ast}^{[1,12];\leq 2} \threePsi \right|$
(which are $\lesssim \varepsilon$ in view of the already proven estimate 
\eqref{E:PSIMIXEDUPTOORDERTWOTRANSVERSALIMPROVED})
on RHS~\eqref{E:LLUNITUPMUFULLSETSTARCOMMUTEDPOINTWISEESTIMATE}
and RHS~\eqref{E:LUNITCOMMUTEDLUNITSMALLIPOINTWISE}.
The conditions on the data imply that the first term 
on RHS~\eqref{E:UPTOTWOTRANSVERSALEIKONALFUNCTIONQUANTITIESGRONWALLREADY}
is $\leq C \varepsilon$. Hence, from Gronwall's inequality, we conclude
$
\left|
		\myarray
			[\Fullset_{\ast \ast}^{[1,11];\leq 1} \upmu]
			{\Fullset_{\ast}^{\leq 11;\leq 2} (\Lunit_{(Small)}^1,\Lunit_{(Small)}^2)}
	\right|(t,u,\vartheta)
	\lesssim 
	\varepsilon \exp(C \TranminusdatasizeWithFactor^{-1})
	\lesssim \varepsilon
$,
which yields 
\eqref{E:ZSTARSTARUPMULINFTY} and \eqref{E:ZSTARLISMALLLISMALLINFTYESTIMATE}.
The estimate \eqref{E:TWORADIALCHICOMMUTEDLINFINITY}
then follows as a consequence of
inequality \eqref{E:POINTWISEESTIMATESFORCHIANDITSDERIVATIVES}
and the estimates
\eqref{E:PSIMIXEDUPTOORDERTWOTRANSVERSALIMPROVED},
\eqref{E:ZSTARSTARUPMULINFTY},
and
\eqref{E:ZSTARLISMALLLISMALLINFTYESTIMATE}.

\medskip
\noindent \textbf{Proof of \eqref{E:UPMULINFINITYALONGP0TBOOT}}:
The estimate \eqref{E:UPMULINFINITYALONGP0TBOOT}
is a trivial consequence of
\eqref{E:SIMPLEUPMUALONGPOESTIMATE}
and \eqref{E:DATAEPSILONVSBOOTSTRAPEPSILON}.

\medskip

\noindent \textbf{Proof of 
	\eqref{E:PSITRANSVERSALLINFINITYBOUNDBOOTSTRAPIMPROVEDSMALL},
	\eqref{E:DENSITYPSITRANSVERSALLINFINITYBOUNDBOOTSTRAPIMPROVEDLARGE},
	\eqref{E:PSITRANSVERSALLINFINITYBOUNDBOOTSTRAPIMPROVEDLARGE},
	\eqref{E:CRUCIALPSITRANSVERSALLINFINITYBOUNDBOOTSTRAPIMPROVEDSMALL},
	\eqref{E:UPTOONETRANSVERSALDERIVATIVEUPMULINFTY},
	\eqref{E:LUNITUPTOONETRANSVERSALUPMULINFINITY},
	and
	\eqref{E:LISMALLLUPTOTWORADIALINFINITYESTIMATE}}:
	We first prove
	\eqref{E:PSITRANSVERSALLINFINITYBOUNDBOOTSTRAPIMPROVEDSMALL},
	\eqref{E:DENSITYPSITRANSVERSALLINFINITYBOUNDBOOTSTRAPIMPROVEDLARGE},
	\eqref{E:PSITRANSVERSALLINFINITYBOUNDBOOTSTRAPIMPROVEDLARGE},
	and
	\eqref{E:CRUCIALPSITRANSVERSALLINFINITYBOUNDBOOTSTRAPIMPROVEDSMALL}.
	We note that a special case of \eqref{E:PSIMIXEDUPTOORDERTWOTRANSVERSALIMPROVED}
	is the estimate 
	$
	\left\|
		\Lunit \Rad^M \threePsi
	\right\|_{L^{\infty}(\Sigma_t^u)}
	\lesssim \varepsilon
	$,
	valid for $M=1,2$.
	Hence, we can integrate along the integral curves of $\Lunit$
	as in \eqref{E:SIMPLEFTCID} and use this 
	estimate to obtain
	$
	\left|
		\Rad^M \threePsi
	\right|(t,u,\vartheta)
	=
	\left|
		\Rad^M \threePsi
	\right|(0,u,\vartheta)
	+ \mathcal{O}(\varepsilon)
	$
	and
	$
	\left|
		\Rad (\Densrenormalized - v^1)
	\right|(t,u,\vartheta)
	=
	\left|
		\Rad (\Densrenormalized - v^1)
	\right|(0,u,\vartheta)
	+ \mathcal{O}(\varepsilon)
	$.
	Using the conditions on the data,
	we arrive at the desired four estimates
	(note that our assumptions on the data and \eqref{E:DATAEPSILONVSBOOTSTRAPEPSILON}
	imply the smallness bounds
	$
	\left\| 
			\Rad (\Densrenormalized - v^1)
	\right\|_{L^{\infty}(\Sigma_0^1)},
		\,
	\left\| 
			\Rad^{[0,2]} v^2
	\right\|_{L^{\infty}(\Sigma_0^1)}
	\lesssim \varepsilon
	$
	and the non-small bounds
	$
	\left\| 
			\Rad^{[0,2]} \Densrenormalized 
	\right\|_{L^{\infty}(\Sigma_0^1)}
	\lesssim 1
	$
	and
	$
	\left\| 
			\Rad^{[0,2]} v^1
	\right\|_{L^{\infty}(\Sigma_0^1)}
	\lesssim 1
	$).

	Inequality \eqref{E:LISMALLLUPTOTWORADIALINFINITYESTIMATE}
	follows in a similar fashion from 
	the already proven estimate \eqref{E:ZSTARLISMALLLISMALLINFTYESTIMATE}.

	We now prove
	\eqref{E:LUNITUPTOONETRANSVERSALUPMULINFINITY}.
	From the evolution equation 
	\eqref{E:UPMUFIRSTTRANSPORT},
	Lemma~\ref{L:SCHEMATICDEPENDENCEOFMANYTENSORFIELDS},
	the commutator estimate
	\eqref{E:LESSPRECISEPERMUTEDVECTORFIELDSACTINGONFUNCTIONSCOMMUTATORESTIMATE} with $f = \upmu$,
	the estimates 
	\eqref{E:PSIMIXEDUPTOORDERTWOTRANSVERSALIMPROVED}
	\eqref{E:ZSTARLISMALLLISMALLINFTYESTIMATE},
	and \eqref{E:ZSTARSTARUPMULINFTY}
	and the bootstrap assumptions,
	we see that for $M=0,1$, we have
	\begin{align} \label{E:LUNITUPMUMAINTERM}
	\Lunit \Rad^M \upmu
	& = 
		\frac{1}{2} 
		\Rad^M 
		\left\lbrace 
			\vec{G}_{\Lunit \Lunit} \contr \Rad \threePsi
		\right\rbrace
		+ 
		\Rad^M \left\lbrace \smoothfunction(\BadVar) \Singletan \threePsi \right\rbrace
		+ 
		[\Lunit,\Rad^M] \upmu
		\\
	& = 
		\frac{1}{2}
		\Rad^M 
		\left\lbrace 
			\vec{G}_{\Lunit \Lunit} \contr \Rad \threePsi
		\right\rbrace
		+ 
		\mathcal{O}(\varepsilon).
		\notag
	\end{align}
	Moreover, 
	from the schematic identity \eqref{E:GFRAMESCALARSDEPENDINGONGOODVARIABLES},
	the estimates 
	\eqref{E:PSIMIXEDUPTOORDERTWOTRANSVERSALIMPROVED}
	and
	\eqref{E:ZSTARLISMALLLISMALLINFTYESTIMATE}
	and the bootstrap assumptions,
	we deduce
	$
	\Lunit \Rad^M 
	\left\lbrace 
		\vec{G}_{\Lunit \Lunit} \contr \Rad \threePsi
	\right\rbrace
	= \mathcal{O}(\varepsilon)
	$.
	Integrating this estimate
	along the integral curves of $\Lunit$
	as in \eqref{E:SIMPLEFTCID},
	we find that
	$
	\left\| 
		\Rad^M 
		\left\lbrace 
			\vec{G}_{\Lunit \Lunit} \contr \Rad \threePsi
		\right\rbrace
	\right\|_{L^{\infty}(\Sigma_t^u)}
	=
	\left\| 
		\Rad^M 
		\left\lbrace 
			\vec{G}_{\Lunit \Lunit} \contr \Rad \threePsi
		\right\rbrace
	\right\|_{L^{\infty}(\Sigma_0^u)}
	+ \mathcal{O}(\varepsilon)
	$.
	From this estimate and \eqref{E:LUNITUPMUMAINTERM},
	we conclude \eqref{E:LUNITUPTOONETRANSVERSALUPMULINFINITY}.

	Finally, we prove \eqref{E:UPTOONETRANSVERSALDERIVATIVEUPMULINFTY}.
	We first integrate along the integral curves of $\Lunit$ as in \eqref{E:SIMPLEFTCID}
	to deduce the following inequality, valid for $M=0,1$:
	$
	\left|
		\Rad^M \upmu
	\right|(t,u,\vartheta)
	\leq
	\left|
		\Rad^M \upmu
	\right|(0,u,\vartheta)
	+
	\int_{s=0}^t
		\left|
			\Lunit \Rad^M \upmu
		\right|(s,u,\vartheta)
	\, ds
	$.
	We now use \eqref{E:LUNITUPTOONETRANSVERSALUPMULINFINITY} to bound
	the time integral in the previous inequality
	by 
	$\leq \TranminusdatasizeWithFactor^{-1} \left\| 
		\Rad^M 
		\left\lbrace 
			\vec{G}_{\Lunit \Lunit} \contr \Rad \threePsi
		\right\rbrace
	\right\|_{L^\infty(\Sigma_0^u)}
	+ C \varepsilon
	$,
	where we have used the assumption 
	$t \leq 2 \TranminusdatasizeWithFactor^{-1}$.
	The desired bound \eqref{E:UPTOONETRANSVERSALDERIVATIVEUPMULINFTY}
	now readily follows from these estimates.
\end{proof}

\section{\texorpdfstring{$L^{\infty}$}{Essential Sup-Norm} Estimates Involving Higher Transversal Derivatives}
\label{S:LINFINITYESTIMATESFORHIGHERTRANSVERSAL}
In Sect.~\ref{S:SHARPESTIMATESFORUPMU},
we derive sharp pointwise estimates for $\upmu$
and some of its derivative, estimates which play a
crucial role in the energy estimates.
The proofs of some of the estimates of Sect.~\ref{S:SHARPESTIMATESFORUPMU} rely on the bound
$\left\|
	\Rad \Rad \upmu
\right\|_{L^{\infty}(\Sigma_t^u)}
\lesssim 1
$.
In this section, we derive this bound
and some related ones,
some of which are needed to prove it.

\subsection{Auxiliary bootstrap assumptions}
\label{SS:BOOTSTRAPFORHIGHERTRANSVERSAL}
We will use auxiliary bootstrap assumptions
to simplify the analysis.
In Prop.~\ref{P:IMPROVEMENTOFHIGHERTRANSVERSALBOOTSTRAP}, 
we derive strict improvements of the assumptions.
\medskip

Our auxiliary bootstrap assumptions are that following inequalities hold on $\mathcal{M}_{\Tboot,U_0}$,
where $\varepsilon$ is the small positive bootstrap parameter from 
Sect.~\ref{SS:PSIBOOTSTRAP}.

\noindent \underline{\textbf{Auxiliary bootstrap assumptions involving three transversal derivatives of}
$\threePsi$.}
For $\Psi \in \lbrace \Densrenormalized ,v^1,v^2 \rbrace$, we have
\begin{align}
	\left\| 
		\Rad \Rad \Rad \Psi
	\right\|_{L^{\infty}(\Sigma_t^u)}
	& \leq 
		\left\| 
			\Rad \Rad \Rad \Psi
		\right\|_{L^{\infty}(\Sigma_0^u)}
		+
		\varepsilon^{1/2}.
		\tag{$\mathbf{AUX}\Rad \Rad \Rad \threePsi$}
	\label{E:WAVEVARIABLEBOOTSTRAPTRANSVERALESTIMATESFORTHREERADLARGE}
\end{align}

\medskip

\noindent \underline{\textbf{Auxiliary bootstrap assumptions involving two transversal derivatives of}
$\upmu$.}

\begin{align}
	\left\| 
		\Lunit \Rad \Rad \upmu
	\right\|_{L^{\infty}(\Sigma_t^u)}
	& \leq 
		\frac{1}{2}
		\left\| 
			\Rad \Rad
			\left\lbrace
				\vec{G}_{\Lunit \Lunit} \contr \Rad \threePsi
			\right\rbrace
		\right\|_{L^{\infty}(\Sigma_0^u)}
		+ \varepsilon^{1/2},
			\label{E:HIGHERLUNITUPMUBOOT}  \tag{$\mathbf{AUX}\Lunit \Rad \Rad \upmu$}  \\
		\left\| 
			\Rad \Rad \upmu
		\right\|_{L^{\infty}(\Sigma_t^u)}
		& \leq
	 	\left\| 
				\Rad \Rad \upmu
		\right\|_{L^{\infty}(\Sigma_0^u)}
		+ 
		2 \TranminusdatasizeWithFactor^{-1} 
		\left\| 
			\Rad \Rad
			\left\lbrace
				\vec{G}_{\Lunit \Lunit} \contr \Rad \threePsi
			\right\rbrace
		\right\|_{L^{\infty}(\Sigma_0^u)}
		+ \varepsilon^{1/2}.
			\label{E:HIGHERUPMUTRANSVERSALBOOT} 
			\tag{$\mathbf{AUX}\Rad \Rad \upmu$}
\end{align}

\subsection{The main estimates involving higher-order transversal derivatives}
\label{SS:HIGHERORDERTRANSVERALMAINESTIMATES}
In the next proposition, we provide the main estimates of 
Sect.~\ref{S:LINFINITYESTIMATESFORHIGHERTRANSVERSAL}.
The proposition yields, in particular, strict improvements of the bootstrap assumptions 
of Sect.~\ref{SS:BOOTSTRAPFORHIGHERTRANSVERSAL}.

\begin{proposition}[\textbf{$L^{\infty}$ estimates involving higher-order transversal derivatives}]
	\label{P:IMPROVEMENTOFHIGHERTRANSVERSALBOOTSTRAP}
 		Under the data-size and bootstrap assumptions 
		of Sects.~\ref{SS:FLUIDVARIABLEDATAASSUMPTIONS}-\ref{SS:PSIBOOTSTRAP} 
		and Sect.~\ref{SS:BOOTSTRAPFORHIGHERTRANSVERSAL}
		and the smallness assumptions of Sect.~\ref{SS:SMALLNESSASSUMPTIONS}, 
		the following estimates hold
		on $\mathcal{M}_{\Tboot,U_0}$.

	\noindent \underline{\textbf{$L^{\infty}$ estimates involving three transversal derivatives of 
	$\threePsi$}.}
	\begin{align} 
	\left\| 
		\Lunit \Rad \Rad \Rad \threePsi
	\right\|_{L^{\infty}(\Sigma_t^u)}
	& \leq 
		C \varepsilon.
	\label{E:IMPROVEDTRANSVERALESTIMATESFORLUNITTHREERAD}
	\end{align}

	Moreover, for $\Psi \in \lbrace \Densrenormalized ,v^1,v^2 \rbrace$, we have
	\begin{subequations}
	\begin{align}
	\left\| 
		\Rad \Rad \Rad \Psi
	\right\|_{L^{\infty}(\Sigma_t^u)}
	& \leq 
		\left\| 
			\Rad \Rad \Rad \Psi
		\right\|_{L^{\infty}(\Sigma_0^u)}
		+
		C \varepsilon.
		\label{E:IMPROVEDTRANSVERALESTIMATESFORTHREERAD}
	\end{align}
	\end{subequations}

	\medskip

	\noindent \underline{\textbf{$L^{\infty}$ estimates involving two transversal derivatives of $\upmu$}}.
	\begin{subequations}
	\begin{align}
		\left\|
			\Lunit \Rad \Rad \upmu
		\right\|_{L^{\infty}(\Sigma_t^u)}
		& \leq
			\frac{1}{2}
			\left\|
				\Rad
				\Rad
				\left\lbrace
					\vec{G}_{\Lunit \Lunit} \contr \Rad \threePsi
				\right\rbrace
			\right\|_{L^{\infty}(\Sigma_0^u)}
			+
			C \varepsilon,
		\label{E:LUNITRADRADUPMULINFTY}
			\\
		\left\|
			\Rad \Rad \upmu
		\right\|_{L^{\infty}(\Sigma_t^u)}
		& \leq
			\left\|
				\Rad \Rad \upmu
			\right\|_{L^{\infty}(\Sigma_0^u)}
			+
			\TranminusdatasizeWithFactor^{-1}
			\left\|
				\Rad
				\Rad
				\left\lbrace
					\vec{G}_{\Lunit \Lunit} \contr \Rad \threePsi
				\right\rbrace
			\right\|_{L^{\infty}(\Sigma_0^u)}
			+
			C \varepsilon.
		\label{E:RADRADUPMULINFTY}
	\end{align}
	\end{subequations}

	\medskip

	\noindent \underline{\textbf{Sharp pointwise estimates involving the critical factor $\vec{G}_{\Lunit \Lunit}$}}.
	Moreover, if $0 \leq M \leq 2$
	and $0 \leq s \leq t < \Tboot$, 
	then we have the following estimates:
	\begin{align}
		\left|
			\Rad^M \vec{G}_{\Lunit \Lunit}(t,u,\vartheta)
			-
			\Rad^M \vec{G}_{\Lunit \Lunit}(s,u,\vartheta)
		\right|
		& \leq C \varepsilon (t - s),
			\label{E:RADDERIVATIVESOFGLLDIFFERENCEBOUND} \\
		\left|
			\left\lbrace
			\Rad^M 
			\left(
				\vec{G}_{\Lunit \Lunit} \contr \Rad \threePsi
			\right)
			\right\rbrace
			(t,u,\vartheta)
			-
			\left\lbrace
				\Rad^M 
				\left(
					\vec{G}_{\Lunit \Lunit} \contr \Rad \threePsi
				\right)
			\right\rbrace
			(s,u,\vartheta)
		\right|
		& \leq C \varepsilon (t - s).
		\label{E:RADDERIVATIVESOFGLLRADPSIDIFFERENCEBOUND}
	\end{align}

	Furthermore, we have
	\begin{align} \label{E:SMALLNESSOFGLL2}
		\left\|
			G_{\Lunit \Lunit}^2
		\right\|_{L^{\infty}(\Sigma_t^u)}
		& \leq C \varepsilon.
	\end{align}

	Finally, with 
	$
	\displaystyle
	\bar{\Speed}'
	:=
	\frac{d}{d \Densrenormalized} \Speed(\Densrenormalized = 0)
	$,  
	we have
	\begin{align} \label{E:LUNITUPMUDOESNOTDEVIATEMUCHFROMTHEDATA}
		\Lunit \upmu(t,u,\vartheta)
		& = 
		- (\bar{\Speed}' + 1)
		\Rad v^1(t,u,\vartheta)
		+ \mathcal{O}(\varepsilon).
	\end{align}


\end{proposition}

\begin{proof}[Proof of Prop.~\ref{P:IMPROVEMENTOFHIGHERTRANSVERSALBOOTSTRAP}]
	See Sect.~\ref{SS:OFTENUSEDESTIMATES} for some comments on the analysis.
	We must derive the estimates in a viable order.
	Throughout this proof, 
	we use the data-size assumptions of
	Sects.~\ref{SS:FLUIDVARIABLEDATAASSUMPTIONS}
	and \ref{SS:DATAFOREIKONALFUNCTIONQUANTITIES}
	and the assumption \eqref{E:DATAEPSILONVSBOOTSTRAPEPSILON}
	without explicitly mentioning them each time. 
	We refer to these as ``conditions on the data.''

	\medskip
	\noindent \textbf{Proof of \eqref{E:IMPROVEDTRANSVERALESTIMATESFORLUNITTHREERAD}-\eqref{E:IMPROVEDTRANSVERALESTIMATESFORTHREERAD}}:
		By \eqref{E:TRANSRENORMALIZEDAPPLIEDWAVEEQUATIONTRANSPORTINTERPRETATION},
		for $\Psi \in \lbrace \Densrenormalized ,v^1,v^2 \rbrace$,
		we have
	\begin{align} \label{E:WAVEEQUATIONTRANSPORTINTERPRETATION}
		\Lunit \Rad \Psi 
		& 
		= 	\smoothfunction(\BadVar,\ginversesphere,\angdiff \vec{x},\Singletan \Psi,\Rad \Psi) 
				\Singletan \Singletan \Psi
			+ \smoothfunction(\BadVar,\ginversesphere,\angdiff \vec{x},\Singletan \Psi,\Rad \Psi) 
				\Singletan \GdVar
					\\
		& \ \
			+ \smoothfunction(\BadVar,\Fullset^{\leq 1} \threePsi) \Tanset^{\leq 1} \Vortrenormalized.
			\notag
	\end{align}
	Commuting \eqref{E:WAVEEQUATIONTRANSPORTINTERPRETATION}
	with $\Rad \Rad$
	and using Lemmas~\ref{L:POINTWISEFORRECTANGULARCOMPONENTSOFVECTORFIELDS}
	and \ref{L:POINTWISEESTIMATESFORGSPHEREANDITSDERIVATIVES},
	the $L^{\infty}$ estimates of Prop.~\ref{P:IMPROVEMENTOFAUX},
	and the auxiliary bootstrap assumptions of Sect.~\ref{SS:BOOTSTRAPFORHIGHERTRANSVERSAL},
	we find that
	\begin{align} \label{E:PROOFWAVEEQNTRANSPORTTWORADCOMMUTED}
		\left|
			\Lunit \Rad \Rad \Rad \Psi
		\right|
		& \leq
		\left|
			\Lunit \Rad \Rad \Rad \Psi
			-
			\Rad \Rad \Lunit \Rad \Psi
		\right| 
		+
		\left|
			\Rad \Rad \Lunit \Rad \Psi
		\right|
		\\
		& \lesssim
			\left|
				\Lunit \Rad \Rad \Rad \Psi
				-
				\Rad \Rad \Lunit \Rad \Psi
			\right|
			+
			\varepsilon.
		\notag
	\end{align}
	Using in addition the commutator estimate
	\eqref{E:LESSPRECISEPERMUTEDVECTORFIELDSACTINGONFUNCTIONSCOMMUTATORESTIMATE}
	with $f = \Rad \Psi$,
	we obtain the bound
	$
	\displaystyle
	\left|
				\Lunit \Rad \Rad \Rad \Psi
				-
				\Rad \Rad \Lunit \Rad \Psi
			\right|
	\lesssim 
	\varepsilon
	$ 
	as well.
	We have therefore proved \eqref{E:IMPROVEDTRANSVERALESTIMATESFORLUNITTHREERAD}.
	The estimates 
	\eqref{E:IMPROVEDTRANSVERALESTIMATESFORTHREERAD}
	then follow from integrating along the integral curves of $\Lunit$
	as in \eqref{E:SIMPLEFTCID}
	and using the estimate \eqref{E:IMPROVEDTRANSVERALESTIMATESFORLUNITTHREERAD}
	and the conditions on the data.

	\medskip
	\noindent \textbf{Proof of 
	\eqref{E:RADDERIVATIVESOFGLLDIFFERENCEBOUND}-\eqref{E:RADDERIVATIVESOFGLLRADPSIDIFFERENCEBOUND}}:
	It suffices to prove that for $M=0,1,2$ and $\imath = 0,1,2$, 
	we have
	\begin{align} \label{E:LDERIVATIVEOFRADIALDERIVATIVESOFCRITICALFACTOR}
		\left|
			\Lunit \Rad^M G_{\Lunit \Lunit}^{\imath}
		\right|,
			\,
		\left|
			\Lunit \Rad^M 
			\Rad \Psi_{\imath}
		\right|
		& \lesssim \varepsilon,
			\\
		\label{E:RADIALDERIVATIVESOFCRITICALFACTOR}
		\left|
			\Rad^M G_{\Lunit \Lunit}^{\imath}
		\right|,
			\,
		\left|
			\Rad^M 
			\Rad \Psi_{\imath}
		\right|
		& \lesssim 1;
	\end{align}
	once we have shown 
	\eqref{E:LDERIVATIVEOFRADIALDERIVATIVESOFCRITICALFACTOR}-\eqref{E:RADIALDERIVATIVESOFCRITICALFACTOR},
	we can obtain the desired estimates by integrating
	along the integral curves of $\Lunit$ from time $s$ to $t$
	(in analogy with \eqref{E:SIMPLEFTCID})
	and using 
	\eqref{E:LDERIVATIVEOFRADIALDERIVATIVESOFCRITICALFACTOR}-\eqref{E:RADIALDERIVATIVESOFCRITICALFACTOR}.
	The estimate for the second term on LHS~\eqref{E:LDERIVATIVEOFRADIALDERIVATIVESOFCRITICALFACTOR}
	follows from
	\eqref{E:PSIMIXEDUPTOORDERTWOTRANSVERSALIMPROVED}
	and 
	\eqref{E:IMPROVEDTRANSVERALESTIMATESFORLUNITTHREERAD}.
	To prove the desired bound \eqref{E:RADIALDERIVATIVESOFCRITICALFACTOR} 
	for
	$
	\left|
		\Lunit \Rad^M G_{\Lunit \Lunit}^{\imath}
	\right|
	$,
	we first use Lemma~\ref{L:SCHEMATICDEPENDENCEOFMANYTENSORFIELDS} to deduce that
	$\vec{G}_{\Lunit \Lunit} = \smoothfunction(\GdVar)$.
	Differentiating this identity
	with $\Lunit \Rad^M$ 
	and using the $L^{\infty}$ estimates of Prop.~\ref{P:IMPROVEMENTOFAUX},
	we obtain the desired bound.
	The estimate \eqref{E:RADIALDERIVATIVESOFCRITICALFACTOR}
	is a simple consequence of
	the relation $\vec{G}_{\Lunit \Lunit} = \smoothfunction(\GdVar)$,
	the $L^{\infty}$ estimates of Prop.~\ref{P:IMPROVEMENTOFAUX},
	and the estimates
	\eqref{E:IMPROVEDTRANSVERALESTIMATESFORTHREERAD}.

\medskip 
\noindent \textbf{Proof of \eqref{E:SMALLNESSOFGLL2}}:	 
From \eqref{E:GLUNITLUNITi},
the fact that $\Radunit^2 = \Radunit_{(Small)}^2$ (see Def.~\ref{D:PERTURBEDPART}),
and Lemma~\ref{L:SCHEMATICDEPENDENCEOFMANYTENSORFIELDS},
we see that 
$G_{\Lunit \Lunit}^2 = \smoothfunction(\GdVar)\GdVar$.
The desired estimate \eqref{E:SMALLNESSOFGLL2}
now follows as a simple consequence of 
the $L^{\infty}$ estimates of Prop.~\ref{P:IMPROVEMENTOFAUX}.

\medskip
\noindent \textbf{Proof of \eqref{E:LUNITUPMUDOESNOTDEVIATEMUCHFROMTHEDATA}}:
From the evolution equation \eqref{E:UPMUFIRSTTRANSPORT},
\eqref{E:PERTURBEDPART},
the identities \eqref{E:GLUNITLUNIT0}-\eqref{E:GLUNITLUNITi},
Lemma~\ref{L:SCHEMATICDEPENDENCEOFMANYTENSORFIELDS},
and the assumption $\bar{\Speed} = 1$,
we deduce that
\begin{align} \label{E:LUPMUPRECISEEXPANSION}
\Lunit \upmu
	& = \frac{1}{2} 
	\sum_{\imath=0}^1 G_{\Lunit \Lunit}^{\imath} \Rad v^1
	+
	\frac{1}{2}
	G_{\Lunit \Lunit}^0 \Rad (\Densrenormalized  - v^1)
	+
	\frac{1}{2} G_{\Lunit \Lunit}^2 \Rad v^2
	+
	\smoothfunction(\BadVar) \Singletan \threePsi
		\\
	& = 
	- (\bar{\Speed}' + 1)\Rad v^1
	+
	\smoothfunction(\GdVar) \Rad (\Densrenormalized  - v^1)
	+
	\smoothfunction(\BadVar,\Rad \threePsi) \GdVar
	+
	\smoothfunction(\BadVar) \Singletan \threePsi.
	\notag
\end{align}
The desired estimate 
\eqref{E:LUNITUPMUDOESNOTDEVIATEMUCHFROMTHEDATA}
now follows from \eqref{E:LUPMUPRECISEEXPANSION}
and the $L^{\infty}$ estimates of Prop.~\ref{P:IMPROVEMENTOFAUX}
(see especially \eqref{E:CRUCIALPSITRANSVERSALLINFINITYBOUNDBOOTSTRAPIMPROVEDSMALL}).


\medskip
\noindent \textbf{Proof of \eqref{E:LUNITRADRADUPMULINFTY}-\eqref{E:RADRADUPMULINFTY}}:
	With the help of 
	the $L^{\infty}$ estimates of Prop.~\ref{P:IMPROVEMENTOFAUX}
	and the bootstrap assumptions,
	we can use the same argument that we used to prove
	\eqref{E:LUNITUPMUMAINTERM}
	in order to conclude that
	\eqref{E:LUNITUPMUMAINTERM}
	also holds with $M=2$.
	The remainder of the proof of \eqref{E:LUNITRADRADUPMULINFTY}-\eqref{E:RADRADUPMULINFTY}
	now proceeds as in the proof
	of \eqref{E:UPTOONETRANSVERSALDERIVATIVEUPMULINFTY}-\eqref{E:LUNITUPTOONETRANSVERSALUPMULINFINITY}
	(which is given just below \eqref{E:LUNITUPMUMAINTERM}),
	thanks to the availability of the already proven estimates
	\eqref{E:RADDERIVATIVESOFGLLDIFFERENCEBOUND}-\eqref{E:RADDERIVATIVESOFGLLRADPSIDIFFERENCEBOUND}
	in the case $M=2$.

\end{proof}

\section{Sharp Estimates for \texorpdfstring{$\upmu$}{the Inverse Foliation Density}}
\label{S:SHARPESTIMATESFORUPMU}
In this section, we derive sharp pointwise estimates for $\upmu$
and some of its derivatives. These estimates provide much more information than
the crude estimates we obtained in
Sects.~\ref{S:PRELIMINARYPOINTWISE} and \ref{S:LINFINITYESTIMATESFORHIGHERTRANSVERSAL}.
The sharp estimates play an essential role in
our derivation a priori energy estimates (see Sect.~\ref{S:ENERGYESTIMATES}). 
The reason is that in order to obtain the energy estimates,
we must know exactly how $\upmu$ vanishes\footnote{It vanishes linearly;
see \eqref{E:MUSTARBOUNDS}. This fact is of fundamental importance for our
a priori energy estimates.
\label{FN:MUVANISHESLINEARLY}}
and how certain ratios with $\upmu$ in the denominator behave.
This is the main information that we derive in this section.

Many results derived in this section are based on a posteriori estimates.
By this,
we mean estimates for quantities at times $0 \leq s \leq t$ 
that depend on the behavior of other quantities at the ``late time'' $t$, 
where $t < \Tboot$.
For this reason, some of our analysis involves functions
$q = q(s,u,\vartheta;t)$, which we view to be 
functions of the geometric coordinates
$(s,u,\vartheta)$ that depend on the  ``late time parameter'' $t$.
When we state and derive estimates for such quantities, 
$s$ is the ``moving'' time variable verifying $0 \leq s \leq t$.

\subsection{Definitions and preliminary ingredients in the analysis}
\label{SS:PRELIMINARYINGREDIENTSFORSHARPMUESTIMATES}

\begin{definition}[\textbf{Auxiliary quantities used to analyze $\upmu$}]
	\label{D:AUXQUANTITIES}
	We define the following quantities, 
	where $0 \leq s \leq t$
	for those quantities that depend on both $s$ and $t:$
	\begin{subequations}
	\begin{align}
	M(s,u,\vartheta;t) 
	& := \int_{s'=s}^{s'=t} 
					\left\lbrace
						\Lunit \upmu(t,u,\vartheta) - \Lunit \upmu(s',u,\vartheta) 
					\right\rbrace
				\, ds',
		\label{E:BIGMDEF} \\
	\mathring{\upmu}(u,\vartheta)
	& := \upmu(s=0,u,\vartheta),
		\\
	\widetilde{M}(s,u,\vartheta;t)
	& := \frac{M(s,u,\vartheta;t)}{\mathring{\upmu}(u,\vartheta) 
		- M(0,u,\vartheta;t)},
		\label{E:WIDETILDEBIGMDEF} \\
	\upmu_{(Approx)}(s,u,\vartheta;t)
		& := 1
			+  \frac{\Lunit \upmu(t,u,\vartheta)}{
				\mathring{\upmu}(u,\vartheta) 
				- M(0,u,\vartheta;t)}s
			+ \widetilde{M}(s,u,\vartheta;t).
		\label{E:MUAPPROXDEF}
	\end{align}
	\end{subequations}
\end{definition}

The following quantity $\upmu_{\star}$ captures the worst-case smallness
of $\upmu$ along $\Sigma_t^u$. Our high-order energies are allowed to
blow up like a positive power of $\upmu_{\star}^{-1}$.

\begin{definition}[\textbf{Definition of} $\upmu_{\star}$]
	\label{D:MUSTARDEF}
	We define
	\begin{align} \label{E:MUSTARDEF}
		\upmu_{\star}(t,u)
		& := \min \lbrace 1, \min_{\Sigma_t^u} \upmu \rbrace.
	\end{align}
\end{definition}

The next lemma provides basic pointwise estimates for the 
auxiliary quantities.

\begin{lemma}[\textbf{First estimates for the auxiliary quantities}]
\label{L:FIRSTESTIMATESFORAUXILIARYUPMUQUANTITIES}
Under the data-size and bootstrap assumptions 
of Sects.~\ref{SS:FLUIDVARIABLEDATAASSUMPTIONS}-\ref{SS:PSIBOOTSTRAP}
and the smallness assumptions of Sect.~\ref{SS:SMALLNESSASSUMPTIONS}, 
the following  
estimates hold for $(t,u,\vartheta) \in [0,\Tboot) \times [0,U_0] \times \mathbb{T}$
and $0 \leq s \leq t$:
\begin{align}
	\mathring{\upmu}(u,\vartheta)
	& = 1 + \mathcal{O}(\varepsilon),
		\label{E:MUINITIALDATAESTIMATE}
		\\
	\mathring{\upmu}(u,\vartheta)
	& = 1 + M(0,u,\vartheta;t) + \mathcal{O}(\varepsilon).
		\label{E:MUAMPLITUDENEARONE}
\end{align}
In addition, the following pointwise estimates hold:
\begin{align} 
	\left|
		\Lunit \upmu(t,u,\vartheta) 
		- 
		\Lunit \upmu(s,u,\vartheta)
	\right|
	& \lesssim \varepsilon(t-s),
		\label{E:LUNITUPMUATTIMETMINUSLUNITUPMUATTIMESPOINTWISEESTIMATE} \\
	|M(s,u,\vartheta;t)|, |\widetilde{M}(s,u,\vartheta;t)|
	& \lesssim 
		\varepsilon (t - s)^2,
		\label{E:BIGMEST} 
		\\
	\upmu(s,u,\vartheta)
	& = (1 + \mathcal{O}(\varepsilon)) 
	\upmu_{(Approx)}(s,u,\vartheta;t).
	\label{E:MUAPPROXMISLIKEMU}
\end{align}
\end{lemma}

\begin{proof}
	\eqref{E:MUINITIALDATAESTIMATE} follows from \eqref{E:UPMUDATATANGENTIALLINFINITYCONSEQUENCES}
	and \eqref{E:DATAEPSILONVSBOOTSTRAPEPSILON}.
	The estimate \eqref{E:LUNITUPMUATTIMETMINUSLUNITUPMUATTIMESPOINTWISEESTIMATE}
	follows from the mean value theorem and
	the estimate $\left| \Lunit \Lunit \upmu \right| \lesssim \varepsilon$,
	which is a special case of \eqref{E:ZSTARSTARUPMULINFTY}.
	The estimate \eqref{E:MUAMPLITUDENEARONE}
	and the estimate \eqref{E:BIGMEST} for $M$
	then follow from definition \eqref{E:BIGMDEF} 
	and the estimates 
	\eqref{E:MUINITIALDATAESTIMATE}
	and
	\eqref{E:LUNITUPMUATTIMETMINUSLUNITUPMUATTIMESPOINTWISEESTIMATE}.
	The estimate \eqref{E:BIGMEST} for $\widetilde{M}$
	follows from definition \eqref{E:WIDETILDEBIGMDEF}, 
	the estimate \eqref{E:BIGMEST} for $M$,
	and \eqref{E:MUAMPLITUDENEARONE}.
To prove \eqref{E:MUAPPROXMISLIKEMU}, we first note the following identity, which is
a straightforward consequence of Def.~\ref{D:AUXQUANTITIES}:
\begin{align} \label{E:MUSPLIT}
	\upmu(s,u,\vartheta) 
	& = \left\lbrace
				\mathring{\upmu}(u,\vartheta) - M(0,u,\vartheta;t)
			\right\rbrace
			\upmu_{(Approx)}(s,u,\vartheta;t).
\end{align}
The desired estimate \eqref{E:MUAPPROXMISLIKEMU} now follows from
\eqref{E:MUSPLIT}	and \eqref{E:MUAMPLITUDENEARONE}.

\end{proof}

To derive some of the estimates of this section,
it is convenient to partition
various subsets of spacetime 
into regions where $\Lunit \upmu < 0$
(and hence $\upmu$ is decaying) 
and regions where $\Lunit \upmu \geq 0$
(and hence $\upmu$ is not decaying).
This motivates the sets given in the next definition.

\begin{definition}[\textbf{Regions of distinct $\upmu$ behavior}]
\label{D:REGIONSOFDISTINCTUPMUBEHAVIOR}
For each 
$t \in [0,\Tboot)$,
$s \in [0,t]$, 
and $u \in [0,U_0]$, 
we partition 
\begin{subequations}
\begin{align}
	[0,u] \times \mathbb{T} 
	& = \Vplus{t}{u} \cup \Vminus{t}{u},
		\label{E:OUINTERVALCROSSS2SPLIT} \\
	\Sigma_s^u
	& = \Sigmaplus{s}{t}{u} \cup \Sigmaminus{s}{t}{u},
	\label{E:SIGMASSPLIT}
\end{align}
\end{subequations}
where
\begin{subequations}
\begin{align}
	\Vplus{t}{u}
	& := 
	\left\lbrace
		(u',\vartheta) \in [0,u] \times \mathbb{T} \ | \  
			\frac{\Lunit \upmu(t,u',\vartheta)}{\mathring{\upmu}(u',\vartheta) - M(0,u',\vartheta;t)}
		\geq 0
	\right\rbrace,
		\label{E:ANGLESANDUWITHNONDECAYINUPMUGBEHAVIOR} \\
	\Vminus{t}{u}
	& := 
	\left\lbrace
		(u',\vartheta) \in [0,u] \times \mathbb{T} \ | \ 
			\frac{\Lunit \upmu(t,u',\vartheta)}{\mathring{\upmu}(u',\vartheta) - M(0,u',\vartheta;t)} < 0
	\right\rbrace,
		\label{E:ANGLESANDUWITHDECAYINUPMUGBEHAVIOR} \\
	\Sigmaplus{s}{t}{u}
	& := 
	\left\lbrace
		(s,u',\vartheta) \in \Sigma_s^u \ | \ (u',\vartheta) \in \Vplus{t}{u}
	\right\rbrace,
		\label{E:SIGMAPLUS} \\
	\Sigmaminus{s}{t}{u}
	& := 
	\left\lbrace
		(s,u',\vartheta) \in \Sigma_s^u \ | \ (u',\vartheta) \in \Vminus{t}{u}
	\right\rbrace.
	\label{E:SIGMAMINUS}
\end{align}
\end{subequations}

\end{definition}

\begin{remark}[\textbf{The role of the denominators in }\eqref{E:ANGLESANDUWITHNONDECAYINUPMUGBEHAVIOR}-\eqref{E:ANGLESANDUWITHDECAYINUPMUGBEHAVIOR}]
\label{R:POSITIVEDENOMINATORS}
Note that \eqref{E:MUAMPLITUDENEARONE} implies 
that the denominator 
$\mathring{\upmu}(u',\vartheta) - M(0,u',\vartheta;t)$
in \eqref{E:ANGLESANDUWITHNONDECAYINUPMUGBEHAVIOR}-\eqref{E:ANGLESANDUWITHDECAYINUPMUGBEHAVIOR}
remains strictly positive all the way up to the shock.
We include the denominator in the definitions
\eqref{E:ANGLESANDUWITHNONDECAYINUPMUGBEHAVIOR}-\eqref{E:ANGLESANDUWITHDECAYINUPMUGBEHAVIOR}
because it helps 
to clarify the connection
between 
the parameter
$\LateTimeLUnitMu$
defined in \eqref{E:CRUCIALLATETIMEDERIVATIVEDEF}
and the sets
$\Vplus{t}{u}$
and
$\Vminus{t}{u}$.
\end{remark}

\subsection{Sharp pointwise estimates for \texorpdfstring{$\upmu$}{the inverse foliation density} and its derivatives}
\label{SS:SHARPPOINTWISEESTIMATESFORUPMU}
In the next proposition, we provide the sharp pointwise estimates for $\upmu$
that we use to close our energy estimates.

\begin{proposition}[\textbf{Sharp pointwise estimates for $\upmu$, $\Lunit \upmu$, and $\Rad \upmu$}]
\label{P:SHARPMU} 
Under the data-size and bootstrap assumptions 
of Sects.~\ref{SS:FLUIDVARIABLEDATAASSUMPTIONS}-\ref{SS:PSIBOOTSTRAP}
and the smallness assumptions of Sect.~\ref{SS:SMALLNESSASSUMPTIONS}, 
the following  
estimates hold for $(t,u,\vartheta) \in [0,\Tboot) \times [0,U_0] \times \mathbb{T}$
and $0 \leq s \leq t$.
\medskip

\noindent \underline{\textbf{Upper bound for $\displaystyle \frac{[\Lunit \upmu]_+}{\upmu}$}.}
\begin{align} \label{E:POSITIVEPARTOFLMUOVERMUISBOUNDED}
	\left\|
		\frac{[\Lunit \upmu]_+}{\upmu}
	\right\|_{L^{\infty}(\Sigma_s^u)}
	& \leq C.
\end{align}

\medskip

\noindent \underline{\textbf{Small $\upmu$ implies $\Lunit \upmu$ is negative}.}
\begin{align} \label{E:SMALLMUIMPLIESLMUISNEGATIVE}
	\upmu(s,u,\vartheta) \leq \frac{1}{4}
	\implies
	\Lunit \upmu(s,u,\vartheta) \leq - \frac{1}{4} \TranminusdatasizeWithFactor,
\end{align}
where $\TranminusdatasizeWithFactor > 0$ is defined in \eqref{E:CRITICALBLOWUPTIMEFACTOR}.

\medskip

\noindent \underline{\textbf{Upper bound for 
$\displaystyle \frac{[\Rad \upmu]_+}{\upmu}$}.}
\begin{align} \label{E:UNIFORMBOUNDFORMRADMUOVERMU}
	\left\|
		\frac{[\Rad \upmu]_+}{\upmu}
	\right\|_{L^{\infty}(\Sigma_s^u)}
	& \leq 
		\frac{C}{\sqrt{\Tboot - s}}.
\end{align}

\medskip

\noindent \underline{\textbf{Sharp spatially uniform estimates}.}
Consider a time interval $s \in [0,t]$ and define
the ($t,u$-dependent) constant $\LateTimeLUnitMu$ by
\begin{align} \label{E:CRUCIALLATETIMEDERIVATIVEDEF}
	\LateTimeLUnitMu 
	& := 
		\sup_{(u',\vartheta) \in [0,u] \times \mathbb{T}} 
		\frac{[\Lunit \upmu]_-(t,u',\vartheta)}{\mathring{\upmu}(u',\vartheta) - M(0,u',\vartheta;t)},
\end{align}
and note that $\LateTimeLUnitMu \geq 0$ in view of the estimate \eqref{E:MUAMPLITUDENEARONE}.
Then
\begin{subequations}
\begin{align}
	\upmu_{\star}(s,u)
	& = \left\lbrace
				1 + \mathcal{O}(\varepsilon)
			\right\rbrace
			\left\lbrace 
				1 - \LateTimeLUnitMu s 
			\right\rbrace,
	\label{E:MUSTARBOUNDS}  
		\\
	\left\| 
		[\Lunit \upmu]_- 
	\right\|_{L^{\infty}(\Sigma_s^u)}
	& =
	\begin{cases}
		\left\lbrace
			1 + \mathcal{O}(\varepsilon^{1/2})
		\right\rbrace
		\LateTimeLUnitMu,
		& \mbox{if } \LateTimeLUnitMu \geq \sqrt{\varepsilon},
		\\
		\mathcal{O}(\varepsilon^{1/2}),
		& \mbox{if } \LateTimeLUnitMu \leq \sqrt{\varepsilon}.
	\end{cases}
	\label{E:LUNITUPMUMINUSBOUND}
\end{align}
\end{subequations}

Furthermore, we have
\begin{subequations}
\begin{align} \label{E:UNOTNECESSARILYEQUALTOONECRUCIALLATETIMEDERIVATIVECOMPAREDTODATAPARAMETER}
	\LateTimeLUnitMu 
	& \leq \left\lbrace
				1 + \mathcal{O}(\varepsilon)
			\right\rbrace
			\TranminusdatasizeWithFactor.
\end{align}

Moreover, when $u = 1$, we have
\begin{align} \label{E:CRUCIALLATETIMEDERIVATIVECOMPAREDTODATAPARAMETER}
	\LateTimeLUnitMu 
	& = \left\lbrace
				1 + \mathcal{O}(\varepsilon)
			\right\rbrace
			\TranminusdatasizeWithFactor.
\end{align}
\end{subequations}

\medskip

\noindent \underline{\textbf{Sharp estimates when $(u',\vartheta) \in \Vplus{t}{u}$}.}
We recall that the set $\Vplus{t}{u}$ is defined in \eqref{E:ANGLESANDUWITHNONDECAYINUPMUGBEHAVIOR}.
If $0 \leq s_1 \leq s_2 \leq t$, then the following estimate holds:
\begin{align} \label{E:LOCALIZEDMUCANTGROWTOOFAST}
	\sup_{(u',\vartheta) \in \Vplus{t}{u}}
	\frac{\upmu(s_2,u',\vartheta)}{\upmu(s_1,u',\vartheta)}
	& \leq C.
\end{align}

In addition, if $s \in [0,t]$ and $\Sigmaplus{s}{t}{u}$ is as defined in \eqref{E:SIGMAPLUS}, then 
\begin{align}  \label{E:KEYMUNOTDECAYBOUND}
		\inf_{\Sigmaplus{s}{t}{u}} \upmu 
	& \geq 1 - C \varepsilon.
\end{align}

In addition, if $s \in [0,t]$ and $\Sigmaplus{s}{t}{u}$ is as defined in \eqref{E:SIGMAPLUS}, then 
\begin{align} 
	\left\| \frac{[\Lunit \upmu]_-}{\upmu} \right\|_{L^{\infty}(\Sigmaplus{s}{t}{u})}
	& \leq C \varepsilon.
		\label{E:KEYMUNOTDECAYINGMINUSPARTLMUOVERMUBOUND}
\end{align}

\medskip

\noindent \underline{\textbf{Sharp estimates when $(u',\vartheta) \in \Vminus{t}{u}$}.}
Assume that the set $\Vminus{t}{u}$ defined in \eqref{E:ANGLESANDUWITHDECAYINUPMUGBEHAVIOR} 
is non-empty, and consider a time interval $s \in [0,t]$. Let $\LateTimeLUnitMu \geq 0$
be as in \eqref{E:CRUCIALLATETIMEDERIVATIVEDEF}.
Then the following estimate holds:
\begin{align} \label{E:LOCALIZEDMUMUSTSHRINK}
	\mathop{\sup_{0 \leq s_1 \leq s_2 \leq t}}_{(u',\vartheta) \in \Vminus{t}{u}}
	\frac{\upmu(s_2,u',\vartheta)}{\upmu(s_1,u',\vartheta)}
	& \leq 1 + C \varepsilon.
\end{align}
Furthermore, if $s \in [0,t]$ and $\Sigmaminus{s}{t}{u}$
is as defined in \eqref{E:SIGMAMINUS}, then 
\begin{align} \label{E:LMUPLUSNEGLIGIBLEINSIGMAMINUS}
	\left\| [\Lunit \upmu]_+ \right\|_{L^{\infty}(\Sigmaminus{s}{t}{u})}
	& \leq C \varepsilon.
\end{align}
Finally, there exists a constant $C > 0$ 
such that if $0 \leq s \leq t$, then
\begin{align}		\label{E:HYPERSURFACELARGETIMEHARDCASEOMEGAMINUSBOUND}
	\left\| 
		[\Lunit \upmu]_- 
	\right\|_{L^{\infty}(\Sigmaminus{s}{t}{u})}
	& \leq 
		\begin{cases}
		\left\lbrace
			1 + C \varepsilon^{1/2}
		\right\rbrace
		\LateTimeLUnitMu,
		& \mbox{if } \LateTimeLUnitMu \geq \sqrt{\varepsilon},
		\\
		C \varepsilon^{1/2},
		& \mbox{if } \LateTimeLUnitMu \leq \sqrt{\varepsilon}.
	\end{cases}
\end{align}

\noindent \underline{\textbf{Approximate time-monotonicity of $\upmu_{\star}^{-1}(s,u)$}.}
There exists a constant $C > 0$ such that if 
$0 \leq s_1 \leq s_2 \leq t$, then
\begin{align} \label{E:MUSTARINVERSEMUSTGROWUPTOACONSTANT}
	\upmu_{\star}^{-1}(s_1,u) & \leq (1 + C \varepsilon) \upmu_{\star}^{-1}(s_2,u).
\end{align}

\end{proposition}

\begin{proof}
See Sect.~\ref{SS:OFTENUSEDESTIMATES} for some comments on the analysis.

\medskip

\noindent \textbf{Proof of \eqref{E:POSITIVEPARTOFLMUOVERMUISBOUNDED}}:
Clearly it suffices to show that
$
\displaystyle
\frac{[\Lunit \upmu(t,u,\vartheta)]_+}{\upmu(t,u,\vartheta)}
\leq C
$
for $0 \leq t < \Tboot$, $u \in [0,U_0]$, and $\vartheta \in \mathbb{T}$.
We may assume that $\Lunit \upmu(t,u,\vartheta) > 0$ since otherwise
\eqref{E:POSITIVEPARTOFLMUOVERMUISBOUNDED} is trivial. 
Then by \eqref{E:LUNITUPMUATTIMETMINUSLUNITUPMUATTIMESPOINTWISEESTIMATE},
for $0 \leq s' \leq s \leq t < \Tboot \leq 2 \TranminusdatasizeWithFactor^{-1}$,
we have that 
$\Lunit \upmu(s',u,\vartheta) 
\geq 
\Lunit \upmu(s,u,\vartheta)
- C \varepsilon(s-s') 
\geq 
- C \varepsilon
$.
Integrating this estimate with respect to $s'$ starting from $s'=0$
and using \eqref{E:MUINITIALDATAESTIMATE},
we find that 
$\upmu(s,u,\vartheta) \geq 1 - C \varepsilon s \geq 1 - C \varepsilon$
and thus $1/\upmu(s,u,\vartheta) \leq 1 + C \varepsilon$.
Also using the bound 
$\left|
	\Lunit \upmu(s,u,\vartheta)
\right|
\leq C
$
proved in \eqref{E:LUNITUPTOONETRANSVERSALUPMULINFINITY},
we conclude the desired estimate.

\medskip

\noindent \textbf{Proof of \eqref{E:SMALLMUIMPLIESLMUISNEGATIVE}}:
By \eqref{E:LUNITUPMUATTIMETMINUSLUNITUPMUATTIMESPOINTWISEESTIMATE},
for $0 \leq s \leq t < \Tboot \leq 2 \TranminusdatasizeWithFactor^{-1}$,
we have that 
$
\Lunit \upmu(s,u,\vartheta) = \Lunit \upmu(0,u,\vartheta) + \mathcal{O}(\varepsilon)
$.
Integrating this estimate with respect to $s$ starting from $s=0$
and using \eqref{E:MUINITIALDATAESTIMATE},
we find that 
$\upmu(s,u,\vartheta) 
	= 1 
	+ 
	\mathcal{O}(\varepsilon) 
	+
	s \Lunit \upmu(0,u,\vartheta) 
$.
It follows that 
whenever $\upmu(s,u,\vartheta) < 1/4$, we have
$
\Lunit \upmu(0,u,\vartheta) 
< 
- \frac{1}{2} \TranminusdatasizeWithFactor(3/4 + \mathcal{O}(\varepsilon))
= - \frac{3}{8} \TranminusdatasizeWithFactor + \mathcal{O}(\varepsilon)
$.
Again using \eqref{E:LUNITUPMUATTIMETMINUSLUNITUPMUATTIMESPOINTWISEESTIMATE}
to deduce that $\Lunit \upmu(s,u,\vartheta) = \Lunit \upmu(0,u,\vartheta) + \mathcal{O}(\varepsilon)$,
we arrive at the desired estimate \eqref{E:SMALLMUIMPLIESLMUISNEGATIVE}.

\medskip

\noindent \textbf{Proof of \eqref{E:UNOTNECESSARILYEQUALTOONECRUCIALLATETIMEDERIVATIVECOMPAREDTODATAPARAMETER} and \eqref{E:CRUCIALLATETIMEDERIVATIVECOMPAREDTODATAPARAMETER}}:
We prove only 
\eqref{E:UNOTNECESSARILYEQUALTOONECRUCIALLATETIMEDERIVATIVECOMPAREDTODATAPARAMETER}
since
\eqref{E:CRUCIALLATETIMEDERIVATIVECOMPAREDTODATAPARAMETER}
 follows from nearly identical arguments. 
From the first line of \eqref{E:LUPMUPRECISEEXPANSION},
\eqref{E:MUAMPLITUDENEARONE},
\eqref{E:LUNITUPMUATTIMETMINUSLUNITUPMUATTIMESPOINTWISEESTIMATE},
and the $L^{\infty}$ estimates of Prop.~\ref{P:IMPROVEMENTOFAUX},
we have 
\[
\frac{\Lunit \upmu(t,u,\vartheta)}{\mathring{\upmu}(u,\vartheta) - M(0,u,\vartheta;t)}
= 
\Lunit \upmu(0,u,\vartheta)
+ \mathcal{O}(\varepsilon)
= \frac{1}{2} \sum_{\imath=0}^1 [G_{\Lunit \Lunit}^{\imath} \Rad v^1](0,u,\vartheta)
+ \mathcal{O}(\varepsilon).
\]
From this estimate
and definitions
\eqref{E:CRITICALBLOWUPTIMEFACTOR}
and
\eqref{E:CRUCIALLATETIMEDERIVATIVEDEF},
we conclude that
$\LateTimeLUnitMu 
\leq 
\TranminusdatasizeWithFactor + \mathcal{O}(\varepsilon) 
= (1 + \mathcal{O}(\varepsilon)) \TranminusdatasizeWithFactor$,
which yields the desired bound \eqref{E:UNOTNECESSARILYEQUALTOONECRUCIALLATETIMEDERIVATIVECOMPAREDTODATAPARAMETER}.

\medskip

\noindent \textbf{Proof of \eqref{E:MUSTARBOUNDS} and \eqref{E:MUSTARINVERSEMUSTGROWUPTOACONSTANT}}:
We first prove \eqref{E:MUSTARBOUNDS}.
We start by establishing the following preliminary estimate
for the crucial quantity $\LateTimeLUnitMu = \LateTimeLUnitMu(t,u)$ 
(see \eqref{E:CRUCIALLATETIMEDERIVATIVEDEF}):
\begin{align} \label{E:LATETIMELMUTIMESTISLESSTHANONE}
	t \LateTimeLUnitMu
	< 1.
\end{align} 
To proceed, we use
\eqref{E:MUAPPROXDEF},
\eqref{E:MUSPLIT},
\eqref{E:MUAMPLITUDENEARONE},
and
\eqref{E:BIGMEST}
to deduce that the following estimate holds
for $(s,u',\vartheta) \in [0,t] \times [0,u] \times \mathbb{T}$:
\begin{align} \label{E:MUFIRSTLOWERBOUND}
	\upmu(s,u',\vartheta)
	& =
		(1 + \mathcal{O}(\varepsilon))
		\left\lbrace
			1 
			+
			\frac{\Lunit \upmu(t,u',\vartheta)}
			{\mathring{\upmu}(u',\vartheta) - M(0,u',\vartheta;t)} s
			+ 
			\mathcal{O}(\varepsilon) (t-s)^2
		\right\rbrace.
\end{align}
Setting $s=t$ in equation \eqref{E:MUFIRSTLOWERBOUND},
taking the min of both sides  
over $(u',\vartheta) \in [0,u] \times \mathbb{T}$,
and appealing to definitions
\eqref{E:MUSTARDEF} and 
\eqref{E:CRUCIALLATETIMEDERIVATIVEDEF},
we deduce that
$\upmu_{\star}(t,u)
= (1 + \mathcal{O}(\varepsilon))(1-\LateTimeLUnitMu t)
$.
Since $\upmu_{\star}(t,u) > 0$ by \eqref{E:BOOTSTRAPMUPOSITIVITY},
we conclude \eqref{E:LATETIMELMUTIMESTISLESSTHANONE}.

Having established the preliminary estimate, 
we now take the min of both sides of \eqref{E:MUFIRSTLOWERBOUND}
over $(u',\vartheta) \in [0,u] \times \mathbb{T}$
and appeal to definitions
\eqref{E:MUSTARDEF} and \eqref{E:CRUCIALLATETIMEDERIVATIVEDEF}
to obtain:
\begin{align} \label{E:HARDERCASEMUFIRSTLOWERBOUND}
	\min_{(u',\vartheta) \in [0,u] \times \mathbb{T}} \upmu(s,u',\vartheta)
	& = 
		(1 + \mathcal{O}(\varepsilon))
		\left\lbrace
			1 
			- \LateTimeLUnitMu s
			+ \mathcal{O}(\varepsilon) (t-s)^2
		\right\rbrace.
\end{align}
We will show that the terms in braces on RHS~\eqref{E:HARDERCASEMUFIRSTLOWERBOUND} verify
\begin{align} \label{E:MUSECONDLOWERBOUND}
	1 
	- \LateTimeLUnitMu s
	+ \mathcal{O}(\varepsilon) (t-s)^2
	& 
	=
	(1 + \smoothfunction(s,u;t))
	\left\lbrace
		1 - \LateTimeLUnitMu s
	\right\rbrace,
\end{align}
where
\begin{align} \label{E:AMPLITUDEDEVIATIONFUNCTIONMUSECONDLOWERBOUND}
	\smoothfunction(s,u;t)
	& = \mathcal{O}(\varepsilon).
\end{align}
The desired estimate \eqref{E:MUSTARBOUNDS}
then follows easily from 
\eqref{E:HARDERCASEMUFIRSTLOWERBOUND}-\eqref{E:AMPLITUDEDEVIATIONFUNCTIONMUSECONDLOWERBOUND}
and 
definition \eqref{E:MUSTARDEF}.
To prove \eqref{E:AMPLITUDEDEVIATIONFUNCTIONMUSECONDLOWERBOUND}, 
we first use \eqref{E:MUSECONDLOWERBOUND} to solve for $\smoothfunction(s,u;t)$: 
\begin{align} \label{E:AMPLITUDEDEVIATIONFUNCTIONEXPRESSION}
	\smoothfunction(s,u;t)
	=
	\frac{\mathcal{O}(\varepsilon) (t-s)^2
				}
				{
				1 - \LateTimeLUnitMu s
				}
	=
	\frac{\mathcal{O}(\varepsilon) (t-s)^2
				}
				{
				1 - \LateTimeLUnitMu t
				+ 
				\LateTimeLUnitMu (t-s)
				}.
\end{align}
We start by considering the case $\LateTimeLUnitMu \leq (1/4) \TranminusdatasizeWithFactor$.
Since $0 \leq s \leq t < \Tboot \leq 2 \TranminusdatasizeWithFactor^{-1}$,
the denominator in the middle expression in \eqref{E:AMPLITUDEDEVIATIONFUNCTIONEXPRESSION}
is $\geq 1/2$, and the desired estimate
\eqref{E:AMPLITUDEDEVIATIONFUNCTIONMUSECONDLOWERBOUND}
follows easily whenever
$\varepsilon$ is sufficiently small.
In remaining case, we have
$\LateTimeLUnitMu > (1/4) \TranminusdatasizeWithFactor$.
Using \eqref{E:LATETIMELMUTIMESTISLESSTHANONE},
we deduce that 
RHS~\eqref{E:AMPLITUDEDEVIATIONFUNCTIONEXPRESSION}
$
\displaystyle
\leq \frac{1}{\LateTimeLUnitMu} \mathcal{O}(\varepsilon) (t-s)
\leq C \varepsilon \TranminusdatasizeWithFactor^{-2} 
\lesssim \varepsilon
$
as desired.

Inequality \eqref{E:MUSTARINVERSEMUSTGROWUPTOACONSTANT} then follows as a simple consequence of 
\eqref{E:MUSTARBOUNDS}.

\medskip

\noindent \textbf{Proof of \eqref{E:LUNITUPMUMINUSBOUND} and \eqref{E:HYPERSURFACELARGETIMEHARDCASEOMEGAMINUSBOUND}}:
To prove \eqref{E:LUNITUPMUMINUSBOUND},
we first use \eqref{E:LUNITUPMUATTIMETMINUSLUNITUPMUATTIMESPOINTWISEESTIMATE}
to deduce that for $0 \leq s \leq t < \Tboot \leq 2 \TranminusdatasizeWithFactor^{-1}$ 
and $(u',\vartheta) \in [0,u] \times \mathbb{T}$,
we have $\Lunit \upmu(s,u',\vartheta) = \Lunit \upmu(t,u',\vartheta) + \mathcal{O}(\varepsilon)$.
Appealing to definition \eqref{E:CRUCIALLATETIMEDERIVATIVEDEF} and
using the estimate \eqref{E:MUAMPLITUDENEARONE},
we find that
$
\displaystyle
\left\| 
	[\Lunit \upmu]_- 
\right\|_{L^{\infty}(\Sigma_s^u)}
=
	\LateTimeLUnitMu
	+ 
	\mathcal{O}(\varepsilon)
$.
If
$
\displaystyle
\sqrt{\varepsilon}
\leq \LateTimeLUnitMu
$,
we see that as long as $\varepsilon$ is sufficiently small, 
we have the desired bound
$
\displaystyle
	\LateTimeLUnitMu
	+ 
	\mathcal{O}(\varepsilon)
	=
	(1 + \mathcal{O}(\varepsilon^{1/2}))
	\LateTimeLUnitMu
$.
On the other hand, if
$
\displaystyle
\LateTimeLUnitMu
\leq 
\sqrt{\varepsilon}
$,
then similar reasoning yields that
$
\displaystyle
\left\| 
	[\Lunit \upmu]_- 
\right\|_{L^{\infty}(\Sigma_s^u)}
=
	\LateTimeLUnitMu
	+ 
	\mathcal{O}(\varepsilon)
= \mathcal{O}(\sqrt{\varepsilon})
$
as desired. We have thus proved \eqref{E:LUNITUPMUMINUSBOUND}.

The estimate \eqref{E:HYPERSURFACELARGETIMEHARDCASEOMEGAMINUSBOUND} 
can be proved via a similar argument
and we omit the details.

\medskip

\noindent \textbf{Proof of \eqref{E:UNIFORMBOUNDFORMRADMUOVERMU}}:
We fix times $s$ and $t$ with $0 \leq s \leq t < \Tboot$
and a point $p \in \Sigma_s^u$ with geometric coordinates 
$(s,\widetilde{u},\widetilde{\vartheta})$.
Let $\iota : [0,u] \rightarrow \Sigma_s^u$ be the integral curve
of $\Rad$ that passes through $p$ and that is parametrized by the values $u'$ of the eikonal function.
We set 
\[
\displaystyle
F(u') := \upmu \circ \iota(u'),
\qquad
\dot{F}(u') := \frac{d}{d u'} F(u') = (\Rad \upmu)\circ \iota(u').
\]
We must bound 
$
\displaystyle
	\frac{[\Rad \upmu]_+}{\upmu}|_p
	= \frac{[\dot{F}(\widetilde{u})]_+}{F(\widetilde{u})}
$.
We may assume that $\dot{F}(\widetilde{u}) > 0$ since otherwise the desired estimate is trivial.
We now set 
\[
H := \sup_{\mathcal{M}_{\Tboot,U_0}} \Rad \Rad \upmu.
\]
If 
$
\displaystyle
F(\widetilde{u}) > \frac{1}{2}
$,
then the desired estimate is a simple consequence of
\eqref{E:UPTOONETRANSVERSALDERIVATIVEUPMULINFTY} with $M=1$.
We may therefore also assume that
$
\displaystyle
F(\widetilde{u}) \leq \frac{1}{2}
$.
Then in view of the estimate
$
\left\|
	\upmu - 1
\right\|_{L^{\infty}\left(\mathcal{P}_0^{\Tboot}\right)}
\lesssim \varepsilon
$
along $\mathcal{P}_0^{\Tboot}$
(see \eqref{E:UPMULINFINITYALONGP0TBOOT}),
we deduce that there exists a $u'' \in [0,\widetilde{u}]$ such that
$\dot{F}(u'') < 0$. Considering also the assumption $\dot{F}(\widetilde{u}) > 0$,
we see that $H > 0$. Moreover, by \eqref{E:RADRADUPMULINFTY}, we have $H \leq C$.
Furthermore, by continuity, there exists a smallest $u_{\ast} \in [0,\widetilde{u}]$
such that $\dot{F}(u') \geq 0$ for $u' \in [u_{\ast},\widetilde{u}]$.
We also set
\begin{align} \label{E:MUMIN}
\upmu_{(Min)}(s,u') 
:= \min_{(u'',\vartheta) \in [0,u'] \times \mathbb{T}} \upmu(s,u'',\vartheta).
\end{align}
The two main steps in the proof are showing that
\begin{align} \label{E:RADMUOVERMUALGEBRAICBOUND}
	\frac{[\Rad \upmu(s,\widetilde{u},\widetilde{\vartheta})]_+}
		{\upmu(s,\widetilde{u},\widetilde{\vartheta})}
	& \leq 
		H^{1/2}\frac{1}{\sqrt{\upmu_{(Min)}}(s,\widetilde{u})}
\end{align}
and that for
$0 \leq s \leq t < \Tboot$, 
we have
\begin{align} \label{E:UPMUMINLOWERBOUND}
	\upmu_{(Min)}(s,u)
	& \geq 
		\max 
		\left\lbrace
			(1 - C \varepsilon)
			\LateTimeLUnitMu 
			(t-s),
			(1 - C \varepsilon) 
			(1 - \LateTimeLUnitMu s)
		\right\rbrace,
\end{align}
where $\LateTimeLUnitMu = \LateTimeLUnitMu(t,u)$
is defined in \eqref{E:CRUCIALLATETIMEDERIVATIVEDEF}.
Once we have obtained \eqref{E:RADMUOVERMUALGEBRAICBOUND}-\eqref{E:UPMUMINLOWERBOUND}
(see below), 
we split the remainder of the proof 
(which is relatively easy)
into the two cases 
$
\displaystyle
\LateTimeLUnitMu \leq \frac{1}{4} \TranminusdatasizeWithFactor
$
and 
$
\displaystyle
\LateTimeLUnitMu > \frac{1}{4} \TranminusdatasizeWithFactor
$.
In the first case 
$
\displaystyle
\LateTimeLUnitMu \leq \frac{1}{4} \TranminusdatasizeWithFactor
$,
we have
$
\displaystyle
1 - \LateTimeLUnitMu s \geq 1 - \frac{1}{4} \TranminusdatasizeWithFactor \Tboot \geq \frac{1}{2}
$,
and the
desired bound 
$
\displaystyle
\frac{[\Rad \upmu(s,\widetilde{u},\widetilde{\vartheta})]_+}
{\upmu(s,\widetilde{u},\widetilde{\vartheta})}
\leq C 
\leq \frac{C}{\Tboot^{1/2}}
\leq \frac{C}{\sqrt{\Tboot-s}}
\leq \mbox{RHS~\eqref{E:UNIFORMBOUNDFORMRADMUOVERMU}}
$
follows easily from 
\eqref{E:RADMUOVERMUALGEBRAICBOUND}
and the second term in the $\min$ on RHS~\eqref{E:UPMUMINLOWERBOUND}.
In the remaining case 
$
\displaystyle
\LateTimeLUnitMu > \frac{1}{4} \TranminusdatasizeWithFactor
$, 
we have 
$
\displaystyle
\frac{1}{\LateTimeLUnitMu} \leq C
$, 
and using the first term in the $\min$ on RHS~\eqref{E:UPMUMINLOWERBOUND}, 
we
deduce that 
$
\displaystyle
\frac{[\Rad \upmu(s,\widetilde{u},\widetilde{\vartheta})]_+}
{\upmu(s,\widetilde{u},\widetilde{\vartheta})}
\leq \frac{C}{\sqrt{t-s}}
$.
Since this estimate holds for all $t < \Tboot$ with a uniform constant $C$,
we conclude \eqref{E:UNIFORMBOUNDFORMRADMUOVERMU} in this case.

We now prove \eqref{E:RADMUOVERMUALGEBRAICBOUND}.
To this end, we will show that 
\begin{align} \label{E:FIRSTRADMUOVERMUALGEBRAICBOUND}
	\frac{[\Rad \upmu(s,\widetilde{u},\widetilde{\vartheta})]_+}
		{\upmu(s,\widetilde{u},\widetilde{\vartheta})}
	& \leq 
		2 H^{1/2} \frac{\sqrt{\upmu(s,\widetilde{u},\widetilde{\vartheta}) - \upmu_{(Min)}(s,\widetilde{u})}}
		{\upmu(s,\widetilde{u},\widetilde{\vartheta})}.
\end{align}
Then viewing RHS~\eqref{E:RADMUOVERMUALGEBRAICBOUND} as a function
of the real variable $\upmu(s,\widetilde{u},\widetilde{\vartheta})$ 
(with all other parameters fixed)
on the domain $[\upmu_{(Min)}(s,\widetilde{u}),\infty)$, we
carry out a simple calculus exercise 
to find that RHS~\eqref{E:RADMUOVERMUALGEBRAICBOUND}
$
\displaystyle
\leq H^{1/2}\frac{1}{\sqrt{\upmu_{(Min)}(s,\widetilde{u})}}
$,
which yields \eqref{E:RADMUOVERMUALGEBRAICBOUND}.
We now prove \eqref{E:FIRSTRADMUOVERMUALGEBRAICBOUND}.
For any $u' \in [u_{\ast},\widetilde{u}]$, 
we use the mean value theorem to obtain 
\begin{align} \label{E:MVTESTIMATES}
\dot{F}(\widetilde{u}) 
- 
\dot{F}(u') 
\leq H (\widetilde{u} - u'),
	\qquad
F(\widetilde{u}) 
- 
F(u') 
\geq 
\min_{u'' \in [u',\widetilde{u}]} \dot{F}(u'')
(\widetilde{u} - u').
\end{align}
Setting 
$
\displaystyle
u_1 := 
\widetilde{u} 
- 
\frac{1}{2} \frac{\dot{F}(\widetilde{u})}{H}
$, 
we find from the first estimate in \eqref{E:MVTESTIMATES} that 
for $u' \in [u_1,\widetilde{u}]$, we have
$
\displaystyle
\dot{F}(u')
\geq 
\frac{1}{2} \dot{F}(\widetilde{u})
$.
Using also the second estimate in \eqref{E:MVTESTIMATES}, we
find that for $u' \in [u_1,\widetilde{u}]$, we have
$
\displaystyle
F(\widetilde{u}) 
- 
F(u_1)
\geq 
\frac{1}{2} \dot{F}(\widetilde{u}) (\widetilde{u}-u_1)
= \frac{1}{4} \frac{\dot{F}^2(\widetilde{u})}{H}
$. 
Noting that 
the definition of $\upmu_{(Min)}$ implies that
$F(u_1) 
\geq \upmu_{(Min)}(s,\widetilde{u})
$,
we deduce that
\begin{align} \label{E:FDIFFERENCELOWERBOUND}
	\upmu (s,\widetilde{u},\widetilde{\vartheta})
	- 
	\upmu_{(Min)}(s,\widetilde{u})
	& \geq 
	\frac{1}{4} \frac{[\Rad \upmu(s,\widetilde{u},\widetilde{\vartheta})]_+^2}{H}.
\end{align}
Taking the square root of \eqref{E:FDIFFERENCELOWERBOUND},
rearranging,
and dividing by
$\upmu(s,\widetilde{u},\widetilde{\vartheta})$,
we conclude the desired estimate \eqref{E:FIRSTRADMUOVERMUALGEBRAICBOUND}.

It remains for us to prove \eqref{E:UPMUMINLOWERBOUND}. 
Reasoning as in the proof of 
\eqref{E:MUFIRSTLOWERBOUND}-\eqref{E:AMPLITUDEDEVIATIONFUNCTIONMUSECONDLOWERBOUND}
and using \eqref{E:LATETIMELMUTIMESTISLESSTHANONE},
we find that for $0 \leq s \leq t < \Tboot$ and $u' \in [0,u]$,
we have
$\upmu_{(Min)}(s,u') 
\geq 
(1 - C \varepsilon)
\left\lbrace
	1 - \LateTimeLUnitMu s
\right\rbrace
\geq
(1 - C \varepsilon)
\LateTimeLUnitMu (t-s)
$.
From these two inequalities, 
we conclude \eqref{E:UPMUMINLOWERBOUND}.

\medskip

\noindent {\textbf{Proof of} \eqref{E:LOCALIZEDMUMUSTSHRINK}}:
A straightforward modification of the proof of \eqref{E:MUSTARBOUNDS},
based on equation 
\eqref{E:MUFIRSTLOWERBOUND}
and on replacing $\LateTimeLUnitMu$ in 
\eqref{E:HARDERCASEMUFIRSTLOWERBOUND}-\eqref{E:MUSECONDLOWERBOUND} 
with $\Lunit \upmu(t,u',\vartheta)$
(without taking the $\min$ on the LHS of the analog of \eqref{E:HARDERCASEMUFIRSTLOWERBOUND}),
yields that for $0 \leq s \leq t < \Tboot$
and $(u',\vartheta) \in \Vminus{t}{u}$,
we have
$
\displaystyle
\upmu(s,u',\vartheta) 
= \left\lbrace
				1 + \mathcal{O}(\varepsilon)
			\right\rbrace
			\left\lbrace 
				1 - \left|\Lunit \upmu(t,u',\vartheta)\right| s 
			\right\rbrace$.
The estimate \eqref{E:LOCALIZEDMUMUSTSHRINK} 
then follows as a simple consequence.

\medskip

\noindent {\textbf{Proof of}
\eqref{E:LOCALIZEDMUCANTGROWTOOFAST}, \eqref{E:KEYMUNOTDECAYBOUND}, \textbf{and} \eqref{E:KEYMUNOTDECAYINGMINUSPARTLMUOVERMUBOUND}}:
By \eqref{E:LUNITUPMUATTIMETMINUSLUNITUPMUATTIMESPOINTWISEESTIMATE},
if $(u',\vartheta) \in \Vplus{t}{u}$ and
$0 \leq s \leq t < \Tboot$,
then
$[\Lunit \upmu]_-(s,u,\vartheta) \leq C \varepsilon$
and
$\Lunit \upmu(s,u,\vartheta) \geq - C \varepsilon$.
Integrating the latter estimate with respect to $s$ from $0$ to $t$
and using \eqref{E:MUINITIALDATAESTIMATE},
we find that  
$\upmu(s,u',\vartheta) 
\geq 1 - C \varepsilon s
\geq 1 - C \varepsilon
$.
Moreover, from \eqref{E:UPTOONETRANSVERSALDERIVATIVEUPMULINFTY} with $M=0$, 
we have the crude bound $\upmu(s,u',\vartheta) \leq C$.
The desired bounds
\eqref{E:LOCALIZEDMUCANTGROWTOOFAST}, \eqref{E:KEYMUNOTDECAYBOUND}, and 
\eqref{E:KEYMUNOTDECAYINGMINUSPARTLMUOVERMUBOUND}
now readily follow
from these estimates.

\medskip

\noindent {\textbf{Proof of} \eqref{E:LMUPLUSNEGLIGIBLEINSIGMAMINUS}}:
By \eqref{E:LUNITUPMUATTIMETMINUSLUNITUPMUATTIMESPOINTWISEESTIMATE},
if $(u',\vartheta) \in \Vminus{t}{u}$ and
$0 \leq s \leq t < \Tboot$,
then
$[\Lunit \upmu]_+(s,u',\vartheta)
=  [\Lunit \upmu]_+(t,u',\vartheta) + \mathcal{O}(\varepsilon) = \mathcal{O}(\varepsilon)$.
The desired bound \eqref{E:LMUPLUSNEGLIGIBLEINSIGMAMINUS} thus follows.

\end{proof}

\subsection{Sharp time-integral estimates involving \texorpdfstring{$\upmu$}{the inverse foliation density}}
\label{SS:SHARPTIMEINTEGRALESTIMATES}
Some error integrals appearing in our top-order energy identities contain a dangerous factor of 
$1/\upmu$. This forces us, in our Gronwall argument for a priori energy estimates, 
to derive estimates for time integrals involving various powers of $\upmu_{\star}^{-1}$.
In Prop.~\ref{P:MUINVERSEINTEGRALESTIMATES}, we derive estimates for these time integrals.
The estimates of Prop.~\ref{P:MUINVERSEINTEGRALESTIMATES}
directly influence the blowup-exponents $P$ 
featured in our high-order energy estimates, 
which are allowed to blow up like $\upmu_{\star}^{-P}$.
In particular, when controlling the size of the blowup-exponents,
we will use the fact that the estimates 
\eqref{E:KEYMUTOAPOWERINTEGRALBOUND}
and
\eqref{E:KEYHYPERSURFACEMUTOAPOWERINTEGRALBOUND}
have coefficient factors $1 + C \sqrt{\varepsilon}$
on the right-hand sides; larger coefficient factors would lead
to larger blowup-exponents.

\begin{proposition}[\textbf{Fundamental estimates for time integrals involving $\upmu^{-1}$}] 
\label{P:MUINVERSEINTEGRALESTIMATES}
	Let $\upmu_{\star}(t,u)$ be as defined in \eqref{E:MUSTARDEF}.
	Let
	\begin{align*}
		\Contwo > 1
	\end{align*}
	be a real number. 
	Under the data-size and bootstrap assumptions 
	of Sects.~\ref{SS:FLUIDVARIABLEDATAASSUMPTIONS}-\ref{SS:PSIBOOTSTRAP}
	and the smallness assumptions of Sect.~\ref{SS:SMALLNESSASSUMPTIONS}, 
	the following  
	estimates hold for $(t,u) \in [0,\Tboot] \times [0,U_0]$.

	\medskip

	\noindent \underline{\textbf{Estimates relevant for borderline top-order spacetime integrals}.}
	There exists a constant $C > 0$ 
	such that if $\Contwo \sqrt{\varepsilon} \leq 1$,
	then
	\begin{align} \label{E:KEYMUTOAPOWERINTEGRALBOUND}
		\int_{s=0}^t 
			\frac{\left\| [\Lunit \upmu]_- \right\|_{L^{\infty}(\Sigma_s^u)}} 
					 {\upmu_{\star}^{\Contwo}(s,u)}
		\, ds 
		& \leq \frac{1 + C \sqrt{\varepsilon}}{\Contwo-1} 
			\upmu_{\star}^{1-\Contwo}(t,u).
	\end{align}

	\noindent \underline{\textbf{Estimates relevant for borderline top-order hypersurface integrals}.}
	There exists a constant $C > 0$ such that
	\begin{align} \label{E:KEYHYPERSURFACEMUTOAPOWERINTEGRALBOUND}
		\left\| 
			\Lunit \upmu 
		\right\|_{L^{\infty}(\Sigmaminus{t}{t}{u})} 
		\int_{s=0}^t 
			\frac{1} 
				{\upmu_{\star}^{\Contwo}(s,u)}
			\, ds 
		& \leq \frac{1 + C \sqrt{\varepsilon}}{\Contwo-1} 
			\upmu_{\star}^{1-\Contwo}(t,u).
	\end{align}

	\medskip

	\noindent \underline{\textbf{Estimates relevant for less dangerous top-order spacetime integrals}.}
	There exists a constant $C > 0$ 
	such that if $\Contwo \sqrt{\varepsilon} \leq 1$,
	then
	\begin{align} \label{E:LOSSKEYMUINTEGRALBOUND}
		\int_{s=0}^t \frac{1} 
			{\upmu_{\star}^{\Contwo}(s,u)}
		\, ds 
		& \leq C \left\lbrace 
			2^{\Contwo}
			+ 
			\frac{1}{\Contwo-1} \right\rbrace \upmu_{\star}^{1-\Contwo}(t,u).
	\end{align} 

	\medskip

	\noindent \underline{\textbf{Estimates for integrals that lead to only $\ln \upmu_{\star}^{-1}$ degeneracy}.} 
	There exists a constant $C > 0$ such that
	\begin{align} \label{E:KEYMUINVERSEINTEGRALBOUND}
		\int_{s=0}^t 
			\frac{\left\| [\Lunit \upmu]_- \right\|_{L^{\infty}(\Sigma_s^u)}} 
					 {\upmu_{\star}(s,u)}
		\, ds 
		& \leq (1 + C \sqrt{\varepsilon}) \ln \upmu_{\star}^{-1}(t,u) + C \sqrt{\varepsilon}.
	\end{align}
	In addition, there exists a constant $C > 0$ such that
	\begin{align} \label{E:LOGLOSSMUINVERSEINTEGRALBOUND}
		\int_{s=0}^t 
			\frac{1}{\upmu_{\star}(s,u)}
		\, ds 
		& \leq  C \left\lbrace \ln \upmu_{\star}^{-1}(t,u) + 1 \right\rbrace.
	\end{align}

	\medskip

	\noindent \underline{\textbf{Estimates for integrals that break the $\upmu_{\star}^{-1}$ degeneracy}.} 
	There exists a constant $C > 0$ such that
	\begin{align} \label{E:LESSSINGULARTERMSMPOINTNINEINTEGRALBOUND}
		\int_{s=0}^t 
			\frac{1} 
			{\upmu_{\star}^{9/10}(s,u)}
		\, ds 
		& \leq C.
	\end{align}

\end{proposition}

\begin{proof}
\noindent {\textbf{Proof of} \eqref{E:KEYMUTOAPOWERINTEGRALBOUND}, 
\eqref{E:KEYHYPERSURFACEMUTOAPOWERINTEGRALBOUND},
and \eqref{E:KEYMUINVERSEINTEGRALBOUND}}:
To prove \eqref{E:KEYMUTOAPOWERINTEGRALBOUND}, we 
first consider the case 
$\LateTimeLUnitMu \geq \sqrt{\varepsilon}$
in \eqref{E:LUNITUPMUMINUSBOUND}.
Using \eqref{E:MUSTARBOUNDS} and \eqref{E:LUNITUPMUMINUSBOUND}, we deduce that
\begin{align} \label{E:PROOFKEYMUTOAPOWERINTEGRALBOUND}
		\int_{s=0}^t 
			\frac{\left\| [\Lunit \upmu]_- \right\|_{L^{\infty}(\Sigma_s^u)}} 
					 {\upmu_{\star}^{\Contwo}(s,u)}
		\, ds 
		& = (1 + \mathcal{O}(\varepsilon^{1/2}))
			\int_{s=0}^t 
				\frac{\LateTimeLUnitMu}{(1 - \LateTimeLUnitMu s)^{\Contwo}}
			\, ds 
				\\
		& \leq \frac{1 + \mathcal{O}(\varepsilon^{1/2})}{\Contwo - 1}
				\frac{1}{(1 - \LateTimeLUnitMu t)^{\Contwo-1}}
			= 	\frac{1 + \mathcal{O}(\varepsilon^{1/2})}{\Contwo - 1}
					\upmu_{\star}^{1-\Contwo}(t,u)
					\notag
	\end{align}
	as desired. We now consider the remaining case
	$\LateTimeLUnitMu \leq \sqrt{\varepsilon}$
	in \eqref{E:LUNITUPMUMINUSBOUND}.
	Using 
	\eqref{E:MUSTARBOUNDS}
	and
	\eqref{E:LUNITUPMUMINUSBOUND}
	and
	the fact that $0 \leq s \leq t < \Tboot \leq 2 \TranminusdatasizeWithFactor^{-1}$,
	we see that for $\varepsilon$ sufficiently small relative to 
	$\TranminusdatasizeWithFactor$, 
	we have
	\begin{align} \label{E:SECONDCASEPROOFKEYMUTOAPOWERINTEGRALBOUND}
		\int_{s=0}^t 
			\frac{\left\| [\Lunit \upmu]_- \right\|_{L^{\infty}(\Sigma_s^u)}} 
					 {\upmu_{\star}^{\Contwo}(s,u)}
		\, ds 
		& \leq
			C \varepsilon^{1/2}
			\int_{s=0}^t 
				\frac{1}{(1 - \LateTimeLUnitMu s)^{\Contwo}}
			\, ds 
				\\
		& \leq
			C \varepsilon^{1/2}
			\frac{1}{(1 - \LateTimeLUnitMu t)^{\Contwo-1}}
			\leq \frac{1}{\Contwo - 1}
			\upmu_{\star}^{1-\Contwo}(t,u)
					\notag
	\end{align}
	as desired.
	We have thus proved \eqref{E:KEYMUTOAPOWERINTEGRALBOUND}.

	Inequality \eqref{E:KEYMUINVERSEINTEGRALBOUND}
	can be proved using similar arguments and we omit the details.

	Inequality \eqref{E:KEYHYPERSURFACEMUTOAPOWERINTEGRALBOUND} 
	can be proved using similar arguments
	with the help of the estimate \eqref{E:HYPERSURFACELARGETIMEHARDCASEOMEGAMINUSBOUND}
	and we omit the details.

	\medskip

\noindent {\textbf{Proof of} \eqref{E:LOSSKEYMUINTEGRALBOUND},
\eqref{E:LOGLOSSMUINVERSEINTEGRALBOUND}, 
and \eqref{E:LESSSINGULARTERMSMPOINTNINEINTEGRALBOUND}}:
	To prove \eqref{E:LOSSKEYMUINTEGRALBOUND}, 
	we first use \eqref{E:MUSTARBOUNDS} to deduce
	\begin{align} \label{E:PROOFLOSSKEYMUINTEGRALBOUND}
		\int_{s=0}^t \frac{1} 
			{\upmu_{\star}^{\Contwo}(s,u)}
		\, ds 
		& \leq 
		C
		\int_{s=0}^t 
			\frac{1}{(1 - \LateTimeLUnitMu s)^{\Contwo}}
		\, ds,
	\end{align}
	where $\LateTimeLUnitMu = \LateTimeLUnitMu(t,u)$
	is defined in \eqref{E:CRUCIALLATETIMEDERIVATIVEDEF}.
	We first assume that 
	$
	\displaystyle
	\LateTimeLUnitMu \leq \frac{1}{4} \TranminusdatasizeWithFactor
	$. 
	Then since $0 \leq t < \Tboot < 2 \TranminusdatasizeWithFactor^{-1}$,
	we see from \eqref{E:MUSTARBOUNDS} that 
	$
	\displaystyle
	\upmu_{\star}(s,u) \geq \frac{1}{4}$ for $0 \leq s \leq t
	$
	and that RHS~\eqref{E:PROOFLOSSKEYMUINTEGRALBOUND} 
	$\leq C 2^{\Contwo} t 
	\leq C 2^{\Contwo} \TranminusdatasizeWithFactor^{-1} 
	\leq 
	C 2^{\Contwo} \leq C 2^{\Contwo} \upmu_{\star}^{1-\Contwo}(t,u)$
	as desired.
	In the remaining case, we have 
	$
	\displaystyle
	\LateTimeLUnitMu > \frac{1}{4} \TranminusdatasizeWithFactor
	$,
	and we can use \eqref{E:MUSTARBOUNDS} 
	and the estimate 
	$
	\displaystyle
	\frac{1}{\LateTimeLUnitMu} 
	\leq C$
	to bound RHS~\eqref{E:PROOFLOSSKEYMUINTEGRALBOUND} by
	\begin{align} \label{E:SECONDPROORELATIONFLOSSKEYMUINTEGRALBOUND}
		&
		\leq
		\frac{C}{\LateTimeLUnitMu}
		\frac{1}{\Contwo - 1}
		\frac{1}{(1 - \LateTimeLUnitMu t)^{\Contwo-1}}
		\leq 
		\frac{C}{\Contwo - 1}
		\upmu_{\star}^{1-\Contwo}(t,u)
	\end{align}
	as desired.

	Inequalities 
	\eqref{E:LOGLOSSMUINVERSEINTEGRALBOUND}
	and
	\eqref{E:LESSSINGULARTERMSMPOINTNINEINTEGRALBOUND} can be proved in a similar
	fashion. We omit the details, aside from remarking that the last step of the proof of 
	\eqref{E:LESSSINGULARTERMSMPOINTNINEINTEGRALBOUND}
	relies on the trivial estimate $(1 - \LateTimeLUnitMu t)^{1/10} \leq 1$.

\end{proof}

\section{The fundamental \texorpdfstring{$L^2$-}{square integral-}controlling quantities}
\label{S:FUNDAMENTALL2CONTROLLINGQUANTITIES}
\setcounter{equation}{0}
In this section, we define the 
``fundamental $L^2$-controlling quantities''
that we use to control $\threePsi$,
$\Vortrenormalized$, and their derivatives in $L^2$.
We also exhibit their coerciveness properties.
These are the quantities that, in Sect.~\ref{S:ENERGYESTIMATES}, 
we will estimate using a Gronwall-type inequality.

\subsection{Definitions of the fundamental \texorpdfstring{$L^2$-}{square integral-}controlling quantities}
\label{SS:DEFSOFL2CONTROLLINGQUANTITIES}

\begin{definition}[\textbf{The main coercive quantities used for controlling 
the solution and its derivatives in} $L^2$]
\label{D:MAINCOERCIVEQUANT}
In terms of the energy-null flux quantities of Def.~\ref{D:ENERGYFLUX}
and the multi-index set $\mathcal{I}^{N;\leq 1}$
of Def.~\ref{E:COMMUTATORMULTIINDICES}, 
we define
\begin{subequations}
\begin{align}
	\totmax{N}(t,u)
	& := 
		\mathop{\max_{\vec{I} \in \mathcal{I}^{N;\leq 1}}}_{
			\Psi \in \lbrace \Densrenormalized  - v^1,v^1,v^2 \rbrace}
		\sup_{(t',u') \in [0,t] \times [0,u]} 
		\left\lbrace
			\enzero[\Fullset^{\vec{I}} \Psi](t',u')
			+ \flzero[\Fullset^{\vec{I}} \Psi](t',u')
		\right\rbrace,
			\label{E:Q0TANNDEF} \\
	\easytotmax{N}(t,u)
	& := 
		\mathop{\max_{\vec{I} \in \mathcal{I}^{N;\leq 1}}}_{\Psi \in 
		\lbrace \Densrenormalized  - v^1,v^2
		\rbrace}
		\sup_{(t',u') \in [0,t] \times [0,u]} 
		\left\lbrace
			\enzero[\Fullset^{\vec{I}} \Psi](t',u')
			+ \flzero[\Fullset^{\vec{I}} \Psi](t',u')
		\right\rbrace,
			\label{E:EASYQ0TANNDEF} \\
	\VorttotTanmax{N}(t,u)
	& :=
		\max_{|\vec{I}|=N}
		\sup_{(t',u') \in [0,t] \times [0,u]} 
		\left\lbrace
			\Vortenzero[\Tanset^{\vec{I}} \Vortrenormalized](t',u')
			+ \Vortflzero[\Tanset^{\vec{I}} \Vortrenormalized](t',u')
		\right\rbrace,
		\label{E:VORTQ0TANNDEF}
			\\
	\totmax{[1,N]}(t,u)
	& := \max_{1 \leq M \leq N} \totmax{M}(t,u),
		\label{E:MAXEDQ0TANLEQNDEF} 
			\\
	\VorttotTanmax{\leq N}(t,u)
	& := \max_{0 \leq M \leq N} \VorttotTanmax{M}(t,u).
		\label{E:VORTMAXEDQ0TANLEQNDEF} 
\end{align}
\end{subequations}
\end{definition}

\begin{remark}[\textbf{Carefully note what is controlled by $\totmax{N}$ and $\easytotmax{N}$}]
	\label{R:WHICHVARIABLESARECONTROLLED}
	Note that $\totmax{N}$ directly controls the derivatives of
	$\Densrenormalized  - v^1$,
	$v^1$,
	and
	$v^2$,
	while control of derivatives of
	$\Densrenormalized $
	in terms of $\totmax{N}$
	(and thus for all three entries of the array $\threePsi$)
	can be achieved via the triangle inequality.
	Similarly, 
	$\easytotmax{N}$
	 directly controls the derivatives of
	$\Densrenormalized  - v^1$
	and
	$v^2$.
	The quantity $\easytotmax{N}$ might seem to be unnecessary, but 
	it plays an important role in our energy estimates.
	The reason is that the most degenerate error terms in the energy estimates
	for $\easytotmax{N}$ are multiplied by a small factor $\varepsilon$;
	see \eqref{E:EASYVARIABLESTOPORDERWAVEENERGYINTEGRALINEQUALITIES}. 
	This is important because $\easytotmax{N}$ appears as a 
	large coefficient (denoted by $C_*$) source term
	in the energy estimates for $\totmax{N}$;
	see \eqref{E:TOPORDERWAVEENERGYINTEGRALINEQUALITIES}.
	The smallness of $\varepsilon$ compensates for the largeness of $C_*$
	and allows us to close our energy estimates; see also the discussion in Subsubsect.~\ref{sec.gd.comp}.
	Similar remarks apply to the integrals
	of Def.~\ref{D:COERCIVEINTEGRAL}.
\end{remark}

The following spacetime integrals are indispensable for controlling the geometric
torus derivatives of $\threePsi$. Their key property is that they
are strong in regions where $\upmu$ is small;
see Lemma~\ref{L:KEYSPACETIMECOERCIVITY}.
Recall that these spacetime integrals
are featured in the basic energy identity for solutions to the 
wave equation; see~\eqref{E:MULTERRORINT}.

\begin{definition}[\textbf{Key coercive spacetime integrals}]
\label{D:COERCIVEINTEGRAL}
	Let $\mathcal{I}^{N;\leq 1}$ be the multi-index set from Def.~\ref{D:COMMUTATORMULTIINDICES}.
	We associate the following integrals to $\Psi$,
	where $[\Lunit \upmu]_- = |\Lunit \upmu|$ when $\Lunit \upmu < 0$
	and $[\Lunit \upmu]_- = 0$ when $\Lunit \upmu \geq 0$:
	\begin{subequations}
	\begin{align} \label{E:COERCIVESPACETIMEDEF} 
		\coercivespacetime[\Psi](t,u)
		& :=
	 	\frac{1}{2}
	 	\int_{\mathcal{M}_{t,u}}
			[\Lunit \upmu]_-
			|\angdiff \Psi|^2
		\, d \vol, 
				\\
		\coercivespacetimemax{N}(t,u) 
		& := 
			\mathop{\max_{\vec{I} \in \mathcal{I}^{N;\leq 1}}}_{\Psi \in 
			\lbrace 
				\Densrenormalized  - v^1,
				v^1,
				v^2
			\rbrace}
			\coercivespacetime[\Fullset^{\vec{I}} \Psi](t,u)
			,
			\\
		\coercivespacetimemax{[1,N]}(t,u) 
		& := \max_{1 \leq M \leq N} \coercivespacetimemax{M}(t,u),
		\label{E:MAXEDCOERCIVESPACETIMEDEF}
			\\
		\easycoercivespacetimemax{N}(t,u) 
		& := 
			\mathop{\max_{\vec{I} \in \mathcal{I}^{N;\leq 1}}}_{\Psi \in 
			\lbrace 
				\Densrenormalized  - v^1,
				v^2
			\rbrace}
			\coercivespacetime[\Fullset^{\vec{I}} \Psi](t,u),
			\\
		\easycoercivespacetimemax{[1,N]}(t,u) 
		& := \max_{1 \leq M \leq N} \easycoercivespacetimemax{M}(t,u).
		\label{E:EASYMAXEDCOERCIVESPACETIMEDEF}
	\end{align}
	\end{subequations}
\end{definition}

\subsection{Comparison of area forms and estimates for the \texorpdfstring{$L^2$-}{square integral} 
norm of time integrals}
\label{SS:PRELIMINARYLEMMAS}
We now provide some preliminary lemmas that we will use
in our $L^2$ analysis.

\begin{lemma}[\textbf{Pointwise estimates for} $\gtancomp$]
\label{L:POINTWISEESTIMATEFORGTANCOMP}
Let $\gtancomp$ be the metric component from
Def.~\ref{D:METRICANGULARCOMPONENT}.
Under the data-size and bootstrap assumptions 
of Sects.~\ref{SS:FLUIDVARIABLEDATAASSUMPTIONS}-\ref{SS:PSIBOOTSTRAP}
and the smallness assumptions of Sect.~\ref{SS:SMALLNESSASSUMPTIONS}, 
the following estimate holds
on $\mathcal{M}_{\Tboot,U_0}$:
\begin{align} \label{E:POINTWISEESTIMATEFORGTANCOMP}
	\gtancomp = 1 + \mathcal{O}(\varepsilon).
\end{align}
\end{lemma}

\begin{proof}
	See Sect.~\ref{SS:OFTENUSEDESTIMATES} for some comments on the analysis.
	Using the identity \eqref{E:LDERIVATIVEOFVOLUMEFORMFACTOR} and the estimate
	\eqref{E:TWORADIALCHICOMMUTEDLINFINITY}, we deduce
	that $\Lunit \ln \gtancomp = \mathcal{O}(\varepsilon)$.
	Integrating this estimate
	along the integral curves of $\Lunit$
	as in \eqref{E:SIMPLEFTCID},
	we find that $\ln \gtancomp(t,u,\vartheta) = \ln \gtancomp(0,u,\vartheta) + \mathcal{O}(\varepsilon)$.
	From this estimate and the small-data bound
	$\gtancomp(0,u,\vartheta) = 1 + \mathcal{O}(\varepsilon)$,
	which we derive just below,
	we conclude the desired bound \eqref{E:POINTWISEESTIMATEFORGTANCOMP}.
	To derive the small-data bound, we use 
	\eqref{E:LITTLEGDECOMPOSED}-\eqref{E:METRICPERTURBATIONFUNCTION},
	the small-data bound \eqref{E:LINFTYCOORADANGILARGEDATAALONGSIGMA1},
	and the bootstrap assumptions \eqref{E:PSIFUNDAMENTALC0BOUNDBOOTSTRAP}
	to obtain
	$
	\gtancomp^2|_{t=0} = g(\CoordAng,\CoordAng)|_{t=0} 
	= g_{22}|_{t=0}
	+ 
	\mathcal{O}(\mathring{\upepsilon})
	=
	1 + \mathcal{O}(\varepsilon)
	$,
	which implies the desired bound.
	\end{proof}

\begin{lemma}[\textbf{Comparison of the forms} $d \spherevol$ \textbf{and} $d \vartheta$]
\label{L:LINEVOLUMEFORMCOMPARISON}
Let $p=p(\vartheta)$ be a non-negative function of $\vartheta$. 
Then the following estimates hold
for $(t,u) \in [0,\Tboot) \times [0,U_0]$:
\begin{align} \label{E:LINEVOLUMEFORMCOMPARISON}
	(1 - C \varepsilon) \int_{\vartheta \in \mathbb{T}} p(\vartheta) d \vartheta
	\leq
	\int_{\ell_{t,u}} p(\vartheta) d \argspherevol{(t,u,\vartheta)}
	& \leq (1 + C \varepsilon) \int_{\vartheta \in \mathbb{T}} p(\vartheta) d \vartheta,
\end{align}
where $d \vartheta$ denotes the standard integration measure on $\mathbb{T}$.

Furthermore, let $p=p(u',\vartheta)$ be a non-negative function of $(u',\vartheta) \in [0,u] \times \mathbb{T}$
that \textbf{does not depend on $t$}.
Then for $s, t \in [0,\Tboot)$ and $u \in [0,U_0]$, 
we have:
\begin{align} \label{E:SIGMATVOLUMEFORMCOMPARISON}
	(1 - C \varepsilon)
	\int_{\Sigma_s^u} p \, d \tvol
	\leq
	\int_{\Sigma_t^u} p \, d \tvol
	& \leq 
	(1 + C \varepsilon)
	\int_{\Sigma_s^u} p \, d \tvol.
\end{align}

\end{lemma}

\begin{proof}
	See Sect.~\ref{SS:OFTENUSEDESTIMATES} for some comments on the analysis.
	From \eqref{E:RESCALEDVOLUMEFORMS} and inequality \eqref{E:POINTWISEESTIMATEFORGTANCOMP},
	we deduce that 
	$d \spherevol = (1 + \mathcal{O}(\varepsilon)) \, d \vartheta$,
	which yields \eqref{E:LINEVOLUMEFORMCOMPARISON}.
	\eqref{E:SIGMATVOLUMEFORMCOMPARISON} then follows from
	\eqref{E:LINEVOLUMEFORMCOMPARISON} and the fact that
	$d \tvol = d \argspherevol{(t,u,\vartheta)} du'$
	along $\Sigma_t^u$.
\end{proof}

\begin{lemma}[\textbf{Estimate for the norm} $\| \cdot \|_{L^2(\Sigma_t^u)}$ \textbf{of time-integrated functions}] 
\label{L:L2NORMSOFTIMEINTEGRATEDFUNCTIONS}
Let $f$ be a scalar-valued function on $\mathcal{M}_{\Tboot,U_0}$ and let
\begin{align} \label{E:BIGFISTIMEINETGRALOFLITTLEF}
	F(t,u,\vartheta) := \int_{t'=0}^t f(t',u,\vartheta) \, dt'.
\end{align}
Under the data-size and bootstrap assumptions 
of Sects.~\ref{SS:FLUIDVARIABLEDATAASSUMPTIONS}-\ref{SS:PSIBOOTSTRAP}
and the smallness assumptions of Sect.~\ref{SS:SMALLNESSASSUMPTIONS}, 
the following estimate holds
for $(t,u) \in [0,\Tboot) \times [0,U_0]$:
\begin{align} \label{E:L2NORMSOFTIMEINTEGRATEDFUNCTIONS}
	\| F \|_{L^2(\Sigma_t^u)} 
	& \leq (1 + C \varepsilon) 
		\int_{t'=0}^t 
			\| f \|_{L^2(\Sigma_{t'}^u)}
		\, dt'.
\end{align}
\end{lemma}

\begin{proof}
	Recall that 
	$\| F \|_{L^2(\Sigma_t^u)}
	: = 
				\left\lbrace
				\int_{u'=0}^u
					\int_{\ell_{t,u'}}
						F^2(t,u',\vartheta) 
					\, d \spherevol
				\, du'
				\right\rbrace^{1/2}
		$.
	Using the estimate \eqref{E:LINEVOLUMEFORMCOMPARISON},
	we may replace $d \spherevol$ in the previous formula 
	with the standard integration measure
	$d \vartheta$ up to an overall multiplicative error factor of $1 + \mathcal{O}(\varepsilon)$.
	The desired estimate \eqref{E:L2NORMSOFTIMEINTEGRATEDFUNCTIONS} follows from this estimate and
	from applying Minkowski's inequality for integrals to equation \eqref{E:BIGFISTIMEINETGRALOFLITTLEF}.
\end{proof}

\subsection{The coerciveness of the fundamental square-integral controlling quantities}
\label{SS:COERCIVENESSOFCONTROLLINGQUANTITIES}
In this section, we quantify the coercive nature of the fundamental 
$L^2$-controlling quantities.

We start by quantifying the coercive nature of the 
spacetime integrals of Def.~\ref{D:COERCIVEINTEGRAL}.

\begin{lemma}[\textbf{Strength of the coercive spacetime integral}]
\label{L:KEYSPACETIMECOERCIVITY}
Let $\mathbf{1}_{\lbrace \upmu \leq 1/4 \rbrace}$
denote the characteristic function of the spacetime
subset 
$
\displaystyle
\lbrace (t,u,\vartheta) \in [1,\infty) \times [0,1] \times \mathbb{T}
	\ | \ 
\upmu(t,u,\vartheta) \leq 1/4 \rbrace
$.
Under the data-size and bootstrap assumptions 
of Sects.~\ref{SS:FLUIDVARIABLEDATAASSUMPTIONS}-\ref{SS:PSIBOOTSTRAP}
and the smallness assumptions of Sect.~\ref{SS:SMALLNESSASSUMPTIONS}, 
the following lower bound holds
for $(t,u) \in [0,\Tboot) \times [0,U_0]$:
\begin{align} \label{E:KEYSPACETIMECOERCIVITY}
		\coercivespacetime[\Psi](t,u) 
		& \geq 
		\frac{1}{8}
		\TranminusdatasizeWithFactor
		\int_{\mathcal{M}_{t,u}}
			\mathbf{1}_{\lbrace \upmu \leq 1/4 \rbrace}
			\left|
				\angdiff \Psi
			\right|^2
		\, d \vol.
	\end{align}
\end{lemma}

\begin{proof}
	Inequality \eqref{E:KEYSPACETIMECOERCIVITY}
	follows from definition \eqref{E:COERCIVESPACETIMEDEF}
	and the estimate \eqref{E:SMALLMUIMPLIESLMUISNEGATIVE}.
\end{proof}

We now quantify the coercivity of the
fundamental $L^2$-controlling quantities from Def.~\ref{D:MAINCOERCIVEQUANT}.

\begin{lemma}[\textbf{The coercivity of the fundamental controlling quantities}]
	\label{L:COERCIVENESSOFCONTROLLING}
	Assume that $1 \leq N \leq 20$ and $0 \leq M \leq \min \lbrace N-1, 1 \rbrace$.
	Under the assumptions of Lemma~\ref{L:KEYSPACETIMECOERCIVITY},
	the following lower bounds hold
	for $(t,u) \in [0,\Tboot) \times [0,U_0]$:
	\begin{subequations}
	\begin{align} \label{E:WAVECOERCIVENESSOFCONTROLLING}
		\totmax{N}(t,u)
		& \geq 
			\max_{\Psi \in \lbrace \Densrenormalized  - v^1,v^1,v^2 \rbrace}
			\Big\lbrace
				\frac{1}{2}
				\left\|
					\sqrt{\upmu} \Lunit \Fullset_{\ast}^{N;M} \Psi
				\right\|_{L^2(\Sigma_t^u)}^2,
					\,
				\left\|
					\Rad \Fullset_{\ast}^{N;M} \Psi
				\right\|_{L^2(\Sigma_t^u)}^2,
					\,
				\frac{1}{2}
				\left\|
					\sqrt{\upmu} \angdiff \Fullset_{\ast}^{N;M} \Psi
				\right\|_{L^2(\Sigma_t^u)}^2,
					\\
			& \ \
				\left\|
					\Lunit \Fullset_{\ast}^{N;M} \Psi
				\right\|_{L^2(\mathcal{P}_u^t)}^2,
					\,
				\left\|
					\sqrt{\upmu} \angdiff \Fullset_{\ast}^{N;M} \Psi
				\right\|_{L^2(\mathcal{P}_u^t)}^2
				\Big\rbrace,
					\notag 
						\\
		\totmax{N}(t,u)
		& \geq 
			\max
			\Big\lbrace
				\frac{1}{8}
				\left\|
					\sqrt{\upmu} \Lunit \Fullset_{\ast}^{N;M} \Densrenormalized
				\right\|_{L^2(\Sigma_t^u)}^2,
					\,
				\frac{1}{4}
				\left\|
					\Rad \Fullset_{\ast}^{N;M} \Densrenormalized
				\right\|_{L^2(\Sigma_t^u)}^2,
					\,
				\frac{1}{8}
				\left\|
					\sqrt{\upmu} \angdiff \Fullset_{\ast}^{N;M} \Densrenormalized
				\right\|_{L^2(\Sigma_t^u)}^2,
					\label{E:DENSITYWAVECOERCIVENESSOFCONTROLLING} 
						\\
			& \ \
				\frac{1}{4}
				\left\|
					\Lunit \Fullset_{\ast}^{N;M} \Densrenormalized
				\right\|_{L^2(\mathcal{P}_u^t)}^2,
					\,
				\frac{1}{4}
				\left\|
					\sqrt{\upmu} \angdiff \Fullset_{\ast}^{N;M} \Densrenormalized
				\right\|_{L^2(\mathcal{P}_u^t)}^2
				\Big\rbrace,
					\notag 
	\end{align}
	\end{subequations}

	\begin{align} \label{E:EASYWAVECOERCIVENESSOFCONTROLLING}
		\easytotmax{N}(t,u)
		& \geq 
			\max_{\Psi \in \lbrace \Densrenormalized  - v^1,v^2 \rbrace}
			\Big\lbrace
				\frac{1}{2}
				\left\|
					\sqrt{\upmu} \Lunit \Fullset_{\ast}^{N;M} \Psi
				\right\|_{L^2(\Sigma_t^u)}^2,
					\,
				\left\|
					\Rad \Fullset_{\ast}^{N;M} \Psi
				\right\|_{L^2(\Sigma_t^u)}^2,
					\,
				\frac{1}{2}
				\left\|
					\sqrt{\upmu} \angdiff \Fullset_{\ast}^{N;M} \Psi
				\right\|_{L^2(\Sigma_t^u)}^2,
					\\
			& \ \
				\left\|
					\Lunit \Fullset_{\ast}^{N;M} \Psi
				\right\|_{L^2(\mathcal{P}_u^t)}^2,
					\,
				\left\|
					\sqrt{\upmu} \angdiff \Fullset_{\ast}^{N;M} \Psi
				\right\|_{L^2(\mathcal{P}_u^t)}^2
				\Big\rbrace.
					\notag 
		\end{align}
	In addition, if $N \leq 21$, then
	\begin{align}
		\VorttotTanmax{N}(t,u)
			& \geq 
			\max
			\Big\lbrace
				\left\|
					\sqrt{\upmu} \Tanset^N \Vortrenormalized
				\right\|_{L^2(\Sigma_t^u)}^2,
					\,
				\left\|
					\Tanset^N \Vortrenormalized
				\right\|_{L^2(\mathcal{P}_u^t)}^2
			\Big\rbrace.
			\label{E:VORTICITYCOERCIVENESSOFCONTROLLING} 
	\end{align}

	Moreover, if $1 \leq N \leq 20$ and $0 \leq M \leq 1$,
	then
	\begin{subequations}
	\begin{align}
		\left\|
			\threePsi
		\right\|_{L^2(\Sigma_t^u)}^2
		& \leq
			C \mathring{\upepsilon}^2
			+ 
			C \totmax{1}(t,u),
				\label{E:LOWESTORDERL2ESTIMATELOSSOFONEDERIVATIVE} \\
		\left\|
			\Fullset_{\ast}^{N;M} \threePsi
		\right\|_{L^2(\Sigma_t^u)}^2,
			\,
		\left\|
			\Fullset_{\ast}^{N;M} \threePsi
		\right\|_{L^2(\ell_{t,u})}^2
		& \leq
			C \mathring{\upepsilon}^2
			+ 
			C \totmax{N}(t,u),
			 \label{E:L2ESTIMATELOSSOFONEDERIVATIVE} \\
		\left\|
			\Rad \threePsi
		\right\|_{L^2(\Sigma_t^u)}^2
		& \leq 
		C
		\left\|
			\Rad \threePsi
		\right\|_{L^2(\Sigma_0^u)}^2
		+ 
		C \mathring{\upepsilon}^2
		+ 
		C \totmax{1}(t,u),
		\label{E:RADPSIL2ESTIMATELOSSOFONEDERIVATIVE}
			\\
		\left\|
			\Rad \Rad \threePsi
		\right\|_{L^2(\Sigma_t^u)}^2
		& \leq 
		C 
		\left\|
			\Rad \Rad \threePsi
		\right\|_{L^2(\Sigma_0^u)}^2
		+ 
		C \mathring{\upepsilon}^2
		+ 
		C \totmax{2}(t,u).
		\label{E:TWORADPSIL2ESTIMATELOSSOFONEDERIVATIVE}
	\end{align}
	\end{subequations}
	Finally, if $N \leq 20$, then
	\begin{align}
		\left\|
			\Tanset^{\leq N} \Vortrenormalized
		\right\|_{L^2(\ell_{t,u})}^2
			& \leq
			C \mathring{\upepsilon}^2
			+ 
			C \Vorttotmax{\leq N+1}(t,u).
			 \label{E:VORTICITYL2ESTIMATEALONGELLTU}
	\end{align}

\end{lemma}

\begin{remark}
	\label{R:SHARPCONSTANTSMATTER}
	The constants 
	$1$ and $\frac{1}{2}$
	in front of the term
	$
	\left\|
		\Rad \Fullset_{\ast}^{N;M} \Psi
	\right\|_{L^2(\Sigma_t^u)}^2
	$
	and the term
	$
	\frac{1}{2}
				\left\|
					\sqrt{\upmu} \angdiff \Fullset_{\ast}^{N;M} \Psi
				\right\|_{L^2(\Sigma_t^u)}^2
	$
	on RHS~\eqref{E:WAVECOERCIVENESSOFCONTROLLING}
	influence the blowup-rate of our top-order energy estimates.
	In turn, this affects the number of derivatives that we need to close
	our estimates.
\end{remark}

\begin{proof}[Proof of Lemma~\ref{L:COERCIVENESSOFCONTROLLING}]
	The estimates stated in 
	\eqref{E:VORTICITYCOERCIVENESSOFCONTROLLING}
	follow easily from Defs.~\ref{D:MAINCOERCIVEQUANT} and \ref{D:ENERGYFLUX}.

	The estimates stated in \eqref{E:WAVECOERCIVENESSOFCONTROLLING}
	follow easily from 
	Lemma~\ref{L:ORDERZEROCOERCIVENESS},
	Def.~\ref{D:MAINCOERCIVEQUANT},
	and Young's inequality.

	The estimates stated in \eqref{E:EASYWAVECOERCIVENESSOFCONTROLLING}
	follow easily from \eqref{E:WAVECOERCIVENESSOFCONTROLLING} and the triangle
	inequality.

	We now prove \eqref{E:L2ESTIMATELOSSOFONEDERIVATIVE}.
	We first note that the estimates for 
	$
	\displaystyle
	\left\|
		\Fullset_{\ast}^{N;M} \threePsi
	\right\|_{L^2(\Sigma_t^u)}^2
	$
	follow easily from integrating the estimates for
	$
	\displaystyle
		\left\|
			\Fullset_{\ast}^{N;M} \threePsi
		\right\|_{L^2(\ell_{t,u})}^2
	$
	with respect to $u$.
	Hence, it suffices to prove the estimates for
	$
	\displaystyle
	\left\|
		\Fullset_{\ast}^{N;M}
	\right\|_{L^2(\ell_{t,u})}^2
	$.
	Our proof is based on the identity \eqref{E:UDERIVATIVEOFLINEINTEGRAL}.
	To proceed, we first use 
	\eqref{E:CONNECTIONBETWEENANGLIEOFGSPHEREANDDEFORMATIONTENSORS},
	the estimate \eqref{E:POINTWISEESTIMATESFORGSPHEREANDITSNOSPECIALSTRUCTUREDERIVATIVES}
	and the $L^{\infty}$ estimates of Prop.~\ref{P:IMPROVEMENTOFAUX}
	to bound the factor 
	$
	\displaystyle
	(1/2) \mytr \angdeform{\Rad}
	$ 
	in \eqref{E:UDERIVATIVEOFLINEINTEGRAL} as follows:
	$
	\displaystyle
	(1/2)
	\left|
		\mytr \angdeform{\Rad}
	\right|
	\lesssim 
	\left|
		\angLie_{\Rad} \gsphere
	\right|
	\lesssim 1
	$.
	Using this estimate,
	the identity \eqref{E:UDERIVATIVEOFLINEINTEGRAL}
	with $f = (\Fullset_{\ast}^{N;M} \Psi)^2$, 
	and Young's inequality,
	we deduce that
	\begin{equation} \label{E:UDERIVATIVEOFSTUINTEGRALGRONWALLREADY}
	\begin{split}
	\left\|
		\Fullset_{\ast}^{N;M} \Psi
	\right\|_{L^2(\ell_{t,u})}^2
	 \leq &
		\left\|
			\Fullset_{\ast}^{N;M} \Psi
		\right\|_{L^2(\ell_{t,0})}^2
		+
		\int_{u'=0}^u
			\left\|
				\Rad \Fullset_{\ast}^{N;M} \Psi
			\right\|_{L^2(\ell_{t,u'})}^2
		\, du'\\
			&+ 
		C
		\int_{u'=0}^u
			\left\|
				\Fullset_{\ast}^{N;M} \Psi
			\right\|_{L^2(\ell_{t,u'})}^2
		\, du'.
\end{split}
\end{equation}
From \eqref{E:UDERIVATIVEOFSTUINTEGRALGRONWALLREADY},
the smallness assumption \eqref{E:SMALLDATAASSUMPTIONSALONGELLT0},
and Gronwall's inequality, 
we deduce
\begin{align} \label{E:UDERIVATIVEOFSTUINTEGRALGRONWALLED}
	\left\|
		\Fullset_{\ast}^{N;M}
	\right\|_{L^2(\ell_{t,u})}^2
	& \leq 
		e^{cu} \mathring{\upepsilon}^2
		+
		e^{cu}
		\int_{u'=0}^u
			\left\|
				\Rad \Fullset_{\ast}^{N;M} \Psi
			\right\|_{L^2(\ell_{t,u'})}^2
		\, du'
			\\
	& = 
			e^{cu} \mathring{\upepsilon}^2
			+
			C e^{cu}
			\left\|
				\Rad \Fullset_{\ast}^{N;M} \Psi
			\right\|_{L^2(\Sigma_t^u)}^2
		\leq 
			C \mathring{\upepsilon}^2
			+
			C 
			\left\|
				\Rad \Fullset_{\ast}^{N;M} \Psi
			\right\|_{L^2(\Sigma_t^u)}^2.
			\notag
\end{align}
The desired bound 
for 
$
\left\|
	\Fullset_{\ast}^{N;M}\Psi
\right\|_{L^2(\ell_{t,u})}^2
$
now follows from
\eqref{E:UDERIVATIVEOFSTUINTEGRALGRONWALLED}
and the already proven estimate
\eqref{E:WAVECOERCIVENESSOFCONTROLLING}
for 
$
\left\|
	\Rad \Fullset_{\ast}^{N;M} \Psi
\right\|_{L^2(\Sigma_t^u)}^2
$.

	We now prove \eqref{E:VORTICITYL2ESTIMATEALONGELLTU}.
	We first use \eqref{E:LDERIVATIVEOFLINEINTEGRAL}
	with $f = (\Tanset^N \Vortrenormalized)^2$,
	the estimate \eqref{E:TWORADIALCHICOMMUTEDLINFINITY},
	and Young's inequality to deduce
	\begin{align} \label{E:VORTICITYLTUESTIMATEDIFFERENTIALVERSION}
		\frac{\partial}{\partial t}
		\left\|
			\Tanset^N \Vortrenormalized
		\right\|_{L^2(\ell_{t,u})}^2
		& \leq
			\left\|
				\Lunit \Tanset^N \Vortrenormalized
			\right\|_{L^2(\ell_{t,u})}^2
			+
			C
			\left\|
				\Tanset^N \Vortrenormalized
			\right\|_{L^2(\ell_{t,u})}^2.
	\end{align}
	Integrating \eqref{E:VORTICITYLTUESTIMATEDIFFERENTIALVERSION}
	from the initial time $0$ to time $t$,
	using Gronwall's inequality, 
	and using the small data assumption \eqref{E:DATAASSUMPTIONSALONGELL1U},
	we obtain the desired bound \eqref{E:VORTICITYL2ESTIMATEALONGELLTU}
	as follows:
	\begin{align} \label{E:VORTICITYLTUESTIMATE}
		\left\|
			\Tanset^N \Vortrenormalized
		\right\|_{L^2(\ell_{t,u})}^2
		& \lesssim
				\left\|
			\Tanset^N \Vortrenormalized
			\right\|_{L^2(\ell_{0,u})}^2
				+
				\int_{s=0}^t
					\left\|
						\Lunit \Tanset^N \Vortrenormalized
					\right\|_{L^2(\ell_{s,u})}^2
				\, ds
			\lesssim
			\mathring{\upepsilon}^2
			+
			\VorttotTanmax{N+1}(t,u).
		\end{align}

	To prove \eqref{E:LOWESTORDERL2ESTIMATELOSSOFONEDERIVATIVE},
	we first use the fundamental theorem of calculus to express
	\begin{align} \label{E:SIMPLEFTCIDENTITY}
		\threePsi(t,u,\vartheta)
		&
		= 
		\threePsi(0,u,\vartheta)
		+
		\int_{t'=0}^t
			\Lunit \threePsi(t',0,\vartheta)
		\, dt'.
	\end{align}
	From \eqref{E:SIMPLEFTCIDENTITY},
	\eqref{E:L2NORMSOFTIMEINTEGRATEDFUNCTIONS},
	and the estimate
	$
	\displaystyle
	\left\| 
		\Lunit \threePsi
	\right\|_{L^2(\Sigma_t^u)} 
	\lesssim
	\sqrt{\totmax{1}}(t,u)
	+
	\mathring{\upepsilon}
	$
	(which is a particular case of the already proven bound \eqref{E:L2ESTIMATELOSSOFONEDERIVATIVE}),
	we deduce
	\begin{align} \label{E:L2NORMOFSIMPLEFTCIDENTITY}
		\| \threePsi \|_{L^2(\Sigma_t^u)} 
		&
		\lesssim 
		\| \threePsi(0,\cdot) \|_{L^2(\Sigma_t^u)} 
		+
		\int_{t'=0}^t
			\left\| 
				\Lunit \threePsi
			\right\|_{L^2(\Sigma_{t'}^u)} 
		\, dt'
			\\
		&
		\lesssim 
		\| \threePsi(0,\cdot) \|_{L^2(\Sigma_t^u)} 
		+
		\int_{t'=0}^t
			\sqrt{\totmax{1}}(t',u)
		\, dt'
		+ 
		\mathring{\upepsilon}.
		\notag
	\end{align}
	From \eqref{E:SIGMATVOLUMEFORMCOMPARISON} with $s=0$ and the 
	smallness assumption \eqref{E:L2SMALLDATAASSUMPTIONSALONGSIGMA0},
	we deduce that 
	$\| \threePsi(0,\cdot) \|_{L^2(\Sigma_t^u)}
	\lesssim
	\| \threePsi \|_{L^2(\Sigma_0^u)}
	\lesssim \mathring{\upepsilon}
	$.
	Moreover, using 
	the fact that $\totmax{1}$
	is increasing in its arguments,
	we deduce that
	$
	\displaystyle
	\int_{t'=0}^t
		\sqrt{\totmax{1}}(t',u)
	\, dt'
	\lesssim
	\totmax{1}^{1/2}(t,u)
	$
	.
	Inserting the above estimates into RHS~\eqref{E:L2NORMOFSIMPLEFTCIDENTITY},
	we arrive at the desired bound \eqref{E:LOWESTORDERL2ESTIMATELOSSOFONEDERIVATIVE}.

	The remaining estimates 
	\eqref{E:RADPSIL2ESTIMATELOSSOFONEDERIVATIVE}
	and
	\eqref{E:TWORADPSIL2ESTIMATELOSSOFONEDERIVATIVE}
	follow from arguments similar to the ones
	we used to prove \eqref{E:LOWESTORDERL2ESTIMATELOSSOFONEDERIVATIVE}
	and we omit those details.
\end{proof}

\section{Sobolev embedding}
\label{SS:SOBOLEVEMBEDDING}
Our main goal in this section is to prove Cor.~\ref{C:PSILINFTYINTERMSOFENERGIES},
which is the Sobolev embedding result that we will use 
to improve the fundamental bootstrap assumptions
\eqref{E:PSIFUNDAMENTALC0BOUNDBOOTSTRAP}-\eqref{E:VORTFUNDAMENTALC0BOUNDBOOTSTRAP}.

\begin{lemma}[\textbf{Sobolev embedding along} $\ell_{t,u}$]
	\label{L:SOBOLEV}
Under the data-size and bootstrap assumptions 
of Sects.~\ref{SS:FLUIDVARIABLEDATAASSUMPTIONS}-\ref{SS:PSIBOOTSTRAP}
and the smallness assumptions of Sect.~\ref{SS:SMALLNESSASSUMPTIONS}, 
the following estimate holds for 
scalar-valued functions $f$ defined on $\ell_{t,u}$
for $(t,u) \in [0,\Tboot) \times [0,U_0]$:
\begin{align} \label{E:SOBOLEV}
		\left\|
			f
		\right\|_{L^{\infty}(\ell_{t,u})}
		& \leq 
			C
		\left\|
			\GeoAng^{\leq 1} f
		\right\|_{L^2(\ell_{t,u})}.
	\end{align}
\end{lemma}

\begin{proof}
	Standard Sobolev embedding yields that
	$\left\|
			f
	\right\|_{L^{\infty}(\mathbb{T})}
	\leq 
	C
	\left\|
		\CoordAng^{\leq 1} f
	\right\|_{L^2(\mathbb{T})}$,
	where the integration measure defining $\| \cdot \|_{L^2(\mathbb{T})}$ is $d \vartheta$.
	From 
	Def.~\ref{D:METRICANGULARCOMPONENT},
	\eqref{E:GEOANGPOINTWISE},
	Lemma~\ref{L:POINTWISEESTIMATEFORGTANCOMP},
	and the $L^{\infty}$ estimates of Prop.~\ref{P:IMPROVEMENTOFAUX},
	we find that
	$
	|\CoordAng^{\leq 1} f| = (1 + \mathcal{O}(\varepsilon))|\GeoAng^{\leq 1} f|
	$.
	From these estimates and Lemma~\ref{L:LINEVOLUMEFORMCOMPARISON},
	we conclude the desired estimate \eqref{E:SOBOLEV}.

\end{proof}

\begin{corollary}[$L^{\infty}$ \textbf{bounds for} $\threePsi$ 
\textbf{and} $\Vortrenormalized$
\textbf{in terms of the fundamental controlling quantities}]
\label{C:PSILINFTYINTERMSOFENERGIES}
Under the assumptions of Lemma~\ref{L:SOBOLEV}, 
the following estimates hold for
$(t,u) \in [0,\Tboot) \times [0,U_0]$:
\begin{subequations}
\begin{align}  \label{E:PSILINFTYINTERMSOFENERGIES}
	\left\|
		\Fullset_{\ast}^{\leq 13;\leq 1} \threePsi
	\right\|_{L^{\infty}(\Sigma_t^u)}
	& \lesssim 
		\totmax{[1,14]}^{1/2}(t,u)
		+ 
		\mathring{\upepsilon},
			\\
	\left\|
		\Tanset^{\leq 13} \Vortrenormalized
	\right\|_{L^{\infty}(\Sigma_t^u)}
	& \lesssim 
		\VorttotTanmax{\leq 15}^{1/2}(t,u)
		+ 
		\mathring{\upepsilon}.
		\label{E:VORTICITYLINFTYINTERMSOFENERGIES}
\end{align}
\end{subequations}

\end{corollary}

\begin{proof}
	See Sect.~\ref{SS:OFTENUSEDESTIMATES} for some comments on the analysis.
	The bound
	$
	\left\|
		\Fullset_{\ast}^{[1,13];\leq 1} \threePsi
	\right\|_{L^{\infty}(\Sigma_t^u)}
		\lesssim 
		\totmax{[1,14]}^{1/2}(t,u) + \mathring{\upepsilon}
	$
	follows from \eqref{E:L2ESTIMATELOSSOFONEDERIVATIVE}
	and Lemma~\ref{L:SOBOLEV}.
	This estimate implies in particular the pointwise bound
	$\left|
		\Lunit \threePsi
	 \right|
		\lesssim 
		\totmax{[1,14]}^{1/2}(t,u)
		+
		\mathring{\upepsilon}
	$.
	Integrating
	along the integral curves of
	$\Lunit$ 
	as in \eqref{E:SIMPLEFTCID}
	and using this bound 
	and the small-data assumption
	$
	\left\|
		\threePsi
	\right\|_{L^{\infty}(\Sigma_0^u)}
	\leq \mathring{\upepsilon}
	$
	(see \eqref{E:LINFTYSMALLDATAASSUMPTIONSALONGSIGMA0}),
	we deduce that
	$
	\left\|
		\threePsi
	\right\|_{L^{\infty}(\Sigma_t^u)}
	\leq 
		C \totmax{[1,14]}^{1/2}(t,u)
		+ 
		C \mathring{\upepsilon}
	$.
	We have thus proved \eqref{E:PSILINFTYINTERMSOFENERGIES}.

	The estimate \eqref{E:VORTICITYLINFTYINTERMSOFENERGIES}
	is
	a direct consequence of Lemma~\ref{L:SOBOLEV}
	and inequality \eqref{E:VORTICITYL2ESTIMATEALONGELLTU}.
\end{proof}

%

\section{Pointwise estimates for the error integrands}
\label{S:POINTWISEESTIMATESFORWAVEEQUATIONERRORTERMS}
\setcounter{equation}{0}
In order to derive a priori estimates 
for the fundamental $L^2$-controlling quantities of
Defs.~\ref{D:MAINCOERCIVEQUANT} and \ref{D:COERCIVEINTEGRAL},
we must first obtain pointwise estimates for the error terms
in the energy identities corresponding to the
wave equations verified by
$\Fullset_*^{[1,20];\leq 1} (\Densrenormalized  - v^1)$,
$\Fullset_*^{[1,20];\leq 1} v^1$,
$\Fullset_*^{[1,20];\leq 1} v^2$,
and the transport equations verified by 
$\Tanset^{\leq 21} \Vortrenormalized$.
By ``energy identities,'' we mean the ones provided by
Props.~\ref{P:DIVTHMWITHCANCELLATIONS}
and \ref{P:ENERGYIDENTITYRENORMALIZEDVORTICITY}.
In this section, we derive these pointwise estimates.
The error terms consist of the following three types, 
ordered in increasing difficulty:
\textbf{i)} Error terms generated by differentiating the inhomogeneous
terms on RHSs \eqref{E:VELOCITYWAVEEQUATION}-\eqref{E:RENORMALIZEDDENSITYWAVEEQUATION};
\textbf{ii)} Error terms corresponding to the last integral on 
RHS~\eqref{E:ENERGYIDENTITYRENORMALIZEDVORTICITY};
\textbf{ii')} Error terms corresponding to the deformation tensor of the
	multiplier vectorfields, which correspond to the 
	last integral on RHS~\eqref{E:E0DIVID};
and \textbf{iii)} Error terms generated by the commutator terms of the form
$[\upmu \square_g,\Fullset_*^{N;\leq 1}] \Psi$ and $[\upmu \Transport,\Tanset^N] \Vortrenormalized$
and their derivatives up to top-order.

We prove the two main propositions 
in Sects.~\ref{SS:PROOFOFPROPWAVEIDOFKEYDIFFICULTENREGYERRORTERMS}
and \ref{SS:PROOFOFPROPVORTICITYIDOFKEYDIFFICULTENREGYERRORTERMS}.
The rest of Sect.~\ref{S:POINTWISEESTIMATESFORWAVEEQUATIONERRORTERMS} consists
of preliminary estimates.

\subsection{Harmless terms}
\label{SS:HARMLESSTERMS}
We start by defining error terms of type
$Harmless_{(Wave)}^{\leq N}$,
which appear in the energy estimates for the wave variables,
and of type
$Harmless_{(Vort)}^{\leq N}$,
which appear in the energy estimates for the specific vorticity.
These terms have a negligible effect on the dynamics,
even near the shock. Most error terms that we encounter 
are of these types.

\begin{definition}[\textbf{Harmless terms}]
	\label{D:HARMLESSTERMS}
	$Harmless_{(Wave)}^{\leq N}$ 
	and
	$Harmless_{(Vort)}^{\leq N}$
	denote any terms such that 
	under the data-size and bootstrap assumptions 
	of Sects.~\ref{SS:FLUIDVARIABLEDATAASSUMPTIONS}-\ref{SS:PSIBOOTSTRAP}
	and the smallness assumptions of Sect.~\ref{SS:SMALLNESSASSUMPTIONS}, 
	the following bound holds on
	$\mathcal{M}_{\Tboot,U_0}$,
	where $1 \leq N \leq 20$ in \eqref{E:WAVEHARMESSTERMPOINTWISEESTIMATE}
	and $N \leq 21$ in \eqref{E:VORTICITYHARMESSTERMPOINTWISEESTIMATE}
	(see Sect.~\ref{SS:STRINGSOFCOMMUTATIONVECTORFIELDS} regarding the vectorfield operator notation):
	\begin{subequations}
	\begin{align} \label{E:WAVEHARMESSTERMPOINTWISEESTIMATE}
		\left| 
			Harmless_{(Wave)}^{\leq N}
		\right|
		& \lesssim 
			\left|
				\Fullset_{\ast}^{[1,N+1];\leq 2} \threePsi
			\right|
			+
			\left|
				\Fullset_{\ast}^{[1,N];\leq 2} \GdVar
			\right|
			+
			\left|
				\Fullset_{\ast \ast}^{[1,N];\leq 1} \BadVar
			\right|
			+
			\left|
				\Fullset^{\leq N;\leq 1} \Vortrenormalized
			\right|,
				\\
		\left| 
			Harmless_{(Vort)}^{\leq N}
		\right|
		& \lesssim 
			\varepsilon
			\left|
				\Fullset_{\ast}^{[1,N];\leq 2} \threePsi
			\right|
			+
			\varepsilon
			\left|
				\Fullset_{\ast}^{[1,N-1];\leq 2} \GdVar
			\right|
			+
			\varepsilon
			\left|
				\Fullset_{\ast \ast}^{[1,N-1];\leq 1} \BadVar
			\right|
			+
			\left|
				\Tanset^{\leq N} \Vortrenormalized
			\right|.
			\label{E:VORTICITYHARMESSTERMPOINTWISEESTIMATE}
	\end{align}
	\end{subequations}
	By definition, the first three terms on RHS~\eqref{E:VORTICITYHARMESSTERMPOINTWISEESTIMATE} are
	absent when $N=0$
	and the second and third terms on RHS~\eqref{E:VORTICITYHARMESSTERMPOINTWISEESTIMATE} 
	are absent when $N=1$.

\end{definition}

\begin{remark}[\textbf{The role of the factors} $\varepsilon$]
	\label{R:ROLEOFEPSIONFACTORINHARMLESSVORTICITYTERMS}
	The smallness factors $\varepsilon$
	on RHS~\eqref{E:VORTICITYHARMESSTERMPOINTWISEESTIMATE} 
	allow us to derive
	(see Lemma~\ref{L:VORTICITYAPRIORIENERGYESTIMATES}),
	via a bootstrap argument,
	energy estimates for 
	$\Vortrenormalized$
	without having to simultaneously derive energy estimates for $\threePsi$.
	This allows us to obtain energy estimates for $\Vortrenormalized$
	in a much simpler way compared to the energy estimates for $\threePsi$.
\end{remark}

\subsection{Identification of the difficult error terms in the commuted equations}
\label{SS:IDOFDIFFICULTERRORTERMS}
In the next proposition, 
which we prove in Sect.~\ref{SS:PROOFOFPROPWAVEIDOFKEYDIFFICULTENREGYERRORTERMS},
we identify the main error terms in the inhomogeneous
wave equations verified by the higher-order versions of
$\Densrenormalized  - v^1$,
$v^1$,
and $v^2$.
The main terms will require careful treatment in the energy estimates,
while the terms denoted by
$Harmless_{(Wave)}^{\leq N}$
will be easy to control.

\begin{proposition}[\textbf{Identification of the key difficult error term factors in the commuted wave equations}]
\label{P:WAVEIDOFKEYDIFFICULTENREGYERRORTERMS}
	Assume that $1 \leq N \leq 20$
	and recall that $\GeoAngFlatRadComponent$ is the scalar-valued
	appearing in Lemma~\ref{L:GEOANGDECOMPOSITION}.
	Under the data-size and bootstrap assumptions 
	of Sects.~\ref{SS:FLUIDVARIABLEDATAASSUMPTIONS}-\ref{SS:PSIBOOTSTRAP}
	and the smallness assumptions of Sect.~\ref{SS:SMALLNESSASSUMPTIONS}, 
	the following pointwise estimates hold
	for $i=1,2$ on $\mathcal{M}_{\Tboot,U_0}$
	(see Sect.~\ref{SS:STRINGSOFCOMMUTATIONVECTORFIELDS} regarding the vectorfield operator notation):
\begin{subequations}
\begin{align}
	\upmu \square_g (\GeoAng^{N-1} \Lunit v^i)
	& = (\angdiffuparg{\#} v^i) \cdot (\upmu \angdiff \GeoAng^{N-1} \mytr \upchi)
			\label{E:WAVELISTHEFIRSTCOMMUTATORIMPORTANTTERMS} \\
	& \ \
		+ [ia] \upmu (\exp \Densrenormalized) \Speed^2 (g_{ab} \Radunit^b) 
					\GeoAng^{N-1} \Lunit \Lunit \Vortrenormalized
					\notag \\
	& \ \
		- [ia] \upmu (\exp \Densrenormalized) \Speed^2
			\left(
				\frac{g_{ab} \GeoAng^b}{g_{cd} \GeoAng^c \GeoAng^d}
			\right)
			\GeoAng^{N-1} \Lunit \GeoAng \Vortrenormalized
				\notag \\
	& \ \
			+ Harmless_{(Wave)}^{\leq N},
			\notag \\
	\upmu \square_g (\GeoAng^N v^i)
	& = (\Rad v^i) \GeoAng^N \mytr \upchi
			+ \GeoAngFlatRadComponent (\angdiffuparg{\#} v^i) \cdot (\upmu \angdiff \GeoAng^{N-1} \mytr \upchi)
			\label{E:WAVEGEOANGISTHEFIRSTCOMMUTATORIMPORTANTTERMS} \\
		& \ \
			+ [ia] \upmu (\exp \Densrenormalized) \Speed^2 (g_{ab} \Radunit^b) 
					\GeoAng^N  \Lunit \Vortrenormalized
					\notag \\
		& \ \
			- [ia] \upmu (\exp \Densrenormalized) \Speed^2
			\left(
				\frac{g_{ab} \GeoAng^b}{g_{cd} \GeoAng^c \GeoAng^d}
			\right)
			\GeoAng^{N+1} \Vortrenormalized
				\notag \\
	& \ \
				+ Harmless_{(Wave)}^{\leq N}.
				 \notag
	\end{align}

	Moreover, if $2 \leq N \leq 20$, then
	\begin{align}
	\upmu \square_g (\GeoAng^{N-1} \Rad v^i)
	& = (\Rad v^i) \GeoAng^{N-1} \Rad \mytr \upchi
			- (\upmu \angdiffuparg{\#} v^i) \cdot (\upmu \angdiff \GeoAng^{N-1} \mytr \upchi)
				\label{E:WAVERADISTHEFIRSTCOMMUTATORIMPORTANTTERMS}  \\
	& \ \
				- 
				[ia] \upmu^2 (\exp \Densrenormalized) \Speed^2 (g_{ab} \Radunit^b) 
					\GeoAng^{N-1} \Lunit \Lunit \Vortrenormalized
				\notag	\\
	& \ \
		+ [ia] \upmu^2 (\exp \Densrenormalized) \Speed^2
			\left(
				\frac{g_{ab} \GeoAng^b}{g_{cd} \GeoAng^c \GeoAng^d}
			\right)
			\GeoAng^N \Lunit \Vortrenormalized
			\notag	\\
	& \ \
			+ Harmless_{(Wave)}^{\leq N}.
			\notag
	\end{align}

	Moreover, if $2 \leq N \leq 20$
	and $\Fullset^{N-1;1}$ contains
	exactly one factor of $\Rad$ with all other factors equal to
	$\GeoAng$, then we have
	\begin{align}
	\upmu \square_g (\Fullset^{N-1;1} \Lunit v^i)
	& = (\angdiffuparg{\#} v^i) \cdot (\upmu \angdiff \GeoAng^{N-2} \Rad \mytr \upchi)
		\label{E:SECONDWAVELISTHEFIRSTCOMMUTATORIMPORTANTTERMS}  \\
	& \ \
		- [ia] \upmu^2 (\exp \Densrenormalized) \Speed^2 (g_{ab} \Radunit^b) 
			\GeoAng^{N-2} \Lunit \Lunit \Lunit \Vortrenormalized
				\notag \\
	& \ \
		+ [ia] \upmu^2 (\exp \Densrenormalized) \Speed^2
			\left(
				\frac{g_{ab} \GeoAng^b}{g_{cd} \GeoAng^c \GeoAng^d}
			\right)
			\GeoAng^{N-1} \Lunit \Lunit \Vortrenormalized
				\notag \\
	& \ \
			+ Harmless_{(Wave)}^{\leq N},
			\notag \\
	\upmu \square_g (\Fullset^{N-1;1} \GeoAng v^i)
	& = (\Rad v^i) \GeoAng^{N-1} \Rad \mytr \upchi
			+ \GeoAngFlatRadComponent (\angdiffuparg{\#} v^i) \cdot (\upmu \angdiff \GeoAng^{N-2} \Rad \mytr \upchi)
				\label{E:SECONDWAVEGEOANGANGISTHEFIRSTCOMMUTATORIMPORTANTTERMS} \\
			& \ \
			- [ia] \upmu^2 (\exp \Densrenormalized) \Speed^2 (g_{ab} \Radunit^b) 
					\GeoAng^{N-1} \Lunit \Lunit \Vortrenormalized
					\notag \\
		& \ \
			+ 
			[ia] \upmu^2 (\exp \Densrenormalized) \Speed^2
			\left(
				\frac{g_{ab} \GeoAng^b}{g_{cd} \GeoAng^c \GeoAng^d}
			\right)
			\GeoAng^N \Lunit \Vortrenormalized
			\notag \\
	& \ \
			+ Harmless_{(Wave)}^{\leq N}.
				 \notag
\end{align}
\end{subequations}

Furthermore, if $1 \leq N \leq 20$
and $\Tanset^{N-1}$ contains a factor of $\Lunit$, 
then
\begin{subequations}
\begin{align} \label{E:WAVEHARMLESSORDERNPURETANGENTIALCOMMUTATORS}
	\upmu \square_g (\Tanset^{N-1} \Lunit v^i)
	& = 
			+ [ia] \upmu (\exp \Densrenormalized) \Speed^2 (g_{ab} \Radunit^b) 
					\Tanset^{N-2} \Tanset^{N-1} \Lunit \Lunit \Vortrenormalized
					\\
	& \ \
			- [ia] \upmu (\exp \Densrenormalized) \Speed^2
			\left(
				\frac{g_{ab} \GeoAng^b}{g_{cd} \GeoAng^c \GeoAng^d}
			\right)
			\Tanset^{N-2} \Tanset^{N-1} \Lunit \GeoAng \Vortrenormalized
				\notag \\
	& \ \
			+ Harmless_{(Wave)}^{\leq N}.
			\notag
\end{align}

Likewise, if $1 \leq N \leq 20$
and $\Fullset^{N-1;1}$ 
contains one or more factors of $\Lunit$, then
\begin{align} \label{E:WAVEHARMLESSORDERNCOMMUTATORS}
	\upmu \square_g (\Fullset^{N-1;1} \Lunit v^i)
	& = 
			- [ia] \upmu^2 (\exp \Densrenormalized) \Speed^2 (g_{ab} \Radunit^b) 
					\Tanset^{N-2} \Lunit \Lunit \Lunit \Vortrenormalized
					\\
	& \ \
			+ [ia] \upmu^2 (\exp \Densrenormalized) \Speed^2
			\left(
				\frac{g_{ab} \GeoAng^b}{g_{cd} \GeoAng^c \GeoAng^d}
			\right)
			\Tanset^{N-2} \Lunit \GeoAng \Lunit \Vortrenormalized
				\notag \\
	& \ \
			+ Harmless_{(Wave)}^{\leq N}.
			\notag
\end{align}
Similarly, if $1 \leq N \leq 20$
and $\Tanset^{N-1}$ 
contains one or more factors of $\Lunit$, then
\begin{align} \label{E:WAVESECONDHARMLESSORDERNCOMMUTATORS}
	\upmu \square_g (\Tanset^{N-1} \Rad v^i)
	& = 
	- [ia] \upmu^2 (\exp \Densrenormalized) \Speed^2 (g_{ab} \Radunit^b) 
					\Tanset^{N-1} \Lunit \Lunit \Vortrenormalized
					\\
	& \ \
		+ [ia] \upmu^2 (\exp \Densrenormalized) \Speed^2
			\left(
				\frac{g_{ab} \GeoAng^b}{g_{cd} \GeoAng^c \GeoAng^d}
			\right)
			\Tanset^{N-1} \GeoAng \Lunit \Vortrenormalized
			\notag	\\
	& \ \
			+ Harmless_{(Wave)}^{\leq N}.
			\notag
\end{align}
In addition, if $1 \leq N \leq 20$
and $\Fullset_{\ast}^{N-1;1}$ 
contains one or more factors of $\Lunit$, then
\begin{align} \label{E:WAVETHIRDHARMLESSORDERNCOMMUTATORS}
	\upmu \square_g (\Fullset_{\ast}^{N-1;1} \GeoAng v^i)
	& = 
			- [ia] \upmu^2 (\exp \Densrenormalized) \Speed^2 (g_{ab} \Radunit^b) 
					\Tanset^{N-2}  \GeoAng \Lunit \Lunit \Vortrenormalized
					\\
	& \ \
			+ 
			[ia] \upmu^2 (\exp \Densrenormalized) \Speed^2
			\left(
				\frac{g_{ab} \GeoAng^b}{g_{cd} \GeoAng^c \GeoAng^d}
			\right)
			\Tanset^{N-2} \GeoAng \GeoAng \Lunit \Vortrenormalized
			\notag	\\
	& \ \
			+ Harmless_{(Wave)}^{\leq N}.
			\notag
\end{align}
Finally, $\Densrenormalized  - v^1$ 
verifies similar estimates according to the following prescription:
\begin{align} \label{E:KEYDIFFICULTFORRENORMALZIEDDENSITY}
	& \mbox{in 
	\eqref{E:WAVELISTHEFIRSTCOMMUTATORIMPORTANTTERMS}-\eqref{E:WAVETHIRDHARMLESSORDERNCOMMUTATORS}
	with $i=1$},
		\\
	& \mbox{we may replace the explicit factors of $v^1$ on the 
	LHSs and RHSs
	with
	$\Densrenormalized  - v^1$
	}
		\notag \\
	& \mbox{as long as we change the sign of all explicit 
		$\Vortrenormalized$-containing products on the RHSs}.
		\notag
\end{align}
\end{subequations}
\end{proposition}

The next proposition, which we prove 
in Sect.~\ref{SS:PROOFOFPROPVORTICITYIDOFKEYDIFFICULTENREGYERRORTERMS}, 
is an analog of Prop.~\ref{P:WAVEIDOFKEYDIFFICULTENREGYERRORTERMS}
for $\Vortrenormalized$. 
Specifically, the proposition identifies the main error terms in the inhomogeneous
transport equations verified by the higher-order versions of
$\Vortrenormalized$.

\begin{proposition}[\textbf{Identification of the key difficult error term factors 
	in the commuted transport equation}]
	\label{P:VORTICITYIDOFKEYDIFFICULTENREGYERRORTERMS}
	Assume that $1 \leq N \leq 20$.
	Under the data-size and bootstrap assumptions 
	of Sects.~\ref{SS:FLUIDVARIABLEDATAASSUMPTIONS}-\ref{SS:PSIBOOTSTRAP}
	and the smallness assumptions of Sect.~\ref{SS:SMALLNESSASSUMPTIONS}, 
	the following pointwise estimates hold
	on $\mathcal{M}_{\Tboot,U_0}$
	(see Sect.~\ref{SS:STRINGSOFCOMMUTATIONVECTORFIELDS} regarding the vectorfield operator notation):
	\begin{subequations}
	\begin{align}
		\upmu \Transport \GeoAng^N \Lunit \Vortrenormalized
		& = 
			(\GeoAng \Vortrenormalized) \GeoAng^{N-1} \Rad \mytr \upchi
			+ Harmless_{(Vort)}^{\leq N+1},
				\label{E:VORTICITYLISTHEFIRSTCOMMUTATORIMPORTANTTERMS} \\
	\upmu \Transport \GeoAng^{N+1} \Vortrenormalized
		& = - g(\GeoAng,\GeoAng) (\Lunit \Vortrenormalized) \GeoAng^{N-1} \Rad \mytr \upchi
				+ \GeoAngFlatRadComponent (\GeoAng \Vortrenormalized) \GeoAng^{N-1} \Rad \mytr \upchi
				+ Harmless_{(Vort)}^{\leq N+1},
				\label{E:VORTICITYGEOANGANGISTHEFIRSTCOMMUTATORIMPORTANTTERMS}
	\end{align}
	\end{subequations}
	where $\GeoAngFlatRadComponent$ is the scalar-valued function from
	\eqref{E:FLATYDERIVATIVERADIALCOMPONENT}.

Furthermore, if $1 \leq N \leq 20$ and $\Tanset^{N+1}$ is any $(N+1)^{st}$ order 
$\mathcal{P}_u$-tangent operator except for 
$\GeoAng^N \Lunit$ or $\GeoAng^{N+1}$,
then
\begin{align} \label{E:VORTICITYHARMLESSORDERNPLUSONECOMMUTATORS}
	\upmu \Transport \Tanset^{N+1} \Vortrenormalized
	& = Harmless_{(Vort)}^{\leq N+1}.
\end{align}

Finally, if $\Singletan \in \Tanset$, then
\begin{align} \label{E:TRIVIALVORTICITYHARMLESSORDERNPLUSONECOMMUTATORS}
	\upmu \Transport \Singletan \Vortrenormalized
	& = Harmless_{(Vort)}^{\leq 1}.
\end{align}

\end{proposition}

\subsection{Technical estimates involving the eikonal function quantities}
\label{L:ESTIMATERELATINGUPMUTOTRCHI}
In this section, we provide two technical lemmas that will allow us to reduce
the analysis of the top-order derivatives of $\upmu$
to those of $\mytr \upchi$.
This is mainly for convenience.

We start with a lemma in which we obtain higher-order analogs of 
Lemma~\ref{L:RADDIRECTIONCHITRANSPORT}.

\begin{lemma}[\textbf{Estimate connecting} $\Fullset^{\vec{I}} \mytr \upchi$ \textbf{to} 
$\angLap \Tanset^{\vec{J}} \upmu$]
	\label{L:TANGENTIALCOMMUTEDRADTRCHIANGLAPUPMUCOMPARISON}
	Assume that $1 \leq N \leq 20$.
	Let $\vec{I} \in \mathcal{I}^{N;1}$,
	where the multi-index set $\mathcal{I}^{N;1}$ 
	is defined in Def.~\ref{D:COMMUTATORMULTIINDICES}.
	Let $\vec{J}$ be any multi-index 
	formed by deleting
	the one entry in $\vec{I}$
	corresponding to the single $\Rad$ differentiation
	and by possibly permuting the remaining entries
	(and thus $|\vec{J}| = N-1$ and the corresponding operator
	is $\Tanset^{\vec{J}}$).
	Under the data-size and bootstrap assumptions 
	of Sects.~\ref{SS:FLUIDVARIABLEDATAASSUMPTIONS}-\ref{SS:PSIBOOTSTRAP}
	and the smallness assumptions of Sect.~\ref{SS:SMALLNESSASSUMPTIONS}, 
	the following estimate holds on
	$\mathcal{M}_{\Tboot,U_0}$
	(see Sect.~\ref{SS:STRINGSOFCOMMUTATIONVECTORFIELDS} regarding the vectorfield operator notation):
	\begin{align} \label{E:TANGENTIALCOMMUTEDRADTRCHIANGLAPUPMUCOMPARISON}
		\left|
			\Fullset^{\vec{I}} \mytr \upchi
			- 
			\angLap \Tanset^{\vec{J}} \upmu
		\right|
		& \lesssim
			\left|
				\Fullset_{\ast}^{[1,N+1];\leq 1} \Psi
			\right|
			+
			\left|
				\myarray
					[\Fullset_{\ast \ast}^{[1,N];0} \BadVar]
					{\Fullset_{\ast}^{[1,N];\leq 1} \GdVar}
			\right|.
	\end{align}
\end{lemma}

\begin{proof}
	See Sect.~\ref{SS:OFTENUSEDESTIMATES} for some comments on the analysis.
	First, using the commutator estimate 
	\eqref{E:LESSPRECISEPERMUTEDVECTORFIELDSACTINGONFUNCTIONSCOMMUTATORESTIMATE} with 
	$f = \mytr \upchi$, $N$ in the role of $N+1$, and $M=1$,
	the estimate \eqref{E:POINTWISEESTIMATESFORCHIANDITSDERIVATIVES},
	and the $L^{\infty}$ estimates of Prop.~\ref{P:IMPROVEMENTOFAUX},
	we write
	$
	\Fullset^{\vec{I}} \mytr \upchi
	=
	\Tanset^{\vec{J}} \Rad \mytr \upchi
	$
	plus error terms with magnitudes
	$\lesssim \mbox{RHS~\eqref{E:TANGENTIALCOMMUTEDRADTRCHIANGLAPUPMUCOMPARISON}}$.
	Next, we apply $\Tanset^{\vec{J}}$
	to \eqref{E:RADDIRECTIONCHITRANSPORT}.
	Using 
	the estimates
	\eqref{E:ANGDIFFXI}-\eqref{E:ANGDIFFXNOSPECIALSTRUCTUREDIFFERENTIATEDPOINTWISE} 
	and
	\eqref{E:POINTWISEESTIMATESFORGSPHEREANDITSSTARDERIVATIVES}
	and the $L^{\infty}$ estimates of Prop.~\ref{P:IMPROVEMENTOFAUX},
	we write
	$
	\Tanset^{\vec{J}} \Rad \mytr \upchi 
	= \Tanset^{\vec{J}} \angLap \upmu 
	$
	plus error terms with magnitudes
	$\lesssim \mbox{RHS~\eqref{E:TANGENTIALCOMMUTEDRADTRCHIANGLAPUPMUCOMPARISON}}$.
	Finally, we use the commutator estimate
	\eqref{E:ANGLAPNOSPECIALSTRUCTUREFUNCTIONCOMMUTATOR} with 
	$f = \upmu$, $N-1$ in the role of $N$, and $M=0$
	and the $L^{\infty}$ estimates of Prop.~\ref{P:IMPROVEMENTOFAUX}
	to write
	$
	\Tanset^{\vec{J}} \angLap \upmu 
	=
	\angLap \Tanset^{\vec{J}} \upmu
	$
	plus error terms with magnitudes
	$\lesssim \mbox{RHS~\eqref{E:TANGENTIALCOMMUTEDRADTRCHIANGLAPUPMUCOMPARISON}}$.
	Combining the above estimates, we conclude \eqref{E:TANGENTIALCOMMUTEDRADTRCHIANGLAPUPMUCOMPARISON}.
\end{proof}

\begin{lemma}[\textbf{Connecting derivatives of $\upmu$ to derivatives of $\mytr \upchi$
up to error terms}] 
\label{L:HIGHGEOANGDERIVATIVESOFUPMUINTERMSOFRADTRCHI}
	Assume that $1 \leq N \leq 20$.
	Under the data-size and bootstrap assumptions 
	of Sects.~\ref{SS:FLUIDVARIABLEDATAASSUMPTIONS}-\ref{SS:PSIBOOTSTRAP}
	and the smallness assumptions of Sect.~\ref{SS:SMALLNESSASSUMPTIONS}, 
	the following estimates hold on
	$\mathcal{M}_{\Tboot,U_0}$
	(see Sect.~\ref{SS:STRINGSOFCOMMUTATIONVECTORFIELDS} regarding the vectorfield operator notation):
	\begin{align} \label{E:HIGHGEOANGDERIVATIVESOFUPMUINTERMSOFRADTRCHI}
	&
	\left|
		\GeoAng^{N+1} \upmu
		-
		g(\GeoAng,\GeoAng) \GeoAng^{N-1} \Rad \mytr \upchi
	\right|,
		\,
	\left|
		\angdiffuparg{\#} \GeoAng^N \upmu
		-
		(\GeoAng^{N-1} \Rad \mytr \upchi) \GeoAng
	\right|
	\\
	& \lesssim 
			\left|
				\Fullset_{\ast}^{[1,N+1];\leq 1} \threePsi
			\right|
		+
		\left|
				\myarray
					[\Fullset_{\ast \ast}^{[1,N];0} \BadVar]
					{\Fullset_{\ast}^{[1,N];\leq 1} \GdVar}
			\right|.
		\notag
\end{align}
\end{lemma}

\begin{proof}
	See Sect.~\ref{SS:OFTENUSEDESTIMATES} for some comments on the analysis.
	We start by proving \eqref{E:HIGHGEOANGDERIVATIVESOFUPMUINTERMSOFRADTRCHI}
	for the first term on the LHS.
	Using \eqref{E:ANGLAPINTERMSOFGEOANGDERIVATIVES}
	with $f = \GeoAng^{N-1} \upmu$,
	we obtain
	\begin{align} \label{E:FIRSTIDHIGHGEOANGDERIVATIVESOFUPMUINTERMSOFRADTRCHI}
	\GeoAng^{N+1} \upmu
	= g(\GeoAng,\GeoAng) \angLap \GeoAng^{N-1} \upmu
		+ 
			\left\lbrace
				\GeoAng 
				\ln g(\GeoAng,\GeoAng) 
			\right\rbrace
			\GeoAng^N \upmu.
	\end{align}
	Using \eqref{E:TANGENTIALCOMMUTEDRADTRCHIANGLAPUPMUCOMPARISON}
	and
	$|\GeoAng| = 1 + \mathcal{O}(\varepsilon)$
	(which is a simple consequence of \eqref{E:GEOANGPOINTWISE} 
	and the $L^{\infty}$ estimates of Prop.~\ref{P:IMPROVEMENTOFAUX}), 
	we deduce that
	$
	g(\GeoAng,\GeoAng) \angLap \GeoAng^{N-1} \upmu
	=
	g(\GeoAng,\GeoAng) \GeoAng^{N-1} \Rad \mytr \upchi
	$
	plus error terms that are bounded in magnitude by 
	$\lesssim$ RHS~\eqref{E:HIGHGEOANGDERIVATIVESOFUPMUINTERMSOFRADTRCHI}.
	Next, we use Lemma~\ref{L:SCHEMATICDEPENDENCEOFMANYTENSORFIELDS}
	and the $L^{\infty}$ estimates of Prop.~\ref{P:IMPROVEMENTOFAUX}
	to deduce that
	$
	|\GeoAng \ln g(\GeoAng,\GeoAng)|
	= |\GeoAng \smoothfunction(\GdVar)|
	\lesssim \varepsilon
	$.
	It follows that the last product on RHS~\eqref{E:FIRSTIDHIGHGEOANGDERIVATIVESOFUPMUINTERMSOFRADTRCHI}
	is bounded in magnitude by 
	$\lesssim |\Fullset_{\ast \ast}^{[1,N];0} \BadVar|$.
	We have thus obtained the desired estimate.

	We now prove \eqref{E:HIGHGEOANGDERIVATIVESOFUPMUINTERMSOFRADTRCHI} for
	$
	\left|
		\angdiffuparg{\#} \GeoAng^N \upmu
		-
		(\GeoAng^{N-1} \Rad \mytr \upchi) \GeoAng
	\right|
	$.
	By \eqref{E:INVERSEANGULARMETRICINTERMSOFGEOANG}, we have
	$
	\displaystyle
	\angdiffuparg{\#} \GeoAng^N \upmu
	= 
	\frac{1}{g(\GeoAng,\GeoAng)}
	(\GeoAng^{N+1} \upmu) \GeoAng
	$.
	The desired estimate is therefore
	a simple consequence of the
	estimates we obtained above for
	$
	\left|
		\GeoAng^{N+1} \upmu
		-
		g(\GeoAng,\GeoAng) \GeoAng^{N-1} \Rad \mytr \upchi
	\right|
	$
	and the estimate $|\GeoAng| = 1 + \mathcal{O}(\varepsilon)$
	noted above.
\end{proof}

\subsection{Pointwise estimates for the deformation tensors of the commutation vectorfields}
\label{SS:IMPORTANTTERMSINDEFTENSORS}
In the next lemma, we identify the main terms in various derivatives
of the frame components of the deformation tensors of the commutation
vectorfields defined in \eqref{E:COMMUTATIONVECTORFIELDS}. 
The main terms are located on the left-hand sides
of the estimates stated in the lemma,
while the right-hand sides of the estimates
contain simple error terms that
will be easy to bound in the energy estimates.
The main terms involve top-order derivatives
of the eikonal function quantities and are difficult
to control in the energy estimates.

\begin{lemma} [\textbf{Identification of the important terms in} 
$\deform{\Lunit}$, 
$\deform{\Rad}$,
\textbf{and} $\deform{\GeoAng}$]
\label{L:IMPORTANTDEFTENSORTERMS}
Assume that $1 \leq N \leq 20$.
Under the data-size and bootstrap assumptions 
of Sects.~\ref{SS:FLUIDVARIABLEDATAASSUMPTIONS}-\ref{SS:PSIBOOTSTRAP}
and the smallness assumptions of Sect.~\ref{SS:SMALLNESSASSUMPTIONS}, 
the following estimates hold on
$\mathcal{M}_{\Tboot,U_0}$
(see Sect.~\ref{SS:STRINGSOFCOMMUTATIONVECTORFIELDS} regarding the vectorfield operator notation).

\medskip

\noindent
\underline{\textbf{Important terms in the derivatives of} $\deform{\Lunit}$.}
For $M=0,1$, we have
	\begin{subequations}
	\begin{align} \label{E:LDEFIMPORTANTRADTRACEANGTERMS}
		\left|
			\Fullset^{N-1;M} \Rad \mytr \angdeform{\Lunit}
			- 2 \angLap \Fullset^{N-1;M} \upmu 
		\right|
		& \lesssim
			\left|
				\Fullset_{\ast}^{[1,N+1];\leq 2} \Psi
			\right|
			+
			\left|
				\myarray
					[\Fullset_{\ast \ast}^{[1,N];\leq 1} \BadVar]
					{\Fullset_{\ast}^{[1,N];\leq 2} \GdVar}
				\right|,
	\end{align}

	\begin{align}  \label{E:LDEFIMPORTANTLIERADSPHERERADANDAGNDIVSPHERETERMS}
		&
		\left| 
			\angLie_{\Fullset}^{N-1;M} \angdiffuparg{\#} \deformarg{\Lunit}{\Lunit}{\Rad}
		\right|,
		\,
		\left|
				\Fullset^{N-1;M} \angdiv \angdeformoneformupsharparg{\Lunit}{\Rad}
				- \angLap \Fullset^{N-1;M} \upmu
		\right|,
		\,
		\left|
				\angLie_{\Fullset}^{N-1;M} \angdiffuparg{\#} \mytr \angdeform{\Lunit}
				-  
				2 \angdiffuparg{\#} \Fullset^{N-1;M} \mytr \upchi
		\right|
			\\
		& \lesssim
			\left|
				\Fullset_{\ast}^{[1,N+1];\leq 2} \Psi
			\right|
			+
			\left|
				\myarray
					[\Fullset_{\ast \ast}^{[1,N];\leq 1} \BadVar]
					{\Fullset_{\ast}^{[1,N];\leq 2} \GdVar}
				\right|.
			\notag
	\end{align}
	\end{subequations}

	\medskip
	\noindent
	\underline{\textbf{Important terms in the derivatives of} $\deform{\Rad}$.}
	We have
	\begin{subequations}
	\begin{align}  \label{E:RADDEFIMPORTANTLIERADSPHERELANDRADTRACESPHERETERMS}
		& \left|
			 \angLie_{\Tanset}^{N-1}	\angLie_{\Rad} \angdeformoneformupsharparg{\Rad}{\Lunit}
				+ \angdiffuparg{\#} \Tanset^{N-1} \Rad \upmu
			\right|,
			\,
		\left|
			\Tanset^{N-1} \Rad \mytr \angdeform{\Rad}
			+ 2 \upmu \Tanset^{N-1} \Rad \mytr \upchi
		\right|
			\\
		& \lesssim
			\left|
				\Fullset_{\ast}^{[1,N+1];\leq 2} \Psi
			\right|
			+
			\left|
				\myarray
					[\Fullset_{\ast \ast}^{[1,N];\leq 1} \BadVar]
					{\Fullset_{\ast}^{[1,N];\leq 2} \GdVar}
				\right|
			+
			1,
			\notag
	\end{align}

	\begin{align} \label{E:RADDEFIMPORTANTANGDIVSPHERELANDANGDIFFPILRADTERMS}
		& 
		\left|
			\Tanset^{N-1} \angdiv \angdeformoneformupsharparg{\Rad}{\Lunit}
				+ \Tanset^{N-1} \Rad \mytr \upchi
		\right|,
			\,
 		 \\
		& \left| 
				\angLie_{\Tanset}^{N-1} \angdiffuparg{\#} \deformarg{\Rad}{\Lunit}{\Rad}
				+ \angdiffuparg{\#} \Tanset^{N-1} \Rad \upmu
			\right|,
			\,
		\left|
			\angLie_{\Tanset}^{N-1} \angdiffuparg{\#} \mytr \angdeform{\Rad}
			+ 2 \upmu \angdiffuparg{\#} \Tanset^{N-1} \mytr \upchi
		\right|
			\notag \\
		& \lesssim
			\left|
				\Fullset_{\ast}^{[1,N+1];\leq 2} \Psi
			\right|
			+
			\left|
				\myarray
					[\Fullset_{\ast \ast}^{[1,N];\leq 1} \BadVar]
					{\Fullset_{\ast}^{[1,N];\leq 2} \GdVar}
				\right|.
			\notag
	\end{align}
	\end{subequations}

	\medskip
	\noindent
	\underline{\textbf{Important terms in the derivatives of} $\deform{\GeoAng}$.}
	For $M=0,1$, we have
	\begin{subequations}
	\begin{align}  \label{E:GEOANGDEFIMPORTANTLIERADSPHERELANDRADTRACESPHERETERMS}
		& \left|
			 \angLie_{\Fullset}^{N-1;M}	\angLie_{\Rad} \angdeformoneformupsharparg{\GeoAng}{\Lunit}
				+ (\angLap \Fullset^{N-1;M} \upmu) \GeoAng 
			\right|,
			\,
		\left|
			\Fullset^{N-1;M} \Rad \mytr \angdeform{\GeoAng}
			- 2 \GeoAngFlatRadComponent \angLap \Fullset^{N-1;M} \upmu
		\right|
			\\
		& \lesssim
			\left|
				\Fullset_{\ast}^{[1,N+1];\leq 2} \Psi
			\right|
			+
			\left|
				\myarray
					[\Fullset_{\ast \ast}^{[1,N];\leq 1} \BadVar]
					{\Fullset_{\ast}^{[1,N];\leq 2} \GdVar}
				\right|,
			\notag
	\end{align}

	\begin{align} \label{E:GEOANGDEFIMPORTANTANGDIVSPHERELANDANGDIFFPILRADTERMS}
		& 
		\left|
			\Fullset^{N-1;M} \angdiv \angdeformoneformupsharparg{\GeoAng}{\Lunit}
				+ \GeoAng \Fullset^{N-1;M} \mytr \upchi
		\right|,
 			\,
		\left|
			\Fullset^{N-1;M} \angdiv \angdeformoneformupsharparg{\GeoAng}{\Rad}
			- \left\lbrace
					\upmu \GeoAng \Fullset^{N-1;M} \mytr \upchi
					+ \GeoAngFlatRadComponent \angLap \Fullset^{N-1;M} \upmu
				\right\rbrace
		\right|,
			\\
		& \left| 
				\angLie_{\Fullset}^{N-1;M} \angdiffuparg{\#} \deformarg{\GeoAng}{\Lunit}{\Rad}
				+ (\angLap \Fullset^{N-1;M} \upmu) \GeoAng
			\right|,
			\,
		\left|
			\angLie_{\Fullset}^{N-1;M} \angdiffuparg{\#} \mytr \angdeform{\GeoAng}
			- 2 \GeoAngFlatRadComponent \angdiffuparg{\#} \Fullset^{N-1;M} \mytr \upchi
		\right|
			\notag \\
		& \lesssim
			\left|
				\Fullset_{\ast}^{[1,N+1];\leq 2} \Psi
			\right|
			+
			\left|
				\myarray
					[\Fullset_{\ast \ast}^{[1,N];\leq 1} \BadVar]
					{\Fullset_{\ast}^{[1,N];\leq 2} \GdVar}
				\right|.
			\notag
	\end{align}
	\end{subequations}

	Moreover, we have
	\begin{align} \label{E:MAINTERMINPUREGEOANGDERIVATIVESOFUPMUANGDEFORMGEOANGLUNITANDGEOANGRAD}
		&
		\left|
			\upmu \angLie_{\GeoAng}^N \angdeformoneformupsharparg{\GeoAng}{\Lunit} 
			+ 
			\upmu (\angLie_{\GeoAng}^N \mytr \upchi) \GeoAng 
		\right|,
			\, 
		\left|
			\angLie_{\GeoAng}^N \angdeformoneformupsharparg{\GeoAng}{\Rad}
			-
			\upmu (\angLie_{\GeoAng}^N \mytr \upchi) \GeoAng 
			-
			\GeoAngFlatRadComponent \angdiffuparg{\#} \GeoAng^N \upmu
		\right|
			\\
		& \lesssim 
			\left|
				\Tanset^{[1,N+1]} \threePsi
			\right|
		+
		\left|
			\myarray
			[\Fullset_{\ast \ast}^{[1,N];0} \BadVar]
			{\Tanset^{[1,N]} \GdVar}
		\right|,
		\notag
	\end{align}

	\begin{align} \label{E:MAINTERMINPUREGEOANGDERIVATIVESOFUPMUANGDEFORMLUNITRAD}
		\left|
			\angLie_{\GeoAng}^N \angdeformoneformupsharparg{\Lunit}{\Rad}
			-
			\angdiffuparg{\#} \GeoAng^N \upmu
		\right|
		& \lesssim 
			\left|
				\Fullset_{\ast}^{[1,N+1];\leq 1} \threePsi
			\right|
		+
		\left|
			\myarray
			[\Fullset_{\ast \ast}^{[1,N];0} \BadVar]
			{\Tanset^{[1,N]} \GdVar}
		\right|.
	\end{align}

	Above, within a given inequality, 
	the symbol 
	$\Fullset^{N-1;M}$
	on the LHS always denotes the same
	order $N-1$ vectorfield operator,
	and similarly for the symbol
	$\Tanset^{N-1}$.

\end{lemma}

\begin{proof}
	See Sect.~\ref{SS:OFTENUSEDESTIMATES} for some comments on the analysis.
	The main point of the proof is to identity the products featuring
	the top-order derivatives of the eikonal function quantities,
	which we place on the LHS of the estimates.
	More precisely, we aim to identify the products containing a factor with
	$N+1$ derivatives on $\upmu$ or $N$ derivatives on $\mytr \upchi$,
	\emph{with none of the derivatives being in the $\Lunit$ direction};
	all other terms are error terms that can be shown to be bounded in magnitude
	by $\lesssim$ the RHSs of the inequalities
	(including top-order derivatives of $\upmu$ or $\mytr \upchi$ involving an $\Lunit$ derivative,
	which we bound with the estimates 
	\eqref{E:LLUNITUPMUFULLSETSTARCOMMUTEDPOINTWISEESTIMATE}
	and
	\eqref{E:LUNITCOMMUTEDLUNITSMALLIPOINTWISE}).

	We start by proving \eqref{E:RADDEFIMPORTANTLIERADSPHERELANDRADTRACESPHERETERMS}
	for the first term 
	$
	\angLie_{\Tanset}^{N-1}	\angLie_{\Rad} \angdeformoneformupsharparg{\Rad}{\Lunit}
	+ 
	\angdiffuparg{\#} \Rad \Tanset^{N-1} \upmu
	$
	on the LHS.
	We apply 
	$\angLie_{\Tanset}^{N-1}	\angLie_{\Rad}$
	to the $\gsphere$-dual of the equation \eqref{E:RADDEFORMSPHERERAD}
	for $\angdeformoneformarg{\Rad}{\Lunit}$.
	By Lemma~\ref{L:SCHEMATICDEPENDENCEOFMANYTENSORFIELDS},
	the $\gsphere$-dual of the terms
	$
	- 2 \upzeta^{(Trans-\threePsi)}
		- 2 \upmu \upzeta^{(Tan-\threePsi)}
	$
	is of the form 
	$
		\smoothfunction(\ginversesphere,\angdiff \vec{x},\Rad \threePsi) \GdVar 
		+
		\smoothfunction(\BadVar,\ginversesphere,\angdiff \vec{x}) \Singletan \threePsi
	$.
	Hence, using
	the estimates 
	\eqref{E:ANGDIFFXI}-\eqref{E:ANGDIFFXNOSPECIALSTRUCTUREDIFFERENTIATEDPOINTWISE}  
	and
	\eqref{E:POINTWISEESTIMATESFORGSPHEREANDITSNOSPECIALSTRUCTUREDERIVATIVES}
	and the $L^{\infty}$ estimates of Prop.~\ref{P:IMPROVEMENTOFAUX},
	we find that the $\angLie_{\Tanset}^{N-1}	\angLie_{\Rad}$ derivative
	of these terms are bounded in magnitude by
	$\lesssim 
	\mbox{RHS~\eqref{E:RADDEFIMPORTANTLIERADSPHERELANDRADTRACESPHERETERMS}}
	$
	as desired.
	It remains for us to consider the terms
	$
	-
	\angLie_{\Tanset}^{N-1}	\angLie_{\Rad}
	\angdiffuparg{\#} \upmu
	=
	-
	\angLie_{\Tanset}^{N-1}	\angLie_{\Rad}
	(\ginversesphere \cdot \angdiff \upmu)
	$
	generated by the first term on the RHS
	of the equation \eqref{E:RADDEFORMSPHERERAD}
	for $\angdeformoneformarg{\Rad}{\Lunit}$.
	We first use
	\eqref{E:POINTWISEESTIMATESFORGSPHEREANDITSNOSPECIALSTRUCTUREDERIVATIVES}
	and the $L^{\infty}$ estimates of Prop.~\ref{P:IMPROVEMENTOFAUX}
	to deduce that all terms in the Leibniz expansion of
	$
	-
	\angLie_{\Tanset}^{N-1}	\angLie_{\Rad}
	(\ginversesphere \cdot \angdiff \upmu)
	$
	are bounded in magnitude by
	$
	\lesssim 
	\mbox{RHS~\eqref{E:RADDEFIMPORTANTLIERADSPHERELANDRADTRACESPHERETERMS}}
	$
	except for the top-order-in-$\upmu$ term
	$
	-
	\ginversesphere
	\cdot
	\angLie_{\Tanset}^{N-1}	\angLie_{\Rad} \angdiff \upmu
	$.
	To handle this term,
	we use the commutator estimate
	\eqref{E:NOSPECIALSTRUCTURETENSORFIELDCOMMUTATORESTIMATE}
	with 
	$\xi = \angdiff \upmu$,
	$N-1$ in the role of $N$,
	and
	$M=1$
	and the $L^{\infty}$ estimates of Prop.~\ref{P:IMPROVEMENTOFAUX}
	to express
	$
	-
	\ginversesphere
	\cdot
	\angLie_{\Tanset}^{N-1}	\angLie_{\Rad} \angdiff \upmu
	=
	- \angdiffuparg{\#} \Rad \Tanset^{N-1} \upmu
	$
	plus error terms that are in magnitude
	$
	\lesssim 
	\mbox{RHS~\eqref{E:RADDEFIMPORTANTLIERADSPHERELANDRADTRACESPHERETERMS}}
	$.
	We then bring the top-order term
	$\angdiffuparg{\#} \Rad \Tanset^{N-1} \upmu$
	over to the left, as is indicated on LHS~\eqref{E:RADDEFIMPORTANTLIERADSPHERELANDRADTRACESPHERETERMS},
	which completes the proof of the desired estimate.

	The proof of \eqref{E:RADDEFIMPORTANTLIERADSPHERELANDRADTRACESPHERETERMS}
	for the second term 
	$
	\Tanset^{N-1} \Rad \mytr \angdeform{\Rad}
			+ 2 \upmu \Tanset^{N-1} \Rad \mytr \upchi
	$
	on the LHS is
	based on the formula \eqref{E:RADDEFORMSPHERE}
	but is otherwise similar.
	We omit the details, 
	noting only that the top-order eikonal function term
	occurs when all derivatives fall on
	the $\mytr \upchi$ factor in the first product
	on RHS~\eqref{E:RADDEFORMSPHERE},
	that we use the estimate \eqref{E:POINTWISEESTIMATESFORCHIANDITSDERIVATIVES}
	to bound the below-top-order derivatives of
	$\mytr \upchi$, 
	and that we use the commutator estimate
	\eqref{E:LESSPRECISEPERMUTEDVECTORFIELDSACTINGONFUNCTIONSCOMMUTATORESTIMATE}
	with $f = \mytr \upchi$ to commute derivatives on
	$\mytr \upchi$  
	(rather than the commutator estimate \eqref{E:NOSPECIALSTRUCTURETENSORFIELDCOMMUTATORESTIMATE}
	used above).

	All three estimates in \eqref{E:RADDEFIMPORTANTANGDIVSPHERELANDANGDIFFPILRADTERMS}
	can be proved using essentially the same ideas
	with the help of the formulas
	\eqref{E:RADDEFORMSCALARS},
	\eqref{E:RADDEFORMSPHERERAD},
	and \eqref{E:RADDEFORMSPHERE}.
	More precisely,
	in the case of the
	first term on LHS~\eqref{E:RADDEFIMPORTANTANGDIVSPHERELANDANGDIFFPILRADTERMS}, 
	we use two new ingredients
	\textbf{i)} 
	the commutator estimate
	\eqref{E:ANGDNOSPECIALSTRUCTURETENSORFIELDCOMMUTATORESTIMATE}
	with $\xi = \angdeformoneformarg{\Rad}{\Lunit}$
	(to commute the operator $\Tanset^{N-1}$
	through the operator
	$\angdiv$
	in the term
	$\Tanset^{N-1} \angdiv \angdeformoneformupsharparg{\Rad}{\Lunit}
	=
	\Tanset^{N-1} 
		\left\lbrace 
			\ginversesphere \cdot \angD \angdeformoneformarg{\Rad}{\Lunit}
		\right\rbrace
	$)
	and \textbf{ii)} 
	we use Lemma~\ref{L:TANGENTIALCOMMUTEDRADTRCHIANGLAPUPMUCOMPARISON}
	to replace, 
	up to error terms bounded in magnitude by 
	$\lesssim \mbox{RHS~\eqref{E:RADDEFIMPORTANTANGDIVSPHERELANDANGDIFFPILRADTERMS}}$,
	the top-order eikonal function quantity term
	$-\angLap \Tanset^{N-1} \upmu$
	(generated by the first term on RHS~\eqref{E:RADDEFORMSPHERERAD})
	with the term
	$- \Tanset^{N-1} \Rad \mytr \upchi$ 
	(which we then bring over to LHS~\eqref{E:RADDEFIMPORTANTANGDIVSPHERELANDANGDIFFPILRADTERMS}).

	The estimates \eqref{E:LDEFIMPORTANTRADTRACEANGTERMS}-\eqref{E:LDEFIMPORTANTLIERADSPHERERADANDAGNDIVSPHERETERMS}
	can be proved
	with the help of 
	\eqref{E:LUNITDEFORMSCALARS}-\eqref{E:LUNITDEFORMSPHERE}
	and are based on the same ideas plus one new ingredient:
	to bound the top derivatives of the quantities $\Lunit \upmu$
	in equation \eqref{E:LUNITDEFORMSCALARS}, we use the estimate
	\eqref{E:LLUNITUPMUFULLSETSTARCOMMUTEDPOINTWISEESTIMATE}
	(and the resulting terms are
	bounded in magnitude by 
	$\lesssim \mbox{RHS~\eqref{E:LDEFIMPORTANTLIERADSPHERERADANDAGNDIVSPHERETERMS}}$)
	as desired; we omit the remaining details.

	The proofs of \eqref{E:GEOANGDEFIMPORTANTLIERADSPHERELANDRADTRACESPHERETERMS}-\eqref{E:GEOANGDEFIMPORTANTANGDIVSPHERELANDANGDIFFPILRADTERMS}
	are based on 
	\eqref{E:GEOANGDEFORMSCALARS}-\eqref{E:GEOANGDEFORMSPHERE}
	and require no new ingredients beyond the ones we used above;
	we therefore omit the details.
	The same remarks apply to the estimates
	\eqref{E:MAINTERMINPUREGEOANGDERIVATIVESOFUPMUANGDEFORMGEOANGLUNITANDGEOANGRAD}-\eqref{E:MAINTERMINPUREGEOANGDERIVATIVESOFUPMUANGDEFORMLUNITRAD}
	(except to obtain \eqref{E:MAINTERMINPUREGEOANGDERIVATIVESOFUPMUANGDEFORMLUNITRAD},
	we rely on \eqref{E:LUNITDEFORMSPHERELUNITANDLUNITDEFORMSPHERERAD}).
	\end{proof}

The next lemma complements Lemma~\ref{L:IMPORTANTDEFTENSORTERMS}
by providing bounds for the derivatives of the deformation tensor frame components 
when an $\Lunit$ differentiation is involved or when the number
of derivatives is below-top-order.
In contrast with Lemma~\ref{L:IMPORTANTDEFTENSORTERMS}, 
no difficult terms appear in the estimates.

\begin{lemma}[\textbf{Pointwise estimates for the negligible derivatives of $\deform{\Lunit}$ and $\deform{\GeoAng}$}]
		\label{L:POINTWISEESTIMATESNEGLIGIBLEDERIVATIVESOFDEFORMATIONTENSORS}
	 Assume that $N=20$ and let
	 	$\Singletan \in \Tanset = \lbrace \Lunit, \GeoAng \rbrace$.
	 	Under the data-size and bootstrap assumptions 
		of Sects.~\ref{SS:FLUIDVARIABLEDATAASSUMPTIONS}-\ref{SS:PSIBOOTSTRAP}
		and the smallness assumptions of Sect.~\ref{SS:SMALLNESSASSUMPTIONS}, 
		the following estimates hold on
		$\mathcal{M}_{\Tboot,U_0}$
		(see Sect.~\ref{SS:STRINGSOFCOMMUTATIONVECTORFIELDS} regarding the vectorfield operator notation).

		First, if $\Fullset^{N;\leq 1}$ contains a factor of $\Lunit$, 
		then
		\begin{align} \label{E:TOPDERIVATIVESOFTANGENTIALDEFORMATIONINVOLVINGONELUNITDERIVATIVE}
		&
		\left|
			\Fullset^{N;\leq 1} \mytr \angdeform{\Singletan}
		\right|,
			\,
		\left|
			\Fullset^{N;\leq 1} \deformarg{\Singletan}{\Lunit}{\Rad}
		\right|,
			\,
		\left|
			\Fullset^{N;\leq 1} \deformarg{\Singletan}{\Rad}{\Radunit}
		\right|,
				\\
		&
		\left|
			\angLie_{\Fullset}^{N;\leq 1} \angdeformoneformupsharparg{\Singletan}{\Lunit}
		\right|,
			\,
		\left|
			\angLie_{\Fullset}^{N;\leq 1} \angdeformoneformupsharparg{\Singletan}{\Rad}
		\right|
			\notag \\
		& \lesssim
			\left|
				\Fullset_{\ast}^{[1,N+1];\leq 2} \Psi
			\right|
			+
			\left|
				\myarray
					[\Fullset_{\ast \ast}^{[1,N];\leq 1}\BadVar]
					{\Fullset_{\ast}^{[1,N];\leq 1} \GdVar}
			\right|.
			\notag
	\end{align}

	In addition, if $\Tanset^N$ contains a factor of $\Lunit$, then
	\begin{subequations}
	\begin{align} \label{E:PARTITOPDERIVATIVESOFRADDEFORMATIONINVOLVINGONELUNITDERIVATIVE}
		\left|
			\Tanset^N \mytr \angdeform{\Rad}
		\right|,
			\,
		\left|
			\angLie_{\Tanset}^N \angdeformoneformupsharparg{\Rad}{\Lunit}
		\right|
		& \lesssim
			\left|
				\Fullset_{\ast}^{[1,N+1];\leq 2} \Psi
			\right|
			+
			\left|
				\myarray
					[\Fullset_{\ast \ast}^{[1,N];\leq 1}\BadVar]
					{\Fullset_{\ast}^{[1,N];\leq 1} \GdVar}
			\right|,
				\\
		\left|
			\Tanset^N \deformarg{\Rad}{\Lunit}{\Rad}
		\right|,
			\,
		\left|
			\Tanset^N  \deformarg{\Rad}{\Rad}{\Radunit}
		\right|
		& \lesssim
			\left|
				\Fullset_{\ast}^{[1,N+1];\leq 2} \Psi
			\right|
			+
			\left|
				\myarray
					[\Fullset_{\ast \ast}^{[1,N];\leq 1} \BadVar]
					{\Fullset_{\ast}^{[1,N];\leq 1} \GdVar}
			\right|.
			\label{E:PARTIITOPDERIVATIVESOFRADDEFORMATIONINVOLVINGONELUNITDERIVATIVE} 
	\end{align}
	\end{subequations}

	Moreover, for $1 \leq N \leq 20$,
	the following below-top-order estimates hold,
	where the operator $\Fullset_{\ast}^{N-1;\leq 1}$
	does not necessarily contain any factor of $\Lunit$:
	\begin{subequations}
	\begin{align} \label{E:BELOWTOPORDERDERIVATIVESOFTANGENTIALDEFORMATION}
		&
		\left|
			\Fullset^{N-1;\leq 1} \mytr \angdeform{\Singletan}
		\right|,
			\,
		\left|
			\Fullset^{N-1;\leq 1} \deformarg{\Singletan}{\Lunit}{\Rad}
		\right|,
			\,
		\left|
			\Fullset^{N-1;\leq 1} \deformarg{\Singletan}{\Rad}{\Radunit}
		\right|,
			\,
		\left|
			\angLie_{\Fullset}^{N-1;\leq 1} \angdeformoneformupsharparg{\Singletan}{\Lunit}
		\right|,
			\,
		\left|
			\angLie_{\Fullset}^{N-1;\leq 1} \angdeformoneformupsharparg{\Singletan}{\Rad}
		\right|
				\\
		& \lesssim
			\left|
				\Fullset_{\ast}^{[1,N];\leq 2} \Psi
			\right|
			+
			\left|
				\myarray
					[\Fullset_{\ast \ast}^{[1,N];\leq 1}\BadVar]
					{\Fullset_{\ast}^{[1,N];\leq 1} \GdVar}
			\right|
			+ 
			1.
			\notag
	\end{align}

	In addition, 
	for $1 \leq N \leq 20$,
	we have the following below-top-order estimates:
	\begin{align}
		&
		\left|
			\Tanset^{N-1} \mytr \angdeform{\Rad}
		\right|,
			\,
		\left|
			\Tanset^{N-1} \Lunit \deformarg{\Rad}{\Lunit}{\Rad}
		\right|,
			\,
		\left|
			\Tanset^{N-1} \Lunit \deformarg{\Rad}{\Rad}{\Radunit}
		\right|,
			\,
		\left|
			\angLie_{\Tanset}^{N-1} \angdeformoneformupsharparg{\Rad}{\Lunit}
		\right|,
			\,
		\left|
			\angLie_{\Tanset}^{N-1} \angdeformoneformupsharparg{\Rad}{\Rad}
		\right|
			\label{E:BELOWTOPORDERDERIVATIVESOFRADDEFORMATION} \\
		& \lesssim
			\left|
				\Fullset_{\ast}^{[1,N];\leq 2} \Psi
			\right|
			+
			\left|
				\myarray
					[\Fullset_{\ast \ast}^{[1,N];\leq 1}\BadVar]
					{\Fullset_{\ast}^{[1,N];\leq 1} \GdVar}
			\right|
			+ 
			1.
			\notag
	\end{align}
	\end{subequations}
\end{lemma}
\begin{proof}
		See Sect.~\ref{SS:OFTENUSEDESTIMATES} for some comments on the analysis.
	We first prove \eqref{E:TOPDERIVATIVESOFTANGENTIALDEFORMATIONINVOLVINGONELUNITDERIVATIVE}.
	From Prop.~\ref{L:DEFORMATIONTENSORFRAMECOMPONENTS},
	equation \eqref{E:UPMUFIRSTTRANSPORT},
	and Lemma~\ref{L:SCHEMATICDEPENDENCEOFMANYTENSORFIELDS},
	we see that the deformation tensor components 
	$\mytr \angdeform{\Singletan}$,
	$\deformarg{\Singletan}{\Lunit}{\Rad}$,
	$\cdots$,
	$\angdeformoneformupsharparg{\Singletan}{\Rad}$
	on LHS~\eqref{E:TOPDERIVATIVESOFTANGENTIALDEFORMATIONINVOLVINGONELUNITDERIVATIVE}
	are schematically of the form
	\[
	\smoothfunction(\BadVar,\ginversesphere,\angdiff \vec{x}) \Singletan \threePsi
	+ \smoothfunction(\GdVar,\ginversesphere,\angdiff \vec{x}) 
		\Rad \threePsi 
	+ \smoothfunction(\GdVar,\ginversesphere,\angdiff \vec{x}) 
		\mytr \upchi
	+ \smoothfunction(\GdVar,\ginversesphere,\angdiff \vec{x}) 
		\angdiff \upmu.
	\]
	We now apply
	$\angLie_{\Fullset_{\ast}}^{\leq N;\leq 1}$.
	Recall that $\Fullset_{\ast}^{N;\leq 1}$
	contains a factor of $\Lunit$ by assumption. Let
	$\Fullset_{\ast}^{N-1;\leq 1}$ denote the factors obtained by
	removing the factor of $\Lunit$.
	If all derivatives fall on 
	$\mytr \upchi$, then we 
	the commutator estimate
	\eqref{E:DETAILEDPERMUTEDVECTORFIELDSACTINGONFUNCTIONSCOMMUTATORESTIMATE}
	with $f = \mytr \upchi$,
	the estimate \eqref{E:POINTWISEESTIMATESFORCHIANDITSDERIVATIVES},
	and the $L^{\infty}$ estimates of Prop.~\ref{P:IMPROVEMENTOFAUX},
	to commute the factor of $\Lunit$ in $\Fullset_{\ast}^{N;\leq 1}$
	so that it hits $\mytr \upchi$ first,
	which implies that
	$
	\left|
		\Fullset_{\ast}^{N;\leq 1} \mytr \upchi
	\right|
	\lesssim
	\left|
		\Fullset_{\ast}^{N-1;\leq 1} \Lunit \mytr \upchi
	\right|
	$
	plus error terms that are bounded by 
	$\lesssim \mbox{RHS~\eqref{E:TOPDERIVATIVESOFTANGENTIALDEFORMATIONINVOLVINGONELUNITDERIVATIVE}}$.
	Moreover, using the
	pointwise estimate \eqref{E:LUNITCOMMUTEDLUNITSMALLIPOINTWISE},
	we obtain that
	$
	\left|
		\Fullset_{\ast}^{N-1;\leq 1} \Lunit \mytr \upchi
	\right|
	\lesssim \mbox{RHS~\eqref{E:TOPDERIVATIVESOFTANGENTIALDEFORMATIONINVOLVINGONELUNITDERIVATIVE}}
	$.
	Moreover,
	we bound the remaining factors multiplying 
	$\Fullset^{N;\leq 1} \mytr \upchi$ 
	in magnitude by $\lesssim 1$ via 
	Lemmas~\ref{L:POINTWISEFORRECTANGULARCOMPONENTSOFVECTORFIELDS}
	and \ref{L:POINTWISEESTIMATESFORGSPHEREANDITSDERIVATIVES}
	and the $L^{\infty}$ estimates of Prop.~\ref{P:IMPROVEMENTOFAUX}.
	In total, we find that the product
	under consideration is
	$
	\lesssim \mbox{RHS~\eqref{E:TOPDERIVATIVESOFTANGENTIALDEFORMATIONINVOLVINGONELUNITDERIVATIVE}}
	$
	as desired.
	Similarly, if all derivatives fall on
	$\angdiff \upmu$, 
	we use the commutator estimate
	\eqref{E:DETAILEDPERMUTEDVECTORFIELDSACTINGONFUNCTIONSCOMMUTATORESTIMATE}
	with $f = \upmu$
	and the $L^{\infty}$ estimates 
	of Prop.~\ref{P:IMPROVEMENTOFAUX},
	we can commute the factor of $\Lunit$
	so that it hits $\upmu$ first,
	which implies that
	$
	\left|
		\angdiff \Fullset_{\ast}^{N;\leq 1} \upmu
	\right|
	\lesssim
	\left|
		\Fullset_{\ast}^{N+1;\leq 1} \upmu
	\right|
	\lesssim
	\left|
		\Fullset_{\ast}^{N;\leq 1} \Lunit \upmu
	\right|
	$
	plus error terms that are bounded by 
	$
	\lesssim \mbox{RHS~\eqref{E:TOPDERIVATIVESOFTANGENTIALDEFORMATIONINVOLVINGONELUNITDERIVATIVE}}
	$.
	Moreover,
	we bound the remaining factors multiplying 
	$\angdiff \Fullset^{N;\leq 1} \upmu$ 
	in magnitude by $\lesssim 1$ via 
	Lemmas~\ref{L:POINTWISEFORRECTANGULARCOMPONENTSOFVECTORFIELDS}
	and \ref{L:POINTWISEESTIMATESFORGSPHEREANDITSDERIVATIVES}
	and the $L^{\infty}$ estimates of Prop.~\ref{P:IMPROVEMENTOFAUX}.
	In total, we find that the product
	under consideration is
	$
	\lesssim \mbox{RHS~\eqref{E:TOPDERIVATIVESOFTANGENTIALDEFORMATIONINVOLVINGONELUNITDERIVATIVE}}
	$
	as desired.
	If most (but not all) derivatives fall on
	$\mytr \upchi$ or $\angdiff \upmu$, then we bound all terms
	using the above arguments.
	If most derivatives fall on 
	$\Singletan \threePsi$, $\Rad \threePsi$,  
	$\BadVar$,
	$\GdVar$,
	$\ginversesphere$,
	or
	$\angdiff \vec{x}$,
	then we bound these factors 
	by
	$
	\lesssim \mbox{RHS~\eqref{E:TOPDERIVATIVESOFTANGENTIALDEFORMATIONINVOLVINGONELUNITDERIVATIVE}}
	$
	with the help of
	Lemmas~\ref{L:POINTWISEFORRECTANGULARCOMPONENTSOFVECTORFIELDS}
	and \ref{L:POINTWISEESTIMATESFORGSPHEREANDITSDERIVATIVES}.
	Moreover, we bound the remaining factors 
	(which multiply the factor with many derivatives on it)
	in magnitude by $\lesssim 1$ via 
	Lemmas~\ref{L:POINTWISEFORRECTANGULARCOMPONENTSOFVECTORFIELDS}
	and \ref{L:POINTWISEESTIMATESFORGSPHEREANDITSDERIVATIVES}
	and the $L^{\infty}$ estimates of Prop.~\ref{P:IMPROVEMENTOFAUX}.
	We have therefore proved \eqref{E:TOPDERIVATIVESOFTANGENTIALDEFORMATIONINVOLVINGONELUNITDERIVATIVE}.

	To prove
	\eqref{E:PARTITOPDERIVATIVESOFRADDEFORMATIONINVOLVINGONELUNITDERIVATIVE}-\eqref{E:PARTIITOPDERIVATIVESOFRADDEFORMATIONINVOLVINGONELUNITDERIVATIVE},
 we first use
	Prop.~\ref{L:DEFORMATIONTENSORFRAMECOMPONENTS},
	equation \eqref{E:UPMUFIRSTTRANSPORT},
	and Lemma~\ref{L:SCHEMATICDEPENDENCEOFMANYTENSORFIELDS}
	to deduce that the deformation tensor components 
	$\mytr \angdeform{\Rad}$,
	$\deformarg{\Rad}{\Lunit}{\Rad}$,
	$\deformarg{\Rad}{\Rad}{\Radunit}$,
	$\angdeformoneformupsharparg{\Rad}{\Lunit}$,
	and 
	$\angdeformoneformupsharparg{\Rad}{\Rad}$
	are schematically of the form
	\[
	\smoothfunction(\GdVar,\ginversesphere,\angdiff \vec{x},\Rad \threePsi) 
	+
	\smoothfunction(\BadVar,\ginversesphere,\angdiff \vec{x}) 
	\Singletan \threePsi
	+
	\smoothfunction(\BadVar,\ginversesphere,\angdiff \vec{x}) 
	\mytr \upchi
	+ \Rad \upmu
	+ \ginversesphere \angdiff \upmu.
	\]
	From this schematic formula 
	and the assumption that $\Tanset^N$ contains a factor of $\Lunit$,
	we see that all terms can be bounded by using the arguments given 
	in the previous paragraph.

	The estimates
	\eqref{E:BELOWTOPORDERDERIVATIVESOFTANGENTIALDEFORMATION}
	and 
	\eqref{E:BELOWTOPORDERDERIVATIVESOFRADDEFORMATION}
	can be proved using similar but simpler arguments and
	we therefore omit the details. This completes our proof of the lemma.

	\end{proof}

\subsection{Pointwise estimates involving the fully modified quantities}
\label{SS:POINTWISEFULLYMODIFIED}
Our main goal in this section is to prove
Prop.~\ref{P:KEYPOINTWISEESTIMATE},
in which we obtain pointwise estimates for the most difficult product that appears in our energy estimates:
$(\Rad v^1) \Fullset_{\ast}^{N;\leq 1} \mytr \upchi$.
As a preliminary step, we prove
Lemma~\ref{L:RENORMALIZEDTOPORDERTRCHIJUNKTRANSPORTINVERTED},
in which we use the transport equation \eqref{E:TOPORDERTRCHIJUNKRENORMALIZEDTRANSPORT}
to derive pointwise estimates for the 
fully modified quantities $\upchifullmodarg{\Fullset_{\ast}^{N;\leq 1}}$
From Def.~\ref{D:TRANSPORTRENORMALIZEDTRCHIJUNK}.

Before proving Lemma~\ref{L:RENORMALIZEDTOPORDERTRCHIJUNKTRANSPORTINVERTED},
we first provide a lemma in which we derive pointwise estimates
for the source term $\upchifullmodinhom$
appearing in the transport equation \eqref{E:TOPORDERTRCHIJUNKRENORMALIZEDTRANSPORT}
satisfied by $\upchifullmodarg{\Fullset_{\ast}^{N;\leq 1}}$.
At the same time, for later use,
we derive pointwise estimates for the terms
$
\upchipartialmodinhom
$
and
$
	\upchipartialmodsourcearg{\Fullset_{\ast}^{N-1;\leq 1}}
$
appearing in the transport equation
\eqref{E:COMMUTEDTRCHIJUNKFIRSTPARTIALRENORMALIZEDTRANSPORTEQUATION}
verified by the partially modified quantity
$\upchipartialmodarg{\Fullset_{\ast}^{N-1;\leq 1}}$.

\begin{lemma}[\textbf{Pointwise estimates for} $\Fullset_{\ast}^{N;\leq 1} \upchifullmodinhom$,  
$\Fullset_{\ast}^{N-1;\leq 1} \upchipartialmodinhom$,
\textbf{and} 
$\upchipartialmodsourcearg{\Fullset_{\ast}^{N-1;\leq 1}}$]
	\label{L:CHIPARTIALMODSOURCETERMPOINTWISE}
	Assume that $1 \leq N \leq 20$. Let 
	$\upchifullmodinhom$ be the quantity defined in 
	\eqref{E:LOWESTORDERTRANSPORTRENORMALIZEDTRCHIJUNKDISCREPANCY},
	let $\upchipartialmodinhom$
	be the quantity defined in \eqref{E:LOWESTORDERTRANSPORTPARTIALRENORMALIZEDTRCHIJUNKDISCREPANCY},
	let $\upchipartialmodinhomarg{\Fullset_{\ast}^{N-1;\leq 1}}$
	be the quantity defined in
	\eqref{E:TRANSPORTPARTIALRENORMALIZEDTRCHIJUNKDISCREPANCY},
	and let
	$
	\upchipartialmodsourcearg{\Fullset_{\ast}^{N-1;\leq 1}}
	$
	be the quantity defined in
	\eqref{E:TRCHIJUNKCOMMUTEDTRANSPORTEQNPARTIALRENORMALIZATIONINHOMOGENEOUSTERM}.
	Under the data-size and bootstrap assumptions 
	of Sects.~\ref{SS:FLUIDVARIABLEDATAASSUMPTIONS}-\ref{SS:PSIBOOTSTRAP}
	and the smallness assumptions of Sect.~\ref{SS:SMALLNESSASSUMPTIONS}, 
	the following pointwise estimates hold
	on $\mathcal{M}_{\Tboot,U_0}$: 
\begin{subequations}
	\begin{align} 
		\left|
			\Fullset_{\ast}^{N;\leq 1} \upchifullmodinhom
			+ 
			\vec{G}_{\Lunit \Lunit}\contr
			\Rad \Fullset_{\ast}^{N;\leq 1} \threePsi
		\right|
		& 
		\lesssim
		\upmu
		\left|
			\Fullset_{\ast}^{[1,N+1];\leq 1} \threePsi
		\right|
		+
		\left|
			\Fullset_{\ast}^{[1,N];\leq 2} \threePsi
		\right|
		+
		\left|
			\Fullset_{\ast}^{[1,N];\leq 1} \GdVar
		\right|
		+
		\left|
			\Fullset_{\ast \ast}^{[1,N];\leq 1} \BadVar
		\right|,
			\label{E:TOPDERIVATIVESOFXPOINTWISEBOUND} \\
		\left|
			\Fullset_{\ast}^{N;\leq 1} \upchifullmodinhom
		\right|
		&
		\lesssim
		\left|
			\Fullset_{\ast}^{[1,N+1];\leq 2} \threePsi
		\right|
		+
		\left|
			\Fullset_{\ast}^{[1,N];\leq 1} \GdVar
		\right|
		+
		\left|
			\Fullset_{\ast \ast}^{[1,N];\leq 1} \BadVar
		\right|,
			\label{E:CRUDEXPOINTWISEBOUND} \\
		\left|
			\Fullset_{\ast}^{N;\leq 1} \upchipartialmodinhom
		\right|
		& \lesssim
			\left|
				\Fullset_{\ast}^{[1,N+1];\leq 1} \threePsi
			\right|
			+
			\left|
				\Fullset_{\ast}^{[1,N];\leq 1} \GdVar
			\right|
			+
		\left|
			\Fullset_{\ast \ast}^{[1,N];0} \BadVar
		\right|,
				\label{E:POINTWISELOWESTORDERTRANSPORTPARTIALRENORMALIZEDTRCHIJUNKDISCREPANCY} \\
		\left|
			\Lunit \upchipartialmodinhomarg{\Fullset_{\ast}^{N-1;\leq 1}}
		\right|,
			\,
		\left|
			\GeoAng \upchipartialmodinhomarg{\Fullset_{\ast}^{N-1;\leq 1}}
		\right|
		& 
		\lesssim
			\left| 
				\Fullset_{\ast}^{[1,N+1];\leq 1} \threePsi 
			\right|,
			\label{E:HARMLESSNATUREOFPARTIALLYMODIFIEDDISCREPANCY}
			\\
		\left|
			\upchipartialmodsourcearg{\Fullset_{\ast}^{N-1;\leq 1}}
		\right|
		& \lesssim 
			\left|
				\Fullset_{\ast}^{[1,N];\leq 1} \GdVar
			\right|
			+
		\left|
			\Fullset_{\ast \ast}^{[1,N];0} \BadVar
		\right|.
			\label{E:CHIPARTIALMODSOURCETERMPOINTWISE}
	\end{align}
	\end{subequations}
	\end{lemma}

	\begin{proof}
	See Sect.~\ref{SS:OFTENUSEDESTIMATES} for some comments on the analysis.
	Throughout this proof, we silently use the $L^{\infty}$ estimates 
	of Prop.~\ref{P:IMPROVEMENTOFAUX}.

To prove \eqref{E:TOPDERIVATIVESOFXPOINTWISEBOUND},
we first use \eqref{E:LOWESTORDERTRANSPORTRENORMALIZEDTRCHIJUNKDISCREPANCY}
and Lemma~\ref{L:SCHEMATICDEPENDENCEOFMANYTENSORFIELDS}
to deduce
$
\upchifullmodinhom 
	= 
	-
	 \vec{G}_{\Lunit \Lunit}\contr\Rad \threePsi
	+
	\upmu \smoothfunction(\GdVar,\ginversesphere,\angdiff \vec{x}) \Singletan \threePsi
$.
We now apply $\Fullset_{\ast}^{N;\leq 1}$ to this identity 
and bring the top-order term 
$\vec{G}_{\Lunit \Lunit}\contr\Rad \Fullset_{\ast}^{N;\leq 1} \threePsi$
over to the left 
(as indicated on LHS~\eqref{E:TOPDERIVATIVESOFXPOINTWISEBOUND}), 
which leaves the commutator terms
$[\vec{G}_{\Lunit \Lunit},\Fullset_{\ast}^{N;\leq 1}]\contr\Rad \threePsi$
and
$\vec{G}_{\Lunit \Lunit}\contr[\Rad,\Fullset_{\ast}^{N;\leq 1}] \threePsi$
on the RHS.
To bound the term
$
\left|
\Fullset_{\ast}^{N;\leq 1}
\left\lbrace 
	\upmu \smoothfunction(\GdVar,\ginversesphere,\angdiff \vec{x}) \Singletan \threePsi
\right\rbrace
\right|
$
by $\leq$ RHS~\eqref{E:TOPDERIVATIVESOFXPOINTWISEBOUND}, 
we use
Lemmas~\ref{L:POINTWISEFORRECTANGULARCOMPONENTSOFVECTORFIELDS}
and \ref{L:POINTWISEESTIMATESFORGSPHEREANDITSDERIVATIVES}.
Note that we have paid special attention to terms in which
all derivatives $\Fullset_{\ast}^{N;\leq 1}$ fall on $\Singletan \threePsi$;
these terms are bounded by 
the first term on RHS~\eqref{E:TOPDERIVATIVESOFXPOINTWISEBOUND}.
To bound 
$
\left|
	[\vec{G}_{\Lunit \Lunit},\Fullset_{\ast}^{N;\leq 1}]\contr\Rad \threePsi
\right|
$
by $\leq$ RHS~\eqref{E:TOPDERIVATIVESOFXPOINTWISEBOUND},
we use the fact that 
$\vec{G}_{\Lunit \Lunit} = \smoothfunction(\GdVar)$ (see Lemma~\ref{L:SCHEMATICDEPENDENCEOFMANYTENSORFIELDS}).
To bound
$
\left|
	\vec{G}_{\Lunit \Lunit}\contr[\Rad,\Fullset_{\ast}^{N;\leq 1}] \threePsi
\right|
$
by $\leq$ RHS~\eqref{E:TOPDERIVATIVESOFXPOINTWISEBOUND}, 
we also use the commutator estimate 
\eqref{E:LESSPRECISEPERMUTEDVECTORFIELDSACTINGONFUNCTIONSCOMMUTATORESTIMATE}
with $f = \threePsi$ and $M \leq 1$.
Combining the above estimates, we arrive at \eqref{E:TOPDERIVATIVESOFXPOINTWISEBOUND}.
The proof of \eqref{E:CRUDEXPOINTWISEBOUND} is similar but simpler and we omit the details.
The same is true for the proof of
\eqref{E:POINTWISELOWESTORDERTRANSPORTPARTIALRENORMALIZEDTRCHIJUNKDISCREPANCY}
since by Lemma~\ref{L:SCHEMATICDEPENDENCEOFMANYTENSORFIELDS},
we have
$
\upchipartialmodinhom 
= 
\smoothfunction(\GdVar,\ginversesphere,\angdiff \vec{x})
\Singletan \threePsi
$.

	We now prove \eqref{E:HARMLESSNATUREOFPARTIALLYMODIFIEDDISCREPANCY}.
	We give the proof only for the second term on the LHS
	since the proof for the first one is similar.
	To proceed, we use \eqref{E:TRANSPORTPARTIALRENORMALIZEDTRCHIJUNKDISCREPANCY} 
	and Lemma~\ref{L:SCHEMATICDEPENDENCEOFMANYTENSORFIELDS}
	to deduce that
	$
	\GeoAng \upchipartialmodinhomarg{\Fullset_{\ast}^{N-1;\leq 1}} 
	= 
	\GeoAng 
	\left\lbrace 
		\smoothfunction(\GdVar,\ginversesphere,\angdiff \vec{x}) 
		\Fullset_{\ast}^{N;\leq 1} \threePsi 
	\right\rbrace$.
	The estimate \eqref{E:HARMLESSNATUREOFPARTIALLYMODIFIEDDISCREPANCY} now follows easily 
	from the previous expression and
	Lemmas~\ref{L:POINTWISEFORRECTANGULARCOMPONENTSOFVECTORFIELDS}
  and \ref{L:POINTWISEESTIMATESFORGSPHEREANDITSDERIVATIVES}.

	We now prove \eqref{E:CHIPARTIALMODSOURCETERMPOINTWISE}.
	We bound the term
	$
	\displaystyle
	\Fullset_{\ast}^{N-1;\leq 1} (\mytr \upchi)^2
	$
	from RHS~\eqref{E:TRCHIJUNKCOMMUTEDTRANSPORTEQNPARTIALRENORMALIZATIONINHOMOGENEOUSTERM}
	by $\leq$ RHS~\eqref{E:CHIPARTIALMODSOURCETERMPOINTWISE}
	with the help of inequality \eqref{E:POINTWISEESTIMATESFORCHIANDITSDERIVATIVES}.
	We bound the term
	$
	\displaystyle
	[\Fullset_{\ast}^{N-1;\leq 1},\vec{G}_{\Lunit \Lunit}]\contr\angLap \threePsi
	$
	with the help of the aforementioned relation $\vec{G}_{\Lunit \Lunit} = \smoothfunction(\GdVar)$ 
	and Cor.~\ref{C:BOUNDSFORDERIVATIVESOFLAPLACIAN}.
	To bound
	$
	\displaystyle
	\vec{G}_{\Lunit \Lunit}\contr[\Fullset_{\ast}^{N-1;\leq 1},\angLap] \threePsi
	$,
	we also use
	the commutator estimate \eqref{E:ANGLAPNOSPECIALSTRUCTUREFUNCTIONCOMMUTATOR}
	with $f = \threePsi$.
	We bound the term
	$
	\displaystyle
	[\Lunit, \Fullset_{\ast}^{N-1;\leq 1}] \mytr \upchi
	$
	with the help of the commutator estimate
	\eqref{E:DETAILEDPERMUTEDVECTORFIELDSACTINGONFUNCTIONSCOMMUTATORESTIMATE}
	with $\mytr \upchi$ in the role of $f$
	and inequality \eqref {E:POINTWISEESTIMATESFORCHIANDITSDERIVATIVES}.
	We bound 
	$
	\displaystyle
	[\Lunit, \Fullset_{\ast}^{N-1;\leq 1}] \upchipartialmodinhom
	$
	with the help of the commutator estimate 
	\eqref{E:DETAILEDPERMUTEDVECTORFIELDSACTINGONFUNCTIONSCOMMUTATORESTIMATE} 
	with 
	$f =\upchipartialmodinhom$
	and \eqref{E:POINTWISELOWESTORDERTRANSPORTPARTIALRENORMALIZEDTRCHIJUNKDISCREPANCY}.
	To bound
	$
	\displaystyle
	\Lunit
				\left\lbrace
					\upchipartialmodinhomarg{\Fullset_{\ast}^{N-1;\leq 1}}
					- 
					\Fullset_{\ast}^{N-1;\leq 1} \upchipartialmodinhom
				\right\rbrace
	$,
	we first note that
	\eqref{E:TRANSPORTPARTIALRENORMALIZEDTRCHIJUNKDISCREPANCY},
	\eqref{E:LOWESTORDERTRANSPORTPARTIALRENORMALIZEDTRCHIJUNKDISCREPANCY},
	and the Leibniz rule imply that the magnitude of this term is
	\[
	\lesssim
	\mathop{\sum_{N_1 + N_2 \leq N}}_{N_1 \geq 1}
	\sum_{M_1 + M_2 \leq 1}
	\left|
		\angLie_{\Fullset}^{N_1;M_1} \vec{G}_{(Frame)}^{\#}
	\right|
	\left|
		\Fullset^{N_2+1;M_2} \threePsi
	\right|
	+
	\left|
		\vec{G}_{(Frame)}^{\#}
	\right|
	\left|
		[\Lunit,\Lunit \Fullset_{\ast}^{N-1;\leq 1}] \threePsi
	\right|.
	\]
	Since Lemma~\ref{L:SCHEMATICDEPENDENCEOFMANYTENSORFIELDS} implies that 
	$\vec{G}_{(Frame)}^{\#} = \smoothfunction(\GdVar,\ginversesphere,\angdiff \vec{x})$,
	the desired bound for the sum follows from
	Lemmas~\ref{L:POINTWISEFORRECTANGULARCOMPONENTSOFVECTORFIELDS}
  and \ref{L:POINTWISEESTIMATESFORGSPHEREANDITSDERIVATIVES}.
	To bound the term 
	$\left|
		[\Lunit,\Lunit \Fullset_{\ast}^{N-1;\leq 1}] \threePsi
	\right|
	$,
	we also use the commutator estimate
	\eqref{E:LESSPRECISEPERMUTEDVECTORFIELDSACTINGONFUNCTIONSCOMMUTATORESTIMATE}
	with $f = \threePsi$ and $M \leq 1$.
	We have therefore proved \eqref{E:CHIPARTIALMODSOURCETERMPOINTWISE},
	which completes the proof of the lemma.

\end{proof}

With the help of the previous lemma, we now derive pointwise
estimates for the fully modified quantities
$\upchifullmodarg{\Fullset_{\ast}^{N;\leq 1}}$.

\begin{lemma}[\textbf{Estimates for solutions to the transport equation verified by} $\upchifullmodarg{\Fullset_{\ast}^{N;\leq 1}}$]
	\label{L:RENORMALIZEDTOPORDERTRCHIJUNKTRANSPORTINVERTED}
	Assume that $N=20$ and 
	let $\upchifullmodarg{\Fullset_{\ast}^{N;\leq 1}}$
	and
	$\upchifullmodinhom$
	be as in
	Prop.~\ref{P:TOPORDERTRCHIJUNKRENORMALIZEDTRANSPORT}.
	Assume first that $\Fullset_{\ast}^{N;\leq 1} = \GeoAng^N$.
	Under the data-size and bootstrap assumptions 
	of Sects.~\ref{SS:FLUIDVARIABLEDATAASSUMPTIONS}-\ref{SS:PSIBOOTSTRAP}
	and the smallness assumptions of Sect.~\ref{SS:SMALLNESSASSUMPTIONS}, 
	the following pointwise estimate holds on
	$\mathcal{M}_{\Tboot,U_0}$:
	\begin{align} \label{E:RENORMALIZEDTOPORDERTRCHIJUNKTRANSPORTINVERTED}
		\left|\upchifullmodarg{\GeoAng^N} \right|(t,u,\vartheta)
		& \leq 
			C
			\left| \upchifullmodarg{\GeoAng^N } \right|(0,u,\vartheta)
				\\
		& \ \ + 	2 (1 + C \varepsilon)
							\int_{s=0}^t 
								\frac{[\Lunit \upmu(s,u,\vartheta)]_-}{\upmu(s,u,\vartheta)}
								\left| \GeoAng^N \upchifullmodinhom \right|(s,u,\vartheta)
							\, ds
					\notag \\
	& \ \ 
				+ 	C
						\int_{s=0}^t 
							\left\lbrace
								\left|
									\Fullset_{\ast}^{[1,N+1];\leq 2} \threePsi
								\right|
								+
								\left|
									\myarray
										[\Fullset_{\ast \ast}^{[1,N];\leq 1} \BadVar]
										{\Fullset_{\ast}^{[1,N];\leq 1} \GdVar}
								\right|
							\right\rbrace
							(s,u,\vartheta)
						\, ds
						\notag 
						\\
			& \ \ 
				+ 	C
						\int_{s=0}^t 
							\left\lbrace
								\upmu
								\left|
									\Fullset^{N+1;\leq 1} \Vortrenormalized
								\right|
							\right\rbrace
							(s,u,\vartheta)
						\, ds
				 + 	C
						\int_{s=0}^t 
							\left|
								\Fullset^{\leq N;\leq 1} \Vortrenormalized
							\right|
							(s,u,\vartheta)
						\, ds.
						\notag
\end{align}

Assume now that
$\Fullset_{\ast}^{N;\leq 1} = \GeoAng^{N-1} \Rad$.
Then
$
\left|\upchifullmodarg{\GeoAng^{N-1} \Rad} \right|(t,u,\vartheta)
$
verifies inequality \eqref{E:RENORMALIZEDTOPORDERTRCHIJUNKTRANSPORTINVERTED},
but with the term
$\left| 
	\upchifullmodarg{\GeoAng^N} 
\right|(0,u,\vartheta)$ on the RHS
replaced by
$\left| 
	\upchifullmodarg{\GeoAng^{N-1} \Rad} 
\right|(0,u,\vartheta) 
+
\left| 
	\upchifullmodarg{\GeoAng^N} 
\right|(0,u,\vartheta)
$,
with the term
$\left| \GeoAng^N \upchifullmodinhom \right|(s,u,\vartheta)$
replaced by
$\left| \GeoAng^{N-1} \Rad \upchifullmodinhom \right|(s,u,\vartheta)$,
and with the following additional double time integral present on the RHS:
\begin{align} \label{E:ADDITIONALTIMEINTEGRALFORFULLYMODIFIEDQUANTITYPOINTWISEESTIMATEONERADIALCASE}
	C 
	\int_{s=0}^t
		\int_{s'=0}^s
			\frac{1}{\upmu_{\star}(s',u)}
			\left\lbrace
				\left|
					\Fullset_{\ast}^{[1,N+1];\leq 2} \Psi
				\right|
			+
			\left|
				\myarray
					[\Fullset_{\ast \ast}^{[1,N];\leq 1} \BadVar]
					{\Fullset_{\ast}^{[1,N];\leq 1} \GdVar}
			\right|
			\right\rbrace
			(s',u,\vartheta)
		\, ds'
	\, ds.
\end{align}
\end{lemma}

\begin{proof}
	See Sect.~\ref{SS:OFTENUSEDESTIMATES} for some comments on the analysis.
	We first prove \eqref{E:RENORMALIZEDTOPORDERTRCHIJUNKTRANSPORTINVERTED}.
	We set $\Fullset_{\ast}^{N;\leq 1} = \GeoAng^N$ in \eqref{E:TOPORDERTRCHIJUNKRENORMALIZEDTRANSPORT}
	and view both sides of the equation 
	as functions of $(s,u,\vartheta)$.
	Noting that
	$
	\displaystyle
	\Lunit = \frac{\partial}{\partial s}
	$
	in the present context,
	we define the integrating factor 
	\begin{align} \label{E:INTFACTFULLYMODIFIEDTRANSPORT}
		\iota(s,u,\vartheta)
		:= \exp
			\left(
				\int_{t'=0}^s
					- 2 \frac{\Lunit \upmu(t',u,\vartheta)}{\upmu(t',u,\vartheta)}
				\, dt'
			\right)
			= \frac{\upmu^2(0,u,\vartheta)}{\upmu^2(s,u,\vartheta)}
	\end{align}
	corresponding to the coefficient of
	$
	\upchifullmodarg{\GeoAng^N}
	$
	on LHS~\eqref{E:TOPORDERTRCHIJUNKRENORMALIZEDTRANSPORT}.
	We then rewrite \eqref{E:TOPORDERTRCHIJUNKRENORMALIZEDTRANSPORT} as
	$
	\displaystyle
		\Lunit
		\left(
			\iota \upchifullmodarg{\GeoAng^N}
		\right)
		= \iota
		\times \mbox{\upshape RHS~\eqref{E:TOPORDERTRCHIJUNKRENORMALIZEDTRANSPORT}}
	$
	and integrate the resulting equation
	with respect to $s$ from time $0$ to time $t$.
	From Def.~\ref{D:REGIONSOFDISTINCTUPMUBEHAVIOR}
	and the estimates
	\eqref{E:LOCALIZEDMUCANTGROWTOOFAST}
	and
	\eqref{E:LOCALIZEDMUMUSTSHRINK},
	we deduce
		\begin{align} \label{E:CRUDEMUOVERMUBOUND}
			\sup_{0 \leq s' \leq t}
				\frac{\upmu(t,u,\vartheta)}{\upmu(s',u,\vartheta)}
				& \leq C.
		\end{align}
		From \eqref{E:INTFACTFULLYMODIFIEDTRANSPORT}
		and \eqref{E:CRUDEMUOVERMUBOUND},
		it is straightforward to see that the desired bound
		\eqref{E:RENORMALIZEDTOPORDERTRCHIJUNKTRANSPORTINVERTED}
		follows once we establish the following bounds for 
		the terms generated by the terms on RHS~\eqref{E:TOPORDERTRCHIJUNKRENORMALIZEDTRANSPORT}:
		\begin{align}
			& \left|
					\upmu [\Lunit, \GeoAng^N] \mytr \upchi
				\right|
			(s,u,\vartheta)
				\label{E:FIRSTRERNORMALIZEDTRANSPORTEQNINHOMOGENEOUSTERMBOUND}  
				\\
			& \leq
			C \varepsilon
			\left|
				\upchifullmodarg{\GeoAng^N} 
			\right|
			(s,u,\vartheta)
			\notag	\\
		& \ \
			+
			C \varepsilon
			\left|
				\Fullset_{\ast}^{[1,N+1];\leq 2} \threePsi
			\right|
			(s,u,\vartheta)
			+
			C \varepsilon
			\left|
				\Fullset_{\ast}^{[1,N];\leq 1} \GdVar
			\right|
			(s,u,\vartheta)
			+
			C \varepsilon
			\left|
				\Fullset_{\ast \ast}^{[1,N];\leq 1} \BadVar
			\right|
			(s,u,\vartheta),
				\notag \\
			&
			2 
			\left|
				\upmu \mytr \upchi \GeoAng^N \mytr \upchi
			\right|
			(s,u,\vartheta)
				 \label{E:SECONDRERNORMALIZEDTRANSPORTEQNINHOMOGENEOUSTERMBOUND} 
				\\
			& \leq
				C \varepsilon
			\left|
				\upchifullmodarg{\GeoAng^N} 
			\right|
			(s,u,\vartheta)
			\notag	\\
		& \ \
			+
			C \varepsilon
			\left|
				\Fullset_{\ast}^{[1,N+1];\leq 2} \threePsi
			\right|
			(s,u,\vartheta)
			+
			C \varepsilon
			\left|
				\Fullset_{\ast}^{[1,N];\leq 1} \GdVar
			\right|
			(s,u,\vartheta)
			+
			C \varepsilon
			\left|
				\Fullset_{\ast \ast}^{[1,N];\leq 1} \BadVar
			\right|
			(s,u,\vartheta),
				\notag \\
			& 2 
			\left(
				\frac{\upmu(t,u,\vartheta)}{\upmu(s,u,\vartheta)}
			\right)^2 
			\left|
				\frac{\Lunit \upmu(s,u,\vartheta)}{\upmu(s,u,\vartheta)}
			\right|
			\left|
				\GeoAng^N \upchifullmodinhom
			\right|
			(s,u,\vartheta)
				\label{E:THIRDRERNORMALIZEDTRANSPORTEQNINHOMOGENEOUSTERMBOUND} \\
			& \leq
				2 (1 + C \varepsilon)
				\left|
					\frac{[\Lunit \upmu]_-(s,u,\vartheta)}{\upmu(s,u,\vartheta)}
				\right|
				\left|
					\GeoAng^N \upchifullmodinhom
				\right|
				(s,u,\vartheta)
				\notag \\
			& \ \
				+ 
		C
		\left|
			\Fullset_{\ast}^{[1,N+1];\leq 2} \threePsi
		\right|
		(s,u,\vartheta)
		+
		C
		\left|
			\Fullset_{\ast}^{[1,N];\leq 1} \GdVar
		\right|
		(s,u,\vartheta)
		+
		C
		\left|
			\Fullset_{\ast \ast}^{[1,N];\leq 1} \BadVar
		\right|
			(s,u,\vartheta),
			 \notag	\\
			& \mbox{all remaining terms on RHS~\eqref{E:TOPORDERTRCHIJUNKRENORMALIZEDTRANSPORT} } 
				\mbox{ are in magnitude}
				\label{E:REMAININGRERNORMALIZEDTRANSPORTEQNINHOMOGENEOUSTERMBOUND} \\
			& \leq 
			C 
			\left|
				\Fullset_{\ast}^{[1,N+1];\leq 2} \threePsi
			\right|
			(s,u,\vartheta)
			+
			C 
			\left|
				\Fullset_{\ast}^{[1,N];\leq 1} \GdVar
			\right|
			(s,u,\vartheta)
			+
			C 
			\left|
				\Fullset_{\ast \ast}^{[1,N];\leq 1} \BadVar
			\right|
			(s,u,\vartheta)
				\notag \\
		& \ \
			+
			\upmu
			\left|
				\Fullset^{N+1;\leq 1} \Vortrenormalized
			\right|
			(s,u,\vartheta)
			+
			\left|
				\Fullset^{[1,N];\leq 1} \Vortrenormalized
			\right|
			(s,u,\vartheta).
			\notag
		\end{align}
		We note that in deriving \eqref{E:RENORMALIZEDTOPORDERTRCHIJUNKTRANSPORTINVERTED},
		the product 
		$ C 
			\varepsilon
			\iota 
			\left|
				\upchifullmodarg{\GeoAng^N} 
			\right|
		$
		arising from the first term on RHS \eqref{E:SECONDRERNORMALIZEDTRANSPORTEQNINHOMOGENEOUSTERMBOUND}
		needs to be treated with Gronwall's inequality.
		However, due to the small factor $\varepsilon$,
		this product has only the negligible effect of 
		contributing to the factors 
		$C \varepsilon$
		on RHS \eqref{E:RENORMALIZEDTOPORDERTRCHIJUNKTRANSPORTINVERTED}.
		We remark that some of the estimates 
		\eqref{E:FIRSTRERNORMALIZEDTRANSPORTEQNINHOMOGENEOUSTERMBOUND}-\eqref{E:REMAININGRERNORMALIZEDTRANSPORTEQNINHOMOGENEOUSTERMBOUND}
	are non-optimal in the sense of the number of $\Rad$ derivatives 
	allowed on the right-hand sides. However, later in the proof,
	when we are analyzing $\upchifullmodarg{\GeoAng^{N-1} \Rad}$, 
	the same number of $\Rad$ derivatives appear on the right-hand sides of 
	the analogous estimates,
	and they are optimal.

		We now prove \eqref{E:FIRSTRERNORMALIZEDTRANSPORTEQNINHOMOGENEOUSTERMBOUND}.
		We first use the commutator estimate \eqref{E:DETAILEDPERMUTEDVECTORFIELDSACTINGONFUNCTIONSCOMMUTATORESTIMATE}
		with $M=0$ and $f = \mytr \upchi$,
		the $L^{\infty}$ estimates of Prop.~\ref{P:IMPROVEMENTOFAUX},
		\eqref{E:LUNITCOMMUTEDLUNITSMALLIPOINTWISE} with $M=0$,
		and \eqref{E:POINTWISEESTIMATESFORCHIANDITSDERIVATIVES}
		to deduce that
		$
		\left|
			\upmu [\Lunit, \GeoAng^N] \mytr \upchi
		\right|
		\lesssim
		\varepsilon 
		\left|
			\upmu \GeoAng^N \mytr \upchi
		\right|
		$
		plus error terms with magnitude
		$\leq$ the sum of the last three terms on
		RHS~\eqref{E:FIRSTRERNORMALIZEDTRANSPORTEQNINHOMOGENEOUSTERMBOUND}.
		We then use definition \eqref{E:TRANSPORTRENORMALIZEDTRCHIJUNK}  and
		the estimate \eqref{E:CRUDEXPOINTWISEBOUND}
		to deduce that
		$
		\varepsilon 
		\left|
			\upmu \GeoAng^N \mytr \upchi
		\right|
		= 
		\varepsilon 
		\left|
			\upchifullmodarg{\GeoAng^N}
		\right|
		$
		plus error terms with magnitude
		$\leq$ the sum of the last three terms on
		RHS~\eqref{E:FIRSTRERNORMALIZEDTRANSPORTEQNINHOMOGENEOUSTERMBOUND}, 
		which yields the desired bound.
		Inequality \eqref{E:SECONDRERNORMALIZEDTRANSPORTEQNINHOMOGENEOUSTERMBOUND} can be proved
		using similar arguments (without the help of a commutator estimate).

		We now prove \eqref{E:THIRDRERNORMALIZEDTRANSPORTEQNINHOMOGENEOUSTERMBOUND}. 
		We first note the simple inequality
		$
		\displaystyle
		\left|
			\frac{\Lunit \upmu(s,u,\vartheta)}{\upmu(s,u,\vartheta)}
		\right|
		\leq
		\left|
			\frac{[\Lunit \upmu]_-(s,u,\vartheta)}{\upmu(s,u,\vartheta)}
		\right|
		+
		\left|
			\frac{[\Lunit \upmu]_+(s,u,\vartheta)}{\upmu(s,u,\vartheta)}
		\right|
		$.
		To bound terms on LHS~\eqref{E:THIRDRERNORMALIZEDTRANSPORTEQNINHOMOGENEOUSTERMBOUND} 
		arising from the factor
		$
		\displaystyle
		\left|
			\frac{[\Lunit \upmu]_+(s,u,\vartheta)}{\upmu(s,u,\vartheta)}
		\right|
		$
		by $\leq$ the terms on the last line of RHS~\eqref{E:THIRDRERNORMALIZEDTRANSPORTEQNINHOMOGENEOUSTERMBOUND},
		we use \eqref{E:POSITIVEPARTOFLMUOVERMUISBOUNDED},
		\eqref{E:CRUDEMUOVERMUBOUND},
		and the estimate \eqref{E:CRUDEXPOINTWISEBOUND}.
		To bound terms on LHS \eqref{E:THIRDRERNORMALIZEDTRANSPORTEQNINHOMOGENEOUSTERMBOUND} 
		arising from the factor
$
\displaystyle
\left|
	\frac{[\Lunit \upmu]_-(s,u,\vartheta)}{\upmu(s,u,\vartheta)}
\right|
$,
we consider the partitions from Def.~\ref{D:REGIONSOFDISTINCTUPMUBEHAVIOR}.
When $(u,\vartheta) \in \Vplus{t}{u}$,
we use the bounds 
\eqref{E:KEYMUNOTDECAYINGMINUSPARTLMUOVERMUBOUND}
and
\eqref{E:CRUDEMUOVERMUBOUND}
to deduce that
$
\displaystyle
\left(
	\frac{\upmu(t,u,\vartheta)}{\upmu(s,u,\vartheta)}
\right)^2
\left|
	\frac{[\Lunit \upmu]_-(s,u,\vartheta)}{\upmu(s,u,\vartheta)}
\right|
\leq C \varepsilon
$.
Combining this bound with \eqref{E:CRUDEXPOINTWISEBOUND}, 
we easily conclude that the terms of interest
are $\leq$ the terms on the last line of
RHS~\eqref{E:THIRDRERNORMALIZEDTRANSPORTEQNINHOMOGENEOUSTERMBOUND}.
Finally, when $(u,\vartheta) \in \Vminus{t}{u}$,
we use \eqref{E:LOCALIZEDMUMUSTSHRINK}
to deduce that
\[
2
\left(
	\frac{\upmu(t,u,\vartheta)}{\upmu(s,u,\vartheta)}
\right)^2
\left|
	\frac{[\Lunit \upmu]_-(s,u,\vartheta)}{\upmu(s,u,\vartheta)}
\right|
\leq 2(1 + C \varepsilon) 
\left|
	\frac{[\Lunit \upmu]_-(s,u,\vartheta)}{\upmu(s,u,\vartheta)}
\right|.
\]
Thus, we conclude that the terms under consideration are
$\leq$ the terms on the first line of
RHS~\eqref{E:THIRDRERNORMALIZEDTRANSPORTEQNINHOMOGENEOUSTERMBOUND},
which completes the proof of \eqref{E:THIRDRERNORMALIZEDTRANSPORTEQNINHOMOGENEOUSTERMBOUND}.

We now prove \eqref{E:REMAININGRERNORMALIZEDTRANSPORTEQNINHOMOGENEOUSTERMBOUND},
starting with the estimate for the term
$
[\Lunit, \GeoAng^N] \upchifullmodinhom
$
on RHS~\eqref{E:TOPORDERTRCHIJUNKRENORMALIZEDTRANSPORT}.
Using the commutator estimate \eqref{E:DETAILEDPERMUTEDVECTORFIELDSACTINGONFUNCTIONSCOMMUTATORESTIMATE}
with $M=0$ and $f = \upchifullmodinhom$,
the $L^{\infty}$ estimates of Prop.~\ref{P:IMPROVEMENTOFAUX},
\eqref{E:LUNITCOMMUTEDLUNITSMALLIPOINTWISE} with $M=0$,
and the estimate \eqref{E:CRUDEXPOINTWISEBOUND}, we deduce that
$
\left|
	[\Lunit, \GeoAng^N] \upchifullmodinhom
\right|
\leq
\mbox{\upshape RHS~\eqref{E:REMAININGRERNORMALIZEDTRANSPORTEQNINHOMOGENEOUSTERMBOUND}}
$
as desired.
We next bound bound the term
$
[\upmu, \GeoAng^N] \Lunit \mytr \upchi
$
on RHS~\eqref{E:TOPORDERTRCHIJUNKRENORMALIZEDTRANSPORT}.
Using the $L^{\infty}$ estimates of Prop.~\ref{P:IMPROVEMENTOFAUX}
and the estimate \eqref{E:LUNITCOMMUTEDLUNITSMALLIPOINTWISE} 
with $M=0$, we deduce that
$
\left|
	[\upmu, \GeoAng^N] \Lunit \mytr \upchi
\right|
\leq
\mbox{\upshape RHS~\eqref{E:REMAININGRERNORMALIZEDTRANSPORTEQNINHOMOGENEOUSTERMBOUND}}
$
as desired.
We next bound the term
$[\GeoAng^N,\Lunit \upmu] \mytr \upchi$
on RHS~\eqref{E:TOPORDERTRCHIJUNKRENORMALIZEDTRANSPORT}.
Using 
the estimate \eqref{E:POINTWISEESTIMATESFORCHIANDITSDERIVATIVES},
the $L^{\infty}$ estimates of Prop.~\ref{P:IMPROVEMENTOFAUX},
and the estimate \eqref{E:LLUNITUPMUFULLSETSTARCOMMUTEDPOINTWISEESTIMATE}
with $M=0$, we deduce that
$
\left|
	[\GeoAng^N,\Lunit \upmu] \mytr \upchi
\right|
\leq
\mbox{\upshape RHS~\eqref{E:REMAININGRERNORMALIZEDTRANSPORTEQNINHOMOGENEOUSTERMBOUND}}
$
as desired.
We now bound the term
$
\GeoAng^N \left(\upmu (\mytr \upchi)^2 \right)
- 
2 \upmu \mytr \upchi \GeoAng^N \mytr \upchi
$
on RHS~\eqref{E:TOPORDERTRCHIJUNKRENORMALIZEDTRANSPORT}.
By the Leibniz rule, we see that the magnitude of this term is
$
\displaystyle
\lesssim
\mathop{\sum_{N_1 + N_2 + N_3 \leq N}}_{N_2,N_3 \leq N-1}
\left|
	\GeoAng^{N_1} \upmu
\right|
\left|
	\GeoAng^{N_2} \mytr \upchi
\right|
\left|
	\GeoAng^{N_3} \mytr \upchi
\right|.
$
Hence, from
the estimate \eqref{E:POINTWISEESTIMATESFORCHIANDITSDERIVATIVES}
and the $L^{\infty}$ estimates of Prop.~\ref{P:IMPROVEMENTOFAUX},
we deduce that all products in the sum are
$
\leq
\mbox{\upshape RHS~\eqref{E:REMAININGRERNORMALIZEDTRANSPORTEQNINHOMOGENEOUSTERMBOUND}}
$
as desired.
Finally, to bound the term
$ 
\GeoAng^N \mathfrak{A}
$
on RHS~\eqref{E:TOPORDERTRCHIJUNKRENORMALIZEDTRANSPORT},
we apply $\GeoAng^N$ to both sides of
\eqref{E:RENORMALIZEDRICLLINHOMOGENEOUSTERM}.
We bound the products of interest in magnitude by
$
\leq
\mbox{\upshape RHS~\eqref{E:REMAININGRERNORMALIZEDTRANSPORTEQNINHOMOGENEOUSTERMBOUND}}
$
with the help of 
the estimates \eqref{E:ANGDIFFXI}-\eqref{E:ANGDIFFXNOSPECIALSTRUCTUREDIFFERENTIATEDPOINTWISE},
\eqref{E:POINTWISEESTIMATESFORGSPHEREANDITSNOSPECIALSTRUCTUREDERIVATIVES},
and the $L^{\infty}$ estimates of Prop.~\ref{P:IMPROVEMENTOFAUX}.
This completes the proof of \eqref{E:REMAININGRERNORMALIZEDTRANSPORTEQNINHOMOGENEOUSTERMBOUND}
and finishes the proof of \eqref{E:RENORMALIZEDTOPORDERTRCHIJUNKTRANSPORTINVERTED}.

We now derive the desired bound for
$
\left|\upchifullmodarg{\GeoAng^{N-1} \Rad} \right|(t,u,\vartheta)
$.
We first note that by using essentially the same arguments used in the proof of \eqref{E:RENORMALIZEDTOPORDERTRCHIJUNKTRANSPORTINVERTED},
we can show that
\eqref{E:FIRSTRERNORMALIZEDTRANSPORTEQNINHOMOGENEOUSTERMBOUND}-\eqref{E:REMAININGRERNORMALIZEDTRANSPORTEQNINHOMOGENEOUSTERMBOUND}
hold with the operator $\GeoAng^N$ on the LHS replaced by $\GeoAng^{N-1} \Rad$,
but with the following changes:
\eqref{E:FIRSTRERNORMALIZEDTRANSPORTEQNINHOMOGENEOUSTERMBOUND}-\eqref{E:SECONDRERNORMALIZEDTRANSPORTEQNINHOMOGENEOUSTERMBOUND}
are replaced with
\begin{align}
			& \left|
					\upmu [\Lunit, \GeoAng^{N-1} \Rad] \mytr \upchi
				\right|
			(s,u,\vartheta)
				\label{E:REVISEDFIRSTRERNORMALIZEDTRANSPORTEQNINHOMOGENEOUSTERMBOUND}  \\
			& \leq
			C \varepsilon
			\left|
				\upchifullmodarg{\GeoAng^{N-1} \Rad} 
			\right|
			(s,u,\vartheta)
			+
			\left|
				\upchifullmodarg{\GeoAng^N} 
			\right|
			(s,u,\vartheta)
				\notag \\
			& 
			\ \
		  +
			C \varepsilon
			\left|
				\Fullset_{\ast}^{[1,N+1];\leq 2} \threePsi
			\right|
			(s,u,\vartheta)
			+
			C \varepsilon
			\left|
				\Fullset_{\ast}^{[1,N];\leq 1} \GdVar
			\right|
			(s,u,\vartheta)
			+
			C \varepsilon
			\left|
				\Fullset_{\ast \ast}^{[1,N];\leq 1} \BadVar
			\right|
			(s,u,\vartheta),
			\notag
\end{align}
\begin{align}
			&
			2 
			\left|
				\upmu \mytr \upchi \GeoAng^{N-1} \Rad \mytr \upchi
			\right|
			(s,u,\vartheta)
				 \label{E:REVISEDSECONDRERNORMALIZEDTRANSPORTEQNINHOMOGENEOUSTERMBOUND} \\
			& \leq
			C \varepsilon
			\left|
				\upchifullmodarg{\GeoAng^{N-1} \Rad} 
			\right|
			(s,u,\vartheta)
			+
			C \varepsilon
			\left|
				\upchifullmodarg{\GeoAng^N} 
			\right|
			(s,u,\vartheta)
				\notag \\
			& \ \
				+
			C \varepsilon
			\left|
				\Fullset_{\ast}^{[1,N+1];\leq 2} \threePsi
			\right|
			(s,u,\vartheta)
			+
			C \varepsilon
			\left|
				\Fullset_{\ast}^{[1,N];\leq 1} \GdVar
			\right|
			(s,u,\vartheta)
			+
			C \varepsilon
			\left|
				\Fullset_{\ast \ast}^{[1,N];\leq 1} \BadVar
			\right|
			(s,u,\vartheta).
				\notag 
\end{align}
The new features are that the second term on RHS~\eqref{E:REVISEDFIRSTRERNORMALIZEDTRANSPORTEQNINHOMOGENEOUSTERMBOUND}
does not contain a small factor $\varepsilon$
and that both RHS~\eqref{E:REVISEDFIRSTRERNORMALIZEDTRANSPORTEQNINHOMOGENEOUSTERMBOUND}
and RHS~\eqref{E:REVISEDSECONDRERNORMALIZEDTRANSPORTEQNINHOMOGENEOUSTERMBOUND}
depend on $\upchifullmodarg{\GeoAng^N}$
(that is, the estimate for 
$\upchifullmodarg{\GeoAng^{N-1} \Rad}$
does not decouple from the one for 
$\upchifullmodarg{\GeoAng^N}$).
To obtain \eqref{E:REVISEDFIRSTRERNORMALIZEDTRANSPORTEQNINHOMOGENEOUSTERMBOUND},
we use the commutator estimate \eqref{E:DETAILEDPERMUTEDVECTORFIELDSACTINGONFUNCTIONSCOMMUTATORESTIMATE}
with $f = \mytr \upchi$ as before, 
but now with $M=1$.
Also using the $L^{\infty}$ estimates of Prop.~\ref{P:IMPROVEMENTOFAUX},
\eqref{E:LUNITCOMMUTEDLUNITSMALLIPOINTWISE} with $M=1$,
and \eqref{E:POINTWISEESTIMATESFORCHIANDITSDERIVATIVES}
we deduce that 
$
		\left|
			\upmu [\Lunit, \GeoAng^{N-1} \Rad] \mytr \upchi
		\right|
		\lesssim
		\varepsilon 
		\left|
			\upmu \GeoAng^{N-1} \Rad \mytr \upchi
		\right|
		+
		\left|
			\upmu \GeoAng^N \mytr \upchi
		\right|
		$
plus error terms that are bounded in magnitude by
$\lesssim \mbox{RHS~\eqref{E:REVISEDFIRSTRERNORMALIZEDTRANSPORTEQNINHOMOGENEOUSTERMBOUND}}$.
We then use definition \eqref{E:TRANSPORTRENORMALIZEDTRCHIJUNK}  and
		the estimate \eqref{E:CRUDEXPOINTWISEBOUND}
		to deduce that
		$
		\varepsilon 
		\left|
			\upmu \GeoAng^{N-1} \Rad \mytr \upchi
		\right|
		= 
		\varepsilon 
		\left|
			\upchifullmodarg{\GeoAng^{N-1} \Rad}
		\right|
		$
		plus error terms with magnitude
		$\leq$ the sum of the last three terms on
		RHS~\eqref{E:REVISEDFIRSTRERNORMALIZEDTRANSPORTEQNINHOMOGENEOUSTERMBOUND}
		and 
		$
		\left|
			\upmu \GeoAng^N \mytr \upchi
		\right|
		= 
		\left|
			\upchifullmodarg{\GeoAng^N}
		\right|
		$
		plus error terms with magnitude
		$\leq$ the sum of the last three terms on
		RHS~\eqref{E:REVISEDFIRSTRERNORMALIZEDTRANSPORTEQNINHOMOGENEOUSTERMBOUND}.
		We have thus proved \eqref{E:REVISEDFIRSTRERNORMALIZEDTRANSPORTEQNINHOMOGENEOUSTERMBOUND}.
		Inequality \eqref{E:REVISEDSECONDRERNORMALIZEDTRANSPORTEQNINHOMOGENEOUSTERMBOUND}
		can be proved using similar arguments (without the help of a commutator estimate).

		We now recall that
		we can rewrite \eqref{E:TOPORDERTRCHIJUNKRENORMALIZEDTRANSPORT} 
		(with $\GeoAng^{N-1} \Rad$ in the role of $\Fullset_{\ast}^{\leq N;\leq 1}$ in that equation)
		in the form
	$
	\displaystyle
		\Lunit
		\left(
			\iota \upchifullmodarg{\GeoAng^{N-1} \Rad}
		\right)
		= \iota
		\times \mbox{\upshape RHS~\eqref{E:TOPORDERTRCHIJUNKRENORMALIZEDTRANSPORT}}
	$
	and integrate the resulting equation
	with respect to $s$ from the initial time $0$ to time $t$.
	With the help of of the estimates obtained in the previous paragraph,
	we can obtain a pointwise
	estimate for $\iota \upchifullmodarg{\GeoAng^{N-1} \Rad}(t,u,\vartheta)$,
	much as in the case of $\iota \upchifullmodarg{\GeoAng^N}(t,u,\vartheta)$,
	where we use Gronwall's inequality to handle the first terms
	$
	\displaystyle
	C \varepsilon
	\left|
		\upchifullmodarg{\GeoAng^{N-1} \Rad} 
	\right|
	$
	on RHS~\eqref{E:REVISEDFIRSTRERNORMALIZEDTRANSPORTEQNINHOMOGENEOUSTERMBOUND}
	and
	RHS~\eqref{E:REVISEDSECONDRERNORMALIZEDTRANSPORTEQNINHOMOGENEOUSTERMBOUND}.
	The new step compared to the argument for $\upchifullmodarg{\GeoAng^N}$
	is that
	we insert the already proven bound \eqref{E:RENORMALIZEDTOPORDERTRCHIJUNKTRANSPORTINVERTED}
	in order to handle the terms
	$
	\displaystyle
	\left|
		\upchifullmodarg{\GeoAng^N} 
	\right|
	(s,u,\vartheta)
	$
	on RHS~\eqref{E:REVISEDFIRSTRERNORMALIZEDTRANSPORTEQNINHOMOGENEOUSTERMBOUND}
	and
	RHS~\eqref{E:REVISEDSECONDRERNORMALIZEDTRANSPORTEQNINHOMOGENEOUSTERMBOUND}.
	In view of the fact that we are using Gronwall's inequality,
	we see that the bound \eqref{E:RENORMALIZEDTOPORDERTRCHIJUNKTRANSPORTINVERTED}
	leads to the presence 
	(on the RHS of the pointwise estimate 
	$
	\displaystyle
	\left|
			\upchifullmodarg{\GeoAng^{N-1} \Rad} 
	\right|
	(t,u,\vartheta)
	\leq \cdots
	$
	)
	of additional integrals of the form
	\begin{align} \label{E:ADDITIONALTIMEINTEGRALS}
	\int_{s=0}^t
		\frac{\upmu^2(t,u,\vartheta)}{\upmu^2(s,u,\vartheta)}
		\times
		\mbox{\upshape RHS~\eqref{E:RENORMALIZEDTOPORDERTRCHIJUNKTRANSPORTINVERTED}}
		(s,u,\vartheta)
	\, ds,
	\end{align}
	which we did not encounter in our proof of \eqref{E:RENORMALIZEDTOPORDERTRCHIJUNKTRANSPORTINVERTED}.
	To handle these additional integrals, we first bound the factor
	$
	\displaystyle
	\frac{\upmu^2(t,u,\vartheta)}{\upmu^2(s,u,\vartheta)}
	$
	in \eqref{E:ADDITIONALTIMEINTEGRALS}
	by $\leq C$ with the help of \eqref{E:CRUDEMUOVERMUBOUND}.
	Next, we use the $L^{\infty}$ estimates of Prop.~\ref{P:IMPROVEMENTOFAUX}
	and the estimate \eqref{E:CRUDEXPOINTWISEBOUND}
	to bound the first integrand
	on RHS~\eqref{E:RENORMALIZEDTOPORDERTRCHIJUNKTRANSPORTINVERTED},
	evaluated at $(s',u,\vartheta)$,
	as follows:
	\begin{align} \label{E:INTEGRANDBOUNDFORALREADYGRONWALLEDTERM}
	\frac{[\Lunit \upmu(s',u,\vartheta)]_-}{\upmu(s',u,\vartheta)}
	\left| \GeoAng^N \upchifullmodinhom \right|(s',u,\vartheta)
	\lesssim 
	\frac{1}{\upmu_{\star}(s',u)}
			\left\lbrace
				\left|
					\Fullset_{\ast}^{[1,N+1];\leq 2} \Psi
				\right|
			+
			\left|
				\myarray
					[\Fullset_{\ast \ast}^{[1,N];\leq 1} \BadVar]
					{\Fullset_{\ast}^{[1,N];\leq 1} \GdVar}
			\right|
			\right\rbrace
			(s',u,\vartheta).
	\end{align}
	Inserting the estimate \eqref{E:INTEGRANDBOUNDFORALREADYGRONWALLEDTERM}
	into the first integrand on RHS~\eqref{E:RENORMALIZEDTOPORDERTRCHIJUNKTRANSPORTINVERTED}
	(with $s$ in the role of $t$ on RHS~\eqref{E:RENORMALIZEDTOPORDERTRCHIJUNKTRANSPORTINVERTED}
	and $s'$ in the role of $t$),
	we generate the double time integral 
	stated in \eqref{E:ADDITIONALTIMEINTEGRALFORFULLYMODIFIEDQUANTITYPOINTWISEESTIMATEONERADIALCASE}.
	The remaining two time integrals on  
	RHS~\eqref{E:INTEGRANDBOUNDFORALREADYGRONWALLEDTERM}
	also generate double time integrals,
	but they are less singular in that they do not involve
	the factor of
	$
	\displaystyle
	\frac{1}{\upmu_{\star}}
	$
	present on RHS~\eqref{E:INTEGRANDBOUNDFORALREADYGRONWALLEDTERM}.
	Hence, these double time integrals are
	$\lesssim$ the single time integrals on RHS~\eqref{E:RENORMALIZEDTOPORDERTRCHIJUNKTRANSPORTINVERTED},
	in view of the following simple bound, 
	which holds for non-negative scalar-valued functions $f$:
	\begin{align} \label{E:SIMPLEDOUBLETIMEINTEGRALINTERMSOFSINGLETIMEINTEGRALBOuND}
	\int_{s=0}^t
		\int_{s'=0}^s
			f(s',u,\vartheta)
		\, ds'
	\, ds
	& \leq
	\int_{s=0}^t
		\int_{s'=0}^t
			f(s',u,\vartheta)
		\, ds'
	\, ds
		\\
	& \leq 
	t
	\int_{s'=0}^t
		f(s',u,\vartheta)
	\, ds'
	\leq 
	C
	\int_{s=0}^t
		f(s,u,\vartheta)
	\, ds.
	\notag
	\end{align}
	Similarly, with the help of \eqref{E:CRUDEMUOVERMUBOUND},
	we bound the time integral 
	generated by the initial data term
	$
	\left| 
		\upchifullmodarg{\GeoAng^N } 
	\right|
	(0,u,\vartheta)
	$
	on RHS~\eqref{E:RENORMALIZEDTOPORDERTRCHIJUNKTRANSPORTINVERTED}
	by 
	$
	\displaystyle
	\leq C 
	\int_{s'=0}^t
		\left| 
			\upchifullmodarg{\GeoAng^N } 
		\right|
		(0,u,\vartheta)
		\, ds'
	\leq C
	\left| 
		\upchifullmodarg{\GeoAng^N } 
	\right|
	(0,u,\vartheta).
	$
	We have thus obtained the desired bound for 
	$
	\left|\upchifullmodarg{\GeoAng^{N-1} \Rad} \right|(t,u,\vartheta)
	$,
	which completes the proof of the lemma.

\end{proof}

Armed with Lemma~\ref{L:RENORMALIZEDTOPORDERTRCHIJUNKTRANSPORTINVERTED}, 
we now derive the main result of this section.

\begin{remark}[\textbf{Boxed constants affect high-order energy blowup-rates}]
	\label{R:BOXEDCONSTANTS}
The ``boxed constants'' such as the $\boxed{2}$ and $\boxed{4}$
appearing on the RHS of inequality \eqref{E:KEYPOINTWISEESTIMATE}
are important because they affect the blowup-rates 
(that is, the powers of $\upmu_{\star}^{-1}$)
featured on the right-hand sides of high-order energy estimates. 
Similar remarks apply to the boxed constants appearing 
on RHSs~\eqref{E:TOPORDERWAVEENERGYINTEGRALINEQUALITIES},
\eqref{E:DIFFICULTTERML2BOUND},
\eqref{E:MOSTDIFFICULTERRORINTEGRALBOUND},
\eqref{E:ANNOYINGLDERIVATIVEBOUNDRYSPATIALINTEGRALFACTORL2ESTIMATE},
and
\eqref{E:ANNOYINGBOUNDRYSPATIALINTEGRALFACTORL2ESTIMATE}.
\end{remark}

\begin{proposition}[\textbf{The key pointwise estimates for}
 $(\Rad v^1) \Fullset_{\ast}^{N;\leq 1} \mytr \upchi$]
	\label{P:KEYPOINTWISEESTIMATE}
	Assume that $N = 20$
	and let
	$\Fullset_{\ast}^{N;\leq 1} \in \lbrace \GeoAng^N, \GeoAng^{N-1} \Rad \rbrace$.
	There exists a constant $C_* > 0$ such that
	under the data-size and bootstrap assumptions 
	of Sects.~\ref{SS:FLUIDVARIABLEDATAASSUMPTIONS}-\ref{SS:PSIBOOTSTRAP}
	and the smallness assumptions of Sect.~\ref{SS:SMALLNESSASSUMPTIONS}, 
	the following pointwise estimate holds on
	$\mathcal{M}_{\Tboot,U_0}$:
	\begin{align} \label{E:KEYPOINTWISEESTIMATE}
		\left|
			(\Rad v^1) \Fullset_{\ast}^{N;\leq 1} \mytr \upchi
		\right|
		& \leq
			\boxed{2} 
			\left\|
				\frac{[\Lunit \upmu]_-}{\upmu}
			\right\|_{L^{\infty}(\Sigma_t^u)}
			\left| 
				\Rad \Fullset_{\ast}^{N;\leq 1} v^1 
			\right|
				\\
		& \ \
				+ C_*
					\frac{1}{\upmu_{\star}(t,u)}
					\left| 
						\Rad \Fullset_{\ast}^{N;\leq 1} (\Densrenormalized  - v^1)
					\right|
			\notag \\
		 &  \ \ + 
						\boxed{4} (1 + C \varepsilon)
						\frac
						{
							\left\|
								[\Lunit \upmu]_-
							\right\|_{L^{\infty}(\Sigma_t^u)}}
						{\upmu_{\star}(t,u)}
					\int_{t'=0}^t 
						\frac
						{
							\left\|
								[\Lunit \upmu]_-
							\right\|_{L^{\infty}(\Sigma_{t'}^u)}}
						{\upmu_{\star}(t',u)}
						\left| 
							\Rad \Fullset_{\ast}^{N;\leq 1} v^1
						\right|(t',u,\vartheta)
				\, dt'
				\notag \\
	&  \ \ + 
					C_*
					\frac{1}{\upmu_{\star}(t,u)}
					\int_{t'=0}^t 
						\frac
						{1}
						{\upmu_{\star}(t',u)}
						\left| 
							\Rad \Fullset_{\ast}^{N;\leq 1} (\Densrenormalized  - v^1)
						\right|(t',u,\vartheta)
				\, dt'
				\notag \\
	& \ \	+ \mbox{\upshape Error},
		\notag 
	\end{align}
	where
\begingroup
\allowdisplaybreaks
	\begin{align}  \label{E:ERRORTERMKEYPOINTWISEESTIMATE}
		\left|
			\mbox{\upshape Error}
		\right|
		(t,u,\vartheta)
		& \lesssim
			\frac{1}{\upmu_{\star}(t,u)}
			\left\lbrace
				\left| 
					\upchifullmodarg{\GeoAng^N} 
				\right|
				+
				\left| 
					\upchifullmodarg{\GeoAng^{N-1} \Rad} 
				\right|
			\right\rbrace
			(0,u,\vartheta)
				\\
	& \ \
			+ 
			\varepsilon
			\frac{1}{\upmu_{\star}}
			\left|
			\Fullset_{\ast}^{[1,N+1];\leq 2} \threePsi
			\right|
			(t,u,\vartheta)
			+ 
			\left|
			\Fullset_{\ast}^{[1,N+1];\leq 2} \threePsi
			\right|
			(t,u,\vartheta)
			\notag	\\
			& 
			\ \
			+
			\frac{1}{\upmu_{\star}(t,u)}
			\left|
				\Fullset_{\ast}^{[1,N];\leq 2} \threePsi
			\right|
			(t,u,\vartheta)
			+
			\frac{1}{\upmu_{\star}(t,u)}
			\left|
				\myarray
					[\Fullset_{\ast \ast}^{[1,N];\leq 1} \BadVar]
					{\Fullset_{\ast}^{[1,N];\leq 1} \GdVar}
			\right|
			(t,u,\vartheta)
				\notag \\
		& \ \
		+
		\frac{1}{\upmu_{\star}(t,u)}
		\int_{t'=0}^t
		\int_{s=0}^{t'}
			\frac{1}{\upmu_{\star}(s,u)}
			\left\lbrace
				\left|
					\Fullset_{\ast}^{[1,N+1];\leq 2} \threePsi
				\right|
			+
			\left|
				\myarray
					[\Fullset_{\ast \ast}^{[1,N];\leq 1} \BadVar]
					{\Fullset_{\ast}^{[1,N];\leq 1} \GdVar}
			\right|
			\right\rbrace
			(s,u,\vartheta)
		\, ds
	\, dt'
		\notag \\
	& \ \ 
			  + 
				\varepsilon
				\frac{1}{\upmu_{\star}(t,u)}
				\int_{t'=0}^t
					\frac{1}{\upmu_{\star}(t',u)}
					\left|
						\Fullset_{\ast}^{[1,N+1];\leq 2} \threePsi
					\right|
					(t',u,\vartheta)
				\, dt'
				\notag 
				\\
	& \ \ 
			  + 
				\frac{1}{\upmu_{\star}(t,u)}
				\int_{t'=0}^t
					\left|
						\Fullset_{\ast}^{[1,N+1];\leq 2} \threePsi
					\right|
					(t',u,\vartheta)
				\, dt'
				\notag 
				\\
		& \ \
				+ 
				\frac{1}{\upmu_{\star}(t,u)}
				\int_{t'=0}^t
					\frac{1}{\upmu_{\star}(t',u)}
					\left\lbrace
						\left|
							\Fullset_{\ast}^{[1,N];\leq 2} \threePsi
						\right|
						+
						\left|
							\myarray
								[\Fullset_{\ast \ast}^{[1,N];\leq 1} \BadVar]
								{\Fullset_{\ast}^{[1,N];\leq 1} \GdVar}
						\right|
					\right\rbrace
					(t',u,\vartheta)
				\, dt'
				\notag
					\\
			& \ \
				+ 	\frac{1}{\upmu_{\star}(t,u)}
						\int_{t'=0}^t 
							\upmu(t',u,\vartheta)
								\left|
									\Fullset_{\ast}^{N+1;\leq 1} \Vortrenormalized
								\right|
							(t',u,\vartheta)
						\, dt'
						\notag \\
			& \ \
						+ 
						\frac{1}{\upmu_{\star}(t,u)}
						\int_{t'=0}^t 
							\left|
								\Fullset^{\leq N;\leq 1} \Vortrenormalized
							\right|
							(t',u,\vartheta)
						\, dt'.
					\notag
		\end{align}
		\endgroup

	Furthermore, we have the following less precise pointwise estimate:
\begingroup
\allowdisplaybreaks
	\begin{align} \label{E:LESSPRECISEKEYPOINTWISEESTIMATE}
		& 
		\left|
			\upmu \Fullset_{\ast}^{N;\leq 1} \mytr \upchi
		\right|
		\\
		& \lesssim
		\left\lbrace
			\left| 
				\upchifullmodarg{\GeoAng^N} 
			\right|
			+
			\left| 
				\upchifullmodarg{\GeoAng^{N-1} \Rad} 
			\right|
		\right\rbrace
		(0,u,\vartheta)
			\notag \\
	& \ \
		+
		\upmu
		\left|
			\Tanset^{N+1} \threePsi
		\right|
		+
		\left|
			\Rad \Fullset_{\ast}^{N;\leq 1} \threePsi
		\right|
		\notag \\
	& \ \ 
		+
		\left|
			\Fullset_{\ast}^{[1,N];\leq 2} \threePsi
		\right|
		+
		\left|
			\myarray
				[\Fullset_{\ast \ast}^{[1,N];\leq 1} \BadVar]
				{\Fullset_{\ast}^{[1,N];\leq 1} \GdVar}
		\right|
		\notag	\\
		& \ \
				+
				\int_{t'=0}^t 
						\frac{1}{\upmu_{\star}(t',u)}
						\left| 
							\Rad \Fullset_{\ast}^{N;\leq 1} \threePsi
						\right|
						(t',u,\vartheta)
				\, dt'
				+ 
				\int_{t'=0}^t 
					\left|
						\Fullset_{\ast}^{N+1;2} \threePsi
					\right|
					(t',u,\vartheta)
				\, dt'
							\notag \\
		& \ \
		+
		\int_{t'=0}^t
		\int_{s=0}^{t'}
			\frac{1}{\upmu_{\star}(s,u)}
			\left\lbrace
				\left|
					\Fullset_{\ast}^{[1,N+1];\leq 2} \threePsi
				\right|
			+
			\left|
				\myarray
					[\Fullset_{\ast \ast}^{[1,N];\leq 1} \BadVar]
					{\Fullset_{\ast}^{[1,N];\leq 1} \GdVar}
			\right|
			\right\rbrace
			(s,u,\vartheta)
		\, ds
	\, dt'
		\notag \\
	& \ \  
					+ 
					\int_{t'=0}^t
					\frac{1}{\upmu_{\star}(t',u)}
					\left\lbrace
						\left|
							\Fullset_{\ast}^{[1,N];\leq 2} \threePsi
						\right|
						+
						\left|
							\myarray
								[\Fullset_{\ast \ast}^{[1,N];\leq 1} \BadVar]
								{\Fullset_{\ast}^{[1,N];\leq 1} \GdVar}
						\right|
					\right\rbrace
					(t',u,\vartheta)
				\, dt'
						\notag
					\\
			& \ \
					+ 
						\int_{t'=0}^t 
								\upmu(t',u,\vartheta)
								\left|
									\Fullset_{\ast}^{N+1;\leq 1} \Vortrenormalized
								\right|
								(t',u,\vartheta)
						\, dt'
							\notag \\
			& \ \
						+ 
						\int_{t'=0}^t 
							\left|
								\Fullset^{\leq N;\leq 1} \Vortrenormalized
							\right|
							(t',u,\vartheta)
						\, dt'.
		\notag
	\end{align}
	\endgroup
\end{proposition}

\begin{proof}
	See Sect.~\ref{SS:OFTENUSEDESTIMATES} for some comments on the analysis.
	Throughout this proof,
	$\mbox{\upshape Error}$ denotes any term verifying the estimate \eqref{E:ERRORTERMKEYPOINTWISEESTIMATE}.
	We first prove \eqref{E:KEYPOINTWISEESTIMATE} in the case
	$\Fullset_{\ast}^{N;\leq 1} = \GeoAng^N$.
	Using \eqref{E:TRANSPORTRENORMALIZEDTRCHIJUNK}-\eqref{E:LOWESTORDERTRANSPORTRENORMALIZEDTRCHIJUNKDISCREPANCY},
	the estimate \eqref{E:TOPDERIVATIVESOFXPOINTWISEBOUND},
	and the simple bound
	$\| \Rad v^1 \|_{L^{\infty}} \lesssim 1$
	(see \eqref{E:PSITRANSVERSALLINFINITYBOUNDBOOTSTRAPIMPROVEDLARGE}),
	we decompose
	\begin{align} \label{E:DIFFICULTPRODUCTDECOMPOSITION}
		(\Rad v^1) \GeoAng^N \mytr \upchi
		& 
		= \frac{1}{\upmu} (\Rad v^1) \upchifullmodarg{\GeoAng^N}
			 +
			\frac{1}{\upmu} (\Rad v^1) \vec{G}_{\Lunit \Lunit}\contr\GeoAng^N \Rad \threePsi
			+ \mbox{\upshape Error}.
	\end{align}
	Next, we use the commutator estimate
	\eqref{E:LESSPRECISEPERMUTEDVECTORFIELDSACTINGONFUNCTIONSCOMMUTATORESTIMATE} with 
	$f = \threePsi$,
	the schematic identity \eqref{E:GFRAMESCALARSDEPENDINGONGOODVARIABLES},
	and the $L^{\infty}$ estimates of Prop.~\ref{P:IMPROVEMENTOFAUX}
	to obtain
	\begin{align} \label{E:SIMPLECOMMESTIMATE}
		\frac{1}{\upmu} 
		(\Rad v^1) 
		\vec{G}_{\Lunit \Lunit}
		\contr \GeoAng^N \Rad \threePsi
		& 
		=
		\frac{1}{\upmu} 
		(\Rad v^1) 
		\vec{G}_{\Lunit \Lunit}
		\contr \Rad \GeoAng^N \threePsi
		+ \mbox{\upshape Error}.
	\end{align}
	Recalling the definition \eqref{E:THREEVECPSI} of $\threePsi$
	and that
	$
	\vec{G}_{\Lunit \Lunit}\contr\Rad \threePsi
	=
	\sum_{\imath= 0}^2
	G_{\Lunit \Lunit}^{\imath} \Rad \Psi_{\imath}
	$,
	and using the transport equation \eqref{E:UPMUFIRSTTRANSPORT},
	we compute that
	\begin{align} \label{E:HARDPOINTWISEESTIMATEFORPRODUCTDECOMPOSITIONOFMAINTERM}
		&
		\frac{1}{\upmu} 
		(\Rad v^1) 
		\vec{G}_{\Lunit \Lunit}\contr\Rad \GeoAng^N \threePsi
			\\
		&
		= 2
			\left(
				\frac{\Lunit \upmu}{\upmu} 
			\right)
			\Rad \GeoAng^N v^1
			\notag \\
		& \ \
			+
			\frac{1}{\upmu} 
			(\Rad v^1) 
			G_{\Lunit \Lunit}^0 
			\Rad \GeoAng^N (\Densrenormalized  - v^1)
			\notag	\\
		& \ \
			-
			\frac{1}{\upmu} 
			(\Rad v^2) 
			G_{\Lunit \Lunit}^2 
			\Rad \GeoAng^N v^1
			+
			\frac{1}{\upmu} 
			(\Rad v^1) 
			G_{\Lunit \Lunit}^2 
			\Rad \GeoAng^N v^2
			+
			\frac{1}{\upmu} 
			\left\lbrace
				\Rad(v^1 - \Densrenormalized) 
			\right\rbrace
			G_{\Lunit \Lunit}^0 
			\Rad \GeoAng^N v^1
			\notag	\\
		& \ \
			+ 
			\left(
				\vec{G}_{\Lunit \Lunit} \contr \Lunit \threePsi
			\right)
			\Rad \GeoAng^N v^1
			+ 
			2
			\left(
				\vec{G}_{\Lunit \Radunit} \contr \Lunit \threePsi
			\right)
			\Rad \GeoAng^N v^1.
			\notag 
	\end{align}
	Using the schematic identity \eqref{E:GFRAMESCALARSDEPENDINGONGOODVARIABLES}
	and the $L^{\infty}$ estimates of Prop.~\ref{P:IMPROVEMENTOFAUX},
	we deduce that the product on the second line of
	RHS~\eqref{E:HARDPOINTWISEESTIMATEFORPRODUCTDECOMPOSITIONOFMAINTERM}
	is bounded in magnitude by the second term $C_* \cdots$ on RHS~\eqref{E:KEYPOINTWISEESTIMATE}.
	Using the $L^{\infty}$ estimates of Prop.~\ref{P:IMPROVEMENTOFAUX}
	(in particular \eqref{E:PSITRANSVERSALLINFINITYBOUNDBOOTSTRAPIMPROVEDSMALL}
	and \eqref{E:CRUCIALPSITRANSVERSALLINFINITYBOUNDBOOTSTRAPIMPROVEDSMALL})
	and the estimate \eqref{E:SMALLNESSOFGLL2},
	we deduce that the terms on the last two lines
	of RHS~\eqref{E:HARDPOINTWISEESTIMATEFORPRODUCTDECOMPOSITIONOFMAINTERM}
	are bounded in magnitude by
	$\lesssim$ the terms on the second line of RHS~\eqref{E:ERRORTERMKEYPOINTWISEESTIMATE}
	and are therefore $\mbox{\upshape Error}$.
	To bound the product
	$
	\displaystyle
	2
			\left(
				\frac{\Lunit \upmu}{\upmu} 
			\right)
			\Rad \GeoAng^N v^1
	$,
	we first split
	$\Lunit \upmu = [\Lunit \upmu]_+ - [\Lunit \upmu]_-$.
	From \eqref{E:POSITIVEPARTOFLMUOVERMUISBOUNDED}, we find that
	the product corresponding to $[\Lunit \upmu]_+$
	is $\mbox{\upshape Error}$, while the product corresponding to $[\Lunit \upmu]_-$
	is clearly bounded in magnitude by the first term
	on RHS~\eqref{E:KEYPOINTWISEESTIMATE}.

	It remains for us to bound the first product 
	$
	\displaystyle
	\frac{1}{\upmu} (\Rad v^1) \upchifullmodarg{\GeoAng^N}
	$
	on RHS~\eqref{E:DIFFICULTPRODUCTDECOMPOSITION}.
	Our argument is based on equation \eqref{E:RENORMALIZEDTOPORDERTRCHIJUNKTRANSPORTINVERTED}.
	To proceed, we multiply both sides of \eqref{E:RENORMALIZEDTOPORDERTRCHIJUNKTRANSPORTINVERTED}
	by 
	$
	\displaystyle
	\frac{1}{\upmu} (\Rad v^1) 
	$.
	We first bound the product generated by the first time integral 
	$2 (1 + C \varepsilon) \cdots$
	on RHS~\eqref{E:RENORMALIZEDTOPORDERTRCHIJUNKTRANSPORTINVERTED}.
	Using \eqref{E:TOPDERIVATIVESOFXPOINTWISEBOUND}
	and the simple bounds 
	$
	\displaystyle
	|[\Lunit \upmu]_-|(s,u,\vartheta) 
	\lesssim |\Lunit \upmu|(s,u,\vartheta) \lesssim 1
	$
	(that is, \eqref{E:LUNITUPTOONETRANSVERSALUPMULINFINITY})
	and 
	$
	\displaystyle
	\left|
		\frac{1}{\upmu} (\Rad v^1) 
	\right|
	(t,u)
	\lesssim 
	\frac{1}{\upmu_{\star}(t,u)}
	$
	(see \eqref{E:PSITRANSVERSALLINFINITYBOUNDBOOTSTRAPIMPROVEDLARGE}),
	we express the product under consideration as
	\begin{align} \label{E:MOSTDIFFICULTTIMEINTEGRALPRODUCTPOINTWISEBOUND}
	- 2 (1 + C \varepsilon)
							\frac{1}{\upmu} (\Rad v^1) 
							\int_{s=0}^t 
								\frac{[\Lunit \upmu]_-(s,u,\vartheta)}{\upmu(s,u,\vartheta)}
								\left| 
									\vec{G}_{\Lunit \Lunit}\contr
									\Rad \GeoAng^N \threePsi
								\right|(s,u,\vartheta)
							\, ds
	+
	\mbox{\upshape Error}.
	\end{align}
	Next, we algebraically decompose the second factor in the integrand 
	in \eqref{E:MOSTDIFFICULTTIMEINTEGRALPRODUCTPOINTWISEBOUND} as
	\begin{align} \label{E:SPLITTINGOFGLLSTARRADGEOANGNVECPSI}
		\vec{G}_{\Lunit \Lunit} \contr \Rad \GeoAng^N \threePsi
		=
		\sum_{\imath=0}^1 G_{\Lunit \Lunit}^{\imath} \Rad \GeoAng^N v^1
		+
		G_{\Lunit \Lunit}^2 \Rad \GeoAng^N v^2
		+
		G_{\Lunit \Lunit}^0 \Rad \GeoAng^N (\Densrenormalized  - v^1).
	\end{align}
	Using the schematic identity \eqref{E:GFRAMESCALARSDEPENDINGONGOODVARIABLES}
	and the $L^{\infty}$ estimates of Prop.~\ref{P:IMPROVEMENTOFAUX},
	we bound the magnitude of the time integral corresponding to
	the product
	$
	G_{\Lunit \Lunit}^0 \Rad \GeoAng^N (\Densrenormalized  - v^1)
	$
	in \eqref{E:SPLITTINGOFGLLSTARRADGEOANGNVECPSI}
	by $\leq$ the fourth term 
	$
	\displaystyle
					C_*
					\frac{1}{\upmu_{\star}(t,u)}
					\int_{t'=0}^t 
						\frac
						{1}
						{\upmu_{\star}(t',u)}
						\left| 
							\Rad \Fullset_{\ast}^{N;\leq 1} (\Densrenormalized  - v^1)
						\right|(t',u,\vartheta)
				\, dt'
	$
	on RHS~\eqref{E:KEYPOINTWISEESTIMATE}.
	Next, using \eqref{E:SMALLNESSOFGLL2},
	we bound the magnitude of the time integral corresponding to
	the product
	$
	G_{\Lunit \Lunit}^2 \Rad \GeoAng^N v^2
	$
	in \eqref{E:SPLITTINGOFGLLSTARRADGEOANGNVECPSI}
	by the time-integral-involving product on RHS~\eqref{E:ERRORTERMKEYPOINTWISEESTIMATE}
	featuring the small coefficient $\varepsilon$
	(and thus the time integral under consideration is of the form $\mbox{\upshape Error}$).
	We now bound the remaining time integral 
	\[
	 2 (1 + C \varepsilon)
							\frac{1}{\upmu} (\Rad v^1) 
							\int_{s=0}^t 
								\frac{[\Lunit \upmu(s,u,\vartheta)]_-}{\upmu(s,u,\vartheta)}
								\left| 
									\sum_{\imath=0}^1 
										G_{\Lunit \Lunit}^{\imath}(s,u,\vartheta) 
										\Rad \GeoAng^N v^1(s,u,\vartheta)
								\right|
							\, ds,
	\]
	which is generated by the sum in \eqref{E:SPLITTINGOFGLLSTARRADGEOANGNVECPSI}.
	We first algebraically decompose
	\begin{align}
	& \sum_{\imath=0}^1 G_{\Lunit \Lunit}^{\imath}(s,u,\vartheta) \Rad \GeoAng^N v^1(s,u,\vartheta)
		\\
	& 
	=
	\sum_{\imath=0}^1 G_{\Lunit \Lunit}^{\imath}(t,u,\vartheta) \Rad \GeoAng^N v^1(s,u,\vartheta)
	+
	\left\lbrace
		\sum_{\imath=0}^1 G_{\Lunit \Lunit}^{\imath}(s,u,\vartheta) 
		-
		\sum_{\imath=0}^1 G_{\Lunit \Lunit}^{\imath}(t,u,\vartheta)
	\right\rbrace
	\Rad \GeoAng^N v^1(s,u,\vartheta).
		\notag
	\end{align}
	Using the estimate \eqref{E:RADDERIVATIVESOFGLLDIFFERENCEBOUND}
	and the bounds
	$
	\displaystyle
	|[\Lunit \upmu]_-|
	\lesssim 1
	$
	and
	$
	\displaystyle
	\left|
		\frac{1}{\upmu} (\Rad v^1) 
	\right|
	(t,u)
	\lesssim 
	\frac{1}{\upmu_{\star}(t,u)}
	$
	noted above,
	we bound the magnitude of the time integral
	featuring the integrand factor
	$
	\sum_{\imath=0}^1 G_{\Lunit \Lunit}^{\imath}(s,u,\vartheta) 
	-
	\sum_{\imath=0}^1 G_{\Lunit \Lunit}^{\imath}(t,u,\vartheta)
	$
	by 
	$\lesssim$
	the time-integral-involving product on RHS~\eqref{E:ERRORTERMKEYPOINTWISEESTIMATE}
	featuring the small coefficient $\varepsilon$
	(and thus the time integral under consideration is of the form $\mbox{\upshape Error}$).
	The remaining time integral that we must estimate
	contains the integrand factor 
	$\sum_{\imath=0}^1 G_{\Lunit \Lunit}^{\imath}(t,u,\vartheta)$,
	which we may pull out of the $\, ds$ integral.
	That is, we must bound
	\begin{align} \label{E:LASTTIMEINTEGRALTOBOUND}
				2 (1 + C \varepsilon)
							\left| 
								\frac{1}{\upmu} 
								(\Rad v^1) 
							\sum_{\imath=0}^1 G_{\Lunit \Lunit}^{\imath}
							\right|
							(t,u,\vartheta)
							\int_{s=0}^t 
								\frac{[\Lunit \upmu]_-(s,u,\vartheta)}{\upmu(s,u,\vartheta)}
								\left| 
									\Rad \GeoAng^N v^1
								\right|(s,u,\vartheta)
							\, ds.
	\end{align}
	Next, using the transport equation \eqref{E:UPMUFIRSTTRANSPORT},
	we algebraically decompose the factor outside of the integral 
	in \eqref{E:LASTTIMEINTEGRALTOBOUND} as follows:
	\begin{align} \label{E:MAINDIFFICULTINTEGRALFACTORDECOMPOSITION}
		\frac{1}{\upmu} 
		(\Rad v^1) 
		\sum_{\imath=0}^1 G_{\Lunit \Lunit}^{\imath}
		& = 2 \frac{\Lunit \upmu}{\upmu}
			\\
		& \ \
				-
				\frac{1}{\upmu} G_{\Lunit \Lunit}^0 \Rad (\Densrenormalized  - v^1)
				-
				\frac{1}{\upmu} G_{\Lunit \Lunit}^2 \Rad v^2
				+ 
				\vec{G}_{\Lunit \Lunit} \contr \Lunit \threePsi
				+ 
				2 \vec{G}_{\Lunit \Radunit} \contr \Lunit \threePsi.
				\notag
	\end{align}
	Substituting the decomposition \eqref{E:MAINDIFFICULTINTEGRALFACTORDECOMPOSITION}
	into \eqref{E:LASTTIMEINTEGRALTOBOUND}
	and using the same arguments given in the lines just below 
	\eqref{E:HARDPOINTWISEESTIMATEFORPRODUCTDECOMPOSITIONOFMAINTERM},
	we bound the term \eqref{E:LASTTIMEINTEGRALTOBOUND} by
	\begin{align} \label{E:MAINPARTLASTTIMEINTEGRALTOBOUND}
	\leq
	\displaystyle
	4 (1 + C \varepsilon)
	\left| 
		\frac{[\Lunit \upmu]_-}{\upmu} (\Rad v^1) 
	\right|
	(t,u,\vartheta)
  \int_{s=0}^t 
								\frac{[\Lunit \upmu]_-(s,u,\vartheta)}{\upmu(s,u,\vartheta)}
								\left| 
									\Rad \GeoAng^N v^1
								\right|(s,u,\vartheta)
							\, ds
	\end{align}
	plus a term that is $\leq$ the time-integral-involving product 
	on RHS~\eqref{E:ERRORTERMKEYPOINTWISEESTIMATE}
	featuring the small coefficient $\varepsilon$
	(and thus is of the form $\mbox{\upshape Error}$).
	Finally, we note that RHS~\eqref{E:MAINPARTLASTTIMEINTEGRALTOBOUND}
	is $\leq$ 
	the third term
	$\boxed{4} (1 + C \varepsilon) \cdots$
	on RHS~\eqref{E:KEYPOINTWISEESTIMATE}
	as desired.
	We have thus proved \eqref{E:KEYPOINTWISEESTIMATE} in the case
	$\Fullset_{\ast}^{N;\leq 1} = \GeoAng^N$.

	We now prove \eqref{E:KEYPOINTWISEESTIMATE} in the remaining case
	$\Fullset_{\ast}^{N;\leq 1} = \GeoAng^{N-1} \Rad$.
	The proof is nearly identical to the case
	$\Fullset_{\ast}^{N;\leq 1} = \GeoAng^N$.
	The only difference is the presence of some additional error terms 
	$\mbox{\upshape Error}$, which 
	appear on RHS~\eqref{E:ERRORTERMKEYPOINTWISEESTIMATE}.
	The additional error terms, 
	namely the second term on the first line of RHS~\eqref{E:ERRORTERMKEYPOINTWISEESTIMATE}
	and the double time integral on RHS~\eqref{E:ERRORTERMKEYPOINTWISEESTIMATE},
	are generated in view of equation
	\eqref{E:ADDITIONALTIMEINTEGRALFORFULLYMODIFIEDQUANTITYPOINTWISEESTIMATEONERADIALCASE}
	and the remarks located just above it.

	The proof of \eqref{E:LESSPRECISEKEYPOINTWISEESTIMATE}
	is based on a subset of the above arguments and is much
	simpler. We therefore omit the details,
	noting only that the main simplification is that we do not have
	to rely on the 
	algebraic decompositions \eqref{E:SPLITTINGOFGLLSTARRADGEOANGNVECPSI} and
	\eqref{E:MAINDIFFICULTINTEGRALFACTORDECOMPOSITION};
	we can instead crudely bound the terms on 
	LHS~\eqref{E:SPLITTINGOFGLLSTARRADGEOANGNVECPSI} and \eqref{E:MAINDIFFICULTINTEGRALFACTORDECOMPOSITION}.

\end{proof}

\subsection{Pointwise estimates for the partially modified quantities}
\label{SS:POINTWISEFORPARTIALLYMODIFIED}

In this section, we derive pointwise estimates for the partially modified quantities
from Def.~\ref{D:TRANSPORTRENORMALIZEDTRCHIJUNK}. 
We also derive pointwise estimates for their $\Lunit$ derivative.

\begin{lemma}[\textbf{Pointwise estimates for the partially modified quantities and their}
$\Lunit$ \textbf{derivative}]
	\label{L:SHARPPOINTWISEPARTIALLYMODIFIED}
	Assume that $N=20$
	and let
	$\Fullset_{\ast}^{N-1;\leq 1} \in \lbrace \GeoAng^{N-1}, \GeoAng^{N-2} \Rad \rbrace$.
	Let $\upchipartialmodarg{\Fullset_{\ast}^{N-1;\leq 1}}$
	be the partially modified quantity defined by \eqref{E:TRANSPORTPARTIALRENORMALIZEDTRCHIJUNK}.
	There exist constants\footnote{For the purpose of the remainder of the proof, there is no need to distinguish between the constants $C$ and $C_{\ast}$. Here, we just use $C_{\ast}$ to denote the (large) constants which would in principle have caused the top-order energy to blow up with a worse rate if it were not for the fact that we have carefully distinguished between the energies $\mathbb Q_N$ and $\mathbb Q_N^{(Partial)}$ (and $\mathbb K_N$ and $\mathbb K^{(Partial)}_N$); see Remark~\ref{R:WHICHVARIABLESARECONTROLLED}. Similar remarks apply to later appearances of $C_*$.} $C > 0$ and $C_* > 0$ 
	such that under the data-size and bootstrap assumptions 
	of Sects.~\ref{SS:FLUIDVARIABLEDATAASSUMPTIONS}-\ref{SS:PSIBOOTSTRAP}
	and the smallness assumptions of Sect.~\ref{SS:SMALLNESSASSUMPTIONS}, 
	the following pointwise estimate holds on
	$\mathcal{M}_{\Tboot,U_0}$:
	\begin{subequations}
	\begin{align} 
		\left|
			\Lunit \upchipartialmodarg{\Fullset_{\ast}^{N-1;\leq 1}}
		\right|
		& 
		\leq
		 \frac{1}{2} 
			\left|
				\sum_{\imath=0}^1 G_{\Lunit \Lunit}^{\imath} 
			\right|
			\left|
				\angLap \Fullset_{\ast}^{N-1;\leq 1} v^1
			\right|
			+
			C_*
			\left|
				\angdiff \Fullset_{\ast}^{N;\leq 1} (\Densrenormalized  - v^1)
			\right|
			\label{E:SHARPPOINTWISELUNITPARTIALLYMODIFIED} \\
	& \ \
			+
			C \varepsilon
			\left|
				\Fullset_{\ast}^{[1,N+1];\leq 1}
			\right|
			+
			C 
			\left|
				\myarray[\Fullset_{\ast \ast}^{[1,N];0} \BadVar]
				{\Fullset_{\ast}^{[1,N];\leq 1} \GdVar}
			\right|,
				\notag \\
		\left|
			\upchipartialmodarg{\Fullset_{\ast}^{N-1;\leq 1}}
		\right|
		(t,u,\vartheta)
		& 
		\leq
		\left|
			\upchipartialmodarg{\Fullset_{\ast}^{N-1;\leq 1}}
		\right|
		(0,u,\vartheta)
			\label{E:SHARPPOINTWISEPARTIALLYMODIFIED}
			\\
	& \ \
		+
		\frac{1}{2} 
		\left|
			\sum_{\imath=0}^1 
			G_{\Lunit \Lunit}^{\imath}
		\right|
		(t,u,\vartheta)
			\int_{t'=0}^t
			\left|
				\angLap \Fullset_{\ast}^{N-1;\leq 1} v^1
			\right|
			(t',u,\vartheta)
			\, dt'
			\notag \\
		& \ \
			+
			C_*
			\int_{t'=0}^t
			\left|
				\angdiff \Fullset_{\ast}^{N;\leq 1} (\Densrenormalized  - v^1)
			\right|
			(t',u,\vartheta)
			\, dt
			\notag \\
		& \ \
		+
			C \varepsilon 
			\int_{t'=0}^t
				\left|
					\Fullset_{\ast}^{[1,N+1];\leq 1} \threePsi
				\right|
				(t',u,\vartheta)
			\, dt'
				\notag \\
	& \ \
			+
			C
			\int_{t'=0}^t
				\left|
					\myarray[\Fullset_{\ast \ast}^{[1,N];0} \BadVar]
					{\Fullset_{\ast}^{[1,N];\leq 1} \GdVar}
				\right|
			(t',u,\vartheta)
			\, dt'.
			\notag
	\end{align}
	\end{subequations}
\end{lemma}

\begin{proof}
	See Sect.~\ref{SS:OFTENUSEDESTIMATES} for some comments on the analysis.
	We first prove \eqref{E:SHARPPOINTWISELUNITPARTIALLYMODIFIED} with the help of
	equation \eqref{E:COMMUTEDTRCHIJUNKFIRSTPARTIALRENORMALIZEDTRANSPORTEQUATION}.
	We begin by algebraically decomposing the 
	first product on RHS~\eqref{E:COMMUTEDTRCHIJUNKFIRSTPARTIALRENORMALIZEDTRANSPORTEQUATION}
	as follows:
	\begin{align}\label{E:PARTIALLYMODQUANTPOINTWISEDECOMPOFKEYTERM}
	\frac{1}{2} \vec{G}_{\Lunit \Lunit}\contr\angLap \Fullset_{\ast}^{N-1;\leq 1} \threePsi
	& = \frac{1}{2} \sum_{\imath=0}^1 G_{\Lunit \Lunit}^{\imath}
				\angLap \Fullset_{\ast}^{N-1;\leq 1} v^1
				\\
	& \ \
		+ \frac{1}{2} G_{\Lunit \Lunit}^2 \angLap \Fullset_{\ast}^{N-1;\leq 1} v^2
		+
		\frac{1}{2} G_{\Lunit \Lunit}^0 \angLap \Fullset_{\ast}^{N-1;\leq 1} 
		(\Densrenormalized  - v^1).
		\notag
	\end{align}
	Clearly the first sum on RHS~\eqref{E:SHARPPOINTWISELUNITPARTIALLYMODIFIED}
	arises from the first sum on RHS~\eqref{E:PARTIALLYMODQUANTPOINTWISEDECOMPOFKEYTERM}.
	Next, using \eqref{E:GFRAMESCALARSDEPENDINGONGOODVARIABLES},
	the $L^{\infty}$ estimates of Prop.~\ref{P:IMPROVEMENTOFAUX},
	and \eqref{E:SMALLNESSOFGLL2},
	we find that
	$
	\|G_{\Lunit \Lunit}^0 \|_{L^{\infty}(\Sigma_t^u)}
	\lesssim 1
	$
	and
	$
	\|G_{\Lunit \Lunit}^2 \|_{L^{\infty}(\Sigma_t^u)}
	\lesssim 
	\varepsilon
	$.
	Hence, we can bound the terms on the last line of
	RHS~\eqref{E:PARTIALLYMODQUANTPOINTWISEDECOMPOFKEYTERM}
	by the second and third terms on
	RHS~\eqref{E:SHARPPOINTWISELUNITPARTIALLYMODIFIED}.
	Finally, 
	to bound the terms on RHS~\eqref{E:TRCHIJUNKCOMMUTEDTRANSPORTEQNPARTIALRENORMALIZATIONINHOMOGENEOUSTERM},
	we simply quote \eqref{E:CHIPARTIALMODSOURCETERMPOINTWISE}.
	We have thus proved \eqref{E:SHARPPOINTWISELUNITPARTIALLYMODIFIED}.

	To derive \eqref{E:SHARPPOINTWISEPARTIALLYMODIFIED}, we integrate
	\eqref{E:SHARPPOINTWISELUNITPARTIALLYMODIFIED} along the integral curves of 
	$\Lunit$ as in \eqref{E:SIMPLEFTCID}.
	The only subtle point is that we 
	bound the time integral of the first sum on RHS~\eqref{E:SHARPPOINTWISELUNITPARTIALLYMODIFIED}
	as follows
	by using \eqref{E:RADDERIVATIVESOFGLLDIFFERENCEBOUND} 
	with $M=0$ and $s=t'$:
	\begin{align} \label{E:PULLGLLOUTOFINTEGRAL}
	 	&
		\int_{t'=0}^t
	 		\left\lbrace
	 		\left| \sum_{\imath=0}^1 G_{\Lunit \Lunit}^{\imath} \right|
				\left|
						\angLap \Fullset_{\ast}^{N-1;\leq 1} v^1
				\right|
			\right\rbrace
			(t',u,\vartheta)
		\, dt'
			\\
		& 
		\leq
		\left| \sum_{\imath=0}^1 G_{\Lunit \Lunit}^{\imath} \right|(t,u,\vartheta)
		\int_{t'=0}^t
	 		\left|
				\angLap \Fullset_{\ast}^{N-1;\leq 1} v^1
			\right|
			(t',u,\vartheta)
		\, dt'
			\notag \\
	& \ \
		+ 
		C \varepsilon
		\int_{t'=0}^t
	 		\left|
				\Fullset_{\ast}^{[1,N+1];\leq 1} v^1
			\right|
			(t',u,\vartheta)
		\, dt'.
		\notag
	\end{align}
	Note that the last term on RHS~\eqref{E:PULLGLLOUTOFINTEGRAL}
	is bounded by the next-to-last term on RHS~\eqref{E:SHARPPOINTWISEPARTIALLYMODIFIED}.
	We have thus proved \eqref{E:SHARPPOINTWISEPARTIALLYMODIFIED}.
\end{proof}

\subsection{Pointwise estimates for the inhomogeneous terms in the wave equations}
In this section, we derive
pointwise estimates for the derivatives 
of the inhomogeneous terms in the geometric wave equations
\eqref{E:VELOCITYWAVEEQUATION}-\eqref{E:RENORMALIZEDDENSITYWAVEEQUATION}.

We start with a lemma in which we decompose the derivatives of the 
$\Vortrenormalized$-involving inhomogeneous terms on RHS~\eqref{E:VELOCITYWAVEEQUATION}
into the main terms and error terms. By ``main terms,'' we mean those products
that involve the top-order derivatives of $\Vortrenormalized$.

\begin{lemma}[\textbf{Identification of the important wave equation inhomogeneous terms involving the 
	top-order derivatives of the vorticity}]
	\label{L:WAVEEQUATIONTOPDERIVATIVESOFVORTICITY}
	Assume that $1 \leq N \leq 20$.
	Under the data-size and bootstrap assumptions 
	of Sects.~\ref{SS:FLUIDVARIABLEDATAASSUMPTIONS}-\ref{SS:PSIBOOTSTRAP}
	and the smallness assumptions of Sect.~\ref{SS:SMALLNESSASSUMPTIONS}, 
	the following pointwise estimate holds on
	$\mathcal{M}_{\Tboot,U_0}$
	(see Sect.~\ref{SS:STRINGSOFCOMMUTATIONVECTORFIELDS} regarding the vectorfield operator notation):
	\begin{subequations}
	\begin{align}
		\Tanset^N
		\left\lbrace
			[ia] \upmu (\exp \Densrenormalized) \Speed^2 g_{ab} \Radunit^b \Lunit \Vortrenormalized
		\right\rbrace
		& =
			[ia] \upmu (\exp \Densrenormalized) \Speed^2 g_{ab} \Radunit^b \Tanset^N \Lunit \Vortrenormalized
			+
			\mbox{\upshape Error},
			\label{E:FIRSTTERMWAVEEQUATIONTOPDERIVATIVESOFVORTICITY} 
			\\
		\Tanset^N
			\left\lbrace
				[ia] \upmu (\exp \Densrenormalized) \Speed^2
				\left(
					\frac{g_{ab} \GeoAng^b}{g_{cd} \GeoAng^c \GeoAng^d}
				\right)
				\GeoAng \Vortrenormalized
			\right\rbrace
		& =
			[ia] \upmu (\exp \Densrenormalized) \Speed^2
				\left(
					\frac{g_{ab} \GeoAng^b}{g_{cd} \GeoAng^c \GeoAng^d}
				\right)
				\Tanset^N \GeoAng \Vortrenormalized
					\label{E:SECONDTERMWAVEEQUATIONTOPDERIVATIVESOFVORTICITY} \\
		& \ \
				+
				\mbox{\upshape Error},
				\notag
	\end{align}
	\end{subequations}
	where
	\begin{align} \label{E:ERRORTERMBOUNDWAVEEQUATIONTOPDERIVATIVESOFVORTICITY}
		\left|
			\mbox{\upshape Error}
		\right|
		& \lesssim 
			\varepsilon
			\left|
				\Fullset_{\ast \ast}^{[1,N];\leq 1} \BadVar
			\right|
			+
			\left|
				\Tanset^{\leq N} \Vortrenormalized
			\right|.
	\end{align}

	Moreover, let $\Fullset_{\ast}^{N;1}$ be a $N^{th}$ order
	vectorfield operator containing exactly one factor of $\Rad$,
	and let $\Tanset^{N-1}$ denote the remaining non-$\Rad$ factors.
	Then we have the following estimates:
	\begin{subequations}
	\begin{align}
		\Fullset_{\ast}^{N;1}
		\left\lbrace
			[ia] \upmu (\exp \Densrenormalized) \Speed^2 g_{ab} \Radunit^b \Lunit \Vortrenormalized
		\right\rbrace
		& =
			-
			[ia] \upmu^2 
			(\exp \Densrenormalized) \Speed^2 g_{ab} \Radunit^b \Tanset^{N-1} \Lunit \Lunit \Vortrenormalized
			+
			\mbox{\upshape Error},
			\label{E:FIRSTTERMWAVEEQUATIONTOPRADCONTAININGDERIVATIVESOFVORTICITY} 
			\\
			\Fullset_{\ast}^{N;1}
			\left\lbrace
				[ia] \upmu (\exp \Densrenormalized) \Speed^2
				\left(
					\frac{g_{ab} \GeoAng^b}{g_{cd} \GeoAng^c \GeoAng^d}
				\right)
				\GeoAng \Vortrenormalized
			\right\rbrace
		& =
			-
			[ia] \upmu^2 
				(\exp \Densrenormalized) \Speed^2
				\left(
					\frac{g_{ab} \GeoAng^b}{g_{cd} \GeoAng^c \GeoAng^d}
				\right)
				\Tanset^{N-1} \GeoAng \Lunit \Vortrenormalized
					\label{E:SECONDTERMWAVEEQUATIONTOPRADCONTAININGDERIVATIVESOFVORTICITY} \\
		& \ \
				+
				\mbox{\upshape Error},
				\notag
	\end{align}
	\end{subequations}
	where $\mbox{\upshape Error}$ satisfies \eqref{E:ERRORTERMBOUNDWAVEEQUATIONTOPDERIVATIVESOFVORTICITY}.
\end{lemma}

\begin{proof}
	See Sect.~\ref{SS:OFTENUSEDESTIMATES} for some comments on the analysis.
	The estimates 
	\eqref{E:FIRSTTERMWAVEEQUATIONTOPDERIVATIVESOFVORTICITY}-\eqref{E:ERRORTERMBOUNDWAVEEQUATIONTOPDERIVATIVESOFVORTICITY}
	are a straightforward consequence
	of Lemma~\ref{L:SCHEMATICDEPENDENCEOFMANYTENSORFIELDS},
	(which implies that
	$
	\upmu (\exp \Densrenormalized) \Speed^2 g_{ab} \Radunit^b
	= \smoothfunction(\BadVar)
	$
	and
	$\upmu (\exp \Densrenormalized) \Speed^2
				\left(
					\frac{g_{ab} \GeoAng^b}{g_{cd} \GeoAng^c \GeoAng^d}
				\right)
	= \smoothfunction(\BadVar)$)
	and the $L^{\infty}$ estimates of Prop.~\ref{P:IMPROVEMENTOFAUX}.

	Similar remarks apply to the estimates
	\eqref{E:FIRSTTERMWAVEEQUATIONTOPRADCONTAININGDERIVATIVESOFVORTICITY}-\eqref{E:SECONDTERMWAVEEQUATIONTOPRADCONTAININGDERIVATIVESOFVORTICITY}.
However, when bounding derivatives of 
$\Vortrenormalized$ that contain the factor of $\Rad$,
we first use the commutator estimate
\eqref{E:LESSPRECISEPERMUTEDVECTORFIELDSACTINGONFUNCTIONSCOMMUTATORESTIMATE} with 
$f = \Vortrenormalized$
and the $L^{\infty}$ estimates of Prop.~\ref{P:IMPROVEMENTOFAUX}
to commute the factor of $\Rad$ so that it hits $\Vortrenormalized$ first;
the commutator terms are of the form $\mbox{\upshape Error}$,
where $\mbox{\upshape Error}$ verifies \eqref{E:ERRORTERMBOUNDWAVEEQUATIONTOPDERIVATIVESOFVORTICITY}.
We then use the transport equation \eqref{E:RENORMALIZEDVORTICTITYTRANSPORTEQUATION}
and the identity \eqref{E:TRANSPORTVECTORFIELDINTERMSOFLUNITANDRADUNIT}
to algebraically replace $\Rad \Vortrenormalized$ with $- \upmu \Lunit \Vortrenormalized$;
this replacement is the origin of the extra factor of $\upmu$
on RHSs~\eqref{E:FIRSTTERMWAVEEQUATIONTOPRADCONTAININGDERIVATIVESOFVORTICITY}-\eqref{E:SECONDTERMWAVEEQUATIONTOPRADCONTAININGDERIVATIVESOFVORTICITY}
compared to \eqref{E:FIRSTTERMWAVEEQUATIONTOPDERIVATIVESOFVORTICITY}-\eqref{E:SECONDTERMWAVEEQUATIONTOPDERIVATIVESOFVORTICITY}. 
Finally, we explicitly place the products containing the top-order (that is, order $N+1$)
derivatives of $\Vortrenormalized$ on
RHSs~\eqref{E:FIRSTTERMWAVEEQUATIONTOPRADCONTAININGDERIVATIVESOFVORTICITY}-\eqref{E:SECONDTERMWAVEEQUATIONTOPRADCONTAININGDERIVATIVESOFVORTICITY}
and again use the $L^{\infty}$ estimates of Prop.~\ref{P:IMPROVEMENTOFAUX}
to conclude that the remaining products are of the form $\mbox{\upshape Error}$,
where $\mbox{\upshape Error}$ verifies \eqref{E:ERRORTERMBOUNDWAVEEQUATIONTOPDERIVATIVESOFVORTICITY}.
This completes the proof of the lemma.
\end{proof}

We now derive estimates for
the derivatives of the null forms
on RHSs~\eqref{E:VELOCITYWAVEEQUATION}-\eqref{E:RENORMALIZEDDENSITYWAVEEQUATION}.

\begin{lemma}[\textbf{Estimates for the null forms}]
	\label{L:POINTWISESTIMATESFORTHENULLFOMRS}
	Assume that $1 \leq N \leq 20$ and 
	let $\mathscr{Q}^i$ and $\mathscr{Q}$
	be the null forms defined by 
	\eqref{E:VELOCITYNULLFORM} and \eqref{E:DENSITYNULLFORM}.
	Under the data-size and bootstrap assumptions 
	of Sects.~\ref{SS:FLUIDVARIABLEDATAASSUMPTIONS}-\ref{SS:PSIBOOTSTRAP}
	and the smallness assumptions of Sect.~\ref{SS:SMALLNESSASSUMPTIONS}, 
	the following estimates hold on $\mathcal{M}_{\Tboot,U_0}$
	(see Sect.~\ref{SS:STRINGSOFCOMMUTATIONVECTORFIELDS} regarding the vectorfield operator notation):
\begin{align} \label{E:POINTWISESTIMATESFORTHENULLFOMRS}
		\left|
			\Fullset^{N;\leq 1} (\upmu \mathscr{Q}^i)
		\right|,
			\,
		\left|
			\Fullset^{N;\leq 1} (\upmu \mathscr{Q})
		\right|
		& \lesssim 
			\left|
				\Fullset_{\ast}^{[1,N+1];\leq 2} \threePsi
			\right|
			+
			\varepsilon
			\left|
				\Fullset_{\ast \ast}^{[1,N];\leq 1} \BadVar
			\right|.
	\end{align}
\end{lemma}
\begin{proof}
	See Sect.~\ref{SS:OFTENUSEDESTIMATES} for some comments on the analysis.
	The estimate \eqref{E:POINTWISESTIMATESFORTHENULLFOMRS} 
	is a straightforward consequence of \eqref{E:UPMUTIMESNULLFORMSSCHEMATIC}
	and the $L^{\infty}$ estimates of Prop.~\ref{P:IMPROVEMENTOFAUX}.
\end{proof}

\subsection{Pointwise estimates for the error terms generated by the multiplier vectorfield}
\label{SS:POINTWISEESTIMATESFORTHEMULTIPLIERVECTORFIELDERRORTERMS}
In this section, we derive pointwise
estimates for the wave equation energy estimate error terms 
generated by the deformation tensor of the multiplier vectorfield $\Mult$.
That is, we obtain pointwise bounds for the terms $\basicenergyerrorarg{\Mult}{i}[\Psi]$
(see \eqref{E:MULTERRORINT})
corresponding the integrand $\upmu \enmomtensor^{\alpha \beta}[\Psi] \deformarg{\Mult}{\alpha}{\beta}$
on RHS~\eqref{E:E0DIVID}.

\begin{lemma}[\textbf{Pointwise bounds for the error terms generated by the deformation tensor of} 
$\Mult$]
\label{L:MULTIPLIERVECTORFIEDERRORTERMPOINTWISEBOUND}
	Let $\Psi$ be a function\footnote{We will eventually apply this estimate
	with the role of $\Psi$ played by
	a derivative of an element of $\lbrace \Densrenormalized - v^1$, $v^1$, $v^2 \rbrace$.}
	and consider the multiplier vectorfield error terms 
	$\basicenergyerrorarg{\Mult}{1}[\Psi],\cdots,\basicenergyerrorarg{\Mult}{5}[\Psi]$
	defined in \eqref{E:MULTERRORINTEG1}-\eqref{E:MULTERRORINTEG5}.
	Let $\varsigma > 0$ be a real number.
	Under the data-size and bootstrap assumptions 
	of Sects.~\ref{SS:FLUIDVARIABLEDATAASSUMPTIONS}-\ref{SS:PSIBOOTSTRAP}
	and the smallness assumptions of Sect.~\ref{SS:SMALLNESSASSUMPTIONS}, 
	the following pointwise inequality holds
	on $\mathcal{M}_{\Tboot,U_0}$ (without any absolute value taken on the left),
	where the implicit constants are independent of $\varsigma$:
	\begin{align} \label{E:MULTIPLIERVECTORFIEDERRORTERMPOINTWISEBOUND}
		\sum_{i=1}^5 \basicenergyerrorarg{\Mult}{i}[\Psi]
		& \lesssim
			(1 + \varsigma^{-1})(\Lunit \Psi)^2
			+ (1 + \varsigma^{-1}) (\Rad \Psi)^2
			+ \upmu |\angdiff \Psi|^2
			+ \varsigma \TranminusdatasizeWithFactor |\angdiff \Psi|^2
				\\
		& \ \
				+ \frac{1}{\sqrt{\Tboot - t}} 
				\upmu |\angdiff \Psi|^2.
				\notag
\end{align}
\end{lemma}

\begin{proof}
	See Sect.~\ref{SS:OFTENUSEDESTIMATES} for some comments on the analysis.
	Only the term $\basicenergyerrorarg{\Mult}{3}[\Psi]$
	is difficult to treat.
	Specifically, using the 
	schematic relations
	\eqref{E:TENSORSDEPENDINGONGOODVARIABLESGOODPSIDERIVATIVES},
	\eqref{E:TENSORSDEPENDINGONGOODVARIABLESBADDERIVATIVES},
	and \eqref{E:ZETABADDERIVATIVESBETTERTHANEXPECTED},
	the estimate \eqref{E:ANGDIFFXI},
	and the $L^{\infty}$ estimates
	of Props.~\ref{P:IMPROVEMENTOFAUX} and \ref{P:IMPROVEMENTOFHIGHERTRANSVERSALBOOTSTRAP},
	it is straightforward to verify that the terms in braces on
	RHS \eqref{E:MULTERRORINTEG1}, \eqref{E:MULTERRORINTEG2}, \eqref{E:MULTERRORINTEG4}, and
	\eqref{E:MULTERRORINTEG5} are bounded in magnitude by $\lesssim 1$.
	It follows that for $i=1,2,4,5$,
	$\left|
		\basicenergyerrorarg{\Mult}{i}[\Psi]
	\right|
	$ 
	is $\lesssim$ the sum of the terms on the first line of
	RHS~\eqref{E:MULTIPLIERVECTORFIEDERRORTERMPOINTWISEBOUND}.
	The quantities $\varsigma$ and $\TranminusdatasizeWithFactor$
	appear on RHS~\eqref{E:MULTIPLIERVECTORFIEDERRORTERMPOINTWISEBOUND} 
	because we use Young's inequality to bound
	$\basicenergyerrorarg{\Mult}{4}[\Psi] \lesssim |\Lunit \Psi||\angdiff \Psi| 
	\leq \varsigma^{-1} \TranminusdatasizeWithFactor^{-1}(\Lunit \Psi)^2 
	+ \varsigma \TranminusdatasizeWithFactor |\angdiff \Psi|^2
	\leq C \varsigma^{-1} (\Lunit \Psi)^2 + C \varsigma \TranminusdatasizeWithFactor |\angdiff \Psi|^2
	$.
	Similar remarks apply to $\basicenergyerrorarg{\Mult}{5}[\Psi]$.

	To bound the difficult term 
	$\basicenergyerrorarg{\Mult}{3}[\Psi]$,
	we also use the estimates
	\eqref{E:POSITIVEPARTOFLMUOVERMUISBOUNDED}
	and 
	\eqref{E:UNIFORMBOUNDFORMRADMUOVERMU},
	which allow us to bound the
	first two terms in braces on RHS~\eqref{E:MULTERRORINTEG3}.
	Note that since no absolute value 
	is taken on LHS~\eqref{E:MULTIPLIERVECTORFIEDERRORTERMPOINTWISEBOUND},
	we may replace the factor 
	$(\Rad \upmu)/\upmu$ 
	from RHS~\eqref{E:MULTERRORINTEG3}
	with the factor $[\Rad \upmu]_+/\upmu$,
	which is bounded by \eqref{E:UNIFORMBOUNDFORMRADMUOVERMU}.
	This completes our proof of \eqref{E:MULTIPLIERVECTORFIEDERRORTERMPOINTWISEBOUND}.
\end{proof}

\subsection{Proof of Prop.~\ref{P:WAVEIDOFKEYDIFFICULTENREGYERRORTERMS}}
\label{SS:PROOFOFPROPWAVEIDOFKEYDIFFICULTENREGYERRORTERMS}
See Sect.~\ref{SS:OFTENUSEDESTIMATES} for some comments on the analysis.
We must derive estimates for the elements
$\Psi \in \lbrace \Densrenormalized  - v^1,v^1,v^2 \rbrace$.
To condense the notation, we
use the following notation for the term in braces on RHS \eqref{E:BOXZCOM}:
\begin{align} \label{E:COMMCURRENTZ}
{\Jcurrent{Z}^{\alpha}[\Psi]}
			:= \deformuparg{Z}{\alpha}{\beta} \D_{\beta} \Psi 
			- \frac{1}{2} \myspacetimetr \deform{Z} \D^{\alpha} \Psi.
\end{align}
Throughout we silently use the
Definition \ref{D:HARMLESSTERMS} of 
$Harmless_{(Wave)}^{\leq N}$ terms.
We prove the estimates 
\eqref{E:SECONDWAVEGEOANGANGISTHEFIRSTCOMMUTATORIMPORTANTTERMS}
and \eqref{E:WAVETHIRDHARMLESSORDERNCOMMUTATORS} 
(corresponding to $\Psi = v^i$),
whose proofs are closely related,
in detail. Later in the proof, we indicate the minor changes needed to obtain
\eqref{E:WAVELISTHEFIRSTCOMMUTATORIMPORTANTTERMS}-\eqref{E:SECONDWAVELISTHEFIRSTCOMMUTATORIMPORTANTTERMS}
and
\eqref{E:WAVEHARMLESSORDERNPURETANGENTIALCOMMUTATORS}-\eqref{E:WAVESECONDHARMLESSORDERNCOMMUTATORS}.
At the very end of the proof, we indicate
the minor changes needed to obtain \eqref{E:KEYDIFFICULTFORRENORMALZIEDDENSITY},
that is, the estimates in the case $\Psi = \Densrenormalized  - v^1$.
%
To proceed, we iterate \eqref{E:BOXZCOM}, 
use the wave equation \eqref{E:VELOCITYWAVEEQUATION},
use the decomposition \eqref{E:VELOCITYWAVEEQUATIONDERIVATIVEOFVORTICITYINHOMOGENEOUSTERMEXPRESSION},
use the estimates 
\eqref{E:FIRSTTERMWAVEEQUATIONTOPRADCONTAININGDERIVATIVESOFVORTICITY}-\eqref{E:SECONDTERMWAVEEQUATIONTOPRADCONTAININGDERIVATIVESOFVORTICITY}
and \eqref{E:POINTWISESTIMATESFORTHENULLFOMRS},
and use the estimates
\begin{align} \label{E:SIMPLEDEFTENSORPOINTWISEBOUNDS}
\left\|
	\Fullset^{\leq 10;\leq 1} \mytr \angdeform{\Lunit}
\right\|_{L^{\infty}(\Sigma_t^u)},
\,
\left\|
	\Fullset^{\leq 10;\leq 1} \mytr \angdeform{\GeoAng}
\right\|_{L^{\infty}(\Sigma_t^u)},
	\,
\left\|
	\Tanset^{\leq 10} \mytr \angdeform{\Rad}
\right\|_{L^{\infty}(\Sigma_t^u)}
\lesssim 1
\end{align}
(which follow from \eqref{E:CONNECTIONBETWEENANGLIEOFGSPHEREANDDEFORMATIONTENSORS}, 
\eqref{E:POINTWISEESTIMATESFORGSPHEREANDITSNOSPECIALSTRUCTUREDERIVATIVES},
and the $L^{\infty}$ estimates of Prop.~\ref{P:IMPROVEMENTOFAUX})
to deduce that
\begin{align} \label{E:COMMUTEDWAVEGEOANGFIRSTMAINTERMPLUSERRORTERM}
		\upmu \square_{g(\threePsi)} (\Fullset^{N-1;1} \GeoAng v^i)
		&  =  \Fullset^{N-1;1}
					\left(
					\upmu 
					\D_{\alpha} 
					\Jcurrent{\GeoAng}^{\alpha}[v^i]
					\right)
						\\
		& \ \
			- [ia] \upmu^2 (\exp \Densrenormalized) \Speed^2 (g_{ab} \Radunit^b) 
				\Tanset^{N-1} \Lunit \Lunit \Vortrenormalized
				\notag \\
		& \ \
			+ [ia] \upmu^2 (\exp \Densrenormalized) \Speed^2
			\left(
				\frac{g_{ab} \GeoAng^b}{g_{cd} \GeoAng^c \GeoAng^d}
			\right)
			\Tanset^{N-1} \GeoAng \Lunit \Vortrenormalized
				\notag \\
		& \ \
			+ \mbox{Error},
			\notag
\end{align}
where
\begin{align} \label{E:GEOANGLOWERORDERRORTERMESTIMATE}
	\left|
		\mbox{Error}
	\right|
	& \lesssim
		\mathop{\sum_{N_1 + N_2 + N_3 \leq N-1}}_{N_2 \leq N-2}
		\sum_{M_1 + M_2 + M_3 \leq 1}
		\sum_{\Singletan_1, \Singletan_2 \in \Tanset}
		\left(
			1 +
			\left|
				\Fullset^{N_1;M_1} \mytr \angdeform{\Singletan_1}
			\right|
		\right)
		\left|
			\Fullset^{N_2;M_2}
			\left(
				\upmu 
				\D_{\alpha} 
				\Jcurrent{\Singletan_2}^{\alpha}[\Fullset^{N_3;M_3} v^i]
			\right)
		\right|
			\\
	& \ \
		+
		\mathop{\sum_{N_1 + N_2 + N_3 \leq N-1}}_{N_2 \leq N-2}
		\sum_{\Singletan \in \Tanset}
		\left(
			1 +
			\left|
				\Tanset^{N_1} \mytr \angdeform{\Singletan}
			\right|
		\right)
		\left|
			\Tanset^{N_2}
			\left(
				\upmu 
				\D_{\alpha} 
				\Jcurrent{\Rad}^{\alpha}[\Tanset^{N_3} v^i]
			\right)
		\right|
			\notag \\
	& \ \
		+
		\mathop{\sum_{N_1 + N_2 + N_3 \leq N-1}}_{N_2 \leq N-2}
		\sum_{\Singletan \in \Tanset}
		\left(
			1 +
			\left|
				\Tanset^{N_1} \mytr \angdeform{\Rad}
			\right|
		\right)
		\left|
			\Tanset^{N_2}
			\left(
				\upmu 
				\D_{\alpha} 
				\Jcurrent{\Singletan}^{\alpha}[\Tanset^{N_3} v^i]
			\right)
		\right|
			\notag \\
	& \ \
			+
			\left|
				\Fullset_{\ast}^{[1,N+1];\leq 2} \threePsi
			\right|
			+
			\myarray[
				\left|
				\Fullset_{\ast \ast}^{[1,N];\leq 1} \BadVar
			\right|]
			{
			\left|
				\Fullset_{\ast}^{[1,N];\leq 2} \GdVar
			\right|
			}
			+
			\left|
				\Tanset^{\leq N} \Vortrenormalized
			\right|.
			\notag
\end{align}
Note that the terms on the last line of
RHS~\eqref{E:COMMUTEDWAVEGEOANGFIRSTMAINTERMPLUSERRORTERM}
are $Harmless_{(Wave)}^{\leq N}$ as desired.

\begin{remark}
	For the purpose of proving \eqref{E:SECONDWAVEGEOANGANGISTHEFIRSTCOMMUTATORIMPORTANTTERMS}
and \eqref{E:WAVETHIRDHARMLESSORDERNCOMMUTATORS},
	the estimate \eqref{E:GEOANGLOWERORDERRORTERMESTIMATE}
	is non-optimal in the sense that some terms on RHS~\eqref{E:GEOANGLOWERORDERRORTERMESTIMATE}
	could be deleted and the inequality would remain true. 
	However, those terms later appear
	when we are deriving the other estimates of the proposition.
	For this reason, we find it convenient to already include them 
	on RHS~\eqref{E:GEOANGLOWERORDERRORTERMESTIMATE}.
\end{remark}

Most of our effort goes towards estimating 
the first term on RHS~\eqref{E:COMMUTEDWAVEGEOANGFIRSTMAINTERMPLUSERRORTERM}.
Equivalently, we may analyze 
the $\Fullset^{N-1;1}$ derivatives of the seven terms
on RHS~\eqref{E:DIVCOMMUTATIONCURRENTDECOMPOSITION}
(with $\Psi = v^i$ and $Z = \GeoAng$ in \eqref{E:DIVCOMMUTATIONCURRENTDECOMPOSITION}).
We will show that if $\Fullset^{N-1;1}$ contains no factor of $\Lunit$, then
\begin{align}
\Fullset^{N-1;1} \mathscr{K}_{(\pi-Danger)}^{(\GeoAng)}[v^i]
& = (\GeoAng^{N-1} \Rad \mytr \upchi) \Rad v^i
	+
	Harmless_{(Wave)}^{\leq N},
	\label{E:PIDANGERANALYSIS} \\
\Fullset^{N-1;1} \mathscr{K}_{(\pi-Less \ Dangerous)}^{(\GeoAng)}[v^i]
& = 
	\upmu \GeoAngFlatRadComponent (\angdiffuparg{\#} \GeoAng^{N-2} \Rad \mytr \upchi)
	\cdot 
	\angdiff v^i
	+ 
	Harmless_{(Wave)}^{\leq N},
	\label{E:PILESSDANGEROUSANALYSIS}
\end{align}
and
\begin{align} \label{E:PIHARMLESSANALYSIS}
&
\Fullset^{N-1;1} \mathscr{K}_{(\pi-Cancel-1)}^{(\GeoAng)}[v^i],
	\,
\Fullset^{N-1;1} \mathscr{K}_{(\pi-Cancel-2)}^{(\GeoAng)}[v^i],
	\,
\Fullset^{N-1;1} \mathscr{K}_{(\pi-Good)}^{(\GeoAng)}[v^i],
	\\
& \Fullset^{N-1;1} \mathscr{K}_{(\Psi)}^{(\GeoAng)}[v^i],
		\,
\Fullset^{N-1;1} \mathscr{K}_{(Low)}^{(\GeoAng)}[v^i]
	\notag
		\\
& = Harmless_{(Wave)}^{\leq N}.
\notag
\end{align}
At the same time,
we will show that if $\Fullset^{N-1;1}$ contains one or more factors of $\Lunit$
(and thus $\Fullset^{N-1;1} = \Fullset_{\ast}^{N-1;1}$), then
\begin{align} \label{E:LCONTAININGPIHARMLESSANALYSIS}
&
\Fullset^{N-1;1} \mathscr{K}_{(\pi-Danger)}^{(\GeoAng)}[v^i],
	\,
\Fullset^{N-1;1} \mathscr{K}_{(\pi-Less \ Dangerous)}^{(\GeoAng)}[v^i],
	\\
&
\Fullset^{N-1;1} \mathscr{K}_{(\pi-Cancel-1)}^{(\GeoAng)}[v^i],
	\,
\Fullset^{N-1;1} \mathscr{K}_{(\pi-Cancel-2)}^{(\GeoAng)}[v^i],
	\,
\Fullset^{N-1;1} \mathscr{K}_{(\pi-Good)}^{(\GeoAng)}[v^i],
\notag	\\
& \Fullset^{N-1;1} \mathscr{K}_{(\Psi)}^{(\GeoAng)}[v^i],
		\,
\Fullset^{N-1;1} \mathscr{K}_{(Low)}^{(\GeoAng)}[v^i]
	\notag
		\\
& = Harmless_{(Wave)}^{\leq N}.
\notag
\end{align}

After establishing \eqref{E:PIHARMLESSANALYSIS}-\eqref{E:LCONTAININGPIHARMLESSANALYSIS},
we will show that
\begin{align} \label{E:WAVEERRORTERMSAREHARMLESS}
	\mbox{RHS~\eqref{E:GEOANGLOWERORDERRORTERMESTIMATE}}
	= Harmless_{(Wave)}^{\leq N}.
\end{align}

Then combining \eqref{E:COMMUTEDWAVEGEOANGFIRSTMAINTERMPLUSERRORTERM},
\eqref{E:PIHARMLESSANALYSIS},
\eqref{E:LCONTAININGPIHARMLESSANALYSIS},
and \eqref{E:WAVEERRORTERMSAREHARMLESS},
we conclude the desired estimates
\eqref{E:SECONDWAVEGEOANGANGISTHEFIRSTCOMMUTATORIMPORTANTTERMS} and 
\eqref{E:WAVETHIRDHARMLESSORDERNCOMMUTATORS}.

We now return to our analysis of the 
first term on RHS~\eqref{E:COMMUTEDWAVEGEOANGFIRSTMAINTERMPLUSERRORTERM}.
We will separately analyze the $\Fullset^{N-1;1}$ derivative
of each of the 
seven terms on RHS~\eqref{E:DIVCOMMUTATIONCURRENTDECOMPOSITION}
(with $Z = \GeoAng$ and $\Psi = v^i$ in \eqref{E:DIVCOMMUTATIONCURRENTDECOMPOSITION}).

\medskip
\noindent
\underline{\textbf{Analysis of} $\Fullset^{N-1;1} \mathscr{K}_{(\pi-Danger)}^{(\GeoAng)}[v^1]$.}
We apply to $\Fullset^{N-1;1}$ to \eqref{E:DIVCURRENTTRANSVERSAL} 
(with $\Psi = v^i$ and $Z = \GeoAng$).
We first analyze the difficult product in which all derivatives fall on 
the factor $\angdiv \angdeformoneformupsharparg{\GeoAng}{\Lunit}$:
\begin{align} \label{E:PIDANGERDIFFICULT}
	- (\Fullset^{N-1;1} \angdiv \angdeformoneformupsharparg{\GeoAng}{\Lunit}) \Rad v^i.
\end{align}
Using the estimate \eqref{E:GEOANGDEFIMPORTANTANGDIVSPHERELANDANGDIFFPILRADTERMS} 
for the first term on the LHS
and the simple bound
$\| \Rad v^i \|_{L^{\infty}(\Sigma_t^u)} \lesssim 1$
(see \eqref{E:PSITRANSVERSALLINFINITYBOUNDBOOTSTRAPIMPROVEDSMALL} and
\eqref{E:PSITRANSVERSALLINFINITYBOUNDBOOTSTRAPIMPROVEDLARGE}),
we deduce from \eqref{E:PIDANGERDIFFICULT} that
\begin{align} \label{E:PIDANGERDIFFICULTDECOMPOSED}
	- (\Fullset^{N-1;1} \angdiv \angdeformoneformupsharparg{\GeoAng}{\Lunit}) \Rad v^i
	= (\GeoAng \Fullset^{N-1;1} \mytr \upchi) \Rad v^i
		+ Harmless_{(Wave)}^{\leq N}.
\end{align}
We first consider the case in which 
$\Fullset^{N-1;1}$ contains no factor of $\Lunit$,
which is relevant for proving \eqref{E:PIDANGERANALYSIS}.
Then $\Fullset^{N-1;1}$ contains $N-2$ factors of $\GeoAng$
and one factor of $\Rad$.
Thus, using using the commutator estimate 
\eqref{E:LESSPRECISEPERMUTEDVECTORFIELDSACTINGONFUNCTIONSCOMMUTATORESTIMATE} with 
$f = \mytr \upchi$, $N$ in the role of $N+1$, and $M=1$,
the estimate \eqref{E:POINTWISEESTIMATESFORCHIANDITSDERIVATIVES},
and the $L^{\infty}$ estimates of Prop.~\ref{P:IMPROVEMENTOFAUX},
we may commute the factor of $\Rad$ so that it hits $\mytr \upchi$ first,
thereby obtaining
$(\GeoAng \Fullset^{N-1;1} \mytr \upchi) 
\Rad v^i
=
(\GeoAng^{N-1} \Rad \mytr \upchi)
\Rad v^i
+ Harmless_{(Wave)}^{\leq N}
$.
The remaining terms obtained from
applying $\Fullset^{N-1;1}$ to \eqref{E:DIVCURRENTTRANSVERSAL}
generate products involving $\leq N-2$
derivatives of 
$\angdiv \angdeformoneformupsharparg{\GeoAng}{\Lunit}$.
We will show that these products are
$Harmless_{(Wave)}^{\leq N}$,
which completes the proof of \eqref{E:PIDANGERANALYSIS}.
To proceed, we again use the $L^{\infty}$ estimates of Prop.~\ref{P:IMPROVEMENTOFAUX},
the estimate \eqref{E:POINTWISEESTIMATESFORCHIANDITSDERIVATIVES},
and \eqref{E:GEOANGDEFIMPORTANTANGDIVSPHERELANDANGDIFFPILRADTERMS}
(with $\leq N-2$ in the role of $N-1$ in \eqref{E:GEOANGDEFIMPORTANTANGDIVSPHERELANDANGDIFFPILRADTERMS})
to deduce that all of the products under consideration
are $Harmless_{(Wave)}^{\leq N}$.
We clarify that
the estimate \eqref{E:GEOANGDEFIMPORTANTANGDIVSPHERELANDANGDIFFPILRADTERMS}
(for the first term on the LHS) generates a factor of 
$\mytr \upchi$ 
with $\leq N-1$ derivatives on it 
(located on LHS~\eqref{E:GEOANGDEFIMPORTANTANGDIVSPHERELANDANGDIFFPILRADTERMS}),
which is in contrast to the factor from \eqref{E:PIDANGERDIFFICULTDECOMPOSED} with $N$ derivatives.
This factor is below-top-order in the sense that we may bound it with
\eqref{E:POINTWISEESTIMATESFORCHIANDITSDERIVATIVES}
and hence the corresponding product of this factor and $\Rad v^i$
contributes only to the $Harmless_{(Wave)}^{\leq N}$ terms.
We have thus proved \eqref{E:PIDANGERANALYSIS}.

We now consider the case in which 
$\Fullset^{N-1;1}$ contains a factor of $\Lunit$,
which is relevant for the estimate \eqref{E:WAVETHIRDHARMLESSORDERNCOMMUTATORS}
Noting that the estimate \eqref{E:PIDANGERDIFFICULTDECOMPOSED} still holds,
we use the same commutator argument given in the previous paragraph
to obtain
$(\GeoAng \Fullset^{N-1;1} \mytr \upchi) 
\Rad v^i
=
(\Fullset^{N-1;1} \Lunit \mytr \upchi)
\Rad v^i
+ Harmless_{(Wave)}^{\leq N}
$, 
where the operators $\Fullset^{N-1;1}$ on the LHS and RHS are not necessarily the same.
Using \eqref{E:LUNITCOMMUTEDLUNITSMALLIPOINTWISE} with $M=1$
and the bound $\| \Rad v^i \|_{L^{\infty}(\Sigma_t^u)} \lesssim 1$ mentioned above,
we deduce that
$(\Fullset^{N-1;1} \Lunit \mytr \upchi)
	= Harmless_{(Wave)}^{\leq N}
$.
The remaining terms obtained from
applying $\Fullset^{N-1;1}$ to \eqref{E:DIVCURRENTTRANSVERSAL}
are $Harmless_{(Wave)}^{\leq N}$ for the same reasons 
given in the previous paragraph.
We have thus proved \eqref{E:LCONTAININGPIHARMLESSANALYSIS} for
the term
$\Fullset^{N-1;1} \mathscr{K}_{(\pi-Danger)}^{(\GeoAng)}[v^i]$.

\medskip
\noindent
\underline{\textbf{Analysis of} $\Fullset^{N-1;1} \mathscr{K}_{(\pi-Cancel-1)}^{(\GeoAng)}[v^i]$.}
We apply $\Fullset^{N-1;1}$ to \eqref{E:DIVCURRENTCANEL1} 
(with $\Psi = v^i$ and $Z = \GeoAng$).
We first analyze the difficult product in which all derivatives fall on 
the deformation tensor components:
\begin{align} \label{E:PICANCEL1DIFFICULT}
	\left\lbrace
		\frac{1}{2} \Fullset^{N-1;1} \Rad \mytr \angdeform{\GeoAng}
		- \Fullset^{N-1;1} \angdiv \angdeformoneformupsharparg{\GeoAng}{\Rad}
		- \upmu \Fullset^{N-1;1} \angdiv \angdeformoneformupsharparg{\GeoAng}{\Lunit}
	\right\rbrace
	\Lunit v^i.
\end{align}
Using the second estimate in
\eqref{E:GEOANGDEFIMPORTANTLIERADSPHERELANDRADTRACESPHERETERMS}
and the first and second estimates in
\eqref{E:GEOANGDEFIMPORTANTANGDIVSPHERELANDANGDIFFPILRADTERMS},
we express the terms in brace in \eqref{E:PICANCEL1DIFFICULT}
as the sum of $Harmless_{(Wave)}^{\leq N}$
terms and terms involving the order $N$ derivatives of
$\upmu$ and $\mytr \upchi$, which \emph{exactly cancel}.
Also using the bound
$\| \Lunit v^i \|_{L^{\infty}(\Sigma_t^u)} \lesssim \varepsilon$ 
(see \eqref{E:PSIMIXEDUPTOORDERTWOTRANSVERSALIMPROVED}),
we conclude that
$\mbox{\eqref{E:PICANCEL1DIFFICULT}}= Harmless_{(Wave)}^{\leq N}$.
The remaining terms obtained from
applying $\Fullset^{N-1;1}$ 
to \eqref{E:DIVCURRENTCANEL1}
can be shown to be 
$Harmless_{(Wave)}^{\leq N}$
by combining essentially the same argument
with the schematic identity \eqref{E:LINEARLYSMALLSCALARSDEPENDINGONGOODVARIABLES}
for $\GeoAngFlatRadComponent$,
the estimates \eqref{E:POINTWISEESTIMATESFORCHIANDITSDERIVATIVES}
and
\eqref{E:GEOANGPOINTWISE},
and the $L^{\infty}$ estimates of Prop.~\ref{P:IMPROVEMENTOFAUX}.
We have thus proved \eqref{E:PIHARMLESSANALYSIS} 
and
\eqref{E:LCONTAININGPIHARMLESSANALYSIS}
for
$\Fullset^{N-1;1} \mathscr{K}_{(\pi-Cancel-1)}^{(\GeoAng)}[v^i]$.

\medskip
\noindent
\underline{\textbf{Analysis of} $\Fullset^{N-1;1} \mathscr{K}_{(\pi-Cancel-2)}^{(\GeoAng)}[v^i]$.}
We apply $\Fullset^{N-1;1}$ to \eqref{E:DIVCURRENTCANEL2} 
(with $\Psi = v^i$ and $Z = \GeoAng$).
We first analyze the difficult product in which all derivatives fall on 
the deformation tensor components:
\begin{align} \label{E:PICANCEL2DIFFICULT}
	\left\lbrace
		- \angLie_{\Fullset}^{N-1;1} \angLie_{\Rad} \angdeformoneformupsharparg{\GeoAng}{\Lunit}
		+ 
		\angdiffuparg{\#} \Fullset^{N-1;1} \deformarg{\GeoAng}{\Lunit}{\Rad}
	\right\rbrace 
				\cdot
	\angdiff v^i.
\end{align}
Using the first estimate in
\eqref{E:GEOANGDEFIMPORTANTLIERADSPHERELANDRADTRACESPHERETERMS}
and the third estimate in
\eqref{E:GEOANGDEFIMPORTANTANGDIVSPHERELANDANGDIFFPILRADTERMS},
we express the terms in braces in 
\eqref{E:PICANCEL2DIFFICULT}
as the sum of $Harmless_{(Wave)}^{\leq N}$
terms and terms involving the order $N$ derivatives of
$\upmu$ and $\mytr \upchi$, which \emph{exactly cancel}.
Also using the bound
$\| \angdiff v^i \|_{L^{\infty}(\Sigma_t^u)} \lesssim \varepsilon$ 
(see \eqref{E:PSIMIXEDUPTOORDERTWOTRANSVERSALIMPROVED}),
we conclude that
$\mbox{\eqref{E:PICANCEL2DIFFICULT}}= Harmless_{(Wave)}^{\leq N}$.
The remaining terms obtained from
applying $\Fullset^{N-1;1}$ 
to \eqref{E:DIVCURRENTCANEL2}
can be shown to be 
$Harmless_{(Wave)}^{\leq N}$
by combining essentially the same argument 
with the estimate \eqref{E:GEOANGPOINTWISE}
and the $L^{\infty}$ estimates of Prop.~\ref{P:IMPROVEMENTOFAUX}.
We have thus proved \eqref{E:PIHARMLESSANALYSIS} 
and \eqref{E:LCONTAININGPIHARMLESSANALYSIS}
for
$\Fullset^{N-1;1} \mathscr{K}_{(\pi-Cancel-2)}^{(\GeoAng)}[v^i]$.

\medskip
\noindent
\underline{\textbf{Analysis of} $\Fullset^{N-1;1} \mathscr{K}_{(\pi-Less \ Dangerous)}^{(\GeoAng)}[v^i]$.}
We apply $\Fullset^{N-1;1}$ to \eqref{E:DIVCURRENTELLIPTIC} 
(with $\Psi = v^i$ and $Z = \GeoAng$).
We first analyze the difficult product in which all derivatives fall on 
the deformation tensor component:
\begin{align} \label{E:PILESSDANGEROUSDIFFICULT}
	\frac{1}{2} \upmu (\angLie_{\Fullset}^{N-1;1} \angdiffuparg{\#} \mytr \angdeform{\GeoAng}) 
	\cdot \angdiff v^i.
\end{align}
Using the fourth estimate in \eqref{E:GEOANGDEFIMPORTANTANGDIVSPHERELANDANGDIFFPILRADTERMS} 
and the simple bounds
$\| \angdiff v^i \|_{L^{\infty}(\Sigma_t^u)} \lesssim \varepsilon$ 
(see \eqref{E:PSIMIXEDUPTOORDERTWOTRANSVERSALIMPROVED}),
$\| \upmu \|_{L^{\infty}(\Sigma_t^u)} \lesssim 1$ 
(see \eqref{E:UPTOONETRANSVERSALDERIVATIVEUPMULINFTY}),
and
$\| \GeoAngFlatRadComponent \|_{L^{\infty}(\Sigma_t^u)} \lesssim \varepsilon$ 
(which follows from \eqref{E:LINEARLYSMALLSCALARSDEPENDINGONGOODVARIABLES},
\eqref{E:PSIMIXEDUPTOORDERTWOTRANSVERSALIMPROVED},
and
\eqref{E:ZSTARLISMALLLISMALLINFTYESTIMATE}),
we deduce that
\begin{align} \label{E:PILESSDIFFICULTDECOMPOSED}
	\frac{1}{2} \upmu (\angLie_{\Fullset}^{N-1;1} \angdiffuparg{\#} \mytr \angdeform{\GeoAng}) 
		\cdot \angdiff v^i
	= \upmu \GeoAngFlatRadComponent (\angdiffuparg{\#} \Fullset_*^{N-1;M} \mytr \upchi)
		\cdot \angdiff v^i
		+ 
		Harmless_{(Wave)}^{\leq N}.
\end{align}
We first consider the case in which 
$\Fullset^{N-1;1}$ contains no factor of $\Lunit$,
which is relevant for proving \eqref{E:PILESSDANGEROUSANALYSIS}.
Then $\Fullset^{N-1;1}$ contains $N-2$ factors of $\GeoAng$
and one factor of $\Rad$. 
We write
$\angdiffuparg{\#} \Fullset^{N-1;M} \mytr \upchi
= 
\ginversesphere
\cdot
\angLie_{\Fullset}^{N-1;M} \angdiff \mytr \upchi
$,
and use the commutator estimate 
\eqref{E:NOSPECIALSTRUCTURETENSORFIELDCOMMUTATORESTIMATE} with 
$\xi = \angdiff \mytr \upchi$, $N-1$ in the role of $N$, and $M=1$,
the estimate \eqref{E:POINTWISEESTIMATESFORCHIANDITSDERIVATIVES},
and the $L^{\infty}$ estimates of Prop.~\ref{P:IMPROVEMENTOFAUX}
to commute the factor of $\Rad$ so that it hits $\mytr \upchi$ first,
thereby obtaining
$
\upmu(\angdiffuparg{\#} \Fullset^{N-1;M} \mytr \upchi)
\cdot \angdiff v^i
=
\upmu (\angdiffuparg{\#} \GeoAng^{N-2} \Rad \mytr \upchi)
\cdot \angdiff v^i
+ 
Harmless_{(Wave)}^{\leq N}
$.
The remaining terms obtained from
applying $\Fullset^{N-1;1}$ to \eqref{E:DIVCURRENTELLIPTIC}
generate products involving $\leq N-2$
derivatives of 
$\angdiff \mytr \angdeform{\GeoAng}$.
We will show that these products are
$Harmless_{(Wave)}^{\leq N}$.
To proceed, we again use the $L^{\infty}$ estimates of Prop.~\ref{P:IMPROVEMENTOFAUX},
the estimate \eqref{E:POINTWISEESTIMATESFORCHIANDITSDERIVATIVES},
the estimate $\| \GeoAngFlatRadComponent \|_{L^{\infty}(\Sigma_t^u)} \lesssim \varepsilon$ mentioned above,
and the fourth estimate in \eqref{E:GEOANGDEFIMPORTANTANGDIVSPHERELANDANGDIFFPILRADTERMS}
(with $\leq N-2$ in the role of $N-1$ in \eqref{E:GEOANGDEFIMPORTANTANGDIVSPHERELANDANGDIFFPILRADTERMS})
to deduce that all of the products under consideration
are $Harmless_{(Wave)}^{\leq N}$.
We clarify that
the estimate \eqref{E:GEOANGDEFIMPORTANTANGDIVSPHERELANDANGDIFFPILRADTERMS}
generates a factor of 
$\mytr \upchi$ 
with $\leq N-1$ derivatives on it 
(located on LHS~\eqref{E:GEOANGDEFIMPORTANTANGDIVSPHERELANDANGDIFFPILRADTERMS}),
which is in contrast to the factor from \eqref{E:PILESSDIFFICULTDECOMPOSED} with $N$ derivatives.
This factor is below-top-order in the sense that we may bound it with
\eqref{E:POINTWISEESTIMATESFORCHIANDITSDERIVATIVES}, 
and the corresponding product 
$
\upmu \GeoAngFlatRadComponent (\angdiffuparg{\#} \Fullset^{\leq N-2;M} \mytr \upchi)
\cdot \angdiff v^i
$
contributes only to the $Harmless_{(Wave)}^{\leq N}$ terms.
The remaining terms obtained from
applying $\Fullset^{N-1;1}$ to \eqref{E:DIVCURRENTELLIPTIC}
can be shown to be 
$Harmless_{(Wave)}^{\leq N}$
by using essentially the same argument and
the $L^{\infty}$ estimates of Prop.~\ref{P:IMPROVEMENTOFAUX}.
We have thus proved \eqref{E:PILESSDANGEROUSANALYSIS}.

We now consider the case in which 
$\Fullset^{N-1;1}$ contains a factor of $\Lunit$,
which is relevant for proving \eqref{E:WAVETHIRDHARMLESSORDERNCOMMUTATORS}.
Noting that the formula \eqref{E:PILESSDIFFICULTDECOMPOSED} still holds,
we use the same commutator argument given in the previous paragraph
to obtain
$
\frac{1}{2} \upmu (\angLie_{\Fullset}^{N-1;1} \angdiffuparg{\#} \mytr \angdeform{\GeoAng}) \cdot \angdiff v^i
=
\upmu \GeoAngFlatRadComponent (\angdiffuparg{\#} \Fullset^{N-1;M} \Lunit \mytr \upchi)
\cdot \angdiff v^i
+ Harmless_{(Wave)}^{\leq N}
$.
Moreover, using the $L^{\infty}$ estimates of Prop.~\ref{P:IMPROVEMENTOFAUX},
\eqref{E:LINEARLYSMALLSCALARSDEPENDINGONGOODVARIABLES},
and
\eqref{E:LUNITCOMMUTEDLUNITSMALLIPOINTWISE},
we deduce that
$
\left|
	\upmu \GeoAngFlatRadComponent (\angdiffuparg{\#} \Fullset^{N-1;M} \Lunit \mytr \upchi)
\cdot \angdiff v^i
\right|
\lesssim
\left|
	\Fullset^{\leq N-1;\leq 1} \Lunit \mytr \upchi)
\right|
= 
Harmless_{(Wave)}^{\leq N}
$.
We have thus proved \eqref{E:LCONTAININGPIHARMLESSANALYSIS} for
the term
$\Fullset^{N-1;1} \mathscr{K}_{(\pi-Less \ Dangerous)}^{(\GeoAng)}[v^i]$.

\medskip
\noindent
\underline{\textbf{Analysis of} $\Fullset^{N-1;1} \mathscr{K}_{(\pi-Good)}^{(\GeoAng)}[v^i]$.}
We apply $\Fullset^{N-1;1}$ to \eqref{E:DIVCURRENTGOOD} 
(with $\Psi = v^i$ and $Z = \GeoAng$).
The main point is that all deformation tensor components on RHS~\eqref{E:DIVCURRENTGOOD}
are hit with an $\Lunit$ derivative. We may therefore bound the products under consideration
using \eqref{E:TOPDERIVATIVESOFTANGENTIALDEFORMATIONINVOLVINGONELUNITDERIVATIVE}-\eqref{E:PARTIITOPDERIVATIVESOFRADDEFORMATIONINVOLVINGONELUNITDERIVATIVE}
and the $L^{\infty}$ estimates of Prop.~\ref{P:IMPROVEMENTOFAUX},
thus concluding that all products are 
$Harmless_{(Wave)}^{\leq N}$.
We have therefore proved
\eqref{E:PIHARMLESSANALYSIS}
and 
\eqref{E:LCONTAININGPIHARMLESSANALYSIS}
for
$\Fullset^{N-1;1} \mathscr{K}_{(\pi-Good)}^{(\GeoAng)}[v^i]$.

\medskip
\noindent
\underline{\textbf{Analysis of} $\Fullset^{N-1;1} \mathscr{K}_{(\Psi)}^{(\GeoAng)}[v^i]$.}
The terms in $\mathscr{K}_{(\Psi)}^{(\GeoAng)}[v^i]$ 
(see \eqref{E:DIVCURRENTPSI}, where $\Psi = v^i$ and $Z = \GeoAng$) 
are of the form
$
\smoothfunction(\BadVar) \pi \Singletan Z v^i
+
\smoothfunction(\BadVar) \angLap v^i
$
where $\Singletan \in \Tanset$, 
$Z \in \Fullset$,
and $\pi$ is one of the following components of $\deform{\GeoAng}$:
$
\pi
\in
\left\lbrace
	\mytr \angdeform{\GeoAng},
	\deformarg{\GeoAng}{\Lunit}{\Rad},
	\deformarg{\GeoAng}{\Rad}{\Radunit},
	\angdeformoneformupsharparg{\GeoAng}{\Lunit},
	\angdeformoneformupsharparg{\GeoAng}{\Rad}
\right\rbrace
$.
We now apply $\Fullset^{N-1;1}$ to the expression
$
\smoothfunction(\BadVar) \pi \Singletan Z v^i
+
\smoothfunction(\BadVar) \angLap v^i
$
and use 
\eqref{E:BOUNDSFORDERIVATIVESOFLAPLACIAN},
\eqref{E:BELOWTOPORDERDERIVATIVESOFTANGENTIALDEFORMATION}-\eqref{E:BELOWTOPORDERDERIVATIVESOFRADDEFORMATION},
and the $L^{\infty}$ estimates of Prop.~\ref{P:IMPROVEMENTOFAUX},
thereby concluding that all products under consideration
are $Harmless_{(Wave)}^{\leq N}$ as desired.
We have thus proved
\eqref{E:PIHARMLESSANALYSIS}
and 
\eqref{E:LCONTAININGPIHARMLESSANALYSIS}
for
$\Fullset^{N-1;1} \mathscr{K}_{(\Psi)}^{(\GeoAng)}[v^i]$.

\medskip
\noindent
\underline{\textbf{Analysis of} $\Fullset^{N-1;1} \mathscr{K}_{(Low)}^{(\GeoAng)}[v^i]$.}
Using Lemma~\ref{L:SCHEMATICDEPENDENCEOFMANYTENSORFIELDS}, 
we see that (see \eqref{E:DIVCURRENTLOW}, where $\Psi = v^i$ and $Z = \GeoAng$) 
$
\mathscr{K}_{(Low)}^{(\GeoAng)}[v^i]
=
\smoothfunction(\Tanset^{\leq 1} \BadVar,\ginversesphere,\angdiff \vec{x},\Rad \threePsi) 
\pi \Singletan \GdVar
$
where 
$\pi$ and $\Singletan$ are as in the previous paragraph.
Hence, we conclude that $\Fullset^{N-1;1} \mathscr{K}_{(Low)}^{(\GeoAng)}[v^i] = Harmless^{\leq N}$
by using the same arguments as in the previous paragraph together with 
Lemmas~\ref{L:POINTWISEFORRECTANGULARCOMPONENTSOFVECTORFIELDS}
and
~\ref{L:POINTWISEESTIMATESFORGSPHEREANDITSDERIVATIVES}
(to bound the derivatives of $\ginversesphere$ and $\angdiff \vec{x}$).
We have thus proved
\eqref{E:PIHARMLESSANALYSIS}
and 
\eqref{E:LCONTAININGPIHARMLESSANALYSIS}
for
$\Fullset^{N-1;1} \mathscr{K}_{(Low)}^{(\GeoAng)}[v^i]$.

\medskip

We have thus established \eqref{E:PIHARMLESSANALYSIS}-\eqref{E:LCONTAININGPIHARMLESSANALYSIS},
which completes the analysis of the desired estimates for $\Psi \in \lbrace v^1, v^2 \rbrace$,
except for the error term bound \eqref{E:GEOANGLOWERORDERRORTERMESTIMATE}
(which we derive below),
in the case that the commutator operator is of the form
$\Fullset^{N-1;1} \GeoAng$,
where $\Fullset^{N-1;1}$ contains exactly one factor of $\Rad$.
We must also establish similar estimates in the remaining cases,
corresponding to the following operators 
on
LHSs~\eqref{E:WAVELISTHEFIRSTCOMMUTATORIMPORTANTTERMS}-\eqref{E:SECONDWAVELISTHEFIRSTCOMMUTATORIMPORTANTTERMS} 
and \eqref{E:WAVEHARMLESSORDERNPURETANGENTIALCOMMUTATORS}-\eqref{E:WAVESECONDHARMLESSORDERNCOMMUTATORS}:
\begin{enumerate}
\item $\GeoAng^{N-1} \Lunit$
\item $\GeoAng^N$
\item $\GeoAng^{N-1} \Rad$
\item $\Fullset^{N-1;1} \Lunit$
	(where $\Fullset^{N-1;1}$ contains exactly one factor of $\Rad$ with all other factors equal to
	$\GeoAng$)
\item $\Tanset^{N-1} \Lunit$ (where $\Tanset^{N-1}$ contains one or more factors of $\Lunit$)
\item $\Fullset^{N-1;1} \Lunit$ (where $\Fullset^{N-1;1}$ contains one or more factors of $\Lunit$)
\item $\Tanset^{N-1} \Rad$ (where $\Tanset^{N-1}$ contains one or more factors of $\Lunit$)
\end{enumerate}
In these remaining seven cases,
we can obtain an analog of the estimate 
\eqref{E:COMMUTEDWAVEGEOANGFIRSTMAINTERMPLUSERRORTERM}
by using the same arguments,
which are based on Lemma~\ref{L:WAVEEQUATIONTOPDERIVATIVESOFVORTICITY}
and the estimates \eqref{E:POINTWISESTIMATESFORTHENULLFOMRS}
and \eqref{E:SIMPLEDEFTENSORPOINTWISEBOUNDS}. 
We note that the error term bound \eqref{E:GEOANGLOWERORDERRORTERMESTIMATE}
remains correct as stated in all of these cases.
We also note that the estimates of Lemma~\ref{L:WAVEEQUATIONTOPDERIVATIVESOFVORTICITY}
yield the explicitly listed terms on 
RHSs~\eqref{E:WAVELISTHEFIRSTCOMMUTATORIMPORTANTTERMS}-\eqref{E:SECONDWAVELISTHEFIRSTCOMMUTATORIMPORTANTTERMS} 
and \eqref{E:WAVEHARMLESSORDERNPURETANGENTIALCOMMUTATORS}-\eqref{E:WAVESECONDHARMLESSORDERNCOMMUTATORS}
that depend on the order $N+1$ derivatives of $\Vortrenormalized$.
Moreover, in these remaining cases, we can use
essentially the same arguments that we used in the case $\Fullset^{N-1;1} \GeoAng$
to establish pointwise estimates for the main term
(that is, the analog of the first one on RHS~\eqref{E:COMMUTEDWAVEGEOANGFIRSTMAINTERMPLUSERRORTERM}).
That is, with the help of 
Lemmas~\ref{L:TANGENTIALCOMMUTEDRADTRCHIANGLAPUPMUCOMPARISON} and \ref{L:IMPORTANTDEFTENSORTERMS},
we can establish analogs of 
\eqref{E:PIHARMLESSANALYSIS}-\eqref{E:LCONTAININGPIHARMLESSANALYSIS}.
The estimates are very similar in nature,  
the only difference being the details of the important terms generated by
Lemma~\ref{L:IMPORTANTDEFTENSORTERMS};
the corresponding important products are precisely
the ones on explicitly listed on 
RHSs~\eqref{E:WAVELISTHEFIRSTCOMMUTATORIMPORTANTTERMS}-\eqref{E:SECONDWAVELISTHEFIRSTCOMMUTATORIMPORTANTTERMS} 
that depend on $N$ derivatives of $\mytr \upchi$.
More precisely, an argument similar to the one that we gave
in the case $\Fullset^{N-1;1} \GeoAng$ yields that in the remaining seven cases stated above,
the estimate \eqref{E:PIDANGERANALYSIS} 
must respectively be replaced with
\begin{align}
\GeoAng^{N-1} \mathscr{K}_{(\pi-Danger)}^{(\Lunit)}[v^i]
& = 0,
	\label{E:LUNITFIRSTPIDANGERANALYSIS} 
		\\
\GeoAng^{N-1} \mathscr{K}_{(\pi-Danger)}^{(\GeoAng)}[v^i]
& = 
	(\Rad v^i) \GeoAng^N \mytr \upchi
	+
	Harmless_{(Wave)}^{\leq N},
	\label{E:ALLGEOANGPIDANGERANALYSIS}
	\\
\GeoAng^{N-1} \mathscr{K}_{(\pi-Danger)}^{(\Rad)}[v^i]
& = (\Rad v^i) \GeoAng^{N-1} \Rad \mytr \upchi
	+
	Harmless_{(Wave)}^{\leq N},
	\label{E:RADFIRSTFOLLOWEDBYALLGEOANGPIDANGERANALYSIS} 
		\\
\Fullset^{N-1;1} \mathscr{K}_{(\pi-Danger)}^{(\Lunit)}[v^i]
& = 0,
	\label{E:ANOTHERLUNITFIRSTPIDANGERANALYSIS} 
		\\
\Tanset^{N-1} \mathscr{K}_{(\pi-Danger)}^{(\Lunit)}[v^i]
& = 0,
	\label{E:LUNITFIRSTFOLLOWEDBYRESTTANPIDANGERANALYSIS} 
		\\
\Fullset^{N-1;1} \mathscr{K}_{(\pi-Danger)}^{(\Lunit)}[v^i]
& = 0,
	\label{E:YETANOTHERLUNITFIRSTPIDANGERANALYSIS} 
		\\
\Tanset^{N-1} \mathscr{K}_{(\pi-Danger)}^{(\Rad)}[v^i]
& = 	Harmless_{(Wave)}^{\leq N},
	\label{E:RADFIRSTFOLLOWEDBYATLEASTONELPIDANGERANALYSIS} 
\end{align}
while the estimate \eqref{E:PILESSDANGEROUSANALYSIS} must
respectively be replaced with
\begin{align}
\GeoAng^{N-1} \mathscr{K}_{(\pi-Less \ Dangerous)}^{(\Lunit)}[v^i]
& =   (\angdiffuparg{\#} v^i) \cdot (\upmu \angdiff \GeoAng^{N-1} \mytr \upchi)
	+
	Harmless_{(Wave)}^{\leq N},
	\label{E:LUNITFIRSTFOLLOWEDBYALLGEOANGPILESSDANGEROUSANALYSIS} 
		\\
\GeoAng^{N-1} \mathscr{K}_{(\pi-Less \ Dangerous)}^{(\GeoAng)}[v^i]
& = \GeoAngFlatRadComponent (\angdiffuparg{\#} v^i) \cdot (\upmu \angdiff \GeoAng^{N-1} \mytr \upchi)
	+
	Harmless_{(Wave)}^{\leq N},
	\label{E:ALLGEOANGLESSDANGEROUSANALYSIS} 
	\\
\GeoAng^{N-1} \mathscr{K}_{(\pi-Less \ Dangerous)}^{(\Rad)}[v^i]
& = - (\upmu \angdiffuparg{\#} v^i) \cdot (\upmu \angdiff \GeoAng^{N-1} \mytr \upchi)
	+
	Harmless_{(Wave)}^{\leq N},
	\label{E:RADFIRSTFOLLOWEDBYALLGEOANGPILESSDANGEROUSANALYSIS} 
	\\
\Fullset^{N-1;1} \mathscr{K}_{(\pi-Less \ Dangerous)}^{(\Lunit)}[v^i]
& = 
	(\angdiffuparg{\#} v^i) \cdot (\upmu \angdiff \GeoAng^{N-2} \Rad \mytr \upchi)
	+ 
	Harmless_{(Wave)}^{\leq N},
	\label{E:LUNITFIRSTPILESSDANGEROUSANALYSIS}
		\\
\Tanset^{N-1;1} \mathscr{K}_{(\pi-Less \ Dangerous)}^{(\Lunit)}[v^i]
& = 
	Harmless_{(Wave)}^{\leq N},
	\label{E:LUNITFIRSTTHENANOTHERLPILESSDANGEROUSANALYSIS}
		\\
\Fullset^{N-1;1} \mathscr{K}_{(\pi-Less \ Dangerous)}^{(\Lunit)}[v^i]
& = 
	Harmless_{(Wave)}^{\leq N},
	\label{E:LUNITFIRSTYETANOTHERPILESSDANGEROUSANALYSIS}
		\\
\Tanset^{N-1} \mathscr{K}_{(\pi-Less \ Dangerous)}^{(\Rad)}[v^i]
& = Harmless_{(Wave)}^{\leq N}
	\label{E:RADFIRSTFOLLOWEDBYALLTANSETWITHONELPILESSDANGEROUSANALYSIS} 
\end{align}
(and all remaining terms are $Harmless_{(Wave)}^{\leq N}$, as in \eqref{E:PIHARMLESSANALYSIS}).

Having treated the difficult main term in all cases, 
we now establish \eqref{E:WAVEERRORTERMSAREHARMLESS}.
We start by bounding the terms
on RHS~\eqref{E:GEOANGLOWERORDERRORTERMESTIMATE}
involving $\leq 10$ derivatives
of the factors 
$\mytr \angdeform{\Lunit}$,
$\mytr \angdeform{\Rad}$,
and
$\mytr \angdeform{\GeoAng}$.
Using \eqref{E:SIMPLEDEFTENSORPOINTWISEBOUNDS}, we see that it suffices
to show that
\begin{align}
& \label{E:TERMSWITHBELOWTOPORDERDEFORMATIONTENSORDERIVATIVESAREHARMLESS}
\mathop{\sum_{N_2 + N_3 \leq N-1}}_{N_2 \leq N-2}
		\sum_{M_2 + M_3 \leq 1}
		\sum_{\Singletan_1, \Singletan_2 \in \Tanset}
		\left|
		\Fullset^{N_2;M_2}
			\left(
				\upmu 
				\D_{\alpha} 
				\Jcurrent{\Singletan_2}^{\alpha}[\Fullset^{N_3;M_3} v^i]
			\right)
		\right|,
			\\
	&
		\mathop{\sum_{N_2 + N_3 \leq N-1}}_{N_2 \leq N-2}
		\sum_{\Singletan \in \Tanset}
		 \left|
			\Tanset^{N_2}
			\left(
				\upmu 
				\D_{\alpha} 
				\Jcurrent{\Rad}^{\alpha}[\Tanset^{N_3} v^i]
			\right)
		\right|,
			\notag \\
	&
		\mathop{\sum_{N_2 + N_3 \leq N-1}}_{N_2 \leq N-2}
		\sum_{\Singletan \in \Tanset}
		\left|
			\Tanset^{N_2}
			\left(
				\upmu 
				\D_{\alpha} 
				\Jcurrent{\Singletan}^{\alpha}[\Tanset^{N_3} v^i]
			\right)
		\right|
			\notag  \\
		& = Harmless_{(Wave)}^{\leq N}.
		\notag
\end{align}
It suffices to again decompose,
with the help of \eqref{E:DIVCOMMUTATIONCURRENTDECOMPOSITION},
the terms $\D_{\alpha} \Jcurrent{\cdot}^{\alpha}$
in \eqref{E:TERMSWITHBELOWTOPORDERDEFORMATIONTENSORDERIVATIVESAREHARMLESS}
and to show that all constituent parts, such as
$
\Fullset^{N_2;M_2} \mathscr{K}_{(\pi-Danger)}^{(\GeoAng)}[\Fullset^{N_3;M_3} v^i]
$
and
$\Tanset^{N_2} \mathscr{K}_{(Low)}^{(\Rad)}[\Tanset^{N_3} v^i]$,
are $Harmless_{(Wave)}^{\leq N}$.
To this end, we repeat the proofs of the above estimates, including
\eqref{E:LUNITFIRSTPIDANGERANALYSIS}-\eqref{E:RADFIRSTFOLLOWEDBYALLTANSETWITHONELPILESSDANGEROUSANALYSIS},
but with $N_2$ (from LHS~\eqref{E:TERMSWITHBELOWTOPORDERDEFORMATIONTENSORDERIVATIVESAREHARMLESS})
in place of $N-1$
and $\Fullset^{N_3;M_3} v^i$
or $\Tanset^{N_3} v^i$
in place of the explicitly written
$v^i$ factors.
The same arguments given above yield that all products
$= Harmless_{(Wave)}^{\leq N_2 + N_3 + 1} \leq Harmless_{(Wave)}^{\leq N}$,
except for the ones corresponding to the 
explicitly written ones on RHSs
\eqref{E:PIDANGERANALYSIS}-\eqref{E:PILESSDANGEROUSANALYSIS},
\eqref{E:ALLGEOANGPIDANGERANALYSIS},
\eqref{E:RADFIRSTFOLLOWEDBYALLGEOANGPIDANGERANALYSIS},
and
\eqref{E:LUNITFIRSTFOLLOWEDBYALLGEOANGPILESSDANGEROUSANALYSIS}-\eqref{E:LUNITFIRSTPILESSDANGEROUSANALYSIS}.
For example, the analog of
the explicitly written term on RHS~\eqref{E:PIDANGERANALYSIS}
is $(\GeoAng^{N_2} \Rad \mytr \upchi) \Rad \Fullset^{N_3;M_3} v^i$
while the analog of
the explicitly written term on RHS~\eqref{E:PILESSDANGEROUSANALYSIS}
is
$
\upmu \GeoAngFlatRadComponent (\angdiffuparg{\#} \GeoAng^{N_2} \Rad \mytr \upchi)
\cdot 
\angdiff \Fullset^{N_3;M_3} v^i$.
We now explain why these explicitly written products are $Harmless_{(Wave)}^{\leq N}$ too.
The important point is that since $N_2 \leq N-2$
on LHS~\eqref{E:TERMSWITHBELOWTOPORDERDEFORMATIONTENSORDERIVATIVESAREHARMLESS},
the factors of $\mytr \upchi$ in these products
are hit with no more than $N-1$ derivatives. We may therefore
pointwise bound these factors using \eqref{E:POINTWISEESTIMATESFORCHIANDITSDERIVATIVES}.
Given this observation, the fact that the products under consideration are
$Harmless_{(Wave)}^{\leq \max \lbrace N_3,N_2+2 \rbrace} \leq Harmless_{(Wave)}^{\leq N}$
follows from the same arguments given
in our prior analysis of $\Fullset^{N-1;1} \mathscr{K}_{(\pi-Danger)}^{(\GeoAng)}[\Psi]$,
$\cdots$,
$\Fullset^{N-1;1} \mathscr{K}_{(Low)}^{(\GeoAng)}[v^i]$.
We have thus shown that
the products
on RHS~\eqref{E:GEOANGLOWERORDERRORTERMESTIMATE}
involving $\leq 10$ derivatives
of the factors 
$\mytr \angdeform{\Lunit}$,
$\mytr \angdeform{\Rad}$,
and
$\mytr \angdeform{\GeoAng}$
are $Harmless_{(Wave)}^{\leq N}$.

To complete the proof of \eqref{E:WAVEERRORTERMSAREHARMLESS},
we must bound the terms on RHS~\eqref{E:GEOANGLOWERORDERRORTERMESTIMATE}
with $11 \leq N_1 \leq N-1 \leq 19$
(which implies that $N_2 + N_3 \leq 8$).
The arguments given in the previous paragraph imply 
that the factors on RHS~\eqref{E:GEOANGLOWERORDERRORTERMESTIMATE} corresponding to $N_2$,
such as
$\Fullset^{N_2;M_2}
			\left(
				\upmu 
				\D_{\alpha} 
				\Jcurrent{\Singletan_2}^{\alpha}[\Fullset^{N_3;M_3} \Psi]
			\right)
$,
are $Harmless_{(Wave)}^{\leq \max \lbrace N_2 + N_3 + 1,N_2 + 2 \rbrace}
\leq 
Harmless_{(Wave)}^{\leq 10}
$.
In particular, 
the $L^{\infty}$ estimates of Prop.~\ref{P:IMPROVEMENTOFAUX}
imply that 
$\| Harmless_{(Wave)}^{\leq 10} \|_{L^{\infty}(\Sigma_t^u)}
\lesssim \varepsilon
$.
Moreover, from
\eqref{E:CONNECTIONBETWEENANGLIEOFGSPHEREANDDEFORMATIONTENSORS}
and 
\eqref{E:POINTWISEESTIMATESFORGSPHEREANDITSSTARDERIVATIVES},
we find that 
the factors 
$\Fullset^{N_1;M_1} \mytr \angdeform{\Singletan_1}$,
$\Tanset^{N_1} \mytr \angdeform{\Singletan}$,
and $\Tanset^{N_1} \mytr \angdeform{\Rad}$
on RHS~\eqref{E:GEOANGLOWERORDERRORTERMESTIMATE}
$= Harmless_{(Wave)}^{\leq N_1 + 1}$.
Combining this bound with the estimate
$
\| Harmless_{(Wave)}^{\leq 10} \|_{L^{\infty}(\Sigma_t^u)}
\lesssim \varepsilon
$,
we conclude that the products under consideration
$= Harmless_{(Wave)}^{\leq N_1 + 1} \leq Harmless_{(Wave)}^{\leq N}$
as desired. We have thus proved
\eqref{E:WAVEERRORTERMSAREHARMLESS},
which completes the proof of Prop.~\ref{P:WAVEIDOFKEYDIFFICULTENREGYERRORTERMS}
except for the estimates \eqref{E:KEYDIFFICULTFORRENORMALZIEDDENSITY}
for the quantity $\Densrenormalized  - v^1$.

To prove \eqref{E:KEYDIFFICULTFORRENORMALZIEDDENSITY},
we first subtract
\eqref{E:VELOCITYWAVEEQUATION} with $i=1$ from \eqref{E:RENORMALIZEDDENSITYWAVEEQUATION} 
to obtain
\begin{align} \label{E:DENSRENORALIZEDMINUSV1WAVEEQUATION}
	\upmu \square_{g(\threePsi)} (\Densrenormalized  - v^1)
	& = \upmu \mathscr{Q}
			- 
			\upmu \mathscr{Q}^1
			+ 
			[1a] (\exp \Densrenormalized) \Speed^2 (\upmu \partial_a \Vortrenormalized)
		  - 
			2 [1a] (\exp \Densrenormalized) \Vortrenormalized (\upmu \Transport v^a).
\end{align}
That is, 
$\Densrenormalized  - v^1$
solves the covariant wave equation with 
inhomogeneous terms equal to RHS~\eqref{E:DENSRENORALIZEDMINUSV1WAVEEQUATION}.
We now repeat the above proofs of the estimates for $v^i$,
with $\Densrenormalized  - v^1$
in the role of $v^i$ and RHS~\eqref{E:DENSRENORALIZEDMINUSV1WAVEEQUATION}
in the role of RHS~\eqref{E:VELOCITYWAVEEQUATION}.
Using nearly identical arguments, we obtain 
\eqref{E:KEYDIFFICULTFORRENORMALZIEDDENSITY}.
For clarity, 
we note that 
$
\Densrenormalized  - v^1
= \Psi_0 - \Psi_1
$
(see Def.~\ref{D:VECPSI})
and hence the derivatives of 
$
\Densrenormalized  - v^1
$
can be controlled, via the triangle inequality,
in terms of the derivatives of $\threePsi$.
This completes the proof of Prop.~\ref{P:WAVEIDOFKEYDIFFICULTENREGYERRORTERMS}.

$\hfill \qed$

\subsection{Proof of Prop.~\ref{P:VORTICITYIDOFKEYDIFFICULTENREGYERRORTERMS}}
\label{SS:PROOFOFPROPVORTICITYIDOFKEYDIFFICULTENREGYERRORTERMS}
See Sect.~\ref{SS:OFTENUSEDESTIMATES} for some comments on the analysis.
Throughout we silently use the
Definition \ref{D:HARMLESSTERMS} of 
$Harmless_{(Vort)}^{\leq N}$ terms.

We first prove \eqref{E:VORTICITYGEOANGANGISTHEFIRSTCOMMUTATORIMPORTANTTERMS}.
	The main point is to identify the products that depend
	on the order $N+1$ derivatives of $\upmu$
	or the order $N$ derivatives of $\mytr \upchi$; all other products will be shown to be
	$Harmless_{(Vort)}^{\leq N+1}$.
	To proceed, we apply $\GeoAng^{N+1}$
	to the transport equation \eqref{E:RENORMALIZEDVORTICTITYTRANSPORTEQUATION}
	and use
	\eqref{E:COMMUTATORFORMULAFORTRANSPORTRENORMANDGEOANG}
	to commute $\GeoAng^{N+1}$ through $\upmu \Transport$.
	Using also
	\eqref{E:NOSPECIALSTRUCTUREDIFFERENTIATEDGEOANGDEFORMSPHERELSHARPPOINTWISE} 
	and
	\eqref{E:STARDIFFERENTIATEDGEOANGDEFORMSPHERERADPOINTWISE}
	with $M=0$
	and
	the $L^{\infty}$ estimates of Prop.~\ref{P:IMPROVEMENTOFAUX},
	we find that
	\begin{align} \label{E:VORTICITYTOPORDERGEOANGDERIVATIVESEQUATIONFIRSTCALCULATION}
		\upmu \Transport \GeoAng^{N+1} \Vortrenormalized
			& = - (\GeoAng^{N+1} \upmu) \Lunit \Vortrenormalized
					+ 
					\left\lbrace
						\upmu \angLie_{\GeoAng}^N \angdeformoneformupsharparg{\GeoAng}{\Lunit} 
						+ \angLie_{\GeoAng}^N \angdeformoneformupsharparg{\GeoAng}{\Rad}
					\right\rbrace
					\cdot \angdiff \Vortrenormalized
					+
					Harmless_{(Vort)}^{\leq N+1}.
	\end{align}
	Note that we have isolated all of the top-order derivatives
	of deformation tensors in the terms in braces on
	RHS~\eqref{E:VORTICITYTOPORDERGEOANGDERIVATIVESEQUATIONFIRSTCALCULATION}.
	Moreover, we clarify that the smallness factors $\varepsilon$
	on RHS~\eqref{E:VORTICITYHARMESSTERMPOINTWISEESTIMATE}
	(with $N+1$ in the role of $N$)
	come from the estimate \eqref{E:VORTICITYUPTOTWOTRANSVERSALLINFTY}
	for the low-order derivatives of $\Vortrenormalized$,
	which appear as factors in quadratic terms that multiply 
	high-order derivatives of deformation tensors.
	In addition, 
	from 
	\eqref{E:HIGHGEOANGDERIVATIVESOFUPMUINTERMSOFRADTRCHI},
	\eqref{E:MAINTERMINPUREGEOANGDERIVATIVESOFUPMUANGDEFORMGEOANGLUNITANDGEOANGRAD},
	and the simple bounds
	$\left\|\Singletan \Vortrenormalized \right\|_{L^{\infty}(\Sigma_t^u)},
		\,
	\left\| \GeoAngFlatRadComponent \right\|_{L^{\infty}(\Sigma_t^u)}
	\lesssim \varepsilon
	$
	(which follow from Lemma~\ref{L:SCHEMATICDEPENDENCEOFMANYTENSORFIELDS} and Prop.~\ref{P:IMPROVEMENTOFAUX}),
	we find that
	\begin{align}
		(\GeoAng^{N+1} \upmu) \Lunit \Vortrenormalized
		& = g(\GeoAng,\GeoAng) (\GeoAng^{N-1} \Rad \mytr \upchi) \Lunit \Vortrenormalized
			+ \mbox{\upshape Error},
			\label{E:FIRSTMAINTERMINPROOFOFVORTICITYLISTHEFIRSTCOMMUTATORIMPORTANTTERMS} \\
		\left\lbrace
			\upmu \angLie_{\GeoAng}^N \angdeformoneformupsharparg{\GeoAng}{\Lunit} 
			+ 
			\angLie_{\GeoAng}^N \angdeformoneformupsharparg{\GeoAng}{\Rad}
		\right\rbrace
		\cdot \angdiff \Vortrenormalized
	& = \GeoAngFlatRadComponent (\angdiffuparg{\#} \GeoAng^N \upmu) \cdot \angdiff \Vortrenormalized
			+ Harmless_{(Vort)}^{\leq N+1}
				\label{E:SECONDMAINTERMINPROOFOFVORTICITYLISTHEFIRSTCOMMUTATORIMPORTANTTERMS} \\
	& = \GeoAngFlatRadComponent (\GeoAng^{N-1} \Rad \mytr \upchi) \GeoAng \Vortrenormalized
		+ Harmless_{(Vort)}^{\leq N+1}.
				\notag 
	\end{align}
	Combining 
	\eqref{E:VORTICITYTOPORDERGEOANGDERIVATIVESEQUATIONFIRSTCALCULATION}
	and
\eqref{E:FIRSTMAINTERMINPROOFOFVORTICITYLISTHEFIRSTCOMMUTATORIMPORTANTTERMS}-\eqref{E:SECONDMAINTERMINPROOFOFVORTICITYLISTHEFIRSTCOMMUTATORIMPORTANTTERMS},
we arrive at the desired estimate \eqref{E:VORTICITYGEOANGANGISTHEFIRSTCOMMUTATORIMPORTANTTERMS}.

The proof of \eqref{E:VORTICITYLISTHEFIRSTCOMMUTATORIMPORTANTTERMS}
is similar but relies on
\eqref{E:COMMUTATORFORMULAFORTRANSPORTRENORMANDLUNIT}
in place of
\eqref{E:COMMUTATORFORMULAFORTRANSPORTRENORMANDGEOANG}
and
\eqref{E:MAINTERMINPUREGEOANGDERIVATIVESOFUPMUANGDEFORMLUNITRAD}
in place of
\eqref{E:MAINTERMINPUREGEOANGDERIVATIVESOFUPMUANGDEFORMGEOANGLUNITANDGEOANGRAD};
we omit the details, noting only that the term
$- (\GeoAng^N \Lunit \upmu) \Lunit \Vortrenormalized$,
which is an analog of the first term 
on RHS~\eqref{E:VORTICITYTOPORDERGEOANGDERIVATIVESEQUATIONFIRSTCALCULATION},
is $Harmless_{(Vort)}^{\leq N+1}$
in view of the bound
$\left\|\Lunit \Vortrenormalized \right\|_{L^{\infty}(\Sigma_t^u)}
	\lesssim \varepsilon
$
mentioned above and inequality \eqref{E:LLUNITUPMUFULLSETSTARCOMMUTEDPOINTWISEESTIMATE} 
with $M=0$.

We now prove \eqref{E:TRIVIALVORTICITYHARMLESSORDERNPLUSONECOMMUTATORS}.
Using the same arguments we used to derive
\eqref{E:VORTICITYLISTHEFIRSTCOMMUTATORIMPORTANTTERMS}
and
\eqref{E:VORTICITYGEOANGANGISTHEFIRSTCOMMUTATORIMPORTANTTERMS},
except now bounding the deformation tensor components 
from the formulas \eqref{E:COMMUTATORFORMULAFORTRANSPORTRENORMANDLUNIT}
and
\eqref{E:COMMUTATORFORMULAFORTRANSPORTRENORMANDGEOANG}
in magnitude by $\lesssim 1$ 
via \eqref{E:BELOWTOPORDERDERIVATIVESOFRADDEFORMATION} with $N=1$
and the $L^{\infty}$ estimates of Prop.~\ref{P:IMPROVEMENTOFAUX},
we deduce that the RHS of the equation
$
\upmu \Transport \Singletan \Vortrenormalized
= \cdots
$
is in magnitude 
$\lesssim 
\left|
	\Lunit \Vortrenormalized
\right|
+
\left|
	\angdiff \Vortrenormalized
\right|
\lesssim 
\left|
	\Tanset^{\leq 1} \Vortrenormalized
\right|
$.
This implies \eqref{E:TRIVIALVORTICITYHARMLESSORDERNPLUSONECOMMUTATORS}.

We now prove \eqref{E:VORTICITYHARMLESSORDERNPLUSONECOMMUTATORS},
starting with the case
$\Tanset^{N+1} = \Tanset^N \GeoAng$.
Then the same arguments we used to prove
\eqref{E:VORTICITYGEOANGANGISTHEFIRSTCOMMUTATORIMPORTANTTERMS}
yield
\begin{align} \label{E:PROOFGEOANGFIRSTVORTICITYHARMLESSORDERNPLUSONECOMMUTATORS}
		\upmu \Transport \Tanset^N \GeoAng \Vortrenormalized
			& = - (\Tanset^N \GeoAng \upmu) \Lunit \Vortrenormalized
					+ 
					\left\lbrace
						\upmu \angLie_{\Tanset}^N \angdeformoneformupsharparg{\GeoAng}{\Lunit} 
						+ \angLie_{\Tanset}^N \angdeformoneformupsharparg{\GeoAng}{\Rad}
					\right\rbrace
					\cdot \angdiff \Vortrenormalized
					+
					Harmless_{(Vort)}^{\leq N+1}.
\end{align}
The key point is that by assumption, 
the operator $\Tanset^N$ contains a factor of $\Lunit$.
Hence, we may use
the commutator estimate \eqref{E:LESSPRECISEPERMUTEDVECTORFIELDSACTINGONFUNCTIONSCOMMUTATORESTIMATE} with 
$f = \upmu$
and the $L^{\infty}$ estimates of Prop.~\ref{P:IMPROVEMENTOFAUX}
to commute the factor of $\Lunit$ in $\Tanset^N \GeoAng \upmu$
so that it hits $\upmu$ first. Also using
the pointwise estimate \eqref{E:LLUNITUPMUFULLSETSTARCOMMUTEDPOINTWISEESTIMATE},
we find that
$
\Tanset^N \GeoAng \upmu
=
Harmless_{(Vort)}^{\leq N+1}
$.
To bound the terms
$
\angLie_{\Tanset}^N \angdeformoneformupsharparg{\GeoAng}{\Lunit}
$
and
$\angLie_{\Tanset}^N \angdeformoneformupsharparg{\GeoAng}{\Rad}$
in magnitude, we use the pointwise estimate
\eqref{E:TOPDERIVATIVESOFTANGENTIALDEFORMATIONINVOLVINGONELUNITDERIVATIVE}.
Combining these pointwise estimates with the $L^{\infty}$ estimates of Prop.~\ref{P:IMPROVEMENTOFAUX},
we conclude that 
all terms on RHS~\eqref{E:PROOFGEOANGFIRSTVORTICITYHARMLESSORDERNPLUSONECOMMUTATORS}
are of the form $Harmless_{(Vort)}^{\leq N+1}$.
We have thus proved \eqref{E:VORTICITYHARMLESSORDERNPLUSONECOMMUTATORS} in the case
$\Tanset^{N+1} = \Tanset^N \GeoAng$.

To finish the proof of \eqref{E:VORTICITYHARMLESSORDERNPLUSONECOMMUTATORS},
it remains only for us to consider the case 
$\Tanset^{N+1} = \Tanset^N \Lunit$
where $\Tanset^N$ contains a factor of $\Lunit$.
Using the same arguments that we used to prove
\eqref{E:VORTICITYLISTHEFIRSTCOMMUTATORIMPORTANTTERMS},
we derive the following analog
of \eqref{E:PROOFGEOANGFIRSTVORTICITYHARMLESSORDERNPLUSONECOMMUTATORS}:
\begin{align} \label{E:PROOFLUNITFIRSTVORTICITYHARMLESSORDERNPLUSONECOMMUTATORS}
		\upmu \Transport \Tanset^N \Lunit \Vortrenormalized
			& = - (\Tanset^N \Lunit \upmu) \Lunit \Vortrenormalized
					+
					\left\lbrace
						\angLie_{\Tanset}^N \angdeformoneformupsharparg{\Lunit}{\Rad}
					\right\rbrace
					\cdot \angdiff \Vortrenormalized
					+
					Harmless_{(Vort)}^{\leq N+1},
\end{align}
where the operator $\Tanset^N$ in \eqref{E:PROOFLUNITFIRSTVORTICITYHARMLESSORDERNPLUSONECOMMUTATORS}
contains a factor of $\Lunit$.
The remainder of the proof now proceeds as in the case
$\Tanset^{N+1} = \Tanset^N \GeoAng$,
but with the estimate
\eqref{E:PARTITOPDERIVATIVESOFRADDEFORMATIONINVOLVINGONELUNITDERIVATIVE}
in place of
\eqref{E:TOPDERIVATIVESOFTANGENTIALDEFORMATIONINVOLVINGONELUNITDERIVATIVE}.
This completes the proof of Prop.~\ref{P:VORTICITYIDOFKEYDIFFICULTENREGYERRORTERMS}.

$\hfill \qed$

\section{Energy estimates}
\label{S:ENERGYESTIMATES}
In this section, we derive the most important estimates of the article:
a priori energy estimates for the solution.
The main result is Prop.~\ref{P:MAINAPRIORIENERGY} (see Sect.~\ref{SS:STATEMENTOFMAINAPRIORIESTIMATES}),
which we prove in Sect.~\ref{SS:PROOFOFPROPMAINAPRIORIENERGY} via a lengthy Gronwall argument,
after deriving many preliminary estimates.
The main preliminary results of this section
are Props.~\ref{P:WAVEENERGYINTEGRALINEQUALITIES} and \ref{P:VORTICITYENERGYINTEGRALINEQUALITIES}
(see Sect.~\ref{SS:STATEMENTOFMAININTEGRALINEQUALITIES}),
in which we derive, with the help of the pointwise estimates 
of Sect.~\ref{S:POINTWISEESTIMATESFORWAVEEQUATIONERRORTERMS} ,
integral inequalities for the fundamental $L^2$-controlling quantities
defined in Sect.~\ref{S:FUNDAMENTALL2CONTROLLINGQUANTITIES}.
To prove Props.~\ref{P:WAVEENERGYINTEGRALINEQUALITIES} and \ref{P:VORTICITYENERGYINTEGRALINEQUALITIES}, 
we must bound the error integrals 
on the right-hand sides of the energy-null flux identities
of Props.~\ref{P:DIVTHMWITHCANCELLATIONS} and \ref{P:ENERGYIDENTITYRENORMALIZEDVORTICITY}
and their higher-order analogs.
We divide the error integrals into various classes,
which we bound in 
Sects.~\ref{SS:L2FOREIKONALNOMODIFIED}-\ref{SS:WAVEEQUATIONENERGYESTIMATESLOSSOFONEDERIVATIVE}.
We combine all of these estimates into proofs of
Props.~\ref{P:WAVEENERGYINTEGRALINEQUALITIES} and \ref{P:VORTICITYENERGYINTEGRALINEQUALITIES}
in Sects.~\ref{SS:PROOFOFPROPWAVEENERGYINTEGRALINEQUALITIES}
and \ref{SS:PROOFOFPROPVORTICITYENERGYINTEGRALINEQUALITIES} respectively.
Compared to the energy estimates in previous works, 
the new feature of the present work
is that we derive a priori energy estimates for the specific vorticity
and control various error integrals that depend on it. In particular, we control
the influence that the specific vorticity has on the variables
$\Densrenormalized$, $v^1$, and $v^2$, which solve covariant wave equations
with vorticity-dependent source terms. We derive estimates
for the specific vorticity itself in Sect.~\ref{SS:MAINVORTICITYAPRIORIENERGYESTIMATES}.
To simplify our proofs, we use an energy bootstrap argument; see Sect.~\ref{SS:L2BOOTSTRAP}
for the bootstrap assumptions.

\subsection{Statement of the main a priori energy estimates}
\label{SS:STATEMENTOFMAINAPRIORIESTIMATES}
We start by stating the proposition featuring our main a priori
energy estimates. Its proof is located 
in Sect.~\ref{SS:PROOFOFPROPMAINAPRIORIENERGY}.

\begin{proposition}[\textbf{The main a priori energy estimates}]
	\label{P:MAINAPRIORIENERGY}
	Consider the $L^2$-controlling quantities 
	$\lbrace \totmax{N}(t,u) \rbrace_{N=1,\cdots,20}$,
	$\lbrace \Vorttotmax{N}(t,u) \rbrace_{N=0,\cdots,21}$,
	and 
	$\lbrace \coercivespacetimemax{N}(t,u) \rbrace_{N=1,\cdots,20}$
	from Defs.~\ref{D:MAINCOERCIVEQUANT}
	and \ref{D:COERCIVEINTEGRAL}.
	There exists a constant $C > 0$ such that
	under the data-size and bootstrap assumptions 
	of Sects.~\ref{SS:FLUIDVARIABLEDATAASSUMPTIONS}-\ref{SS:PSIBOOTSTRAP}
	and the smallness assumptions of Sect.~\ref{SS:SMALLNESSASSUMPTIONS}, 
	the following estimates hold
	for $(t,u) \in [0,\Tboot) \times [0,U_0]$:
	\begin{subequations}
	\begin{align}
		\sqrt{\totmax{15+K}}(t,u)
		+ \sqrt{\coercivespacetimemax{15+K}}(t,u)
		& \leq C \mathring{\upepsilon} \upmu_{\star}^{-(K+.9)}(t,u),
			&& (0 \leq K \leq 5),
				\label{E:MULOSSMAINAPRIORIENERGYESTIMATES} \\
		\sqrt{\totmax{[1,14]}}(t,u)
		+ \sqrt{\coercivespacetimemax{[1,14]}}(t,u)
		& \leq C \mathring{\upepsilon},
		&& 
			\label{E:NOMULOSSMAINAPRIORIENERGYESTIMATES}
				\\
		\sqrt{\VorttotTanmax{21}}(t,u)
		& \leq C \mathring{\upepsilon} \upmu_{\star}^{-6.4}(t,u),
		 && \label{E:TOPVORTMULOSSMAINAPRIORIENERGYESTIMATES} \\
		\sqrt{\VorttotTanmax{16+K}}(t,u)
		& \leq C \mathring{\upepsilon} \upmu_{\star}^{-(K+.9)}(t,u),
			&& (0 \leq M \leq 4),
				\label{E:VORTMULOSSMAINAPRIORIENERGYESTIMATES} \\
		\sqrt{\VorttotTanmax{\leq 15}}(t,u)
		& \leq C \mathring{\upepsilon}.
		&&
			\label{E:VORTNOMULOSSMAINAPRIORIENERGYESTIMATES}
	\end{align}
	\end{subequations}
\end{proposition}

To initiate the proof of Prop.~\ref{P:MAINAPRIORIENERGY},
we provide the following simple lemma,
which shows that the fundamental $L^2$-controlling quantities are
initially $\lesssim \mathring{\upepsilon}^2$.
\begin{lemma}[\textbf{The fundamental controlling quantities are initially small}]
\label{L:INITIALSIZEOFL2CONTROLLING}
	Under the data-size assumptions
	of Sects.~\ref{SS:FLUIDVARIABLEDATAASSUMPTIONS} and \ref{SS:DATAFOREIKONALFUNCTIONQUANTITIES},
	the following estimates hold for 
	$t \in [0,2 \TranminusdatasizeWithFactor^{-1}]$
	and
	$u \in [0,U_0]$: 
	\begin{subequations}
	\begin{align} \label{E:WAVEINITIALSIZEOFL2CONTROLLING}
		\totmax{[0,20]}(0,u),
			\, 
		\totmax{[0,20]}(t,0)
		&
		\leq C \mathring{\upepsilon}^2,
			\\
		\VorttotTanmax{\leq 21}(0,u),
			\,
		\VorttotTanmax{\leq 21}(t,0)
		& \leq C \mathring{\upepsilon}^2.
		\label{E:VORTICITYINITIALSIZEOFL2CONTROLLING}
	\end{align}
	\end{subequations}
\end{lemma}

\begin{proof}
	We first note that by
	\eqref{E:UPMUDATATANGENTIALLINFINITYCONSEQUENCES}
	and \eqref{E:SIMPLEUPMUALONGPOESTIMATE},
	we have $\upmu \approx 1$ along $\Sigma_0^1$
	and along $\mathcal{P}_0^{2 \TranminusdatasizeWithFactor^{-1}}$.
	Using these estimates,
	Def.~\ref{D:MAINCOERCIVEQUANT},
	and
	Lemma~\ref{L:ORDERZEROCOERCIVENESS},
	we see that
	$
	\displaystyle
		\totmax{[0,20]}(0,u)
		\lesssim 
		\sum_{\Psi \in \lbrace \Densrenormalized - v^1,v^1,v^2 \rbrace}
		\left\|
			\Fullset_*^{[1,21];\leq 2} \Psi
		\right\|_{L^2(\Sigma_0^u)}^2
	$
	and
	$
	\displaystyle
		\VorttotTanmax{\leq 21}(0,u)
		\lesssim 
		\left\|
			\Tanset^{\leq 21} \Vortrenormalized
		\right\|_{L^2(\Sigma_0^u)}^2
	$.
	Similarly,
	for $0 \leq t \leq 2 \TranminusdatasizeWithFactor^{-1}$,
	we have
	$
	\displaystyle
		\totmax{[0,20]}(t,0)
		\lesssim 
		\sum_{\Psi \in \lbrace \Densrenormalized  - v^1,v^1,v^2 \rbrace}
		\left\|
			\Fullset_*^{[1,21];\leq 1} \Psi
		\right\|_{L^2(\mathcal{P}_0^t)}^2
	$
	and
	$
	\displaystyle
		\VorttotTanmax{\leq 21}(t,0)
		\lesssim 
		\left\|
			\Tanset^{\leq 21} \Vortrenormalized
		\right\|_{L^2(\mathcal{P}_0^t)}^2
	$.
	The estimates \eqref{E:WAVEINITIALSIZEOFL2CONTROLLING}-\eqref{E:VORTICITYINITIALSIZEOFL2CONTROLLING} 
	follow from these estimates
	and the initial data assumptions 
	\eqref{E:L2SMALLDATAASSUMPTIONSALONGSIGMA0}
	and \eqref{E:PSIL2SMALLDATAASSUMPTIONSALONGP0}.
\end{proof}

\subsection{Statement of the integral inequalities that we use to derive a priori estimates}
\label{SS:STATEMENTOFMAININTEGRALINEQUALITIES}
We prove Prop.~\ref{P:MAINAPRIORIENERGY} using a lengthy Gronwall
argument based on the sharp estimates for $\upmu$ derived in Sect.~\ref{S:SHARPESTIMATESFORUPMU}
and the energy inequalities provided by the next two proposition, 
Props.~\ref{P:WAVEENERGYINTEGRALINEQUALITIES}
and \ref{P:VORTICITYENERGYINTEGRALINEQUALITIES},
which we prove 
in Sects.~\ref{SS:PROOFOFPROPWAVEENERGYINTEGRALINEQUALITIES}
and \ref{SS:PROOFOFPROPVORTICITYENERGYINTEGRALINEQUALITIES} respectively.
See Remark~\ref{R:BOXEDCONSTANTS}
regarding the boxed constants on RHS~\eqref{E:TOPORDERWAVEENERGYINTEGRALINEQUALITIES}.

\begin{proposition}[\textbf{Integral inequalities for the wave variable controlling quantities}]
	\label{P:WAVEENERGYINTEGRALINEQUALITIES}
	Consider the $L^2$-controlling quantities 
	$\lbrace \totmax{N}(t,u) \rbrace_{N=1,\cdots,20}$,
	$\lbrace \easytotmax{N}(t,u) \rbrace_{N=1,\cdots,20}$,
	$\lbrace \Vorttotmax{N}(t,u) \rbrace_{N=0,\cdots,21}$,
	and 
	$\lbrace \coercivespacetimemax{N}(t,u) \rbrace_{N=1,\cdots,20}$
	from Defs.~\ref{D:MAINCOERCIVEQUANT}
	and \ref{D:COERCIVEINTEGRAL}.
	Assume that $N =20$ and $\varsigma > 0$.
	There exist constants $C > 0$ and $C_{\ast} > 0$,
	independent of $\varsigma$, such that
	under the data-size and bootstrap assumptions 
	of Sects.~\ref{SS:FLUIDVARIABLEDATAASSUMPTIONS}-\ref{SS:PSIBOOTSTRAP}
	and the smallness assumptions of Sect.~\ref{SS:SMALLNESSASSUMPTIONS}, 
	the following estimates hold
	for $(t,u) \in [0,\Tboot) \times [0,U_0]$:
\begingroup
\allowdisplaybreaks
	\begin{subequations}
	\begin{align} \label{E:TOPORDERWAVEENERGYINTEGRALINEQUALITIES}
		& 
		\max
		\left\lbrace
			\totmax{N}(t,u),
			\coercivespacetimemax{N}(t,u)
		\right\rbrace
		\\
		& \leq C (1 + \varsigma^{-1}) \mathring{\upepsilon}^2 \upmu_{\star}^{-3/2}(t,u)
				\notag \\
		& \ \ 
			+ 
			\boxed{6}
			\int_{t'=0}^t
					\frac{\left\| [\Lunit \upmu]_- \right\|_{L^{\infty}(\Sigma_{t'}^u)}} 
							 {\upmu_{\star}(t',u)} 
				  \totmax{N}(t',u)
				\, dt'
				\notag \\
		& \ \
			+ 
			\boxed{8.1}
			\int_{t'=0}^t
				\frac{\left\| [\Lunit \upmu]_- \right\|_{L^{\infty}(\Sigma_{t'}^u)}} 
								 {\upmu_{\star}(t',u)} 
						\sqrt{\totmax{N}}(t',u) 
						\int_{s=0}^{t'}
							\frac{\left\| [\Lunit \upmu]_- \right\|_{L^{\infty}(\Sigma_s^u)}} 
									{\upmu_{\star}(s,u)} 
							\sqrt{\totmax{N}}(s,u) 
						\, ds
				\, dt'
			\notag	\\
		& \ \
			+ 
			\boxed{2}
			\frac{1}{\upmu_{\star}^{1/2}(t,u)}
			\sqrt{\totmax{N}}(t,u)
			\left\| 
				\Lunit \upmu 
			\right\|_{L^{\infty}(\Sigmaminus{t}{t}{u})}
			\int_{t'=0}^t
				\frac{1}{\upmu_{\star}^{1/2}(t',u)} \sqrt{\totmax{N}}(t',u)
			\, dt'
			\notag \\
		& \ \ 
			+ 
			C_{\ast}
			\int_{t'=0}^t
					\frac{1} 
							 {\upmu_{\star}(t',u)} 
				  \sqrt{\totmax{N}}(t',u) 
					\sqrt{\easytotmax{N}}(t',u) 
				\, dt'
				\notag \\
		& \ \
			+ 
			C_{\ast}
			\int_{t'=0}^t
				\frac{1}
					{\upmu_{\star}(t',u)} 
						\sqrt{\totmax{N}}(t',u) 
						\int_{s=0}^{t'}
							\frac{1} 
									{\upmu_{\star}(s,u)} 
							\sqrt{\easytotmax{N}}(s,u) 
						\, ds
				\, dt'
			\notag	\\
		& \ \
			+ 
			C_{\ast}
			\frac{1}{\upmu_{\star}^{1/2}(t,u)}
			\sqrt{\totmax{N}}(t,u)
			\int_{t'=0}^t
				\frac{1}{\upmu_{\star}^{1/2}(t',u)} \sqrt{\easytotmax{N}}(t',u)
			\, dt'
			\notag \\
		& \ \
			+ 
			C \varepsilon
			\int_{t'=0}^t
					\frac{1} 
						{\upmu_{\star}(t',u)} 
				  \totmax{N}(t',u)
				\, dt'
				\notag
				\\
		& \ \
			+ 
			C \varepsilon
			\int_{t'=0}^t
				\frac{1} 
						{\upmu_{\star}(t',u)} 
						\sqrt{\totmax{N}}(t',u) 
						\int_{s=0}^{t'}
							\frac{1} 
									{\upmu_{\star}(s,u)} 
							\sqrt{\totmax{N}}(s,u) 
						\, ds
				\, dt'
			\notag	\\
		& \ \
			+ 
			C \varepsilon
			\frac{1}{\upmu_{\star}^{1/2}(t,u)}
			\sqrt{\totmax{N}}(t,u)
			\int_{t'=0}^t
				\frac{1}{\upmu_{\star}^{1/2}(t',u)} \sqrt{\totmax{N}}(t',u)
			\, dt'
				\notag \\
		& \ \
			+ 
			C
			\sqrt{\totmax{N}}(t,u)
			\int_{t'=0}^t
				\frac{1}{\upmu_{\star}^{1/2}(t',u)} \sqrt{\totmax{N}}(t',u)
			\, dt'
				\notag \\
		& \ \
			+
			C
			\int_{t'=0}^t
				\frac{1}{\sqrt{\Tboot - t'}} \totmax{N}(t',u)
			 \, dt'
		+ C (1 + \varsigma^{-1})
					\int_{t'=0}^t
					\frac{1} 
							 {\upmu_{\star}^{1/2}(t',u)} 
				  \totmax{N}(t',u)
				\, dt'
				\notag	\\
		& \ \
			+ C
					\int_{t'=0}^t
					\frac{1} 
							 {\upmu_{\star}(t',u)} 
				  \sqrt{\totmax{N}}(t',u)
				  \int_{s = 0}^{t'}
				  	\frac{1} 
							 {\upmu_{\star}^{1/2}(s,u)} 
							 \sqrt{\totmax{N}}(s,u)
				  \, ds
				\, dt'
				\notag	\\
		& \ \
			+ C
					\int_{t'=0}^t
					\frac{1} 
							 {\upmu_{\star}(t',u)} 
				  \sqrt{\totmax{N}}(t',u)
				  \int_{s = 0}^{t'}
				  	\frac{1}{\upmu_{\star}(s,u)}
				  	\int_{s' = 0}^s
				  	\frac{1} 
							 {\upmu_{\star}^{1/2}(s',u)} 
							 \sqrt{\totmax{N}}(s',u)
						\, ds'	 
				  \, ds
				\, dt'
				\notag	\\
		& \ \
			+ C (1 + \varsigma^{-1})
					\int_{u'=0}^u
						\totmax{N}(t,u')
					\, du'
			+ C \varepsilon \totmax{N}(t,u)
			+ C \varsigma \totmax{N}(t,u)
			+ C \varsigma \coercivespacetimemax{N}(t,u)
				\notag \\
		& \ \
			+ C
				(1 + \varsigma^{-1})
				\int_{t'=0}^t
					\frac{1} 
							{\upmu_{\star}^{5/2}(t',u)} 
				  \totmax{[1,N-1]}(t',u)
				\, dt'
			+ 
			C (1 + \varsigma^{-1})
					\int_{u'=0}^u
						\totmax{[1,N-1]}(t,u')
					\, du'
					\notag \\
		& \ \
			+
			C \varepsilon \totmax{[1,N-1]}(t,u)
			+ 
			C \varsigma \totmax{[1,N-1]}(t,u)
			+ 
			C \varsigma \coercivespacetimemax{[1,N-1]}(t,u)
				\notag 
				\\
		& \ \ 
			+
		C 
		\int_{t'=0}^t
			\frac{1}{\upmu_{\star}^{3/2}(t',u)}
			\left\lbrace
				\int_{s=0}^{t'}
					\sqrt{\VorttotTanmax{N+1}}(s,u)
				\, ds
			\right\rbrace^2
		\, dt'
				\notag \\
		& \ \
			+
			C 
			\int_{t'=0}^t
				\frac{1}{\upmu_{\star}^{3/2}(t',u)}
				\left\lbrace
					\int_{s=0}^{t'}
						\frac{1} 
						{\upmu_{\star}^{1/2}(s,u)}
					\sqrt{\VorttotTanmax{\leq N}}(s,u)
					\, ds
				\right\rbrace^2
			\, dt'
			\notag \\
		& \ \
			+ C
				\int_{t'=0}^t
					\VorttotTanmax{\leq N+1}(t',u)
				\, dt'
			+ C
				\int_{u'=0}^u
					\VorttotTanmax{\leq N}(t,u')
				\, du',
				\notag
	\end{align}
	\end{subequations}
	Moreover,\footnote{Instead of the ``large coefficient'' terms, one only has corresponding ``small coefficient'' terms 
	\begin{align}
	& C \varepsilon
			\int_{t'=0}^t
				\frac{1} {\upmu_{\star}(t',u)} 
						\sqrt{\totmax{N}}(t',u) 
						\int_{s=0}^{t'}
							\frac{1}{\upmu_{\star}(s,u)} 
							\sqrt{\totmax{N}}(s,u) 
						\, ds
				\, dt'
			\notag	\\
		& \ \
			+ 
			C \varepsilon
			\int_{t'=0}^t
					\frac{1} 
						{\upmu_{\star}(t',u)} 
				  \totmax{N}(t',u)
				\, dt'
				\notag
				\\
		& \ \
			+ 
			C \varepsilon
			\frac{1}{\upmu_{\star}^{1/2}(t,u)}
			\sqrt{\totmax{N}}(t,u)
			\int_{t'=0}^t
				\frac{1}{\upmu_{\star}^{1/2}(t',u)} \sqrt{\totmax{N}}(t',u)
			\, dt',
				\notag 
	\end{align}
	which are also featured on RHS~\eqref{E:TOPORDERWAVEENERGYINTEGRALINEQUALITIES}.} 

	\begin{align} \label{E:EASYVARIABLESTOPORDERWAVEENERGYINTEGRALINEQUALITIES}
		& \mbox{inequality \eqref{E:TOPORDERWAVEENERGYINTEGRALINEQUALITIES} holds
		with the LHS replaced with \ }
			\\
	& 
	\max
	\left\lbrace
		\easytotmax{N}(t,u),
		\easycoercivespacetimemax{N}(t,u)
	\right\rbrace
			\notag \\
	& \mbox{and \underline{without} 
		the six ``large-coefficient'' terms \ }
		\boxed{6} \cdots, 
		\boxed{2} \cdots,
		\boxed{8.1} \cdots, 
		C_{\ast} \cdots
			\notag \\
	&	\mbox{on the RHS.} \notag 
	\end{align}
	\endgroup
	In addition if $2 \leq N \leq 20$ and $\varsigma > 0$,
	then
	\begin{align} \label{E:BELOWTOPORDERWAVEENERGYINTEGRALINEQUALITIES}
		& 
		\max
		\left\lbrace
			\totmax{[1,N-1]}(t,u),
			\coercivespacetimemax{[1,N-1]}(t,u)
		\right\rbrace
		\\
		& \leq C \mathring{\upepsilon}^2 
			\notag \\
		& \ \
			+
			C
			\int_{t'=0}^t
				\frac{1}{\upmu_{\star}^{1/2}(t',u)} 
						\sqrt{\totmax{[1,N-1]}}(t',u) 
						\int_{s=0}^{t'}
							\frac{1}{\upmu_{\star}^{1/2}(s,u)} 
							\sqrt{\totmax{N}}(s,u) 
						\, ds
				\, dt'
			\notag	\\
		& \ \
			+
			C
			\int_{t'=0}^t
				\frac{1}{\sqrt{\Tboot - t'}} \totmax{[1,N-1]}(t',u)
			 \, dt'
			 \notag \\
		& \ \
			+ C (1 + \varsigma^{-1})
					\int_{t'=0}^t
					\frac{1} 
							 {\upmu_{\star}^{1/2}(t',u)} 
				  \totmax{[1,N-1]}(t',u)
				\, dt'
				\notag \\
	 & \ \
			+ C (1 + \varsigma^{-1})
					\int_{u'=0}^u
						\totmax{[1,N-1]}(t,u')
					\, du'
				\notag	\\
		& \ \
			+ C \varsigma \coercivespacetimemax{[1,N-1]}(t,u)
				\notag
				\\
		& \ \
			+ C
				\int_{t'=0}^t
					\VorttotTanmax{\leq N}(t',u)
				\, dt'
			+ C
				\int_{u'=0}^u
					\VorttotTanmax{\leq N-1}(t,u')
				\, du'.
				\notag
	\end{align}

\end{proposition}

\begin{proposition}[\textbf{Integral inequalities for the specific vorticity-controlling quantities}]
	\label{P:VORTICITYENERGYINTEGRALINEQUALITIES}
	Consider the $L^2$-controlling quantities 
	$\lbrace \totmax{N}(t,u) \rbrace_{N=1,\cdots,20}$,
	$\lbrace \easytotmax{N}(t,u) \rbrace_{N=1,\cdots,20}$,
	$\lbrace \Vorttotmax{N}(t,u) \rbrace_{N=0,\cdots,21}$,
	and 
	$\lbrace \coercivespacetimemax{N}(t,u) \rbrace_{N=1,\cdots,20}$
	from Defs.~\ref{D:MAINCOERCIVEQUANT}
	and \ref{D:COERCIVEINTEGRAL}.
	Assume that $N=20$.
	There exists a constant $C > 0$ such that
	under the data-size and bootstrap assumptions 
	of Sects.~\ref{SS:FLUIDVARIABLEDATAASSUMPTIONS}-\ref{SS:PSIBOOTSTRAP}
	and the smallness assumptions of Sect.~\ref{SS:SMALLNESSASSUMPTIONS}, 
	the following estimates hold
	for $(t,u) \in [0,\Tboot) \times [0,U_0]$:
\begin{subequations}
\begin{align}
	\VorttotTanmax{N+1}(t,u)
	& \leq
	C
	\mathring{\upepsilon}^2
	\frac{1}
	{\upmu_{\star}^{3/2}(t,u)}
	+
	C
	\varepsilon^2 
	\frac{1}{\upmu_{\star}(t,u)}
	\totmax{N}(t,u)
		\label{E:TOPORDERVORTICITYENERGYINTEGRALINEQUALITIES} \\
 & \ \
	+
	C
	\varepsilon^2
	\int_{t'=0}^t
				\frac{1}{\upmu_{\star}^2(t',u)} 
				\left\lbrace
					\int_{s=0}^t
						\frac{1}{\upmu_{\star}(s,u)} 
						\sqrt{\totmax{N}}(s,u)
					\, ds
				\right\rbrace^2
			\, dt'
		\notag \\
	& \ \
			+ 
			C \varepsilon^2
			\frac{1}{\upmu_{\star}^2(t,u)}
			\totmax{[1,N-1]}(t,u)
			+
			C \varepsilon^2 \coercivespacetimemax{[1,N]}(t,u)
				\notag \\
	& \ \
		+
		C \varepsilon^2
		\int_{t'=0}^t
			\frac{1}{\upmu_{\star}^2(t',u)}
			\left\lbrace
				\int_{s=0}^t
					\sqrt{\VorttotTanmax{N+1}}(s,u)
				\, ds
			\right\rbrace^2
		\, dt'
			\notag \\
	& \ \
		+
		C \varepsilon^2
		\int_{t'=0}^t
			\frac{1}{\upmu_{\star}^2(t',u)}
			\left\lbrace
				\int_{s=0}^t
					\frac{1}{\upmu_{\star}^{1/2}(s,u)}
					\sqrt{\VorttotTanmax{[1,N]}}(s,u)
				\, ds
			\right\rbrace^2
		\, dt'
		\notag \\
	& \ \
		+
		C
		\int_{u'=0}^u
			\VorttotTanmax{\leq N+1}(t,u')
		\, du'.
			\notag
\end{align}
Similarly, if $N \leq 20$, then
\begin{align}
	\VorttotTanmax{\leq N}(t,u)
	& \leq 
	C \mathring{\upepsilon}^2
	+
	C \varepsilon^2
		\int_{t'=0}^t
			\left\lbrace
				\int_{s=0}^{t'}
					\frac{1}{\upmu_{\star}^{1/2}(s,u)}
					\sqrt{\totmax{[1,N]}}(s,u)
				\, ds
			\right\rbrace^2
		\, dt'
	\label{E:BELOWTOPORDERVORTICITYENERGYINTEGRALINEQUALITIES} 
		\\
&  \ \
	+
	\underbrace{
	C \varepsilon^2 \totmax{[1,N-1]}(t,u)
	+ 
	C \varepsilon^2 \coercivespacetimemax{[1,N-1]}(t,u)
	}_{\mbox{Absent if $N=0$}}
		\notag \\
& \ \
	+
	C
	\int_{u'=0}^u
		\VorttotTanmax{\leq N}(t,u')
	\, du'.
			\notag 
\end{align}
\end{subequations}

\end{proposition}

\subsection{Bootstrap assumptions for the fundamental \texorpdfstring{$L^2$-}{square integral}
controlling quantities of the wave variables}
\label{SS:L2BOOTSTRAP}
	To facilitate our proof of Prop.~\ref{P:MAINAPRIORIENERGY},
	it is convenient to make $L^2$-type bootstrap assumptions,
	which we state in this section.
	Specifically, let
	$\lbrace \totmax{N}(t,u) \rbrace_{N=1,\cdots,20}$,
	$\lbrace \Vorttotmax{N}(t,u) \rbrace_{N=0,\cdots,21}$,
	and 
	$\lbrace \coercivespacetimemax{N}(t,u) \rbrace_{N=1,\cdots,20}$
	be the $L^2$-controlling quantities
	from Defs.~\ref{D:MAINCOERCIVEQUANT}
	and \ref{D:COERCIVEINTEGRAL}.
	We assume that the following inequalities hold
	for $(t,u) \in [0,\Tboot) \times [0,U_0]$,
	where $\varepsilon$ is the small bootstrap parameter appearing in 
	Sect.~\ref{SS:PSIBOOTSTRAP}:
	\begin{subequations}
	\begin{align}
		\sqrt{\totmax{15+M}}(t,u)
		+ \sqrt{\coercivespacetimemax{15+M}}(t,u)
		& \leq \sqrt{\varepsilon} \upmu_{\star}^{-(M+.9)}(t,u),
			&& (0 \leq M \leq 5),
				\label{E:BAMULOSSMAINAPRIORIENERGYESTIMATES} \\
		\sqrt{\totmax{[1,14]}}(t,u)
		+ \sqrt{\coercivespacetimemax{[1,14]}}(t,u)
		& \leq \sqrt{\varepsilon},
		&& 
			\label{E:BANOMULOSSMAINAPRIORIENERGYESTIMATES}
				\\
		\sqrt{\VorttotTanmax{21}}(t,u)
		& \leq \sqrt{\varepsilon} \upmu_{\star}^{-6.4}(t,u),
		 && \label{E:BATOPVORTMULOSSMAINAPRIORIENERGYESTIMATES} \\
		\sqrt{\VorttotTanmax{16+M}}(t,u)
		& \leq \sqrt{\varepsilon} \upmu_{\star}^{-(M+.9)}(t,u),
			&& (0 \leq M \leq 4),
				\label{E:BAVORTMULOSSMAINAPRIORIENERGYESTIMATES} \\
		\sqrt{\VorttotTanmax{\leq 15}}(t,u)
		& \leq \sqrt{\varepsilon}.
		&&
			\label{E:BAVORTNOMULOSSMAINAPRIORIENERGYESTIMATES}
		\end{align}
	\end{subequations}

\subsection{Preliminary \texorpdfstring{$L^2$}{square integral} estimates for the eikonal function quantities that do not require modified quantities}
\label{SS:L2FOREIKONALNOMODIFIED}
In Lemma~\ref{L:EASYL2BOUNDSFOREIKONALFUNCTIONQUANTITIES},
we derive a priori estimates for the below-top-order
derivatives of the eikonal function quantities
and, in the case that at least one $\Lunit$-differentiation is involved,
their top-order derivatives. These estimates are simple consequence of the transport inequalities
derived in Prop.~\ref{P:IMPROVEMENTOFAUX} and can be derived 
without using the modified quantities.

We start with a simple commutator lemma.

\begin{lemma}[\textbf{Simple commutator lemma}]
	\label{L:PUTALLRADFACTORSUPFRONT}
	Assume that $1 \leq N \leq 21$
	and $0 \leq M \leq \min \lbrace 2, N-1 \rbrace$.
	Under the data-size and bootstrap assumptions 
	of Sects.~\ref{SS:FLUIDVARIABLEDATAASSUMPTIONS}-\ref{SS:PSIBOOTSTRAP}
	and the smallness assumptions of Sect.~\ref{SS:SMALLNESSASSUMPTIONS}, 
	the following pointwise estimates hold
	on $\mathcal{M}_{\Tboot,U_0}$: 
	\begin{align} \label{E:PUTALLRADFACTORSUPFRONT}
	\left|
		\Fullset_{\ast}^{[1,N];\leq M} \threePsi
	\right|
	& 
	\lesssim
	\sum_{\widetilde{M}=0}^M 
	\left|
		\Rad^{\widetilde{M}} \Tanset^{[1,N-\widetilde{M}]} \threePsi
	\right|
	+
	\varepsilon 
	\left|
		\myarray
			[\Fullset_{\ast \ast}^{[1,N-1];(M-1)_+} \BadVar]
			{\Fullset_{\ast}^{[1,N-1];\leq M} \GdVar}
	\right|,
	\end{align}
	where  $(M-1)_+ = \max \lbrace 0, M-1 \rbrace$
	and when $N=1$,
	only
	$
	\left|
		\Tanset \threePsi
	\right|
	$
	appears on RHS~\eqref{E:PUTALLRADFACTORSUPFRONT}.
\end{lemma}

\begin{proof}
		We repeatedly use
		\eqref{E:DETAILEDPERMUTEDVECTORFIELDSACTINGONFUNCTIONSCOMMUTATORESTIMATE}
		with $f = \threePsi$ and the $L^{\infty}$ estimates of Prop.~\ref{P:IMPROVEMENTOFAUX}
		to commute the factors of $\Rad$ acting on $\threePsi$
		on LHS~\eqref{E:PUTALLRADFACTORSUPFRONT}
		to the front (so they are the last to hit $\threePsi$).
\end{proof}

We now provide the main estimates of this section.

\begin{lemma}[$L^2$ \textbf{bounds for the eikonal function quantities that do not require modified quantities}]
	\label{L:EASYL2BOUNDSFOREIKONALFUNCTIONQUANTITIES}
	Assume that $1 \leq N \leq 20$.
	Under the data-size and bootstrap assumptions 
	of Sects.~\ref{SS:FLUIDVARIABLEDATAASSUMPTIONS}-\ref{SS:PSIBOOTSTRAP}
	and the smallness assumptions of Sect.~\ref{SS:SMALLNESSASSUMPTIONS}, 
	the following $L^2$ estimates hold for $(t,u) \in [0,\Tboot) \times [0,U_0]$
	(see Sect.~\ref{SS:STRINGSOFCOMMUTATIONVECTORFIELDS} regarding the vectorfield operator notation):
	\begin{subequations}
	\begin{align}
		\threemyarray
		[
			\left\|
				\Lunit \Fullset_{\ast}^{[1,N];\leq 1} \upmu
			\right\|_{L^2(\Sigma_t^u)}
		]
		{
			\left\|
				\Lunit \Fullset_{\ast}^{\leq N;\leq 2} \Lunit_{(Small)}^i
			\right\|_{L^2(\Sigma_t^u)}
		}
		{
			\left\|
				\Lunit \Fullset^{\leq N-1;\leq 2} \mytr \upchi
			\right\|_{L^2(\Sigma_t^u)}
		}
		& \lesssim 
				\mathring{\upepsilon}
				+
				\frac{\sqrt{\totmax{[1,N]}}(t,u)}{\upmu_{\star}^{1/2}(t,u)},
			\label{E:LUNITSTARDIFFERENTIATEDEIKONALINTERMSOFCONTROLLING}
				 \\
		\threemyarray
		[
			\left\|
				\Fullset_{\ast \ast}^{[1,N];\leq 1} \upmu
			\right\|_{L^2(\Sigma_t^u)}
		]
		{
			\left\|
				\Fullset_{\ast}^{[1,N];\leq 2} \Lunit_{(Small)}^i
			\right\|_{L^2(\Sigma_t^u)}
		}
		{
			\left\|
				\Fullset^{\leq N-1;\leq 2} \mytr \upchi
			\right\|_{L^2(\Sigma_t^u)}
		}
		& \lesssim 
			\mathring{\upepsilon}
			+ 
			\int_{s=0}^t
				\frac{\sqrt{\totmax{[1,N]}}(s,u)}{\upmu_{\star}^{1/2}(s,u)}
			\, ds.
				 \label{E:TANGENGITALEIKONALINTERMSOFCONTROLLING}
	\end{align}
	\end{subequations}

\end{lemma}

\begin{proof}
	See Sect.~\ref{SS:OFTENUSEDESTIMATES} for some comments on the analysis.
	We set $q_N(t) := \mbox{LHS~\eqref{E:TANGENGITALEIKONALINTERMSOFCONTROLLING}}$.
	From \eqref{E:LLUNITUPMUFULLSETSTARCOMMUTEDPOINTWISEESTIMATE},
	\eqref{E:LUNITCOMMUTEDLUNITSMALLIPOINTWISE},
	Lemma~\ref{L:COERCIVENESSOFCONTROLLING},
	\eqref{E:PUTALLRADFACTORSUPFRONT},
	\eqref{E:SIGMATVOLUMEFORMCOMPARISON},
	and Lemma~\ref{L:L2NORMSOFTIMEINTEGRATEDFUNCTIONS},
	we deduce 
	\begin{align} \label{E:GRONWALLREADYEIKONALFUNCTIONQUANTITIES}
		q_N(t)
		& \leq 
		C q_N(0)
		+ C 
			\int_{s=0}^t
				q_N(s)
			\, ds
			+
			C
			\int_{s=0}^t
				\frac{\sqrt{\totmax{[1,N]}}(s,u)}{\upmu_{\star}^{1/2}(s,u)}
			\, ds.
	\end{align}
	Next, we note that $q_N(0) \lesssim \mathring{\upepsilon}$,
	an estimate that follows from 
	the estimate 
	\eqref{E:POINTWISEESTIMATESFORCHIANDITSDERIVATIVES} for $\mytr \upchi$
	and our data-size assumptions.
	We now apply Gronwall's inequality to \eqref{E:GRONWALLREADYEIKONALFUNCTIONQUANTITIES}
	to conclude that $q_N(t) \lesssim$ RHS \eqref{E:TANGENGITALEIKONALINTERMSOFCONTROLLING} as desired.
	We have thus proved \eqref{E:TANGENGITALEIKONALINTERMSOFCONTROLLING}.

	To obtain the estimates \eqref{E:LUNITSTARDIFFERENTIATEDEIKONALINTERMSOFCONTROLLING},
	we take the norm 
	$\left\|
	\cdot
	\right\|_{L^2(\Sigma_t^u)} 
	$
	of the inequalities 
	\eqref{E:LLUNITUPMUFULLSETSTARCOMMUTEDPOINTWISEESTIMATE},
	\eqref{E:LUNITCOMMUTEDLUNITSMALLIPOINTWISE}
	and argue as above using the already proven estimates
	\eqref{E:TANGENGITALEIKONALINTERMSOFCONTROLLING}.
	In these estimates, we encounter the integrals
	$
	\displaystyle
	\int_{s=0}^t
		\frac{\sqrt{\totmax{[1,N]}}(s,u)}{\upmu_{\star}^{1/2}(s,u)}
	\, ds
	$,
	which we (inefficiently) bound by 
	$
	\lesssim \sqrt{\totmax{[1,N]}}(t,u) 
	\leq \upmu_{\star}^{-1/2}(t,u) \sqrt{\totmax{[1,N]}}(t,u)
	$
	with the help of inequality \eqref{E:LESSSINGULARTERMSMPOINTNINEINTEGRALBOUND}.

\end{proof}

\subsection{Estimates for the easiest error integrals}
\label{SS:ENERGYESTIMATESEASIESTERRORINTEGRALS}
In this section, we derive estimates for the easiest
error integrals that we encounter in our energy estimates
for $\threePsi$ and $\Vortrenormalized$. These
error integrals do not contribute to the blowup
featured in our high-order energy estimates.

We start with a lemma relevant for bounding the specific vorticity. 
Specifically, we bound the 
error integrals corresponding to 
the last integral on RHS~\eqref{E:ENERGYIDENTITYRENORMALIZEDVORTICITY}.

\begin{lemma}[\textbf{The simplest transport equation error integrals}]
	\label{L:SIMPLESTTRANSPORTEQUATIONERROR}
		Assume that $N \leq 21$.
		Under the data-size and bootstrap assumptions 
	of Sects.~\ref{SS:FLUIDVARIABLEDATAASSUMPTIONS}-\ref{SS:PSIBOOTSTRAP}
	and the smallness assumptions of Sect.~\ref{SS:SMALLNESSASSUMPTIONS}, 
	the following integral estimates hold for $(t,u) \in [0,\Tboot) \times [0,U_0]$
	(see Sect.~\ref{SS:STRINGSOFCOMMUTATIONVECTORFIELDS} regarding the vectorfield operator notation):
		\begin{align} \label{E:SIMPLESTTRANSPORTEQUATIONERROR}
			\left|
				\int_{\mathcal{M}_{t,u}}
					\left\lbrace
						\Lunit \upmu
						+
						\upmu \mytr \angk
					\right\rbrace
					(\Tanset^{N} \Vortrenormalized)^2 
				\, d \vol
			\right|
			& 
			\lesssim 
			\int_{u'=0}^u
				\VorttotTanmax{N}(t,u')
			\, du'.
		\end{align}
\end{lemma}
\begin{proof}
	From 
	equation \eqref{E:ANGKDECOMPOSED},
	Lemma~\ref{L:SCHEMATICDEPENDENCEOFMANYTENSORFIELDS},
	Lemma~\ref{L:POINTWISEFORRECTANGULARCOMPONENTSOFVECTORFIELDS},
	and the $L^{\infty}$ estimates of Prop.~\ref{P:IMPROVEMENTOFAUX},
	we deduce
	$\| \Lunit \upmu + \upmu \mytr \angk \|_{L^{\infty}(\Sigma_t^u)}
	\lesssim 1
	$.
	Hence, using Lemma~\ref{L:COERCIVENESSOFCONTROLLING}, we conclude that
	\begin{align}
		\mbox{LHS~\eqref{E:SIMPLESTTRANSPORTEQUATIONERROR}}
		\lesssim
			\int_{\mathcal{M}_{t,u}}
				(\Tanset^{N} \Vortrenormalized)^2
			\, d \vol
		& =
		\int_{u'=0}^u
			\| \Tanset^{N} \Vortrenormalized \|_{L^2(\mathcal{P}_{u'}^t)}^2
		\, du'
		\lesssim
		\int_{u'=0}^u
			\VorttotTanmax{N}(t,u')
		\, du'.
	\end{align}
\end{proof}

The next lemma is a more complicated analog of Lemma~\ref{L:SIMPLESTTRANSPORTEQUATIONERROR}
for the wave variables
$\lbrace \Densrenormalized  - v^1,v^1,v^2 \rbrace$.
We recall that in \eqref{E:MULTERRORINT}, we decomposed
the integrand corresponding
to the last term on the right-hand side 
of the wave equation energy identity \eqref{E:E0DIVID}.
Moreover, we recall that one of the integrand pieces
is coercive and is critically important
for controlling geometric torus derivatives; 
we isolated it in Def.~\ref{D:COERCIVEINTEGRAL}.
In the next lemma,
we bound the error integrals 
corresponding to the remaining terms in the
decomposition \eqref{E:MULTERRORINT}.

\begin{lemma}[\textbf{Error integrals involving the deformation tensor of the multiplier vectorfield}]
	\label{L:MULTIPLIERVECTORFIELDERRORINTEGRALS}
	Let $\Psi \in \lbrace \Densrenormalized  - v^1,v^1,v^2 \rbrace$.
	Assume that $1 \leq N \leq 20$ and $\varsigma > 0$.
	Let $\basicenergyerrorarg{\Mult}{i}[\Fullset_{\ast}^{N;\leq 1} \Psi]$,
	$(i=1,\cdots,5)$,
	be the quantities defined by \eqref{E:MULTERRORINTEG1}-\eqref{E:MULTERRORINTEG5}
	(with $\Fullset_{\ast}^{N;\leq 1} \Psi$ in the role of $\Psi$).
	Under the data-size and bootstrap assumptions 
	of Sects.~\ref{SS:FLUIDVARIABLEDATAASSUMPTIONS}-\ref{SS:PSIBOOTSTRAP}
	and the smallness assumptions of Sect.~\ref{SS:SMALLNESSASSUMPTIONS}, 
	the following integral estimates hold for $(t,u) \in [0,\Tboot) \times [0,U_0]$,
	where the implicit constants are independent of $\varsigma$
	(and without any absolute value taken on the left):
	\begin{align}
		\int_{\mathcal{M}_{t,u}}
			\sum_{i=1}^5 \basicenergyerrorarg{\Mult}{i}[\Fullset_{\ast}^{N;\leq 1} \Psi]
		\, d \vol
		& \lesssim
		 	\int_{t'=0}^t
				\frac{1}{\sqrt{\Tboot - t'}} \totmax{[1,N]}(t',u)
			 \, dt'
			 \\
		& \ \
			+
		 	(1 + \varsigma^{-1})
		 		\int_{t'=0}^t
					\totmax{[1,N]}(t',u)
				\, dt'
			 \notag	\\
		 & \ \ 
		 	+
			(1 + \varsigma^{-1})
			\int_{u'=0}^u
				\totmax{[1,N]}(t,u')
			\, du'
			+ \varsigma 
				\coercivespacetimemax{[1,N]}(t,u).
				\notag
	\end{align}
\end{lemma}

\begin{proof}
	We integrate inequality \eqref{E:MULTIPLIERVECTORFIEDERRORTERMPOINTWISEBOUND} 
	(with $\Fullset_{\ast}^{N;\leq 1} \Psi$ in the role of $\Psi$)
	over the domain $\mathcal{M}_{t,u}$ and use 
	Lemmas~\ref{L:KEYSPACETIMECOERCIVITY} and \ref{L:COERCIVENESSOFCONTROLLING}.
\end{proof}

In the next lemma, we bound
$
\left\|
	\Fullset^{N;1} \Vortrenormalized
\right\|_{L^2(\Sigma_{t,u})}
$
in terms of the fundamental $L^2$-controlling quantities.
These estimates play a role in bounding some of the inhomogeneous
terms in the wave equations that depend on the derivatives of
the vorticity up to top-order.

\begin{lemma}[$L^2$ \textbf{estimates involving one transversal derivative of the specific vorticity}]
\label{L:L2NORMOFVORTICITYTRANSVERSALDERIVATIVESINTERMSOFTANGENTIALDERIVATIVES}
	Assume that $1 \leq N \leq 21$.
	Under the data-size and bootstrap assumptions 
	of Sects.~\ref{SS:FLUIDVARIABLEDATAASSUMPTIONS}-\ref{SS:PSIBOOTSTRAP}
	and the smallness assumptions of Sect.~\ref{SS:SMALLNESSASSUMPTIONS}, 
	the following $L^2$ estimates hold for $(t,u) \in [0,\Tboot) \times [0,U_0]$:
	\begin{align} \label{E:L2NORMOFVORTICITYTRANSVERSALDERIVATIVESINTERMSOFTANGENTIALDERIVATIVES}
		\left\|
			\Fullset^{N;1} \Vortrenormalized
		\right\|_{L^2(\Sigma_{t,u})}
		& 
		\lesssim 
		\mathring{\upepsilon}
		+
		\sqrt{\VorttotTanmax{\leq N}}(t,u)
		+
		\varepsilon
		\int_{s=0}^t
			\frac{\sqrt{\totmax{[1,N-1]}}(t,u)}{\upmu_{\star}^{1/2}(s,u)}
		\, ds,
	\end{align}
	where the last term on 
	RHS~\eqref{E:L2NORMOFVORTICITYTRANSVERSALDERIVATIVESINTERMSOFTANGENTIALDERIVATIVES}
	is absent when $N=1$.
	\end{lemma}

\begin{proof}
	We use
	\eqref{E:DETAILEDPERMUTEDVECTORFIELDSACTINGONFUNCTIONSCOMMUTATORESTIMATE}
	with $f = \Vortrenormalized$
	and the $L^{\infty}$ estimates of Prop.~\ref{P:IMPROVEMENTOFAUX}
		to commute the factors of $\Rad$ acting on $\Vortrenormalized$
		on LHS~\eqref{E:L2NORMOFVORTICITYTRANSVERSALDERIVATIVESINTERMSOFTANGENTIALDERIVATIVES}
		so that it is the first to hit $\Vortrenormalized$, 
		thereby obtaining
	\begin{align} \label{E:VORTICITYSIMPLYCOMMUTATOR}
	\left|
		\Fullset^{N;1} \Vortrenormalized
	\right|
	\lesssim
	\left|
		\Tanset^{\leq N-1} \Rad \Vortrenormalized
	\right|
	+
	\left|
		\Tanset^{\leq N-1} \Vortrenormalized
	\right|
	+
	\varepsilon 
	|\Fullset_{\ast}^{[1,N-1];\leq 1} \GdVar|
	+ 
	\varepsilon 
	|\Fullset_{\ast \ast}^{[1,N-1]} \BadVar|,
	\end{align}
	where the last three terms on 
	RHS~\eqref{E:VORTICITYSIMPLYCOMMUTATOR}
	are absent when $N=1$.
	Using Lemma~\ref{L:COERCIVENESSOFCONTROLLING},
	and the estimate \eqref{E:TANGENGITALEIKONALINTERMSOFCONTROLLING},
	we see that the 
	norms $\| \cdot \|_{L^2(\Sigma_{t,u})}$ of the last two terms
	on RHS~\eqref{E:VORTICITYSIMPLYCOMMUTATOR}
	are $\lesssim$ 
	RHS~\eqref{E:L2NORMOFVORTICITYTRANSVERSALDERIVATIVESINTERMSOFTANGENTIALDERIVATIVES}.
	Next, we use the fundamental theorem of calculus to obtain the following pointwise
	estimate for the second term on RHS~\eqref{E:VORTICITYSIMPLYCOMMUTATOR}:
	\begin{align} \label{E:FTCESTIMATEVORTICITYSIMPLYCOMMUTATOR}
	\left|
		\Tanset^{\leq N-1} \Vortrenormalized
	\right|(t,u,\vartheta)
	\lesssim
		\left|
			\Tanset^{\leq N-1} \Vortrenormalized
		\right|(0,u,\vartheta)
	+ \int_{s=0}^t
			\left|
				\Lunit \Tanset^{\leq N-1} \Vortrenormalized
			\right|(s,u,\vartheta)
		\, ds.
	\end{align}
	Using Lemma~\ref{L:L2NORMSOFTIMEINTEGRATEDFUNCTIONS}
	and Lemma~\ref{L:COERCIVENESSOFCONTROLLING},
	we see that the norm $\| \cdot \|_{L^2(\Sigma_{t,u})}$
	of the last term on RHS~\eqref{E:FTCESTIMATEVORTICITYSIMPLYCOMMUTATOR} is 
	$
	\displaystyle
	\lesssim
	\mathring{\upepsilon} 
	+
	\int_{s=0}^t
		\frac{\sqrt{\VorttotTanmax{[1,N]}}(s,u)}{\upmu_{\star}^{1/2}(s,u)}
	\, ds
	$.
	Manifestly, we have 
	$
	\mathring{\upepsilon} 
	\lesssim 
	$
	RHS~\eqref{E:L2NORMOFVORTICITYTRANSVERSALDERIVATIVESINTERMSOFTANGENTIALDERIVATIVES}.
	Moreover, using inequality \eqref{E:LESSSINGULARTERMSMPOINTNINEINTEGRALBOUND}
	and the fact that $\VorttotTanmax{[1,N]}$ is increasing in its arguments,
	we see that the previous time integral is 
	$\lesssim \sqrt{\VorttotTanmax{[1,N]}}(t,u)$, which is
	$\lesssim$ RHS~\eqref{E:L2NORMOFVORTICITYTRANSVERSALDERIVATIVESINTERMSOFTANGENTIALDERIVATIVES}
	as desired.
	To bound the norm $\| \cdot \|_{L^2(\Sigma_{t,u})}$
	of the first term on RHS~\eqref{E:FTCESTIMATEVORTICITYSIMPLYCOMMUTATOR},
	we use \eqref{E:SIGMATVOLUMEFORMCOMPARISON}
	with $s=0$ 
	and the small-data assumption \eqref{E:L2SMALLDATAASSUMPTIONSALONGSIGMA0}
	to obtain 
	$\| \Tanset^{\leq N-1} \Vortrenormalized(1,\cdot) \|_{L^2(\Sigma_{t,u})}
	\lesssim
	\| \Tanset^{\leq N-1} \Vortrenormalized \|_{L^2(\Sigma_{1,u})}
	\leq \mathring{\upepsilon}
	$
	as desired.

	It remains for us to bound the norm $\| \cdot \|_{L^2(\Sigma_{t,u})}$
	of the first term on RHS~\eqref{E:VORTICITYSIMPLYCOMMUTATOR}.
	Using equations \eqref{E:RENORMALIZEDVORTICTITYTRANSPORTEQUATION}
	and \eqref{E:TRANSPORTVECTORFIELDINTERMSOFLUNITANDRADUNIT}
	to algebraically express $\Rad \Vortrenormalized = - \upmu \Lunit \Vortrenormalized$
	and using the estimates 
	\eqref{E:UPTOONETRANSVERSALDERIVATIVEUPMULINFTY},
	\eqref{E:LUNITUPTOONETRANSVERSALUPMULINFINITY},
	\eqref{E:ZSTARSTARUPMULINFTY},
	and \eqref{E:VORTICITYUPTOTWOTRANSVERSALLINFTY},
	we find that
	$
	\left|
		\Tanset^{\leq N-1} \Rad \Vortrenormalized
	\right|
	\lesssim
	\left|
		\sqrt{\upmu} \Tanset^N \Vortrenormalized
	\right|
	+
	\left|
		\Tanset^{\leq N-1} \Vortrenormalized
	\right|
	+
	\varepsilon
	\left|
		\Fullset_{\ast \ast}^{[1,N-1]} \BadVar
	\right|
	$,
	where the last two terms on the RHS are absent when $N=1$.
	Lemma~\ref{L:COERCIVENESSOFCONTROLLING} immediately yields
	that 
	$
	\left\| 
		\sqrt{\upmu} \Tanset^N \Vortrenormalized 
	\right\|_{L^2(\Sigma_{t,u})}
	\lesssim \sqrt{\VorttotTanmax{\leq N}}(t,u)
	$,
	while the arguments given below \eqref{E:FTCESTIMATEVORTICITYSIMPLYCOMMUTATOR}
	imply that 
	$
	\left\| 
		\Tanset^{\leq N-1} \Vortrenormalized 
	\right\|_{L^2(\Sigma_{t,u})}
	\lesssim 
	\mathring{\upepsilon}
	+ 
	\sqrt{\VorttotTanmax{\leq N}}(t,u)
	$
	as desired. 
	Moreover, above we showed that
	$
	\varepsilon
	\left\|
		\Fullset_{\ast \ast}^{[1,N-1]} \BadVar
	\right\|_{L^2(\Sigma_{t,u})}
	\lesssim
	\mbox{RHS~\eqref{E:L2NORMOFVORTICITYTRANSVERSALDERIVATIVESINTERMSOFTANGENTIALDERIVATIVES}}.
	$
	This completes the proof of 
	\eqref{E:L2NORMOFVORTICITYTRANSVERSALDERIVATIVESINTERMSOFTANGENTIALDERIVATIVES}.
\end{proof}

In the next lemma, we bound the error integrals corresponding
to the $Harmless_{(Wave)}^{\leq N}$
and the 
$Harmless_{(Vort)}^{\leq N}$
terms.

\begin{lemma}[$L^2$ \textbf{bounds for error integrals involving} $Harmless_{(Wave)}^{\leq N}$ 
\textbf{or}
$Harmless_{(Vort)}^{\leq N}$
\textbf{terms}]
	\label{L:STANDARDPSISPACETIMEINTEGRALS}
	Let $\Psi \in \lbrace \Densrenormalized  - v^1,v^1,v^2 \rbrace$.
	Assume that $1 \leq N \leq 20$ and $\varsigma > 0$.
	Recall that the terms $Harmless_{(Wave)}^{\leq N}$ 
	and
	$Harmless_{(Vort)}^{\leq N}$
	are defined in
	Def.~\ref{D:HARMLESSTERMS}.
	Under the data-size and bootstrap assumptions 
	of Sects.~\ref{SS:FLUIDVARIABLEDATAASSUMPTIONS}-\ref{SS:PSIBOOTSTRAP}
	and the smallness assumptions of Sect.~\ref{SS:SMALLNESSASSUMPTIONS}, 
	the following integral estimates hold for $(t,u) \in [0,\Tboot) \times [0,U_0]$,
	where the implicit constants are independent of $\varsigma$:
	\begin{align}  \label{E:WAVESHARMLESSSPACETIMEINTEGRALS}
		\int_{\mathcal{M}_{t,u}}
		 	& \left|
					\myarray[(1 + \upmu) \Lunit \Fullset_{\ast}^{N;\leq 1} \Psi]
						{\Rad \Fullset_{\ast}^{N;\leq 1} \Psi}
				\right|
				\left|
					Harmless_{(Wave)}^{\leq N}
				\right|
		 \, d \vol
		   \\
		& \lesssim
		 	(1 + \varsigma^{-1})
			\int_{t'=0}^t
				\frac{\totmax{[1,N]}(t',u)}{\upmu_{\star}^{1/2}(t',u)}
			\, dt'
			+
			(1 + \varsigma^{-1})
			\int_{u'=0}^u
				\totmax{[1,N]}(t,u')
			\, du'
				\notag \\
		& \ \
			+ \varsigma
			  \coercivespacetimemax{[1,N]}(t,u)
			+
			\int_{t'=0}^t
				\Vorttotmax{\leq N}(t',u)
			\, dt'
			+
			\int_{u'=0}^u
				\VorttotTanmax{\leq N}(t,u')
			\, du'
			+ \mathring{\upepsilon}^2.
			\notag 
	\end{align}
	Moreover, if $N \leq 21$, then
	\begin{align}
		\int_{\mathcal{M}_{t,u}}
			&
			\left|
				\Tanset^{N} \Vortrenormalized
			\right|
			\left|
				Harmless_{(Vort)}^{\leq N}
			\right|
		\, d \vol
			\label{E:VORTICITYHARMLESSSPACETIMEINTEGRALS} \\
		& \lesssim
			\underbrace{
			\varepsilon^2 \totmax{[1,N-1]}(t,u)
			+ 
			\varepsilon^2 \coercivespacetimemax{[1,N-1]}(t,u)
			}_{\mbox{Absent if $N=0,1$}}
			+
			\int_{u'=0}^u
				\VorttotTanmax{\leq N}(t,u')
			\, du'
			+ 
			\varepsilon^2 \mathring{\upepsilon}^2.
			\notag
	\end{align}
\end{lemma}
\begin{proof}
	See Sect.~\ref{SS:OFTENUSEDESTIMATES} for some comments on the analysis.
	To prove 
	\eqref{E:WAVESHARMLESSSPACETIMEINTEGRALS}
	and \eqref{E:VORTICITYHARMLESSSPACETIMEINTEGRALS}
	we must estimate the spacetime integrals of various quadratic terms.
	We derive the desired estimates for five representative
	quadratic terms: four in the case of
	\eqref{E:WAVESHARMLESSSPACETIMEINTEGRALS}
	and one in the case of
	\eqref{E:VORTICITYHARMLESSSPACETIMEINTEGRALS}.
	The remaining terms can be bounded using similar or simpler arguments
	and we omit those details.
	As our first example, we bound the spacetime integral of
	$\left|
			\Lunit \Fullset_{\ast}^{N;\leq 1} \Psi
		\right|
	\left|
		\GeoAng^{N+1} \Psi
	\right|
	$
	(note that $\GeoAng^{N+1} \Psi = Harmless_{(Wave)}^{\leq N}$).
	Using
	spacetime Cauchy-Schwarz,
	Lemmas~\ref{L:KEYSPACETIMECOERCIVITY} and \ref{L:COERCIVENESSOFCONTROLLING},
	and simple estimates of the form
	$ab \lesssim a^2 + b^2$,
	and separately treating the regions
	$\lbrace \upmu \geq 1/4 \rbrace$ 
	and $\lbrace \upmu < 1/4 \rbrace$
	when bounding the integral of
	$\left|
		\GeoAng^{N+1} \Psi
	\right|^2$,
	we deduce the desired estimate as follows:
	\begin{align}
		& \int_{\mathcal{M}_{t,u}}
		 		\left|
					\Lunit \Fullset_{\ast}^{N;\leq 1} \Psi
				\right|
				\left|
					\GeoAng^{N+1} \Psi
				\right|
		 \, d \vol
		 	\label{E:FIRSTHARMLESSEXAMPLEINTEGRAL} \\
		 & \lesssim
		 \left\lbrace
		 \int_{\mathcal{M}_{t,u}}
		 		\left|
					\Lunit \Fullset_{\ast}^{N;\leq 1} \Psi
				\right|^2
			\, d \vol
			\right\rbrace^{1/2}
			\left\lbrace
			\int_{\mathcal{M}_{t,u}}
		 		\left|
					\GeoAng^{N+1} \Psi
				\right|^2
			\, d \vol
			\right\rbrace^{1/2}
				\notag \\
			& \lesssim
				(1 + \varsigma^{-1})
				\int_{u'=0}^u
					\int_{\mathcal{P}_{u'}^t}
						\left|
							\Lunit \Fullset_{\ast}^{N;\leq 1} \Psi
						\right|^2
					\, d \conevol
				\, du'
				\notag \\
		& \ \
				+
				\int_{u'=0}^u
					\int_{\mathcal{P}_{u'}^t}
						\upmu
						\left|
							\angdiff \GeoAng^N \Psi
						\right|^2
					\, d \conevol
				\, du'
				+
				\varsigma \TranminusdatasizeWithFactor
				\int_{\mathcal{M}_{t,u}}
				\textbf{1}_{\lbrace \upmu < 1/4 \rbrace}
		 		\left|
					\angdiff \GeoAng^N \Psi
				\right|^2
			\, d \vol
				\notag \\
			& \lesssim
				(1 + \varsigma^{-1})
				\int_{u'=0}^u
					\totmax{[1,N]}(t,u')
				\, du'
				+
				\varsigma
				\coercivespacetimemax{[1,N]}(t,u).
				\notag
	\end{align}

	As our second example, we bound the spacetime integral of
	$\left|
		\Lunit \Fullset_{\ast}^{N;\leq 1} \Psi
	 \right|
	 \left|
	 		\Fullset_{\ast \ast}^{[1,N];\leq 1} \upmu
	 \right|$.
	Using spacetime Cauchy-Schwarz,
	Lemmas~\ref{L:KEYSPACETIMECOERCIVITY} and \ref{L:COERCIVENESSOFCONTROLLING},
	inequalities
	\eqref{E:LESSSINGULARTERMSMPOINTNINEINTEGRALBOUND}
	and 
	\eqref{E:TANGENGITALEIKONALINTERMSOFCONTROLLING},
	simple estimates of the form
	$ab \lesssim a^2 + b^2$,
	and the fact that
	$\totmax{[1,N]}$
	is increasing in its arguments,
	we derive the desired estimate as follows: 
	\begin{align}
		& \int_{\mathcal{M}_{t,u}}
		 		\left|
					\Lunit \Fullset_{\ast}^{N;\leq 1} \Psi
				\right|
				\left|
					\Fullset_{\ast \ast}^{[1,N];\leq 1} \upmu
				\right|
		 \, d \vol
		 	\label{E:SECONDHARMLESSEXAMPLEINTEGRAL} \\
		 & \lesssim
		 	\int_{u'=0}^u
					\int_{\mathcal{P}_{u'}^t}
						\left|
							\Lunit \Fullset_{\ast}^{N;\leq 1} \Psi
						\right|^2
					\, d \conevol
			\, du'
			+
			\int_{t'=0}^t
		 			\int_{\Sigma_{t'}^u}
		 				\left|
							\Fullset_{\ast \ast}^{[1,N];\leq 1} \upmu
						\right|^2
					\, d \tvol
			\, dt'
				\notag
					\\ 
		& \lesssim
				\int_{u'=0}^u
					\totmax{[1,N]}(t,u')
				\, du'
				+
				\int_{t'=0}^t
					\left\lbrace
						\int_{s=0}^{t'}
							\frac{\sqrt{\totmax{[1,N]}}(s,u)}{\upmu_{\star}^{1/2}(s,u)}
						\, ds
					\right\rbrace^2
					+ 
					\mathring{\upepsilon}^2
				\, dt'
				\notag
					\\
		& \lesssim
				\int_{u'=0}^u
					\totmax{[1,N]}(t,u')
				\, du'
				+
				\int_{t'=0}^t
					\totmax{[1,N]}(t',u)
				\, dt'
				+ \mathring{\upepsilon}^2.
				\notag
	\end{align}

	As our third example, 
	we bound the spacetime integral of 
	$\left|
			\Lunit \Fullset_{\ast}^{N;\leq 1} \Psi
	 \right|
	 \left|
	 		\Fullset_{\ast}^{N+1;2} \Psi
	 \right|
	$.
	From Lemma~\ref{L:PUTALLRADFACTORSUPFRONT} with $M=2$
	and $N+1$ in the role of $N$, we obtain
	$|\Fullset_{\ast}^{N+1;2} \Psi|
	\lesssim
	|\Rad \Rad \Tanset^{[1,N-1]} \Psi|
	+
	|\Rad \Tanset^{[1,N]} \Psi|
	+ 
	|\Tanset^{[1,N+1]} \Psi|
	+
	\varepsilon |\Fullset_{\ast}^{[1,N];\leq 2} \GdVar|
	+ \varepsilon |\Fullset_{\ast \ast}^{[1,N];\leq 1} \BadVar|
	$. 
	Thus, we must bound the integral of the five corresponding products
	from the RHS of the previous inequality.
	To bound the integral of the first product,
	we argue as in the proof of 
	\eqref{E:FIRSTHARMLESSEXAMPLEINTEGRAL}
	to deduce that
	\begin{align}
		& \int_{\mathcal{M}_{t,u}}
		 		\left|
					\Lunit \Fullset_{\ast}^{N;\leq 1} \Psi
				\right|
				\left|
				 \Rad \Rad \Tanset^{[1,N-1]} \Psi
				\right|
		 \, d \vol
		 	\label{E:FINALHARMLESSEXAMPLEINTEGRAL} \\
		& \lesssim
		 	\int_{u'=0}^u
					\int_{\mathcal{P}_{u'}^t}
						\left|
							\Lunit \Fullset_{\ast}^{N;\leq 1} \Psi
						\right|^2
					\, d \conevol
			\, du'
			+
			\int_{t'=0}^t
		 			\int_{\Sigma_{t'}^u}
		 				\left|
							\Rad \Rad \Tanset^{[1,N-1]} \Psi
						\right|^2
					\, d \tvol
			\, dt'
				\notag
					\\ 
		& \lesssim
				\int_{u'=0}^u
					\totmax{[1,N]}(t,u')
				\, du'
				+
				\int_{t'=0}^t
					\totmax{[1,N]}(t',u)
				\, dt',
				\notag
	\end{align}
	which is $\lesssim$ RHS~\eqref{E:WAVESHARMLESSSPACETIMEINTEGRALS} as desired.
	The second product
	$\left|
			\Lunit \Fullset_{\ast}^{N;\leq 1} \Psi
	 \right|
	 \left|
			\Rad \Tanset^{[1,N]} \Psi
	 \right|
	$
	can be bounded in the same way.
	Similar reasoning yields that the integral
	of the third product
	$|\Lunit \Fullset_{\ast}^{N;\leq 1} \Psi|
	 |\Tanset^{[1,N+1]} \Psi|
	$
	is $\lesssim$ RHS~\eqref{E:FIRSTHARMLESSEXAMPLEINTEGRAL} 
	plus RHS~\eqref{E:SECONDHARMLESSEXAMPLEINTEGRAL} as desired.
	We clarify that the factor $\mathring{\upepsilon}^2$ is generated by the square of 
	RHS~\eqref{E:L2ESTIMATELOSSOFONEDERIVATIVE},
	which is needed to bound $\Singletan \Psi$.
	Similar reasoning, together with inequality
	\eqref{E:TANGENGITALEIKONALINTERMSOFCONTROLLING},
	yields that
	the integral
	of the third product
	$	\varepsilon
		|\Lunit \Fullset_{\ast}^{N;\leq 1} \Psi|
		|\Fullset_{\ast}^{[1,N];\leq 2} \GdVar|
	$
	and
	the integral
	of the fourth product
	$	\varepsilon
		|\Lunit \Fullset_{\ast}^{N;\leq 1} \Psi|
	 	|\Fullset_{\ast \ast}^{[1,N];\leq 1} \BadVar|
	$
	are $\lesssim$ RHS~\eqref{E:FIRSTHARMLESSEXAMPLEINTEGRAL}
	plus RHS~\eqref{E:SECONDHARMLESSEXAMPLEINTEGRAL} as desired.
	We clarify that we have used
	the fact that $\totmax{[1,N]}$ is increasing in its arguments
	and the estimate \eqref{E:LESSSINGULARTERMSMPOINTNINEINTEGRALBOUND}
	to bound the time integrals on RHSs
	\eqref{E:TANGENGITALEIKONALINTERMSOFCONTROLLING} 
	by $\lesssim \totmax{[1,N]}(t,u)$,
	as we did in passing to the last line of \eqref{E:SECONDHARMLESSEXAMPLEINTEGRAL}.

	As our fourth example, we bound the integral of
	$
	\left|
		\Lunit \Fullset_{\ast}^{N;\leq 1} \Psi
	 \right|
	 \left|
	 		\Fullset^{N;1} \Vortrenormalized
	 \right|
	$.
	We argue as in the proof of 
	\eqref{E:FIRSTHARMLESSEXAMPLEINTEGRAL}
	and use Lemmas~\ref{L:COERCIVENESSOFCONTROLLING}
	and \ref{L:L2NORMOFVORTICITYTRANSVERSALDERIVATIVESINTERMSOFTANGENTIALDERIVATIVES}
	to deduce that
	\begin{align}
		& \int_{\mathcal{M}_{t,u}}
		 		\left|
					\Lunit \Fullset_{\ast}^{N;\leq 1} \Psi
				\right|
				\left|
				 \Fullset^{N;1} \Vortrenormalized
				\right|
		 \, d \vol
		 	\label{E:FOURTHHARMLESSEXAMPLEINTEGRAL} \\
		& \lesssim
		 	\int_{u'=0}^u
					\int_{\mathcal{P}_{u'}^t}
						\left|
							\Lunit \Fullset_{\ast}^{N;\leq 1} \Psi
						\right|^2
					\, d \conevol
			\, du'
			+
			\int_{t'=0}^t
		 			\int_{\Sigma_{t'}^u}
		 				\left|
							\Fullset^{N;1} \Vortrenormalized
						\right|^2
					\, d \tvol
			\, dt'
				\notag
					\\ 
		& \lesssim
				\int_{u'=0}^u
					\totmax{[1,N]}(t,u')
				\, du'
				+
				\int_{t'=0}^t
					\left\lbrace
					\int_{s=0}^{t'}
						\frac{\sqrt{\totmax{[1,N-1]}}(t,u)}{\upmu_{\star}^{1/2}(s,u)}
					\, ds
					\right\rbrace^2
				\, dt'
					\notag \\
	& \ \
				+
				\int_{t'=0}^t
					\Vorttotmax{\leq N}(t',u)
				\, dt'
				+
				\mathring{\upepsilon}^2.
				\notag
	\end{align}
	Using
	inequality \eqref{E:LESSSINGULARTERMSMPOINTNINEINTEGRALBOUND}
	and the fact that the $\totmax{[1,N-1]}$ is increasing in its arguments,
	we bound the double time integral on RHS~\eqref{E:FOURTHHARMLESSEXAMPLEINTEGRAL}
	by 
	$
	\displaystyle
	\lesssim 
	\int_{t'=0}^t
		\totmax{[1,N-1]}(t',u)
	\, dt'
	\leq
	\int_{t'=0}^t
		\totmax{[1,N]}(t',u)
	\, dt'
	$.
	We conclude that
	RHS~\eqref{E:FOURTHHARMLESSEXAMPLEINTEGRAL}
	is
	$\lesssim$ RHS~\eqref{E:WAVESHARMLESSSPACETIMEINTEGRALS} as desired.
	This completes our proof of the representative estimates
	from \eqref{E:WAVESHARMLESSSPACETIMEINTEGRALS}.

	We now prove one representative estimate
	from \eqref{E:VORTICITYHARMLESSSPACETIMEINTEGRALS}.
	Specifically, we bound the integral of the product
	$
	\varepsilon
	|\Tanset^{\leq N} \Vortrenormalized|
	|\Fullset_{\ast}^{N;\leq 2} \threePsi|
	$
	in the cases $1 \leq N \leq 21$.
	Using Lemma~\ref{L:PUTALLRADFACTORSUPFRONT},
	we bound the last factor as follows:
	\begin{align} \label{E:SIMPLECOMMEST}
	|\Fullset_{\ast}^{N;\leq 2} \threePsi|
	& \lesssim
	|\Rad \Fullset_{\ast}^{[1,N-1];\leq 1} \threePsi|
	+ 
	|\Lunit \Tanset^{\leq N-1} \threePsi|
	+
	|\angdiff \Tanset^{\leq N-1} \threePsi|
		\\
	& \ \
		+
		\sum_{i=1}^2|\Fullset_{\ast}^{[1,N-1];\leq 2} \Lunit_{(Small)}^i|
		+ 
		|\Fullset_{\ast \ast}^{[1,N-1];\leq 1} \upmu|,
		\notag
	\end{align}
	where the first, fourth, and fifth terms on RHS~\eqref{E:SIMPLECOMMEST}
	are absent when $N=1$.
	Thus, we must bound the integral of 
	$\varepsilon$ times the five corresponding products
	from the RHS of the previous inequality.
	To this end, we first use Young's inequality to obtain
	\begin{align}
		& 
			\varepsilon
			\int_{\mathcal{M}_{t,u}}
		 		\left|
					\Tanset^{\leq N} \Vortrenormalized
				\right|
				\left|
					\Fullset_{\ast}^{N;\leq 2} \threePsi 
				\right|
		 \, d \vol
		 	\label{E:HARMLESSVORTICITYEXAMPLEINTEGRAL} \\
		& \lesssim
		 	\int_{u'=0}^u
				\left\|
					\Tanset^{\leq N} \Vortrenormalized
				\right\|_{L^2(\mathcal{P}_{u'}^t)}^2
			\, du'
				+
			\varepsilon^2
			\int_{t'=0}^t
				\left\|
					\Rad \Fullset_{\ast}^{[1,N-1];\leq 1} \threePsi
				\right\|_{L^2(\Sigma_{t'}^u)}^2
			\, dt'
				\notag \\
		& \ \
			+
			\varepsilon^2
			\int_{u'=0}^u
				\left\|
					\Lunit \Tanset^{\leq N-1} \threePsi
				\right\|_{L^2(\mathcal{P}_{u'}^t)}^2
			\, du'
			+
			\varepsilon^2
			\int_{u'=0}^u
				\left\|
					\angdiff \Tanset^{\leq N-1} \threePsi
				\right\|_{L^2(\mathcal{P}_{u'}^t)}^2
			\, du'
			\notag \\
		& \ \
			+
			\varepsilon^2
			\int_{t'=0}^t
				\sum_{i=1,2}
				\left\|
					\Fullset_{\ast}^{[1,N-1];\leq 2} \Lunit_{(Small)}^i
				\right\|_{L^2(\Sigma_{t'}^u)}^2
			\, dt'
			+
			\varepsilon^2
			\int_{t'=0}^t
				\left\|
					\Fullset_{\ast \ast}^{[1,N-1];\leq 1} \upmu
				\right\|_{L^2(\Sigma_{t'}^u)}^2
			\, dt'.
				\notag 
	\end{align}
	By Lemma~\ref{L:COERCIVENESSOFCONTROLLING},
	the integrals on the first line of RHS~\eqref{E:HARMLESSVORTICITYEXAMPLEINTEGRAL}
	and the $\Lunit \Tanset^{\leq N-1} \threePsi$ integral 
	on the second line
	are $\lesssim$ RHS~\eqref{E:VORTICITYHARMLESSSPACETIMEINTEGRALS}
	as desired. 
	Moreover, using Lemma~\ref{L:KEYSPACETIMECOERCIVITY},
	we find that the $\angdiff \Tanset^{\leq N-1} \threePsi$ integral on
	RHS~\eqref{E:HARMLESSVORTICITYEXAMPLEINTEGRAL}
	is, for $N \geq 2$, $\lesssim$ 
	the $\varepsilon^2 \coercivespacetimemax{[1,N-1]}(t,u)$ term
	on RHS~\eqref{E:VORTICITYHARMLESSSPACETIMEINTEGRALS} as desired.
	In the case $N=1$, to bound the 
	$\angdiff \threePsi$ integral
	on RHS~\eqref{E:HARMLESSVORTICITYEXAMPLEINTEGRAL} by
	$\lesssim$ RHS~\eqref{E:VORTICITYHARMLESSSPACETIMEINTEGRALS}, 
	we again use Lemma~\ref{L:COERCIVENESSOFCONTROLLING}.
	Finally, the argument given in our second example above
	implies that the terms on the last line 
	of RHS~\eqref{E:HARMLESSVORTICITYEXAMPLEINTEGRAL}
	are $\lesssim$ RHS~\eqref{E:VORTICITYHARMLESSSPACETIMEINTEGRALS} 
	as desired. This completes our proof of the representative estimate
	from \eqref{E:VORTICITYHARMLESSSPACETIMEINTEGRALS}.

\end{proof}

\subsection{Estimates for wave equation error integrals involving top-order vorticity terms}
\label{SS:WAVEENERGYESTIMTATESINVOLVINGTOPORDERVORTICITY}
Recall that the geometric wave equation
\eqref{E:VELOCITYWAVEEQUATION}
has inhomogeneous terms depending on $\Vortrenormalized$.
In the next lemma, we derive 
simple estimates for the corresponding error integrals that
depend on the top-order derivatives of $\Vortrenormalized$.
The precise form of the terms that we bound corresponds to the explicit
vorticity-involving terms on  
RHS~\eqref{E:WAVELISTHEFIRSTCOMMUTATORIMPORTANTTERMS}-\eqref{E:KEYDIFFICULTFORRENORMALZIEDDENSITY}.

\begin{lemma}[\textbf{Estimates for wave equation integrals involving top-order vorticity terms}]
\label{L:WAVEERRORINTEGRALSINVOLVINGTOPORDERVORTICITY}
	Assume that 
	$\Psi 
	\in 
	\lbrace 
		\Densrenormalized  - v^1,v^1,v^2 
	\rbrace$ 
	and that $1 \leq N \leq 20$.
	Under the data-size and bootstrap assumptions 
	of Sects.~\ref{SS:FLUIDVARIABLEDATAASSUMPTIONS}-\ref{SS:PSIBOOTSTRAP}
	and the smallness assumptions of Sect.~\ref{SS:SMALLNESSASSUMPTIONS}, 
	the following integral estimates hold for $(t,u) \in [0,\Tboot) \times [0,U_0]$
	(see Sect.~\ref{SS:STRINGSOFCOMMUTATIONVECTORFIELDS} regarding the vectorfield operator notation):
	\begin{align}
	 & 
	\int_{\mathcal{M}_{t,u}}
		\left|
			\myarray[\Rad \Fullset_{\ast}^{N;\leq 1} \Psi]
				{(1 + \upmu) \Lunit \Fullset_{\ast}^{N;\leq 1} \Psi}
		\right|
		\left|
			\myarray[{[ia]} \upmu (\exp \Densrenormalized) \Speed^2 (g_{ab} \Radunit^b) 
					\Tanset^{N+1} \Vortrenormalized]
				{{[ia]} \upmu (\exp \Densrenormalized) \Speed^2
			\left(
				\frac{g_{ab} \GeoAng^b}{g_{cd} \GeoAng^c \GeoAng^d}
			\right)
			\Tanset^{N+1} \Vortrenormalized}
		\right|
	\, d \vol
			\label{E:WAVEERRORINTEGRALSINVOLVINGTOPORDERVORTICITY} \\
		& \lesssim
		\notag
		\int_{t'=0}^t
			\totmax{N}(t',u)
		\, dt'
		+
		\int_{t'=0}^t
			\VorttotTanmax{N+1}(t',u)
		\, dt'.
	\end{align}
\end{lemma}

\begin{proof}
	Using 
	the schematic relations \eqref{E:SCALARSDEPENDINGONGOODVARIABLES},
	the $L^{\infty}$ estimates of Prop.~\ref{P:IMPROVEMENTOFAUX},
	Young's inequality,
	and Lemma~\ref{L:COERCIVENESSOFCONTROLLING},
	we bound LHS~\eqref{E:WAVEERRORINTEGRALSINVOLVINGTOPORDERVORTICITY}
	as follows:
	\begin{align}
	& 
	\int_{\mathcal{M}_{t,u}}
		\left|
			\myarray[\Rad \Fullset_{\ast}^{N;\leq 1} \Psi]
				{(1 + \upmu) \Lunit \Fullset_{\ast}^{N;\leq 1} \Psi}
		\right|
		\left|
			\myarray[{[ia]} \upmu (\exp \Densrenormalized) \Speed^2 (g_{ab} \Radunit^b) 
					\Tanset^{N+1} \Vortrenormalized]
				{{[ia]} \upmu (\exp \Densrenormalized) \Speed^2
			\left(
				\frac{g_{ab} \GeoAng^b}{g_{cd} \GeoAng^c \GeoAng^d}
			\right)
			\Tanset^{N+1} \Vortrenormalized}
		\right|
	\, d \vol
			\\
	& \lesssim
		\int_{\mathcal{M}_{t,u}}
			\left|
				\Rad \Fullset_{\ast}^{N;\leq 1} \Psi
			\right|^2
		\, d \vol
		+
		\int_{\mathcal{M}_{t,u}}
			\upmu
			\left|
				\Lunit \Fullset_{\ast}^{N;\leq 1} \Psi
			\right|^2
		\, d \vol
		+
		\int_{\mathcal{M}_{t,u}}
			\upmu
			\left|
				\Tanset^{N+1} \Vortrenormalized
			\right|^2
		\, d \vol
			\notag \\
	& \lesssim
		\int_{t'=0}^t
			\left\|
				\Rad \Fullset_{\ast}^{N;\leq 1} \Psi
			\right\|_{L^2(\Sigma_{t'}^u	)}^2
		\, dt'
		+
		\int_{t'=0}^t
			\left\|
				\sqrt{\upmu} \Lunit \Fullset_{\ast}^{N;\leq 1} \Psi
			\right\|_{L^2(\Sigma_{t'}^u	)}^2
		\, dt'
			\notag \\
	& \ \
		+
		\int_{t'=0}^t
			\left\|
				\sqrt{\upmu} \Tanset^{N+1} \Vortrenormalized
			\right\|_{L^2(\Sigma_{t'}^u	)}^2
		\, dt'
			\notag 
				\\
	& \lesssim
		\int_{t'=0}^t
			\totmax{N}(t',u)
		\, dt'
		+
		\int_{t'=0}^t
			\VorttotTanmax{N+1}(t',u)
		\, dt'.
		\notag
	\end{align}
	\end{proof}

\subsection{\texorpdfstring{$L^2$}{Square integral} bounds for the difficult top-order error integrals in terms of 
\texorpdfstring{$\totmax{[1,N]}$}{the fundamental controlling quantities}}
\label{SS:MOSTDIFFICULTENERGYESTIMATEINTEGRALS}
In this section, 
we derive estimates for the difficult error integrals
that we encounter in our energy estimates for
$\threePsi$. 
These error integrals 
would cause derivative loss 
if they were not treated carefully and moreover, they
make a substantial contribution to the blowup-rates featured 
in our high-order energy estimates.
Our arguments here rely on the fully modified
quantities defined in Sect.~\ref{S:MODIFIED}.

The main result is Lemma~\ref{L:MOSTDIFFICULTERRORINTEGRALS}.
We start with a preliminary lemma in which we estimate
the most difficult product that appears in our wave equation energy estimates.

\begin{lemma}[$L^2$ \textbf{bound for the most difficult product}]
	\label{L:DIFFICULTTERML2BOUND}
		Assume that 
		$N=20$ and that
		$\Fullset_{\ast}^{N;\leq 1} \in \lbrace \GeoAng^N, \GeoAng^{N-1} \Rad \rbrace$.
		There exist constants $C > 0$ and $C_* > 0$ such that
		under the data-size and bootstrap assumptions 
		of Sects.~\ref{SS:FLUIDVARIABLEDATAASSUMPTIONS}-\ref{SS:PSIBOOTSTRAP}
		and the smallness assumptions of Sect.~\ref{SS:SMALLNESSASSUMPTIONS},
		the following $L^2$ estimate holds for the difficult product 
		$(\Rad v^1) \Fullset_{\ast}^{N;\leq 1} \mytr \upchi$ 
		from Prop.~\ref{P:KEYPOINTWISEESTIMATE}
		whenever $(t,u) \in [0,\Tboot) \times [0,U_0]$:
	\begingroup
\allowdisplaybreaks
	\begin{align} \label{E:DIFFICULTTERML2BOUND}
		\left\|
			(\Rad v^1) \Fullset_{\ast}^{N;\leq 1} \mytr \upchi
		\right\|_{L^2(\Sigma_t^u)}
		& \leq
			\boxed{2} 
			\frac{\left\| [\Lunit \upmu]_- \right\|_{L^{\infty}(\Sigma_t^u)}} 
						{\upmu_{\star}(t,u)} 
				\sqrt{\totmax{N}}(t,u)
			\\
	& \ \ +
			\boxed{4.05}
			\frac{\left\| [\Lunit \upmu]_- \right\|_{L^{\infty}(\Sigma_t^u)}} 
						{\upmu_{\star}(t,u)} 
			\int_{s=0}^t
				\frac{\left\| [\Lunit \upmu]_- \right\|_{L^{\infty}(\Sigma_s^u)}} 
				{\upmu_{\star}(s,u)} 
				\sqrt{\totmax{N}}(s,u) 
			\, ds
			\notag \\
	& \ \
		+
		C_* \frac{1}{\upmu_{\star}(t,u)} \sqrt{\easytotmax{N}}(t,u)
			\notag \\
	& \ \ +
			C_* \frac{1}{\upmu_{\star}(t,u)} 
			\int_{s=0}^t
				\frac{1} 
				{\upmu_{\star}(s,u)} 
				\sqrt{\easytotmax{N}}(s,u) 
			\, ds
			\notag \\
	& \ \ +
			C \varepsilon
			\frac{1} {\upmu_{\star}(t,u)} 
			\int_{s=0}^t
				\frac{1}{\upmu_{\star}(s,u)} 
				\sqrt{\totmax{N}}(s,u) 
			\, ds
			\notag \\
	& \ \
			+ C
				\frac{1}{\upmu_{\star}(t,u)}
				\int_{s'=0}^t
					\frac{1}{\upmu_{\star}(s',u)}
					\int_{s=0}^{s'}
						\frac{1}{\upmu_{\star}^{1/2}(s,u)}
						\sqrt{\totmax{N}}(s,u)
					\ ds
				\, ds'
			\notag  \\
	& \ \
		+ C
			\frac{1}{\upmu_{\star}(t,u)}
			\int_{s=0}^t
					\frac{1}{\upmu_{\star}^{1/2}(s,u)}
					\sqrt{\totmax{N}}(s,u)
			\, ds
				\notag \\
		& \ \
			+ C \frac{1}{\upmu_{\star}^{1/2}(t,u)} \sqrt{\totmax{N}}(t,u)
			+ C \frac{1}{\upmu_{\star}^{3/2}(t,u)} \sqrt{\totmax{[1,N-1]}}(t,u)
				\notag  \\
		& \ \ 
			+
			C \frac{1}{\upmu_{\star}(t,u)}
			\int_{s=0}^t
				\sqrt{\VorttotTanmax{N+1}}(s,u)
			\, ds
				\notag \\
		& \ \
			+
			C \frac{1}{\upmu_{\star}(t,u)}
			\int_{s=0}^t
				\frac{1} 
				{\upmu_{\star}^{1/2}(s,u)}
				\sqrt{\VorttotTanmax{\leq N}}(s,u)
			\, ds
			\notag \\
		&  \ \
			+ C \frac{1}{\upmu_{\star}^{3/2}(t,u)} \mathring{\upepsilon}.
			\notag
	\end{align}
	\endgroup

	Moreover, we have the following less degenerate estimates:
	\begingroup
\allowdisplaybreaks
	\begin{align} \label{E:LESSDEGENERATEDIFFICULTTERML2BOUND}
		\myarray[
		\left\|
			(\Rad v^2) \Fullset_{\ast}^{N;\leq 1} \mytr \upchi
		\right\|_{L^2(\Sigma_t^u)}
		]
		{\left\|
			\left\lbrace
				\Rad(\Densrenormalized  - v^1) 
			\right\rbrace
			\Fullset_{\ast}^{N;\leq 1} \mytr \upchi
		\right\|_{L^2(\Sigma_t^u)}
		}
	& \lesssim
			\varepsilon
			\frac{1}{\upmu_{\star}(t,u)}
			\sqrt{\totmax{N}}(t,u)
				\\
	& \ \ 
		+
			\varepsilon
			\frac{1}{\upmu_{\star}(t,u)}
			\int_{s=0}^t
				\frac{1} 
				{\upmu_{\star}(s,u)} 
					\sqrt{\totmax{N}}(s,u)
			\, ds
				\notag
				\\
	& \ \ 
			+
			\varepsilon
			\frac{1}{\upmu_{\star}(t,u)}
			\int_{s=0}^t
				\sqrt{\VorttotTanmax{N+1}}(s,u)
			\, ds
			\notag \\
	& \ \
			+
			\varepsilon
			\frac{1}{\upmu_{\star}(t,u)}
			\int_{s=0}^t
				\frac{1} 
				{\upmu_{\star}^{1/2}(s,u)}
				\sqrt{\VorttotTanmax{\leq N}}(s,u)
			\, ds
				\notag \\
		& \ \
			+ 
			\varepsilon
			\frac{1}{\upmu_{\star}^{3/2}(t,u)} \sqrt{\totmax{[1,N-1]}}(t,u)
			+ \mathring{\upepsilon}
				\frac{1}{\upmu_{\star}^{3/2}(t,u)}.
				\notag
	\end{align}
	\endgroup

	Furthermore, we have the following less precise estimate:
	\begin{align} \label{E:LESSPRECISEDIFFICULTTERML2BOUND}
		\left\|
			\upmu \Fullset_{\ast}^{N;\leq 1} \mytr \upchi
		\right\|_{L^2(\Sigma_t^u)}
		& \lesssim
			\sqrt{\totmax{[1,N]}}(t,u)
			+
			\int_{s=0}^t
				\frac{1} 
				{\upmu_{\star}(s,u)} 
					\sqrt{\totmax{[1,N]}}(s,u)
			\, ds
					\\
	& \ \ 
			+
			\int_{s=0}^t
				\sqrt{\VorttotTanmax{N+1}}(s,u)
			\, ds
			+
			\int_{s=0}^t
				\frac{1} 
				{\upmu_{\star}^{1/2}(s,u)}
				\sqrt{\VorttotTanmax{\leq N}}(s,u)
			\, ds
				\notag
				\\
	& \ \
		+ \mathring{\upepsilon} 
				\left\lbrace 
					\ln \upmu_{\star}^{-1}(t,u) 
					+ 
					1 
				\right\rbrace.
				\notag
	\end{align}
\end{lemma}
\begin{proof}
	We first prove \eqref{E:DIFFICULTTERML2BOUND}. We take the norm
	$\| \cdot \|_{L^2(\Sigma_t^u)}$
	of both sides of inequality \eqref{E:KEYPOINTWISEESTIMATE}.
	Using Lemma~\ref{L:COERCIVENESSOFCONTROLLING}, we see
	that the norm $\| \cdot \|_{L^2(\Sigma_t^u)}$ of the first term 
	on RHS~\eqref{E:KEYPOINTWISEESTIMATE} is bounded by the first term 
	on RHS~\eqref{E:DIFFICULTTERML2BOUND}.
	Similarly, using Lemma~\ref{L:COERCIVENESSOFCONTROLLING},
	we see that
	the norm $\| \cdot \|_{L^2(\Sigma_t^u)}$
	of the second term on RHS~\eqref{E:KEYPOINTWISEESTIMATE} is bounded by
	the term
	$
	\displaystyle
	C_* \frac{1}{\upmu_{\star}(t,u)} \sqrt{\easytotmax{N}}(t,u)
	$
	on RHS~\eqref{E:DIFFICULTTERML2BOUND}.
	Next we use Lemmas~\ref{L:L2NORMSOFTIMEINTEGRATEDFUNCTIONS} and \ref{L:COERCIVENESSOFCONTROLLING} 
	to bound
	the norm $\| \cdot \|_{L^2(\Sigma_t^u)}$
	of the third term on RHS~\eqref{E:KEYPOINTWISEESTIMATE} 
	by the term $\boxed{4.05} \cdots$
	on RHS~\eqref{E:DIFFICULTTERML2BOUND}.
	Similarly, using Lemmas~\ref{L:L2NORMSOFTIMEINTEGRATEDFUNCTIONS} and \ref{L:COERCIVENESSOFCONTROLLING} 
	we see that
	the norm $\| \cdot \|_{L^2(\Sigma_t^u)}$
	of the fourth term on RHS~\eqref{E:KEYPOINTWISEESTIMATE} is bounded by
	the 
	$\sqrt{\easytotmax{N}}$-involving time integral term
	on RHS~\eqref{E:DIFFICULTTERML2BOUND}
	(which is multiplied by $C_*$).

	It remains for us to explain why the norm $\| \cdot \|_{L^2(\Sigma_t^u)}$
	of the terms $\mbox{\upshape Error}$
	on RHS~\eqref{E:KEYPOINTWISEESTIMATE}
	are $\leq$ the sum of the 
	terms on lines five to eleven of RHS~\eqref{E:DIFFICULTTERML2BOUND}.
	With the exception of the bound for the terms
	on the first line RHS~\eqref{E:ERRORTERMKEYPOINTWISEESTIMATE},
	the desired bounds follow from
	the same estimates used above
	together with those of Lemmas~\ref{L:EASYL2BOUNDSFOREIKONALFUNCTIONQUANTITIES}
	and \ref{L:L2NORMOFVORTICITYTRANSVERSALDERIVATIVESINTERMSOFTANGENTIALDERIVATIVES},
	inequalities
	\eqref{E:LOSSKEYMUINTEGRALBOUND},
	\eqref{E:LOGLOSSMUINVERSEINTEGRALBOUND},
	and
	\eqref{E:LESSSINGULARTERMSMPOINTNINEINTEGRALBOUND},
	the fact that the $\totmax{M}$ are increasing in their arguments, 
	and simple inequalities of the form
	$ab \lesssim a^2 + b^2$.
	Finally, we must bound 
the norm $\| \cdot \|_{L^2(\Sigma_t^u)}$
	of the terms 
	on the first line of RHS~\eqref{E:ERRORTERMKEYPOINTWISEESTIMATE}.
	To this end, we first use 
	\eqref{E:SIGMATVOLUMEFORMCOMPARISON} with $s=0$
	to deduce
	$
	\displaystyle
	\left\|
		\left| 
			\upchifullmodarg{\GeoAng^N} 
		\right|
		(1,\cdot)
		+
		\left| 
			\upchifullmodarg{\GeoAng^{N-1} \Rad} 
		\right|
		(1,\cdot)
\right\|_{L^2(\Sigma_t^u)}
\lesssim
\left\|
	\left| 
		\upchifullmodarg{\GeoAng^N} 
	\right|
	+
	\left| 
		\upchifullmodarg{\GeoAng^{N-1} \Rad} 
	\right|
\right\|_{L^2(\Sigma_0^u)}
$.
Next, from definition \eqref{E:TRANSPORTPARTIALRENORMALIZEDTRCHIJUNK},
the simple inequality
$|\vec{G}_{(Frame)}| = |\smoothfunction(\GdVar,\angdiff \vec{x})| \lesssim 1$
(which follows from Lemmas~\ref{L:SCHEMATICDEPENDENCEOFMANYTENSORFIELDS}
and \ref{L:POINTWISEFORRECTANGULARCOMPONENTSOFVECTORFIELDS}
and the $L^{\infty}$ estimates of Prop.~\ref{P:IMPROVEMENTOFAUX}),
the estimate \eqref{E:POINTWISEESTIMATESFORCHIANDITSDERIVATIVES},
and our assumptions on the data,
we find that 
$
\displaystyle
\left\|
	\left| 
		\upchifullmodarg{\GeoAng^N} 
	\right|
	+
	\left| 
		\upchifullmodarg{\GeoAng^{N-1} \Rad} 
	\right|
\right\|_{L^2(\Sigma_0^u)}
\lesssim \mathring{\upepsilon}
$.
It follows that
the norm $\| \cdot \|_{L^2(\Sigma_t^u)}$
of the terms on the first line of RHS~\eqref{E:ERRORTERMKEYPOINTWISEESTIMATE}
is 
$
\displaystyle
\lesssim 
\mathring{\upepsilon}
\frac{1} 
				{\upmu_{\star}(t,u)}
$
as desired. This completes the proof of \eqref{E:DIFFICULTTERML2BOUND}.

The proof of \eqref{E:LESSPRECISEDIFFICULTTERML2BOUND} 
is based on inequality \eqref{E:LESSPRECISEKEYPOINTWISEESTIMATE}
and is similar but much simpler; we omit the details, noting
only that inequality \eqref{E:LOGLOSSMUINVERSEINTEGRALBOUND}
leads to the presence of the factor 
$
	\ln \upmu_{\star}^{-1}(t,u) + 1 
$.

The estimate \eqref{E:LESSDEGENERATEDIFFICULTTERML2BOUND}
then follows from \eqref{E:LESSPRECISEDIFFICULTTERML2BOUND} 
and the estimates
$
\| \Rad v^2 \|_{L^{\infty}(\Sigma_t^u)}
\lesssim 
\varepsilon
$
and
$
	\| \Rad (\Densrenormalized  - v^1) \|_{L^{\infty}(\Sigma_t^u)} 
	\lesssim
	\varepsilon
$
(that is, 
\eqref{E:PSITRANSVERSALLINFINITYBOUNDBOOTSTRAPIMPROVEDSMALL} and \eqref{E:CRUCIALPSITRANSVERSALLINFINITYBOUNDBOOTSTRAPIMPROVEDSMALL}).

\end{proof}

Armed with Lemma~\ref{L:DIFFICULTTERML2BOUND},
we now derive the main result of this section.

\begin{lemma}[\textbf{Bound for the most difficult error integrals}]
	\label{L:MOSTDIFFICULTERRORINTEGRALS}
		Assume that $N=20$.
		There exist constants $C > 0$ and $C_* > 0$ such that
		under the data-size and bootstrap assumptions 
		of Sects.~\ref{SS:FLUIDVARIABLEDATAASSUMPTIONS}-\ref{SS:PSIBOOTSTRAP}
		and the smallness assumptions of Sect.~\ref{SS:SMALLNESSASSUMPTIONS},
		the following integral inequalities hold
		for $(t,u) \in [0,\Tboot) \times [0,U_0]$:
	\begin{align} \label{E:MOSTDIFFICULTERRORINTEGRALBOUND}
		&
		2
		\left|
			\int_{\mathcal{M}_{t,u}}
				(\Rad \Fullset_{\ast}^{N;\leq 1} v^1)	 
				(\Rad v^1) \GeoAng^N \mytr \upchi
			\, d \vol
		\right|,
			\,
		2
		\left|
			\int_{\mathcal{M}_{t,u}}
				(\Rad \Fullset_{\ast}^{N;\leq 1} v^1)	 
				(\Rad v^1) \GeoAng^{N-1} \Rad \mytr \upchi
			\, d \vol
		\right|
			\\
		& \leq
			\boxed{4} 
				\int_{t'=0}^t
					\frac{\left\| [\Lunit \upmu]_- \right\|_{L^{\infty}(\Sigma_{t'}^u)}} 
							{\upmu_{\star}(t',u)} 
							\totmax{N}(t',u)
							\, dt'
			\notag \\
	& \ \
			+ 
			\boxed{8.1}
			\int_{t'=0}^t
				\frac{\left\| [\Lunit \upmu]_- \right\|_{L^{\infty}(\Sigma_{t'}^u)}} 
								 {\upmu_{\star}(t',u)} 
						\sqrt{\totmax{N}}(t',u) 
						\int_{s=0}^{t'}
							\frac{\left\| [\Lunit \upmu]_- \right\|_{L^{\infty}(\Sigma_s^u)}} 
									{\upmu_{\star}(s,u)} 
							\sqrt{\totmax{N}}(s,u) 
						\, ds
				\, dt'
			\notag	\\
		& \ \ 
			+ 
			C_{\ast}
			\int_{t'=0}^t
					\frac{1} 
							 {\upmu_{\star}(t',u)} 
				  \sqrt{\totmax{N}}(t',u) 
					\sqrt{\easytotmax{N}}(t',u) 
				\, dt'
				\notag \\
	& \ \
			+ 
			C_{\ast}
			\int_{t'=0}^t
				\frac{1}
					{\upmu_{\star}(t',u)} 
						\sqrt{\totmax{N}}(t',u) 
						\int_{s=0}^{t'}
							\frac{1} 
									{\upmu_{\star}(s,u)} 
							\sqrt{\easytotmax{N}}(s,u) 
						\, ds
				\, dt'
			\notag	\\
	& \ \
			+ 
			C \varepsilon
			\int_{t'=0}^t
				\frac{1} 
						{\upmu_{\star}(t',u)} 
						\sqrt{\totmax{N}}(t',u) 
						\int_{s=0}^{t'}
							\frac{1} 
									{\upmu_{\star}(s,u)} 
							\sqrt{\totmax{N}}(s,u) 
						\, ds
				\, dt'
			\notag	\\
	& \ \
			+ C
			\int_{t'=0}^t
				\frac{1}{\upmu_{\star}(t',u)}
				\sqrt{\totmax{N}}(t',u)
				\int_{s'=0}^{t'}
					\frac{1}{\upmu_{\star}(s',u)}
					\int_{s=0}^{s'}
						\frac{1}{\upmu_{\star}^{1/2}(s,u)}
						\sqrt{\totmax{N}}(s,u)
					\ ds
				\, ds'
			\, dt'
			\notag  \\
	& \ \
		+ C
		\int_{t'=0}^t
			\frac{1}{\upmu_{\star}(t',u)}
			\sqrt{\totmax{N}}(t',u)
			\int_{s=0}^{t'}
					\frac{1}{\upmu_{\star}^{1/2}(s,u)}
					\sqrt{\totmax{N}}(s,u)
			\, ds
		\, dt'
				\notag \\
		& \ \
			+ C 
				\int_{t'=0}^t
					\frac{1}{\upmu_{\star}^{1/2}(t,u)} \totmax{N}(t',u)
				\, dt'
			+ C 
				\int_{t'=0}^t
					\frac{1}{\upmu_{\star}^{5/2}(t',u)} \totmax{[1,N-1]}(t',u)
				\, dt'
				\notag  \\
		& \ \ 
			+
		C 
		\int_{t'=0}^t
			\frac{1}{\upmu_{\star}^{3/2}(t',u)}
			\left\lbrace
				\int_{s=0}^{t'}
					\sqrt{\VorttotTanmax{N+1}}(s,u)
				\, ds
			\right\rbrace^2
		\, dt'
				\notag \\
		& \ \
			+
			C 
			\int_{t'=0}^t
				\frac{1}{\upmu_{\star}^{3/2}(t',u)}
				\left\lbrace
					\int_{s=0}^{t'}
						\frac{1} 
						{\upmu_{\star}^{1/2}(s,u)}
						\sqrt{\VorttotTanmax{\leq N}}(s,u)
					\, ds
					\right\rbrace^2
			\, dt'
			\notag \\
		&  \ \
			+ 
			C 
				\mathring{\upepsilon}
				\frac{1}{\upmu_{\star}^{3/2}(t,u)}.
			\notag
	\end{align}

	Moreover, we have the following less degenerate estimates:
	\begin{align} \label{E:LESSDEGENERATEDIFFICULTERRORINTEGAL}
		&
		2
		\left|
			\int_{\mathcal{M}_{t,u}}
				(\Rad \Fullset_{\ast}^{N;\leq 1} v^2)	 
				(\Rad v^2) \GeoAng^N \mytr \upchi
			\, d \vol
		\right|,
			\,
		2
		\left|
			\int_{\mathcal{M}_{t,u}}
				(\Rad \Fullset_{\ast}^{N;\leq 1} v^2)	 
				(\Rad v^2) \GeoAng^{N-1} \Rad \mytr \upchi
			\, d \vol
		\right|,
			\\
		&
		2
		\left|
			\int_{\mathcal{M}_{t,u}}
				\left\lbrace
					\Rad \Fullset_{\ast}^{N;\leq 1} (\Densrenormalized  - v^1)	 
				\right\rbrace
				\left\lbrace
					\Rad (\Densrenormalized  - v^1)
				\right\rbrace
				\GeoAng^N \mytr \upchi
			\, d \vol
		\right|,
			\,
				\notag \\
		& 2
		\left|
			\int_{\mathcal{M}_{t,u}}
				\left\lbrace
					\Rad \Fullset_{\ast}^{N;\leq 1} (\Densrenormalized  - v^1)	 
				\right\rbrace
				\left\lbrace
					\Rad (\Densrenormalized  - v^1)
				\right\rbrace 
				\GeoAng^{N-1} \Rad \mytr \upchi
			\, d \vol
		\right|
		\notag	\\
		& \lesssim
			\varepsilon
			\int_{t'=0}^t
				\frac{1}{\upmu_{\star}(t',u)}
				\totmax{N}(t',u)
			\, dt'
			\notag	\\
	& \ \ 
		+
			\varepsilon
			\int_{t'=0}^t
				\frac{1}{\upmu_{\star}(t',u)}
				\sqrt{\totmax{N}}(t',u)
				\int_{s=0}^t
					\frac{1} 
					{\upmu_{\star}(s,u)} 
					\sqrt{\totmax{N}}(s,u)
				\, ds
			\, ds
					\notag
				\\
	& \ \ 
			+
			\varepsilon
			\int_{t'=0}^t
				\frac{1}{\upmu_{\star}(t',u)}
				\sqrt{\totmax{N}}(t',u)
				\int_{s=0}^{t'}
					\sqrt{\VorttotTanmax{N+1}}(s,u)
				\, ds
			\, dt'
			\notag \\
	& \ \
			+
			\varepsilon
			\int_{t'=0}^t
				\frac{1}{\upmu_{\star}(t',u)}
				\sqrt{\totmax{N}}(t',u)
				\int_{s=0}^{t'}
					\frac{1} 
					{\upmu_{\star}^{1/2}(s,u)}
					\sqrt{\VorttotTanmax{\leq N}}(s,u)
				\, ds
			\, dt'
				\notag \\
		& \ \
			+ 
			\varepsilon
			\int_{t'=0}^t
				\frac{1}{\upmu_{\star}^{5/2}(t',u)} \totmax{[1,N-1]}(t',u)
			\, dt'
			+ 
				\mathring{\upepsilon}
				\frac{1}{\upmu_{\star}^{3/2}(t,u)}.
				\notag
	\end{align}

\end{lemma}

\begin{proof}
		We first prove \eqref{E:MOSTDIFFICULTERRORINTEGRALBOUND}.
		We treat only the first integral on the LHS since the second one can be
		treated using identical arguments.
		To proceed, we first use Cauchy-Schwarz and \eqref{E:WAVECOERCIVENESSOFCONTROLLING} 
		to bound it by
\begin{align} \label{E:NEWPROOFFIRSTSTEPDIFFICULTINTEGRALBOUND}
		\leq 2 
		\int_{t' = 0}^t 
			\sqrt{\totmax{N}}(t',u) 
			\left\| 
				(\Rad v^1) \GeoAng^N \mytr \upchi 
			\right\|_{L^2(\Sigma_{t'}^u)} 
		\, dt'.
\end{align}
We now substitute the estimate \eqref{E:DIFFICULTTERML2BOUND}
(with $t$ in \eqref{E:DIFFICULTTERML2BOUND} replaced by $t'$)
for the second factor in the integrand \eqref{E:NEWPROOFFIRSTSTEPDIFFICULTINTEGRALBOUND}.
Following this substitution,
the desired bound of 
$\mbox{RHS~\eqref{E:NEWPROOFFIRSTSTEPDIFFICULTINTEGRALBOUND}}$
by $\leq \mbox{RHS~\eqref{E:MOSTDIFFICULTERRORINTEGRALBOUND}}$
follows easily with the help of simple estimates of the form
$ab \lesssim a^2 + b^2$, 
the fact that $\totmax{N}$
is increasing in its arguments,
and the estimate
\eqref{E:LOSSKEYMUINTEGRALBOUND},
which we use to 
bound the error integral
$
\displaystyle
\mathring{\upepsilon}
				\int_{t'=0}^t
					\frac{1}{\upmu_{\star}^{3/2}(t',u)}
					\sqrt{\totmax{N}}(t',u)
				\, dt'
$
by 
$
\displaystyle
\lesssim
\mathring{\upepsilon}^2
				\int_{t'=0}^t
					\frac{1}{\upmu_{\star}^{5/2}(t',u)}
				\, dt'
+
			\int_{t'=0}^t
					\frac{1}{\upmu_{\star}^{1/2}(t',u)}
					\totmax{N}(t',u)
	\, dt'
\lesssim
\mathring{\upepsilon}^2
				\frac{1}{\upmu_{\star}^{3/2}(t,u)}
	+
			\int_{t'=0}^t
					\frac{1}{\upmu_{\star}^{1/2}(t',u)}
				\totmax{N}(t',u)
	\, dt'
$.

The proof of \eqref{E:LESSDEGENERATEDIFFICULTERRORINTEGAL}
is similar but simpler. To bound the first integral 
on LHS~\eqref{E:LESSDEGENERATEDIFFICULTERRORINTEGAL},
we first argue as above to deduce that it is bounded by
RHS~\eqref{E:NEWPROOFFIRSTSTEPDIFFICULTINTEGRALBOUND},
but with $\Rad v^1$ on the RHS replaced by $\Rad v^2$.
We then use the estimate \eqref{E:LESSDEGENERATEDIFFICULTTERML2BOUND}
in place of the estimate \eqref{E:DIFFICULTTERML2BOUND}
used in the proof of \eqref{E:MOSTDIFFICULTERRORINTEGRALBOUND}.
The remainder of the proof now proceeds as in the proof 
of \eqref{E:MOSTDIFFICULTERRORINTEGRALBOUND}.
The remaining three integrals on LHS~\eqref{E:LESSDEGENERATEDIFFICULTERRORINTEGAL}
can be bounded in the same way.
\end{proof}

\subsection{\texorpdfstring{$L^2$}{Square integral} bounds for less degenerate top-order error integrals in terms of 
\texorpdfstring{$\totmax{[1,N]}$}{the fundamental controlling quantities}}
\label{SS:LESSDEGENERATEENERGYESTIMATEINTEGRALS}
In this section, we bound some top-order error integrals
for which we need the fully modified
quantities defined in Sect.~\ref{S:MODIFIED}
to avoid losing a derivative. 
However, the error integrals
contain a helpful factor of $\upmu$.
For this reason, the estimates are easier to derive and less degenerate.

\begin{lemma}[\textbf{Bounds for less degenerate top-order error integrals}]
	\label{L:LESSDEGENERATEENERGYESTIMATEINTEGRALS}
	Assume that $\Psi \in \lbrace \Densrenormalized  - v^1,v^1,v^2 \rbrace$
	and that $N = 20$.
	Recall that $\GeoAngFlatRadComponent$
	is the scalar-valued function appearing in Lemma~\ref{L:GEOANGDECOMPOSITION}.
	Under the data-size and bootstrap assumptions 
	of Sects.~\ref{SS:FLUIDVARIABLEDATAASSUMPTIONS}-\ref{SS:PSIBOOTSTRAP}
	and the smallness assumptions of Sect.~\ref{SS:SMALLNESSASSUMPTIONS}, 
	the following integral estimates hold for $(t,u) \in [0,\Tboot) \times [0,U_0]$:
	\begin{subequations}
	\begin{align} \label{E:FIRSTLESSDEGENERATEENERGYESTIMATEINTEGRALS}
		&
		\left|
			\int_{\mathcal{M}_{t,u}}
				\GeoAngFlatRadComponent
				(\Rad \Fullset_{\ast}^{N;\leq 1} \Psi)	 
				(\Rad \Psi) 
				(\angdiffuparg{\#} \Psi)
				\cdot
				\myarray[\upmu \angdiff \GeoAng^{N-1} \mytr \upchi]
					{\upmu \angdiff \GeoAng^{N-2} \Rad \mytr \upchi}
			\, d \vol
		\right|
			\\
		& \lesssim
			\int_{t'=0}^t
				\left\lbrace 
					\ln \upmu_{\star}^{-1}(t',u) 
					+ 
					1 
				\right\rbrace^2
				\totmax{[1,N]}(t',u)
			\, dt'
			+
			\int_{t'=0}^t
				\Vorttotmax{\leq N+1}(t',u)
			\, dt'
			+ 
			\mathring{\upepsilon}^2,
			\notag
\end{align}
\begin{align}
	&
		\left|
			\int_{\mathcal{M}_{t,u}}
				(1 + 2 \upmu)
				\GeoAngFlatRadComponent 
				(\Lunit \Fullset_{\ast}^{N;\leq 1} \Psi)	 
				(\Rad \Psi)
				(\angdiffuparg{\#} \Psi)
				\cdot
				\myarray[\upmu \angdiff \GeoAng^{N-1} \mytr \upchi]
					{\upmu \angdiff \GeoAng^{N-2} \Rad \mytr \upchi}
			\, d \vol
		\right|
			\label{E:SECONDLESSDEGENERATEENERGYESTIMATEINTEGRALS} \\
		& \lesssim
			\int_{t'=0}^t
				\left\lbrace 
					\ln \upmu_{\star}^{-1}(t',u) 
					+ 
					1 
				\right\rbrace^2
				\totmax{[1,N]}(t',u)
			\, dt'
			+
			\int_{u'=0}^u
				\totmax{[1,N]}(t,u')
			\, du'
				\notag \\
	& \ \
		+
		\int_{t'=0}^t
			\Vorttotmax{\leq N+1}(t',u)
		\, dt'
		+ \mathring{\upepsilon}^2.
			\notag
	\end{align}
	\end{subequations}
\end{lemma}
\begin{proof}
	We prove \eqref{E:SECONDLESSDEGENERATEENERGYESTIMATEINTEGRALS}
	only for the first product on the LHS since the proof for the second term is identical.
	To proceed,
	we first use the schematic identity \eqref{E:LINEARLYSMALLSCALARSDEPENDINGONGOODVARIABLES}
	for $\GeoAngFlatRadComponent$,
	the $L^{\infty}$ estimates of Prop.~\ref{P:IMPROVEMENTOFAUX},
	Young's inequality,
	Lemma~\ref{L:COERCIVENESSOFCONTROLLING},
	and inequality \eqref{E:LESSPRECISEDIFFICULTTERML2BOUND}
	to obtain
	\begin{align} \label{E:FIRSTPROOFSTEPSECONDLESSDEGENERATEENERGYESTIMATEINTEGRALS}
		&
		\left|
			\int_{\mathcal{M}_{t,u}}
				(1 + 2 \upmu)
				\GeoAngFlatRadComponent 
				(\Lunit \Fullset_{\ast}^{N;\leq 1} \Psi)	 
				(\Rad \Psi)
				(\angdiffuparg{\#} \Psi)
				\cdot
				(\upmu \angdiff \GeoAng^{N-1} \mytr \upchi)
			\, d \vol
		\right|
			\\
		& \lesssim 
			\int_{u'=0}^u
				\left\|
					\Lunit \Fullset_{\ast}^{N;\leq 1} \Psi 
				\right\|_{L^2(\mathcal{P}_{u'}^t)}^2
			\, du'
			+
			\int_{t'=0}^t
				\left\|
					\upmu \GeoAng^N \mytr \upchi
				\right\|_{L^2(\Sigma_{t'}^u)}^2
			\, dt'
				\notag \\
		& \lesssim
			\int_{u'=0}^u
				\totmax{[1,N]}(t,u')
			\, du'
		+
		\int_{t'=0}^t
			\totmax{[1,N]}(t',u)
		\, dt'
		+
		\int_{t'=0}^t
			\left\lbrace
				\int_{s=0}^{t'}
					\frac{1} 
					{\upmu_{\star}(s,u)} 
					\sqrt{\totmax{N}}(s,u)
				\, ds
			\right\rbrace^2
		\, dt'
			\notag \\
		& \ \
			+
			\int_{t'=0}^t
				\left\lbrace
					\int_{s=0}^{t'}
						\sqrt{\VorttotTanmax{N+1}}(s,u)
					\, ds
				\right\rbrace^2
			\, dt'
			+
			\int_{t'=0}^t
				\left\lbrace
					\int_{s=0}^{t'}
						\frac{1} 
						{\upmu_{\star}^{1/2}(s,u)}
						\sqrt{\VorttotTanmax{\leq N}}(s,u)
					\, ds
				\right\rbrace^2
			\, dt'
			\notag	\\
		& \ \
			+
			\mathring{\upepsilon}^2
			\int_{t'=0}^t
				\left\lbrace 
					\ln \upmu_{\star}^{-1}(t,u) 
					+ 
					1 
				\right\rbrace^2
			\, dt'.
			\notag
	\end{align}
	Using inequalities
	\eqref{E:LOGLOSSMUINVERSEINTEGRALBOUND}
	and 
	\eqref{E:LESSSINGULARTERMSMPOINTNINEINTEGRALBOUND}
	and the fact that
	$\totmax{M}$
	and $\Vorttotmax{M}$ are increasing in their arguments,
	we conclude that
	RHS~\eqref{E:FIRSTPROOFSTEPSECONDLESSDEGENERATEENERGYESTIMATEINTEGRALS}
	is $\lesssim$ RHS~\eqref{E:SECONDLESSDEGENERATEENERGYESTIMATEINTEGRALS}
	as desired.

	The proof of \eqref{E:FIRSTLESSDEGENERATEENERGYESTIMATEINTEGRALS}
	is similar, except in the first step we bound
	LHS~\eqref{E:FIRSTLESSDEGENERATEENERGYESTIMATEINTEGRALS} by 
	$
	\displaystyle
	\lesssim
	\int_{t'=0}^t
				\left\|
					\Rad \Fullset_{\ast}^{N;\leq 1} \Psi
				\right\|_{L^2(\Sigma_{t'}^u)}^2
			\, dt'
			+
			\int_{t'=0}^t
				\left\|
					\upmu \GeoAng^N \mytr \upchi
				\right\|_{L^2(\Sigma_{t'}^u)}^2
			\, dt'
	$.
\end{proof}

\subsection{Error integrals requiring integration by parts with respect to \texorpdfstring{$\Lunit$}{the rescaled null vectorfield}}
\label{SS:ERROINTEGRALSINVOLVINGIBPL}
In deriving our top-order energy estimates for the wave variables
$\lbrace \Densrenormalized  - v^1,v^1,v^2 \rbrace$,
we encounter some difficult error integrals that we 
can control only by integrating by parts with respect to $\Lunit$;
see the error integral \eqref{E:LUNITINVOLVINGDIFFICULTERRORINTEGRAL}
and the discussion below it.
It turns out that in carrying out this procedure,
we must use the partially modified quantities 
of Sect.~\ref{S:MODIFIED}
in order to avoid generating error terms that are too large to control.
This results in the presence of two types of error integrals,
which we bound in this section: spacetime error integrals,
some of which involve the partially modified
quantities, and $\Sigma_t^u$ ``boundary'' error integrals, 
some of which also involve the partially modified quantities.
We remark that we treat the most difficult of these error integrals 
in Lemmas~\ref{L:ANNOYINGBOUNDRYSPATIALINTEGRALFACTORL2ESTIMATE} and
\ref{L:BOUNDSFORDIFFICULTTOPORDERINTEGRALSINVOLVINGLUNITIBP}.

\begin{lemma}[\textbf{A difficult top-order hypersurface} $L^2$ \textbf{estimate}]
	\label{L:ANNOYINGBOUNDRYSPATIALINTEGRALFACTORL2ESTIMATE}
	Assume that $N=20$
	and let $\Fullset_{\ast}^{N-1;\leq 1} \in \lbrace \GeoAng^{N-1}, \GeoAng^{N-2} \Rad \rbrace$.
	Let $\upchipartialmodarg{\Fullset_{\ast}^{N-1;\leq 1}}$
	be the corresponding partially modified quantity defined by \eqref{E:TRANSPORTPARTIALRENORMALIZEDTRCHIJUNK}.
	There exist constants $C > 0$ and $C_* > 0$ such that
	under the data-size and bootstrap assumptions 
	of Sects.~\ref{SS:FLUIDVARIABLEDATAASSUMPTIONS}-\ref{SS:PSIBOOTSTRAP}
	and the smallness assumptions of Sect.~\ref{SS:SMALLNESSASSUMPTIONS},
	the following $L^2$ estimate holds for $(t,u) \in [0,\Tboot) \times [0,U_0]$:
	\begin{subequations}
	\begin{align} \label{E:ANNOYINGLDERIVATIVEBOUNDRYSPATIALINTEGRALFACTORL2ESTIMATE}
		\left\|
			\frac{1}{\sqrt{\upmu}} (\Rad v^1) \Lunit \upchipartialmodarg{\Fullset_{\ast}^{N-1;\leq 1}}
		\right\|_{L^2(\Sigma_t^u)}
		& \leq
			\boxed{\sqrt{2}}
			\frac{\left\| [\Lunit \upmu]_- \right\|_{C^0(\Sigma_t^u)}}
			{\upmu_{\star}(t,u)}
			\sqrt{\totmax{N}}(t,u) 
				\\
	& \ \
			+
			C_*
			\frac{1}
			{\upmu_{\star}(t,u)}
			\sqrt{\easytotmax{N}}(t,u) 
			\notag \\
		& \ \
			+ C \frac{1}{\upmu_{\star}^{1/2}(t,u)} \sqrt{\totmax{[1,N]}}(t,u)
			+ C \varepsilon \frac{1}{\upmu_{\star}(t,u)} \sqrt{\totmax{[1,N]}}(t,u)
				\notag \\
		& \ \
			+
			C 
			\frac{1}
			{\upmu_{\star}(t,u)}
			\sqrt{\totmax{[1,N-1]}}(t,u) 
			+ C \mathring{\upepsilon} \frac{1}{\upmu_{\star}^{1/2}(t,u)},
			\notag
	\end{align}
	\begin{align}
		\left\|
			\frac{1}{\sqrt{\upmu}} (\Rad v^1) \upchipartialmodarg{\Fullset_{\ast}^{N-1;\leq 1}}
		\right\|_{L^2(\Sigma_t^u)}
		& \leq
			\boxed{\sqrt{2}}
			\left\| \Lunit \upmu \right\|_{L^{\infty}(\Sigmaminus{t}{t}{u})}
			\frac{1}{\upmu_{\star}^{1/2}(t,u)}
			\int_{t'=0}^t
				\frac{1}{\upmu_{\star}^{1/2}(t',u)} \sqrt{\totmax{[1,N]}}(t',u)
			\, dt'
				\label{E:ANNOYINGBOUNDRYSPATIALINTEGRALFACTORL2ESTIMATE} \\
		& \ \
			+
			C_*
			\frac{1}{\upmu_{\star}^{1/2}(t,u)}
			\int_{t'=0}^t
				\frac{1}{\upmu_{\star}^{1/2}(t',u)} \sqrt{\easytotmax{[1,N]}}(t',u)
			\, dt'
			\notag \\
		& \ \
			+ 
			C 
			\int_{t'=0}^t
				\frac{1}{\upmu_{\star}^{1/2}(t',u)} \sqrt{\totmax{[1,N]}}(t',u)
			\, dt'
				\notag \\
		& \ \
			+ 
			C \varepsilon
			\frac{1}{\upmu_{\star}^{1/2}(t,u)}
			\int_{t'=0}^t
				\frac{1}{\upmu_{\star}^{1/2}(t',u)} \sqrt{\totmax{[1,N]}}(t',u)
			\, dt'
			\notag \\
	& \ \
		+ C \mathring{\upepsilon}
			\frac{1}{\upmu_{\star}^{1/2}(t,u)}.
			\notag
	\end{align}
	\end{subequations}

	Moreover, we have the following less precise estimates:
	\begin{subequations}
	\begin{align}  \label{E:NOTTOOHARDLUNITAPPLIEDTOBOUNDRYSPATIALINTEGRALFACTORL2ESTIMATE}
		\left\|
			\Lunit \upchipartialmodarg{\Fullset_{\ast}^{N-1;\leq 1}}
		\right\|_{L^2(\Sigma_t^u)}
		& \lesssim
			\frac{1}{\upmu_{\star}^{1/2}(t,u)} 
			\sqrt{\totmax{[1,N]}}(t,u)
			+
			\mathring{\upepsilon},
				\\
		\left\|
			\upchipartialmodarg{\Fullset_{\ast}^{N-1;\leq 1}}
		\right\|_{L^2(\Sigma_t^u)}
		& \lesssim
			\int_{t'=0}^t
				\frac{1}{\upmu_{\star}^{1/2}(t',u)} \sqrt{\totmax{[1,N]}}(t',u)
			\, dt'
			+
			\mathring{\upepsilon}.
		\label{E:MUCHEASIERBOUNDRYSPATIALINTEGRALFACTORL2ESTIMATE}
	\end{align}
	\end{subequations}
\end{lemma}
\begin{proof}
	We start by proving \eqref{E:ANNOYINGLDERIVATIVEBOUNDRYSPATIALINTEGRALFACTORL2ESTIMATE}.
	We multiply inequality \eqref{E:SHARPPOINTWISELUNITPARTIALLYMODIFIED} by
	$
	\displaystyle
	\frac{1}{\sqrt{\upmu}} \Rad v^1.
	$
	We first consider the difficult product generated by the first term on 
	RHS~\eqref{E:SHARPPOINTWISELUNITPARTIALLYMODIFIED}.
	To proceed, we multiply the identity \eqref{E:MAINDIFFICULTINTEGRALFACTORDECOMPOSITION}
	by $\sqrt{\upmu}$
	and use the schematic identity \eqref{E:GFRAMESCALARSDEPENDINGONGOODVARIABLES}
	and the $L^{\infty}$ estimates of Prop.~\ref{P:IMPROVEMENTOFAUX}
	(in particular \eqref{E:PSITRANSVERSALLINFINITYBOUNDBOOTSTRAPIMPROVEDSMALL}
	and \eqref{E:CRUCIALPSITRANSVERSALLINFINITYBOUNDBOOTSTRAPIMPROVEDSMALL})
	to obtain
	\begin{align} \label{E:CRUCIALPRODUCTMAINTERMPLUSERRORS}
		\frac{1}{2}
		\frac{1}{\sqrt{\upmu}}
		\left|
			\sum_{\imath=0}^1
			G_{\Lunit \Lunit}^{\imath}
			\Rad v^1
		\right|
		& 
			=
			\frac{\Lunit \upmu}{\sqrt{\upmu}}
			+
			\mathcal{O}(\varepsilon) \frac{1}{\sqrt{\upmu}}.
	\end{align}
	Inserting \eqref{E:CRUCIALPRODUCTMAINTERMPLUSERRORS} into the product of 
	$
	\displaystyle
	\frac{1}{\sqrt{\upmu}} \Rad v^1
	$
	and the first product on RHS~\eqref{E:SHARPPOINTWISELUNITPARTIALLYMODIFIED},
	we obtain the terms
	\begin{align} \label{E:DIFFICULTPRODUCTANNOYINGLDERIVATIVEBOUNDRYSPATIALINTEGRALFACTORL2ESTIMATE}
			\frac{[\Lunit \upmu]_-}{\upmu}
			\left|
				\sqrt{\upmu} \angLap \Fullset_{\ast}^{N-1;\leq 1} v^1
			\right|
			+
			\frac{[\Lunit \upmu]_+}{\upmu}
			\left|
				\sqrt{\upmu} \angLap \Fullset_{\ast}^{N-1;\leq 1} v^1
			\right|
			+
			\frac{\mathcal{O}(\varepsilon)}{\upmu}
			\left|
				\sqrt{\upmu} \angLap \Fullset_{\ast}^{N-1;\leq 1} v^1
			\right|.
	\end{align}
	Using Lemma~\ref{L:COERCIVENESSOFCONTROLLING}, we see that
	the norm $\| \cdot \|_{L^2(\Sigma_t^u)}$
	of the first product in
	\eqref{E:DIFFICULTPRODUCTANNOYINGLDERIVATIVEBOUNDRYSPATIALINTEGRALFACTORL2ESTIMATE}
	is bounded by the first term on
	RHS~\eqref{E:ANNOYINGLDERIVATIVEBOUNDRYSPATIALINTEGRALFACTORL2ESTIMATE}
	as desired.
	Next, we again use Lemma~\ref{L:COERCIVENESSOFCONTROLLING}
	and the estimate \eqref{E:POSITIVEPARTOFLMUOVERMUISBOUNDED}
	to bound the norm $\| \cdot \|_{L^2(\Sigma_t^u)}$
	of the second and third products on 
	RHS~\eqref{E:DIFFICULTPRODUCTANNOYINGLDERIVATIVEBOUNDRYSPATIALINTEGRALFACTORL2ESTIMATE}
	by the terms on the third line of
	RHS~\eqref{E:ANNOYINGLDERIVATIVEBOUNDRYSPATIALINTEGRALFACTORL2ESTIMATE}
	as desired.
	In proving the remaining estimates,
	we use \eqref{E:PSITRANSVERSALLINFINITYBOUNDBOOTSTRAPIMPROVEDLARGE} to bound
	$
	\displaystyle
	\left\|
		\frac{1}{\sqrt{\upmu}} \Rad v^1
	\right\|_{L^{\infty}(\Sigma_t^u)}
	\leq
	C
	\frac{1}{\sqrt{\upmu_{\star}(t,u)}}
	$
	and thus it remains for us to bound
	the norm $\| \cdot \|_{L^2(\Sigma_t^u)}$
	of the remaining two terms on RHS~\eqref{E:SHARPPOINTWISELUNITPARTIALLYMODIFIED}
	and to multiply those bounds by
	$
	\displaystyle
	C \frac{1}{\sqrt{\upmu_{\star}(t,u)}}
	$.
	To handle the product generated by the
	second term on RHS~\eqref{E:SHARPPOINTWISELUNITPARTIALLYMODIFIED},
	we use Lemma~\ref{L:COERCIVENESSOFCONTROLLING},
	which implies that its 
	norm $\| \cdot \|_{L^2(\Sigma_t^u)}$
	is bounded by the $\easytotmax{N}$-involving term
	on RHS~\eqref{E:ANNOYINGLDERIVATIVEBOUNDRYSPATIALINTEGRALFACTORL2ESTIMATE}
	(which has the coefficient $C_*$).
	To bound the product generated by the
	next-to-last term $C \varepsilon \cdots$ on RHS~\eqref{E:SHARPPOINTWISELUNITPARTIALLYMODIFIED},
	we use Lemma~\ref{L:COERCIVENESSOFCONTROLLING}
	(the product under consideration is bounded by $\leq$ the 
	term
	$
	\displaystyle
	C \varepsilon \frac{1}{\upmu_{\star}(t,u)} \sqrt{\totmax{[1,N]}}(t,u)
	$ on RHS~\eqref{E:ANNOYINGLDERIVATIVEBOUNDRYSPATIALINTEGRALFACTORL2ESTIMATE}).
	To bound the product generated by the
	last term on RHS~\eqref{E:SHARPPOINTWISELUNITPARTIALLYMODIFIED},
	we use the estimate \eqref{E:TANGENGITALEIKONALINTERMSOFCONTROLLING},
	inequality \eqref{E:LESSSINGULARTERMSMPOINTNINEINTEGRALBOUND},
	and the fact that
	$\totmax{[1,N]}$
	is increasing in its arguments.
	We have thus proved \eqref{E:ANNOYINGLDERIVATIVEBOUNDRYSPATIALINTEGRALFACTORL2ESTIMATE}.

	We now prove \eqref{E:ANNOYINGBOUNDRYSPATIALINTEGRALFACTORL2ESTIMATE}.
	We multiply inequality \eqref{E:SHARPPOINTWISEPARTIALLYMODIFIED} by
	$
	\displaystyle
	\frac{1}{\sqrt{\upmu}} \Rad v^1.
	$
	The most difficult product is generated by the second term on 
	RHS~\eqref{E:SHARPPOINTWISEPARTIALLYMODIFIED}:
	\begin{align} \label{E:DIFFICULTTERMSHARPPOINTWISEPARTIALLYMODIFIED}
		\frac{1}{2} 
		\frac{1}{\sqrt{\upmu}}
		\left|
			\sum_{\imath=0}^1
			G_{\Lunit \Lunit}^{\imath}
			\Rad v^1
		\right|
		(t,u,\vartheta)
			\int_{t'=0}^t
			\left|
				\angLap \Fullset_{\ast}^{N-1;\leq 1} v^1
			\right|
			(t',u,\vartheta)
			\, dt'.
	\end{align}
We now substitute RHS~\eqref{E:CRUCIALPRODUCTMAINTERMPLUSERRORS}
for the product
$
\displaystyle
\frac{1}{\sqrt{\upmu}}
\left|
			\sum_{\imath=0}^1
			G_{\Lunit \Lunit}^{\imath}
			\Rad v^1
		\right|
$
in \eqref{E:DIFFICULTTERMSHARPPOINTWISEPARTIALLYMODIFIED}
and take the norm $\| \cdot \|_{L^2(\Sigma_t^u)}$
of the resulting expression.
With the help of 
Lemma~\ref{L:L2NORMSOFTIMEINTEGRATEDFUNCTIONS} and
Lemma~\ref{L:COERCIVENESSOFCONTROLLING},
we see that the norm
$\| \cdot \|_{L^2(\Sigma_t^u)}$
of the product generated by the second product
on RHS~\eqref{E:CRUCIALPRODUCTMAINTERMPLUSERRORS}
is bounded by
the next-to-last term
$
C \varepsilon
	\cdots
$
on RHS~\eqref{E:ANNOYINGBOUNDRYSPATIALINTEGRALFACTORL2ESTIMATE}.
To handle the remaining product
(corresponding to the term
$
\displaystyle
\frac{\Lunit \upmu}{\sqrt{\upmu}}
$
on RHS~\eqref{E:CRUCIALPRODUCTMAINTERMPLUSERRORS}),
we first decompose $\Sigma_t^u = \Sigmaplus{s}{t}{u} \cup \Sigmaminus{s}{t}{u}$
as in \eqref{E:SIGMASSPLIT} 
and again use Lemmas~\ref{L:L2NORMSOFTIMEINTEGRATEDFUNCTIONS} and \ref{L:COERCIVENESSOFCONTROLLING} 
as well as the simple estimate
$\| \Lunit \upmu \|_{L^2(\Sigma_t^u)}
\lesssim 1
$
(see \eqref{E:LUNITUPTOONETRANSVERSALUPMULINFINITY})
to bound it by
\begin{align} \label{E:DIFFICULTPARTANNOYINGBOUNDRYSPATIALINTEGRALFACTORL2ESTIMATE}
	& \leq 
			\sqrt{2}
			\left\| \Lunit \upmu \right\|_{L^{\infty}(\Sigmaminus{t}{t}{u})}
			\frac{1}{\upmu_{\star}^{1/2}(t,u)}
			\int_{t'=0}^t
				\frac{1}{\upmu_{\star}^{1/2}(t',u)} \sqrt{\totmax{[1,N]}}(t',u)
			\, dt'
				\\
& \ \
			+
			\sqrt{2}
			\left\| \frac{\Lunit \upmu}{\sqrt{\upmu}} \right\|_{L^{\infty}(\Sigmaplus{t}{t}{u})}
			\int_{t'=0}^t
				\frac{1}{\upmu_{\star}^{1/2}(t',u)} \sqrt{\totmax{[1,N]}}(t',u)
			\, dt'
			\notag	\\
	& \ \
			+ 
			C \varepsilon
			\frac{1}{\upmu_{\star}^{1/2}(t,u)}
			\int_{t'=0}^t
				\frac{1}{\upmu_{\star}^{1/2}(t',u)} \sqrt{\totmax{[1,N]}}(t',u)
			\, dt'.
			\notag
\end{align}
The first and third products on RHS~\eqref{E:DIFFICULTPARTANNOYINGBOUNDRYSPATIALINTEGRALFACTORL2ESTIMATE} 
are manifestly bounded by RHS~\eqref{E:ANNOYINGBOUNDRYSPATIALINTEGRALFACTORL2ESTIMATE}.
To bound the second product on RHS~\eqref{E:DIFFICULTPARTANNOYINGBOUNDRYSPATIALINTEGRALFACTORL2ESTIMATE}
by RHS~\eqref{E:ANNOYINGBOUNDRYSPATIALINTEGRALFACTORL2ESTIMATE},
we need only to use the following estimate to bound the factor multiplying the time integral:
\begin{align} \label{E:LMUSQRTMUOVERSQRTMUSIGMAPLUSBOUND}
	\left\| \frac{\Lunit \upmu}{\sqrt{\upmu}} \right\|_{L^{\infty}(\Sigmaplus{t}{t}{u})}
	\leq
	C
	\left\| \frac{[\Lunit \upmu]_+}{\upmu} \right\|_{L^{\infty}(\Sigmaplus{t}{t}{u})}
	+
	C
	\left\| \frac{[\Lunit \upmu]_-}{\upmu} \right\|_{L^{\infty}(\Sigmaplus{t}{t}{u})}
	\leq
	C.
\end{align}
The estimate \eqref{E:LMUSQRTMUOVERSQRTMUSIGMAPLUSBOUND}
is a straightforward consequence of 
the estimates
\eqref{E:UPTOONETRANSVERSALDERIVATIVEUPMULINFTY} and
\eqref{E:LUNITUPTOONETRANSVERSALUPMULINFINITY} (with $M=0$),
\eqref{E:POSITIVEPARTOFLMUOVERMUISBOUNDED},
and
\eqref{E:KEYMUNOTDECAYINGMINUSPARTLMUOVERMUBOUND}.
	We now bound the norm $\| \cdot \|_{L^2(\Sigma_t^u)}$
	of the product of
	$
	\displaystyle
	\frac{1}{\sqrt{\upmu}} \Rad v^1
	$
	and the remaining four terms on
	RHS~\eqref{E:SHARPPOINTWISEPARTIALLYMODIFIED}.
	In all of the remaining estimates,
	we rely on the bound
	$
	\displaystyle
	\left\|
		\frac{1}{\sqrt{\upmu}} \Rad v^1
	\right\|_{L^{\infty}(\Sigma_t^u)}
	\leq
	C
	\frac{1}{\sqrt{\upmu_{\star}(t,u)}}
	$
	noted in the proof of \eqref{E:ANNOYINGLDERIVATIVEBOUNDRYSPATIALINTEGRALFACTORL2ESTIMATE};
	it therefore remains for us to bound
	the norm $\| \cdot \|_{L^2(\Sigma_t^u)}$
	of the remaining four terms on RHS~\eqref{E:SHARPPOINTWISEPARTIALLYMODIFIED}
	and to multiply those bounds by
	$
	\displaystyle
	C \frac{1}{\sqrt{\upmu_{\star}(t,u)}}
	$.
	To bound the product corresponding to 
	the first term 
	$
	\displaystyle
	\left|
		\upchipartialmodarg{\Fullset_{\ast}^{N;\leq 1}}
	\right|
		(0,u,\vartheta)
	$
	on RHS~\eqref{E:SHARPPOINTWISEPARTIALLYMODIFIED},
	we first use 
	\eqref{E:SIGMATVOLUMEFORMCOMPARISON} with $s=0$
	to deduce
	$
	\displaystyle
	\left\|
	\left|
		\upchipartialmodarg{\GeoAng^{N-1}}
	\right|
	(1,\cdot)
\right\|_{L^2(\Sigma_t^u)}
\lesssim
\left\|
	\upchipartialmodarg{\GeoAng^{N-1}}
\right\|_{L^2(\Sigma_0^u)}
$.
Next, from definition \eqref{E:TRANSPORTPARTIALRENORMALIZEDTRCHIJUNK},
the simple inequality
$|G_{(Frame)}| = |\smoothfunction(\GdVar,\angdiff \vec{x})| \lesssim 1$
(which follows from Lemmas~\ref{L:SCHEMATICDEPENDENCEOFMANYTENSORFIELDS}
and \ref{L:POINTWISEFORRECTANGULARCOMPONENTSOFVECTORFIELDS}
and the $L^{\infty}$ estimates of Prop.~\ref{P:IMPROVEMENTOFAUX}),
the estimate \eqref{E:POINTWISEESTIMATESFORCHIANDITSDERIVATIVES},
and our assumptions on the data,
we find that 
$
\left\|
	\upchipartialmodarg{\Fullset_{\ast}^{N;\leq 1}}
\right\|_{L^2(\Sigma_0^u)}
\lesssim \mathring{\upepsilon}
$.
In total, we conclude that the product under consideration 
is bounded in the norm
$
\left\|
	\cdot
\right\|_{L^2(\Sigma_t^u)}
$
by the last term on RHS~\eqref{E:ANNOYINGBOUNDRYSPATIALINTEGRALFACTORL2ESTIMATE} 
as desired.
To bound the norm $\| \cdot \|_{L^2(\Sigma_t^u)}$ of the second time integral
	$C_* \cdots$
	on RHS~\eqref{E:SHARPPOINTWISEPARTIALLYMODIFIED},
	we use Lemmas~\ref{L:L2NORMSOFTIMEINTEGRATEDFUNCTIONS} and \ref{L:COERCIVENESSOFCONTROLLING}.
	Multiplying by 
	$
	\displaystyle
	C \frac{1}{\sqrt{\upmu_{\star}(t,u)}}
	$,
	we find that the term of interest 
	is bounded by the second term $C_* \cdots$ 
	on RHS~\eqref{E:ANNOYINGBOUNDRYSPATIALINTEGRALFACTORL2ESTIMATE}.
	Similarly, we see that the product generated by the time integral
	$C \varepsilon \cdots$
	on RHS~\eqref{E:SHARPPOINTWISEPARTIALLYMODIFIED}
	is bounded by $\leq$
	the $C \varepsilon \cdots$ term on
	RHS~\eqref{E:ANNOYINGBOUNDRYSPATIALINTEGRALFACTORL2ESTIMATE}.
	To bound the product generated by the last time integral
	on RHS~\eqref{E:SHARPPOINTWISEPARTIALLYMODIFIED},
	we use a similar argument together with 
	\eqref{E:TANGENGITALEIKONALINTERMSOFCONTROLLING},
	except that as a preliminary step, we bound
	the time integral on RHS~\eqref{E:TANGENGITALEIKONALINTERMSOFCONTROLLING}
	by $\lesssim \sqrt{\totmax{[1,N]}}(t,u)$
	with the help of \eqref{E:LESSSINGULARTERMSMPOINTNINEINTEGRALBOUND}.
	We have thus proved \eqref{E:ANNOYINGBOUNDRYSPATIALINTEGRALFACTORL2ESTIMATE}.

	The proofs of 
	\eqref{E:NOTTOOHARDLUNITAPPLIEDTOBOUNDRYSPATIALINTEGRALFACTORL2ESTIMATE}
	and
	\eqref{E:MUCHEASIERBOUNDRYSPATIALINTEGRALFACTORL2ESTIMATE}
	are based on a subset of the above arguments and are much
	simpler; we therefore omit the details,
	noting only that the main simplification is that we 
	do \emph{not} have to rely on the estimate
	\eqref{E:CRUCIALPRODUCTMAINTERMPLUSERRORS},
	which played a fundamental role in our proofs of
	\eqref{E:ANNOYINGLDERIVATIVEBOUNDRYSPATIALINTEGRALFACTORL2ESTIMATE}
	and 
	\eqref{E:ANNOYINGBOUNDRYSPATIALINTEGRALFACTORL2ESTIMATE}.
	\end{proof}

\begin{lemma}[\textbf{Bounds connected to easy top-order 
error integrals requiring integration by parts with respect to} $\Lunit$]
	\label{L:HARMLESSIBPERRORINTEGRALS}
	Let $\Psi \in \lbrace \Densrenormalized  - v^1,v^1,v^2 \rbrace$.
	Assume that $N = 20$ and $\varsigma > 0$.
	For $i=1,2$, 
	let
	$\mbox{\upshape Error}_i[\Fullset_{\ast}^{N;\leq 1} \Psi; 
		\upchipartialmodarg{\Fullset_{\ast}^{N-1;\leq 1}}]$
	be the error integrands defined in
	\eqref{E:LUNITIBPSPACETIMEERROR}
	and
	\eqref{E:LUNITIBPHYPERSURFACEERROR},
	where the partially modified quantity
	$\upchipartialmodarg{\Fullset_{\ast}^{N-1;\leq 1}}$ 
	defined in \eqref{E:TRANSPORTPARTIALRENORMALIZEDTRCHIJUNK} is in role of $\ThirdSmoothFunction$
	and we are assuming no relationship between the operators
	$\Fullset_{\ast}^{N;\leq 1}$ and $\Fullset_{\ast}^{N-1;\leq 1}$
	(see Sect.~\ref{SS:STRINGSOFCOMMUTATIONVECTORFIELDS} regarding the vectorfield operator notation).
	Under the data-size and bootstrap assumptions 
	of Sects.~\ref{SS:FLUIDVARIABLEDATAASSUMPTIONS}-\ref{SS:PSIBOOTSTRAP}
	and the smallness assumptions of Sect.~\ref{SS:SMALLNESSASSUMPTIONS},
	the following estimates hold for $(t,u) \in [0,\Tboot) \times [0,U_0]$,
	where the implicit constants are independent of $\varsigma$:
	\begin{subequations}
	\begin{align} \label{E:HARMLESSIBPSPACETIMEERRORINTEGRALS}
		\int_{\mathcal{M}_{t,u}}
			& 
			\left|
				\mbox{\upshape Error}_1[\Fullset_{\ast}^{N;\leq 1} \Psi; \upchipartialmodarg{\Fullset_{\ast}^{N-1;\leq 1} }]
			\right|
		\, d \vol
			\\
		& \lesssim 
			(1 + \varsigma^{-1})
			\int_{s=0}^t
				\frac{1}{\upmu_{\star}^{1/2}(s,u)} 
				\totmax{[1,N]}(s,u)
			\, ds
			+
			\int_{s=0}^t
				\frac{1}{\upmu_{\star}^{3/2}(s,u)} 
				\totmax{[1,N-1]}(s,u)
			\, ds
			\notag \\
		& \ \
			+
			\varsigma \coercivespacetimemax{[1,N]}(t,u)
			+ (1 + \varsigma^{-1}) \mathring{\upepsilon}^2,
				\notag 
		\end{align}
		\begin{align}
		\int_{\Sigma_t^u}
		\left|
			\mbox{\upshape Error}_2[\Fullset_{\ast}^{N;\leq 1} \Psi; \upchipartialmodarg{\Fullset_{\ast}^{N-1;\leq 1} }]
		\right|
		\, d \vol
		& \lesssim 
			\mathring{\upepsilon}^2
			+ \varepsilon \totmax{[1,N]}(t,u)
			\label{E:HARMLESSIBPHYPERSURFACEERRORINTEGRALS},
			\\
		\int_{\Sigma_0^u}
		\left|
			\mbox{\upshape Error}_2[\Fullset_{\ast}^{N;\leq 1} \Psi; \upchipartialmodarg{\Fullset_{\ast}^{N-1;\leq 1} }]
		\right|
		\, d \vol
		& \lesssim \mathring{\upepsilon}^2,
			\label{E:DATAHARMLESSIBPHYPERSURFACEERRORINTEGRALS}
				\\
		\int_{\Sigma_0^u}
			\left|
				(1 + 2 \upmu) (\Rad \Psi) (\GeoAng \Tanset^N \Psi) \upchipartialmodarg{\Fullset_{\ast}^{N-1;\leq 1} }
			\right|
		\, d \vol
		& \lesssim
			\mathring{\upepsilon}^2.
			\label{E:QUADRATICDATAHARMLESSIBPHYPERSURFACEERRORINTEGRALS}
	\end{align}
	\end{subequations}
\end{lemma}

\begin{proof}
	See Sect.~\ref{SS:OFTENUSEDESTIMATES} for some comments on the analysis.
	We first prove \eqref{E:HARMLESSIBPSPACETIMEERRORINTEGRALS}.
	All products on RHS~\eqref{E:LUNITIBPSPACETIMEERROR}
	contain a quadratic factor of
	$(\angdiff \Fullset_{\ast}^{N;\leq 1} \Psi) \upchipartialmodarg{\Fullset_{\ast}^{N-1;\leq 1}}$,
	$(\GeoAng \Fullset_{\ast}^{N;\leq 1} \Psi) \upchipartialmodarg{\Fullset_{\ast}^{N-1;\leq 1}}$,
	$(\Fullset_{\ast}^{N;\leq 1} \Psi) \upchipartialmodarg{\Fullset_{\ast}^{N-1;\leq 1} }$,
	or 
	$(\Fullset_{\ast}^{N;\leq 1} \Psi) \Lunit \upchipartialmodarg{\Fullset_{\ast}^{N-1;\leq 1} }$.
	Using inequalities
	\eqref{E:NOSPECIALSTRUCTUREDIFFERENTIATEDGEOANGDEFORMSPHERELSHARPPOINTWISE}
	and
	\eqref{E:NOSPECIALSTRUCTUREDIFFERNTIATEDANGDEFORMTANGENTPOINTWISE}
	and the $L^{\infty}$ estimates of Prop.~\ref{P:IMPROVEMENTOFAUX},
	we find that the remaining factors in the products are bounded in 
	the norm
	$\| \cdot \|_{L^{\infty}(\Sigma_t^u)}$ by $\lesssim 1$.
	Hence, it suffices to bound the 
	magnitude of the spacetime integrals of the four quadratic terms
	by $\lesssim$ RHS~\eqref{E:HARMLESSIBPSPACETIMEERRORINTEGRALS}.
	To bound the spacetime integral of 
	$
	\left|
		(\GeoAng \Fullset_{\ast}^{N;\leq 1} \Psi) \upchipartialmodarg{\Fullset_{\ast}^{N-1;\leq 1}} 
	\right|
	$, 
	we use spacetime Cauchy-Schwarz,
	Lemmas~\ref{L:KEYSPACETIMECOERCIVITY} and \ref{L:COERCIVENESSOFCONTROLLING},
	inequalities
	\eqref{E:LESSSINGULARTERMSMPOINTNINEINTEGRALBOUND}
	and 
	\eqref{E:MUCHEASIERBOUNDRYSPATIALINTEGRALFACTORL2ESTIMATE},
	simple estimates of the form
	$ab \lesssim a^2 + b^2$,
	and the fact that
	$\totmax{[1,N]}$
	is increasing in its arguments
	to deduce
	\begin{align} \label{E:SPACETIMEINTEGRALBOUNDFORINTEGRATIONBYPARTISINL}
		\int_{\mathcal{M}_{t,u}}
			&
			\left|
				(\GeoAng \Fullset_{\ast}^{N;\leq 1} \Psi) \upchipartialmodarg{\Fullset_{\ast}^{N-1;\leq 1} }
			\right|
		\, d \vol
			\\  
		& 
		\lesssim
		\varsigma 
		\TranminusdatasizeWithFactor
		\int_{\mathcal{M}_{t,u}}
			\left|
				\angdiff \Fullset_{\ast}^{N;\leq 1} \Psi
			\right|^2
		\, d \vol
		+
		\varsigma^{-1} \TranminusdatasizeWithFactor^{-1}
		\int_{s=0}^t
			\left\|
				\upchipartialmodarg{\Fullset_{\ast}^{N-1;\leq 1} }
			\right\|_{L^2(\Sigma_s^u)}^2
		\, ds
			\notag \\
		& \lesssim
			\varsigma
			\coercivespacetimemax{[1,N]}(t,u)
			+
			\varsigma^{-1}
			\int_{s=0}^t
				\left\lbrace
				\int_{t'=0}^s
					\frac{1}{\upmu_{\star}^{1/2}(t',u)} \totmax{[1,N]}^{1/2}(t',u)
				\, dt'
			  \right\rbrace^2
			  +
			  \varsigma^{-1}
			  \mathring{\upepsilon}^2
			\, ds
			\notag
			\\
		& \lesssim
			\varsigma
			\coercivespacetimemax{[1,N]}(t,u)
			+
			\varsigma^{-1}
			\int_{s=0}^t
				\totmax{[1,N]}(s,u)
			\, ds
			+
			\varsigma^{-1} \mathring{\upepsilon}^2,
			\notag
	\end{align}
	which is $\leq$ RHS~\eqref{E:HARMLESSIBPSPACETIMEERRORINTEGRALS} as desired.
	We clarify that in passing to the last inequality in 
	\eqref{E:SPACETIMEINTEGRALBOUNDFORINTEGRATIONBYPARTISINL},
	we have used the fact that $\totmax{[1,N]}$ is increasing in its arguments
	and the estimate \eqref{E:LESSSINGULARTERMSMPOINTNINEINTEGRALBOUND}
	to deduce that
	$
	\displaystyle
	\int_{t'=0}^s
		\frac{1}{\upmu_{\star}^{1/2}(t',u)} \totmax{[1,N]}^{1/2}(t',u)
	\, dt'
	\lesssim 
	\totmax{[1,N]}^{1/2}(s,u)
	$,
	as we did in passing to the last line of \eqref{E:SECONDHARMLESSEXAMPLEINTEGRAL}.

	The spacetime integral of
	$
	\left|
		(\angdiff \Fullset_{\ast}^{N;\leq 1} \Psi) \upchipartialmodarg{\GeoAng^{N-1}} 
	\right|
	$
	can be bounded in the same way.

	The spacetime integral of
	$\left| (\Fullset_{\ast}^{N;\leq 1} \Psi) \upchipartialmodarg{\Fullset_{\ast}^{N-1;\leq 1} } \right|$
	can be bounded by $\leq$ RHS~\eqref{E:HARMLESSIBPSPACETIMEERRORINTEGRALS}
	by using essentially the same arguments; we omit the details.

	To bound the spacetime integral of
	$
	\left|
		(\Fullset_{\ast}^{N;\leq 1} \Psi) \Lunit \upchipartialmodarg{\Fullset_{\ast}^{N-1;\leq 1} }
	\right|
	$, 
	by $\leq$ RHS~\eqref{E:HARMLESSIBPSPACETIMEERRORINTEGRALS},
	we first use Cauchy-Schwarz,
	Lemma~\ref{L:COERCIVENESSOFCONTROLLING},
	and inequality \eqref{E:NOTTOOHARDLUNITAPPLIEDTOBOUNDRYSPATIALINTEGRALFACTORL2ESTIMATE}
	to deduce
	\begin{align} \label{E:GEOANGNLUNITPARTIALLMODIFIEDSPACETIMEINTEGRAL}
		\int_{\mathcal{M}_{t,u}}
			\left|
				(\Fullset_{\ast}^{N;\leq 1} \Psi) \Lunit \upchipartialmodarg{\Fullset_{\ast}^{N-1;\leq 1} }
			\right|
		\, d \vol
		& \lesssim
			\int_{s=0}^t
				\left\|
					\Fullset_{\ast}^{N;\leq 1} \Psi
				\right\|_{L^2(\Sigma_s^u)}
				\left\|
					\Lunit \upchipartialmodarg{\Fullset_{\ast}^{N-1;\leq 1} }
				\right\|_{L^2(\Sigma_s^u)}
		\, ds
			\\
		& \lesssim
			\int_{s=0}^t
				\frac{1}{\upmu_{\star}(s,u)} 
				\totmax{[1,N-1]}^{1/2}(s,u)
				\totmax{[1,N]}^{1/2}(s,u)
			\, ds
			\notag \\
		& \ \
			+
			\mathring{\upepsilon}
			\int_{s=0}^t
				\frac{1}{\upmu_{\star}^{1/2}(s,u)} 
				\totmax{[1,N-1]}^{1/2}(s,u)
			\, ds.
			\notag
		\end{align}
	Finally, using simple estimates of the form $ab \lesssim a^2 + b^2$,
	the estimate \eqref{E:LESSSINGULARTERMSMPOINTNINEINTEGRALBOUND},
	and the fact that
	$\totmax{[1,N]}$
	is increasing in its arguments,
	we bound RHS~\eqref{E:GEOANGNLUNITPARTIALLMODIFIEDSPACETIMEINTEGRAL} by $\lesssim$ RHS
	\eqref{E:HARMLESSIBPSPACETIMEERRORINTEGRALS} as desired. This concludes the proof of 
	\eqref{E:HARMLESSIBPSPACETIMEERRORINTEGRALS}.

	We now prove \eqref{E:HARMLESSIBPHYPERSURFACEERRORINTEGRALS}
	and \eqref{E:DATAHARMLESSIBPHYPERSURFACEERRORINTEGRALS}.
	Using the estimate \eqref{E:NOSPECIALSTRUCTUREDIFFERNTIATEDANGDEFORMTANGENTPOINTWISE}
	and the $L^{\infty}$ estimates of Prop.~\ref{P:IMPROVEMENTOFAUX},
	we see that RHS~\eqref{E:LUNITIBPHYPERSURFACEERROR}
	is bounded in magnitude by 
	$\lesssim 
	\varepsilon 
		\left|
			\Fullset_{\ast}^{N;\leq 1} \Psi 
		\right|
		\left|
			\upchipartialmodarg{\Fullset_{\ast}^{N-1;\leq 1} } 
		\right|
	$.
	Next, we use Cauchy-Schwarz on $\Sigma_t^u$,
	Lemma~\ref{L:COERCIVENESSOFCONTROLLING}, 
	\eqref{E:MUCHEASIERBOUNDRYSPATIALINTEGRALFACTORL2ESTIMATE},
	and the estimate \eqref{E:LESSSINGULARTERMSMPOINTNINEINTEGRALBOUND}
	to deduce that
	\begin{align} \label{E:NOTSODIFFICULTSIGMATUERRORINTEGRAL}
	\varepsilon 
		\int_{\Sigma_t^u}
			\left|
				\Fullset_{\ast}^{N;\leq 1} \Psi 
			\right|
			\left|
				\upchipartialmodarg{\Fullset_{\ast}^{N-1;\leq 1} } 
			\right|
		\, d \tvol
		& \lesssim
		\varepsilon
		\left\|
			\Fullset_{\ast}^{N;\leq 1} \Psi 
		\right\|_{L^2(\Sigma_t^u)}
		\left\|
			\upchipartialmodarg{\Fullset_{\ast}^{N-1;\leq 1} }  
		\right\|_{L^2(\Sigma_t^u)}
		 \\
		&
		\lesssim
		\varepsilon
		\left\lbrace
			\totmax{[1,N]}^{1/2}(t,u)
			+
			\mathring{\upepsilon}
		\right\rbrace^2
		\lesssim \mbox{RHS~\eqref{E:HARMLESSIBPHYPERSURFACEERRORINTEGRALS}}
		\notag
	\end{align}
	as desired. 
	We clarify that in passing to the second inequality of \eqref{E:NOTSODIFFICULTSIGMATUERRORINTEGRAL},
	we have used \eqref{E:LESSSINGULARTERMSMPOINTNINEINTEGRALBOUND}
	and the fact that $\totmax{[1,N]}$ is increasing in its arguments to bound
	the time integral on RHS~\eqref{E:MUCHEASIERBOUNDRYSPATIALINTEGRALFACTORL2ESTIMATE} by
	$\lesssim \totmax{[1,N]}^{1/2}(t,u)$.
	\eqref{E:DATAHARMLESSIBPHYPERSURFACEERRORINTEGRALS} then follows from
	\eqref{E:HARMLESSIBPHYPERSURFACEERRORINTEGRALS} with $t=0$
	and Lemma~\ref{L:INITIALSIZEOFL2CONTROLLING}.

	The proof of \eqref{E:QUADRATICDATAHARMLESSIBPHYPERSURFACEERRORINTEGRALS} is similar.
	The difference is that the
	$L^{\infty}$ estimates of Prop.~\ref{P:IMPROVEMENTOFAUX}  
	imply only that LHS~\eqref{E:QUADRATICDATAHARMLESSIBPHYPERSURFACEERRORINTEGRALS}
	is 
	$
	\lesssim
	\int_{\Sigma_0^u}
			\left|
				\Tanset^{N+1} \Psi 
			\right|
			\left|
				\upchipartialmodarg{\Fullset_{\ast}^{N-1;\leq 1} } 
			\right|
		\, d \tvol
	$,
	without a gain of a factor $\varepsilon$.
	However, this integral is quadratically small in the data-size
	parameter $\mathring{\upepsilon}$,
	as is easy to verify using the arguments given in the previous paragraph. 
	We have thus proved \eqref{E:QUADRATICDATAHARMLESSIBPHYPERSURFACEERRORINTEGRALS}.
\end{proof}

\begin{lemma}[\textbf{Bounds for difficult top-order spacetime error integrals connected to integration by parts involving $\Lunit$}]
	\label{L:BOUNDSFORDIFFICULTTOPORDERINTEGRALSINVOLVINGLUNITIBP}
	Assume that $N =20$ and $\varsigma > 0$.
	Let 
	$\upchipartialmodarg{\GeoAng^{N-1}}$
	and
	$\upchipartialmodarg{\GeoAng^{N-2} \Rad}$ be the
	partially modified quantities
	defined in \eqref{E:TRANSPORTPARTIALRENORMALIZEDTRCHIJUNK}.
	There exists a constant $C > 0$,
	independent of $\varsigma$,
	such that under the data-size and bootstrap assumptions 
	of Sects.~\ref{SS:FLUIDVARIABLEDATAASSUMPTIONS}-\ref{SS:PSIBOOTSTRAP}
	and the smallness assumptions of Sect.~\ref{SS:SMALLNESSASSUMPTIONS}, 
	the following integral estimates hold for $(t,u) \in [0,\Tboot) \times [0,U_0]$
	(see Sect.~\ref{SS:STRINGSOFCOMMUTATIONVECTORFIELDS} regarding the vectorfield operator notation):
	\begin{align} \label{E:DIFFICULTLUNITSPACETIMEIBPINTEGRALBOUND}
		&
		\left|
			\int_{\mathcal{M}_{t,u}}
				(1 + 2 \upmu) (\Rad v^1) (\GeoAng \Fullset_{\ast}^{N;\leq 1} v^1) 
				\Lunit \upchipartialmodarg{\GeoAng^{N-1}}
			\, d \vol
		\right|,
			\\
		&
		\left|
			\int_{\mathcal{M}_{t,u}}
				(1 + 2 \upmu) (\Rad v^1) (\GeoAng \Fullset_{\ast}^{N;\leq 1} v^1) 
				\Lunit \upchipartialmodarg{\GeoAng^{N-2} \Rad}
			\, d \vol
		\right|
			\notag \\
		& \leq
			\boxed{2}
			\int_{t'=0}^t
					\frac{\left\| [\Lunit \upmu]_- \right\|_{L^{\infty}(\Sigma_{t'}^u)}} 
							 {\upmu_{\star}(t',u)} 
				  \totmax{[1,N]}(t',u)
				\, dt'
			\notag	\\
		& \ \
			+
			C \varepsilon
			\int_{t'=0}^t
				\frac{1} 
					{\upmu_{\star}(t',u)} 
				 \totmax{[1,N]}(t',u)
			\, dt'
			+
			C
			\int_{t'=0}^t
				\frac{1} 
							 {\upmu_{\star}^{1/2}(t',u)} 
				 \totmax{[1,N]}(t',u)
			\, dt'
			\notag	\\
		& \ \
			+ C \mathring{\upepsilon}^2,
			\notag
		\end{align}

		\begin{align}
		&
		\left|
		\int_{\Sigma_t^u}
				(1 + 2 \upmu) (\Rad v^1) (\GeoAng \Fullset_{\ast}^{N;\leq 1} v^1) 
				\upchipartialmodarg{\GeoAng^{N-1}}
		\, d \vol
		\right|,
			\label{E:DIFFICULTLUNITHYPERSURFACEIBPINTEGRALBOUND} \\
		&
		\left|
		\int_{\Sigma_t^u}
				(1 + 2 \upmu) (\Rad v^1) (\GeoAng \Fullset_{\ast}^{N;\leq 1} v^1) 
				\upchipartialmodarg{\GeoAng^{N-2} \Rad}
		\, d \vol
		\right|
		\notag \\
		& \leq 
			\boxed{2}
			\frac{1}{\upmu_{\star}^{1/2}(t,u)}
			\sqrt{\totmax{[1,N]}}(t,u)
			\left\| \Lunit \upmu \right\|_{L^{\infty}(\Sigmaminus{t}{t}{u})}
			\int_{t'=0}^t
				\frac{1}{\upmu_{\star}^{1/2}(t',u)} \sqrt{\totmax{[1,N]}}(t',u)
			\, dt'
				\notag \\
		& \ \
			+ 
			C \varepsilon
			\frac{1}{\upmu_{\star}^{1/2}(t,u)}
			\sqrt{\totmax{[1,N]}}(t,u)
			\int_{t'=0}^t
				\frac{1}{\upmu_{\star}^{1/2}(t',u)} \sqrt{\totmax{[1,N]}}(t',u)
			\, dt'
			\notag
				\\
		& \ \
			+ 
			C 
			\sqrt{\totmax{[1,N]}}(t,u)
			\int_{t'=0}^t
				\frac{1}{\upmu_{\star}^{1/2}(t',u)} \sqrt{\totmax{[1,N]}}(t',u)
			\, dt'
				\notag \\
		&  \ \
			+ C \varsigma \totmax{[1,N]}(t,u)
			+ C \varsigma^{-1} \mathring{\upepsilon}^2 \frac{1}{\upmu_{\star}(t,u)}.
			\notag
	\end{align}

	Moreover, we have the following less degenerate estimates:
	\begin{align} \label{E:LESSDEGENERATEDIFFICULTLUNITSPACETIMEIBPINTEGRALBOUND}
		&
		\left|
			\int_{\mathcal{M}_{t,u}}
				(1 + 2 \upmu) 
				\myarray[(\Rad v^2) (\GeoAng \Fullset_{\ast}^{N;\leq 1} v^2) ]
					{\left\lbrace
					\Rad (\Densrenormalized  - v^1)
				\right\rbrace
				\left\lbrace
					\GeoAng \Fullset_{\ast}^{N;\leq 1} (\Densrenormalized  - v^1)
				\right\rbrace}
				\Lunit \upchipartialmodarg{\GeoAng^{N-1}}
			\, d \vol
		\right|,
			\\
		&
		\left|
			\int_{\mathcal{M}_{t,u}}
				(1 + 2 \upmu) 
				\myarray[(\Rad v^2) (\GeoAng \Fullset_{\ast}^{N;\leq 1} v^2) ]
					{\left\lbrace
					\Rad (\Densrenormalized  - v^1)
				\right\rbrace
				\left\lbrace
					\GeoAng \Fullset_{\ast}^{N;\leq 1} (\Densrenormalized  - v^1)
				\right\rbrace}
				\Lunit \upchipartialmodarg{\GeoAng^{N-2} \Rad}
			\, d \vol
		\right|
			\notag \\
		& \leq
			C \varepsilon
			\int_{t'=0}^t
				\frac{1} 
					{\upmu_{\star}(t',u)} 
				 \totmax{[1,N]}(t',u)
			\, dt'
			+ C \mathring{\upepsilon}^2,
			\notag
	\end{align}

	\begin{align}
		&
		\left|
		\int_{\Sigma_t^u}
				(1 + 2 \upmu) 
				\myarray[(\Rad v^2) (\GeoAng \Fullset_{\ast}^{N;\leq 1} v^2) ]
					{\left\lbrace
					\Rad (\Densrenormalized  - v^1)
				\right\rbrace
				\left\lbrace
					\GeoAng \Fullset_{\ast}^{N;\leq 1} (\Densrenormalized  - v^1)
				\right\rbrace}
				\upchipartialmodarg{\GeoAng^{N-1}}
		\, d \vol
		\right|,
			\label{E:LESSDEGENERATEDIFFICULTLUNITHYPERSURFACEIBPINTEGRALBOUND} \\
			&
		\left|
		\int_{\Sigma_t^u}
				(1 + 2 \upmu) 
				\myarray[(\Rad v^2) (\GeoAng \Fullset_{\ast}^{N;\leq 1} v^2) ]
					{\left\lbrace
					\Rad (\Densrenormalized  - v^1)
				\right\rbrace
				\left\lbrace
					\GeoAng \Fullset_{\ast}^{N;\leq 1} (\Densrenormalized  - v^1)
				\right\rbrace}
				\upchipartialmodarg{\GeoAng^{N-2} \Rad}
		\, d \vol
		\right| 
			\notag \\
		& \leq 
			C \varepsilon
			\frac{1}{\upmu_{\star}^{1/2}(t,u)}
			\sqrt{\totmax{[1,N]}}(t,u)
			\int_{t'=0}^t
				\frac{1}{\upmu_{\star}^{1/2}(t',u)} \sqrt{\totmax{[1,N]}}(t',u)
			\, dt'
			\notag
				\\
		& + C \varepsilon \totmax{[1,N]}(t,u)
			+ C \mathring{\upepsilon}^2 \frac{1}{\upmu_{\star}(t,u)}.
			\notag
	\end{align}

	\end{lemma}

	\begin{proof}
	We prove \eqref{E:DIFFICULTLUNITSPACETIMEIBPINTEGRALBOUND}
	only for the first term on the LHS since the second term
	can be treated in an identical fashion.
	To proceed, we first use Cauchy-Schwarz, 
	and the estimates
		$
		\left|
			\GeoAng
		\right|
		\leq 1 + C \varepsilon
		$,
		$
		\left\| 
			\Rad \Psi
		\right\|_{L^{\infty}(\Sigma_t^u)} 
		\lesssim 1
		$,
		and
		$
		\left\| 
			\upmu
		\right\|_{L^{\infty}(\Sigma_t^u)} 
		\lesssim 1
		$
		(which follow from \eqref{E:GEOANGPOINTWISE} and
		the $L^{\infty}$ estimates of Prop.~\ref{P:IMPROVEMENTOFAUX})
		to bound the LHS by 
		\begin{align} \label{E:ANGDIFFTIMESPSIPARTIALMODDIFFICULTTOPORDERFIRSTESTIMATE}
		&	\leq
		(1 + C \varepsilon)
		\int_{t'=0}^t
		\left\|
			\sqrt{\upmu} \angdiff \Fullset_{\ast}^{N;\leq 1} \Psi
		\right\|_{L^2(\Sigma_{t'}^u)}
		\left\|
			\frac{1}{\sqrt{\upmu}} (\Rad \Psi) \Lunit \upchipartialmodarg{\GeoAng^{N-1}}
		\right\|_{L^2(\Sigma_{t'}^u)}
		\, dt'
			\\
		& \	\
			+ 
		C 
		\int_{t'=0}^t
			\left\|
				\sqrt{\upmu} \angdiff \Fullset_{\ast}^{N;\leq 1} \Psi
			\right\|_{L^2(\Sigma_{t'}^u)}
			\left\|
				\Lunit \upchipartialmodarg{\GeoAng^{N-1}}
			\right\|_{L^2(\Sigma_{t'}^u)}
		\, dt'.
	\notag
	\end{align}
	The desired estimate \eqref{E:DIFFICULTLUNITSPACETIMEIBPINTEGRALBOUND}
	now follows from \eqref{E:ANGDIFFTIMESPSIPARTIALMODDIFFICULTTOPORDERFIRSTESTIMATE},
	Lemma~\ref{L:COERCIVENESSOFCONTROLLING},
	and inequalities \eqref{E:ANNOYINGLDERIVATIVEBOUNDRYSPATIALINTEGRALFACTORL2ESTIMATE}
	and \eqref{E:NOTTOOHARDLUNITAPPLIEDTOBOUNDRYSPATIALINTEGRALFACTORL2ESTIMATE}.
	We clarify that to bound the integral
	$
	\displaystyle
	\int_{t'=0}^t
		C \mathring{\upepsilon}
		\frac{1}{\upmu_{\star}^{1/2}(t',u)}
		\totmax{[1,N]}^{1/2}(t',u)
	\, dt'
	$,
	which is generated by the last term 
	on RHS~\eqref{E:ANNOYINGLDERIVATIVEBOUNDRYSPATIALINTEGRALFACTORL2ESTIMATE},
	we first use Young's inequality to bound the integrand by
	$
	\displaystyle
	\lesssim
	\frac{\mathring{\upepsilon}^2}{\upmu_{\star}^{1/2}(t',u)}
	+ 
	\frac{\totmax{[1,N]}(t',u)}{\upmu_{\star}^{1/2}(t',u)}
	$.
	We then bound the time integral of the first term in the previous expression by 
	$\lesssim \mathring{\upepsilon}^2
	$ 
	with the help of the estimate \eqref{E:LESSSINGULARTERMSMPOINTNINEINTEGRALBOUND}
	and the time integral of the second by 
	$\leq$ the third term on RHS~\eqref{E:DIFFICULTLUNITSPACETIMEIBPINTEGRALBOUND}.

	The proof of \eqref{E:LESSDEGENERATEDIFFICULTLUNITSPACETIMEIBPINTEGRALBOUND}
	is similar but simpler and is based on the estimates
	$
	\left\| 
		\Rad v^2 \right 
	\|_{L^{\infty}(\Sigma_t^u)}
	\lesssim \varepsilon
	$
	and
	$
	\left\| 
		\Rad (\Densrenormalized  - v^1) 
	\right\|_{L^{\infty}(\Sigma_t^u)}
	\lesssim \varepsilon
	$
	(see \eqref{E:PSITRANSVERSALLINFINITYBOUNDBOOTSTRAPIMPROVEDSMALL} and \eqref{E:CRUCIALPSITRANSVERSALLINFINITYBOUNDBOOTSTRAPIMPROVEDSMALL})
	and the estimate \eqref{E:NOTTOOHARDLUNITAPPLIEDTOBOUNDRYSPATIALINTEGRALFACTORL2ESTIMATE};
	we omit the details.

	The proof of \eqref{E:DIFFICULTLUNITHYPERSURFACEIBPINTEGRALBOUND}
	is similar to the proof of \eqref{E:DIFFICULTLUNITSPACETIMEIBPINTEGRALBOUND} 
	but relies on 
	\eqref{E:ANNOYINGBOUNDRYSPATIALINTEGRALFACTORL2ESTIMATE}
	and \eqref{E:MUCHEASIERBOUNDRYSPATIALINTEGRALFACTORL2ESTIMATE}
	in place of 
	\eqref{E:ANNOYINGLDERIVATIVEBOUNDRYSPATIALINTEGRALFACTORL2ESTIMATE}
	and
	\eqref{E:NOTTOOHARDLUNITAPPLIEDTOBOUNDRYSPATIALINTEGRALFACTORL2ESTIMATE};
	we omit the details, noting only that we encounter the term
	$
	\displaystyle
	C \mathring{\upepsilon}
				\frac{1}{\upmu_{\star}^{1/2}(t,u)}
				\totmax{[1,N]}^{1/2}(t,u)
	$
	generated by the last term on RHS~\eqref{E:ANNOYINGBOUNDRYSPATIALINTEGRALFACTORL2ESTIMATE}.
	We bound this term by using Young's inequality as follows:
	$
	\displaystyle
	C \mathring{\upepsilon}
				\frac{1}{\upmu_{\star}^{1/2}(t,u)}
				\totmax{[1,N]}^{1/2}(t,u)
	\leq 
	C \varsigma^{-1} \mathring{\upepsilon}^2 \frac{1}{\upmu_{\star}(t,u)}
	+ C \varsigma \totmax{[1,N]}(t,u)
	$.

	The proof of \eqref{E:LESSDEGENERATEDIFFICULTLUNITHYPERSURFACEIBPINTEGRALBOUND}
	is similar to the proof of \eqref{E:DIFFICULTLUNITHYPERSURFACEIBPINTEGRALBOUND}
	but is simpler. It is based on the estimates
	$
	\left\| 
		\Rad v^2 \right 
	\|_{L^{\infty}(\Sigma_t^u)}
	\lesssim \varepsilon
	$
	and
	$
	\left\| 
		\Rad (\Densrenormalized  - v^1)
	\right\|_{L^{\infty}(\Sigma_t^u)}
	\lesssim \varepsilon
	$
	noted above
	and the estimate \eqref{E:MUCHEASIERBOUNDRYSPATIALINTEGRALFACTORL2ESTIMATE};
	we omit the details.

	\end{proof}

\subsection{Estimates for the most degenerate top-order transport equation error integrals}
\label{SS:MOSTDEGENERATETRANSPORTEQUATIONERRORINTEGRALS}
In the next lemma, we bound the most degenerate error integrals appearing
in the top-order energy estimates for the specific vorticity,
which are generated by the main terms from
Prop.~\ref{P:VORTICITYIDOFKEYDIFFICULTENREGYERRORTERMS}.
These error integrals are responsible for 
the large blowup-exponent $6.4$ in the factor
$\upmu_{\star}^{-6.4}(t,u)$
on RHS~\eqref{E:TOPVORTMULOSSMAINAPRIORIENERGYESTIMATES}.

\begin{lemma}[\textbf{Estimates for the most degenerate top-order transport equation error integrals}]
	\label{L:MOSTDEGENERATETOPORDERTRANSPORTEQUATIONERRORINTEGRALS}
	Assume that $N = 20$
	and recall that $\GeoAngFlatRadComponent$ is the scalar-valued
	appearing in Lemma~\ref{L:GEOANGDECOMPOSITION}.
	Under the data-size and bootstrap assumptions 
	of Sects.~\ref{SS:FLUIDVARIABLEDATAASSUMPTIONS}-\ref{SS:PSIBOOTSTRAP}
	and the smallness assumptions of Sect.~\ref{SS:SMALLNESSASSUMPTIONS}, 
	the following integral estimates hold for $(t,u) \in [0,\Tboot) \times [0,U_0]$:
	\begin{align} \label{E:MOSTDEGENERATETOPORDERTRANSPORTEQUATIONERRORINTEGRALS}
		&
		\left|
			\int_{\mathcal{M}_{t,u}}
				\threemyarray[(\GeoAng \Vortrenormalized) \GeoAng^{N-1} \Rad \mytr \upchi]
					{g(\GeoAng,\GeoAng) (\Lunit \Vortrenormalized) \GeoAng^{N-1} \Rad  \mytr \upchi}
					{\GeoAngFlatRadComponent (\GeoAng \Vortrenormalized) \GeoAng^{N-1} \Rad \mytr \upchi}
				\Tanset^{N+1} \Vortrenormalized
			\, d \vol
		\right|
			\\
	& \lesssim
	\varepsilon^2 
	\frac{1}{\upmu_{\star}(t,u)}
	\totmax{N}(t,u)
	+
	\varepsilon^2
	\int_{t'=0}^t
				\frac{1}{\upmu_{\star}^2(t',u)} 
				\left\lbrace
					\int_{s=0}^{t'}
						\frac{1}{\upmu_{\star}(s,u)} 
						\sqrt{\totmax{N}}(s,u)
					\, ds
				\right\rbrace^2
	\, dt'
		\notag \\
	& \ \
			+
			\varepsilon^2
			\int_{t'=0}^t
				\frac{1}{\upmu_{\star}^2(t',u)}
				\left\lbrace
					\int_{s=0}^{t'}
						\sqrt{\VorttotTanmax{N+1}}(s,u)
					\, ds
				\right\rbrace^2
				\, dt'
			\notag
				\\
	 & \ \
			+
			\varepsilon^2
			\int_{t'=0}^t
				\frac{1}{\upmu_{\star}^2(t',u)}
				\left\lbrace
					\int_{s=0}^{t'}
						\frac{1} 
						{\upmu_{\star}^{1/2}(s,u)}
						\sqrt{\VorttotTanmax{\leq N}}(s,u)
					\, ds
				\right\rbrace^2
				\, dt'
				\notag \\
	& \ \
		+
		\int_{u'=0}^u
			\VorttotTanmax{N+1}(t,u')
		\, du'
		+
			\mathring{\upepsilon}^2
			\frac{1}{\upmu_{\star}^{3/2}(t,u)}.
	\notag
\end{align}
\end{lemma}

\begin{proof}
We prove \eqref{E:MOSTDEGENERATETOPORDERTRANSPORTEQUATIONERRORINTEGRALS}
for the term
$
	\left|
		\int_{\mathcal{M}_{t,u}}
			(\GeoAng \Vortrenormalized) (\GeoAng^{N-1} \Rad \mytr \upchi)
			\Tanset^{N+1} \Vortrenormalized
		\, d \vol
	\right|
$
in detail.
The other two error integrals
on LHS~\eqref{E:MOSTDEGENERATETOPORDERTRANSPORTEQUATIONERRORINTEGRALS}
can be handled using nearly identical arguments,
the schematic relations \eqref{E:SCALARSDEPENDINGONGOODVARIABLES}
and \eqref{E:LINEARLYSMALLSCALARSDEPENDINGONGOODVARIABLES},
and the $L^{\infty}$ estimates of Prop.~\ref{P:IMPROVEMENTOFAUX};
we omit those details.
To proceed, we use the bound
$
\| 
	\GeoAng \Vortrenormalized
\|_{L^{\infty}(\Sigma_t^u)}
\lesssim \varepsilon
$
(see \eqref{E:VORTICITYUPTOTWOTRANSVERSALLINFTY})
and Young's inequality
to deduce that the error integral under consideration is
\begin{align} \label{E:FIRSTBOUNDMOSTDEGENERATETOPORDERTRANSPORTEQUATIONERRORINTEGRALS}
	& \lesssim 
		\varepsilon^2
		\int_{\mathcal{M}_{t,u}}
			(\GeoAng^{N-1} \Rad \mytr \upchi)^2
		\, d \vol
		+
		\int_{\mathcal{M}_{t,u}}
			(\Tanset^{N+1} \Vortrenormalized)^2
		\, d \vol
			\\
	& \lesssim 
		\varepsilon^2
		\int_{t'=0}^t
			\| 
				\GeoAng^{N-1} \Rad \mytr \upchi 
			\|_{L^2(\Sigma_{t'}^u)}^2
		\, dt'
		+
		\int_{u'=0}^u
			\| 
				\Tanset^{N+1} \Vortrenormalized
			\|_{L^2(\mathcal{P}_{u'}^t)}^2
		\, du'.
			\notag 
\end{align}
Using Lemma~\ref{L:COERCIVENESSOFCONTROLLING}, we bound
the last integral on RHS~\eqref{E:FIRSTBOUNDMOSTDEGENERATETOPORDERTRANSPORTEQUATIONERRORINTEGRALS}
by 
$
\lesssim
\int_{u'=0}^u
	\VorttotTanmax{N+1}(t,u')
\, du'
$
as desired.
To handle the remaining time integral 
on RHS~\eqref{E:FIRSTBOUNDMOSTDEGENERATETOPORDERTRANSPORTEQUATIONERRORINTEGRALS},
we use the estimate \eqref{E:LESSPRECISEDIFFICULTTERML2BOUND}
to bound it as follows:
\begin{align} \label{E:SECONDBOUNDMOSTDEGENERATETOPORDERTRANSPORTEQUATIONERRORINTEGRALS}
		\varepsilon^2
		\int_{t'=0}^t
			\| 
				\GeoAng^{N-1} \Rad \mytr \upchi 
			\|_{L^2(\Sigma_{t'}^u)}^2
		\, dt'
	& \lesssim 
		\varepsilon^2
		\int_{t'=0}^t
			\frac{1}{\upmu_{\star}^2(t',u)}
			\| 
				\upmu \GeoAng^{N-1} \Rad  \mytr \upchi 
			\|_{L^2(\Sigma_{t'}^u)}^2
		\, dt'
			\\
	& \lesssim 
			\varepsilon^2
			\int_{t'=0}^t
				\frac{1}{\upmu_{\star}^2(t',u)}
				\totmax{N}(t',u)
			\, dt'
				\notag \\
	& \ \
	+
	\varepsilon^2
	\int_{t'=0}^t
				\frac{1}{\upmu_{\star}^2(t',u)} 
				\left\lbrace
					\int_{s=0}^{t'}
						\frac{1}{\upmu_{\star}(s,u)} 
						\sqrt{\totmax{N}}(s,u)
					\, ds
				\right\rbrace^2
	\, dt'
		\notag \\
	& \ \
			+
			\varepsilon^2
			\int_{t'=0}^t
				\frac{1}{\upmu_{\star}^2(t',u)}
				\left\lbrace
					\int_{s=0}^{t'}
						\sqrt{\VorttotTanmax{N+1}}(s,u)
					\, ds
				\right\rbrace^2
				\, dt'
			\notag
				\\
	 & \ \
			+
			\varepsilon^2
			\int_{t'=0}^t
				\frac{1}{\upmu_{\star}^2(t',u)}
				\left\lbrace
					\int_{s=0}^{t'}
						\frac{1} 
						{\upmu_{\star}^{1/2}(s,u)}
						\sqrt{\VorttotTanmax{\leq N}}(s,u)
					\, ds
				\right\rbrace^2
				\, dt'
			\notag
				\\
			&  \ \
				+ 
				\varepsilon^2
				\mathring{\upepsilon}^2 
				\int_{t'=0}^t
					\frac{1}{\upmu_{\star}^2(t',u)}
					\left\lbrace 
						\ln \upmu_{\star}^{-1}(t',u) 
						+ 
						1 
					\right\rbrace^2
				\, dt'.
				\notag
\end{align}
Using the fact that $\totmax{N}$ is increasing in its arguments
and the estimate \eqref{E:LOSSKEYMUINTEGRALBOUND}, we
find that
RHS~\eqref{E:SECONDBOUNDMOSTDEGENERATETOPORDERTRANSPORTEQUATIONERRORINTEGRALS}
$\lesssim$
RHS~\eqref{E:MOSTDEGENERATETOPORDERTRANSPORTEQUATIONERRORINTEGRALS}
as desired.
\end{proof}

\subsection{Estimates for transport equation error integrals involving a loss of one derivative}
\label{SS:TRANSPORTEQUATIONENERGYESTIMATESLOSSOFONEDERIVATIVE}
In the next lemma, we estimate some error integrals that arise
when bounding the below-top-order derivatives of the specific vorticity.
We allow the estimates to lose one derivative. The advantage is that
the right-hand sides of the estimates are much less singular with
respect to powers of $\upmu_{\star}^{-1}$ compared to the estimates we would
obtain in an approach that avoids derivative loss. 
This fact is crucially important for our energy estimate descent scheme,
in which the below-top-order energy estimates become successively less singular
with respect to powers of $\upmu_{\star}^{-1}$.

\begin{lemma}[\textbf{Estimates for transport equation error integrals involving a loss of one derivative}]
\label{L:TRANSPORTEQUATIONERRORINTEGRALSLOSEONEDERIVATIVE}
	Assume that $2 \leq N \leq 20$
	and recall that $\GeoAngFlatRadComponent$ is the scalar-valued
	appearing in Lemma~\ref{L:GEOANGDECOMPOSITION}.
	Under the data-size and bootstrap assumptions 
	of Sects.~\ref{SS:FLUIDVARIABLEDATAASSUMPTIONS}-\ref{SS:PSIBOOTSTRAP}
	and the smallness assumptions of Sect.~\ref{SS:SMALLNESSASSUMPTIONS}, 
	the following integral estimates hold for $(t,u) \in [0,\Tboot) \times [0,U_0]$
	(see Sect.~\ref{SS:STRINGSOFCOMMUTATIONVECTORFIELDS} regarding the vectorfield operator notation):
\begin{align}
	& 
	\left|
		\int_{\mathcal{M}_{t,u}}
			\threemyarray[(\GeoAng \Vortrenormalized) \GeoAng^{N-2} \Rad \mytr \upchi]
				{g(\GeoAng,\GeoAng) (\Lunit \Vortrenormalized) \GeoAng^{N-2} \Rad \mytr \upchi}
				{\GeoAngFlatRadComponent (\GeoAng \Vortrenormalized) \GeoAng^{N-2} \Rad \mytr \upchi}
			\Tanset^N \Vortrenormalized
		\, d \vol
	\right|
		\label{E:VORTTRCHIDERIVATIVELOSINGERRORINTEGRALS} \\
	& \lesssim
		\varepsilon^2
		\int_{t'=0}^t
			\left\lbrace
				\int_{s=0}^{t'}
					\frac{\sqrt{\totmax{[1,N]}}(s,u)}{\upmu_{\star}^{1/2}(s,u)}
				\, ds
			\right\rbrace^2
		\, dt'
		+
		\int_{u'=0}^u
			\VorttotTanmax{\leq N}(t,u')
		\, du'
		+
		\varepsilon^2
		\mathring{\upepsilon}^2.
		\notag
\end{align}
\end{lemma}

\begin{proof}
The proof is the same as the proof of Lemma~\ref{L:MOSTDEGENERATETOPORDERTRANSPORTEQUATIONERRORINTEGRALS}
except for one key difference:
we use the estimate \eqref{E:TANGENGITALEIKONALINTERMSOFCONTROLLING}
to bound the term
$
\displaystyle
	\| 
		\GeoAng^{N-2} \Rad \mytr \upchi 
	\|_{L^2(\Sigma_{t'}^u)}^2
$,
in place of the estimate
\eqref{E:LESSPRECISEDIFFICULTTERML2BOUND}
used in bounding the term
$
\displaystyle
	\| 
		\GeoAng^{N-1} \Rad \mytr \upchi 
	\|_{L^2(\Sigma_{t'}^u)}^2
$
on the first line of RHS~\eqref{E:SECONDBOUNDMOSTDEGENERATETOPORDERTRANSPORTEQUATIONERRORINTEGRALS}.

\end{proof}

\subsection{Estimates for wave equation error integrals involving a loss of one derivative}
\label{SS:WAVEEQUATIONENERGYESTIMATESLOSSOFONEDERIVATIVE}
We now provide an analog of Lemma~\ref{L:TRANSPORTEQUATIONERRORINTEGRALSLOSEONEDERIVATIVE} 
for the wave equations.
Specifically, in the next lemma, we estimate some error integrals that arise
when bounding the below-top-order derivatives of the elements of
$\lbrace \Densrenormalized  - v^1,v^1,v^2 \rbrace$.
As in Lemma~\ref{L:TRANSPORTEQUATIONERRORINTEGRALSLOSEONEDERIVATIVE},
we allow the estimates to lose one derivative, and the gain is that
the right-hand sides of the estimates are much less singular with
respect to powers of $\upmu_{\star}^{-1}$ compared to the estimates we would
obtain in an approach that avoids derivative loss.

\begin{lemma}[\textbf{Estimates for wave equation error integrals involving a loss of one derivative}]
\label{L:WAVEEQUATIONERRORINTEGRALSLOSEONEDERIVATIVE}
	Let $\Psi \in \lbrace \Densrenormalized  - v^1,v^1,v^2 \rbrace$
	and assume that $2 \leq N \leq 20$.
	Recall that $\GeoAngFlatRadComponent$
	is the scalar-valued function appearing in Lemma~\ref{L:GEOANGDECOMPOSITION}.
	Under the data-size and bootstrap assumptions 
	of Sects.~\ref{SS:FLUIDVARIABLEDATAASSUMPTIONS}-\ref{SS:PSIBOOTSTRAP}
	and the smallness assumptions of Sect.~\ref{SS:SMALLNESSASSUMPTIONS},
	the following integral estimates hold for $(t,u) \in [0,\Tboot) \times [0,U_0]$
	(see Sect.~\ref{SS:STRINGSOFCOMMUTATIONVECTORFIELDS} regarding the vectorfield operator notation):
\begin{align}
	& 
		\int_{\mathcal{M}_{t,u}}
			\left|
				\myarray[\Rad \Fullset_{\ast}^{N-1;\leq 1} \Psi]
					{(1 + 2 \upmu) \Lunit \Fullset_{\ast}^{N-1;\leq 1} \Psi}
			\right|
			\left|
				\fourmyarray[
					(\Rad \Psi) \Fullset_{\ast}^{N-1; \leq 1} \mytr \upchi
					]
					{
					- (\upmu \angdiffuparg{\#} \Psi) \cdot (\upmu \angdiff \Fullset^{N-2;\leq 1} \mytr \upchi)
					}
					{
					\GeoAngFlatRadComponent (\angdiffuparg{\#} \Psi) \cdot (\upmu \angdiff \Fullset^{N-2;\leq 1} 
						\mytr \upchi)
					}
					{
					(\angdiffuparg{\#} \Psi) \cdot (\upmu \angdiff \Fullset^{N-2;\leq 1} \mytr \upchi)
					}
			\right|
		\, d \vol
			\label{E:WAVETRCHIDERIVATIVELOSINGERRORINTEGRALS} 
			\\
	& \lesssim
		\int_{t'=0}^t
			\frac{\sqrt{\totmax{[1,N-1]}}(t',u)}{\upmu_{\star}^{1/2}(t',u)}
			\left\lbrace
				\int_{s=0}^{t'}
					\frac{\sqrt{\totmax{N}}(s,u)}{\upmu_{\star}^{1/2}(s,u)}
				\, ds
			\right\rbrace
		\, dt'
		+
		\int_{t'=0}^t
			\frac{\totmax{[1,N-1]}(t',u)}{\upmu_{\star}^{1/2}(t',u)}
		\, dt'
		+
		\mathring{\upepsilon}^2.
		\notag
\end{align}


\end{lemma}

\begin{proof}
	It suffices to consider
	only the first term $(\Rad \Psi) \Fullset_{\ast}^{N-1; \leq 1} \mytr \upchi$
	in the second array on LHS~\eqref{E:WAVETRCHIDERIVATIVELOSINGERRORINTEGRALS}
	since the other three terms in the array can be bounded using the same arguments.
	They are in fact smaller in view of the estimates
	$\left\|
		\GeoAngFlatRadComponent
	\right\|_{L^{\infty}(\Sigma_t^u)}
	\lesssim \varepsilon
	$,
	$\left\|
		\angdiff \Psi
	\right\|_{L^{\infty}(\Sigma_t^u)}
	\lesssim \varepsilon
	$,
	and
	$\left\|
		\upmu
	\right\|_{L^{\infty}(\Sigma_t^u)}
	\lesssim 1
	$,
	which are simple consequences of \eqref{E:LINEARLYSMALLSCALARSDEPENDINGONGOODVARIABLES}
	and the $L^{\infty}$ estimates of Prop.~\ref{P:IMPROVEMENTOFAUX}.
	To proceed, we use Cauchy-Schwarz along $\Sigma_{t'}^u$,
	the $L^{\infty}$ estimates of Prop.~\ref{P:IMPROVEMENTOFAUX},
	Lemma~\ref{L:COERCIVENESSOFCONTROLLING},
	the estimate \eqref{E:TANGENGITALEIKONALINTERMSOFCONTROLLING},
	the simple estimate 
	$\mathring{\upepsilon} \sqrt{\totmax{[1,N-1]}}(t',u)
	\leq 
	\mathring{\upepsilon}^2 
	+
	\totmax{[1,N-1]}(t',u)
	$,
	and inequality \eqref{E:LESSSINGULARTERMSMPOINTNINEINTEGRALBOUND}
	to bound the spacetime integral under consideration as follows:
	\begin{align}
		&
		\int_{\mathcal{M}_{t,u}}
			\left|
				\myarray[\Rad \Fullset_{\ast}^{N-1;\leq 1} \Psi]
				{(1 + 2 \upmu) \Lunit \Fullset_{\ast}^{N-1;\leq 1} \Psi}
			\right|
			\left|
				(\Rad \Psi) \Fullset_{\ast}^{N-1; \leq 1} \mytr \upchi
			\right|
		\, d \vol
			\\
	& \lesssim
		\int_{t'=0}^t
			\left\lbrace
				\left\|
					\Rad \Fullset_{\ast}^{N-1;\leq 1} \Psi
				\right\|_{L^2(\Sigma_{t'}^u	)}
				+
				\left\|
					\Lunit \Fullset_{\ast}^{N-1;\leq 1} \Psi
				\right\|_{L^2(\Sigma_{t'}^u	)}
			\right\rbrace
			\left\|
					\Fullset_{\ast}^{N-1; \leq 1} \mytr \upchi
				\right\|_{L^2(\Sigma_{t'}^u	)}
		\, dt'
			\notag \\
& \lesssim
		\int_{t'=0}^t
			\frac{\sqrt{\totmax{[1,N-1]}}(t',u)}{\upmu_{\star}^{1/2}(t',u)}
			\left\lbrace
				\mathring{\upepsilon}
				+
				\int_{s=0}^{t'}
					\frac{\sqrt{\totmax{N}}(s,u)}{\upmu_{\star}^{1/2}(s,u)}
				\, ds
			\right\rbrace
		\, dt'
		\notag
			\\
	& \lesssim
		\int_{t'=0}^t
			\frac{\sqrt{\totmax{[1,N-1]}}(t',u)}{\upmu_{\star}^{1/2}(t',u)}
			\left\lbrace
				\int_{s=0}^{t'}
					\frac{\sqrt{\totmax{N}}(s,u)}{\upmu_{\star}^{1/2}(s,u)}
				\, ds
			\right\rbrace
			\, dt'
			+
		\int_{t'=0}^t
			\frac{\totmax{[1,N-1]}(t',u)}{\upmu_{\star}^{1/2}(t',u)}
		\, dt'
			\notag \\
	& \ \
		+
		\mathring{\upepsilon}^2
		\int_{t'=0}^t
			\frac{1}{\upmu_{\star}^{1/2}(t',u)}
		\, dt'
		\notag
			\\
	& \lesssim
		\int_{t'=0}^t
			\frac{\sqrt{\totmax{[1,N-1]}}(t',u)}{\upmu_{\star}^{1/2}(t',u)}
			\left\lbrace
				\int_{s=0}^{t'}
					\frac{\sqrt{\totmax{N}}(s,u)}{\upmu_{\star}^{1/2}(s,u)}
				\, ds
			\right\rbrace
			\, dt'
			+
		\int_{t'=0}^t
			\frac{\totmax{[1,N-1]}(t',u)}{\upmu_{\star}^{1/2}(t',u)}
		\, dt'
			+
			\mathring{\upepsilon}^2.
			\notag
	\end{align}

\end{proof}

\subsection{Proof of Prop.~\ref{P:VORTICITYENERGYINTEGRALINEQUALITIES}}
\label{SS:PROOFOFPROPVORTICITYENERGYINTEGRALINEQUALITIES}
We first prove \eqref{E:TOPORDERVORTICITYENERGYINTEGRALINEQUALITIES}.
Let $\vec{I}$ be a $\Tanset$ multi-index with $|\vec{I}|=21$.
From \eqref{E:ENERGYIDENTITYRENORMALIZEDVORTICITY},
we deduce that
\begin{align} \label{E:TANGENTCOMMUTEDENERGYIDENTITYRENORMALIZEDVORTICITY}
		\Vortenzero[\Tanset^{\vec{I}} \Vortrenormalized](t,u)
		+ 
		\Vortflzero[\Tanset^{\vec{I}} \Vortrenormalized](t,u)
		& 
		=
		\Vortenzero[\Tanset^{\vec{I}} \Vortrenormalized](0,u)
		+ 
		\Vortflzero[\Tanset^{\vec{I}} \Vortrenormalized](t,0)
			\\
		& \ \
			+
			\int_{\mathcal{M}_{t,u}}
				\left\lbrace
					\Lunit \upmu
					+
					\upmu \mytr \angk
				\right\rbrace
				(\Tanset^{\vec{I}} \Vortrenormalized)^2 
			\, d \vol
				\notag \\
	& \ \
			+
			2
			\int_{\mathcal{M}_{t,u}}
				(\Tanset^{\vec{I}} \Vortrenormalized) 
				\upmu \Transport \Tanset^{\vec{I}} \Vortrenormalized
			\, d \vol.
			\notag
\end{align}
We will show that the magnitude of RHS~\eqref{E:TANGENTCOMMUTEDENERGYIDENTITYRENORMALIZEDVORTICITY}
is $\leq \mbox{RHS~\eqref{E:TOPORDERVORTICITYENERGYINTEGRALINEQUALITIES}}$.
Then taking the max of that inequality over all $\vec{I}$ with $|\vec{I}|=21$
and appealing to Def.~\ref{D:MAINCOERCIVEQUANT}, we arrive at
\eqref{E:TOPORDERVORTICITYENERGYINTEGRALINEQUALITIES}.
The first integral on RHS~\eqref{E:TANGENTCOMMUTEDENERGYIDENTITYRENORMALIZEDVORTICITY} was
treated in Lemma~\ref{L:SIMPLESTTRANSPORTEQUATIONERROR}.
To bound the last integral on RHS~\eqref{E:TANGENTCOMMUTEDENERGYIDENTITYRENORMALIZEDVORTICITY},
we first use Prop.~\ref{P:VORTICITYIDOFKEYDIFFICULTENREGYERRORTERMS} to express 
the integrand factor
$
\upmu \Transport \Tanset^{\vec{I}} \Vortrenormalized
$
as the products explicitly indicated
on either RHS~\eqref{E:VORTICITYLISTHEFIRSTCOMMUTATORIMPORTANTTERMS}
or RHS~\eqref{E:VORTICITYGEOANGANGISTHEFIRSTCOMMUTATORIMPORTANTTERMS}
plus $Harmless_{(Vort)}^{\leq 21}$ error terms.
The error integrals
$
\int_{\mathcal{M}_{t,u}}
	(\Tanset^{\vec{I}} \Vortrenormalized) 
	Harmless_{(Vort)}^{\leq 21}
\, d \vol
$
were treated in Lemma~\ref{L:STANDARDPSISPACETIMEINTEGRALS}.
The remaining three error integrals, 
which correspond to the products explicitly indicated
on either RHS~\eqref{E:VORTICITYLISTHEFIRSTCOMMUTATORIMPORTANTTERMS}
and RHS~\eqref{E:VORTICITYGEOANGANGISTHEFIRSTCOMMUTATORIMPORTANTTERMS},
were treated in 
Lemma~\ref{L:MOSTDEGENERATETOPORDERTRANSPORTEQUATIONERRORINTEGRALS}.
We have thus proved \eqref{E:TOPORDERVORTICITYENERGYINTEGRALINEQUALITIES}.

The proof of \eqref{E:BELOWTOPORDERVORTICITYENERGYINTEGRALINEQUALITIES} 
in the cases $2 \leq N \leq 20$ is similar.
The only difference is that we bound the explicitly listed products 
on RHS~\eqref{E:VORTICITYLISTHEFIRSTCOMMUTATORIMPORTANTTERMS}
and RHS~\eqref{E:VORTICITYGEOANGANGISTHEFIRSTCOMMUTATORIMPORTANTTERMS}
(with $N-1$ in the role of $N$ in 
\eqref{E:VORTICITYLISTHEFIRSTCOMMUTATORIMPORTANTTERMS}-\eqref{E:VORTICITYGEOANGANGISTHEFIRSTCOMMUTATORIMPORTANTTERMS}) in a different way: by using the derivative-losing 
Lemma~\ref{L:TRANSPORTEQUATIONERRORINTEGRALSLOSEONEDERIVATIVE}
in place of Lemma~\ref{L:MOSTDEGENERATETOPORDERTRANSPORTEQUATIONERRORINTEGRALS}.
The proof of \eqref{E:BELOWTOPORDERVORTICITYENERGYINTEGRALINEQUALITIES} 
in the case $N=1$ is similar but simpler and
relies on equation
\eqref{E:TRIVIALVORTICITYHARMLESSORDERNPLUSONECOMMUTATORS}.
The proof of \eqref{E:BELOWTOPORDERVORTICITYENERGYINTEGRALINEQUALITIES} 
when $N=0$ is even simpler since,
by \eqref{E:RENORMALIZEDVORTICTITYTRANSPORTEQUATION},
the last integral on RHS~\eqref{E:TANGENTCOMMUTEDENERGYIDENTITYRENORMALIZEDVORTICITY}
completely vanishes.

\hfill $\qed$

\subsection{Proof of Prop.~\ref{P:WAVEENERGYINTEGRALINEQUALITIES}}
\label{SS:PROOFOFPROPWAVEENERGYINTEGRALINEQUALITIES}
\ \\

\noindent \textbf{Proof of \eqref{E:TOPORDERWAVEENERGYINTEGRALINEQUALITIES}}:
We set $N=20$ (which corresponds to the top-order number of commutations of the wave equations
\eqref{E:VELOCITYWAVEEQUATION}-\eqref{E:RENORMALIZEDDENSITYWAVEEQUATION}).
Let $\Fullset_{\ast}^{N;\leq 1}$ be an $N^{th}$-order 
vectorfield operator involving at most one $\Rad$ factor
and let $\Psi \in \lbrace \Densrenormalized  - v^1,v^1,v^2 \rbrace$.
From \eqref{E:E0DIVID} with $\Fullset_{\ast}^{N;\leq 1} \Psi$ in the role of $\Psi$, 
the decomposition \eqref{E:MULTERRORINT} with 
$\Fullset_{\ast}^{N;\leq 1} \Psi$ in the role of $\Psi$,
and definition \eqref{E:COERCIVESPACETIMEDEF},
we have
\begin{align} \label{E:E0DIVIDMAINESTIMATES}
	\enzero[\Fullset_{\ast}^{N;\leq 1} \Psi](t,u)
	& 
	+ 
	\flzero[\Fullset_{\ast}^{N;\leq 1} \Psi](t,u)
	+
	\coercivespacetime[\Fullset_{\ast}^{N;\leq 1} \Psi](t,u)
	\\
	& 
		=
		\enzero[\Fullset_{\ast}^{N;\leq 1} \Psi](0,u)
		+ 
		\flzero[\Fullset_{\ast}^{N;\leq 1} \Psi](t,0)
			\notag \\
	& \ \
		- 
		\int_{\mathcal{M}_{t,u}}
			\left\lbrace
				(1 + 2 \upmu) (\Lunit \Fullset_{\ast}^{N;\leq 1} \Psi)
				+ 
				2 \Rad \Fullset_{\ast}^{N;\leq 1} \Psi 
			\right\rbrace
			\upmu \square_g(\Fullset_{\ast}^{N;\leq 1} \Psi) 
		\, d \vol
							\notag \\
				& \ \ 
						+
						\sum_{i=1}^5 
						\int_{\mathcal{M}_{t,u}}
							\basicenergyerrorarg{\Mult}{i}[\Fullset_{\ast}^{N;\leq 1} \Psi]
						\, d \vol.
							\notag
\end{align}
We will show that 
$
\mbox{RHS~\eqref{E:E0DIVIDMAINESTIMATES}} 
\leq \mbox{RHS~\eqref{E:TOPORDERWAVEENERGYINTEGRALINEQUALITIES}}
$.
Then, taking the max over that estimate for all such operators of order precisely $N$ 
and over 
$\Psi \in \lbrace \Densrenormalized  - v^1,v^1,v^2 \rbrace$
and appealing to Defs.~\ref{D:MAINCOERCIVEQUANT} and \ref{D:COERCIVEINTEGRAL},
we conclude \eqref{E:TOPORDERWAVEENERGYINTEGRALINEQUALITIES}.

To show that 
$
\mbox{RHS~\eqref{E:E0DIVIDMAINESTIMATES}} 
\leq \mbox{RHS~\eqref{E:TOPORDERWAVEENERGYINTEGRALINEQUALITIES}}
$,
we first use Lemma~\ref{L:INITIALSIZEOFL2CONTROLLING} to deduce that
$
\enzero[\Fullset_{\ast}^{N;\leq 1} \Psi](0,u)
+ 
\flzero[\Fullset_{\ast}^{N;\leq 1} \Psi](t,0)
\lesssim \mathring{\upepsilon}^2
$,
which is $\leq$ the first term on 
RHS~\eqref{E:TOPORDERWAVEENERGYINTEGRALINEQUALITIES} 
as desired.

To bound the last integral
$
\sum_{i=1}^5
\int_{\mathcal{M}_{t,u}} 
	\basicenergyerrorarg{\Mult}{i}[\cdots]
$
on RHS~\eqref{E:E0DIVIDMAINESTIMATES} 
by 
$
\leq \mbox{RHS~\eqref{E:TOPORDERWAVEENERGYINTEGRALINEQUALITIES}}
$,
we use Lemma~\ref{L:MULTIPLIERVECTORFIELDERRORINTEGRALS}.

We now address the first integral 
$
- \int_{\mathcal{M}_{t,u}} \cdots
$
on RHS~\eqref{E:E0DIVIDMAINESTIMATES}.
If $\Fullset_{\ast}^{N;\leq 1}$ is \emph{not} of the form 
$\GeoAng^{N-1} \Lunit$, 
$\GeoAng^N$,
$\GeoAng^{N-1} \Rad$,
$\Fullset_*^{N-1;1} \Lunit$,
or
$\Fullset_*^{N-1;1} \GeoAng$,
where $\Fullset_*^{N-1;1}$ contains
exactly one factor of $\Rad$ and $N-2$ factors
of $\GeoAng$,
then 
it is easy to see that
$\Fullset_{\ast}^{N;\leq 1}$ 
must be of the form
of one of the operators on LHSs \eqref{E:WAVEHARMLESSORDERNPURETANGENTIALCOMMUTATORS}-\eqref{E:WAVETHIRDHARMLESSORDERNCOMMUTATORS}.
The desired bound thus follows from 
\eqref{E:WAVEHARMLESSORDERNPURETANGENTIALCOMMUTATORS}-\eqref{E:WAVETHIRDHARMLESSORDERNCOMMUTATORS},
\eqref{E:KEYDIFFICULTFORRENORMALZIEDDENSITY},
\eqref{E:WAVESHARMLESSSPACETIMEINTEGRALS},
and
\eqref{E:WAVEERRORINTEGRALSINVOLVINGTOPORDERVORTICITY}.
Note that these bounds do not produce any of 
the difficult ``boxed-constant-involving'' terms on 
RHS~\eqref{E:TOPORDERWAVEENERGYINTEGRALINEQUALITIES}.

We now address the first integral 
$
- \int_{\mathcal{M}_{t,u}} \cdots
$
on RHS~\eqref{E:E0DIVIDMAINESTIMATES}
when $\Psi = v^1$
and $\Fullset_*^{N;\leq 1}$ is one of the 
five operators not treated in the previous paragraph,
that is, when $\Fullset_*^{N;\leq 1}$ is one of
$\GeoAng^{N-1} \Lunit$, 
$\GeoAng^N$,
$\GeoAng^{N-1} \Rad$,
$\Fullset^{N-1;1} \Lunit$,
or
$\Fullset^{N-1;1} \GeoAng$,
where $\Fullset^{N-1;1} = \Fullset_*^{N-1;1}$ contains
exactly one factor of $\Rad$ and $N-2$ factors
of $\GeoAng$.
We consider in detail only the case 
$\Fullset_*^{N;\leq 1} = \GeoAng^N$;
the other four cases
can be treated in an identical fashion
(with the help of Prop.~\ref{P:WAVEIDOFKEYDIFFICULTENREGYERRORTERMS})
and we omit those details.
Moreover, the estimates for the wave variables
$\Psi = v^2$ 
and $\Psi = \Densrenormalized  - v^1$
are less degenerate and easier to derive; we will briefly comment on them below.
To proceed, we substitute RHS~\eqref{E:WAVEGEOANGISTHEFIRSTCOMMUTATORIMPORTANTTERMS}
(in the case $i=1$)
for the integrand factor $\upmu \square_g(\GeoAng^N v^1)$ 
on RHS~\eqref{E:E0DIVIDMAINESTIMATES}.
It suffices for us to bound the
integrals corresponding to the terms
$(\Rad v^1) \GeoAng^N \mytr \upchi$
and
$\GeoAngFlatRadComponent (\angdiffuparg{\#} v^1) \cdot (\upmu \angdiff \GeoAng^{N-1} \mytr \upchi)$
from RHS~\eqref{E:WAVEGEOANGISTHEFIRSTCOMMUTATORIMPORTANTTERMS};
the integrals generated by the $\Vortrenormalized$-involving terms
on RHS~\eqref{E:WAVEGEOANGISTHEFIRSTCOMMUTATORIMPORTANTTERMS} 
were suitably bounded in Lemma~\ref{L:WAVEERRORINTEGRALSINVOLVINGTOPORDERVORTICITY},
while the above argument has already addressed how to
bound the integrals generated by $Harmless_{(Wave)}^{\leq N}$
terms (via \eqref{E:WAVESHARMLESSSPACETIMEINTEGRALS}).
To bound the difficult integral
\begin{align} \label{E:DIFFICULTINTEGRALFORV1MAINPROOF}
	-           2
							\int_{\mathcal{M}_{t,u}}
							(\Rad \GeoAng^N v^1)	 
							(\Rad v^1) \GeoAng^N \mytr \upchi
							\, d \vol
\end{align}
in magnitude by $\leq \mbox{RHS~\eqref{E:TOPORDERWAVEENERGYINTEGRALINEQUALITIES}}$,
we use the estimate \eqref{E:MOSTDIFFICULTERRORINTEGRALBOUND},
which accounts for
the portion $\boxed{4} \cdots$ of the first boxed constant integral
$\boxed{6} \cdots$ on RHS~\eqref{E:TOPORDERWAVEENERGYINTEGRALINEQUALITIES}
and the full portion of 
the boxed constant integral $\boxed{8.1} \cdots$ 
on RHS~\eqref{E:TOPORDERWAVEENERGYINTEGRALINEQUALITIES}.

We now bound the magnitude of the error integral
\begin{align} \label{E:LUNITINVOLVINGDIFFICULTERRORINTEGRAL}
	-           \int_{\mathcal{M}_{t,u}}
							(1 + 2 \upmu)
							(\Lunit \GeoAng^N v^1)	 
							(\Rad v^1) \GeoAng^N \mytr \upchi
							\, d \vol.
\end{align}
To proceed, we use 
\eqref{E:TRANSPORTPARTIALRENORMALIZEDTRCHIJUNK}-\eqref{E:TRANSPORTPARTIALRENORMALIZEDTRCHIJUNKDISCREPANCY} 
to decompose
$ \GeoAng^N \mytr \upchi 
= \GeoAng \upchipartialmodarg{\GeoAng^{N-1}}
	- \GeoAng \upchipartialmodinhomarg{\GeoAng^{N-1}}$.
Since RHS~\eqref{E:HARMLESSNATUREOFPARTIALLYMODIFIEDDISCREPANCY} $= Harmless^{\leq N}$,
we have already suitably bounded the error integrals generated by
$\GeoAng \upchipartialmodinhomarg{\GeoAng^{N-1}}$.
We therefore must bound the magnitude of
\begin{align} \label{E:TOPORDERLUNITIBPDIFFICULTTERMS}
	-
	\int_{\mathcal{M}_{t,u}}
		(1 + 2 \upmu)
		(\Lunit \GeoAng^N v^1)	 
		(\Rad v^1) 
		\GeoAng \upchipartialmodarg{\GeoAng^{N-1}}
	\, d \vol
\end{align}
by $\leq \mbox{RHS~\eqref{E:TOPORDERWAVEENERGYINTEGRALINEQUALITIES}}$.
To this end, we integrate by parts using \eqref{E:LUNITANDANGULARIBPIDENTITY}
with 
$\ThirdSmoothFunction := \upchipartialmodarg{\GeoAng^{N-1}}$.
We bound the error integrals on the last line of RHS \eqref{E:LUNITANDANGULARIBPIDENTITY}
and the $\int_{\Sigma_0^u} \cdots$ integral on the second line
using Lemma~\ref{L:HARMLESSIBPERRORINTEGRALS}.
It remains for us to bound the first two (difficult) integrals on
RHS \eqref{E:LUNITANDANGULARIBPIDENTITY} 
in magnitude by $\leq \mbox{RHS~\eqref{E:TOPORDERWAVEENERGYINTEGRALINEQUALITIES}}$.
The desired bounds have been derived in the estimates
\eqref{E:DIFFICULTLUNITSPACETIMEIBPINTEGRALBOUND}-\eqref{E:DIFFICULTLUNITHYPERSURFACEIBPINTEGRALBOUND}
of Lemma~\ref{L:BOUNDSFORDIFFICULTTOPORDERINTEGRALSINVOLVINGLUNITIBP}.
Note that these estimates account for
the remaining portion $\boxed{2} \cdots$ of the first boxed constant integral
$\boxed{6} \cdots$ on RHS~\eqref{E:TOPORDERWAVEENERGYINTEGRALINEQUALITIES}
and the full portion of 
the boxed constant integral $\boxed{2} \cdots$ 
on RHS~\eqref{E:TOPORDERWAVEENERGYINTEGRALINEQUALITIES}.

To finish deriving the desired estimates in the case $\Psi = v^1$,
it remains for us to bound the two error integrals generated by the term
$\GeoAngFlatRadComponent (\angdiffuparg{\#} v^1) \cdot (\upmu \angdiff \GeoAng^{N-1} \mytr \upchi)$
from RHS~\eqref{E:WAVEGEOANGISTHEFIRSTCOMMUTATORIMPORTANTTERMS}.
These two integrals were suitably bounded in magnitude
by $\leq \mbox{RHS~\eqref{E:TOPORDERWAVEENERGYINTEGRALINEQUALITIES}}$
in Lemma~\ref{L:LESSDEGENERATEENERGYESTIMATEINTEGRALS}
(note that we are using the simple bound
$\left\lbrace 
					\ln \upmu_{\star}^{-1}(t',u) 
					+ 
					1 
				\right\rbrace^2
\lesssim 
\upmu_{\star}^{-1/2}(t',u)
$
in order to bound the 
integrand factors in the first integrals
on RHS \eqref{E:FIRSTLESSDEGENERATEENERGYESTIMATEINTEGRALS}
and RHS \eqref{E:SECONDLESSDEGENERATEENERGYESTIMATEINTEGRALS}).
Note also
that these estimates do not contribute to the difficult ``boxed-constant-involving'' 
products on RHS~\eqref{E:TOPORDERWAVEENERGYINTEGRALINEQUALITIES}.
We have thus shown that when $\Psi = v^1$,
the desired inequality 
$\mbox{RHS\eqref{E:E0DIVIDMAINESTIMATES}}
\leq 
\mbox{RHS~\eqref{E:TOPORDERWAVEENERGYINTEGRALINEQUALITIES}}
$
holds.

We now comment on the cases
$\Psi = v^2$ 
and $\Psi =  \Densrenormalized  - v^1$.
The proofs that
$
\mbox{RHS~\eqref{E:E0DIVIDMAINESTIMATES}} 
\leq \mbox{RHS~\eqref{E:TOPORDERWAVEENERGYINTEGRALINEQUALITIES}}
$
in these cases
are essentially the same as in the case $\Psi = v^1$,
except that in bounding the analog of the error integral \eqref{E:DIFFICULTINTEGRALFORV1MAINPROOF},
we now use the less degenerate estimate
\eqref{E:LESSDEGENERATEDIFFICULTERRORINTEGAL}
in place of
 \eqref{E:MOSTDIFFICULTERRORINTEGRALBOUND}
and, in bounding the analog of the error integral \eqref{E:TOPORDERLUNITIBPDIFFICULTTERMS},
we use the less degenerate estimates
\eqref{E:LESSDEGENERATEDIFFICULTLUNITSPACETIMEIBPINTEGRALBOUND}-\eqref{E:LESSDEGENERATEDIFFICULTLUNITHYPERSURFACEIBPINTEGRALBOUND}
in place of
\eqref{E:DIFFICULTLUNITSPACETIMEIBPINTEGRALBOUND}-\eqref{E:DIFFICULTLUNITHYPERSURFACEIBPINTEGRALBOUND}.
These less degenerate estimates do not produce any of the
``boxed-constant-involving'' products
on RHS~\eqref{E:TOPORDERWAVEENERGYINTEGRALINEQUALITIES}
because they all gain a smallness factor of $\varepsilon$
via the factors $\Rad v^2$ and $\Rad (\Densrenormalized  - v^1))$
(which verify the smallness estimates 
\eqref{E:PSITRANSVERSALLINFINITYBOUNDBOOTSTRAPIMPROVEDSMALL}
and
\eqref{E:CRUCIALPSITRANSVERSALLINFINITYBOUNDBOOTSTRAPIMPROVEDSMALL}).
In total, we have proved \eqref{E:TOPORDERWAVEENERGYINTEGRALINEQUALITIES}.

\medskip
\noindent \textbf{Proof of \eqref{E:EASYVARIABLESTOPORDERWAVEENERGYINTEGRALINEQUALITIES}}:
The argument given in the previous paragraph yields 
\eqref{E:EASYVARIABLESTOPORDERWAVEENERGYINTEGRALINEQUALITIES}.

\medskip
\noindent \textbf{Proof of \eqref{E:BELOWTOPORDERWAVEENERGYINTEGRALINEQUALITIES}}:
We repeat the proof of \eqref{E:TOPORDERWAVEENERGYINTEGRALINEQUALITIES}
with $M$ in the role of $N$, where $1 \leq M \leq N-1$,
and with one important change:
we bound the difficult error integrals
such as
\[
	-           2
							\int_{\mathcal{M}_{t,u}}
							(\Rad \GeoAng^N \Psi)	 
							(\Rad \Psi) \GeoAng^N \mytr \upchi
							\, d \vol,
							\qquad
								-           \int_{\mathcal{M}_{t,u}}
							(1 + 2 \upmu)
							(\Lunit \GeoAng^N \Psi)	 
							(\Rad \Psi) \GeoAng^N \mytr \upchi
							\, d \vol
\]
in a different way: 
by using Lemma~\ref{L:WAVEEQUATIONERRORINTEGRALSLOSEONEDERIVATIVE}.
More precisely, we replace $N$ with $M$ in 
\eqref{E:WAVELISTHEFIRSTCOMMUTATORIMPORTANTTERMS}-\eqref{E:SECONDWAVEGEOANGANGISTHEFIRSTCOMMUTATORIMPORTANTTERMS}
and consider the explicitly listed products on the RHSs that involve
the derivatives of $\mytr \upchi$ 
(see also \eqref{E:KEYDIFFICULTFORRENORMALZIEDDENSITY} in the case 
$\Psi = \Densrenormalized  - v^1$).
We bound the corresponding
error integrals by using the derivative-losing Lemma~\ref{L:WAVEEQUATIONERRORINTEGRALSLOSEONEDERIVATIVE}
in place of the arguments used in proving \eqref{E:TOPORDERWAVEENERGYINTEGRALINEQUALITIES}.

\hfill $\qed$

\subsection{The main vorticity a priori energy estimates}
\label{SS:MAINVORTICITYAPRIORIENERGYESTIMATES}
The energy estimates for the specific vorticity are easy to derive
with the help of the bootstrap assumptions. We provide
them in the next lemma.

\begin{lemma}[\textbf{The main a priori energy estimates for the specific vorticity}]
	\label{L:VORTICITYAPRIORIENERGYESTIMATES}
		Under the data-size and bootstrap assumptions 
		of Sects.~\ref{SS:FLUIDVARIABLEDATAASSUMPTIONS}-\ref{SS:PSIBOOTSTRAP}
		and Sect.~\ref{SS:L2BOOTSTRAP}
		and the smallness assumptions of Sect.~\ref{SS:SMALLNESSASSUMPTIONS}, 
		the a priori energy estimates
		\eqref{E:TOPVORTMULOSSMAINAPRIORIENERGYESTIMATES}-\eqref{E:VORTNOMULOSSMAINAPRIORIENERGYESTIMATES}
		for the vorticity hold
		on $\mathcal{M}_{\Tboot,U_0}$.
\end{lemma}

\begin{proof}
	We start by deriving the desired estimates 
	\eqref{E:TOPVORTMULOSSMAINAPRIORIENERGYESTIMATES}-\eqref{E:VORTNOMULOSSMAINAPRIORIENERGYESTIMATES}
	for $\VorttotTanmax{21}(t,u)$ and $\VorttotTanmax{\leq 20}(t,u)$.
	Below we will use inequalities
	\eqref{E:TOPORDERVORTICITYENERGYINTEGRALINEQUALITIES}-\eqref{E:BELOWTOPORDERVORTICITYENERGYINTEGRALINEQUALITIES}
	and the bootstrap assumptions 
	\eqref{E:BAMULOSSMAINAPRIORIENERGYESTIMATES}-\eqref{E:BAVORTNOMULOSSMAINAPRIORIENERGYESTIMATES}
	to obtain the following inequalities:
	\begin{align}
		\VorttotTanmax{21}(t,u)
		& \leq
		C 
		\left(
			\mathring{\upepsilon}^2 
			+
			\varepsilon^3 
		\right)
		\upmu_{\star}^{-12.8}(t,u)
		+
		C
		\int_{u'=0}^u
			\VorttotTanmax{21}(t,u')
		\, du'
		+
		C
		\int_{u'=0}^u
			\VorttotTanmax{\leq 20}(t,u')
		\, du',
		 	\label{E:TOPORDERVORTICITYGRONWALLREADY} \\
		\VorttotTanmax{\leq 20}(t,u)
		& \leq
		C 
		\left(
			\mathring{\upepsilon}^2 
			+
			\varepsilon^3 
		\right)
		\upmu_{\star}^{-9.8}(t,u)
		+
		C
		\int_{u'=0}^u
			\VorttotTanmax{\leq 20}(t,u')
		\, du'.
		\label{E:JUSTBELOWTOPORDERVORTICITYGRONWALLREADY}
	\end{align}
	Then from \eqref{E:JUSTBELOWTOPORDERVORTICITYGRONWALLREADY}
	and Gronwall's inequality in $u$, we obtain 
	$
	\displaystyle
	\VorttotTanmax{\leq 20}(t,u) 
	\leq C  
	\left(
			\mathring{\upepsilon}^2 
			+
			\varepsilon^3 
		\right)
		\upmu_{\star}^{-9.8}
	$.
	Inserting this estimate into the last integral on RHS~\eqref{E:TOPORDERVORTICITYGRONWALLREADY},
	we find that $\VorttotTanmax{21}(t,u)$ obeys inequality \eqref{E:TOPORDERVORTICITYGRONWALLREADY}
	but with the last integral deleted. 
	Hence, from Gronwall's inequality in $u$, 
	we obtain 
	$
	\displaystyle
	\VorttotTanmax{21}(t,u) 
	\leq C 
	\left(
			\mathring{\upepsilon}^2 
			+
			\varepsilon^3 
		\right)
	\upmu_{\star}^{-12.8}(t,u)
	$.
	Recalling the assumption $\varepsilon^{3/2} \leq \mathring{\upepsilon}$
	(see \eqref{E:DATAEPSILONVSBOOTSTRAPEPSILON}),
	we see that we have shown \eqref{E:TOPVORTMULOSSMAINAPRIORIENERGYESTIMATES}
	and the estimate \eqref{E:VORTMULOSSMAINAPRIORIENERGYESTIMATES}  
	for $\sqrt{\VorttotTanmax{20}}(t,u)$.

	It remains for us 
	to derive \eqref{E:TOPORDERVORTICITYGRONWALLREADY}-\eqref{E:JUSTBELOWTOPORDERVORTICITYGRONWALLREADY}.
	To derive \eqref{E:TOPORDERVORTICITYGRONWALLREADY}
	we set $N=20$ in \eqref{E:TOPORDERVORTICITYENERGYINTEGRALINEQUALITIES},
	which yields an integral inequality for $\VorttotTanmax{21}(t,u)$.
	We then insert the bootstrap assumptions
	\eqref{E:BAMULOSSMAINAPRIORIENERGYESTIMATES}-\eqref{E:BAVORTNOMULOSSMAINAPRIORIENERGYESTIMATES}
	into all terms on
	RHS~\eqref{E:TOPORDERVORTICITYENERGYINTEGRALINEQUALITIES}
	except for the last integral 
	$
	\displaystyle
	C
		\int_{u'=0}^u
			\VorttotTanmax{\leq N+1}(t,u')
		\, du'
	$.
	It immediately follows 
	that all of the terms generated by the bootstrap assumptions, 
	except for the ones involving time integrals, 
	are $\leq$ the  
	$
	\displaystyle
	C 
		\left(
			\mathring{\upepsilon}^2 
			+
			\varepsilon^3 
		\right)
	\upmu_{\star}^{-12.8}(t,u)
	$ 
	term on RHS~\eqref{E:TOPORDERVORTICITYGRONWALLREADY} as desired.
	We now explain how to handle the terms generated by the 
	time integrals on RHS~\eqref{E:TOPORDERVORTICITYENERGYINTEGRALINEQUALITIES}.
	We consider in detail only the term 
	$
	\displaystyle
	C
	\varepsilon^2
	\int_{t'=0}^t
				\frac{1}{\upmu_{\star}^2(t',u)} 
				\left\lbrace
					\int_{s=0}^t
						\frac{1}{\upmu_{\star}(s,u)} 
						\sqrt{\totmax{20}}(s,u)
					\, ds
				\right\rbrace^2
			\, dt'
	$;
	the remaining time integrals on RHS~\eqref{E:TOPORDERVORTICITYENERGYINTEGRALINEQUALITIES}
	can be bounded in a similar fashion and we omit the details.
	To proceed, we use the bootstrap assumptions,
	the estimate \eqref{E:LOSSKEYMUINTEGRALBOUND},
	and the assumption \eqref{E:DATAEPSILONVSBOOTSTRAPEPSILON}
	to deduce that the double time integral under consideration is
	\begin{align}
	&
	\leq
	C
	\varepsilon^3
	\int_{t'=0}^t
				\frac{1}{\upmu_{\star}^2(t',u)} 
				\left\lbrace
					\int_{s=0}^t
						\frac{1}{\upmu_{\star}^{6.9}(s,u)} 
					\, ds
				\right\rbrace^2
			\, dt'
				\\
	& \leq 
	C
	\varepsilon^3
	\int_{t'=0}^t
		\frac{1}{\upmu_{\star}^{13.8}(t',u)} 
	\, dt'
	\leq 
	C \varepsilon^3 \upmu_{\star}^{-12.8}(t,u)
	\leq C \mathring{\upepsilon}^2 \upmu_{\star}^{-12.8}(t,u)
	\notag
	\end{align}
	as desired.
	We have thus proved \eqref{E:TOPORDERVORTICITYGRONWALLREADY}.
	The proof of \eqref{E:JUSTBELOWTOPORDERVORTICITYGRONWALLREADY} is 
	based on inequality \eqref{E:BELOWTOPORDERVORTICITYENERGYINTEGRALINEQUALITIES} 
	with $N=20$
	but is otherwise similar to the proof of \eqref{E:TOPORDERVORTICITYGRONWALLREADY}; 
	we omit the details.
	We have thus obtained the desired estimates for
	$\VorttotTanmax{21}(t,u)$ and $\VorttotTanmax{20}(t,u)$.

	We now explain how to derive the estimates
	\eqref{E:VORTMULOSSMAINAPRIORIENERGYESTIMATES}-\eqref{E:VORTNOMULOSSMAINAPRIORIENERGYESTIMATES}
	for 
	$\sqrt{\VorttotTanmax{\leq 19}}$,
	$\sqrt{\VorttotTanmax{\leq 18}}$,
	$\cdots$,
	$\sqrt{\VorttotTanmax{0}}$.
	The desired estimates can be derived 
	from inequality \eqref{E:BELOWTOPORDERVORTICITYENERGYINTEGRALINEQUALITIES}, 
	the bootstrap assumptions
	\eqref{E:BANOMULOSSMAINAPRIORIENERGYESTIMATES}-\eqref{E:BAVORTNOMULOSSMAINAPRIORIENERGYESTIMATES},
	Gronwall's inequality in $u$,
	and the assumption \eqref{E:DATAEPSILONVSBOOTSTRAPEPSILON}
	by using essentially the same arguments
	that we used to derive 
	the estimates for $\sqrt{\VorttotTanmax{20}}$.
	However, there is one minor new feature that is needed to
	obtain the estimates \eqref{E:VORTNOMULOSSMAINAPRIORIENERGYESTIMATES}
	for $\sqrt{\VorttotTanmax{\leq 15}}(t,u)$:
	in carrying out the above procedure, 
	we encounter a term that needs to be treated
	using a slightly different argument:
	the term
	\begin{align} \label{E:TIMEINTEGRALTHATBREAKSVORTICITYMUDEGENERACY}
	C \varepsilon^3
		\int_{t'=0}^t
			\left\lbrace
				\int_{s=0}^{t'}
					\frac{1}{\upmu_{\star}^{1.4}(s,u)}
				\, ds
			\right\rbrace^2
		\, dt'
	\end{align}
	generated by the double time integral on RHS~\eqref{E:BELOWTOPORDERVORTICITYENERGYINTEGRALINEQUALITIES}.
	The new part of the argument is that in
	addition to inequality \eqref{E:LOSSKEYMUINTEGRALBOUND},
	we must also use inequality
	\eqref{E:LESSSINGULARTERMSMPOINTNINEINTEGRALBOUND};
	\emph{inequality \eqref{E:LESSSINGULARTERMSMPOINTNINEINTEGRALBOUND}
	is what allows us to break the $\upmu_{\star}^{-1}$ degeneracy}.
	More precisely, to bound
	the term \eqref{E:TIMEINTEGRALTHATBREAKSVORTICITYMUDEGENERACY},
	we use
	inequalities
	\eqref{E:LOSSKEYMUINTEGRALBOUND}
	and
	\eqref{E:LESSSINGULARTERMSMPOINTNINEINTEGRALBOUND}
	and the assumption \eqref{E:DATAEPSILONVSBOOTSTRAPEPSILON}
	to deduce that
	\begin{align} \label{E:VORTICITYMUDEGENERACYBROKEN}
		\mbox{\eqref{E:TIMEINTEGRALTHATBREAKSVORTICITYMUDEGENERACY}}
		\leq
		C \varepsilon^3
		\int_{t'=0}^t
			\frac{1}{\upmu_{\star}^{.8}(s,u)}
		\, dt'
	\leq
		C \varepsilon^3
	\leq C \mathring{\upepsilon}^2
	\end{align}
	as desired, where RHS~\eqref{E:VORTICITYMUDEGENERACYBROKEN}
	\emph{does not involve the singular factor $\upmu_{\star}^{-1}$}!
	We have thus 
	obtained the desired estimates
	\eqref{E:TOPVORTMULOSSMAINAPRIORIENERGYESTIMATES}-\eqref{E:BAVORTNOMULOSSMAINAPRIORIENERGYESTIMATES},
	which completes the proof of the lemma.
\end{proof}

\subsection{Proof of Prop.~\ref{P:MAINAPRIORIENERGY}}
\label{SS:PROOFOFPROPMAINAPRIORIENERGY}
To simplify the proof, we assume that the energy bootstrap assumptions
\eqref{E:BAMULOSSMAINAPRIORIENERGYESTIMATES}-\eqref{E:BAVORTNOMULOSSMAINAPRIORIENERGYESTIMATES}
hold for $(t,u) \in [0,\Tboot) \times [0,U_0]$.
To prove the proposition, it 
suffices to derive, under the energy bootstrap assumptions,
the estimates \eqref{E:MULOSSMAINAPRIORIENERGYESTIMATES}-\eqref{E:VORTNOMULOSSMAINAPRIORIENERGYESTIMATES}
for $(t,u) \in [0,\Tboot) \times [0,U_0]$.
We can then use a standard continuity-in-$t$ argument 
for the fundamental $L^2$-controlling quantities
to deduce that the estimates \eqref{E:MULOSSMAINAPRIORIENERGYESTIMATES}-\eqref{E:VORTNOMULOSSMAINAPRIORIENERGYESTIMATES}
do in fact hold for $(t,u) \in [0,\Tboot) \times [0,U_0]$
and, in view of our assumption
$\mathring{\upepsilon} \leq \varepsilon$, 
that the bootstrap assumptions are never saturated
(for $\mathring{\upepsilon}$ sufficiently small).
Note that this argument relies on
Lemma~\ref{L:INITIALSIZEOFL2CONTROLLING},
which implies that the fundamental $L^2$-controlling quantities
do not saturate inequalities
\eqref{E:BAMULOSSMAINAPRIORIENERGYESTIMATES}-\eqref{E:BAVORTNOMULOSSMAINAPRIORIENERGYESTIMATES}
at the initial time $0$.

We now recall that in Lemma~\ref{L:VORTICITYAPRIORIENERGYESTIMATES}, 
we derived,
with the help of the energy bootstrap assumptions,
the a priori vorticity energy estimates 
\eqref{E:TOPVORTMULOSSMAINAPRIORIENERGYESTIMATES}-\eqref{E:VORTNOMULOSSMAINAPRIORIENERGYESTIMATES}.
Hence, it remains only for us to derive
the wave variable energy estimates
\eqref{E:MULOSSMAINAPRIORIENERGYESTIMATES}-\eqref{E:NOMULOSSMAINAPRIORIENERGYESTIMATES}.
We are of course free to use
the vorticity energy estimates
\eqref{E:TOPVORTMULOSSMAINAPRIORIENERGYESTIMATES}-\eqref{E:VORTNOMULOSSMAINAPRIORIENERGYESTIMATES}
in the remainder of the proof.

\medskip

\noindent \textbf{Estimates for} 
	$\totmax{20}$, 
	$\coercivespacetimemax{20}$, 
	$\totmax{[1,19]}$, 
	\textbf{and} 
	$\coercivespacetimemax{[1,19]}$: 
	These estimates are highly coupled and must be treated as a system
	featuring also
	$\easytotmax{20}$
	and
	$\easycoercivespacetimemax{20}$. 
	To proceed, we set
	\begin{align}
		F(t,u) 
		& := \sup_{(\hat{t},\hat{u}) \in [0,t] \times [0,u]} 
			\iota_F^{-1}(\hat{t},\hat{u})
			\max
			\left\lbrace
				\totmax{20}(\hat{t},\hat{u}),
				\coercivespacetimemax{20}(\hat{t},\hat{u})
			\right\rbrace, 
			\label{E:TOPORDERENERGYRESCALED} \\
		G(t,u) 
		& := \sup_{(\hat{t},\hat{u}) \in [0,t] \times [0,u]} 
			\iota_G^{-1}(\hat{t},\hat{u})
			\max
			\left\lbrace
				\easytotmax{20}(\hat{t},\hat{u}),
				\easycoercivespacetimemax{20}(\hat{t},\hat{u})
			\right\rbrace, 
			\label{E:EASYTOPORDERENERGYRESCALED} \\
		H(t,u) 
		& := \sup_{(\hat{t},\hat{u}) \in [0,t] \times [0,u]} 
		\iota_H^{-1}(\hat{t},\hat{u})
		\max
		\left\lbrace
			\totmax{[1,19]}(\hat{t},\hat{u}),
			\coercivespacetimemax{[1,19]}(\hat{t},\hat{u})
		\right\rbrace,
		\label{E:JUSTBELOWTOPORDERENERGYRESCALED}
	\end{align}
	where for $1 \leq t' \leq \hat{t} \leq t$ 
	and $0 \leq u' \leq \hat{u} \leq U_0 \leq 1$,
	we define
	\begin{align} \label{E:INTEGRATINGFACTOR}
		\iota_1(t')
		&:=
		\exp
		\left(
			\int_{s=0}^{t'}
				\frac{1}{\sqrt{\Tboot - s}}
			\, ds
		\right)
		= \exp
			\left(
				2 \sqrt{\Tboot}
				- 2 \sqrt{\Tboot - t'}
			\right),
				\\
	\iota_2(t',u')
		&:=
		\exp
		\left(
			\int_{s=0}^{t'}
				\frac{1}{\upmu_{\star}^{9/10}(s,u')}
			\, ds
		\right),
			\label{E:SECONDINTEGRATINGFACTOR} \\
	\iota_F(t',u')
	= \iota_G(t',u') 
	& :=
		\upmu_{\star}^{-11.8}(t',u')
		\iota_1^c(t')
		\iota_2^c(t',u')
		e^{c t'}
		e^{c u'},
			\\
	\iota_H(t',u')
	& :=
		\upmu_{\star}^{-9.8}(t',u')
		\iota_1^c(t')
		\iota_2^c(t',u')
		e^{c t'}
		e^{c u'},
		\label{E:JUSTBELOWTOPORDERTOTALINTEGRATINGFACTOR}
	\end{align}
	and $c$ is a sufficiently large positive constant that we choose below.
	The functions \eqref{E:INTEGRATINGFACTOR}-\eqref{E:JUSTBELOWTOPORDERTOTALINTEGRATINGFACTOR}
	are approximate integrating factors that
	we will use to absorb all error terms 
	on the RHSs of the inequalities of Prop.~\ref{P:WAVEENERGYINTEGRALINEQUALITIES}
	back into the LHSs.
	We claim that to obtain the desired estimates for 
	$\totmax{20}$, 
	$\coercivespacetimemax{20}$,
	$\easytotmax{20}$,
	$\easycoercivespacetimemax{20}$,
	$\totmax{[1,19]}$,
	and 
	$\coercivespacetimemax{[1,19]}$,
	it suffices to prove 
	\begin{align} \label{E:DESIREDRESCALEDAPRIORIBOUND}
		F(t,u) 
		& \leq C \mathring{\upepsilon}^2,
		\qquad
		G(t,u) 
		\leq C \mathring{\upepsilon}^2,
		\qquad
		H(t,u) 
		\leq C \mathring{\upepsilon}^2,
	\end{align}
	where $C$ in \eqref{E:DESIREDRESCALEDAPRIORIBOUND} is allowed to depend on $c$.
	To justify the claim, we use the fact that for a fixed $c$,
	the functions
	$\iota_1^c(t)$,
	$\iota_2^c(t,u)$,
	$e^{c t}$,
	and
	$e^{c u}$
	are uniformly bounded from above by a positive constant
	for $(t,u) \in [0,\Tboot) \times [0,U_0]$;
	all of these estimates are simple to derive,
	except for \eqref{E:SECONDINTEGRATINGFACTOR},
	which relies on \eqref{E:LESSSINGULARTERMSMPOINTNINEINTEGRALBOUND}.

	To prove \eqref{E:DESIREDRESCALEDAPRIORIBOUND},
	it suffices to
	show that
	there exist positive constants 
	$\upalpha_1$,
	$\upalpha_2$, 
	$\upalpha_3$,
	$\upbeta_1$,
	$\upbeta_3$,
	$\upgamma_1$,
	and 
	$\upgamma_3$
	with
	\begin{align} \label{E:CONSTSBOUND}
		&
		\upalpha_1 
		+ 
		\upalpha_2 
		\sqrt{
					\upbeta_1
					+
					\upbeta_3 \frac{\upgamma_1}{1 - \upgamma_3}
				}
		+ 
		\frac{\upalpha_3 \upgamma_1}{1 - \gamma_3} 
		< 1,
		\qquad
		\upgamma_3 < 1
	\end{align}
	such that if $c$ is sufficiently large, then
	\begin{align} \label{E:TOPORDERALMOSTESTIMATED}
		F(t,u) 
		& \leq C \mathring{\upepsilon}^2 
			+ \upalpha_1 F(t,u) 
			+ \upalpha_2 F^{1/2}(t,u) G^{1/2}(t,u) 
			+ \upalpha_3 H(t,u),
			\\
		G(t,u) 
		& \leq C \mathring{\upepsilon}^2 
			+ \upbeta_1 F(t,u) 
			+ \upbeta_3 H(t,u),\\			
		H(t,u) \label{E:BELOWTOPORDERALMOSTESTIMATED}
		& \leq C \mathring{\upepsilon}^2 
			+ \upgamma_1 F(t,u) 
		  + \upgamma_3 H(t,u).
	\end{align}
One key reason that we will be able to obtain \eqref{E:CONSTSBOUND}
is that we will be able to make
  $\upalpha_3$,
	$\upbeta_1$,
	and
	$\upbeta_3$
	as small as we want 
 by choosing $\varsigma$ and $\varepsilon$ to be sufficiently small.
 Once we have obtained \eqref{E:TOPORDERALMOSTESTIMATED}-\eqref{E:BELOWTOPORDERALMOSTESTIMATED},
 we easily deduce from those estimates that
	\begin{align} \label{E:MOREBELOWTOPORDERALMOSTESTIMATED}
		H(t,u)
		& \leq 
			C \mathring{\upepsilon}^2 
			+ 
			\frac{\upgamma_1}{1 - \upgamma_3} 
			F(t,u),
			\\
		G(t,u)
		& \leq 
			C \mathring{\upepsilon}^2 
			+ 
			\left\lbrace
				\upbeta_1
				+ 
				\upbeta_3 \frac{\upgamma_1}{1 - \upgamma_3}
			\right\rbrace
			F(t,u),
			 \label{E:EASYTOPORDERALMOSTESTIMATED} \\
		F(t,u)
		& \leq
			C \mathring{\upepsilon}^2
			+
			C \mathring{\upepsilon} F^{1/2}(t,u)
			+ 
			\left\lbrace
				\upalpha_1
				+
				\upalpha_2 
				\sqrt{
					\upbeta_1
					+
					\upbeta_3 \frac{\upgamma_1}{1 - \upgamma_3}
				}
				+ 
				\frac{\upalpha_3 \upgamma_1}{1 - \upgamma_3} 
			\right\rbrace
			F(t,u).
		\label{E:MORETOPORDERALMOSTESTIMATED}
	\end{align}
	The desired bounds \eqref{E:DESIREDRESCALEDAPRIORIBOUND}
	(for $\mathring{\upepsilon}$ sufficiently small)
	now follow easily from
	\eqref{E:CONSTSBOUND}
	and
	\eqref{E:MOREBELOWTOPORDERALMOSTESTIMATED}-\eqref{E:MORETOPORDERALMOSTESTIMATED}.

	It remains for us to derive \eqref{E:TOPORDERALMOSTESTIMATED}-\eqref{E:BELOWTOPORDERALMOSTESTIMATED}.
	To this end, we will use the critically important estimates of Prop.~\ref{P:MUINVERSEINTEGRALESTIMATES}
	as well as the following estimates, which are easy to derive:
	\begin{align}
		\int_{t'=0}^{\hat{t}}
				\iota_1^c(t') \frac{1}{\sqrt{\Tboot - t'}} 
		\, dt'
		=
		\frac{1}{c}
		\int_{t'=0}^{\hat{t}}
				\frac{d}{dt'}
				\left\lbrace
					\iota_1^c(t') 
				\right\rbrace
		\, dt'
		& \leq \frac{1}{c} \iota_1^c(\hat{t}),
			\label{E:IF1SMALLNESS} \\
		\int_{t'=0}^{\hat{t}}
				\iota_2^c(t',\hat{u}) \frac{1}{\upmu_{\star}^{9/10}(t',\hat{u})} 
		\, dt'
		=
		\frac{1}{c}
		\int_{t'=0}^{\hat{t}}
				\frac{d}{dt'}
				\left\lbrace
					\iota_2^c(t',\hat{u}) 
				\right\rbrace
		\, dt'
		& \leq \frac{1}{c} \iota_2^c(\hat{t},\hat{u}),
			\label{E:IF2SMALLNESS} \\
		\int_{t'=0}^{\hat{t}}
				e^{c t'}
		\, dt'
		& \leq \frac{1}{c} e^{c \hat{t}},
			\label{E:IF3SMALLNESS} \\
		\int_{u'=0}^{\hat{u}}
				e^{c u'} 
		\, du'
		& \leq \frac{1}{c} e^{c \hat{u}}.
		\label{E:IF4SMALLNESS}
		\end{align}
	We will close the estimates 
	by taking $c$ to be large and $\varepsilon$ to be small.

	\begin{remark}
	\label{R:CONSTANTSINDEPENDENTOFLITTLEC}
	\emph{We stress that from now through
	inequality \eqref{E:NEARLYBELOWTOPORDERALMOSTESTIMATED},
	the constants $C$ can be chosen to be independent of $c$.}
	\end{remark}

	In our analysis, we will often use the fact that
	$\iota_1^c(\cdot)$,
	$\iota_2^c(\cdot)$
	$e^{c \cdot}$,
	and $e^{c \cdot}$ are non-decreasing in their arguments.
	Also, we will often use the estimate \eqref{E:MUSTARINVERSEMUSTGROWUPTOACONSTANT}, 
	which implies that for
	$t' \leq \hat{t}$
	and $u' \leq \hat{u}$, we have the approximate monotonicity inequality
	\begin{align} \label{E:APPROXIMATETIMEMONOTONICITYGRONWALLPROOF}
		(1 + C \varepsilon) \upmu_{\star}(t',u') 
		& \geq \upmu_{\star}(\hat{t},\hat{u}).
	\end{align}
	We use these monotonicity properties below without
	explicitly mentioning them each time.

	We now set $N=20$, 
	multiply both sides of inequality
	\eqref{E:TOPORDERWAVEENERGYINTEGRALINEQUALITIES}
	by $\iota_F^{-1}(t,u)$ and then set $(t,u) = (\hat{t},\hat{u})$.
	Similarly, we multiply both sides of the inequality described in
	\eqref{E:EASYVARIABLESTOPORDERWAVEENERGYINTEGRALINEQUALITIES}
	by $\iota_G^{-1}(t,u)$
	and the inequality
	\eqref{E:BELOWTOPORDERWAVEENERGYINTEGRALINEQUALITIES}
	by $\iota_H^{-1}(t,u)$ and, in both cases,
	set $(t,u) = (\hat{t},\hat{u})$.
	To deduce \eqref{E:TOPORDERALMOSTESTIMATED}-\eqref{E:BELOWTOPORDERALMOSTESTIMATED},
	the difficult step is to obtain suitable bounds
	for the terms generated by the terms on 
	RHSs~\eqref{E:TOPORDERWAVEENERGYINTEGRALINEQUALITIES}-\eqref{E:BELOWTOPORDERWAVEENERGYINTEGRALINEQUALITIES}.
	Once we have obtained suitable bounds, 
	we can then take 
	$\sup_{(\hat{t},\hat{u}) \in [0,t] \times [0,u]}$
	of both sides of the resulting inequalities, and by virtue of
	definitions
	\eqref{E:TOPORDERENERGYRESCALED}-\eqref{E:JUSTBELOWTOPORDERENERGYRESCALED},
	we will easily
	conclude \eqref{E:TOPORDERALMOSTESTIMATED}-\eqref{E:BELOWTOPORDERALMOSTESTIMATED}.

	We start by showing how to obtain suitable bounds
	for the terms on RHS~\eqref{E:TOPORDERWAVEENERGYINTEGRALINEQUALITIES} 
	that involve the vorticity energies.
	These estimates are easy to derive because we have already derived suitable
	estimates for the vorticity energies. Specifically, we must handle the terms
	\begin{align}
		&
		C 
		\sup_{(\hat{t},\hat{u}) \in [0,t] \times [0,u]}
		\iota_F^{-1}(\hat{t},\hat{u})
		\int_{t'=0}^{\hat{t}}
			\frac{1}{\upmu_{\star}^{3/2}(t',\hat{u})}
			\left\lbrace
				\int_{s=0}^{t'}
					\sqrt{\VorttotTanmax{21}}(s,\hat{u})
				\, ds
			\right\rbrace^2
		\, dt',
			\label{E:MAINPROOFTOPORDERVORTICITYTERMAINWAVEENERGIESDOUBLETIMEINGEGRATED}	 
			\\
		& 
			C 
			\sup_{(\hat{t},\hat{u}) \in [0,t] \times [0,u]}
			\iota_F^{-1}(\hat{t},\hat{u})
			\int_{t'=0}^{\hat{t}}
				\frac{1}{\upmu_{\star}^{3/2}(t',\hat{u})}
				\left\lbrace
					\int_{s=0}^{t'}
						\frac{1} 
						{\upmu_{\star}^{1/2}(s,\hat{u})}
						\sqrt{\VorttotTanmax{\leq 20}}(s,\hat{u})
					\, ds
				\right\rbrace^2
			\, dt',
			\label{E:MAINPROOFBELOWTOPORDERVORTICITYTERMAINWAVEENERGIESDOUBLETIMEINTEGRATED}
				\\
		&
		C
		\sup_{(\hat{t},\hat{u}) \in [0,t] \times [0,u]}
		\iota_F^{-1}(\hat{t},\hat{u})
		\int_{t'=0}^{\hat{t}}
		\VorttotTanmax{\leq 21}(t',u)
		\, dt',
			\label{E:MAINPROOFTOPORDERVORTICITYTERMAINWAVEENERGIESTIMEINGEGRATED} 
			\\
		&
		C
		\sup_{(\hat{t},\hat{u}) \in [0,t] \times [0,u]}
		\iota_F^{-1}(\hat{t},\hat{u})
		\int_{u'=0}^{\hat{u}}
			\VorttotTanmax{\leq 20}(\hat{t},u')
		\, du'
		\label{E:MAINPROOFBELOWTOPORDERVORTICITYTERMAINWAVEENERGIESUINGEGRATED}
	\end{align}
	generated by the integrals on
	the last three lines of RHS~\eqref{E:TOPORDERWAVEENERGYINTEGRALINEQUALITIES}.
	To proceed, we insert the already proven vorticity estimates
	\eqref{E:TOPVORTMULOSSMAINAPRIORIENERGYESTIMATES}-\eqref{E:VORTNOMULOSSMAINAPRIORIENERGYESTIMATES}
	into the integrands in
	\eqref{E:MAINPROOFTOPORDERVORTICITYTERMAINWAVEENERGIESDOUBLETIMEINGEGRATED}-\eqref{E:MAINPROOFBELOWTOPORDERVORTICITYTERMAINWAVEENERGIESUINGEGRATED}.
	With the help of inequality \eqref{E:LOSSKEYMUINTEGRALBOUND},
	we obtain
	\begin{align}  
		&
		C 
		\int_{t'=0}^{\hat{t}}
			\frac{1}{\upmu_{\star}^{3/2}(t',\hat{u})}
			\left\lbrace
				\int_{s=0}^{t'}
					\sqrt{\VorttotTanmax{21}}(s,\hat{u})
				\, ds
			\right\rbrace^2
		\, dt'
			\label{E:BOUNDEDMAINPROOFTOPORDERVORTICITYTERMAINWAVEENERGIESDOUBLETIMEINGEGRATED}
				\\
		& \leq 
		C \mathring{\upepsilon}^2
		\int_{t'=0}^{\hat{t}}
			\upmu_{\star}^{-3/2}(t',\hat{u})
			\left\lbrace
				\int_{s=0}^{t'}
					\upmu_{\star}^{-6.4}(s,\hat{u})
				\, ds
			\right\rbrace^2
		\, dt'
			 \notag \\
		& \leq 
		C \mathring{\upepsilon}^2
		\int_{t'=0}^{\hat{t}}
			\frac{1}{\upmu_{\star}^{12.3}(t',\hat{u})}
		\, dt'
		\leq C \mathring{\upepsilon}^2 \upmu_{\star}^{-11.3}(\hat{t},\hat{u}),
			\notag \\
		&
		\int_{t'=0}^{\hat{t}}
				\frac{1}{\upmu_{\star}^{3/2}(t',\hat{u})}
				\left\lbrace
					\int_{s=0}^{t'}
					\frac{1} 
					{\upmu_{\star}^{1/2}(s,\hat{u})}
					\sqrt{\VorttotTanmax{\leq 20}}(s,\hat{u}) \, ds
				\right\rbrace^2
				\, dt'
				\label{E:BOUNDEDMAINPROOFBELOWTOPORDERVORTICITYTERMAINWAVEENERGIESDOUBLETIMEINTEGRATED} 
				\\
	& \leq 
		C \mathring{\upepsilon}^2
		\int_{t'=0}^{\hat{t}}
			\upmu_{\star}^{-3/2}(t',\hat{u})
			\left\lbrace
				\int_{s=0}^{t'}
					\upmu_{\star}^{-5.4}(s,\hat{u})
				\, ds
			\right\rbrace^2
		\, dt'
			 \notag \\
		& \leq 
		C \mathring{\upepsilon}^2
		\int_{t'=0}^{\hat{t}}
			\upmu_{\star}^{-10.3}(t',\hat{u})
		\, dt'
		\leq C \mathring{\upepsilon}^2 \upmu_{\star}^{-9.3}(\hat{t},\hat{u}),
			\notag 
	\end{align}
	\begin{align}
		C
		\int_{t'=0}^{\hat{t}}
			\VorttotTanmax{\leq 21}(t',\hat{u})
		\, dt'
		& \leq
		C \mathring{\upepsilon}^2
		\int_{t'=0}^{\hat{t}}
			\upmu_{\star}^{-12.8}(t',\hat{u})
		\, dt'
		\leq
		C \mathring{\upepsilon}^2 \upmu_{\star}^{-11.8}(\hat{t},\hat{u}),
			\label{E:BOUNDEDMAINPROOFTOPORDERVORTICITYTERMAINWAVEENERGIESTIMEINGEGRATED} \\
		C
		\int_{u'=0}^{\hat{u}}
			\VorttotTanmax{\leq 20}(\hat{t},u')
		\, du'
		& \leq 
			C \VorttotTanmax{\leq 20}(\hat{t},\hat{u})
		\leq C \mathring{\upepsilon}^2 \upmu_{\star}^{-9.8}(\hat{t},\hat{u}).
		\label{E:BOUNDEDMAINPROOFBELOWTOPORDERVORTICITYTERMAINWAVEENERGIESUINGEGRATED}
		\end{align}
		Multiplying
		\eqref{E:BOUNDEDMAINPROOFTOPORDERVORTICITYTERMAINWAVEENERGIESDOUBLETIMEINGEGRATED}-\eqref{E:BOUNDEDMAINPROOFBELOWTOPORDERVORTICITYTERMAINWAVEENERGIESUINGEGRATED} 
		by $\iota_F^{-1}(\hat{t},\hat{u})$
		and then taking $\sup_{(\hat{t},\hat{u}) \in [0,t] \times [0,u]}$, we conclude
		\begin{align}
			\mbox{\eqref{E:MAINPROOFTOPORDERVORTICITYTERMAINWAVEENERGIESDOUBLETIMEINGEGRATED}}
			& \leq 
			C \mathring{\upepsilon}^2
			\sup_{(\hat{t},\hat{u}) \in [0,t] \times [0,u]}
			\upmu_{\star}^{.5}(\hat{t},\hat{u})
			\iota_1^{-c}(\hat{t})
			\iota_2^{-c}(\hat{t},\hat{u})
			e^{-c \hat{t}}
			e^{-c \hat{u}}
			\leq C \mathring{\upepsilon}^2,
				\label{E:SUPEDBOUNDEDMAINPROOFTOPORDERVORTICITYTERMAINWAVEENERGIESDOUBLETIMEINGEGRATED} 
				\\
			\mbox{\eqref{E:MAINPROOFBELOWTOPORDERVORTICITYTERMAINWAVEENERGIESDOUBLETIMEINTEGRATED}}
			& \leq 
			C \mathring{\upepsilon}^2
			\sup_{(\hat{t},\hat{u}) \in [0,t] \times [0,u]}
			\upmu_{\star}^{2.5}(\hat{t},\hat{u})
			\iota_1^{-c}(\hat{t})
			\iota_2^{-c}(\hat{t},\hat{u})
			e^{-c \hat{t}}
			e^{-c \hat{u}}
			\leq C \mathring{\upepsilon}^2,
			\label{E:SUPEDBOUNDEDMAINPROOFBELOWTOPORDERVORTICITYTERMAINWAVEENERGIESDOUBLETIMEINTEGRATED}
				\\
			\mbox{\eqref{E:MAINPROOFTOPORDERVORTICITYTERMAINWAVEENERGIESTIMEINGEGRATED}}
			& \leq 
			C \mathring{\upepsilon}^2
			\sup_{(\hat{t},\hat{u}) \in [0,t] \times [0,u]}
			\iota_1^{-c}(\hat{t})
			\iota_2^{-c}(\hat{t},\hat{u})
			e^{-c \hat{t}}
			e^{-c \hat{u}}
			\leq C \mathring{\upepsilon}^2,
				\label{E:SUPEDBOUNDEDMAINPROOFTOPORDERVORTICITYTERMAINWAVEENERGIESTIMEINGEGRATED} 
				\\
			\mbox{\eqref{E:MAINPROOFBELOWTOPORDERVORTICITYTERMAINWAVEENERGIESUINGEGRATED}}
			& \leq 
			C \mathring{\upepsilon}^2
			\sup_{(\hat{t},\hat{u}) \in [0,t] \times [0,u]}
			\upmu_{\star}^2(\hat{t},\hat{u})
			\iota_1^{-c}(\hat{t})
			\iota_2^{-c}(\hat{t},\hat{u})
			e^{-c \hat{t}}
			e^{-c \hat{u}}
			\leq C \mathring{\upepsilon}^2
			\label{E:SUPEDBOUNDEDMAINPROOFBELOWTOPORDERVORTICITYTERMAINWAVEENERGIESUINGEGRATED}
		\end{align}
		as desired.
		We have thus accounted for the influence of the vorticity 
		in the top-order wave energies.


	We now show how to obtain suitable bounds for 
	the terms generated by
	the ``borderline'' terms
	$\boxed{6} \int \cdots$,
	$\boxed{8.1} \int \cdots$, 
	and 
	$
	\displaystyle
	\boxed{2} \frac{1}{\upmu_{\star}^{1/2}(t,u)} 
	\sqrt{\totmax{20}}(t,u)
	\left\| 
				\Lunit \upmu 
			\right\|_{L^{\infty}(\Sigmaminus{t}{t}{u})}
	\int \cdots
	$
	on RHS~\eqref{E:TOPORDERWAVEENERGYINTEGRALINEQUALITIES}
	(where we recall that $N=20$ in this part of the proof).
	The terms generated by the
	remaining ``non-borderline'' terms on RHS~\eqref{E:TOPORDERWAVEENERGYINTEGRALINEQUALITIES} 
	are easier to treat.
	We start with the term 
	$\boxed{6} 
		\iota_F^{-1}(\hat{t},\hat{u}) \int_{t'=0}^{\hat{t}} \cdots
	$.
	Multiplying and 
	dividing by
	$\upmu_{\star}^{11.8}(t', \hat{u})$ in the integrand, 
	taking 
	$
	\displaystyle
	\sup_{t' \in [0,\hat{t}]} \upmu_{\star}^{11.8}(t',\hat{u}) \totmax{20}(t',\hat{u})
	$,
	pulling the $\sup$-ed quantity out of the integral,
	and using the critically important integral estimate
	\eqref{E:KEYMUTOAPOWERINTEGRALBOUND} with $\Contwo = 12.8$, 
	we find that
	\begin{align} \label{E:MOSTDANGEROUSINTEGRALESTIMATE}
	& \boxed{6}
			\iota_F^{-1}(\hat{t},\hat{u})
			\int_{t'=0}^{\hat{t}}
					\frac{\left\| [\Lunit \upmu]_- \right\|_{L^{\infty}(\Sigma_{t'}^{\hat{u}})}} 
							 {\upmu_{\star}(t',\hat{u})} 
				  \totmax{20}(t',\hat{u})
			\, dt'
				\\
		& \leq 
				\boxed{6}
				\iota_F^{-1}(\hat{t},\hat{u})
				\sup_{t' \in [0,\hat{t}]}
				\left\lbrace 
					\upmu_{\star}^{11.8}(t',\hat{u}) \totmax{20}(t',\hat{u})
				\right\rbrace
				\int_{t'=0}^{\hat{t}}
					\left\| [\Lunit \upmu]_- \right\|_{L^{\infty}(\Sigma_{t'}^{\hat{u}})}
					\upmu_{\star}^{-12.8}(t',\hat{u})
				\, dt'
				\notag \\
		& \leq 
				\boxed{6}
				\upmu_{\star}^{11.8}(\hat{t},\hat{u})
				\sup_{t' \in [0,\hat{t}]}
				\left\lbrace 
					\iota_1^{-c}(t')
					\iota_2^{-c}(t',\hat{u})
					e^{-c t'} 
					e^{-c \hat{u}}
					\upmu_{\star}^{11.8}(t',\hat{u}) 
					\totmax{20}(t',\hat{u})
				\right\rbrace
					\notag \\
			& \ \ \ \
					\times
					\int_{t'=0}^{\hat{t}}
					\left\| [\Lunit \upmu]_- \right\|_{L^{\infty}(\Sigma_{t'}^{\hat{u}})}
					\upmu_{\star}^{-12.8}(t',\hat{u})
				\, dt'
				\notag \\
		&   \leq
				\frac{6 + C \sqrt{\varepsilon}}{11.8}
				F(\hat{t},\hat{u})
				\leq
				\frac{6 + C \sqrt{\varepsilon}}{11.8}
				F(t,u).
				\notag
	\end{align}

To handle the integral 
$\boxed{8.1} \iota_F^{-1}(\hat{t},\hat{u}) \int \cdots$,
we use a similar argument, but this time taking into account
that there are two time integrations. We find that
\begin{align} \label{E:SECONDMOSTDANGEROUSINTEGRALESTIMATE}
	&
	\boxed{8.1}
			\iota_F^{-1}(\hat{t},\hat{u})
			\int_{t'=0}^{\hat{t}}
				\frac{\left\| [\Lunit \upmu]_- \right\|_{L^{\infty}(\Sigma_{t'}^{\hat{u}})}} 
								 {\upmu_{\star}(t',\hat{u})} 
				\sqrt{\totmax{20}}(t',\hat{u}) 
						\int_{s=0}^{t'}
							\frac{\left\| [\Lunit \upmu]_- \right\|_{L^{\infty}(\Sigma_s^{\hat{u}})}} 
									{\upmu_{\star}(s,\hat{u})} 
							\sqrt{\totmax{20}}(s,\hat{u}) 
						\, ds
				\, dt'
				\\
		& \leq
				\frac{8.1 + C \sqrt{\varepsilon}}{5.9 \times 11.8}
				F(t,u).
				\notag
\end{align}

To handle the integral $\boxed{2} \iota_F^{-1}(\hat{t},\hat{u}) \int \cdots$,
we use a similar argument based on the critically important estimate 
\eqref{E:KEYHYPERSURFACEMUTOAPOWERINTEGRALBOUND}.
We find that
\begin{align} \label{E:THIRDMOSTDANGEROUSINTEGRALESTIMATE}
	&
	\boxed{2}
	\iota_F^{-1}(\hat{t},\hat{u})
	\frac{1}{\upmu_{\star}^{1/2}(\hat{t},\hat{u})}
	\sqrt{\totmax{20}}(\hat{t},\hat{u})
	\left\| 
			\Lunit \upmu 
	\right\|_{L^{\infty}(\Sigmaminus{\hat{t}}{\hat{t}}{\hat{u}})}
	\int_{t'=0}^t
			\frac{1}{\upmu_{\star}^{1/2}(t',\hat{u})} \sqrt{\totmax{20}}(t',\hat{u})
	\, dt'
		\\
	& \leq \frac{2 + C \sqrt{\varepsilon}}{5.4} F(t,u).
	\notag
\end{align}

The important point is that for small $\varepsilon$,
the factors 
$\displaystyle \frac{6 + C \sqrt{\varepsilon}}{11.8}$ on RHS
\eqref{E:MOSTDANGEROUSINTEGRALESTIMATE},
$\displaystyle \frac{8.1 + C \sqrt{\varepsilon}}{5.9 \times 11.8}$
on RHS \eqref{E:SECONDMOSTDANGEROUSINTEGRALESTIMATE},
and $\displaystyle\frac{2 + C \sqrt{\varepsilon}}{5.4}$
on RHS \eqref{E:THIRDMOSTDANGEROUSINTEGRALESTIMATE}
sum to 
$\displaystyle 
\frac{6}{11.8} 
+ \frac{8.1}{5.9 \times 11.8} 
+ \frac{2}{5.4} 
+ C \sqrt{\varepsilon} 
< 1
$.
This sum is the main contributor to the constant
$\upalpha_1$ on RHS \eqref{E:TOPORDERALMOSTESTIMATED}.

	We now derive suitable bounds for the three terms
	on RHS~\eqref{E:TOPORDERWAVEENERGYINTEGRALINEQUALITIES} 
	that are multiplied by the large constant $C_{\ast}$.
	We bound these terms using essentially the same reasoning
	that we used in proving 
	\eqref{E:MOSTDANGEROUSINTEGRALESTIMATE},
	\eqref{E:SECONDMOSTDANGEROUSINTEGRALESTIMATE},
	and \eqref{E:THIRDMOSTDANGEROUSINTEGRALESTIMATE},
	but we use only the crude inequality \eqref{E:LOSSKEYMUINTEGRALBOUND}
	in place of the delicate inequalities \eqref{E:KEYMUTOAPOWERINTEGRALBOUND} 
	and \eqref{E:KEYHYPERSURFACEMUTOAPOWERINTEGRALBOUND}.
	We find that
	\begin{align}
		 C_{\ast}
		\iota_F^{-1}(\hat{t},\hat{u})
		&
		\int_{t'=0}^{\hat{t}}
					\frac{1} 
							 {\upmu_{\star}(t',\hat{u})} 
				  \sqrt{\totmax{20}}(t',\hat{u}) 
					\sqrt{\easytotmax{20}}(t',\hat{u}) 
		\, dt'
			\label{E:FIRSTLARGECONSTANTTERM} \\
		& \leq C F^{1/2}(t,u)G^{1/2}(t,u),
				\notag \\
			C_{\ast}
			\iota_F^{-1}(\hat{t},\hat{u})
			&
			\int_{t'=0}^{\hat{t}}
				\frac{1}
					{\upmu_{\star}(t',\hat{u})} 
						\sqrt{\totmax{20}}(t',\hat{u}) 
						\int_{s=0}^{t'}
							\frac{1} 
									{\upmu_{\star}(s,\hat{u})} 
							\sqrt{\easytotmax{20}}(s,\hat{u}) 
						\, ds
				\, dt'
				\label{E:SECONDLARGECONSTANTTERM}
					\\
			& \leq C F^{1/2}(t,u)G^{1/2}(t,u),
			\notag	\\
			C_{\ast}
			\iota_F^{-1}(\hat{t},\hat{u})
			&
			\frac{1}{\upmu_{\star}^{1/2}(\hat{t},\hat{u})}
			\sqrt{\totmax{20}}(\hat{t},\hat{u})
			\int_{t'=0}^{\hat{t}}
				\frac{1}{\upmu_{\star}^{1/2}(t',\hat{u})} \sqrt{\easytotmax{20}}(t',\hat{u})
			\, dt'
				\label{E:THIRDLARGECONSTANTTERM} \\
			& \leq C F^{1/2}(t,u)G^{1/2}(t,u),
			\notag
	\end{align}
	where the constants $C$ on RHSs\eqref{E:FIRSTLARGECONSTANTTERM}-\eqref{E:THIRDLARGECONSTANTTERM}
	are large and sum to the large constant $\upalpha_2$ on RHS~\eqref{E:TOPORDERALMOSTESTIMATED}.
	We remark that the largeness of $\upalpha_2$ 
	will not preclude us from closing the estimates because we will
	gain smallness in $G(t,u)$ by using a separate argument given below.

The remaining integrals 
on RHS~\eqref{E:TOPORDERWAVEENERGYINTEGRALINEQUALITIES}
are easier to treat.
We now show how to bound the term arising from the integral
on the $14^{th}$ line of RHS \eqref{E:TOPORDERWAVEENERGYINTEGRALINEQUALITIES}, 
which involves three time integrations.
The term arising from the integrals on the $13^{th}$ line of 
RHS~\eqref{E:TOPORDERWAVEENERGYINTEGRALINEQUALITIES}
can be handled using similar arguments, so we do not provide those details.
We claim that the following sequence of inequalities holds for the term of interest,
which yields the desired bound:
\begin{align} \label{E:ESTIMATEFORTRIPLETIMEINTEGRATEDERRORINTEGRAL}
					& 
					C
					\iota_F^{-1}(\hat{t},\hat{u})
					\int_{t'=0}^{\hat{t}}
					\frac{1} 
							 {\upmu_{\star}(t',\hat{u})} 
				  \sqrt{\totmax{20}}(t',\hat{u})
				  \int_{s = 0}^{t'}
				  	\frac{1}{\upmu_{\star}(s,\hat{u})}
				  	\int_{s' = 0}^s
				  	\frac{1} 
							 {\upmu_{\star}^{1/2}(s',\hat{u})} 
							 \sqrt{\totmax{20}}(s',\hat{u})
						\, ds'	 
				  \, ds
				\, dt'
					\\
				& \leq
					\frac{C}{c} 
					\iota_F^{-1}(\hat{t},\hat{u})
					\iota_2^{c/2}(\hat{t},\hat{u})
					\int_{t'=0}^{\hat{t}}
					\frac{1} 
							 {\upmu_{\star}(t',\hat{u})} 
				  \sqrt{\totmax{20}}(t',\hat{u})
				  	\notag \\
				 & \ \
				 		\times
				 		\int_{s = 0}^{t'}
				  	\frac{1}{\upmu_{\star}(s,\hat{u})}
				  	\sup_{(s',u') \in [1,s] \times [0,\hat{u}]}
						\left\lbrace 
							\iota_2^{-c/2}(s',u') \sqrt{\totmax{20}}(s',u')
						\right\rbrace
				  	\, ds
				\, dt'
					\notag \\
			& \leq
					\frac{C}{c} 
					\iota_F^{-1}(\hat{t},\hat{u})
					\iota_2^{c/2}(\hat{t},\hat{u})
					\int_{t'=0}^{\hat{t}}
					\frac{1} 
							 {\upmu_{\star}(t',\hat{u})} 
				  \sqrt{\totmax{20}}(t',\hat{u})
				  \notag \\
			& \ \ 
						\times
						\sup_{(s',u') \in [0,t'] \times [0,\hat{u}]}
						\left\lbrace 
							\upmu_{\star}(s',u') \iota_2^{-c/2}(s',u') \sqrt{\totmax{20}}(s',u')
						\right\rbrace
				  \int_{s = 0}^{t'}
				  	\frac{1}{\upmu_{\star}^2(s,\hat{u})}
				  \, ds
				\, dt'
					\notag \\
		& \leq
					\frac{C}{c} 
					\iota_F^{-1}(\hat{t},\hat{u})
					\iota_2^{c/2}(\hat{t},\hat{u})
					\sup_{(s',u') \in [0,\hat{t}] \times [0,\hat{u}]}
						\left\lbrace 
							\upmu_{\star}(s',u') \iota_2^{-c/2}(s',u') \sqrt{\totmax{20}}(s',u')
						\right\rbrace
						\notag \\
		& \times 
					\sup_{(s',u') \in [0,\hat{t}] \times [0,\hat{u}]}
						\left\lbrace 
							\sqrt{\totmax{20}}(s',u')
						\right\rbrace
					\int_{t'=0}^{\hat{t}}
						\frac{1} 
							 {\upmu_{\star}(t',\hat{u})} 
				  	\int_{s = 0}^{t'}
				  		\frac{1}{\upmu_{\star}^2(s,\hat{u})}
				  	\, ds
				\, dt'
				\notag \\
		& \leq
					\frac{C}{c} 
					\upmu_{\star}(\hat{t},\hat{u})
					\sup_{(s',u') \in [0,\hat{t}] \times [0,\hat{u}]}
						\left\lbrace 
							\iota_F^{-1}(s',u') \totmax{20}(s',u')
						\right\rbrace
				\times 
					\int_{t'=0}^{\hat{t}}
						\frac{1} 
							 {\upmu_{\star}(t',\hat{u})} 
				  	\int_{s = 0}^{t'}
				  		\frac{1}{\upmu_{\star}^2(s,\hat{u})}
				  	\, ds
				\, dt'
				\notag \\
		& \leq \frac{C}{c} F(\hat{t},\hat{u})
			\leq \frac{C}{c} F(t,u)
			\notag,
\end{align}
which yields the desired smallness factor $\displaystyle \frac{1}{c}$.
We now explain how to derive 
\eqref{E:ESTIMATEFORTRIPLETIMEINTEGRATEDERRORINTEGRAL}.
To deduce the first inequality,
we multiplied and divided by
$\iota_2^{c/2}(t',\hat{u})$ in the integral 
$\int \cdots ds'$,
then pulled
$\displaystyle 
\sup_{(s',u') \in [0,s] \times [0,\hat{u}]}
						\left\lbrace 
							\iota_2^{-c/2}(s',u') \sqrt{\totmax{20}}(s',u')
						\right\rbrace
$
out of the integral, and finally used \eqref{E:IF2SMALLNESS}
to gain the smallness factor $\displaystyle \frac{1}{c}$
from the remaining terms
$
\displaystyle
\int_{s'=0}^s
	\frac{1} 
		{\upmu_{\star}^{1/2}(s',\hat{u})} 
		\iota_2^{c/2}(s',\hat{u})
\, ds'
$.
To derive the second
inequality in \eqref{E:ESTIMATEFORTRIPLETIMEINTEGRATEDERRORINTEGRAL},
we multiplied and divided by
$\upmu_{\star}(s,\hat{u})$ in the integral 
$\int \cdots ds$,
and used the approximate monotonicity property
\eqref{E:APPROXIMATETIMEMONOTONICITYGRONWALLPROOF}
to pull the factor 
$
\displaystyle
\sup_{(s',u') \in [0,t'] \times [0,\hat{u}]}
						\left\lbrace 
							\upmu_{\star}(s',u') \iota_2^{-c/2}(s',u') \sqrt{\totmax{20}}(s',u')
						\right\rbrace
$
out of the $ds$ integral, which costs us a harmless multiplicative factor of $1 + C \varepsilon$.
The third inequality in \eqref{E:ESTIMATEFORTRIPLETIMEINTEGRATEDERRORINTEGRAL}
follows easily. To derive the fourth inequality,
we use the monotonicity of
$\iota_1^c(\cdot)$,
$\iota_2^c(\cdot)$
and $e^{c \cdot}$,
and the approximate monotonicity property
\eqref{E:APPROXIMATETIMEMONOTONICITYGRONWALLPROOF}.
To derive the fifth inequality, 
we use inequality \eqref{E:LOSSKEYMUINTEGRALBOUND} twice.
The final inequality follows easily.

Similarly, we claim that we can bound the terms on the 
$8^{th}$ through $11^{th}$ lines of RHS~\eqref{E:TOPORDERWAVEENERGYINTEGRALINEQUALITIES} 
and the second term on the $12^{th}$ line of RHS~\eqref{E:TOPORDERWAVEENERGYINTEGRALINEQUALITIES}
as follows:
\begin{align} \label{E:FIRSTEASYHYPERSUFRACEINTEGRALENERGYESTIMATETERM}
			\iota_F^{-1}(\hat{t},\hat{u})
			C
			\varepsilon
			\int_{t'=0}^{\hat{t}}
					\frac{1} 
						{\upmu_{\star}(t',\hat{u})} 
				  \totmax{20}(t',\hat{u})
				\, dt'
			& \leq 
				C \varepsilon F(\hat{t},\hat{u})
				\leq C \varepsilon F(t,u), 
				\\
		\iota_F^{-1}(\hat{t},\hat{u})
		C \varepsilon
			\int_{t'=0}^{\hat{t}}
				\frac{1} 
						{\upmu_{\star}(t',\hat{u})} 
						\sqrt{\totmax{20}}(t',\hat{u}) 
						\int_{s=0}^{t'}
							\frac{1} 
									{\upmu_{\star}(s,\hat{u})} 
							\sqrt{\totmax{20}}(s,\hat{u}) 
						\, ds
				\, dt'
				& \leq 
				C \varepsilon F(\hat{t},\hat{u})
				\leq C \varepsilon F(t,u), 
					\label{E:ANOTHERFIRSTEASYHYPERSUFRACEINTEGRALENERGYESTIMATETERM} 
					\\
			\iota_F^{-1}(\hat{t},\hat{u})
			C
			\varepsilon
			\frac{1}{\upmu_{\star}^{1/2}(\hat{t},\hat{u})}
			\sqrt{\totmax{20}}(\hat{t},\hat{u})
			\int_{t'=0}^{\hat{t}}
				\frac{1}{\upmu_{\star}^{1/2}(t',\hat{u})} \sqrt{\totmax{20}}(t',\hat{u})
			\, dt'
			& \leq 
				C \varepsilon F(\hat{t},\hat{u})
				\leq C \varepsilon F(t,u), 
				\label{E:SECONDEASYHYPERSUFRACEINTEGRALENERGYESTIMATETERM} 
				\\
			\iota_F^{-1}(\hat{t},\hat{u})
			C
			\sqrt{\totmax{20}}(\hat{t},\hat{u})
			\int_{t'=0}^{\hat{t}}
				\frac{1}{\upmu_{\star}^{1/2}(t',\hat{u})} \sqrt{\totmax{20}}(t',\hat{u})
			\, dt'
			& 
			\leq \frac{C}{c} F(\hat{t},\hat{u})
			\leq \frac{C}{c} F(t,u),
			\label{E:THIRDEASYHYPERSUFRACEINTEGRALENERGYESTIMATETERM}
	\end{align}
	\begin{align}
			\iota_F^{-1}(\hat{t},\hat{u})
			C (1 + \varsigma^{-1})
					\int_{t'=0}^{\hat{t}}
					\frac{1} 
							 {\upmu_{\star}^{1/2}(t',\hat{u})} 
				  \totmax{20}(t',\hat{u})
				\, dt'
			& \leq 
				\frac{C}{c} (1 + \varsigma^{-1}) F(\hat{t},\hat{u})
				\leq \frac{C}{c} (1 + \varsigma^{-1}) F(t,u).
				\label{E:ANOTHEREASYHYPERSUFRACEINTEGRALENERGYESTIMATETERM}
		\end{align}
To derive \eqref{E:FIRSTEASYHYPERSUFRACEINTEGRALENERGYESTIMATETERM},
we use arguments similar to the ones we used 
in deriving \eqref{E:MOSTDANGEROUSINTEGRALESTIMATE}, 
but in place of the delicate estimate \eqref{E:KEYMUTOAPOWERINTEGRALBOUND},
we use the estimate \eqref{E:LOSSKEYMUINTEGRALBOUND},
whose imprecision is compensated for by the availability of the smallness 
factor $\varepsilon$.
Similar remarks apply to \eqref{E:ANOTHERFIRSTEASYHYPERSUFRACEINTEGRALENERGYESTIMATETERM},
but we rely on the fact that there are two time integrations.
The proof of \eqref{E:ANOTHEREASYHYPERSUFRACEINTEGRALENERGYESTIMATETERM}
is similar, but we multiply and divide by 
by $\iota_2^{-c}(t',\hat{u})$ in the integrand
and use the estimate
\eqref{E:IF2SMALLNESS} to gain the smallness factor
$
\displaystyle \frac{1}{c}
$.
To derive \eqref{E:SECONDEASYHYPERSUFRACEINTEGRALENERGYESTIMATETERM},
we use arguments similar to the ones we used above,
but we now multiply and divide by 
$\upmu_{\star}^{5.9}(t',\hat{u})$
in the time integral on LHS \eqref{E:SECONDEASYHYPERSUFRACEINTEGRALENERGYESTIMATETERM}
and use \eqref{E:LOSSKEYMUINTEGRALBOUND}.
To derive \eqref{E:THIRDEASYHYPERSUFRACEINTEGRALENERGYESTIMATETERM},
we use similar arguments based on multiplying and dividing by
$\iota_2^{c/2}(t',\hat{u})$ in the time integral
and using \eqref{E:IF2SMALLNESS}.

Similarly, we derive the bound
\begin{align}
	&
	C
	\iota_F^{-1}(\hat{t},\hat{u})
	\int_{t'=0}^{\hat{t}}
		\frac{1}{\sqrt{\Tboot - t'}} \totmax{20}(t',\hat{u})
	\, dt'
	\leq \frac{C}{c} F(\hat{t},\hat{u})
	\leq \frac{C}{c} F(t,u)
\end{align}
for the first term on the $12^{th}$ line of RHS \eqref{E:TOPORDERWAVEENERGYINTEGRALINEQUALITIES}
by multiplying and dividing by $\iota_1^{c}(t')$
in the integrand and using
\eqref{E:IF1SMALLNESS} to gain the smallness factor
$
\displaystyle
\frac{1}{c}
$.

Similarly, we derive the bound
\begin{align}
	&
	C
	(1 + \varsigma^{-1})
	\iota_F^{-1}(\hat{t},\hat{u})
	\int_{u'=0}^{\hat{u}}
		\totmax{20}(\hat{t},u')
	\, du'
	\leq \frac{C}{c} F(\hat{t},\hat{u})
	\leq \frac{C}{c} F(t,u)
\end{align}
for the first term on the $15^{th}$ line of RHS \eqref{E:TOPORDERWAVEENERGYINTEGRALINEQUALITIES}
by multiplying and dividing by $e^{cu'}$
in the integrand and using \eqref{E:IF4SMALLNESS} to 
gain the smallness factor
$
\displaystyle
\frac{1}{c}
$.

It is easy to see that the terms 
arising from the term 
on the 
first line and the last three terms on the $15^{th}$
line of RHS \eqref{E:TOPORDERWAVEENERGYINTEGRALINEQUALITIES},
namely
$
C (1 + \varsigma^{-1}) \mathring{\upepsilon}^2 \upmu_{\star}^{-3/2}(t,u)
$,
$
C \varepsilon \totmax{20}(t,u)
$,
$
C \varsigma \totmax{20}(t,u)
$,
and
$
C \varsigma \coercivespacetimemax{20}(t,u)
$,
are respectively bounded 
(after multiplying by $\iota_F^{-1}$ and taking the relevant sup)
by 
$\leq C (1 + \varsigma^{-1}) \mathring{\upepsilon}^2
$,
$
\leq C \varepsilon F(t,u)
$,
$
\leq
C \varsigma F(t,u)
$,
and 
$
\leq
C \varsigma F(t,u)
$.

To bound the term arising from the first term on 
the $16^{th}$ line of RHS~\eqref{E:TOPORDERWAVEENERGYINTEGRALINEQUALITIES}
(where we recall that $N=20$), 
we argue as follows with the help of
\eqref{E:IF2SMALLNESS} and \eqref{E:APPROXIMATETIMEMONOTONICITYGRONWALLPROOF}:
\begin{align} \label{E:BELOWTOPORDERTERMCOUPLEDINTOTOPORDERESTIMATE}
			& C (1 + \varsigma^{-1}) 
			\iota_F^{-1}(\hat{t},\hat{u})
			\int_{t'=0}^{\hat{t}}
					\frac{1} 
							 {\upmu_{\star}^{5/2}(t',\hat{u})} 
				  \totmax{[1,19]}(t',\hat{u})
			\, dt'
				\\
			& 
			\leq
			C (1 + \varsigma^{-1}) 
			\upmu_{\star}^{9.8}(\hat{t},\hat{u})
			\iota_1^{-c}(\hat{t})
			\iota_2^{-c}(\hat{t},\hat{u})
			e^{-c \hat{t}} 
			e^{-c \hat{u}}
			\sup_{t' \in [0,\hat{t}]}
				\left(
					\frac{\upmu_{\star}(\hat{t},\hat{u})} 
					{\upmu_{\star}(t',\hat{u})} 
				\right)^2
					\notag \\
	& \ \
			\times
			\left\lbrace
				\sup_{t' \in [0,\hat{t}]}
				\iota_2^{-c}(t') \totmax{[1,19]}(t',\hat{u})
			\right\rbrace
			\times
			\int_{t'=0}^{\hat{t}}
				\frac{\iota_2^c(t')} 
				{\upmu_{\star}^{1/2}(t',\hat{u})} 
			\, dt'
				\notag 
				\\
		& 
			\leq
			C (1 + \varsigma^{-1}) 
			\iota_2^{-c}(\hat{t})
			\left\lbrace 
				\sup_{t' \in [0,\hat{t}]}
				\iota_H^{-1}(t',\hat{u})\totmax{[1,19]}(t',\hat{u})
			\right\rbrace
			\times
			\int_{t'=0}^{\hat{t}}
				\frac{\iota_2^c(t')} 
				{\upmu_{\star}^{1/2}(t',\hat{u})} 
			\, dt'
				\notag 
				\\
		& \leq \frac{C}{c} (1 + \varsigma^{-1}) 
						H(\hat{t},\hat{u})
			\leq \frac{C}{c} (1 + \varsigma^{-1})  H(t,u).
			\notag
\end{align}
Using a similar argument based on \eqref{E:IF4SMALLNESS},
we bound the term arising from the second term on 
the $16^{th}$ line of RHS~\eqref{E:TOPORDERWAVEENERGYINTEGRALINEQUALITIES}
as follows:
\begin{align} \label{E:BELOWTOPORDERUNITEGRATEDTERMCOUPLEDINTOTOPORDERESTIMATE}
			& C (1 + \varsigma^{-1}) 
			\iota_F^{-1}(\hat{t},\hat{u})
			\int_{u'=1}^{u}
				\totmax{[1,19]}(\hat{t},u')
			\, du'
				\\
			& 
			\leq
			C (1 + \varsigma^{-1}) 
			\upmu_{\star}^{11.8}(\hat{t},\hat{u})
			\iota_1^{-c}(\hat{t})
			\iota_2^{-c}(\hat{t},\hat{u})
			e^{-c \hat{t}} 
			e^{-c \hat{u}}
			\left\lbrace
				\sup_{u' \in [0,\hat{u}]}
				e^{-c u'} \totmax{[1,19]}(\hat{t},u')
			\right\rbrace
			\times
			\int_{u'=0}^{\hat{u}}
				e^{c u'}
			\, du'
					\notag \\
	& 	\leq
			C (1 + \varsigma^{-1}) 
			\upmu_{\star}^2(\hat{t},\hat{u})
			e^{-c \hat{u}}
				\left\lbrace
					\sup_{u' \in [0,\hat{u}]}
					\iota_H^{-1}(\hat{t},\hat{u'}) \totmax{[1,19]}(\hat{t},u')
				\right\rbrace
			\times
			\int_{u'=0}^{\hat{u}}
				e^{c u'}
			\, du'
				\notag \\
		& \leq \frac{C}{c} (1 + \varsigma^{-1}) 
						H(\hat{t},\hat{u})
			\leq 
				\frac{C}{c} (1 + \varsigma^{-1})  
				H(t,u).
			\notag
\end{align}
To bound the terms arising from the three terms on the $17^{th}$ line of
RHS~\eqref{E:TOPORDERWAVEENERGYINTEGRALINEQUALITIES},
 we argue as follows (again recalling that $N=20$):
\begin{align} 
			C \varepsilon
			\iota_F^{-1}(\hat{t},\hat{u})
			\totmax{[1,19]}(\hat{t},\hat{u})
			& 
			= 
			C \varepsilon 
			\upmu_{\star}^2(\hat{t},\hat{u})
			\iota_H^{-1}(\hat{t},\hat{u})
			\coercivespacetimemax{[1,N-1]}(\hat{t},\hat{u})
				\label{E:BELOWTOPORDERNONINTEGRATEDWAVEENERGYCOUPLEDINTOTOPORDERESTIMATE} \\
			& 
			\leq 
			C \varepsilon \upmu_{\star}^2(\hat{t},\hat{u}) H(\hat{t},\hat{u})
			\leq C \varepsilon H(\hat{t},\hat{u})
			\leq C \varepsilon H(t,u),
			\notag
				\\
			C \varsigma 
			\iota_F^{-1}(\hat{t},\hat{u})
			\totmax{[1,19]}(\hat{t},\hat{u})
			& 
			= 
			C \varsigma 
			\upmu_{\star}^2(\hat{t},\hat{u}) 
			\iota_H^{-1}(\hat{t},\hat{u})
			\coercivespacetimemax{[1,N-1]}(\hat{t},\hat{u})
				\label{E:VARSIGMABELOWTOPORDERNONINTEGRATEDWAVEENERGYCOUPLEDINTOTOPORDERESTIMATE} 
				\\
			& 
			\leq 
			C \varsigma \upmu_{\star}^2(\hat{t},\hat{u}) H(\hat{t},\hat{u})
			\leq C \varsigma H(\hat{t},\hat{u})
			\leq C \varsigma H(t,u),
			\notag
				\\
			C \varsigma
			\iota_F^{-1}(\hat{t},\hat{u})
			\coercivespacetimemax{[1,N-1]}(\hat{t},\hat{u})
			& 
			= C \varsigma \upmu_{\star}^2(\hat{t},\hat{u})
			\iota_H^{-1}(\hat{t},\hat{u})
			\coercivespacetimemax{[1,N-1]}(\hat{t},\hat{u})
					\label{E:BELOWTOPORDERCOERCIVESPACETIMETERMCOUPLEDINTOTOPORDERESTIMATE}  \\
			& 
			\leq
			C \varsigma \upmu_{\star}^2(\hat{t},\hat{u})
			H(\hat{t},\hat{u})
		 \leq C \varsigma H(\hat{t},\hat{u})
			\leq 
			C \varsigma H(t,u).
		\notag
	\end{align}
Inserting all of these estimates into the RHS of
$\iota_F^{-1}(\hat{t},\hat{u}) \times$ \eqref{E:TOPORDERWAVEENERGYINTEGRALINEQUALITIES}$(\hat{t},\hat{u})$
and taking $\sup_{(\hat{t},\hat{u}) \in [0,t] \times [0,u]}$ of both sides,
we deduce that
\begin{align} \label{E:NEARLYTOPORDERALMOSTESTIMATED}
		F(t,u) 
		& \leq 
				C 
				(1 + \varsigma^{-1})
				\mathring{\upepsilon}^2 
			+ \left\lbrace
					\frac{6}{11.8} 
					+ \frac{8.1}{5.9 \times 11.8} 
					+ \frac{2}{5.4} 
					+ C \sqrt{\varepsilon} 
					+ C \varsigma
					+ \frac{C}{c} (1 + \varsigma^{-1})
				\right\rbrace
				F(t,u) 
					\\
		& \ \
			+ 
			C
			\left\lbrace
				 \varepsilon
				+ 
				\varsigma
				+
				\frac{1}{c} (1 + \varsigma^{-1})
			\right\rbrace
			H(t,u)
			+ C 
			F^{1/2}(t,u) G^{1/2}(t,u).
			\notag
	\end{align}

We now bound the terms 
$\displaystyle \iota_G^{-1}(\hat{t},\hat{u}) \times \cdots$ arising from the terms
described in \eqref{E:EASYVARIABLESTOPORDERWAVEENERGYINTEGRALINEQUALITIES}.
We claim that the following analog of 
\eqref{E:NEARLYTOPORDERALMOSTESTIMATED}
holds:
	\begin{align} \label{E:EASYNEARLYTOPORDERALMOSTESTIMATED}
		G(t,u) 
		& \leq 
			C (1 + \varsigma^{-1})\mathring{\upepsilon}^2 
			+
			C
			\left\lbrace
				\varepsilon
				+ 
				\varsigma
				+ \frac{1}{c} (1 + \varsigma^{-1})
			\right\rbrace
			F(t,u) 
			+ 
			C
			\left\lbrace
					\varepsilon
				+ \varsigma
				+ \frac{1}{c} (1 + \varsigma^{-1})
			\right\rbrace
			H(t,u).
	\end{align}
	The proof of \eqref{E:EASYNEARLYTOPORDERALMOSTESTIMATED}
	is similar to the proof of \eqref{E:NEARLYTOPORDERALMOSTESTIMATED},
	but with the following key changes:
	\textbf{i)}: in view of \eqref{E:EASYVARIABLESTOPORDERWAVEENERGYINTEGRALINEQUALITIES},
	the terms corresponding to
	\eqref{E:MOSTDANGEROUSINTEGRALESTIMATE}-\eqref{E:THIRDMOSTDANGEROUSINTEGRALESTIMATE}
	and
	\eqref{E:FIRSTLARGECONSTANTTERM}-\eqref{E:THIRDLARGECONSTANTTERM}
	are absent from RHS~\eqref{E:EASYNEARLYTOPORDERALMOSTESTIMATED}.
	To treat the terms
	$\displaystyle \iota_G^{-1}(\hat{t},\hat{u}) \times \cdots$
	corresponding to the three new terms
	explicitly listed in
	\eqref{E:EASYVARIABLESTOPORDERWAVEENERGYINTEGRALINEQUALITIES},
	we argue as in the proof of
	\eqref{E:MOSTDANGEROUSINTEGRALESTIMATE},
	\eqref{E:SECONDMOSTDANGEROUSINTEGRALESTIMATE},
	and
	\eqref{E:THIRDMOSTDANGEROUSINTEGRALESTIMATE},
	but using the cruder inequality \eqref{E:LOSSKEYMUINTEGRALBOUND}
	in place of the delicate inequalities
	\eqref{E:KEYMUTOAPOWERINTEGRALBOUND} 
	and \eqref{E:KEYHYPERSURFACEMUTOAPOWERINTEGRALBOUND}.
	We find that
	\begin{align}
		C \varepsilon
		\iota_G^{-1}(\hat{t},\hat{u})& 
		\int_{t'=0}^{\hat{t}}
				\frac{1} {\upmu_{\star}(t',\hat{t})} 
						\sqrt{\totmax{N}}(t',\hat{t}) 
						\int_{s=0}^{t'}
							\frac{1}{\upmu_{\star}(s,\hat{t})} 
							\sqrt{\totmax{N}}(s,\hat{t}) 
						\, ds
				\, dt'
			\label{E:FIRSTDEGENERATEEASYTOPORDERESTIMATEWEAKCOUPLING} \\
		& \leq
			C \varepsilon  F(t,u),
			\notag \\
		& C \varepsilon
		\iota_G^{-1}(\hat{t},\hat{u})
			\int_{t'=0}^{\hat{t}}
					\frac{1} 
						{\upmu_{\star}(t',\hat{u})} 
				  \totmax{N}(t',\hat{u})
				\, dt'
				\label{E:SECONDDEGENERATEEASYTOPORDERESTIMATEWEAKCOUPLING} \\
		& \leq  C \varepsilon  F(t,u),
			\notag \\
		&
			C \varepsilon
			\iota_G^{-1}(\hat{t},\hat{u})
			\frac{1}{\upmu_{\star}^{1/2}(\hat{t},\hat{u})}
			\sqrt{\totmax{N}}(\hat{t},\hat{u})
			\int_{t'=0}^{\hat{t}}
				\frac{1}{\upmu_{\star}^{1/2}(t',\hat{u})} \sqrt{\totmax{N}}(t',\hat{u})
			\, dt'
				\label{E:THIRDDEGENERATEEASYTOPORDERESTIMATEWEAKCOUPLING} \\
		& \leq C \varepsilon  F(t,u).
		\notag 
	\end{align}
	The desired estimate \eqref{E:EASYNEARLYTOPORDERALMOSTESTIMATED} now follows
	from the same arguments used to prove \eqref{E:NEARLYTOPORDERALMOSTESTIMATED},
	where
	\eqref{E:FIRSTDEGENERATEEASYTOPORDERESTIMATEWEAKCOUPLING}-\eqref{E:THIRDDEGENERATEEASYTOPORDERESTIMATEWEAKCOUPLING}
	contribute to the second product on RHS~\eqref{E:EASYNEARLYTOPORDERALMOSTESTIMATED}.

We now bound the terms 
$\displaystyle \iota_H^{-1}(\hat{t},\hat{u}) \times \cdots$ arising from the terms on
RHS~\eqref{E:BELOWTOPORDERWAVEENERGYINTEGRALINEQUALITIES}.
All terms except the one arising from the
integral involving the top-order factor $\sqrt{\totmax{20}}$
(featured in the $ds$ integral on RHS \eqref{E:BELOWTOPORDERWAVEENERGYINTEGRALINEQUALITIES})
can be bounded by
$
\displaystyle 
\leq C \mathring{\upepsilon}^2 
+ \frac{C}{c}(1 + \varsigma^{-1})G(t,u) 
+ C \varsigma H(t,u)
$
by using essentially the same arguments given above.
In particular, we use the
already proven specific vorticity energy estimates
\eqref{E:VORTMULOSSMAINAPRIORIENERGYESTIMATES}-\eqref{E:VORTNOMULOSSMAINAPRIORIENERGYESTIMATES}
to handle the terms
generated by the integrals on
the last line of RHS~\eqref{E:BELOWTOPORDERWAVEENERGYINTEGRALINEQUALITIES}.
To handle the remaining term involving the top-order factor 
$\sqrt{\totmax{20}}$, 
we use arguments similar to the ones
we used to prove \eqref{E:ESTIMATEFORTRIPLETIMEINTEGRATEDERRORINTEGRAL}
(in particular, we use inequality \eqref{E:LOSSKEYMUINTEGRALBOUND} twice)
to bound it as follows:
\begin{align} \label{E:TOPORDERTERMCOUPLEDINTOBELOWTOPORDERESTIMATE}
			& C
			\iota_H^{-1}(\hat{t},\hat{u})
			\int_{t'=0}^{\hat{t}}
				\frac{1}{\upmu_{\star}^{1/2}(t',\hat{u})} 
						\sqrt{\totmax{[1,19]}}(t',\hat{u}) 
						\int_{s=0}^{t'}
							\frac{1}{\upmu_{\star}^{1/2}(s,\hat{u})} 
							\sqrt{\totmax{20}}(s,\hat{u}) 
						\, ds
				\, dt'
				\\
			& 
			\leq
			C
			\iota_H^{-1}(\hat{t},\hat{u})
			\sup_{(t',u') \in [0,\hat{t}] \times [0,\hat{u}]}
			\left\lbrace 
				\upmu_{\star}^{4.9}(t',\hat{u}) \sqrt{\totmax{[1,19]}}(t',u')
			\right\rbrace
			\times
			\sup_{(t',u') \in [0,\hat{t}] \times [0,\hat{u}]}
			\left\lbrace 
				\upmu_{\star}^{5.9}(t',\hat{u}) \sqrt{\totmax{20}}(t',u')
			\right\rbrace
				\notag \\
			& \times
				\int_{t'=0}^{\hat{t}}
						\frac{1}{\upmu_{\star}^{5.4}(t',\hat{u})} 
						\int_{s=0}^{t'}
							\frac{1}{\upmu_{\star}^{6.4}(s,\hat{u})}
						\, ds
				\, dt'
					\notag
					\\
		& \leq  C
						F^{1/2}(\hat{t},\hat{u})
						H^{1/2}(\hat{t},\hat{u})
			\leq C F(t,u)
				+ \frac{1}{2}H(t,u).
			\notag
\end{align}
Inserting all of these estimates into the RHS of
$\iota_H^{-1}(\hat{t},\hat{u}) \times$ 
\eqref{E:BELOWTOPORDERWAVEENERGYINTEGRALINEQUALITIES}$(\hat{t},\hat{u})$
and taking $\sup_{(\hat{t},\hat{u}) \in [0,t] \times [0,u]}$ of both sides,
we deduce that
\begin{align} \label{E:NEARLYBELOWTOPORDERALMOSTESTIMATED}
	H(t,u) 
		& \leq C \mathring{\upepsilon}^2 
				+ 
				C
				\left\lbrace
					1
					+
					\frac{1}{c}
				\right\rbrace
				F(t,u)
				+ 
				C
				\left\lbrace
					\frac{1}{2}
					+
					\varsigma
					+
					(1 + \varsigma^{-1}) \frac{1}{c}
				\right\rbrace
				H(t,u).
\end{align}

We now consider the system of three inequalities
\eqref{E:NEARLYTOPORDERALMOSTESTIMATED},
\eqref{E:EASYNEARLYTOPORDERALMOSTESTIMATED},
and
\eqref{E:NEARLYBELOWTOPORDERALMOSTESTIMATED},
and we remind the reader that the constants $C$ 
in these inequalities can be chosen to be independent of $c$.
The desired estimates \eqref{E:TOPORDERALMOSTESTIMATED}-\eqref{E:BELOWTOPORDERALMOSTESTIMATED}
now follow from first choosing $\varsigma$ to be sufficiently small,
then choosing $c$ to be sufficiently large, then choosing
$\varepsilon$ to be sufficiently small,
and using the aforementioned fact that
$\displaystyle 
\frac{6}{11.8} 
+ \frac{8.1}{5.9 \times 11.8} 
+ \frac{2}{5.4} 
+ C \sqrt{\varepsilon} 
< 1
$.


\noindent \textbf{Estimates for} 
$
\max
\left\lbrace
	\totmax{[1,18]},
	\coercivespacetimemax{[1,18]}
\right\rbrace
$,
$
\max
\left\lbrace
	\totmax{[1,17]},
	\coercivespacetimemax{[1,17]}
\right\rbrace
$,
$\cdots$,
$
\max
\left\lbrace
	\totmax{1},
	\coercivespacetimemax{1}
\right\rbrace
$ 
\textbf{via a descent scheme}:
We now explain how to use inequality
\eqref{E:BELOWTOPORDERWAVEENERGYINTEGRALINEQUALITIES}
to derive the estimates
for 
$
\max
\left\lbrace
	\totmax{[1,18]},
	\coercivespacetimemax{[1,18]}
\right\rbrace
$,
$\max
\left\lbrace
	\totmax{[1,17]},
	\coercivespacetimemax{[1,17]}
\right\rbrace$,
$\cdots$,
$
\max
\left\lbrace
	\totmax{1},
	\coercivespacetimemax{1}
\right\rbrace
$ 
by downward induction.
Unlike our analysis of the strongly coupled pair
$
\max
\left\lbrace
	\totmax{20},
	\coercivespacetimemax{20}
\right\rbrace
$
and 
$
\max
\left\lbrace
	\totmax{[1,19]},
	\coercivespacetimemax{[1,19]}
\right\rbrace
$,
we can derive the desired estimates
for 
$
\max
\left\lbrace
	\totmax{[1,18]},
	\coercivespacetimemax{[1,18]}
\right\rbrace
$
by using only inequality \eqref{E:BELOWTOPORDERWAVEENERGYINTEGRALINEQUALITIES}
and the already derived estimates
for 
$
\max
\left\lbrace
	\totmax{[1,19]},
	\coercivespacetimemax{[1,19]}
\right\rbrace
$.
At the end of the proof, we will describe
the minor changes needed to derive
the desired estimates for
$
\max
\left\lbrace
	\totmax{[1,17]},
	\coercivespacetimemax{[1,17]}
\right\rbrace
$, 
$\cdots$,
$
\max
\left\lbrace
\totmax{1} + \coercivespacetimemax{1}
\right\rbrace
$.

To begin, we define the following analogs of 
\eqref{E:JUSTBELOWTOPORDERTOTALINTEGRATINGFACTOR}
and
\eqref{E:JUSTBELOWTOPORDERENERGYRESCALED}:
\begin{align} 	\label{E:MORETHANONEBELOWTOPORDERTOTALINTEGRATINGFACTOR}
	\iota_{\widetilde{H}}(t',u')
	& :=
		\upmu_{\star}^{-7.8}(t',u')
		\iota_1^c(t')
		\iota_2^c(t',u')
		e^{c t'}
		e^{c u'},
\end{align}
\begin{align}
		\widetilde{H}(t,u) 
		& 
		:= 
		\sup_{(\hat{t},\hat{u}) \in [0,t] \times [0,u]} 
		\iota_{\widetilde{H}}^{-1}(\hat{t},\hat{u})
		\max
		\left\lbrace
			\totmax{[1,18]}(\hat{t},\hat{u}),
			\coercivespacetimemax{[1,18]}(\hat{t},\hat{u})
		\right\rbrace.
		\label{E:MORETHANONEBELOWTOPORDERENERGYRESCALED}
	\end{align}
Note that the power of $\upmu_{\star}^{-1}$
in the factor $\upmu_{\star}^{-7.8}$ has been reduced
by two in \eqref{E:MORETHANONEBELOWTOPORDERTOTALINTEGRATINGFACTOR}
compared to \eqref{E:JUSTBELOWTOPORDERTOTALINTEGRATINGFACTOR},
which corresponds to less singular behavior of
$
\max
\left\lbrace
	\totmax{[1,18]},
	\coercivespacetimemax{[1,18]}
\right\rbrace
$
near the shock. As before, to prove the desired estimate
\eqref{E:MULOSSMAINAPRIORIENERGYESTIMATES} (now with $K=3$),
it suffices to prove
\begin{align} \label{E:MORETHANONEBELOWTOPORDERDESIREDBOUND}
	\widetilde{H}(t,u)
	\leq C \mathring{\upepsilon}^2.
\end{align}

We now set $N=19$, 
multiply both sides of inequality
\eqref{E:BELOWTOPORDERWAVEENERGYINTEGRALINEQUALITIES}
by $\iota_{\widetilde{H}}^{-1}(t,u)$ and then set $(t,u) = (\hat{t},\hat{u})$.
With one exception,
we can bound
all terms arising from the integrals on RHS \eqref{E:BELOWTOPORDERWAVEENERGYINTEGRALINEQUALITIES}
by 
$\displaystyle 
\leq C \mathring{\upepsilon}^2 
+ \frac{C}{c}(1 + \varsigma^{-1})\widetilde{H}(t,u) 
+ \varsigma \widetilde{H}(t,u)
$
(where $C$ is independent of $c$)
by using the same arguments
that we used in deriving the estimate for
$
\max
\left\lbrace
	\totmax{[1,19]},
	\coercivespacetimemax{[1,19]}
\right\rbrace
$.
The exceptional term is the one arising from the 
integral involving the above-present-order factor $\sqrt{\totmax{[1,19]}}$.
We bound the exceptional term as follows by using
inequality \eqref{E:LOSSKEYMUINTEGRALBOUND}, 
the approximate monotonicity of $\iota_{\widetilde{H}}$,
and the estimate
$
\sqrt{\totmax{[1,19]}} \leq C_c \mathring{\upepsilon} \upmu_{\star}^{-4.9}(t,u)
$
(which follows from the already proven estimate \eqref{E:DESIREDRESCALEDAPRIORIBOUND} for $H(t,u)$):
\begin{align} \label{E:BELOWTOPORDERTERMCOUPLEDINTOMOREBELOWTOPORDERESTIMATE}
			& C
			\iota_{\widetilde{H}}^{-1}(\hat{t},\hat{u})
			\int_{t'=0}^{\hat{t}}
				\frac{1}{\upmu_{\star}^{1/2}(t',\hat{u})} 
						\sqrt{\totmax{[1,18]}}(t',\hat{u}) 
						\int_{s=0}^{t'}
							\frac{1}{\upmu_{\star}^{1/2}(s,\hat{u})} 
							\sqrt{\totmax{[1,19]}}(s,\hat{u}) 
						\, ds
				\, dt'
				\\
			& 
			\leq
			C_c
			\mathring{\upepsilon}
			\iota_{\widetilde{H}}^{-1/2}(\hat{t},\hat{u})
			\sup_{(t',u') \in [0,\hat{t}] \times [0,\hat{u}]}
			\left\lbrace 
				\iota_{\widetilde{H}}^{-1/2}(t',u')
				\sqrt{\totmax{[1,18]}}(t',u')
			\right\rbrace
			\times
			\int_{t'=0}^{\hat{t}}
						\frac{1}{\upmu_{\star}^{1/2}(t',\hat{u})} 
						\int_{s=0}^{t'}
							\frac{1}{\upmu_{\star}^{5.4}(s,\hat{u})}
						\, ds
				\, dt'
					\notag
					\\
		& \leq 
			C_c
			\mathring{\upepsilon}
			\iota_{\widetilde{H}}^{-1/2}(\hat{t},\hat{u})
			\upmu_{\star}^{-3.9}(\hat{t},\hat{u})
			\sup_{(t',u') \in [0,\hat{t}] \times [0,\hat{u}]}
			\left\lbrace 
				\iota_{\widetilde{H}}^{-1/2}(t',u')
				\sqrt{\totmax{[1,18]}}(t',u')
			\right\rbrace
			\notag
				\\
		& \leq 
			C_c
			\mathring{\upepsilon}
			\widetilde{H}^{1/2}(\hat{t},\hat{u})
			\leq 
			C_c \mathring{\upepsilon}^2
			+
			\frac{1}{2}
			\widetilde{H}(t,u).
			\notag
\end{align}
In total, we have obtained the following
analog of \eqref{E:NEARLYBELOWTOPORDERALMOSTESTIMATED}:
\begin{align}
	\widetilde{H}(t,u) \label{E:NEARLYMORETHANONEBELOWTOPORDERALMOSTESTIMATED}
		& \leq C_c \mathring{\upepsilon}^2 
				+ \frac{C}{c} (1 + \varsigma^{-1}) \widetilde{H}(t,u)
				+ \frac{1}{2} \widetilde{H}(t,u)
				+ C \varsigma \widetilde{H}(t,u),
\end{align}
where $C_c$ is the only constant that depends on $c$.
The desired bound \eqref{E:MORETHANONEBELOWTOPORDERDESIREDBOUND}
easily follows from \eqref{E:NEARLYMORETHANONEBELOWTOPORDERALMOSTESTIMATED}
by first choosing $\varsigma$ to be sufficiently small
and then $c$ to be sufficiently large
so that we can absorb all factors of $\widetilde{H}$ 
on RHS \eqref{E:NEARLYMORETHANONEBELOWTOPORDERALMOSTESTIMATED}
into the LHS.

The desired bounds \eqref{E:NOMULOSSMAINAPRIORIENERGYESTIMATES} for
$
\max
\left\lbrace
	\totmax{[1,17]},
	\coercivespacetimemax{[1,17]}
\right\rbrace
$,
$
\max
\left\lbrace
	\totmax{[1,16]},
	\coercivespacetimemax{[1,16]}
\right\rbrace
$
$\cdots$
can be (downward) inductively derived
by using an argument similar to the one
we used to bound
$
\max
\left\lbrace
	\totmax{[1,18]},
	\coercivespacetimemax{[1,18]}
\right\rbrace
$,
which relied on the already proven bounds for 
$
\max
\left\lbrace
	\totmax{[1,19]},
	\coercivespacetimemax{[1,19]}
\right\rbrace
$.
The only difference is that
we define the analog of the approximating integrating factor 
\eqref{E:MORETHANONEBELOWTOPORDERTOTALINTEGRATINGFACTOR}
to be
$ \upmu_{\star}^{-P}
	\iota_1^c(t')
	\iota_2^c(t',u')
	e^{c t'}
	e^{c u'}
$,
where $P = 5.8$ for the 
$
\max
\left\lbrace
	\totmax{[1,17]},
	\coercivespacetimemax{[1,17]}
\right\rbrace
$ 
estimate,
$P = 3.8$ for the 
$
\max
\left\lbrace
	\totmax{[1,16]},
	\coercivespacetimemax{[1,16]}
\right\rbrace
$ 
estimate,
$P = 1.8$ for the 
$
\max
\left\lbrace
	\totmax{[1,15]},
	\coercivespacetimemax{[1,15]}
\right\rbrace
$ 
estimate,
and $P=0$ for the  
$
\max
\left\lbrace
	\totmax{[1,\leq 14]},
	\coercivespacetimemax{[1,\leq 14]}
\right\rbrace
$
estimates;
these latter estimates \emph{do not involve any singular factor} of $\upmu_{\star}^{-1}$.
There is one important new detail relevant
for these estimates: in deriving the analog of the inequalities
\eqref{E:BELOWTOPORDERTERMCOUPLEDINTOMOREBELOWTOPORDERESTIMATE}
for 
$
\max
\left\lbrace
	\totmax{[1,\leq 14]},
	\coercivespacetimemax{[1,\leq 14]}
\right\rbrace
$,
we use the estimate \eqref{E:LESSSINGULARTERMSMPOINTNINEINTEGRALBOUND}
in place of the estimate \eqref{E:LOSSKEYMUINTEGRALBOUND};
as in our proof of Lemma~\ref{L:VORTICITYAPRIORIENERGYESTIMATES},
the estimate \eqref{E:LESSSINGULARTERMSMPOINTNINEINTEGRALBOUND}
allows us to break the $\upmu_{\star}^{-1}$ degeneracy.
This completes the proof of Prop.~\ref{P:MAINAPRIORIENERGY}.

\section{The Main Theorem}
\label{S:STABLESHOCKFORMATION}
We now state and prove the main theorem.

\begin{theorem}[\textbf{Stable shock formation}]
\label{T:MAINTHEOREM}
Let $(\Densrenormalized, v^1, v^2, \Vortrenormalized)$
be a solution to the $2D$ compressible Euler equations 
in the form 
\eqref{E:VELOCITYWAVEEQUATION}-\eqref{E:RENORMALIZEDVORTICTITYTRANSPORTEQUATION}
under any physical\footnote{Physical in the sense described below equation \eqref{E:BARATROPICEOS}.} 
barotropic equation of state
except for that of a Chaplygin gas (see \eqref{E:EOSCHAPLYGINGAS})
and let $u$ be a solution to the eikonal equation \eqref{E:INTROEIKONAL}.
Let $\threePsi = (\Psi_0,\Psi_1,\Psi_2):= (\Densrenormalized , v^1,v^2)$
denote the array of the difference between the wave variables and the constant state
solution $(0,0,0)$.
Assume that the solution verifies the size assumptions
on $\Sigma_0^1$ and $\mathcal{P}_0^{2 \TranminusdatasizeWithFactor^{-1}}$
stated in Sects.~\ref{SS:FLUIDVARIABLEDATAASSUMPTIONS} and \ref{SS:DATAFOREIKONALFUNCTIONQUANTITIES}
as well as the smallness assumptions of Sect.~\ref{SS:SMALLNESSASSUMPTIONS}.
In particular, let
$\mathring{\upepsilon}$,
$\mathring{\updelta}$,
and $\TranminusdatasizeWithFactor$
be the data-size parameters from
\eqref{E:CRITICALBLOWUPTIMEFACTOR},
\eqref{E:L2SMALLDATAASSUMPTIONSALONGSIGMA0}-\eqref{E:SMALLDATAASSUMPTIONSALONGELLT0},
and
\eqref{E:LUNITISMALLDATA}-\eqref{E:SIMPLEUPMUALONGPOESTIMATE}.
Assume the genericity condition\footnote{For any barotropic equation of state
except for that of the Chaplygin gas (see \eqref{E:EOSCHAPLYGINGAS}),
there exist choices of the background density $\bar{\rho}$ 
such that the condition \eqref{E:NONDEGENERACYCONDITION} holds.}
\begin{align} \label{E:NONDEGENERACYCONDITION}
	\bar{\Speed}' + 1 
	\neq 0,
\end{align}
where $\bar{\Speed}' := \Speed'(\Densrenormalized = 0)$
denotes the value of $\Speed'$ corresponding to the background constant state.
Let $\Upsilon$ be the change of variables map
from geometric to Cartesian coordinates (see Def.~\ref{D:CHOVMAP}).
For each $U_0 \in [0,1]$, let 
\begin{align*}
	& T_{(Lifespan);U_0}
		\\
	& := 
	\sup 
	\Big\lbrace 
		t \in [1,\infty) \ | \ \mbox{the solution exists classically on } \mathcal{M}_{t;U_0}
			\\
	& \ \ \ \ \ \ \ \ \ \ \ \ \ \ \ \ \ \ \ \ \ \ \ \
		\mbox{ and $\Upsilon$ is a diffeomorphism from } [1,t) \times [0,U_0] \times \mathbb{T} 
		\mbox{ onto its image}
		\Big\rbrace
\end{align*}
(see Figure~\ref{F:SOLIDREGION} on pg.~\pageref{F:SOLIDREGION}).
If $\mathring{\upepsilon}$ is sufficiently small\footnote{Recall that in Subsect.~\ref{SS:EXISTENCEOFDATA}, 
we show that there exists an open set of solutions satisfying the desired
smallness conditions.}
relative to 
$\mathring{\updelta}^{-1}$
and 
$\TranminusdatasizeWithFactor$
(in the sense explained in Sect.~\ref{SS:SMALLNESSASSUMPTIONS}),
then the following conclusions hold,
where all constants can be chosen to be independent of $U_0$.

\noindent \underline{\textbf{Dichotomy of possibilities.}}
One of the following mutually disjoint possibilities must occur,
where $\upmu_{\star}(t,u)$ is defined in \eqref{E:MUSTARDEF}.
\begin{enumerate}
	\renewcommand{\labelenumi}{\textbf{\Roman{enumi})}}
	\item $T_{(Lifespan);U_0} > 2 \TranminusdatasizeWithFactor^{-1}$. 
		In particular, the solution exists classically on the spacetime
		region $\mbox{\upshape cl} \mathcal{M}_{2 \TranminusdatasizeWithFactor^{-1},U_0}$,
		where $\mbox{\upshape cl}$ denotes closure.
		Furthermore, $\inf \lbrace \upmu_{\star}(s,U_0) \ | \ s \in [0,2 \TranminusdatasizeWithFactor^{-1}] \rbrace > 0$.
	\item $0 <  T_{(Lifespan);U_0} \leq 2 \TranminusdatasizeWithFactor^{-1}$,
		and 
		\begin{align} \label{E:MAINTHEOREMLIFESPANCRITERION}
		T_{(Lifespan);U_0} 
		= \sup 
			\left\lbrace 
			t \in [1, 2 \TranminusdatasizeWithFactor^{-1}) \ | \
				\inf \lbrace \upmu_{\star}(s,U_0) \ | \ s \in [1,t) \rbrace > 0
			\right\rbrace.
		\end{align}
\end{enumerate}
In addition, case $\textbf{II)}$ occurs when $U_0 = 1$. In this case, we have
\begin{align} \label{E:CLASSICALLIFESPANASYMPTOTICESTIMATE}
	T_{(Lifespan);1} 
	= 
	\left\lbrace 
		1 + \mathcal{O}(\mathring{\upepsilon}) 
	\right\rbrace
	\TranminusdatasizeWithFactor^{-1}.
\end{align}

\noindent \underline{\textbf{What happens in Case I).}}
In case $\textbf{I)}$, 
all bootstrap assumptions,
the estimates of Props.~\ref{P:IMPROVEMENTOFAUX} and \ref{P:IMPROVEMENTOFHIGHERTRANSVERSALBOOTSTRAP},
and the energy estimates of Prop.~\ref{P:MAINAPRIORIENERGY}
hold on $\mbox{\upshape cl} \mathcal{M}_{2 \TranminusdatasizeWithFactor^{-1},U_0}$
with all factors $\varepsilon$ on the RHS of all inequalities 
replaced by $C \mathring{\upepsilon}$.
Moreover, for $0 \leq K \leq 5$, the following estimates hold
for $(t,u) \in [0,2 \TranminusdatasizeWithFactor^{-1}] \times [0,U_0]$:
\begin{subequations}
	\begin{align}
		\left\|
			\Fullset_{\ast}^{[1,14];\leq 1} \upmu
		\right\|_{L^2(\Sigma_t^u)},
			\,
		\left\|
			\Fullset_{\ast}^{\leq 14;\leq 2}
			\Lunit_{(Small)}^i
		\right\|_{L^2(\Sigma_t^u)},
			\,
		\left\|
			\Fullset_{\ast}^{\leq 13;\leq 2}
			\mytr \upchi
		\right\|_{L^2(\Sigma_t^u)}
		& \leq 
			C \mathring{\upepsilon},
				 \label{E:LOWORDERTANGENTIALEIKONALL2MAINTHEOREM}
				 \\
		\left\|
			\Fullset_{\ast}^{15 + K;\leq 1} \upmu
		\right\|_{L^2(\Sigma_t^u)},
			\,
		\left\|
			\Fullset_{\ast}^{15 + K;\leq 2} 
			\Lunit_{(Small)}^i
		\right\|_{L^2(\Sigma_t^u)},
			\,
		\left\|
			\Fullset_{\ast}^{14 + K;\leq 2} 
			\mytr \upchi
		\right\|_{L^2(\Sigma_t^u)}
		& \leq 
			C \mathring{\upepsilon} \upmu_{\star}^{-(K+.4)}(t,u),
				 \label{E:MIDORDERTANGENTIALEIKONALL2MAINTHEOREM} 
				 	\\
		\left\|
			\Lunit \Fullset_{\ast}^{20; \leq 1} \upmu
		\right\|_{L^2(\Sigma_t^u)},
			\,
		\left\|
			\Lunit \Fullset_{\ast}^{20; \leq 2} \Lunit_{(Small)}^i
		\right\|_{L^2(\Sigma_t^u)},
			\,
		\left\|
			\Lunit \Fullset_{\ast}^{19; \leq 2} \mytr \upchi
		\right\|_{L^2(\Sigma_t^u)}
		& \leq 
			C \mathring{\upepsilon} \upmu_{\star}^{-6.4}(t,u),
				 \label{E:HIGHORDERLUNITTANGENTIALEIKONALL2MAINTHEOREM} 
		\\
		\left\|
			\upmu \GeoAng^{20} \mytr \upchi
		\right\|_{L^2(\Sigma_t^u)},
			\,
		\left\|
			\upmu \GeoAng^{19} \Rad \mytr \upchi
		\right\|_{L^2(\Sigma_t^u)}
		& \leq 
			C \mathring{\upepsilon} \upmu_{\star}^{-5.9}(t,u).
				 \label{E:HIGHORDERANGULARTRCHIL2MAINTHEOREM} 
	\end{align}
\end{subequations}

\noindent \underline{\textbf{What happens in Case II).}}
In case $\textbf{II)}$, 
all bootstrap assumptions,
the estimates of Props.~\ref{P:IMPROVEMENTOFAUX} and \ref{P:IMPROVEMENTOFHIGHERTRANSVERSALBOOTSTRAP},
and the energy estimates of Prop.~\ref{P:MAINAPRIORIENERGY}
hold on $\mathcal{M}_{T_{(Lifespan);U_0},U_0}$
with all factors $\varepsilon$ on the RHS of all inequalities 
replaced by $C \mathring{\upepsilon}$.
Moreover, for $0 \leq K \leq 5$, the estimates 
\eqref{E:LOWORDERTANGENTIALEIKONALL2MAINTHEOREM}-\eqref{E:HIGHORDERANGULARTRCHIL2MAINTHEOREM}
hold for $(t,u) \in [1,T_{(Lifespan);U_0}) \times [0,U_0]$.
In addition, the scalar-valued functions
$\Fullset^{\leq 11;\leq 2} \threePsi$,
$\Rad \Rad \Rad \threePsi$,
$\Fullset^{\leq 10;\leq 2} \Lunit^i$,
$\Fullset^{\leq 11;\leq 1} \upmu$,
and
$\Rad \Rad \upmu$
extend to 
$\Sigma_{T_{(Lifespan);U_0}}^{U_0}$ 
as functions of 
the geometric coordinates $(t,u,\vartheta)$ 
that are uniformly bounded in $L^{\infty}$.
Furthermore, the Cartesian component functions
$g_{\alpha \beta}(\threePsi)$ 
verify the estimate
$g_{\alpha \beta} = m_{\alpha \beta} + \mathcal{O}(\mathring{\upepsilon})$ 
(where $m_{\alpha \beta} = \mbox{\upshape diag}(-1,1,1)$ is the standard Minkowski metric)
and have the same extension properties as $\threePsi$ and its derivatives with
respect to the vectorfields mentioned above.

Moreover,	let $\Sigma_{T_{(Lifespan);U_0}}^{U_0;(Blowup)}$
be the (non-empty) subset of $\Sigma_{T_{(Lifespan);U_0}}^{U_0}$ 
defined by
\begin{align} \label{E:BLOWUPPOINTS}
	\Sigma_{T_{(Lifespan);U_0}}^{U_0;(Blowup)}
	:= \left\lbrace
			(T_{(Lifespan);U_0},u,\vartheta)
			\ | \
			\upmu(T_{(Lifespan);U_0},u,\vartheta)
			= 0
		\right\rbrace.
\end{align}
Then for each point $(T_{(Lifespan);U_0},u,\vartheta) \in \Sigma_{T_{(Lifespan);U_0}}^{U_0;(Blowup)}$,
there exists a past neighborhood containing it such that the following lower bound holds in
the neighborhood:
\begin{align} \label{E:BLOWUPPOINTINFINITE}
	\left| \Radunit \Densrenormalized (t,u,\vartheta) \right|,
		\,
	\left| \Radunit v^1 (t,u,\vartheta) \right|
	\geq \frac{\TranminusdatasizeWithFactor}
	{8 \left|\bar{\Speed}' + 1 \right|} 
	\frac{1}{\upmu(t,u,\vartheta)}.
\end{align}
In \eqref{E:BLOWUPPOINTINFINITE}, 
$
\displaystyle
\frac{\TranminusdatasizeWithFactor}
{8 \left|\bar{\Speed}' + 1 \right|}
$
is a \textbf{positive} data-dependent constant
(see \eqref{E:NONDEGENERACYCONDITION}),
and the $\ell_{t,u}$-transversal vectorfield $\Radunit$ is near-Euclidean-unit length:
$\delta_{ab} \Radunit^a \Radunit^b = 1 + \mathcal{O}(\mathring{\upepsilon})$.
In particular, 
$\Radunit \Densrenormalized$
and
$\Radunit v^1$ blow up like $1/\upmu$ at all points in $\Sigma_{T_{(Lifespan);U_0}}^{U_0;(Blowup)}$.
Conversely, at all points in
$(T_{(Lifespan);U_0},u,\vartheta) \in \Sigma_{T_{(Lifespan);U_0}}^{U_0} \backslash \Sigma_{T_{(Lifespan);U_0}}^{U_0;(Blowup)}$,
we have
\begin{align} \label{E:NONBLOWUPPOINTBOUND}
	\left| \Radunit \Densrenormalized (T_{(Lifespan);U_0},u,\vartheta) \right|,
		\,
	\left| \Radunit v^1 (T_{(Lifespan);U_0},u,\vartheta) \right|
	< \infty.
\end{align}

\end{theorem}

\begin{proof}
	Let $\BigConst > 1$ be a constant 
	(we will enlarge it as needed throughout the proof).
	We define
	\begin{align}
		T_{(Max);U_0} &:= \mbox{ The supremum of the set of times } 
		\Tboot \in [0,2 \TranminusdatasizeWithFactor^{-1}] \mbox{ such that:} 
			\label{E:LIFESPANPROOF}	\\
			& \bullet \mbox{$\threePsi$, $\Vortrenormalized$, $u$, $\upmu$, $\Lunit_{(Small)}^i$,
				and all of the other quantities} 
					\notag \\
			& \ \ \mbox{defined throughout the article
				  exist classically on } \mathcal{M}_{\Tboot,U_0}.
					\notag \\
		& \bullet 
			\mbox{The change of variables map $\Upsilon$ from Def.~\ref{D:CHOVMAP} 
			is a (global) $C^{1,1}$ }
				\notag \\
		& \ \ 
			\mbox{diffeomorphism from $[0,\Tboot) \times [0,U_0] \times \mathbb{T}$
			onto its image $\mathcal{M}_{\Tboot,U_0}$}.
				\notag \\
		& \bullet \inf \left\lbrace \upmu_{\star}(t,U_0) \ | \ t \in [0,\Tboot) \right\rbrace > 0.
				\notag \\
		& \bullet \mbox{The fundamental } L^{\infty} \mbox{ bootstrap assumptions } 
						\eqref{E:PSIFUNDAMENTALC0BOUNDBOOTSTRAP} \mbox{ and } \eqref{E:VORTFUNDAMENTALC0BOUNDBOOTSTRAP}
						\notag \\
		& \ \ \mbox{ hold with } \varepsilon := \BigConst \mathring{\upepsilon}
								\mbox{ for  } (t,u) \in \times [0,\Tboot) \times [0,U_0].
					\notag \\
			& \bullet \mbox{The following $L^2$-type energy bounds hold for }
				(t,u) \in \times [0,\Tboot) \times [0,U_0]:
				\notag \\
		& \ \ 
		\totmax{15+K}^{1/2}(t,u)
		+ 
		\coercivespacetimemax{15+K}^{1/2}(t,u)
		\leq \BigConst \mathring{\upepsilon} \upmu_{\star}^{-(K+.9)}(t,u),
			\qquad (0 \leq K \leq 5),
				\label{E:PROOFMULOSSMAINAPRIORIENERGYESTIMATES} \\
		& \ \ 
		\totmax{[1,14]}^{1/2}(t,u)
		+ 
		\coercivespacetimemax{[1,14]}^{1/2}(t,u)
		\leq \BigConst \mathring{\upepsilon},
		\label{E:PROOFNOMULOSSMAINAPRIORIENERGYESTIMATES}
								\\
		& \ \ 
		\Vorttotmax{21}^{1/2}(t,u)
		\leq \BigConst \mathring{\upepsilon} \upmu_{\star}^{-6.4}(t,u),
			\label{E:VORTICITYPROOFMULOSSMAINAPRIORIENERGYESTIMATES} \\
		& \ \ 
		\Vorttotmax{16+K}^{1/2}(t,u)
		\leq \BigConst \mathring{\upepsilon} \upmu_{\star}^{-(K+.9)}(t,u),
			\qquad (0 \leq K \leq 4),
				\label{E:SECONDVORTICITYPROOFMULOSSMAINAPRIORIENERGYESTIMATES} \\
		& \ \ 
		\Vorttotmax{\leq 15}^{1/2}(t,u)
		\leq \BigConst \mathring{\upepsilon}.
		\label{E:VORTICITYPROOFNOMULOSSMAINAPRIORIENERGYESTIMATES}
		\end{align}
	It is a standard result that
	if $\mathring{\upepsilon}$ is sufficiently small and $\BigConst$ 
	is sufficiently large, then $T_{(Max);U_0}> 0$
	(this is a standard local well-posedness result
	combined with the initial smallness of the
	$L^2$-controlling quantities obtained in Lemma~\ref{L:INITIALSIZEOFL2CONTROLLING}).

	We now show that the energy bounds
	\eqref{E:PROOFMULOSSMAINAPRIORIENERGYESTIMATES}-\eqref{E:VORTICITYPROOFNOMULOSSMAINAPRIORIENERGYESTIMATES}
	and the fundamental $L^{\infty}$ bootstrap assumption 
	\eqref{E:PSIFUNDAMENTALC0BOUNDBOOTSTRAP} and \eqref{E:VORTFUNDAMENTALC0BOUNDBOOTSTRAP}
	are not saturated for $(t,u) \in [1,T_{(Max);U_0}) \times [0,U_0]$.
	The non-saturation of the energy bounds (for $\BigConst$ sufficiently large)
	is provided by Prop.~\ref{P:MAINAPRIORIENERGY}.
	The non-saturation of the fundamental $L^{\infty}$ bootstrap assumptions
	\eqref{E:PSIFUNDAMENTALC0BOUNDBOOTSTRAP} then follows from
	Cor.~\ref{C:PSILINFTYINTERMSOFENERGIES}.
	Consequently, we conclude that all of the estimates proved throughout the article
	hold on $\mathcal{M}_{\Tboot,U_0}$ with
	the smallness parameter $\varepsilon$ replaced by $C \mathring{\upepsilon}$.
	We use this fact throughout the remainder of the proof without further remark. 

	Next, we show that \eqref{E:LOWORDERTANGENTIALEIKONALL2MAINTHEOREM}-\eqref{E:HIGHORDERANGULARTRCHIL2MAINTHEOREM}
	hold for $(t,u) \in [1,T_{(Max);U_0}) \times [0,U_0]$.
	To obtain \eqref{E:LOWORDERTANGENTIALEIKONALL2MAINTHEOREM}-\eqref{E:HIGHORDERLUNITTANGENTIALEIKONALL2MAINTHEOREM},
	we insert the energy estimates of Prop.~\ref{P:MAINAPRIORIENERGY}
	into the RHS of the inequalities of Lemma~\ref{L:EASYL2BOUNDSFOREIKONALFUNCTIONQUANTITIES}
	and use inequalities \eqref{E:LOSSKEYMUINTEGRALBOUND} and \eqref{E:LESSSINGULARTERMSMPOINTNINEINTEGRALBOUND}
	as well as the fact that $\totmax{[1,M]}$ is increasing in its arguments.
	To obtain inequality \eqref{E:HIGHORDERANGULARTRCHIL2MAINTHEOREM}, 
	we also insert the energy estimates of Prop.~\ref{P:MAINAPRIORIENERGY}
	into RHS~\eqref{E:LESSPRECISEDIFFICULTTERML2BOUND}
	and use inequality
	\eqref{E:LOSSKEYMUINTEGRALBOUND}.

	We now establish the dichotomy of possibilities.
	We first show that
	if 
	\[
		\inf \left\lbrace \upmu_{\star}(t,U_0) \ | \ t \in [1,T_{(Max);U_0}) \right\rbrace > 0,
	\]
	then $T_{(Max);U_0} = 2 \TranminusdatasizeWithFactor^{-1}$.
	This fact can be established 
	using the same arguments given
	in the proof of \cite{jSgHjLwW2016}*{Theorem 15.1}
	(for $\mathring{\upepsilon}$ sufficiently small),
	which were based on analogs of the fundamental $L^{\infty}$ bootstrap assumptions
	(now known to be non-saturated)
	and the $L^{\infty}$ estimates of Props.~\ref{P:IMPROVEMENTOFAUX} and 
	\ref{P:IMPROVEMENTOFHIGHERTRANSVERSALBOOTSTRAP}.
	We will not repeat the (straightforward but tedious) proof here; we note only that
	the above assumption for $\upmu_{\star}$ can be combined
	with other simple estimates to yield that
	$\Upsilon$ extends as a
	global $C^{1,1}$ diffeomorphism from
	$[1,T_{(Max);U_0}] \times [0,U_0] \times \mathbb{T}$
	onto its image and moreover, that
	neither the solution nor its derivatives can blow up
	with respect to geometric or Cartesian coordinates for
	times in $[0,2 \TranminusdatasizeWithFactor^{-1}]$.
	We have thus shown that \textbf{I)} $T_{(Max);U_0} = 2 \TranminusdatasizeWithFactor^{-1}$
	or \textbf{II)} $\inf \left\lbrace \upmu_{\star}(t,U_0) \ | \ t \in [1,T_{(Max);U_0}) \right\rbrace = 0$.

	We now show that case \textbf{II)} corresponds to a singularity
	and that the classical lifespan is characterized by \eqref{E:MAINTHEOREMLIFESPANCRITERION}.
	To this end, we first use
	\eqref{E:CRUCIALPSITRANSVERSALLINFINITYBOUNDBOOTSTRAPIMPROVEDSMALL},
	\eqref{E:LUNITUPMUDOESNOTDEVIATEMUCHFROMTHEDATA},
	\eqref{E:SMALLMUIMPLIESLMUISNEGATIVE},
	and the identity $\Rad = \upmu \Radunit$
	to deduce that inequality \eqref{E:BLOWUPPOINTINFINITE} holds.
	Furthermore, from 
	\eqref{E:LITTLEGDECOMPOSED}-\eqref{E:METRICPERTURBATIONFUNCTION},
	\eqref{E:LINEARLYSMALLSCALARSDEPENDINGONGOODVARIABLES},
	and the $L^{\infty}$ estimates of Prop.~\ref{P:IMPROVEMENTOFAUX},
	we deduce that
	$|\Radunit| 
	: =
	\sqrt{g_{ab} \Radunit^a \Radunit^b} 
	= 1 + \smoothfunction(\GdVar) \GdVar 
	=  
	1 
	+ \mathcal{O}(\mathring{\upepsilon})
	$.
	From this estimate and \eqref{E:BLOWUPPOINTINFINITE}, 
	we deduce that at points in $\Sigma_{T_{(Max);U_0},U_0}$
	where $\upmu$ vanishes, $|\Radunit \Psi|$ must blow up like $1/\upmu$.
	Hence, $T_{(Max);U_0}$ is the classical lifespan. That is, 
	we have $T_{(Max);U_0} = T_{(Lifespan);U_0}$
	as well as the characterization \eqref{E:MAINTHEOREMLIFESPANCRITERION} of the 
	classical lifespan. The estimate \eqref{E:NONBLOWUPPOINTBOUND}
	is an immediate consequence of
	the estimates
	\eqref{E:DENSITYPSITRANSVERSALLINFINITYBOUNDBOOTSTRAPIMPROVEDLARGE}-\eqref{E:PSITRANSVERSALLINFINITYBOUNDBOOTSTRAPIMPROVEDLARGE}
	and the identity $\Rad = \upmu \Radunit$.

	To obtain \eqref{E:CLASSICALLIFESPANASYMPTOTICESTIMATE},
	we use
	\eqref{E:MUSTARBOUNDS}
	and
	\eqref{E:CRUCIALLATETIMEDERIVATIVECOMPAREDTODATAPARAMETER}
	to deduce that $\upmu_{\star}(t,1)$
	vanishes for the first time when
	$t = \TranminusdatasizeWithFactor^{-1} 
	+ \mathcal{O}(\mathring{\upepsilon})
	$.

	We now derive the statements regarding the quantities that extend to 
	$\Sigma_{T_{(Lifespan);U_0}}^{U_0}$ as $L^{\infty}$ functions.
	Let $q$ denote any of the quantities
	$\Fullset^{\leq 11;\leq 2} \Psi$,
	$\cdots$,
	$\Rad \Rad \upmu$
	that, in the theorem, are stated 
	to extend to $\Sigma_{T_{(Lifespan);U_0}}^{U_0}$ as an $L^{\infty}$ function
	of the geometric coordinates.
	The $L^{\infty}$ estimates of 
	Props.~\ref{P:IMPROVEMENTOFAUX}
	and \ref{P:IMPROVEMENTOFHIGHERTRANSVERSALBOOTSTRAP}
	imply that $\| \Lunit q \|_{L^{\infty}(\Sigma_t^{U_0})}$
	is uniformly bounded
	for $0 \leq t < T_{(Lifespan);U_0}$.
	Recalling that $\Lunit = \frac{\partial}{\partial t}$, we
	conclude that $q$ extends to $\Sigma_{T_{(Lifespan);U_0}}^{U_0}$
	as an element of $L^{\infty}(\Sigma_{T_{(Lifespan);U_0}}^{U_0})$
	as desired.
	The estimate 
	$g_{\alpha \beta}(\threePsi) = m_{\alpha \beta} + \mathcal{O}(\mathring{\upepsilon})$
	and the extension properties of the $\Fullset$-derivatives of
	the scalar-valued functions $g_{\alpha \beta}(\threePsi)$ then follow from
	\eqref{E:LITTLEGDECOMPOSED},
	the already proven bound $\| \threePsi \|_{L^{\infty}(\Sigma_t^{U_0})} \lesssim \mathring{\upepsilon}$,
	and the extension properties of the $\Fullset$-derivatives of $\threePsi$ obtained just above.
	This completes the proof of the theorem.
\end{proof}

\bibliographystyle{amsalpha}
\bibliography{JBib}

\end{document}